\documentclass{smfbook}
\usepackage{sga1-smf}

\setboolean{orig}{false}

\begin{document}
\frontmatter
\title[Rev\^etements \'etales et groupe fondamental (SGA 1)]{%
S\'eminaire de G\'eom\'etrie Alg\'ebrique du~Bois~Marie
\\
1960--61
\\
Rev\^etements \'etales et groupe fondamental
\\(SGA 1)
}
\alttitle{Rev\^etements \'etales et groupe fondamental (SGA 1)}

\author{Un s\'eminaire dirig\'e par A.\ Grothendieck\\
\smaller Augment\'e de deux expos\'es de Mme M.\ Raynaud}

\begin{abstract}
Ce volume est une \'edition recompos\'ee et annot\'ee du livre \og Rev\^etements \'Etales et Groupe Fondamental\fg, Lecture Notes in Mathematics, 224, Springer-Verlag, Berlin-Heidelberg-New York, 1971, par Alexander Grothendieck et al.\\
Le texte pr\'esente les fondements d'une th\'eorie du groupe fondamental en G\'eom\'etrie Alg\'ebrique, dans le point de vue \og kroneckerien\fg permettant de traiter sur le m\^eme pied le cas d'une vari\'et\'e alg\'ebrique au sens habituel, et celui d'un anneau des entiers d'un corps de nombres, par exemple.
\end{abstract}

\begin{altabstract}
This volume is a new updated edition of the book ``Rev\^etements \'Etales et Groupe Fondamental'', Lecture Notes in Mathematics, 224, Springer-Verlag, Berlin-Heidelberg-New York, 1971, by Alexander Grothendieck et al.\\
The text presents the foundations of a theory of the fundamental group in Algebraic Geometry from the Kronecker point of view, allowing one to treat on an equal footing the case of an algebraic variety in the usual sense, and that of the ring of integers in a number field, for instance.
\end{altabstract}



\altkeywords{\'Etale morphism, smooth morphism, flat morphism, scheme, fundamental group, covering, theory of descent, specialization}

\maketitle

\cleardoublepage
\renewcommand{\baselinestretch}{1.05}\normalfont
\makeschapterhead{Pr\'eface}
\thispagestyle{empty}
\vfill
Le pr\'esent texte est une \'edition recompos\'ee
\ifthenelse{\boolean{orig}}{\ignorespaces}{et annot\'ee}
du livre \og Rev\^etements
\'Etales et Groupe Fondamental\fg, Lecture Notes in Mathematics, 224,
Springer-Verlag, Berlin-Heidelberg-New York, 1971, par Alexander
Grothendieck et al.

La composition en \LaTeX 2e a \'et\'e r\'ealis\'ee par des volontaires, participant \`a un projet dirig\'e par Bas Edixhoven, dont on trouvera plus de d\'etails \`a l'adresse \url{http://www.math.leidenuniv.nl/~edix/}. La mise en page a \'et\'e termin\'ee par la \hbox{Soci\'et\'e} math\'ematique de France.

\ifthenelse{\boolean{orig}}
{Cette version est cens\'ee reproduire le \emph{texte
original}. Dans une autre version, des erreurs typographiques ont
\'et\'e corrig\'ees et quelques mises \`a jour ont \'et\'e effectu\'ees par Michel Raynaud. Cette derni\`ere est \'edit\'ee par la Soci\'et\'e math\'ematique de France (\url{http://smf.emath.fr/}) dans la s\'erie \og Documents Math\'ematiques\fg.}
{Ceci est une \emph{version l\'eg\`erement corrig\'ee} du texte
original. Elle est \'edit\'ee par la \hbox{Soci\'et\'e} math\'ematique de France (\url{http://smf.emath.fr/}) dans la s\'erie \og Documents \hbox{Math\'ematiques}\fg. Quelques mises \`a jour ont \'et\'e effectu\'ees par Michel Raynaud. Celles-ci sont d\'elimit\'ees par des crochets \lcrochetbf\,\rcrochetbf\ et indiqu\'ees par le symbole (MR) (remarques pages \pageref{X.2.14}, \pageref{XI.1.4}, \pageref{XII.5.6}, \pageref{XIII.2.13} et note \Ref{III.6.6.p24} page \pageref{III.6.6.p24}). Afin d'uniformiser les notations, le corps r\'esiduel d'un point  $x$  est not\'e  $\kres(x)$,
et celui d'un anneau local  $A$  est not\'e  $\kres(A)$.\par
Il existe \'egalement une version \'electronique cens\'ee reproduire le texte original.}

Les deux versions ont un seul fichier source commun, se trouvant sur
le serveur des archives arXiv.org e-Print \`a
\url{http://arxiv.org/}. Les diff\'erences entre les deux versions
sont document\'ees dans le fichier source.

L'ancienne num\'erotation des pages est incorpor\'ee dans la marge
(le nombre~$n$ indiquant le d\'ebut de la page~$n$).

\vfill
\newpage
The present text is a new
\ifthenelse{\boolean{orig}}{\ignorespaces}{updated}
edition of the book ``Rev\^etements \'Etales et
Groupe Fondamental'', Lecture Notes in Mathematics, 224,
Springer-Verlag, Berlin-Heidelberg-New York, 1971, by Alexander
Grothendieck et al.

Typesetting in \LaTeX 2e was done by volunteers, participating in a project directed by Bas Edixhoven, on which one can find more details at
\url{http://www.math.leidenuniv.nl/~edix/}. It has been completed by the Soci\'et\'e \hbox{math\'ematique} de France.

\ifthenelse{\boolean{orig}}
{This version is meant to reproduce the \emph{original text}. There is
another version, in which some typograpical mistakes have been
corrected and some notes have been added by Michel Raynaud. The latter has been published by the Soci\'et\'e math\'ematique de France (\url{http://smf.emath.fr/}) as a volume of the series ``Documents Math\'ematiques''.}
{This is a slightly \emph{corrected version} of the original
text. It is published by the \hbox{Soci\'et\'e} \hbox{math\'ematique} de France (\url{http://smf.emath.fr/}) as a volume of the series ``\hbox{Documents} Math\'ematiques''. Some updating remarks have been added by Michel Raynaud. These are bounded by brackets \lcrochetbf\,\rcrochetbf\ and indicated by the symbol (MR) (Remarks on pages \pageref{X.2.14}, \pageref{XI.1.4}, \pageref{XII.5.6}, \pageref{XIII.2.13} and footnote \Ref{III.6.6.p24} on page \pageref{III.6.6.p24}). In order to unify notation, the residue field of a point $x$ is denoted by $\kres(x)$ and the residue field of a local ring $A$ is denoted by $\kres(A)$.\par
There also exists an electronic version that is meant to reproduce the original text.}

Both versions are produced from the same source file, that can be
downloaded from the arXiv.org e-Print server at
\url{http://arxiv.org/}. All differences between the two versions are
documented in the source file.

The old page numbering is incorporated in the margin (the number~$n$
marking the beginning of page~$n$).

\vfill

\cleardoublepage
\renewcommand{\baselinestretch}{1.1}\normalfont
\makeschapterhead{Introduction}
\thispagestyle{empty}
\label{I.0}

Dans la premi\`ere partie de cette introduction, nous donnons des
pr\'ecisions sur le contenu du pr\'esent volume; dans la
deuxi\`eme, sur l'ensemble du \og \emph{S\'eminaire de
G\'eom\'etrie Alg\'ebrique du Bois-Marie}\fg, dont le pr\'esent
volume constitue le tome premier.

\let\oldthesubsection\thesubsection
\def\thesubsection{\arabic{subsection}}
\subsection{}
\label{I.introduction.1}
Le pr\'esent volume pr\'esente les fondements d'une
th\'eorie du groupe fondamental en G\'eom\'etrie Alg\'ebrique,
dans le point de vue \og kroneckerien\fg permettant de traiter sur le
m\^eme pied le cas d'une vari\'et\'e alg\'ebrique au sens
habituel, et celui d'un anneau des entiers d'un corps de nombres, par
exemple. Ce point de vue ne s'exprime d'une fa\c con satisfaisante
que dans le langage des sch\'emas, et nous utiliserons librement ce
langage, ainsi que les r\'esultats principaux expos\'es dans les
trois premiers chapitres des \emph{\'El\'ements de G\'eom\'etrie
Alg\'ebrique} de J\ptbl \textsc{Dieudonn\'e} et A\ptbl \textsc{Grothendieck}, (cit\'e EGA dans
la suite). L'\'etude du pr\'esent volume du \og \emph{S\'eminaire
de G\'eom\'etrie Alg\'ebrique du Bois-Marie}\fg ne demande pas
d'autres connaissances de la G\'eom\'etrie Alg\'ebrique, et peut
donc servir d'introduction aux techniques actuelles de
G\'eom\'etrie Alg\'ebrique, \`a un lecteur d\'esireux de se
familiariser avec ces techniques.

Les expos\'es \Ref{I} \`a \Ref{XI} de ce livre sont une reproduction
textuelle, pratiquement inchang\'ee, des notes
mim\'eographi\'ees du S\'eminaire oral, qui \'etaient
distribu\'ees par les soins de l'\emph{Institut des Hautes \'Etudes
Scientifiques\footnote{Ainsi que les notes des s\'eminaires faisant
suite \`a celui-ci. Ce mode de distribution s'\'etant
av\'er\'e impraticable et insuffisant \`a la longue, tous les
\og \emph{S\'eminaire de G\'eom\'etrie Alg\'ebrique du
Bois-Marie}\fg para\^itront d\'esormais sous forme de livre comme
le pr\'esent volume.}}. Nous nous sommes born\'es \`a rajouter
quelques notes de bas de page au texte primitif, de corriger quelques
erreurs de frappe, et de faire un ajustage terminologique, le mot
\og morphisme simple\fg ayant notamment \'et\'e remplac\'e
entre-temps par celui de \og morphisme lisse\fg, qui ne pr\^ete pas aux
m\^emes confusions.

Les expos\'es \Ref{I} \`a \Ref{IV} pr\'esentent les notions locales de
morphisme \emph{\'etale} et de morphisme \emph{lisse}; ils
n'utilisent gu\`ere le langage des sch\'emas, expos\'e dans le
Chapitre~I des \emph{\'El\'ements}\footnote{Une \'etude plus
compl\`ete est maintenant disponible dans les \emph{\'El\'ements},
Chap~IV, \oldS\S 17 et 18.}. L'expos\'e~\Ref{V} pr\'esente la description
axiomatique du groupe fondamental d'un sch\'ema, utile m\^eme dans
le cas classique o\`u ce sch\'ema se r\'eduit au spectre d'un
corps, o\`u on trouve une reformulation fort commode de la
th\'eorie de Galois habituelle. Les expos\'es \Ref{VI} et \Ref{VIII}
pr\'esentent la \emph{th\'eorie de la descente}, qui a pris une
importance croissance en G\'eom\'etrie Alg\'ebrique dans ces
derni\`eres ann\'ees, et qui pourrait rendre des services
analogues en G\'eom\'etrie Analytique et en Topologie. Il convient
de signaler que l'expos\'e VII n'avait pas \'et\'e
r\'edig\'e, et sa substance se trouve
\ifthenelse{\boolean{orig}}
{incorpor\'e}
{incorpor\'ee}
dans un
travail de J\ptbl Giraud (M\'ethode de la Descente, Bull.\ Soc.\ Math.\
France, M\'emoire~2, 1964, viii + 150~p.). Dans l'expos\'e \Ref{IX}, on
\'etudie plus sp\'ecifiquement la descente des morphismes
\'etales, obtenant une approche syst\'ematique pour des
th\'eor\`emes du type de \textsc{Van Kampen} pour le groupe fondamental,
qui apparaissent ici comme de simples traductions de th\'eor\`emes
de descente. Il s'agit essentiellement d'un proc\'ed\'e de calcul
du groupe fondamental d'un sch\'ema connexe~$X$, muni d'un morphisme
surjectif et propre, disons $X'\to X$, en termes des groupes
fondamentaux des composantes connexes de~$X'$ et des produits
fibr\'es~$X'\times_X X'$, $X'\times_X X'\times_X X'$, et des
homomorphismes induits entre ces groupes par les morphismes
simpliciaux canoniques entre les sch\'emas
pr\'ec\'edents. L'expos\'e~\Ref{X} donne la th\'eorie de la
\emph{sp\'ecialisation du groupe fondamental}, pour un morphisme
propre et lisse, dont le r\'esultat le plus frappant consiste en la
d\'etermination (\`a peu de chose pr\`es) du groupe fondamental
d'une courbe alg\'ebrique lisse en caract\'eristique $p>0$,
gr\^ace au r\'esultat connu par voie transcendante en
caract\'eristique nulle. L'expos\'e~\Ref{XI} donne quelques exemples et
compl\'ements, en explicitant notamment sous forme cohomologique la
th\'eorie des rev\^etements de \emph{\textsc{Kummer}}, et celle
d'\emph{\textsc{Artin}-\textsc{Schreier}}. Pour d'autres commentaires sur le texte,
voir l'\emph{Avertissement} \`a la version multigraphi\'ee, qui
fait suite \`a la pr\'esente Introduction.

Depuis la r\'edaction en 1961 du pr\'esent S\'eminaire a
\'et\'e d\'evelopp\'e, en collaboration par M\ptbl \textsc{Artin} et
moi-m\^eme, le langage de la \emph{topologie \'etale} et une
th\'eorie cohomologique correspondante, expos\'ee dans la partie
SGA~4 \og Cohomologie \'etale des sch\'emas\fg du \emph{S\'eminaire
de G\'eom\'etrie Alg\'ebrique}, \`a para\^itre dans la
m\^eme s\'erie que le pr\'esent volume. Ce langage, et les
r\'esultats dont il dispose d\`es \`a pr\'esent, fournissent
un outil particuli\`erement souple pour l'\'etude du groupe
fondamental, permettant de mieux comprendre et de d\'epasser
certains des r\'esultats expos\'es ici. Il y aurait donc lieu de
reprendre enti\`erement la th\'eorie du groupe fondamental de ce
point de vue (tous les r\'esultats-clefs figurant en fait d\`es
\`a pr\'esent dans \loccit). C'est ce qui \'etait projet\'e
pour le chapitre des \emph{\'El\'ements} consacr\'e au groupe
fondamental, qui devait contenir \'egalement plusieurs autres
d\'eveloppements qui n'ont pu trouver leur place ici (s'appuyant sur
la technique de r\'esolution des singularit\'es): calcul du
\og groupe fondamental local\fg d'un anneau local complet en termes d'une
r\'esolution des singularit\'es convenable de cet anneau, formules
de K\"unneth locales et globales pour le groupe fondamental sans
hypoth\`ese de propret\'e (\cf Exp\ptbl \Ref{XIII}), les r\'esultats de
M\ptbl \textsc{Artin} sur la comparaison des groupes fondamentaux locaux d'un
anneau local hens\'elien excellent et de son compl\'et\'e
(SGA~4~XIX). Signalons \'egalement la n\'ecessit\'e de
d\'evelopper une th\'eorie du groupe fondamental d'un topos, qui
englobera \`a la fois la th\'eorie topologique ordinaire, sa
version semi-simpliciale, la variante \og profinie\fg
\ifthenelse{\boolean{orig}}
{d\'evelopp\'e}
{d\'evelopp\'ee}
dans l'expos\'e~\Ref{V} du pr\'esent volume, et la variante
pro-discr\`ete un peu plus g\'en\'erale de SGA~3~X~7
(adapt\'ee au cas de sch\'emas non normaux et non unibranches). En
attendant une refonte d'ensemble de la th\'eorie dans cette optique,
l'expos\'e \Ref{XIII} de Mme \textsc{Raynaud}, utilisant le langage et les
r\'esultats de SGA~4, est destin\'e \`a montrer le parti qu'on
peut tirer dans quelques questions typiques, en g\'en\'eralisant
notamment certains r\'esultats de l'expos\'e~\Ref{X} \`a des
sch\'emas relatifs non propres. On y donne en particulier la
structure du groupe fondamental \og premier \`a $p$\fg d'une courbe
alg\'ebrique non compl\`ete en
\ifthenelse{\boolean{orig}}
{car.}
{caract\'eristique}
quelconque (que j'avais
\ifthenelse{\boolean{orig}}
{annonc\'e}
{annonc\'ee}
en 1959, mais dont une d\'emonstration n'avait pas
\'et\'e publi\'ee \`a ce jour).

Malgr\'e ces nombreuses lacunes et imperfections (d'autres diront:
\`a cause de ces lacunes et imperfections), je pense que le
pr\'esent volume pourra \^etre utile au lecteur qui d\'esire se
familiariser avec la th\'eorie du groupe fondamental, ainsi que
comme ouvrage de r\'ef\'erence, en attendant la r\'edaction et
la parution d'un texte \'echappant aux critiques que je viens
d'\'enum\'erer.

\subsection{}
\label{I.introduction.2}
Le pr\'esent volume constitue le tome 1 du
\og \emph{S\'eminaire de G\'eom\'etrie Alg\'ebrique du
Bois-Marie}\fg, dont les volumes suivants sont pr\'evus pour
para\^itre dans la m\^eme s\'erie que celui-ci. Le but que se
propose le \emph{S\'eminaire}, parall\`element au trait\'e
\og \emph{\'El\'ements de G\'eom\'etrie Alg\'ebrique}\fg de J\ptbl \textsc{Dieudonn\'e} et
A\ptbl \textsc{Grothendieck}, est de jeter les fondements de la G\'eom\'etrie
Alg\'ebrique, suivant les points de vue dans ce dernier ouvrage. La
r\'ef\'erence standard pour tous les volumes du
\emph{S\'eminaire} est constitu\'ee par les Chapitres~I, II, III
des \og \emph{\'El\'ements de G\'eom\'etrie Alg\'ebrique}\fg
(cit\'es EGA~I, II, III), et on suppose le lecteur en possession du
bagage d'alg\`ebre commutative et l'alg\`ebre homologique que ces
chapitres impliquent\footnote{Voir l'Introduction \`a EGA~I pour des
pr\'ecisions \`a ce sujet.}. De plus, dans chaque volume du
\emph{S\'eminaire} il sera r\'ef\'er\'e librement, dans le
mesure des besoins, \`a des volumes ant\'erieurs du m\^eme
\emph{S\'eminaire}, ou \`a d'autres chapitres publi\'es ou sur
le point de para\^itre des \og \'El\'ements\fg.

Chaque \emph{partie} du \emph{S\'eminaire} est centr\'ee sur un
sujet principal, indiqu\'e dans le titre du ou des volumes
correspondants; le s\'eminaire oral porte g\'en\'eralement sur
une ann\'ee acad\'emique, parfois plus. Les expos\'es \`a
l'int\'erieur de chaque partie du \emph{S\'eminaire} sont
g\'en\'eralement dans une d\'ependance logique \'etroite les
uns par rapport aux autres; par contre, les diff\'erentes parties du
\emph{S\'eminaire} sont dans une large mesure logiquement
ind\'ependants les uns par rapport aux autres. Ainsi, la partie
\og Sch\'emas en Groupes\fg est \`a peu pr\`es enti\`erement
ind\'ependante des deux parties du \emph{S\'eminaire} qui la
pr\'ec\`edent chronologiquement; par contre, elle fait un
fr\'equent appel aux r\'esultats de EGA~IV. Voici la liste des
parties du \emph{S\'eminaire} qui doivent para\^itre
prochainement (cit\'es SGA~1 \`a SGA~7 dans la suite):
\begin{enumerateb}
\item[SGA 1.] Rev\^etements \'etales et groupe fondamental, 1960
et 1961.
\item[SGA 2.] Cohomologie locale des faisceaux coh\'erents et
th\'eor\`emes de Lefschetz locaux et globaux, 1961/62.
\item[SGA 3.] Sch\'emas en groupes, 1963 et 1964 (3~volumes),
en~coll. avec M\ptbl \textsc{Demazure}.
\item[SGA 4.] Th\'eorie des topos et cohomologie \'etale des
sch\'emas, 1963/64 (3~volumes), (en~coll. avec M\ptbl \textsc{Artin} et
J.L\ptbl \textsc{Verdier}).
\item[SGA 5.] Cohomologie $l$-adique et fonctions~$L$, 1964 et 1965
(2~volumes).
\item[SGA 6.] Th\'eorie des intersections et th\'eor\`eme de
Riemann-Roch, 1966/67 (2~volumes) (en~coll. avec P\ptbl \textsc{Berthelot} et
L\ptbl \textsc{Illusie}).
\item[SGA 7.] Groupes de monodromie locale en g\'eom\'etrie
alg\'ebrique.
\end{enumerateb}

Trois parmi ces \emph{S\'eminaires} partiels ont \'et\'e
dirig\'es en \emph{collaboration} avec d'autres
math\'e\-maticiens, qui figureront comme co-auteurs sur la
couverture des volumes correspondants. Quant aux autres participants
actifs du \emph{S\'eminaire}, dont le r\^ole (tant au point de vue
r\'edactionnel que de celui du travail de mise au point
math\'ematique) est all\'e croissant d'ann\'ee en ann\'ee, le
nom de chaque participant figure en t\^ete des expos\'es dont il
est responsable comme conf\'erencier ou comme r\'edacteur, et la
liste de ceux qui figurent dans un volume d\'etermin\'e se trouve
\ifthenelse{\boolean{orig}}
{indiqu\'e}
{indiqu\'ee}
sur la page de garde dudit volume.

Il convient de donner quelques pr\'ecisions sur les rapports entre
le \emph{S\'eminaire} et les \emph{\'El\'ements}. Ces derniers
\'etaient destin\'es en principe \`a donner un expos\'e
d'ensemble des notions et techniques jug\'ees les plus fondamentales
dans la G\'eom\'etrie Alg\'ebrique, \`a mesure que ces notions
et techniques elles-m\^emes se d\'egagent, par le jeu naturel
d'exigences de coh\'erence logique et esth\'etique. Dans cette
optique, il \'etait naturel de consid\'erer le
\emph{S\'eminaire} comme une version pr\'eliminaire des
\emph{\'El\'ements}, destin\'ee \`a \^etre englob\'ee \`a
peu pr\`es totalement, t\^ot ou tard, dans ces derniers. Ce
processus avait d\'ej\`a commenc\'e dans une certaine mesure il
y a quelques ann\'ees, puisque les expos\'es I \`a IV du
pr\'esent volume SGA~1 sont enti\`erement englob\'es par EGA~IV,
et que les expos\'es VI \`a VIII devaient l'\^etre d'ici
quelques ann\'ees dans EGA~VI. Cependant, \`a mesure que se
d\'eveloppe le travail d'\'edification entrepris dans les
\emph{\'El\'ements} et dans le \emph{S\'eminaire}, et que les
proportions d'ensemble se pr\'ecisent, le principe initial
(d'apr\`es lequel le \emph{S\'eminaire} ne constituerait qu'une
version pr\'eliminaire et provisoire) appara\^it de moins en
moins r\'ealiste en raison (entre autres) des limites impos\'ees
par la pr\'evoyante nature \`a la dur\'ee de la vie
humaine. Compte tenu du soin g\'en\'eralement apport\'e dans la
r\'edaction des diff\'erentes parties du \emph{S\'eminaire}, il
n'y aura lieu sans doute de reprendre une telle partie dans les
\emph{\'El\'ements} (ou des trait\'es qui en prendraient la
rel\`eve) que lorsque des progr\`es ult\'erieurs \`a la
r\'edaction permettront d'y apporter des am\'eliorations tr\`es
substantielles, aux prix de modifications assez profondes. C'est le
cas d\`es \`a pr\'esent pour le pr\'esent s\'eminaire SGA~1,
comme on l'a dit plus haut, et \'egalement pour SGA~2 (gr\^ace aux
r\'esultats r\'ecents de Mme~\textsc{Raynaud}). Par contre, rien n'indique
actuellement qu'il en sera ainsi dans un proche avenir pour aucune des
parties cit\'ees plus haut SGA~3 \`a SGA~7.

Les r\'ef\'erences \`a l'int\'erieur du \og S\'eminaire de
G\'eom\'etrie Alg\'ebrique du
\ifthenelse{\boolean{orig}}
{Bois-Marie}
{Bois Marie}\fg sont donn\'ees
ainsi. Une \emph{r\'ef\'erence int\'erieure} \`a une des
parties SGA~1 \`a SGA~7 du S\'eminaire est donn\'ee dans le
style III~9.7, o\`u le chiffre III d\'esigne le num\'ero de
l'expos\'e, qui figure en haut de chaque page de l'expos\'e en
question, et~9.7 le num\'ero de l'\'enonc\'e (ou de la
d\'efinition, remarque, etc.) \`a l'int\'erieur de
l'expos\'e. Le cas \'ech\'eant, des nombres d\'ecimaux plus
longs peuvent \^etre utilis\'es, par exemple 9.7.1, 9.7.2 pour
d\'esigner par exemple les diverses \'etapes dans la
d\'emonstration d'une proposition~9.7. La r\'ef\'erence III~9
d\'esigne la paragraphe 9 de l'expos\'e~III. Le num\'ero de
l'expos\'e est omis pour les r\'ef\'erences internes \`a un
expos\'e. Pour une \emph{r\'ef\'erence \`a une autre} des
\emph{parties} du \emph{S\'eminaire}, on utilise les m\^emes
sigles, mais pr\'ec\'ed\'es de la mention de la partie en
question des SGA, SGA~1~III~9.7. De m\^eme, la r\'ef\'erence
EGA~IV~11.5.7 signifie: \'El\'ements de G\'eom\'etrie
Alg\'ebrique, Chap\ptbl IV, \'enonc\'e (ou d\'efinition
\ifthenelse{\boolean{orig}}{etc...}{etc.})~11.5.7; ici, le premier chiffre arabe d\'esigne encore le
num\'ero du paragraphe. \`A part ces conventions en vigueur dans
l'ensemble des SGA, la bibliographie relative \`a un expos\'e sera
g\'en\'eralement rassembl\'ee \`a la fin de celui-ci, et il y
sera r\'ef\'er\'e \`a l'int\'erieur de l'expos\'e par des
num\'eros entre crochets comme
\ifthenelse{\boolean{orig}}{[3]}{[\textbf{3}]}, suivant l'usage.

Enfin, pour la commodit\'e du lecteur, chaque fois que cela semblera
n\'ecessaire, nous joindrons \`a la fin des volumes des SGA un
index des notations, et un index terminologique contenant s'il y a
lieu une traduction anglaise des termes fran\c cais utilis\'es.

Je tiens \`a joindre \`a cette introduction un commentaire
extra-math\'ematique. Au mois de novembre 1969 j'ai eu connaissance
du fait que l'Institut des Hautes \'Etudes Scientifiques, dont j'ai
\'et\'e professeur essentiellement depuis sa fondation, recevait
depuis trois ans des subventions du Minist\`ere des
Arm\'ees. D\'ej\`a comme chercheur d\'ebutant j'ai trouv\'e
extr\^emement regrettable le peu de scrupules que se font la plupart
des scientifiques pour accepter de collaborer sous une forme ou une
autre avec les appareils militaires. Mes motivations \`a ce moment
\'etaient essentiellement de nature morale, donc peu susceptibles
d'\^etres prises au s\'erieux. Aujourd'hui elles acqui\`erent
une force et une dimension nouvelle, vu le danger de destruction de
l'esp\`ece humaine dont nous menace la prolif\'eration des
appareils militaires et des moyens de destruction massives dont ces
appareils disposent. Je me suis expliqu\'e ailleurs de fa\c con
plus d\'etaill\'ee sur ces questions, beaucoup plus importantes
que l'avancement de n'importe quelle science (y compris la
math\'ematique); on pourra par exemple consulter \`a ce sujet
l'article de G\ptbl Edwards dans le \no 1 du journal Survivre (Ao\^ut
1970), r\'esumant un expos\'e plus d\'etaill\'e de ces
questions que j'avais fait ailleurs. Ainsi, je me suis trouv\'e
travailler pendant trois ans dans une institution alors qu'elle participait \`a
mon insu \`a un mode de financement que je consid\`ere comme
immoral et dangereux\footnote{Il va de soi que l'opinion que je viens
d'exprimer n'engage que ma propre responsabilit\'e, et non pas celle
de la maison d'\'edition Springer qui \'edite le pr\'esent
volume.}. \'Etant \`a pr\'esent seul \`a avoir cette opinion parmi
mes coll\`egues \`a l'IHES, ce qui a condamn\'e \`a
l'\'echec mes efforts pour obtenir la suppression des subventions
militaires du budget de l'IHES, j'ai pris la d\'ecision qui
s'imposait et quitte l'IHES le 30 septembre 1970 et suspends
\'egalement toute collaboration scientifique avec cette institution,
aussi longtemps qu'elle continuera \`a accepter de telles
subventions. J'ai demand\'e \`a M\ptbl Motchane, directeur de l'IHES,
que l'IHES s'abstienne \`a partir du 1er octobre 1970 de diffuser des textes
math\'ematiques dont je suis auteur, ou faisant partie du
S\'eminaire de G\'eom\'etrie Alg\'ebrique du Bois Marie. Comme
il a \'et\'e dit plus haut, la diffusion de ce s\'eminaire va
\^etre assur\'ee par la maison Julius Springer, dans le s\'erie
des Lecture Notes. Je suis heureux de remercier ici la maison Springer
et Monsieur K\ptbl Peters pour l'aide efficace et courtoise qu'ils m'ont
apport\'ee pour rendre possible cette publication, en se chargeant
en particulier de la frappe pour la photooffset des nouveaux
expos\'es rajout\'es aux anciens s\'eminaires, et des
expos\'es manquants des s\'eminaires incomplets.

Je remercie \'egalement M\ptbl J.P\ptbl Delale, qui s'est charg\'e du
travail ingrat de compiler l'index des notations et l'index
terminologique.

\bigskip\hfill Massy, ao\^ut 1970.

\cleardoublepage
\makeschapterhead{Avertissement}
\thispagestyle{empty}
\label{I.avertissement}

Chacun des expos\'es r\'edig\'es donne la substance de plusieurs
expos\'es oraux cons\'ecutifs. Il n'a pas sembl\'e utile d'en
pr\'eciser les dates.

L'expos\'e~VII, auquel il est r\'ef\'er\'e \`a diverses
reprises au cours de l'expos\'e~\Ref{VIII}, n'a pas \'et\'e
r\'edig\'e par le conf\'erencier, qui dans les conf\'erences
orales s'\'etait born\'e \`a esquisser le langage de la descente
dans les cat\'egories g\'en\'erales, en se pla\c cant \`a un
point de vue strictement utilitaire et sans entrer dans les
difficult\'es logiques soulev\'ees par ce langage. Il est apparu
qu'un expos\'e correct de ce langage sortirait des limites des
pr\'esentes notes, ne serait-ce que par sa longueur. Pour un
expos\'e en forme de la th\'eorie de la descente, je renvoie \`a
un article en pr\'eparation de Jean \textsc{Giraud}. En attendant sa
parution\footnote{\label{footnotegiraud}Actuellement paru: J\ptbl \textsc{Giraud},
\emph{M\'ethodes de la descente}, M\'emoire \no 2 de la
Soci\'et\'e math\'ematiques de France, 1964.}, je pense qu'un
lecteur attentif n'aura pas de peine \`a suppl\'eer par ses
propres moyens aux r\'ef\'erences fant\^omes de
l'Expos\'e~VIII.

D'autres expos\'es oraux, se pla\c cant apr\`es l'Expos\'e~\Ref{XI},
et auxquels il est fait allusion \`a certains endroits du texte,
n'ont pas non plus \'et\'e r\'edig\'es, et \'etaient
destin\'es \`a former la substance d'un Expos\'e~\Ref{XII} et d'un
Expos\'e~\Ref{XIII}. Les premiers de ces expos\'es oraux reprenaient,
dans le cadre des sch\'emas et des espaces analytiques avec
\'el\'ements nilpotents (tels qu'ils sont introduits dans le
S\'eminaire Cartan 1960/61) la construction de l'espace analytique associ\'e
\`a un pr\'esch\'ema localement de type fini sur un corps
valu\'e complet~$k$, les th\'eor\`emes du type GAGA dans le cas
o\`u $k$ est le corps des complexes, et l'application \`a la
comparaison du groupe fondamental d\'efini par voie transcendante et
le groupe fondamental \'etudi\'e dans ces notes (comparer
A\ptbl Grothendieck, Fondements de la G\'eom\'etrie Alg\'ebrique,
S\'eminaire Bourbaki \no 190, page~10, d\'ecembre 1959). Les
derniers expos\'es oraux esquissaient la g\'en\'eralisation des
m\'ethodes d\'evelopp\'ees dans le texte pour l'\'etude des
rev\^etements admettant de la ramification mod\'er\'ee, et de la
structure du groupe fondamental d'une courbe compl\`ete priv\'ee
d'un nombre fini de points (comparer \loccit, \no 182, page 27,
th\'eor\`eme~14). Ces expos\'es n'introduisent aucune id\'ee
essentiellement nouvelle, c'est pourquoi il n'a pas sembl\'e
indispensable d'en donner une r\'edaction en forme avant la parution
des chapitres correspondants des \'El\'ements de G\'eom\'etrie
Alg\'ebrique\footnote{Ils sont inclus dans le pr\'esent volume
dans l'Exp\ptbl \Ref{XII} de Mme~Raynaud avec une d\'emonstration
diff\'erente de la d\'emonstration originale expos\'ee dans le
S\'eminaire oral (\cf introduction).}.

Par contre, les th\'eor\`emes du type Lefschetz pour le groupe
fondamental et le groupe de Picard, tant au point de vue local que
global, ont fait l'objet d'un S\'eminaire s\'epar\'e en 1962,
qui a \'et\'e compl\`etement r\'edig\'e et est \`a la
disposition des usagers\footnote{\emph{Cohomologie \'etale des
faisceaux coh\'erents et th\'eor\`emes de Lefschetz locaux et
globaux} (SGA~2), paru dans North Holland Pub.\ Cie.}. Signalons que
les r\'esultats d\'evelopp\'es tant dans le pr\'esent
S\'eminaire que dans celui de 1962 seront utilis\'es de fa\c con
essentielle dans la parution de plusieurs r\'esultats clefs dans la
cohomologie \'etale des pr\'esch\'emas, qui feront l'objet d'un
S\'eminaire (conduit par M\ptbl Artin et moi-m\^eme) en 1963/64,
actuellement en pr\'eparation\footnote{\emph{Cohomologie \'etale
des sch\'emas} (cit\'e SGA~4), \`a para\^itre dans cette
m\^eme s\'erie.}.

Les expos\'es \Ref{I} \`a~\Ref{IV}, de nature essentiellement locale et
tr\`es \'el\'ementaire, seront absorb\'es enti\`erement par
le chapitre IV des \emph{\'El\'ements de G\'eom\'etrie
Alg\'ebrique}, dont la premi\`ere partie est \`a l'impression et
qui sera sans doute publi\'e vers fin~64. Ils pourront n\'eanmoins
\^etre utiles \`a un lecteur qui d\'esirerait se mettre au
courant des propri\'et\'es essentielles des morphismes lisses,
\'etales ou plats, avant d'entrer dans les arcanes d'un trait\'e
syst\'ematique. Quant aux autres expos\'es, ils seront
absorb\'es dans le chapitre VIII\footnote{En fait, par suite de
modification du plan initialement pr\'evu pour les
\emph{\'El\'ements}, l'\'etude du groupe fondamental y est
report\'ee \`a un chapitre ult\'erieur \`a celui qu'on vient
d'indiquer. Comparer l'introduction qui pr\'ec\`ede le pr\'esent
\og Avertissement\fg.} des \og \'El\'ements\fg, dont la publication ne
pourra gu\`ere \^etre envisag\'ee avant plusieurs ann\'ees.

\bigskip\hfill Bures, juin 1963.

{\def\footnotemark{}\let\\\relax
\chapterspace{-4}
\tableofcontents}


\mainmatter
\renewcommand{\chaptername}{Expos\'e}
\let\thesubsection\oldthesubsection

\chapterspace{-2}
\chapter{Morphismes \'etales}
\label{I}

\marginpar{1}

Pour simplifier l'exposition, on suppose que tous les
pr\'esch\'emas envisag\'es sont localement noeth\'eriens, du
moins apr\`es le num\'ero~\Ref{I.2}.

\section{Notions de calcul diff\'erentiel}
\label{I.1}
Soit $X$ un pr\'esch\'ema sur~$Y$, soit $\Delta_{X/Y}$ ou $\Delta$
\label{indnot:ab}\oldindexnot{$\Delta_{X/Y}$ ou simplement $\Delta$|hyperpage}%
le morphisme diagonal $X\to X\times_Y X$. C'est un morphisme
d'immersion, donc un morphisme d'immersion ferm\'ee de~$X$ dans un
ouvert $V$ de~$X\times_Y X$. Soit $\cal{I}_X$ l'id\'eal du
sous-pr\'esch\'ema ferm\'e correspondant \`a la diagonale
dans~$V$ (N.B. si on veut faire les choses intrins\`equement, sans
supposer $X$ s\'epar\'e sur~$Y$ --- hypoth\`ese qui serait
canularesque --- on devrait consid\'erer l'image inverse ensembliste
de~$\cal{O}_{X\times X}$ dans~$X$, et d\'esigner par $\cal{I}_X$
l'id\'eal d'augmentation dans ce dernier...). Le faisceau
$\cal{I}_X/\cal{I}_X^2$ peut \^etre regard\'e comme un faisceau
quasi-coh\'erent sur~$X$, on le d\'enote
par~$\mathit{\Omega}^1_{X/Y}$.
\label{indnot:ac}\oldindexnot{$\mathit{\Omega}^1_{X/Y}$|hyperpage}%
Il est de type fini si $X\to Y$ est de type fini. Il se comporte bien
par rapport \`a extension de la base $Y'\to Y$. On introduit aussi
les faisceaux $\cal{O}_{X\times_Y X}/\cal{I}_X^{n+1}=\cal{P}^n_{X/Y}$,
\label{indnot:ad}\oldindexnot{$\cal{P}_{X/Y}^n$|hyperpage}%
ce sont des faisceaux d'\emph{anneaux} sur~$X$, faisant de~$X$ un
pr\'esch\'ema qu'on peut d\'enoter par $\Delta_{X/Y}^n$
\label{indnot:ae}\oldindexnot{$\Delta_{X/Y}^n$|hyperpage}%
et appeler le \emph{\nieme voisinage infinit\'esimal de}~$X/Y$.
\index{voisinage infinit\'esimal de~$X/Y$|hyperpage}%
Le sorite en est d'une trivialit\'e totale, bien qu'il soit assez
long\footnote{\cf EGA IV~16.3.}; il serait prudent de n'en parler
qu'au moment o\`u on en dit quelque chose de serviable, avec les
morphismes lisses.

\section{Morphismes quasi-finis}
\label{I.2}
\begin{proposition}
\label{I.2.1}
Soit $A\to B$ un homomorphisme local (N.B. les anneaux sont maintenant
noeth\'eriens), \ifthenelse{\boolean{orig}}{$m$}{$\goth{m}$} l'id\'eal maximal de~$A$. Conditions
\'equivalentes:
\begin{enumerate}
\item[(i)] \ifthenelse{\boolean{orig}}
{$B/mB$ est de dimension finie sur~$k=A/m$}
{$B/\goth{m}B$ est de dimension finie sur~$k=A/\goth{m}$}.
\item[(ii)] \ifthenelse{\boolean{orig}}
{$mB$ est un id\'eal de d\'efinition, et $B/r(B)=\kappa(B)$ est une extension
de~$k=\kappa(A)$}
{$\goth{m}B$ est un id\'eal de d\'efinition, et $B/\goth{r}(B)=\kres(B)$ est une extension de~$k=\kres(A)$}
\item[(iii)] Le compl\'et\'e $\widehat{B}$ est fini sur celui
$\widehat{A}$ de~$A$.
\end{enumerate}
\end{proposition}
On dit alors
\marginpar{2}
que $B$ \emph{est quasi-fini}
\index{quasi-fini (morphisme, alg\`ebre)|hyperpage}%
sur~$A$. Un morphisme $f\colon X\to Y$
est dit quasi-fini en~$x$ (ou le $Y$-pr\'esch\'ema $f$ est dit
quasi-fini en~$x$) si $\cal{O}_x$ est quasi-fini
sur~$\cal{O}_{f(x)}$. Cela revient aussi \`a dire que $x$ est
\emph{isol\'e dans sa fibre}~$f^{-1}(x)$. Un morphisme est dit
quasi-fini s'il l'est en tout point\footnote{Dans EGA~II~6.2.3 on
suppose de plus $f$ de type fini.}.
\begin{corollaire}
\label{I.2.2}
Si $A$ est complet, quasi-fini \'equivaut \`a fini.
\ifthenelse{\boolean{orig}}
{On pourrait donner le sorite \textup{(i) (ii) (iii) (iv) (v)} habituel pour les morphismes quasi-finis, mais ce ne semble pas indispensable ici.}
{}
\end{corollaire}

\ifthenelse{\boolean{orig}}
{}
{On pourrait donner le sorite \textup{(i), (ii), (iii), (iv), (v),} habituel pour les morphismes quasi-finis, mais ce ne semble pas indispensable ici.}

\section{Morphismes non ramifi\'es ou nets}
\label{I.3}
\begin{proposition}
\label{I.3.1}
Soit $f\colon X\to Y$ un morphisme de type fini, $x\in X$,
$y=f(x)$. Conditions \'equivalentes:
\begin{enumerate}
\item[(i)] \ifthenelse{\boolean{orig}}
{$\cal{O}_x/m_y\cal{O}_x$ est une extension finie
s\'eparable de~$\kappa(y)$}
{$\cal{O}_x/\goth{m}_y\cal{O}_x$ est une extension finie
s\'eparable de~$\kres(y)$}.
\item[(ii)] $\mathit{\Omega}^1_{X/Y}$ est nul en~$x$.
\item[(iii)] Le morphisme diagonal $\Delta_{X/Y}$ est une immersion
ouverte au voisinage de~$x$.
\end{enumerate}
\end{proposition}
Pour l'implication (i)$\To$(ii), on est ramen\'e
aussit\^ot par Nakayama au cas o\`u $Y=\Spec(k)$,
$X=\Spec(k')$, o\`u c'est bien connu et d'ailleurs trivial
sur la d\'efinition de s\'eparable; (ii)$\To$(iii)
d'apr\`es une caract\'erisation agr\'eable et facile des
immersions ouvertes, utilisant Krull; (iii)$\To$(i) car on est
encore ramen\'e au cas o\`u $Y=\Spec(k)$ et o\`u le
morphisme diagonal est une immersion ouverte partout. Il faut alors
prouver que $X$ est fini d'anneau s\'eparable sur~$k$, on est
ramen\'e pour ceci au cas o\`u $k$ est alg\'ebriquement
clos. Mais alors tout point ferm\'e de~$X$ est isol\'e (car
identique \`a l'image inverse de la diagonale par le morphisme $X\to
X\times_k X$ d\'efini par~$x$), d'o\`u le fait que $X$ est
fini. On peut supposer alors $X$ r\'eduit \`a un point,
d'anneau~$A$, donc $A\otimes_k A\to A$ est un isomorphisme, d'o\`u
$A=k$ cqfd.

\begin{definition}
\label{I.3.2}
\begin{enumerate}
\item[a)] On dit alors que $f$ est \emph{net},
\index{net (morphisme, alg\`ebre)|hyperpage}%
ou encore \emph{non ramifi\'e},
\index{non ramifi\'e (morphisme, alg\`ebre)|hyperpage}%
en~$x$, ou que $X$ est net, ou encore non ramifi\'e, en $x$ sur~$Y$.
\item[b)] Soit $A\to B$ un homomorphisme local, on dit qu'il est
\emph{net}, ou \emph{non ramifi\'e}, ou que $B$ est une alg\`ebre
locale \emph{nette}, ou \emph{non ramifi\'ee} sur~$A$, si 
\ifthenelse{\boolean{orig}}
{$B/r(A)B$ est une extension finie s\'eparable de~$A/r(A)$ \ie si $r(A)B=r(B)$}
{$B/\goth{r}(A)B$ est une extension finie s\'eparable de~$A/\goth{r}(A)$, \ie si $\goth{r}(A)B=\goth{r}(B)$}
et $\kres(B)$ est une extension s\'eparable
de~$\kres(A)$\kern1pt\footnote{\Cf remords dans III~\Ref{III.1.2}.}.
\end{enumerate}
\end{definition}

\begin{remarquesstar}
Le fait
\marginpar{3}
que $B$ soit net sur~$A$ se reconna\^it d\'ej\`a sur
les compl\'et\'es de~$A$ et de~$B$. Net implique quasi-fini.
\end{remarquesstar}

\begin{corollaire}
\label{I.3.3}
L'ensemble des points o\`u $f$ est net est ouvert.
\end{corollaire}

\begin{corollaire}
\label{I.3.4}
Soient $X'$, $X$ deux pr\'esch\'emas de type fini sur~$Y$, et
$g\colon X'\to X$ un $Y$-morphisme. Si $X$ est net sur~$Y$, le
morphisme graphe $\Gamma_g\colon X'\to X\times_Y X$ est une immersion
ouverte.
\end{corollaire}

En effet, c'est l'image inverse du morphisme diagonal $X\to X\times_Y
X$ par
$$
g\times_Y\id_{X'}\colon X'\times_Y X\to X\times_Y X.
$$

On peut aussi introduire l'id\'eal annulateur $\goth{d}_{X/Y}$
\label{indnot:af}\oldindexnot{$\goth{d}_{X/Y}$|hyperpage}%
de~$\mathit{\Omega}^1_{X/Y}$, appel\'e id\'eal diff\'erente
\index{diff\'erente (id\'eal)|hyperpage}%
\index{ideal differente@id\'eal diff\'erente|hyperpage}%
de~$X/Y$; il d\'efinit un sous-pr\'esch\'ema ferm\'e de~$X$
qui, ensemblistement, est l'ensemble des points o\`u $X/Y$ est
ramifi\'e, \ie non net.

\begin{proposition}
\label{I.3.5}
\begin{enumerate}
\item[(i)] Une immersion est nette.
\item[(ii)] Le compos\'e de deux morphismes nets l'est.
\item[(iii)] Extension de base dans un morphisme net en est un autre.
\end{enumerate}
\end{proposition}

Se voit indiff\'eremment sur (ii) ou (iii) (le deuxi\`eme me
semble plus amusant). On peut bien entendu aussi pr\'eciser, en
donnant des \'enonc\'es ponctuels; ce n'est plus g\'en\'eral
qu'en apparence (sauf dans le cadre de la d\'efinition~b)), et
barbant. On obtient comme d'habitude des corollaires:

\begin{corollaires}
\label{I.3.6}
\begin{enumerate}
\item[(iv)] Produit cart\'esien de deux morphismes nets en est un
autre.
\item[(v)] Si $gf$ est net, $f$ est net.
\item[(vi)] Si $f$ est net, $f_{\text{\textup{r\'ed}}}$ est net.
\end{enumerate}
\end{corollaires}

\begin{proposition}
\label{I.3.7}
Soit $A\to B$ un homomorphisme local, on suppose l'extension
r\'esiduelle $\kres(B)/\kres(A)$ triviale ou $\kres(A)$ alg\'ebriquement
clos. Pour que $B/A$ soit net, il faut et il suffit que $\widehat{B}$
soit (comme $\widehat{A}$-alg\`ebre) un quotient de~$\widehat{A}$.
\end{proposition}

\begin{remarquesstar}
\begin{enumerate}
\item[--] Dans le cas o\`u on ne suppose pas l'extension
r\'esiduelle triviale, on peut se ramener \`a ce cas en faisant
une extension finie plate convenable sur~$A$ qui d\'etruise ladite
extension.
\item[--] Donner l'exemple o\`u $A$ est l'anneau local d'un point
double ordinaire d'une courbe, $B$ d'un point du normalis\'e: alors
$A\subset B$, $B$ est net sur
\marginpar{4}
$A$ \`a extension r\'esiduelle triviale, et
$\widehat{A}\to\widehat{B}$ est surjectif mais \emph{non injectif}. On
va donc renforcer la notion de nettet\'e.
\end{enumerate}
\end{remarquesstar}

\section{Morphismes \'etales. Rev\^etements \'etales}
\label{I.4}
On va admettre tout ce qui nous sera n\'ecessaire sur les morphismes
plats; ces faits seront d\'emontr\'es ult\'erieurement, s'il y a
lieu\footnote{\Cf Exp\ptbl \Ref{IV}.}.

\begin{definition}
\label{I.4.1}
\begin{enumerate}
\item[a)] Soit $f\colon X\to Y$ un morphisme de type fini. On dit que
$f$ est \emph{\'etale}
\index{etale (morphisme, algebre)@\'etale (morphisme, alg\`ebre)|hyperpage}%
en $x$ si $f$ est plat en $x$ et net en~$x$. On dit que $f$ est
\'etale s'il l'est en tous les points. On dit que $X$ est \'etale
en $x$ sur~$Y$, ou que c'est un $Y$-pr\'esch\'ema \'etale en~$x$
\ifthenelse{\boolean{orig}}{etc...}{etc.}
\item[b)] Soit $f\colon A\to B$ un homomorphisme local. On dit que $f$
est \'etale, ou que $B$ est \'etale sur~$A$, si $B$ est plat et
non ramifi\'e sur~$A$.\footnote{\Cf remords dans III~\Ref{III.1.2}.}
\end{enumerate}
\end{definition}

\begin{proposition}
\label{I.4.2}
Pour que $B/A$ soit \'etale, il faut et il suffit que
$\widehat{B}/\widehat{A}$ le soit.
\end{proposition}

En
\marginpar{5}
effet, c'est vrai s\'epar\'ement pour \og net\fg et pour \og plat\fg.

\begin{corollaire}
\label{I.4.3}
Soit $f\colon X\to Y$ de type fini, et $x\in X$. Le fait que $f$ soit
\'etale en $x$ ne d\'epend que de l'homomorphisme local
$\cal{O}_{f(x)}\to \cal{O}_x$, et m\^eme seulement de
l'homomorphisme correspondant pour les compl\'et\'es.
\end{corollaire}

\begin{corollaire}
\label{I.4.4}
Supposons que l'extension r\'esiduelle $\kres(A)\to \kres(B)$ soit triviale,
ou que $\kres(A)$ soit alg\'ebriquement clos. Alors $B/A$ est \'etale
\sss $\widehat{A}\to\widehat{B}$ est un isomorphisme.
\end{corollaire}

On conjugue la platitude et~\Ref{I.3.7}.
\begin{proposition}
\label{I.4.5}
Soit $f\colon X\to Y$ un morphisme de type fini. Alors l'ensemble des
points o\`u il est \'etale est ouvert.
\end{proposition}

En effet, c'est vrai s\'epar\'ement pour \og net\fg et \og plat\fg.

Cette proposition montre qu'on peut laisser tomber les \'enonc\'es
\og ponctuels\fg dans l'\'etude des morphismes de type fini \'etales
quelque part.
\begin{proposition}
\label{I.4.6}
\begin{enumerate}
\item[(i)] Une immersion ouverte est \'etale.
\item[(ii)] Le compos\'e de deux morphismes \'etales est
\'etale.
\item[(iii)] Extension de la base.
\end{enumerate}
\end{proposition}

En effet, (i) est trivial, et pour (ii) et (iii) il suffit de noter
que c'est vrai pour \og net\fg et pour \og plat\fg. \`A vrai dire, il y a
aussi des \'enonc\'es correspondants pour les homomorphismes
locaux (sans condition de finitude), qui en tout \'etat de cause
devront figurer au multiplodoque (\`a commencer par le cas: net).

\begin{corollaire}
\label{I.4.7}
Un produit cart\'esien de deux morphismes \'etales est itou.
\end{corollaire}

\begin{corollaire}
\label{I.4.8}
Soient $X$ et $X'$ de type fini sur~$Y$, $g\colon X\to X'$ un
$Y$-morphisme. Si $X'$ est non ramifi\'e sur~$Y$ et $X$ \'etale
sur~$Y$, alors $g$ est \'etale.
\end{corollaire}

En effet, $g$ est le compos\'e du morphisme graphe $\Gamma_g\colon
X\to X\times_Y X'$ qui est une immersion ouverte par~\Ref{I.3.4}, et
du morphisme de projection qui est \'etale car d\'eduit du
morphisme \'etale $X\to Y$ par le \og changement de base\fg $X'\to Y$.

\begin{definition}
\label{I.4.9}
On appelle rev\^etement \'etale (\resp net)
\index{etale (revetement)@\'etale (rev\^etement)|hyperpage}%
\index{net (revetement)@net (rev\^etement)|hyperpage}%
\index{revetement etale (\resp net)@rev\^etement \'etale (\resp net)|hyperpage}%
de~$Y$ un $Y$-sch\'ema $X$ qui est fini sur~$Y$ et \'etale
(\resp net) sur~$Y$.
\end{definition}

\ifthenelse{\boolean{orig}}{}
{\enlargethispage{.5cm}}%
La premi\`ere condition signifie que $X$ est d\'efini par un
faisceau coh\'erent d'alg\`ebres $\cal{B}$ sur~$Y$. La
deuxi\`eme signifie alors que $\cal{B}$ est localement libre sur~$Y$
(\resp rien du tout), \emph{et} que de plus, pour tout $y\in Y$, la
fibre $\cal{B}(y)=\cal{B}_y\otimes_{\cal{O}_y}\kres(y)$ soit une
alg\`ebre s\'eparable (= compos\'e fini d'extensions finies
s\'eparables) sur~$\kres(y)$.

\begin{proposition}
\label{I.4.10}
Soit $X$ un rev\^etement plat de~$Y$ de degr\'e $n$ (la
d\'efinition de ce terme m\'eritait de figurer dans~\Ref{I.4.9})
d\'efini par un faisceau coh\'erent localement libre $\cal{B}$
d'alg\`ebres. On d\'efinit de fa\c con bien connue
l'homomorphisme trace $\cal{B}\to\cal{A}$ (qui est un homomorphisme de
$\cal{A}$-modules, o\`u $\cal{A}=\cal{O}_Y$). Pour que $X$ soit
\'etale, il \fets que la forme bilin\'eaire
$\trace_{\cal{B}/\cal{A}}xy$ correspondante d\'efinisse un
isomorphisme de~$\cal{B}$ sur~$\cal{B}$, ou ce qui revient au
m\^eme, que la \emph{section discriminant}
$$
d_{X/Y} = d_{\cal{B}/\cal{A}} \in
\Gamma(Y,\tbigwedge^n\check{\cal{B}}\otimes_{\cal{A}}
\tbigwedge^n\check{\cal{B}})
$$
soit inversible, ou enfin que l'id\'eal discriminant d\'efini par
cette section soit l'id\'eal unit\'e.
\end{proposition}
En effet, on est ramen\'e au cas o\`u $Y=\Spec(k)$, et alors c'est
un crit\`ere de s\'eparabilit\'e bien connue, et trivial par
passage \`a la cl\^oture alg\'ebrique de~$k$.

\begin{remarquestar}
On
\marginpar{6}
aura un \'enonc\'e moins trivial plus bas, quand on ne suppose
pas a priori que $X$ est plat sur~$Y$ mais qu'on fait une
hypoth\`ese de normalit\'e.
\end{remarquestar}

\section{La propri\'et\'e fondamentale des morphismes \'etales}
\label{I.5}

\begin{theoreme}
\label{I.5.1}
Soit $f\colon X \to Y$ un morphisme de type fini. Pour que $f$ soit
une immersion ouverte, il \fets que ce soit un
morphisme \emph{\'etale et radiciel}.
\end{theoreme}

Rappelons que radiciel
\index{radiciel (morphisme)|hyperpage}%
signifie: injectif, \`a extensions
r\'esiduelles radicielles (et
\ifthenelse{\boolean{orig}}{en}{on}
peut aussi rappeler que cela signifie que le morphisme reste injectif
par toute extension de la base). La n\'ecessit\'e est triviale,
reste la suffisance. On va donner deux d\'emonstrations
diff\'erentes, la premi\`ere plus courte, la deuxi\`eme plus
\'el\'ementaire.

1) Un morphisme plat est ouvert, donc on peut supposer (rempla\c cant $Y$ par~$f(X)$) que $f$ est un \emph{hom\'eomorphisme}
\ifthenelse{\boolean{orig}}
{s\^ur}
{sur}.
Par toute extension de base, il restera vrai que
\ifthenelse{\boolean{orig}}{$f^\prime$}{$f$}
est plat, radiciel, surjectif, donc un hom\'eomorphisme, a fortiori
\ifthenelse{\boolean{orig}}{ferm\'ee}{ferm\'e}
Donc $f$ est \emph{propre}. Donc $f$ est \emph{fini}
(r\'ef\'erence: th\'eor\`eme de Chevalley) d\'efini par un
faisceau coh\'erent $\mathcal{B}$ d'alg\`ebres. $\mathcal{B}$ est
localement libre, de plus en vertu de l'hypoth\`ese il est partout
de rang~$1$, donc $X=Y$, cqfd.

2) On peut supposer $Y$ et $X$ \emph{affines}. On se ram\`ene de
plus facilement \`a prouver ceci: si $Y = \Spec(A)$, $A$ local, et
si $f^{-1}(y)$ est non vide ($y$ \'etant le point ferm\'e de~$Y$)
alors $X=Y$ (en effet, cela impliquera que tout $y \in f(X)$ a un
voisinage ouvert $U$ tel que $X|U = U$). On aura $X = \Spec(B)$, on
veut prouver $A=B$. Mais pour ceci, on est ramen\'e \`a prouver
l'assertion analogue en rempla\c cant $A$ par~$\hat A$, $B$ par $B
\otimes_A \hat A$ (compte tenu que $\hat A$ est fid\`element plat
sur~$A$). On peut donc supposer $A$ \emph{complet}. Soit $x$ le point
au-dessus de~$y$, d'apr\`es le corollaire~\Ref{I.2.2}
$\mathcal{O}_x$ est fini sur~$A$ donc (\'etant plat et radiciel
sur~$A$) est identique \`a~$A$. Donc on a $X = Y \amalg X^\prime$
(somme disjointe). Comme $X$ est radiciel sur~$Y$, $X^\prime$ est
vide. On a fini.
\ifthenelse{\boolean{orig}}{}
{\enlargethispage{.5cm}}%
\begin{corollaire}
\label{I.5.2}
Soit $f\colon X \to Y$ un morphisme d'\emph{immersion ferm\'ee} et
\emph{\'etale}. Si $X$ est connexe, $f$ est un isomorphisme de~$X$
sur une composante connexe de~$Y$.
\end{corollaire}

En effet, $f$ est aussi une immersion ouverte. On en d\'eduit:

\begin{corollaire}
\label{I.5.3}
Soit
\marginpar{7}
$X$ un $Y$-sch\'ema net, $Y$ connexe. Alors toute section
de~$X$ sur~$Y$ est un isomorphisme de~$Y$ sur une composante connexe
de~$X$. Il y a donc correspondance biunivoque entre l'ensemble de ces
sections, et l'ensemble des composantes connexes $X_i$ de~$X$ telles
que la projection $X_i \to Y$ soit un isomorphisme, (ou, ce qui
revient au m\^eme par~\Ref{I.5.1}, surjectif et radiciel). En
particulier, une section est connue quand en conna\^it sa valeur en
un point.
\end{corollaire}

Seule la premi\`ere assertion demande une d\'emonstration;
d'apr\`es~\Ref{I.5.2} il suffit de remarquer qu'une section est une
immersion ferm\'ee (car $X$ est s\'epar\'e sur~$Y$) et \'etale
en vertu de~\Ref{I.4.8}.

\begin{corollaire}
\label{I.5.4}
Soient $X$ et $Y$ deux pr\'esch\'emas sur~$S$, $X$ net
s\'epar\'e sur~$S$ et $Y$ connexe. Soient $f$, $g$ deux
$S$-morphismes de~$Y$ dans~$X$, $y$ un point de~$Y$, on suppose $f(y)
= g(y) = x$ et les homomorphismes r\'esiduels $\kres(x) \to \kres(y)$
d\'efinis par $f$ et~$g$ identiques (\og $f$ et~$g$ co\"incident
g\'eom\'etriquement en~$y$\fg). Alors $f$ et~$g$ sont identiques.
\end{corollaire}

R\'esulte de~\Ref{I.5.3} en se ramenant au cas o\`u $Y=S$, en
rempla\c cant $X$ par $X \times_S Y$.

Voici une variante particuli\`erement importante de~\Ref{I.5.3}.

\begin{theoreme}
\label{I.5.5}
Soient $S$ un pr\'esch\'ema, $X$ et $Y$ deux
$S$-pr\'esch\'emas, $S_0$ un sous-pr\'esch\'ema ferm\'e
de~$S$ ayant m\^eme espace sous-jacent que~$S$, $X_0 = X \times_S
S_0$ et $Y_0 = Y \times_S S_0$ les \og restrictions\fg de~$X$ et~$Y$
sur~$S_0$. On suppose $X$ \'etale sur~$S$. Alors l'application
naturelle
$$
\Hom_S(Y,X) \to \Hom_{S_0}(Y_0,X_0)
$$
est \emph{bijective}.
\end{theoreme}

On est encore ramen\'e au cas o\`u $Y=S$, et alors cela
r\'esulte de la description \og topologique\fg des sections de~$X/Y$
donn\'ee dans~\Ref{I.5.3}.

\begin{scholiestar}
Ce r\'esultat comporte une assertion d'\emph{unicit\'e} et
d'\emph{existence} de \emph{morphismes}. Il peut aussi s'exprimer
(lorsque $X$ et~$Y$ sont tous deux pris \'etales sur~$S$) que le
foncteur $X \mto X_0$ de la cat\'egorie des $S$-sch\'emas
\'etales dans la cat\'egorie des $S_0$-sch\'emas \'etales est
\emph{pleinement fid\`ele}, \ie \'etablit une
\emph{\'equivalence} de la premi\`ere avec une
\emph{sous-cat\'egorie pleine} de la seconde. Nous verrons plus bas
que c'est m\^eme une \'equivalence de la premi\`ere et de la
seconde (ce qui sera un th\'eor\`eme d'\emph{existence de
$S$-sch\'emas \'etales}).
\end{scholiestar}

La
\marginpar{8}
forme suivante, plus g\'en\'erale en apparence, de~\Ref{I.5.5}.\
est souvent commode:

\begin{corollaire}[\og Th\'eor\`eme de prolongement des rel\`evements\fg]
\label{I.5.6}
\index{prolongement des rel\`evements (th\'eor\`eme de)|hyperpage}%
\index{theoreme de prolongement des relevements@th\'eor\`eme de prolongement des rel\`evements|hyperpage}%
Consid\'erons un diagramme commutatif
$$
\xymatrix{
X\ar[d]&\ar[l]Y_0\ar[d]\\
S&\ar[l]Y
}
$$
de morphismes, o\`u $X \to S$ est \'etale et $Y_0 \to Y$ est une
immersion ferm\'ee bijective. Alors on peut trouver un morphisme
unique $Y \to X$ qui rende les deux triangles correspondants
commutatifs.
\end{corollaire}

En effet, rempla\c cant $S$ par~$Y$ et $X$ par $X \times_S Y$, on
est ramen\'e au cas o\`u $Y=S$, et alors c'est le cas particulier
de~\Ref{I.5.5} pour $Y=S$.

Signalons aussi la cons\'equence imm\'ediate suivante
de~\Ref{I.5.1} (que nous n'avons pas
\ifthenelse{\boolean{orig}}{donn\'e}{donn\'ee}
en corollaire~1 pour ne pas interrompre la ligne d'id\'ees
d\'evelopp\'ee \`a la suite de~\Ref{I.5.1}):

\begin{proposition}
\label{I.5.7}
Soient $X$, $X^\prime$ deux pr\'esch\'emas de type finis et plats
sur~$Y$, et soit $g\colon X \to X^\prime$ un $Y$-morphisme. Pour que
$g$ soit une immersion ouverte (\resp un isomorphisme) il faut et il
suffit que pour tout $y\in Y$, le morphisme induit sur les fibres
$$
g \otimes_Y \kres(y)\colon X \otimes_Y \kres(y) \to X^\prime
\otimes_Y \kres(y)
$$
le soit.
\end{proposition}

Il suffit de prouver la suffisance; comme c'est vrai pour la notion de
surjection, on est ramen\'e au cas d'une immersion
ouverte. D'apr\`es~\Ref{I.5.1}, il faut v\'erifier que $g$ est
\emph{radiciel}, ce qui est trivial, et qu'il est \emph{\'etale}, ce
qui r\'esulte du corollaire~\Ref{I.5.9} ci-dessous.

\begin{corollaire}
\label{I.5.8}
(devrait passer au \No\Ref{I.3}) Soient $X$ et $X^\prime$ deux
$Y$-pr\'esch\'emas, $g\colon X \to X^\prime$ un $Y$-morphisme, $x$
un point de~$X$ et $y$ sa projection sur~$Y$. Pour que $g$ soit
quasi-fini (\resp net) en~$x$, il \fets qu'il en soit
de m\^eme de~$g\otimes_Y \kres(y)$.
\end{corollaire}

En effet, les deux alg\`ebres sur~$k\big(g(x)\big)$ qu'il faut
regarder pour s'assurer que l'on a bien un morphisme quasi-fini \resp net en~$x$ sont les m\^emes pour~$g$ et $g \otimes_Y \kres(y)$.

\begin{corollaire}
\label{I.5.9}
Avec
\marginpar{9}
les notations de~\Ref{I.5.8}, supposons $X$ et~$X^\prime$ plats
et de type fini sur~$Y$. Pour que $g$ soit plat (\resp \'etale)
en~$x$, il \fets que $g \otimes_Y \kres(y)$ le soit.
\end{corollaire}

Pour \og plat\fg l'\'enonc\'e n'est mis que pour m\'emoire, c'est
un des crit\`eres fondamentaux de platitude\footnote{\Cf IV~\Ref{IV.5.9}.}. Pour \'etale, cela en r\'esulte; compte tenu
de~\Ref{I.5.8}.

\section{Application aux extensions \'etales des anneaux locaux complets}
\label{I.6}

Ce num\'ero est un cas particulier de r\'esultats sur les
pr\'esch\'emas formels, qui devront figurer dans le
multiplodoque. N\'eanmoins, on s'en tire ici \`a meilleur compte,
\ie sans la d\'etermination locale explicite des morphismes
\'etales au \No\Ref{I.7} (utilisant le Main Theorem). C'est peut-\^etre
une raison suffisante de garder le pr\'esent num\'ero (m\^eme
dans le multiplodoque) \`a cette place.

\begin{theoreme}
\label{I.6.1}
Soit $A$ un anneau local complet (noeth\'erien bien s\^ur), de
corps r\'esiduel~$k$. Pour toute $A$-alg\`ebre~$B$, soit $R(B) = B
\otimes_A k$ consid\'er\'e comme $k$-alg\`ebre, elle d\'epend
donc fonctoriellement de~$B$. Alors $R$ d\'efinit une
\emph{\'equivalence} de la cat\'egorie des $A$-alg\`ebres
\emph{finies et \'etales sur~$A$} avec la cat\'egorie des
alg\`ebres \emph{de rang fini s\'eparables} sur~$k$.
\end{theoreme}

Tout d'abord, le foncteur en question est pleinement fid\`ele, comme
il r\'esulte du fait plus g\'en\'eral:

\begin{corollaire}
\label{I.6.2}
Soient $B$, $B^\prime$ deux $A$-alg\`ebres finies sur~$A$. Si $B$
est \'etale sur~$A$, alors l'application canonique
$$
\Hom_{\textup{$A$-alg}}(B,B^\prime) \to
\Hom_{\textup{$k$-alg}}\big(R(B),R(B^\prime)\big)
$$
est bijective.
\end{corollaire}

On est ramen\'e au cas o\`u $A$ est artinien (en rempla\c cant
$A$ par~$A/\goth{m}^n$), et alors c'est un cas particulier
de~\Ref{I.5.5}.

Il reste \`a prouver que pour toute $k$-alg\`ebre finie et
s\'eparable (pourquoi ne pas dire: \'etale, c'est plus court) $L$,
il existe un $B$ \'etale sur~$A$ tel que $R(B)$ soit isomorphe
\`a~$L$. On peut supposer que $L$ est une extension s\'eparable
de~$k$, comme telle elle admet un g\'en\'erateur~$x$, \ie est
isomorphe \`a une alg\`ebre $k[t]/Fk[t]$ o\`u $F \in k[t]$ est
un polyn\^ome unitaire. On rel\`eve~$F$ en un polyn\^ome
unitaire $F_1$ dans $A[t]$, et on prend $B = A[t]/F_1A[t]$.

\section{Construction locale des morphismes non ramifi\'es et \'etales}
\label{I.7}
\marginpar{10}

\begin{proposition}
\label{I.7.1}
Soient $A$ un anneau noeth\'erien, $B$ une alg\`ebre finie
sur~$A$, $u$ un g\'en\'erateur de~$B$ sur~$A$, $F \in A[t]$ tel
que $F(u) = 0$ (on ne suppose pas $F$ unitaire), $u^\prime =
F^\prime(u)$ (o\`u $F^\prime$ est le polyn\^ome d\'eriv\'e),
$\goth{q}$ un id\'eal premier de~$B$ ne contenant pas~$u^\prime$,
$\goth{p}$ sa trace sur~$A$. Alors $B_\goth{q}$ est net sur~$A_\goth{p}$.
\end{proposition}

En d'autres termes, posant $Y = \Spec(A)$, $X= \Spec(B)$,
$X_{u^\prime} = \Spec(B_{u^\prime})$, $X_{u^\prime}$ est non
ramifi\'e sur~$Y$. L'\'enonce r\'esulte du suivant, plus
pr\'ecis:

\begin{corollaire}
\label{I.7.2}
L'id\'eal diff\'erente de~$B/A$ contient $u^\prime B$, et lui est
\'egal si l'homo\-morphisme naturel $A[t]/FA[t] \to B$ (appliquant
$t$ dans~$u$) est un isomorphisme.
\end{corollaire}

Soit $J$ le noyau de l'homomorphisme $C = A[t] \to B$, ce noyau
contient $FA[t]$, et lui est \'egal dans le deuxi\`eme cas
envisag\'e dans~\Ref{I.7.2}. Comme il est surjectif,
$\Omega^1_{B/A}$ s'identifie au quotient de~$\Omega^1_{C/A}$ par le
sous-module engendr\'e par $J\Omega^1_{C/A}$ et $d(J)$ (il aurait
fallu expliciter au \No\Ref{I.1} la d\'efinition de l'homomorphisme~$d$, et
le calcul de~$\Omega^1$ pour une alg\`ebre de polyn\^omes).
Identifiant $\Omega^1_{C/A}$ \`a~$C$ gr\^ace \`a la base~$dt$,
on trouve $B/B\cdot J^\prime$ donc la diff\'erente est engendr\'ee
par l'ensemble~$J^\prime$ des images dans~$B$ des d\'eriv\'es des
$G \in J$, (et il suffit de prendre des~$G$ engendrant~$J$). Comme
$F\in J$, \resp $F$ est un g\'en\'erateur de~$J$, on a
fini. (N.B.\ On devrait mettre~\Ref{I.7.2} en prop.\
et~\Ref{I.7.1} en corollaire). On trouve:

\begin{corollaire}
\label{I.7.3}
Sous les conditions de~\Ref{I.7.1}, supposant $F$ unitaire et que
$A[t]/FA[t] \to B$ est un isomorphisme, pour que $B_\goth{q}$ soit
\'etale sur~$A_\goth{p}$, il
\ifthenelse{\boolean{orig}}{f et.\ s}{faut et il suffit}
que~$\goth{q}$ ne contienne
pas~$u^\prime$.
\end{corollaire}

En effet, comme $B$ est plat sur~$A$, \'etale \'equivaut \`a
net, et on peut appliquer~\Ref{I.7.2}.

\begin{corollaire}
\label{I.7.4}
Sous les conditions de~\Ref{I.7.3} pour que $B$ soit \'etale
sur~$A$ il
\ifthenelse{\boolean{orig}}{f et.\ s.\ }{faut et il suffit }%
que~$u^\prime$ soit inversible, ou encore que l'id\'eal engendr\'e
par $F$, $F^\prime$ dans $A[t]$ soit l'id\'eal unit\'e.
\end{corollaire}

Le dernier crit\`ere r\'esulte du premier et de Nakayama
(dans~$B$).

Un polyn\^ome unitaire $F\in A[t]$ ayant la propri\'et\'e
\'enonc\'ee dans le corollaire~\Ref{I.7.4}
\marginpar{11}
est dit \emph{polyn\^ome s\'eparable}
\index{separable (polynome)@s\'eparable (polyn\^ome)|hyperpage}%
(si $F$ n'est pas unitaire, il faudrait au moins exiger que le
coefficient de son terme dominant soit inversible; dans le cas o\`u
$A$ est un corps, on retrouve la d\'efinition usuelle).


\begin{corollaire}
\label{I.7.5}
Soit $B$ une alg\`ebre \emph{finie} sur l'anneau
\emph{local}~$A$. On suppose que $K(A)$ est infini ou que $B$ soit
local. Soit~$n$ le rang de~$L = B \otimes_A K(A)$ sur~$K(A) = k$. Pour
que $B$ soit net (\resp \'etale) sur~$A$, il faut et il suffit que
$B$ soit isomorphe \`a un quotient de (\resp isomorphe \`a)
$A[t]/FA[t]$, o\`u $F$ est un polyn\^ome unitaire s\'eparable,
qu'on peut supposer (\resp qui est n\'ecessairement) de
degr\'e~$n$.
\end{corollaire}

Il n'y a qu'\`a prouver la n\'ecessit\'e. Supposons $B$ net
sur~$A$, donc $L$ s\'eparable sur~$k$, il r\'esulte alors de
l'hypoth\`ese faite que $L/k$ admet un g\'en\'erateur~$\xi$,
donc les $\xi^i$ ($0 \leq i < n$) forment une base de~$L$
sur~$k$. Soit $u\in B$ relevant~$\xi$, alors par Nakayama les $u^i$
($0\leq i < n$) engendrent le $A$-module~$B$ (\resp en forment une
base), en particulier on peut trouver un polyn\^ome unitaire $F \in
A[t]$ tel que $F(u)= 0$, et $B$ sera isomorphe \`a un quotient de
(\resp isomorphe \`a) $A[t]/FA[t]$. Enfin, en vertu
de~\Ref{I.7.4}.\ appliqu\'e \`a $L/k$, $F$ et~$F^\prime$
engendrent $A[t]$
\ifthenelse{\boolean{orig}}{module}{modulo}
$\goth{m} A[t]$, donc (d'apr\`es Nakayama dans $A[t]/FA[t]$) $F$
et~$F^\prime$ engendrent~$A[t]$, on a fini.

\begin{theoreme}
\label{I.7.6}
Soient $A$ un anneau local, $A \to \cal{O}$ un homomorphisme local tel
que~$\cal{O}$ soit isomorphe \`a une alg\`ebre localis\'ee d'une
alg\`ebre de type fini sur~$A$. Supposons
\ifthenelse{\boolean{orig}}{$O$}{$\cal{O}$}
\emph{net} sur~$A$. Alors on peut trouver une $A$-alg\`ebre $B$,
enti\`ere sur~$A$, un id\'eal maximal $\goth{n}$ de~$B$, un
g\'en\'erateur $u$ de~$B$ sur~$A$, un polyn\^ome unitaire $F\in
A[t]$, tels que $\goth{n} \not\ni F^\prime(u)$ et que $\cal{O}$ soit
isomorphe (comme $A$-alg\`ebre) \`a~$B_{\goth{n}}$. Si $\cal{O}$
est \'etale sur~$A$, on peut prendre $B=A[t]/FA[t]$.
\end{theoreme}

(Bien entendu, on a l\`a des conditions aussi suffisantes...)

Signalons d'abord les agr\'eables corollaires:

\begin{corollaire}
\label{I.7.7}
Pour que $\cal{O}$ soit net sur~$A$, il \fets que $\cal{O}$ soit
isomorphe au quotient d'une alg\`ebre analogue et \emph{\'etale}
sur~$A$.
\end{corollaire}

En effet, on prendra $\cal{O}^\prime = B^\prime_{\goth{n}^\prime}$,
o\`u $B^\prime = A[t]/FA[t]$ et o\`u $\goth{n}^\prime$ est l'image
inverse de~$\goth{n}$ dans~$B^\prime$.

\begin{corollaire}
\label{I.7.8}
Soit
\marginpar{12}
$f\colon X \to Y$ un morphisme de type fini, $x\in X$. Pour que
$f$ soit net en~$x$, il \fets qu'il existe un voisinage ouvert
$U$ de~$x$ tel que $f|U$ se factorise en $U \to X^\prime \to Y$,
o\`u la premi\`ere fl\`eche est une immersion ferm\'ee et la
seconde un morphisme \'etale.
\end{corollaire}

C'est une simple traduction de~\Ref{I.7.7}.

Montrons comment le jargon de~\Ref{I.7.6} r\'esulte de
l'\'enonc\'e principal: en effet, il existe par~\Ref{I.7.7} un
\'epimorphisme $\cal{O}^\prime \to \cal{O}$, o\`u $\cal{O}$ a les
propri\'et\'es voulues; mais comme $\cal{O}^\prime$ et~$\cal{O}$
sont \'etales sur~$A$, le morphisme $\cal{O}^\prime \to \cal{O}$ est
\'etale par~\Ref{I.4.8} donc un isomorphisme.

\subsubsection*{D\'emonstration de~\Ref{I.7.6}}
Elle reprend une d\'emonstration du s\'eminaire
Chevalley. D'apr\`es le \emph{Main Theorem} on aura $\cal{O} =
B_{\goth{n}}$, o\`u $B$ est une alg\`ebre finie sur~$A$ et
$\goth{n}$ en est un id\'eal maximal. Alors $B/\goth{n} =
K(\cal{O})$ est une extension s\'eparable donc monog\`ene de~$k$;
si $\goth{n}_i$ ($1 \leq i \leq r$) sont les id\'eaux maximaux
de~$B$ distincts de~$\goth{n}$, il existe donc un \'el\'ement $u$
de~$B$ qui appartient \`a tous les~$\goth{n}_i$, et dont l'image
dans~$B/\goth{n}$ en est un g\'en\'erateur. Or $B/\goth{n} =
B_\goth{n}/\goth{n} B_\goth{n} = B_\goth{n}/\goth{m} B_\goth{n}$
(o\`u $\goth{m}$ est l'id\'eal maximal de~$A$). Admettons un
instant le

\begin{lemme}
\label{I.7.9}
Soient $A$ un anneau local, $B$ une alg\`ebre finie sur~$A$,
$\goth{n}$ un id\'eal maximal de~$B$, $u$ un \'el\'ement de~$B$
dont l'image dans $B_\goth{n}/\goth{m} B_\goth{n}$ l'engendre comme
alg\`ebre sur~$k = A/\goth{m}$, et qui se trouve dans tous les
id\'eaux maximaux de~$B$ distincts de~$\goth{n}$. Soit $B^\prime =
B[u]$, $\goth{n}^\prime = \goth{n} B^\prime$. Alors l'homomorphisme
canonique $B^\prime_{\goth{n}^\prime} \to B_\goth{n}$ est un
isomorphisme.
\end{lemme}

\begin{lemme}
\label{I.7.10}
(aurait d\^u figurer en corollaire \`a~\Ref{I.7.1}
avant~\Ref{I.7.5} qu'il implique). Soit $B$ une alg\`ebre finie
sur~$A$ engendr\'ee par un \'el\'ement~$u$, soit $\goth{n}$ un
id\'eal maximal de~$B$ tel que $B_\goth{n}$ soit non ramifi\'e
sur~$A$. Alors il existe un polyn\^ome unitaire $F \in A[t]$ tel que
$F(u) = 0$ et $F^\prime(u) \notin \goth{n}$.
\end{lemme}

Soit en effet $n$ le rang de la $k$-alg\`ebre $L = B \otimes_A k$,
d'apr\`es Nakayama il existe un polyn\^ome unitaire de
degr\'e~$n$ dans $A[t]$, tel que $F(u) = 0$. Soit $f$ le
polyn\^ome d\'eduit de~$F$ par r\'eduction mod~$\goth{m}$, alors
$L$ est $k$-isomorphe \`a $k[t]/fk[t]$, donc par~\Ref{I.7.3}
$f^\prime(\xi)$ n'est pas contenu dans l'id\'eal maximal de~$L$ qui
correspond \`a~$\goth{n}$ ($\xi$ d\'esignant l'image de~$t$
dans~$L$, \ie l'image de~$u$ dans~$L$). Comme $f^\prime(\xi)$ est
l'image de~$F^\prime(u)$, on a fini.

Le
\marginpar{13}
th\'eor\`eme~\Ref{I.7.6} r\'esulte maintenant de la conjonction
de~\Ref{I.7.9} et~\Ref{I.7.10}. Reste \`a
prouver~\Ref{I.7.9}. Posons $S^\prime = B^\prime - \goth{n}^\prime$,
donc $B^\prime {S^\prime}^{-1} = B^\prime_{\goth{n}^\prime}$.
$$
\begin{array}{ccccc}
B&\longrightarrow&B{S^\prime}^{-1}&\longrightarrow&BS^{-1} =
B_\goth{n}\\ \Big\uparrow&&\Big\uparrow&&\\
B^\prime&\longrightarrow&B^\prime{S^\prime}^{-1}
\rlap{\hbox{$\displaystyle \ = B^\prime_{\goth{n}^\prime}$}}&&
\end{array}
$$
Soit de m\^eme $S=B-\goth{n}$, donc $BS^{-1} = B_\goth{n}$, on a
donc un homomorphisme naturel $B{S^\prime}^{-1} \to BS^{-1} =
B_\goth{n}$, prouvons que c'est un isomorphisme, \ie que les
\'el\'ements de~$S$ sont inversibles dans $B{S^\prime}^{-1}$,
\ie que tout id\'eal maximal $\goth{p}$ de ce dernier ne rencontre
pas~$S$, \ie induit~$\goth{n}$ sur~$B$. En effet, comme
$B{S^\prime}^{-1}$ est fini sur~$B^\prime{S^\prime}^{-1} =
B^\prime_{\goth{n}^\prime}$, $\goth{p}$ induit l'unique id\'eal
maximal $\goth{n}^\prime B_{\goth{n}^\prime}$
de~$B^\prime_{\goth{n}^\prime}$, donc induit l'id\'eal maximal
$\goth{n}^\prime$ de~$B^\prime$; comme $B$ est fini sur~$B^\prime$,
l'id\'eal $\goth{q}$ de~$B$ induit par~$\goth{p}$ \'etant
au-dessus de~$\goth{n}^\prime$, est n\'ecessairement maximal, et ne
contient pas~$u$, donc est identique \`a~$\goth{n}$. (On vient
d'utiliser que $u$ appartient \`a tout id\'eal maximal de~$B$
distinct de~$\goth{n}$). Prouvons maintenant que $B{S^\prime}^{-1}$
\'egale $B^\prime{S^\prime}^{-1}$: comme il est fini sur ce dernier,
on est ramen\'e par Nakayama \`a prouver l'\'egalit\'e mod
$\goth{n}^\prime B {S^\prime}^{-1}$ et a fortiori il suffit de prouver
l'\'egalit\'e mod $\goth{m} B{S^\prime}^{-1}$; or
$B{S^\prime}^{-1}/\goth{m} B{S^\prime}^{-1} = B_\goth{n}/\goth{m}
B_\goth{n}$ est engendr\'e sur~$k$ par~$u$ (on utilise ici l'autre
propri\'et\'e de~$u$) donc l'image de~$B^\prime$ (et a fortiori de
$B^\prime {S^\prime}^{-1}$) dedans est tout (comme sous-anneau
contenant~$k$ et l'image de~$u$).

\begin{remarquestar}
On doit pouvoir \'enoncer le th\'eor\`eme~\Ref{I.7.6} pour un
anneau $\cal{O}$ qui est seulement semi-local, de fa\c con \`a
coiffer aussi~\Ref{I.7.5}: on fera l'hypoth\`ese que
$\cal{O}/\goth{m}\cal{O}$ est une $k$-alg\`ebre \emph{monog\`ene};
on pourra donc trouver un $u\in B$ dont l'image dans
\ifthenelse{\boolean{orig}}{ $B/mB$}{$B/\goth{m}B$}
est un g\'en\'erateur, et appartenant \`a tous les id\'eaux
maximaux de~$B$ ne provenant pas de~$\cal{O}$. Les lemmes~\Ref{I.7.9}
et~\Ref{I.7.10} doivent s'adapter sans difficult\'e. Plus
g\'en\'eralement, ...
\end{remarquestar}

\section[Rel\`evement infinit\'esimal des sch\'emas \'etales]{Rel\`evement infinit\'esimal des sch\'emas \'etales.
Application aux sch\'emas formels}
\label{I.8}

\begin{proposition}
\label{I.8.1}
Soient $Y$ un pr\'esch\'ema, $Y_0$ un sous-pr\'esch\'ema,
$X_0$ un $Y_0$-sch\'ema \'etale, $x$ un point de~$X_0$. Alors il
existe un $Y$-sch\'ema \'etale~$X$, un voisinage $U_0$ de~$x$
dans~$X_0$, et un $Y_0$-isomorphisme $U_0 \isomto X \times_Y Y_0$.
\end{proposition}

Soit en effet $y$ la projection de~$x$ dans~$Y_0$,
appliquant~\Ref{I.7.6} au homomorphisme local \'etale $A_0 \to B_0$
des anneaux locaux de~$y$ et~$x$ dans $Y_0$ et~$X_0$: on trouve un
isomorphisme%
\marginpar{14}%
$$
B_0 = {C_0}_{\goth{n}_0}\qquad C_0 = A_0[t]/F_0A_0[t]
$$
o\`u $F_0$ est un polyn\^ome unitaire et $\goth{n}_0$ est un
id\'eal maximal de~$C_0$ ne contenant pas la classe de
$F^\prime_0(t)$ dans~$C_0$. Soit $A$ l'anneau local de~$y$ dans~$Y$,
soit $F$ un polyn\^ome unitaire dans $A[t]$ donnant~$F_0$ par
l'homomorphisme surjectif $A \to A_0$ (on rel\`eve les coefficients
de~$F_0$), soit enfin
\ifthenelse{\boolean{orig}}{$\cal{C} = A[t]/FA[t]$}{$C = A[t]/FA[t]$}
et $\goth{n}$ l'id\'eal maximal
\ifthenelse{\boolean{orig}}{de~$\cal{C}$}{de~$C$}
image inverse de
\ifthenelse{\boolean{orig}}{$\goth{n}_0 x$}{$\goth{n}_0$}
par l'\'epimorphisme naturel $C \to C \otimes_A A_0 = C_0$. Posons
$$
\ifthenelse{\boolean{orig}}{B = C_\goth{n}}{B = C_\goth{n}.}
$$
Il est imm\'ediat par construction et~\Ref{I.7.1} que $B$ est
\'etale sur~$A$, et qu'on a un isomorphisme $B \otimes_A A_0 =
A_0$. On sait (Chap\ptbl I)
qu'il existe un $Y$-sch\'ema de type fini~$X$ et un point~$z$ de~$X$
au-dessus de~$y$ tel que $\cal{O}_z$ soit $A$-isomorphe \`a~$C$;
comme ce dernier est \'etale sur~$A = \cal{O}_y$, on peut (en
prenant $X$ assez petit) supposer que $X$ est \'etale sur~$Y$. Soit
$X_0^\prime = X \times_Y Y_0$, alors l'anneau local de~$z$
dans~$X_0^\prime$ s'identifie \`a $\cal{O}_z \otimes_A A_0 = B
\otimes_A A_0$, donc est isomorphe \`a~$B_0$. Cet isomorphisme est
d\'efini pas un isomorphisme d'un voisinage~$U_0$ de~$x$ dans~$X_0$
sur un voisinage de~$z$ dans~$X_0^\prime$ (loc.\ cit\'e), qu'on peut
supposer identique \`a~$X_0^\prime$ en prenant $X$ assez petit. On a
fini.

\begin{corollaire}
\label{I.8.2}
\'Enonc\'e analogue pour des \emph{rev\^etements} \'etales, en
supposant le corps r\'esiduel $\kres(y)$ infini.
\end{corollaire}

La d\'emonstration est la m\^eme, \Ref{I.7.5} rempla\c cant~\Ref{I.7.6}.

\begin{theoreme}
\label{I.8.3}
Le foncteur envisag\'e dans~\Ref{I.5.5} est une
\emph{\'equivalence} de \emph{cat\'egories}.
\end{theoreme}

En vertu du th\'eor\`eme~\Ref{I.5.5}, il reste \`a montrer que
tout $S_0$-sch\'ema \'etale~$X_0$ est isomorphe \`a un
$S_0$-sch\'ema $X \times_S S_0$, o\`u $X$ est un $S$-sch\'ema
\'etale. L'espace topologique sous-jacent \`a~$X$ devra \^etre
n\'ecessairement identique \`a celui de~$X_0$, $X_0$ s'identifiant
de plus \`a un sous-pr\'esch\'ema ferm\'e de~$X$. Le
probl\`eme est donc \'equivalent au suivant: trouver sur l'espace
topologique sous-jacent $|X_0|$ \`a~$X_0$ un faisceau
d'alg\`ebres~$\cal{O}_X$ sur~$f_0^*(\cal{O}_S)$ (o\`u $f_0$ est
la projection $X_0 \to S_0$, regard\'ee ici comme application
continue des espaces sous-jacents), faisant de~$|X_0|$ un
$S$-pr\'esch\'ema \'etale~$X$, et un homomorphisme
d'alg\`ebres $\cal{O}_X \to \cal{O}_{X_0}$, compatible avec
l'homomorphisme $f_0^*(\cal{O}_S) \to f_0^*(\cal{O}_{S_0})$ sur
les faisceaux de scalaires, induisant un isomorphisme $\cal{O}_X
\otimes_{f_0^*(\cal{O}_S)} f_0^*(\cal{O}_{S_0}) \isomto
\cal{O}_{X_0}$. (Alors $X$ sera un $S$-pr\'esch\'ema \'etale se
r\'eduisant suivant~$X_0$,
\marginpar{15}
donc sera s\'epar\'e sur~$S$ puisque $X_0$ l'est sur~$S_0$, et $X$
r\'epond \`a la question). Si d'ailleurs $(U_i)$ est un
recouvrement de~$X_0$ par des ouverts, et si on a trouv\'e une
solution du probl\`eme dans chacun des~$U_i$, il r\'esulte du
th\'eor\`eme d'unicit\'e~\Ref{I.5.5} que ces solutions se
recollent (\ie les faisceaux d'alg\`ebres qui les d\'efinissent,
munis de leurs homomorphismes d'augmentation, se recollent), et on
constate aussit\^ot que l'espace annel\'e ainsi construit
au-dessus de~$S$ est un $S$-pr\'esch\'ema \'etale~$X$ muni d'un
isomorphisme $X \times_S S_0 \isomfrom X_0$. Il suffit donc de trouver
une solution localement, ce qui est assur\'e par~\Ref{I.8.1}.

\begin{corollaire}
\label{I.8.4}
Soient $S$ un pr\'esch\'ema formel localement noeth\'erien, muni
d'un id\'eal de d\'efinition~$J$, soit $S_0 =
\big(|S|,\cal{O}_S/\cal{J}\big)$ le pr\'esch\'ema ordinaire
correspondant. Alors le foncteur $\goth{X} \mto \goth{X} \times_S
S_0$ de la cat\'egorie des rev\^etements \'etales de~$S$ dans la
cat\'egorie des rev\^etements \'etales de~$S_0$ est une
\'equivalence de cat\'egories.
\end{corollaire}

Bien entendu, on appelera rev\^etement \'etale d'un
pr\'esch\'ema \emph{formel}~$S$ un rev\^etement de~$S$, \ie un
pr\'esch\'ema formel sur~$S$ d\'efini \`a l'aide d'un faisceau
coh\'erent d'alg\`ebres~$\cal{B}$, tel que $\cal{B}$ soit
\emph{localement libre} et que les fibres r\'esiduelles $\cal{B}_s
\otimes_{\cal{O}_s} \kres(s)$ de~$\cal{B}$ soient des alg\`ebres
\emph{s\'eparables} sur~$\kres(s)$. Si on d\'esigne par~$S_n$ le
pr\'esch\'ema ordinaire
\ifthenelse{\boolean{orig}}{$\big(|S|\ \cal{O}_S/J^{n+1}\big)$,}
{$\big(|S|,\cal{O}_S/J^{n+1}\big)$,}
la donn\'ee d'un faisceau coh\'erent d'alg\`ebres $\cal{B}$
sur~$S$ \'equivaut \`a la donn\'ee d'une suite de faisceaux
coh\'erents d'alg\`ebres
\ifthenelse{\boolean{orig}}{ $B_n$}{$\cal{B}_n$}
sur
\ifthenelse{\boolean{orig}}{les~$X_n$,}{les~$S_n$,}
munis d'un syst\`eme transitif d'homomorphismes
\ifthenelse{\boolean{orig}}{$B_m\to B_n$}{$\cal{B}_m\to \cal{B}_n$}
($m \geq n$) d\'efinissant des isomorphismes
\ifthenelse{\boolean{orig}}
{$B_m \otimes_{\cal{O}_{S_m}} \cal{O}_{S_n} \isomto B_n$.}
{$\cal{B}_m \otimes_{\cal{O}_{S_m}} \cal{O}_{S_n} \isomto \cal{B}_n$.}
Il est imm\'ediat que $\cal{B}$ est localement libre si et seulement
si les $\cal{B}_n$ sur les~$S_n$ le sont, et que la condition de
s\'eparabilit\'e est v\'erifi\'ee si et seulement si elle
l'est pour~$\cal{B}_0$, ou encore pour tous
\ifthenelse{\boolean{orig}}{les~$B_n$.}{les~$\cal{B}_n$.}
Ainsi, $\cal{B}$ est \'etale sur~$S$ si et seulement si les
\ifthenelse{\boolean{orig}}{$B_n$}{$\cal{B}_n$}
sur les~$S_n$ le sont. Compte tenu de cela, \Ref{I.8.4} r\'esulte
aussit\^ot de~\Ref{I.8.3}.

\begin{remarquestar}
Il n'\'etait pas n\'ecessaire dans~\Ref{I.8.4} de se borner au cas
des \emph{rev\^etements}. C'est cependant le seul utilis\'e pour
l'instant.
\end{remarquestar}

\section{Propri\'et\'es de permanence}
\label{I.9}
Soit $A\to B$ un homomorphisme local et \'etale, nous examinons ici
quelques cas o\`u une certaine propri\'et\'e pour $A$
entra\^ine la m\^eme propri\'et\'e pour~$B$, ou
r\'eciproquement.
\marginpar{16}
Un certain nombre de telles propositions sont d\'ej\`a
cons\'equences du simple fait que $B$ est \emph{quasi-fini} et
\emph{plat} sur~$A$, et nous nous bornerons \`a en \og rappeler\fg
quelques-unes: $A$~\emph{et} $B$ \emph{ont m\^eme dimension de
Krull, et m\^eme profondeur} (\og codimension cohomologique\fg de
Serre, dans la terminologie encore courante). Il en r\'esulte par
exemple que $A$ \emph{est Cohen-Macaulay si et seulement si $B$ l'est}.
D'ailleurs, pour tout id\'eal premier $\goth{q}$ de~$B$, induisant
$\goth{p}$ sur~$A$, $B_{\goth{q}}$ sera encore quasi-fini et plat
sur~$A_{\goth{p}}$, pourvu qu'on suppose que $B$ soit localis\'ee
d'une alg\`ebre de type fini sur~$A$ (cela r\'esulte du fait que
l'ensemble des points o\`u un morphisme de type fini est quasi-fini
\resp plat est ouvert); et d'ailleurs \emph{tout} id\'eal premier
$\goth{p}$ de~$A$ est induit par un id\'eal premier $\goth{q}$
de~$B$ (car $B$ est \emph{fid\`element} plat sur~$A$). Il en
r\'esulte par exemple que \emph{$\goth{p}$ et~$\goth{q}$ ont
m\^eme rang}; et encore que \emph{$A$ est sans id\'eal premier
immerg\'e si et seulement si~$B$ l'est}.

Nous allons nous borner donc aux propositions plus sp\'eciales au
cas des morphismes \'etales.

\begin{proposition}
\label{I.9.1}
Soit $A\to B$ un homomorphisme local \'etale. Pour que $A$ soit
r\'egulier, il faut et il suffit que $B$ le soit.
\end{proposition}
En effet, soit $k$ le corps r\'esiduel de~$A$, $L$ celui de~$B$.
Comme $B$ est plat sur~$A$ et que
\ifthenelse{\boolean{orig}}
{$L=B\otimes_A k$}
{$L=B\otimes_A k$,}
\ie $\goth{n}=\goth{m}B$ (o\`u $\goth{m}$,~$\goth{n}$ sont les
id\'eaux maximaux de~$A$,~$B$) la filtration $\goth{m}$-adique sur~$B$ est identique \`a sa filtration $\goth{n}$-adique et on aura
$$
\gr^*(B)=\gr^*(A)\otimes_k L.
$$
Il s'ensuit que $\gr^*(B)$ est une alg\`ebre de polyn\^omes sur~$L$ si et seulement si $\gr^*(A)$ est une alg\`ebre de polyn\^omes
sur~$k$. cqfd. (N.B.\ on n'a pas utilis\'e le fait que $L/k$ est
s\'eparable).

\begin{corollaire}
\label{I.9.2}
Soit $f\colon X\to Y$ un morphisme \'etale. Si $Y$ est
r\'egulier, $X$ l'est, la r\'eciproque \'etant vraie si $f$ est
\ifthenelse{\boolean{orig}}{surjective.}{surjectif.}
\end{corollaire}

\setcounter{subsection}{1}

\begin{proposition}
\label{prop:I.9.2}
Soit $f\colon X\to Y$ un morphisme \'etale. Si $Y$ est r\'eduit,
il en est de m\^eme de~$X$, la r\'eciproque \'etant vraie si $f$
est surjectif.
\end{proposition}

Cela \'equivaut au
\begin{corollaire}
\label{I.9.3}
Soit $f\colon A\to B$ un homomorphisme local \'etale, $B$ \'etant
isomorphe \`a une
\marginpar{17}
$A$-alg\`ebre localis\'ee d'une $A$-alg\`ebre de type fini.
Pour que $A$ soit r\'eduit, il faut et il suffit que $B$ le soit.
\end{corollaire}

La n\'ecessit\'e est triviale, puisque $A\to B$ est injectif ($B$
\'etant fid\`element plat sur~$A$). Suffisance: soient
\ifthenelse{\boolean{orig}}
{$p_i$ les id\'eaux premiers minimaux de~$A$, par
hypoth\`ese l'application naturelle $A\to\prod_i A/p_i$ est
injective, donc tensorisant avec le $A$-module plat $B$, on trouve que
$B\to\prod_i B/p_iB$ est injective, et on est ramen\'e \`a
prouver que les $B/p_iB$ sont r\'eduits. Comme
$B/p_iB$ est \'etale sur~$A/p_i$}
{$\goth{p}_i$ les id\'eaux premiers minimaux de~$A$, par
hypoth\`ese l'application naturelle $A\to\prod_i A/\goth{p}_i$ est
injective, donc tensorisant avec le $A$-module plat $B$, on trouve que
$B\to\prod_i B/\goth{p}_iB$ est injective, et on est ramen\'e \`a
prouver que les $B/\goth{p}_iB$ sont r\'eduits. Comme
$B/\goth{p}_iB$ est \'etale sur~$A/\goth{p}_i$},
on est ramen\'e au
cas $A$ int\`egre. Soit~$K$ son corps des fractions, alors $A\to K$
\'etant injectif, il en est de m\^eme ($B$~\'etant $A$-plat) de
$B\to B\otimes_A K$, on est ramener \`a prouver que ce dernier
anneau est r\'eduit. Or, $B$ \'etant localis\'ee d'une
$A$-alg\`ebre de type fini sur~$A$, est l'anneau local d'un point~$x$ d'un sch\'ema de type fini et \emph{\'etale} $X=\Spec(C)$ sur~$Y=\Spec(A)$, donc $B\otimes_A K$ est un anneau localis\'e (par
rapport \`a un ensemble multiplicativement stable convenable) de
l'anneau $C\otimes_A K$ de~$X\otimes_A K$. Comme $X\otimes_A K$ est
\'etale sur~$K$, son anneau est un produit fini de corps (extensions
s\'eparables de~$K$), il en est donc de m\^eme de~$B\otimes_A
K$, cqfd.

\begin{corollaire}
\label{I.9.4}
Soit $f\colon A\to B$ un homomorphisme local \'etale, supposons $A$
analytiquement r\'eduit
\index{analytiquement r\'eduit|hyperpage}%
(\ie le compl\'et\'e $\hat{A}$ de~$A$ sans \'el\'ements
nilpotents). Alors $B$ est analytiquement r\'eduit, et a~fortiori
r\'eduit.
\end{corollaire}

En effet,
\ifthenelse{\boolean{orig}}{$B$}{$\hat{B}$}
est \emph{fini} et \'etale sur~$\hat{A}$, et on applique~\Ref{I.9.3}.

\begin{theoreme}
\label{I.9.5}
Soit $f\colon A\to B$ un homomorphisme local, $B$ \'etant isomorphe
\`a une alg\`ebre localis\'ee d'une $A$-alg\`ebre de type
fini. Alors:
\begin{enumerate}
\item[(i)] Si $f$ est \'etale, $A$ est normal si et seulement si $B$
l'est.
\item[(ii)] Si $A$ est normal, $f$ est \'etale si et seulement si
$f$ est injectif et net (et alors $B$ est normal par \textup{(i)}).
\end{enumerate}
\end{theoreme}

Nous allons donner deux d\'emonstrations diff\'erentes de (i), la
premi\`ere utilise certaines des propri\'et\'es des morphismes
plats quasi-finis (rappel\'es au d\'ebut du num\'ero) sans
utiliser~\Ref{I.7.6} (et par l\`a, le Main Theorem); c'est l'inverse pour la
deuxi\`eme d\'emonstration. Enfin, pour~(ii) il semble qu'on ait
besoin du Main Theorem en tous cas.

\subsubsection*{Premi\`ere d\'emonstration}
On utilise la condition n\'ecessaire et suffisante suivante
\marginpar{18}
de normalit\'e d'un anneau local noeth\'erien $A$ de dimension
$\neq0$.

\begin{enonce*}{Crit\`ere de Serre}
\index{Serre (crit\`ere de)|hyperpage}%
\textup{(i)}~Pour tout id\'eal premier $\goth{p}$ de~$A$ de rang~$1$,
$A_{\goth{p}}$ est normal (ou ce qui revient au m\^eme,
r\'egulier); \textup{(ii)}~Pour tout id\'eal premier $\goth{p}$ de~$A$ de
rang~$\geq 2$, on a profondeur $A_{\goth{p}}$ $\geq2$.\footnote{\Cf EGA~IV~5.8.6.}
\end{enonce*}

Nous admettrons ici ce crit\`ere, qui est cens\'e figurer au
\ifthenelse{\boolean{orig}}{par.}{paragraphe}
des plats. Son principal avantage est qu'il ne suppose pas a priori
$A$ r\'eduit, ni a~fortiori int\`egre. Ici, on peut d\'ej\`a
supposer $\dim{A}=\dim{B}\neq0$.

D'apr\`es les rappels du d\'ebut du num\'ero, les id\'eaux
premiers $\goth{p}$ de~$A$ qui sont de rang~$1$ (\resp de
rang~$\ge2$) sont exactement les traces sur~$A$ des id\'eaux
premiers~$\goth{q}$ de~$B$ qui sont de rang~$1$ (\resp de
rang~$\ge2$). Enfin, si $\goth{p}$ et $\goth{q}$ se correspondent,
$B_{\goth{q}}$ est \'etale sur~$A_{\goth{p}}$, donc a m\^eme
profondeur que $A_{\goth{p}}$, et est r\'egulier si et seulement si
$A_{\goth{p}}$ l'est (\Ref{I.9.1}). Appliquant le crit\`ere de Serre, on
trouve que $A$ est normal si et seulement si $B$ l'est.

\subsubsection*{Deuxi\`eme d\'emonstration}
Supposons $B$ normal, soit $L$ son corps des fractions, $K$ celui de
$A$ ($A$~est int\`egre puisque $B$ l'est). On a vu dans la
d\'emonstration de~\Ref{I.9.3} que $B\otimes_A K$ est un compos\'e fini de
corps, comme il est contenu dans $L$ c'est un corps, et comme il
contient $B$ c'est $L$. Un \'el\'ement de~$K$ entier sur~$A$ est
entier sur~$B$, donc est dans $B$ puisque
\ifthenelse{\boolean{orig}}{$b$}{$B$}
est normal, donc dans $A$ car $B\cap K=A$ (comme il r\'esulte du
fait que $B$ est fid\`element plat sur~$A$).

Supposons maintenant $A$ normal, prouvons que $B$ l'est. En vertu
de~\Ref{I.7.6} on aura $B=B^\prime_{\goth{n}}$, o\`u $B^\prime=A[t]/FA[t]$,
$F$ et $\goth{n}$ \'etant comme dans~\Ref{I.7.6}. Donc $L=B\otimes_A K$
sera un localis\'e de~$B^\prime\otimes_A K=K[t]/FK[t]$, et un
produit de corps; extensions finies s\'eparables de~$K$ ce dernier
produit (comme chaque fois qu'on localise un anneau artinien, ici
$B^\prime_K$ par rapport \`a un ensemble multiplicativement stable)
est un facteur direct de~$B^\prime_K$, correspondant donc \`a une
d\'ecomposition $F=F_1F_2$ dans $K[t]$, le g\'en\'erateur de~$L$
correspondant \`a $t$ \'etant annul\'e d\'ej\`a par~$F_1$.
Or, \emph{$A$ \'etant normal, les $F_i$ sont dans $A[t]$} (supposant
qu'ils sont unitaires). Remarquant que $B\to L=B\otimes_A K$ est
injectif ($A\to K$ l'\'etant et $B$ \'etant plat sur~$A$) il
s'ensuit qu'on aura d\'ej\`a
\ifthenelse{\boolean{orig}}{$F_1(u)=0$.}{$F_1(u)=0$, avec $u$ la
classe de~$t$ dans~$L$.}
Supposant qu'on ait pris $F$ de degr\'e minimum, il s'ensuivra que
$F_2=1$ (N.B.\ on aura
$F^\prime(u)=F^\prime_1(u)F_2(u)+F_1(u)F^\prime_2(u)=F^\prime_1(u)F_2(u)$
puisque $F_1(u)=0$, d'o\`u $F^\prime_1(u)\neq0$ puisque
\ifthenelse{\boolean{orig}}{$F^\prime(u)\neq0$.}{$F^\prime(u)\neq0$.)}

Donc
\marginpar{19}
on a
\begin{equation*}
\label{eq:I.9.5.*}
\tag{$*$} L=B\otimes_A K=K[t]/FK[t]
\end{equation*}
$F$ \'etant par suite un polyn\^ome s\'eparable dans $K[t]$
(mais \'evidemment par n\'ecessairement dans $A[t]$). (N.B.\ Pour
l'instant, on a seulement montr\'e, essentiellement, que dans~\Ref{I.7.6} on
peut choisir $F$ et
\ifthenelse{\boolean{orig}}{\ \ }{$\goth{n}$}
de telle fa\c con que --- avec les notations prises ici --- \ifthenelse{\boolean{orig}}
{$B^\prime\to B^\prime_{\goth{n}}=\cal{B}$}
{$B^\prime\to B^\prime_{\goth{n}}=B$}
soit \emph{injectif}; on s'est servi pour cela de la normalit\'e de
$A$; je ne sais pas si cela reste vrai sans hypoth\`ese de
normalit\'e).

Rappelons maintenant le lemme bien connu, extrait du Cours de Serre de
l'an dernier:

\begin{lemme}
\label{I.9.6}
Soient $K$ un anneau, $F\in K[t]$ un polyn\^ome unitaire
s\'eparable, $L=K[t]/FK[t]$, $u$~la classe de~$t$ dans~$L$ (de sorte
que $F^\prime(u)$ est un \'el\'ement inversible de~$L$). Alors on
a les formules (o\`u $n=\deg{F}$):
\begin{enumerate}
\item[]$\trace_{L/K} u^i/F^\prime(u)=0$\ \ si\ \
\ifthenelse{\boolean{orig}}{$0\le i< N-1$,}{$0\le i< n-1$,}
\item[]$\trace_{L/K} u^{n-1}/F^\prime(u)=1$.
\end{enumerate}
\end{lemme}

\begin{corollaire}
\label{I.9.7}
Le d\'eterminant de la matrice $(u^j\cdot u^i/F^\prime(u))_{0\le i,j\le
n-1}$ est \'egal \`a
\ifthenelse{\boolean{orig}}{$(-1)^n$}{$(-1)^{n(n-1)/2}$},
donc inversible dans tout
sous-anneau $A$ de~$K$.
\end{corollaire}

\begin{corollaire}
\label{I.9.8}
\ifthenelse{\boolean{orig}}{Soit}{Soient}
$A$ un sous-anneau de~$K$, $V$ le $A$-module engendr\'e par les
$u^i$~($0\le i\le n-1$) dans $L$, $V^\prime$ le sous-$A$-module de~$L$
form\'e des $x\in L$ tels que $\trace_{L/K}(xy)\in A$ pour tout
$y\in V$ (\ie pour $y$ de la forme $u^i$, $0\le i\le n-1$). Alors
$V^\prime$ est le $A$-module ayant pour base les $u^i/F^\prime(u)$
($0\le i\le n-1$).
\end{corollaire}

\begin{corollaire}
\label{I.9.9}
Supposons que $K$ soit le corps des fractions d'un anneau int\`egre
normal $A$, $F$ ayant ses coefficients dans $A$. Alors avec les
notations de \Ref{I.9.8}, $V^\prime$ contient la cl\^oture normale
$A^\prime$ de~$A$ dans $L$, qui est donc contenue dans
$A[u]/F^\prime(u)$ et a~fortiori dans $A[u][F^\prime(u)^{-1}]$.
\end{corollaire}

Appliquons ce dernier corollaire \`a la situation que nous avions
obtenue dans la d\'e\-mon\-stra\-tion: comme $F^\prime(u)$ est
inversible dans $B$ qui contient $A[u]$, $B$ contient $A^\prime$.
D'apr\`es le Main Theorem, (ou \`a partir du fait que
$B=A[u]_{\goth{n}}$) $B$ est une alg\`ebre localis\'ee de
$A^\prime$. Comme $A^\prime$ est normal, il en est de m\^eme
de~$B$.

\subsubsection*{D\'emonstration de \textup{(ii)}}
On
\marginpar{20}
proc\`ede comme dans la d\'emonstration qui pr\'ec\`ede pour
prouver qu'on peut, dans~\Ref{I.7.6}, choisir $F$ de telle fa\c con que l'on
ait encore ($*$). Le seul obstacle a priori est que, $B$ n'\'etant
plus suppos\'e plat sur~$A$, on ne peut plus affirmer que $B\to L$
est injectif, de sorte que le raisonnement ne s'appliquera a priori
qu'\`a l'image $B_1$ de~$B$ par ledit homomorphisme. Il s'ensuit
aussit\^ot que $B_1$ est plat sur~$A$ (comme localis\'ee d'une
alg\`ebre libre sur~$A$). En vertu de \Ref{I.4.8} le morphisme $B\to B_1$
est \'etale, donc un isomorphisme, ce qui ach\`eve la
d\'emonstration.

(Du point de vue r\'edaction, il faudrait intervertir les deux
derni\`eres d\'emonstrations, et mettre dans un num\'ero \`a
part les calculs formels du lemme et de ses corollaires).

\begin{corollaire}
\label{I.9.10}
Soit $f\colon X\to Y$ un morphisme \'etale. Si $Y$ est normal, $X$
l'est, la r\'eciproque est vraie si $f$ est
\ifthenelse{\boolean{orig}}{surjective.}{surjectif.}
\end{corollaire}

\begin{corollaire}
\label{I.9.11}
Soit $f\colon X\to Y$ un morphisme dominant, $Y$ \'etant normal et
$X$ connexe. Si $f$ est net, $f$ est \'etale, donc $X$ est normal
et par suite (\'etant connexe) irr\'eductible.
\end{corollaire}

Soit $U$ l'ensemble des points o\`u $f$ est \'etale, il est
ouvert, et il suffit de montrer qu'il est aussi ferm\'e et non vide.
$U$ contient l'image inverse du point g\'en\'erique de~$Y$ (car
pour une alg\`ebre sur un corps, non ramifi\'e $=$ \'etale) donc
($X$ dominant $Y$) est non vide. Si $x$ appartient \`a
l'adh\'erence de~$U$, alors il appartient \`a l'adh\'erence
d'une composante irr\'eductible $U_i$ de~$U$, donc \`a une
composante irr\'eductible $X_i\overset{v}{=}\bar{U}_i$
de~$X$ qui rencontre $U$, et par suite domine $Y$ (car toute
composante de~$U$, plat sur~$Y$, domine $Y$). Par suite, si $y$ est
la projection de~$x$ sur~$Y$, $\cal{O}_y\to\cal{O}_x$ est
\emph{injectif} (compte tenu que $\cal{O}_y$ est int\`egre). Comme
$\cal{O}_y$ est normal et $\cal{O}_y\to\cal{O}_x$ net, on conclut
\`a l'aide de \Ref{I.9.5}(ii).

\begin{corollaire}
\label{I.9.12}
Soit $f\colon X\to Y$ un morphisme de type fini dominant, avec $Y$
normal et $X$ irr\'eductible. Alors l'ensemble des points o\`u
$f$ est \'etale est identique au compl\'ementaire du support de
$\Omega^1_{X/Y}$, \ie au compl\'ementaire du
sous-pr\'esch\'ema de~$X$ d\'efini par l'id\'eal
diff\'erente $\goth{d}_{X/Y}$.
\end{corollaire}

(C'est
\marginpar{21}
cela l'\'enonc\'e \og moins trivial\fg auquel il \'etait fait
allusion dans la remarque du~\No\Ref{I.4}.)

\begin{remarquestar}
On se gardera de croire qu'un rev\^etement \'etale connexe d'un
sch\'ema irr\'eductible soit lui-m\^eme irr\'eductible, quand
on ne suppose pas la base normale. Cette question sera
\'etudi\'ee au~\No \Ref{I.11}.
\end{remarquestar}

\section{Rev\^etements \'etales d'un sch\'ema normal}
\label{I.10}

\begin{proposition}
\label{I.10.1}
Soit $X$ un pr\'esch\'ema \'etale s\'epar\'e sur~$Y$ normal
connexe de corps~$K$. Alors les composantes connexes
\ifthenelse{\boolean{orig}}{de}{}
$X_i$ de~$X$ sont int\`egres, leurs corps $K_i$ sont des extensions
finies s\'eparables de~$K$, $X_i$ s'identifie \`a une partie
ouverte non vide du normalis\'e de~$X$ dans $K_i$ (donc $X$ \`a
une partie ouverte dense du normalis\'e de~$Y$ dans $R(X)=L=\prod
K_i$).
\end{proposition}

D'apr\`es \Ref{I.9.10} $X$ est normal, a~fortiori ses anneaux locaux sont
int\`egres, donc les composantes connexes de~$X$ sont
irr\'eductibles. Comme $X_i$ est normal, et fini et dominant
au-dessus de~$Y$, il r\'esulte d'un cas particulier (\`a peu
pr\`es trivial d'ailleurs) du Main Theorem que $X_i$ est un ouvert
du normalis\'e de~$X$ dans le corps $K_i$ de~$X$.

\begin{corollaire}
\label{I.10.2}
Sous les conditions \Ref{I.10.1}, $X$ est fini sur~$Y$ (\ie un
rev\^etement \'etale de~$Y$) si et seulement si $X$ est isomorphe
au normalis\'e $Y'$ de~$Y$ dans $L=R(X)$ (anneau des fonctions
rationnelles sur~$X$).
\end{corollaire}

En effet, on sait que ce normalis\'e est fini sur~$Y$ ($Y$ \'etant
normal et $R/K$ s\'eparable), inversement si $X$ est fini sur~$Y$ il
l'est sur~$Y'$, donc son image dans $Y'$ est ferm\'ee, d'autre part
elle est dense.

Une alg\`ebre $L$ de rang fini sur~$K$ sera dite \emph{non
ramifi\'ee sur~$X$} (ou simplement non ramifi\'ee sur~$K$, si $X$
est sous-entendu) si $L$ est une alg\`ebre s\'eparable sur~$K$,
\ie compos\'ee directe d'extensions s\'eparables $K_i$, et si le
normalis\'e $Y'$ de~$Y$ dans $L$ (somme disjointe des normalis\'es
de~$Y$ dans les $K_i$) est non ramifi\'e ($=$ \'etale par \Ref{I.9.11})
sur~$Y$. Donc:

\begin{corollaire}
\label{I.10.3}
Pour tout $X$ fini sur~$Y$ et dont toute composante irr\'eductible
domine $Y$, soit $R(X)$ l'anneau des fonctions rationnelles sur~$X$
(produit des anneaux
\marginpar{22}
locaux des points g\'en\'eriques des composantes irr\'eductibles
de~$X$), de sorte que $X\mto R(X)$ est un foncteur, \`a valeurs
dans les alg\`ebres de rang fini sur~$K=R(Y)$. Ce foncteur
\'etablit une \'equivalence de la cat\'egorie des rev\^etement
\'etales connexes de~$Y$ avec la cat\'egorie des extensions $L$ de
$K$ non ramifi\'ees sur~$Y$.
\end{corollaire}

Le foncteur inverse est le foncteur normalisation.

Supposons $Y$ affine, donc d\'efini par un anneau normal $A$ de corps
des fractions~$K$. Soit $L$ une extension finie de~$K$ compos\'e
directe de corps, alors par d\'efinition la normalis\'ee $Y'$ de
$Y$ dans $L$ est isomorphe \`a $\Spec(A')$, o\`u $A'$ est le
normalis\'e de~$A$ dans~$L$. Dire que $L$ est non ramifi\'e sur~$Y$ signifie que $A'$ est non ramifi\'e (ou encore: \'etale)
sur~$A$.
\ifthenelse{\boolean{orig}}{si}{Si}
$A$ est local, il revient au m\^eme de dire que les anneaux locaux
$A^\prime_{\goth{n}}$ (o\`u
\ifthenelse{\boolean{orig}}{$n$}{$\goth{n}$}
parcourt l'ensemble fini des id\'eaux maximaux de~$A'$, \ie de ses
id\'eaux premiers induisant l'id\'eal maximal $\goth{m}$ de~$A$)
soient non ramifi\'es ($=$ \'etales) sur
\ifthenelse{\boolean{orig}}{$X$}{}
l'anneau local~$A$. Enfin, notons aussi que le crit\`ere par le
discriminant (\Ref{I.4.10}) peut aussi s'appliquer dans cette situation (plus
g\'en\'eralement, une variante dudit crit\`ere devrait
s'\'enoncer ainsi, sans condition pr\'eliminaire de platitude
lorsque $X$ domine $Y$, $Y$ \'etant n\'eanmoins suppos\'e
localement int\`egre: $A\to B$ et $B\to B\otimes_A K =L$ sont
injectifs --- alors $\trace_{L/K}$ est d\'efinie --- et
$\trace_{L/K}(xy)$ induit une \emph{forme bilin\'eaire fondamentale}
$B\times B\to A$, \ie il existe des $x_i\in B$ ($1\le i\le n$, $n=$
rang de~$L$ sur~$K$) tels que $\trace(x_ix_j)\in A$ pour tout $i$,~$j$
et
\ifthenelse{\boolean{orig}}{de~$(\trace(x_ix_j))$}
{$\det(\trace(x_ix_j))$}
est inversible dans~$A$).

Le sorite (4.6) implique aussit\^ot le sorite de la non ramification
dans le cadre classique:

\begin{proposition}
\label{I.10.4}
\ifthenelse{\boolean{orig}}{Soient}{Soit}
$Y$ un pr\'esch\'ema normal int\`egre, de corps~$K$.
\textup{(i)} $K$ est non ramifi\'e
\ifthenelse{\boolean{orig}}{sur~$Y$}{sur~$Y$.}
\textup{(ii)} Si $L$ est une extension de~$K$ non ramifi\'ee sur~$Y$, si $Y'$
est un pr\'esch\'ema normal de corps $L$ et dominant $Y$ (par
exemple le normalis\'e de~$Y$ dans~$L$) et $M$ une extension de~$L$
non ramifi\'ee sur~$Y'$, alors $M/K$ est non ramifi\'ee sur~$X$
(\emph{transitivit\'e} de la non ramification). \textup{(iii)} Soit $Y'$ un
pr\'esch\'ema normal int\`egre dominant $Y$, de corps $K'/K$; si
$L$ est une extension de~$K$ non ramifi\'ee sur~$Y$, alors
$L\otimes_K K'$ est une extension de~$K'$ non ramifi\'ee sur~$Y'$
(propri\'et\'e de \emph{translation}).
\end{proposition}

De
\marginpar{23}
plus:

\begin{corollaire}
\label{I.10.5}
Sous les conditions de \textup{(iii)}, si $Y=\Spec(A)$, $Y'=\Spec(A')$, alors
le normalis\'e
\ifthenelse{\boolean{orig}}
{$A'$}
{$\bar{A'}$}
de~$Y'$ dans $L'=L\otimes_K K'$ s'identifie
\`a $\bar{A}\otimes_A A'$, o\`u $\bar{A}$ est le normalis\'e de
$A$ dans $L$.
\end{corollaire}

Habituellement, les gens (qui r\'epugnent \`a la consid\'eration
d'anneaux non int\`egres, fussent-ils compos\'es directs de corps)
\'enoncent la propri\'et\'e de translation sous la forme (plus
faible) suivante:

\begin{corollaire}
\label{I.10.6}
Sous les conditions de \textup{(iii)}, soit $L_1$ une \emph{extension
compos\'ee} de~$L/K$ (non ramifi\'ee sur~$Y$) et de~$K'/K$. Alors
$L_1/K'$ est non ramifi\'ee sur~$Y'$. Dans le cas o\`u
$Y=\Spec(A)$, $Y'=\Spec(A')$, on aura de plus
$$
\overline{A'}=A[\overline{A},A']
$$
\ie l'anneau $\overline{A'}$ normalis\'e de~$A'$ dans $L_1$ est la
$A$-alg\`ebre engendr\'ee par $A'$ et le normalis\'e
$\overline{A}$ de~$A$ dans $L$.
\end{corollaire}

Ce dernier fait est d'ailleurs faux sans hypoth\`ese de non
ramification, m\^eme dans le cas d'extensions compos\'ees de corps
de nombres...

Pour terminer ce num\'ero, nous allons donner l'interpr\'etation
de la notion de rev\^etement \'etale correspondant \`a l'image
intuitive de cette notion: il doit y avoir le \og nombre maximum\fg de
points au-dessus du point consid\'er\'e $y\in Y$, et en
particulier il ne doit pas y avoir \og plusieurs points confondus\fg
au-dessus de~$y$. Pour d\'emontrer les r\'esultats dans ce sens
avec toute la g\'en\'eralit\'e d\'esirable, nous allons
admettre ici la proposition~\Ref{I.10.7} plus bas (dont la d\'emonstration
sera dans le multiplodoque, Chap\ptbl IV, par\ptbl 15, et utilise la technique
des ensembles constructibles de Chevalley, et un petit peu de
th\'eorie de descente...).

Un morphisme de type fini $f\colon X\to Y$ est dit
\emph{universellement ouvert} si pour toute extension de la base
$Y'\to Y$ (avec $Y'$ localement noeth\'erien) le morphisme $f'\colon
X'=X\times_Y Y'\to Y'$ est ouvert, \ie transforme ouverts en
ouverts. On peut d'ailleurs se borner au cas o\`u $Y'$ est de type
fini sur~$Y$ (et m\^eme o\`u $Y'$ est de la forme
$Y[t_1,\dots,t_r]$, o\`u les $t_i$ sont des ind\'etermin\'ees).
Un morphisme universellement ouvert est a~fortiori ouvert (la
r\'eciproque \'etant fausse), d'autre part si $f$ est ouvert, $X$
et $Y$ \'etant irr\'eductibles, alors toutes les composantes de
toutes les fibres de~$f$ ont m\^eme dimension (savoir la dimension de
la fibre g\'en\'erique
\marginpar{24}
$f^{-1}(z)$, $z$ le point g\'en\'erique
de~$Y$). Enfin si $Y$ est normal, cette derni\`ere condition
implique d\'ej\`a que $f$ est \emph{universellement} ouvert
(th\'eor\`eme de Chevalley). Il s'ensuit par exemple que si
$f\colon X\to Y$ est un morphisme \emph{quasi-fini}, avec $Y$ normal
irr\'eductible, alors $f$ est universellement
\ifthenelse{\boolean{orig}}{ouverte}{ouvert}
(ou encore:
\ifthenelse{\boolean{orig}}{ouverte}{ouvert})
si et seulement si toute composante irr\'eductible de~$X$
domine~$Y$. Rappelons aussi qu'un morphisme plat (de type fini)
\'etant ouvert, est aussi universellement ouvert. Ces
pr\'eliminaires pos\'es, \og rappelons\fg la

\begin{proposition}
\label{I.10.7}
Soit $f\colon X\to Y$ un morphisme quasi-fini s\'epar\'e
universellement ouvert. Pour tout $y\in Y$, soit $n(y)$ le \og nombre
g\'eom\'etrique de points de la fibre $f^{-1}(y)$\fg, \'egal
\`a la somme des degr\'es s\'eparables des extensions
r\'esiduelles
\ifthenelse{\boolean{orig}}
{$K(x)/K(y)$}
{$\kres(x)/\kres(y)$}, pour les points $x\in f^{-1}(y)$. Alors
la fonction $y\mto n(y)$ sur~$Y$ est semi-continue
sup\'erieurement. Pour qu'elle soit constante au voisinage du point
$y$ (\ie pour qu'on ait $n(y)=n(z_i)$, o\`u les $z_i$ sont les
points g\'en\'eriques des composantes irr\'eductibles de~$Y$ qui
contiennent $y$) il faut qu'il existe un voisinage $U$ de~$y$ tel que
$X|U$ soit \emph{fini} sur~$U$.\footnote{\Cf EGA IV 15.5.1}
\end{proposition}

\begin{corollaire}
\label{I.10.8}
Si $y\mto n(y)$ est constante et $Y$ g\'eom\'etriquement
unibranche\footnote{Pour la d\'efinition, \cf ci-dessous
\no \Ref{I.11}, p\ptbl \pageref{I.11}}, les composantes irr\'eductibles
de~$X$ sont disjointes.
\end{corollaire}

\begin{proposition}
\label{I.10.9}
Soit $f\colon X\to Y$ un morphisme \emph{\'etale} s\'epar\'e.
Avec les notations~\Ref{I.10.7} la fonction $y\mto n(y)$ est
semi-continue sup\'erieurement. Pour qu'elle soit constante au
voisinage du point $y$, (\ie pour qu'on ait $n(y)=n(z_i)$, o\`u
les $z_i$ sont les points g\'en\'eriques des composantes
irr\'eductibles de~$Y$ qui contiennent $y$) il faut et il suffit
qu'il existe un voisinage ouvert $U$ de~$y$ tel que $X|U$ soit fini
sur U, \ie soit un \emph{rev\^etement \'etale} de~$U$.
\end{proposition}

\begin{corollaire}
\label{I.10.10}
Pour qu'un morphisme \'etale s\'epar\'e $f\colon X\to Y$, $Y$
connexe, soit fini (\ie fasse de~$X$ un \emph{rev\^etement
\'etale} de~$Y$) il faut et il suffit que toutes les fibres de~$f$
aient m\^eme nombre g\'eom\'etrique de points.
\end{corollaire}

\ifthenelse{\boolean{orig}}{Dans~\Ref{I.10.8}}{Dans~\Ref{I.10.7}}
et son corollaire, il n'y avait pas d'hypoth\`ese de normalit\'e
sur~$Y$. Si on fait une telle hypoth\`ese, on trouve
l'\'enonc\'e plus fort (pris le plus souvent comme
d\'efinition de la non ramification d'un rev\^etement):

\begin{theoreme}
\label{I.10.11}
Soit
\marginpar{25}
$f\colon X\to Y$ un morphisme \emph{quasi-fini} s\'epar\'e.
On suppose que~$Y$ est \emph{irr\'eductible}, que toute composante
de~$X$ domine $Y$, que $X$ soit r\'eduit (\ie $\cal{O}_X$ sans
\'el\'ements nilpotents). Soit $n$ le degr\'e de~$X$ sur~$Y$
(somme des degr\'es, sur le corps $K$ de~$Y$, des corps $K_i$ des
composantes irr\'eductibles $X_i$ de~$X$). Soit $y$ un point normal
de~$Y$. Alors le nombre g\'eom\'etrique $n(y)$ de points de~$X$
au-dessus de~$y$ est $\le n$, l'\'egalit\'e ayant lieu si et
seulement si il existe un voisinage ouvert $U$ de~$y$ tel que $X|U$
soit un \emph{rev\^etement \'etale} de~$U$.
\end{theoreme}

Le \og seulement si\fg \'etant trivial, prouvons le \og si\fg. Soit $z$ le
point g\'en\'erique de~$Y$, on a $n(z)=$ (somme des degr\'es
s\'eparables des $K_i/K$) $\le n$ et par \Ref{I.10.7} on a $n(y)\le n(z)\le
n$, l'\'egalit\'e impliquant que $X|U$ est \emph{fini} sur~$U$
pour un voisinage~$U$ convenable de~$y$. On peut donc supposer $X$
fini sur~$Y$ et la fonction $n(y')$ sur~$Y$ constante. Enfin par~\Ref{I.10.8}
$X$ est alors r\'eunion disjointe de ses composantes
irr\'eductibles et pour prouver qu'il est non ramifi\'e en $y$, on
est ramen\'e au cas o\`u $X$ est irr\'eductible, donc
int\`egre. Enfin on peut supposer $Y=\Spec(\cal{O}_y)$. Le
th\'eor\`eme se r\'eduit alors \`a l'\'enonc\'e classique
suivant:

\begin{corollaire}
\label{I.10.12}
Soient $A$ un anneau local normal (noeth\'erien comme toujours) de
corps $K$, $L$ une extension finie de~$K$ de degr\'e $n$, degr\'e
s\'eparable $n_s$, $B$ un sous-anneau de~$L$ fini sur~$A$, de corps
des fractions
\ifthenelse{\boolean{orig}}{$L$.,}{$L$,}
$\goth{m}$ l'id\'eal maximal de~$A$ et $n'$ le
degr\'e s\'eparable de~$B/\goth{m}B$ sur~$A/\goth{m}A=k$ ($=$
somme des degr\'es s\'eparables des extensions r\'esiduelles de
cet anneau). On a $n'\le n_s$ et a~fortiori $n'\le n$. Cette
derni\`ere in\'egalit\'e est une \'egalit\'e si et seulement
si $B$ est non~ramifi\'e ($=$ \'etale) sur~$A$.
\end{corollaire}

Il reste seulement \`a montrer que $n'=n$ implique que $B$ est
\'etale sur~$A$. Rappelons la d\'emonstration quand $k$ est
infini: on doit seulement montrer que $R=B/\goth{m}B$ est
s\'eparable sur~$k$; s'il n'en \'etait pas ainsi il en
r\'esulterait (par un lemme connu) qu'il existe un \'el\'ement
$a$ de~$R$ dont le polyn\^ome minimal sur~$k$ est de degr\'e
$>n'$. Cet \'el\'ement provient d'un \'el\'ement $x$ de~$B$,
dont le polyn\^ome minimal sur~$K$ (en tant qu'\'el\'ement de
$L$) est de degr\'e $\le n$; d'autre part ce dernier a ses
coefficients dans~$A$ puisque $A$ est normal, et donne donc par
r\'eduction mod $\goth{m}$ un polyn\^ome unitaire $F\in k[t]$, de
degr\'e $\le n=n'$, tel que $F(a)=0$, absurde.

Dans
\marginpar{26}
le cas g\'en\'eral ($k$ pouvant \^etre fini), reprenant le
langage g\'eom\'etrique, on consid\`ere $Y'=\Spec(A[t])$ qui est
fid\`element plat sur~$Y$, et le point g\'en\'erique $y'$ de la
fibre $\Spec(k[t])$ de~$Y'$ sur~$y$. Alors $X$ est net sur~$Y$ en $y$
si et seulement si $X'=X\times_Y Y'=\Spec(B[t])$ est net en $y'$ sur~$Y'$, comme on constate aussit\^ot. D'autre part, d'apr\`es le
choix de~$y'$, son corps r\'esiduel est $k(t)$ donc infini. Comme
$y'$ est un point normal de~$Y'$, on est ramen\'e au cas
pr\'ec\'edent.







\section{Quelques compl\'ements}
\label{I.11}

Nous avons d\'ej\`a dit qu'un rev\^etement \'etale connexe
d'un sch\'ema int\`egre n'est pas n\'e\-ces\-sai\-rement
int\`egre. Voici deux exemples de ce fait.

a) Soit $C$ une courbe alg\'ebrique \`a point double
ordinaire~$x$, $C^\prime$ sa normalis\'ee, $a$ et~$b$ les deux
points de~$C^\prime$ au-dessus de~$x$. Soient $C_i^\prime$ ($i=1$,
$2$) deux copies de~$C^\prime$, $a_i$ et~$b_i$ le point de
$C_i^\prime$ qui correspond \`a $a$ \resp $b$. Dans la courbe somme
$C_1^\prime \amalg C_2^\prime$, identifions
\ifthenelse{\boolean{orig}}{$a_1^x$}{$a_1$}
et~$b_2$ d'une part, $a_2$ et~$b_1$ d'autre part (on laisse au lecteur
le soin de pr\'eciser le processus d'identification; il sera
expliqu\'e au Chap\ptbl VI du multiplodoque mais, dans le cas des
courbes sur un corps alg\'ebriquement clos, est trait\'e dans le
livre de Serre sur les courbes alg\'ebriques). On trouve une courbe
$C^{\prime\prime}$ \emph{connexe} et \emph{r\'eductible}, qui est un
rev\^etement \'etale de degr\'e~$2$ de~$C$. Le lecteur
v\'erifiera que de fa\c con g\'en\'erale, les rev\^etements
\'etales connexes \og galoisiens\fg $C^{\prime\prime}$ de~$C$ dont
l'image inverse $C^{\prime\prime} \times_C C^\prime$ est un
rev\^etement \emph{trivial} de~$C^\prime$ (\ie isomorphe \`a la
somme d'un certain nombre de copies de~$C^\prime$) sont \og cycliques\fg
de degr\'e~$n$, et pour tout entier~$n>0$, on peut construire un
rev\^etement \'etale connexe cyclique de degr\'e~$n$. Dans le
langage du groupe fondamental qui sera d\'evelopp\'e plus tard,
cela signifie que le quotient de~$\pi_1(C)$ par le sous-groupe
invariant ferm\'e engendr\'e par l'image de~$\pi_1(C^\prime) \to
\pi_1(C)$ (homomorphisme induit par la projection) est isomorphe au
compactifi\'e de~$\ZZ$. De fa\c con plus pr\'ecise, on doit
pouvoir montrer que le groupe fondamental de~$C$ est isomorphe au
produit libre (topologique) du groupe fondamental
\ifthenelse{\boolean{orig}}{de~$C$}{de~$C^\prime$}
par le compactifi\'e de~$\ZZ$. Notons que ce sont des questions
de ce genre qui ont donn\'e naissance \`a la \og th\'eorie de la
descente\fg pour les sch\'emas.

b) Soit $A$ un anneau local complet int\`egre, on sait que son
normalis\'e $A^\prime$ est fini
\marginpar{27}
sur~$A$ (Nagata), donc c'est un anneau semi-local complet, donc local
puisqu'il est int\`egre. Supposons que l'extension r\'esiduelle
$L/k$ qu'il d\'efinit soit non radicielle (dans le cas contraire, on
dira que $A$ est g\'eom\'etriquement unibranche, \cf plus bas). Ce
sera le cas par exemple pour l'anneau
\ifthenelse{\boolean{orig}}
{$\RR[s,t]/(s^2+t^2)\RR[s,t]$,}
{$\RR[[s,t]]/(s^2+t^2)\RR[[s,t]]$,}
o\`u $\RR$ est le corps
des r\'eels. Soit alors $k^\prime$ une extension galoisienne finie
de~$k$ telle que $L \otimes_k k^\prime$ se d\'ecompose; et soit~$B$
une alg\`ebre finie \'etale sur~$A$ correspondant \`a
l'extension r\'esiduelle~$k^\prime$ (rappelons que $B$ est
essentiellement unique). Alors $B^\prime = A^\prime \otimes_A B$
sur~$B$ a l'alg\`ebre r\'esiduelle $L \otimes_k k^\prime$ qui
n'est pas locale, donc $B^\prime$ n'est pas un anneau local donc
(\'etant complet) a des diviseurs de~$0$. Comme il est contenu dans
l'anneau total des fractions de~$B$ (car libre sur~$A^\prime$ donc
sans torsion sur~$A^\prime$ donc sans torsion sur~$A$, donc contenu
dans $B^\prime \otimes_A K = B^\prime_{(K)} = A^\prime_{(K)} \otimes_K
B_{(K)} = B_{(K)}$ puisque $A^\prime_{(K)} = K$) il s'ensuit que $B$
n'est pas int\`egre. Dans le cas de l'anneau
$\RR[s,t]/(s^2+t^2)\RR[s,t]$, prenant $k^\prime/k=
\CC/\RR$, on trouve pour~$B$ l'anneau local de deux droites
s\'ecantes du plan en leur point d'intersection.

Notons d'ailleurs que s'il existe un rev\^etement connexe \'etale
$X$ de~$Y$ int\`egre qui ne soit pas irr\'eductible, alors toute
composante irr\'eductible de~$X$ donne un exemple d'un
rev\^etement non ramifi\'e $X^\prime$ de~$Y$, dominant~$Y$, qui
n'est pas \'etale sur~$Y$. Dans le cas de l'exemple~a), on obtient
ainsi que $C^\prime$ est non ramifi\'e sur~$C$, sans \^etre
\'etale en les deux points $a$ et~$b$ (comme on constate d'ailleurs
directement par inspection des compl\'et\'es des anneaux locaux
de~$x$ et~$a$: du point de vue \og formel\fg, $C^\prime$ au point~$a$
s'identifie \`a un sous-sch\'ema ferm\'e de~$C$ au point~$x$,
savoir l'une des deux \og branches\fg de~$C$ passant par~$x$).

Dans a) et~b), on voit que la non validit\'e des conclusions
de~\Ref{I.9.5}\ptbl (i) et~(ii) est li\'ee directement au fait qu'un
point de~$Y$ \og \'eclate\fg en des points \emph{distincts} du
normalis\'e (\ifthenelse{\boolean{orig}}{dans~b,}{dans~b),}
le fait que l'extension r\'esiduelle soit non radicielle
doit \^etre interpr\'et\'ee g\'eom\'etriquement de cette
fa\c con). De fa\c con pr\'ecise, nous dirons qu'un anneau local
int\`egre~$A$ est \emph{g\'eom\'etriquement unibranche}
\index{geometriquement unibranche@g\'eom\'etriquement unibranche|hyperpage}%
\index{unibranche (geometriquement)@unibranche (g\'eom\'etriquement)|hyperpage}%
si son normalis\'e n'a qu'un seul id\'eal maximal, l'extension
r\'esiduelle correspondante \'etant radicielle; un point~$y$ d'un
pr\'esch\'ema int\`egre est dit g\'eom\'etriquement
unibranche si son anneau local l'est. Exemples: un point normal, un
point de rebroussement ordinaire d'une courbe, \ifthenelse{\boolean{orig}}{etc...}{etc.} Il semble que si
$Y$ admet un point qui
\marginpar{28}
n'est pas unibranche, il existe toujours
\ifthenelse{\boolean{orig}}{de}{un}
rev\^etement \'etale
connexe non irr\'eductible de~$Y$; c'est du moins ce que nous avons
montr\'e dans le cas~b), lorsque $Y$ est le spectre d'un anneau
local complet. On peut montrer par contre que \emph{si tous les points
de~$Y$ sont g\'eom\'etriquement unibranches, alors tout
$Y$-pr\'esch\'ema non ramifi\'e connexe dominant~$Y$ est
\'etale} et irr\'eductible. La d\'emonstration reprend celle
de~\Ref{I.9.5}, en utilisant la g\'en\'eralisation suivante du
th\'eor\`eme~\Ref{I.8.3}, qui sera d\'emontr\'ee plus tard
\`a l'aide de la technique de
descente\footnote{\Cf IX~\Ref{IX.4.10}. Pour une d\'emonstration
plus directe, \cf EGA~IV 18.10.3, utilisant une variante
de~\Ref{I.9.5} pour des anneaux locaux g\'eom\'etriquement
unibranches.}:

\emph{Soit $Y^\prime \to Y$ un morphisme fini, radiciel, surjectif
(\ie ce qu'on pourrait appeler un \og hom\'eomorphisme
universel\fg). Consid\'erons le foncteur $X \mto X \times_Y
Y^\prime = X^\prime$ des $Y$-pr\'esch\'emas dans les
$Y^\prime$-pr\'esch\'emas. Ce foncteur induit une \'equivalence
de la cat\'egorie des $Y$-sch\'emas \'etales avec la
cat\'egorie des $Y^\prime$-sch\'emas \'etales.} On pourra
appliquer par exemple ce r\'esultat dans le cas o\`u $Y^\prime$
est le normalis\'e de~$Y$, $Y$ \'etant suppos\'e unibranche (et
$Y^\prime$ fini sur~$Y$, ce qui est vrai dans tous les cas qu'on
rencontre en pratique), ou au cas d'un $Y^{\prime\prime}$ \og en
sandwich\fg entre~$Y$ et son normalis\'e (qui n'a plus besoin
d'\^etre fini sur~$Y$).

\chapter{Morphismes lisses: g\'en\'eralit\'es,~propri\'et\'es~diff\'erentielles}
\label{II}
\marginpar{29}
Les renvois \`a l'expos\'e~\Ref{I} sont indiqu\'es par~I. On
rappelle que les anneaux sont noeth\'eriens, et les
pr\'esch\'emas localement noeth\'eriens.

\section{G\'en\'eralit\'es}
\label{II.1}

Soit $Y$ un pr\'esch\'ema, soient $t_1,\dots,t_n$ des
ind\'etermin\'ees, on pose
\begin{equation}
\label{eq:II.1.1}
Y[t_1,\dots,t_n]=Y\otimes_\ZZ \ZZ[t_1,\dots,t_n]\quoi.
\end{equation}
Donc $Y[t_1,\dots,t_n]$ est un $Y$-sch\'ema, affine au-dessus
de~$Y$, d\'efini par le faisceau quasi-coh\'erent d'alg\`ebres
$\cal{O}_Y[t_1,\dots,t_n]$. La donn\'ee d'une section de ce
pr\'esch\'ema au-dessus de~$Y$ \'equivaut donc \`a la
donn\'ee de~$n$ sections de~$\cal{O}_Y$ (correspondant aux images
des~$t_i$ par l'homomorphisme correspondant). Si $Y^\prime$ est
au-dessus de~$Y$, on a
\begin{equation}
\label{eq:II.1.2}
Y[t_1,\dots,t_n]\times_Y Y'=Y'[t_1,\dots,t_n]\quoi,
\end{equation}
(ce qui implique que la donn\'ee d'un $Y$-morphisme de~$Y^\prime$
dans $Y[t_1,\dots,t_n]$ \'equivaut \`a la donn\'ee de~$n$
sections de~$\cal{O}_{Y^\prime}$), d'autre part on a
\begin{equation}
\label{eq:II.1.3}
\big(Y[t_1,\dots,t_n]\big)[t_{n+1},\dots,t_m]=Y[t_1,\dots,t_m]\quoi,
\end{equation}
en vertu de la formule analogue pour les anneaux de polyn\^omes
sur~$\ZZ$. La formule~\eqref{eq:II.1.2} implique que
$Y[t_1,\dots,t_n]$ varie fonctoriellement avec~$Y$.

$Y[t_1,\dots,t_n]$ est de type fini et plat au-dessus de~$Y$.

\begin{definition}
\label{II.1.1}
Soit $f\colon X \to Y$ un morphisme, faisant de~$X$ un
$Y$-pr\'esch\'ema. On dit que $f$ est
\emph{lisse}\footnote{Ancienne terminologie: $f$ est \emph{simple}
\index{simple: synonyme d\'esuet de lisse|hyperpage}%
en~$x$, ou $x$ est un point \emph{simple} pour~$f$. Cette
terminologie pr\^etait \`a confusion dans divers contextes
(alg\`ebres simples, groupes simples) et a d\^u \^etre
abandonn\'ee.}
\index{lisse (morphisme, alg\`ebre)|hyperpage}%
en $x \in X$, ou que $X$ est \emph{lisse sur~$Y$ en~$x$}, s'il existe
un entier $n \geq 0$, un voisinage ouvert $U$ de~$x$, et un
$Y$-morphisme \'etale de~$U$ dans $Y[t_1,\dots,t_n]$. On dit
que~$f$
\ifthenelse{\boolean{orig}}{\emph{(\resp $X$)}}{(\resp $X$)}
est \emph{lisse} s'il est lisse en tous les points de~$X$. Une
alg\`ebre~$B$ sur un anneau~$A$ est dite lisse en un id\'eal
premier $\goth{p}$ de~$B$, si $\Spec(B)$ est lisse sur~$\Spec(A)$
\ifthenelse{\boolean{orig}}{\emph{au}}{au}
point~$\goth{p}$;
\marginpar{30}
$B$ est dite lisse sur~$A$ si elle est lisse sur~$A$ en tout id\'eal
premier~$\goth{p}$ de~$B$. Enfin, un homomorphisme local $A \to B$
d'anneaux locaux est dit lisse (ou $B$ est dite lisse sur~$A$)\footnote{Il vaut mieux dire alors, comme dans EGA IV 18.6.1, que $B$
est \og \emph{essentiellement lisse}\fg sur~$A$.}
\index{essentiellement lisse|hyperpage}%
si $B$ est localis\'ee d'une alg\`ebre de type fini $B_1$ lisse
sur~$A$.
\end{definition}

On note que la notion de lissit\'e de~$X$ sur~$Y$ est locale sur~$X$
et sur~$Y$; si $X$ est lisse sur~$Y$, il est localement de type fini
sur~$Y$.

\setcounter{subsection}{0}
\begin{proposition}
\label{prop:II.1.1}
L'ensemble des points $x$ de~$X$ en lesquels $f$ est lisse est ouvert.
\end{proposition}

C'est trivial sur la d\'efinition.

\begin{corollaire}
\label{II.1.2}
Si $B$ est lisse sur~$A$ en~$\goth{p}$, alors il est lisse sur~$A$
en~$\goth{q}$ pour tout id\'eal premier~$\goth{q}$ de~$B$ contenu
dans~$\goth{p}$.
\end{corollaire}

\Ref{prop:II.1.1} implique aussi que les deux derni\`eres
d\'efinitions~\Ref{II.1.1} co\"incident dans leur domaine commun
d'existence.

\begin{proposition}
\label{II.1.3}
\textup{(i)}~Un morphisme \'etale, en particulier une immersion
ouverte, un morphisme identique, est lisse. \textup{(ii)}~Une extension
de la base dans un morphisme lisse donne un morphisme
lisse. \textup{(iii)}~Le compos\'e de deux morphismes lisses est
lisse.
\end{proposition}

(i)~est trivial sur la d\'efinition, on a plus pr\'ecis\'ement:

\begin{corollaire}
\label{II.1.4}
\'etale $=$ quasi-fini $+$ lisse.
\end{corollaire}

(ii)~r\'esulte aussit\^ot du fait analogue pour les morphismes
\'etales (I~\Ref{I.4.6}) et pour les projections $Y[t_1,\dots,t_n]
\to Y$ (\cf \eqref{eq:II.1.2}). Pour (iii), cela r\'esulte
formellement du fait que c'est vrai s\'epar\'ement pour
\og \'etale\fg (I~\Ref{I.4.6}) et des projections du type
\ifthenelse{\boolean{orig}}{$Y[t_1,\dots,t_n]$}
{$Y[t_1,\dots,t_n]\to Y$}
(\cf \eqref{eq:II.1.3}), et des deux faits cit\'es pour~(ii):
Supposons $Y$ lisse sur~$Z$ et $X$ lisse sur~$Y$, prouvons que $X$ est
lisse sur~$Z$; on peut supposer $Y$ \'etale sur~$Z[t_1,\dots,t_n]$
et $X$ \'etale sur~$Y[s_1,\dots,s_m]$, la premi\`ere
hypoth\`ese implique donc que $Y[s_1,\dots,s_m]$ est \'etale sur~$Z[t_1,\dots,t_n][s_1,\dots,s_m] = Z[t_1,\dots,s_m]$, donc $X$ est
\'etale sur~$Z[t_1,\dots,s_m]$, cqfd.

\begin{remarque}
\label{II.1.5}
L'entier $n$ qui figure dans d\'ef\ptbl \Ref{II.1.1} est bien
d\'etermin\'e, car on constate
\marginpar{31}
aussit\^ot que c'est la dimension de l'anneau local de~$x$ dans sa
fibre $f^{-1}\big(f(x)\big)$. On l'appelle \og dimension relative\fg
de~$X$ sur~$Y$. Elle se comporte additivement pour la composition des
morphismes.
\end{remarque}

\section{Quelques crit\`eres de lissit\'e d'un morphisme}
\label{II.2}

\begin{theoreme}
\label{II.2.1}
Soit $f\colon X\to Y$ un morphisme localement de type fini, soit $x\in
X$ et $y=f(x)$. Pour que $f$ soit lisse en $x$, il faut et il suffit
que \textup{(a)} $f$ soit plat en $x$, et \textup{(b)} $f^{-1}(y)$ soit lisse sur~$\kres(y)$
en $x$.
\end{theoreme}

Le compos\'e de deux morphismes plats \'etant plat, et
$Y[t_{1},\dots,t_{n}]\to Y$ \'etant un morphisme plat, on voit que
lisse implique plat; compte tenu de~\Ref{II.1.3} (ii) cela prouve la
n\'ecessit\'e. Supposons (a) et (b) v\'erifi\'ees, soient $V$
un voisinage affine de~$y$ d'anneau~$A$, $U$ un voisinage affine de
$x$ au-dessus de~$V$, d'anneau $B$. Prenant $U$ assez petit, on peut
supposer par (b) qu'il existe un $\kres(y)$-morphisme \emph{\'etale}
$$
g\colon U|f^{-1}(y)\to\Spec k[t_{1},\dots,t_{n}]\quad (k=\kres(y))
$$
d\'efini par $n$ sections $g_{i}$ du faisceau structural de
$U|f^{-1}(y)$. On constate facilement qu'on peut supposer que les
$g_{i}$ (qui a priori sont des \'el\'ements de
$B\otimes_{A}k=BS^{-1}$, o\`u $S=A-\goth{p}$, $\goth{p}$ l'id\'eal
premier de~$A$ correspondant \`a $y$) proviennent de sections du
faisceau structural de
\ifthenelse{\boolean{orig}}{$V$,}{$U$,}
donc que $g$ est induit par un morphisme,
encore not\'e~$g$
$$
\ifthenelse{\boolean{orig}}{g\colon V\to Y[t_{1},\dots,t_{n}]}
{g\colon U\to Y[t_{1},\dots,t_{n}]}
$$
(quitte \`a multiplier les $g_{i}$ par un m\^eme \'el\'ement
non nul de~$k$). Or
\ifthenelse{\boolean{orig}}{$V$}{$U$}
est plat sur~$Y$ par (a), il en est de m\^eme de
$Y[t_{1},\dots,t_{n}]$, d'autre part $g$ induit un morphisme
\'etale entre les fibres au-dessus de~$y$, donc $g$ est \'etale en
$x$ par~(I~\Ref{I.5.8}), cqfd.

\begin{corollaire}
\label{II.2.2}
Soient $S$ un pr\'esch\'ema, $f\colon X\to Y$ un $S$-morphisme de
type fini, $Y$ \'etant de type fini et plat sur~$S$, $x\in X$, $s$
la projection de~$x$ sur~$S$. Pour que $f$ soit lisse en $x$, il faut
et il suffit que $X$ soit plat (ou encore: lisse) sur~$S$ en $x$, et
que le morphisme $f_{s}\colon X_{s}\to Y_{s}$ induit sur les fibres de
$s$ soit lisse en $x$.
\end{corollaire}

Seule la suffisance demande une d\'emonstration, et r\'esulte du
crit\`ere~\Ref{II.2.1}, joint au crit\`ere de
platitude~(I~\Ref{I.5.9}).

Pour
\marginpar{32}
\'enoncer le r\'esultat suivant, \og rappelons\fg qu'un morphisme
$f\colon X\to Y$ localement de type fini est dit
\emph{\'equidimensionnel}
\index{equidimensionnel@\'equidimensionnel (morphisme)|hyperpage}%
en le point $x\in X$ si (posant $y=f(x)$) on peut trouver un voisinage
ouvert $U$ de~$x$, dont toute composante domine une composante de~$Y$
tel que, pour tout $y'\in Y$, les composantes irr\'eductibles de
$f^{-1}(y')\cap U$ aient toutes une m\^eme dimension
ind\'ependante de~$y'$. Il suffit d'ailleurs dans cette condition de
prendre pour $y'$ les points g\'en\'eriques des composantes
irr\'eductibles de~$Y$ passant par $y$, \emph{et} le point $y$. Si
par exemple $X$ et $Y$ sont int\`egres et $f$ dominant, la condition
signifie que les composantes des $f^{-1}(y)$ passant par $x$ ont \og la
bonne\fg dimension, \ie la dimension de la fibre g\'en\'erique
(rappelons qu'elles sont toujours $\ge$ la dimension de la fibre
g\'en\'erique). Si $f$ est \'equidimensionnel en $x$, la
dimension de sa fibre en $x$ \'etant $n$, et si $g\colon U\to
Y'=Y[t_{1},\dots,t_{n}]$ est un $Y$-morphisme d'un voisinage~$U$ de
$x$, induisant un morphisme sur les fibres de~$y$ qui est quasi-fini
en $x$ (ou encore, ce qui revient au m\^eme, si $g$ est quasi-fini
en~$x$), alors on montre que toute composante irr\'eductible de~$U$
passant par $x$ domine une composante irr\'eductible de~$Y'$. D'ailleurs en vertu du \og lemme de normalisation\fg, un tel $g$
existe toujours (et r\'eciproquement, s'il existe un $Y$-morphisme
quasi-fini $g$ d'un voisinage ouvert $U$ de~$x$ dans un $Y$-sch\'ema
de la forme $Y'=Y[t_{1},\dots,t_{n}]$, tel que toute composante de~$U$ passant par $x$ domine une composante de~$Y'$, alors $f$ est
\'equidimensionnel en~$x$). Ceci pos\'e:

\begin{proposition}
\label{II.2.3}
Soient $f\colon X\to Y$ un morphisme localement de type fini, $x$ un
point de~$X$, $y=f(x)$, on suppose $\cal{O}_{y}$ normal. Pour que $f$
soit lisse en $x$, il faut et il suffit que $f$ soit
\'equidimensionnel en $x$, et que $f^{-1}(y)$ soit lisse sur~$\kres(y)$
en $x$.
\end{proposition}

On voit aussit\^ot sur la d\'efinition qu'un morphisme lisse est
\'equidimensionnel (N.B. un morphisme plat de type fini n'est pas
n\'ecessairement \'equidimensionnel en $x$, m\^eme si sa fibre
en $x$ est irr\'eductible). Prouvons la r\'eciproque. Comme
$f^{-1}(y)$ est lisse sur~$\kres(y)$ en $x$, on peut supposer (rempla\c cant au besoin $X$ par un voisinage convenable de~$x$) qu'il existe
un $Y$-morphisme
$$
g\colon X\to Y[t_{1},\dots,t_{n}]=Y'
$$
induisant un morphisme \'etale sur les fibres de~$y$, et a fortiori
quasi-fini en $x$. Donc
\marginpar{33}
$g$ est non ramifi\'e, et ($f$ \'etant \'equidimensionnel en
$x$) les composantes irr\'eductibles de~$X$ passant par $x$ dominent
chacun une composante de~$Y'$, a fortiori l'homomorphisme
$\cal{O}_{y'}\to\cal{O}_{x}$ d\'eduit de~$g$ (o\`u $y'=g(x)$) est
\emph{injectif}. Cet homomorphisme est de plus non ramifi\'e, et
$\cal{O}_{y'}$ est normal puisque localis\'e de l'anneau
$\cal{O}_{y}[t_{1},\dots,t_{n}]$, qui est normal puisque
$\cal{O}_{y}$ l'est. Donc l'homomorphisme $\cal{O}_{y'}\to\cal{O}_{x}$
est \'etale (I~\Ref{I.9.5}~(ii)).

\begin{remarques}
\label{II.2.4}
L'\'enonc\'e pr\'ec\'edent vaut encore en rempla\c cant
l'hypoth\`ese que $\cal{O}_{y}$ est normal par l'hypoth\`ese plus
faible que $Y$ est \emph{g\'eom\'etriquement unibranche} en $y$,
(\cf I~\Ref{I.11}) --- puisque (I~\Ref{I.9.5}) vaut sous cette
hypoth\`ese. Profitons de l'occasion pour signaler en m\^eme temps
que si le corps r\'esiduel d'un anneau local int\`egre $A$ est
alg\'ebriquement clos, alors analytiquement int\`egre
\index{analytiquement int\`egre|hyperpage}%
(\ie $\widehat{A}$ est int\`egre) implique g\'eom\'etriquement
unibranche, la r\'eciproque \'etant vraie de plus dans toute
cat\'egorie de \og bons anneaux\fg, de fa\c con pr\'ecise dans une
cat\'egorie d'anneaux stable par les op\'erations usuelles, et
o\`u la compl\'etion d'un anneau local normal est normale
(condition remplie, en vertu du \og th\'eor\`eme de normalit\'e
analytique\fg de Zariski, dans la cat\'egorie des alg\`ebres
affines et leurs localis\'ees)\footnote{\Cf EGA IV~7.8.}.

\og Rappelons\fg enfin dans le contexte actuel le r\'esultat suivant,
d\^u \`a Hironaka\footnote{\Cf EGA IV 5.12.10.} qui permet parfois
de s'assurer que $f^{-1}(y)$ est un sch\'ema r\'eduit, \ie que
c'est aussi ce que de nombreux g\'eom\`etres alg\'ebristes
consid\'eraient abusivement comme la fibre (sans multiplicit\'e)
de~$f$ au-dessus de~$x$ (savoir $f^{-1}(y)_{\red}$):
\end{remarques}

\begin{proposition}
\label{II.2.5}
Soient $f\colon X\to Y$ un morphisme dominant de type fini de
pr\'esch\'emas r\'eduits, $y$ un point de~$Y$ tel que
$\cal{O}_{y}$ soit r\'egulier. On suppose que toutes les composantes
de~$f^{-1}(y)$ sont de multiplicit\'e~$1$
\index{composante de multiplicit\'e~$1$|hyperpage}%
(\cf d\'efinition plus bas), et que $f^{-1}(y)_{\red}$ est
normal. Alors $f^{-1}(y)$ est r\'eduit donc normal, $X$ est normal
en tous les points de~$f^{-1}(y)$, enfin $X$ est plat sur~$Y$ en tous
les points de~$f^{-1}(y)$.
\end{proposition}
On
\marginpar{34}
dit qu'une composante $Z$ de~$f^{-1}(y)$ est de
\emph{multiplicit\'e}~$1$
\index{multiplicit\'e d'une composante|hyperpage}%
si, $x$ d\'esignant le point g\'en\'erique de~$Z$, on a (i)
$\dim\cal{O}_{x}=\dim\cal{O}_{y}$ (\ie $Z$ n'est pas \og composante
exc\'edentaire\fg, c'est-\`a-dire n'est pas \og de dimension trop
grande\fg; (ii) l'id\'eal maximal de~$\cal{O}_{x}$ est engendr\'e
par l'id\'eal maximal de~$\cal{O}_{y}$ (qui a priori, en vertu du
choix de~$x$, engendre un id\'eal de d\'efinition
de~$\cal{O}_{x}$).

Compte tenu de~\Ref{II.2.3} ou de~\Ref{II.2.1} on trouve donc:

\begin{corollaire}
\label{II.2.6}
Soient $f\colon X\to Y$ un morphisme dominant de type fini de
pr\'esch\'emas r\'eduits, $y$ un point de~$Y$ tel que
$\cal{O}_{y}$ soit r\'egulier. Pour que $f$ soit lisse aux points de
$X$ au-dessus de~$y$, il faut et il suffit que les composantes de
$f^{-1}(y)$ soient de multiplicit\'es $1$, et que $f^{-1}(y)_{\red}$
soit lisse sur~$\kres(y)$.
\end{corollaire}

Cette situation \'etait surtout consid\'er\'ee par le pass\'e
quand $Y$ \'etait le spectre d'un anneau de valuation discr\`ete
$A$, et \'etait d\'esign\'ee commun\'ement sous des vocables
tels que: \og si la r\'eduction de~$X$ par rapport \`a la valuation
donn\'ee est jolie\fg... De plus, $X$ d\'esignait alors un
sous-sch\'ema (si on peut dire) ferm\'e d'un $\PP_{K}^{n}$ ($K$
\'etant le corps des fractions de~$A$) et faute d'un
\ifthenelse{\boolean{orig}}{language}{langage}
ad\'equat, le r\^ole plus intrins\`eque d'un objet \og d\'efini
sur~$A$\fg (et non seulement sur~$K$) n'apparaissait gu\`ere.

\section{Propri\'et\'es de permanence}
\label{II.3}

\begin{proposition}
\label{II.3.1}
Soit $f\colon X\to Y$ un morphisme, soit $x\in X$ et
$y=f(x)$. Supposons $f$ lisse en~$x$. Pour que $\cal{O}_{x}$ soit
r\'eduit (\resp r\'egulier, \resp normal) il faut et il suffit que
$\cal{O}_{y}$ le soit.
\end{proposition}

Cet \'enonc\'e est en effet connu quand $X$ est de la forme
$Y[t_{1},\dots,t_{n}]$, et il est d\'emontr\'e dans (I, \no \Ref{I.9})
pour un morphisme \'etale; le cas g\'en\'eral s'en d\'eduit
aussit\^ot gr\^ace \`a la d\'efinition~\Ref{II.1.1}.

Nous ne d\'etaillons pas ici les autres propri\'et\'es de
permanence, r\'esultant d\'ej\`a de la seule platitude, ou du
fait que $X$ est localement quasi-fini et plat au-dessus d'un
$Y$-pr\'esch\'ema de la forme $Y[t_{1},\dots,t_{n}]$ (ou, comme
\marginpar{35}
nous dirons, que $X$ est
\ifthenelse{\boolean{orig}}{Cohen-Macauley}{Cohen-Macaulay}
%
\index{Cohen-Macaulay (sch\'ema de)|hyperpage}%
au-dessus de~$Y$). Signalons seulement que de ce dernier fait
r\'esulte que
\begin{equation}
\label{eq:II.3.1}
{\dim\cal{O}_{x}=\dim\cal{O}_{y}+n-d,\quad
\prof\cal{O}_{x}=\prof\cal{O}_{y}+n-d},
\end{equation}
o\`u $n$ est la dimension de la fibre de~$f$ en $x$, et $d$ le
degr\'e de transcendance de~$\kres(x)$ sur~$\kres(y)$, d'o\`u (posant
$\coprof=\dim-\prof$)\footnote{Pour ces formules, \cf EGA~IV 6.1
et~6.3.}
\begin{equation}
\label{eq:II.3.2}
{\coprof\cal{O}_{x}=\coprof\cal{O}_{y}}.
\end{equation}
Il en r\'esulte par exemple que $\cal{O}_{x}$ est Cohen-Macaulay
(\resp sans composantes immerg\'ees) si et seulement si il en est de
m\^eme de~$\cal{O}_{y}$.

\section{Propri\'et\'es diff\'erentielles des morphismes lisses}
\label{II.4}
Pour simplifier, nous nous restreindrons pour l'essentiel au calcul
diff\'erentiel d'ordre~$1$, nous bornant \`a de rapides
indications pour l'ordre sup\'erieur (o\`u les r\'esultats sont
tout aussi simples).

Pour la d\'efinition du faisceau des $1$-diff\'erentielles
$\it{\Omega}^1_{X/Y}$ d'un $Y$-pr\'esch\'ema $X$,
\cf (I~\No\Ref{I.1}). Supposons que $X$ et $Y$ soient des
$S$-pr\'esch\'emas, le morphisme structural $f\colon X\to Y$
\'etant un $S$-morphisme. Alors $f$ d\'efinit un homomorphisme de
Modules (compatible avec $f$)
\begin{equation}
\label{eq:II.4.1}
{f^*\colon \it{\Omega}^1_{Y/S}\to\it{\Omega}^1_{X/S}}
\end{equation}
en d'autres termes, $\it{\Omega}^1_{X/S}$ est \emph{contravariant} en
le $S$-pr\'esch\'ema $X$. D'ailleurs~\eqref{eq:II.4.1}
\'equivaut \`a un homomorphisme de Modules sur~$X$
\begin{equation*}
\label{eq:II.4.1bis}
\tag{4.1\textrm{bis}}
{f^*\big(\it{\Omega}^1_{Y/S}\big)\to\it{\Omega}^1_{X/S}}
\end{equation*}
\'egalement d\'enot\'e par $f^*$ \`a d\'efaut de mieux, et
qui s'ins\`ere dans une suite exacte canonique d'homomorphismes de
Modules
\begin{equation}
\label{eq:II.4.2}
{f^*\big(\it{\Omega}^1_{Y/S}\big)\to
\it{\Omega}^1_{X/S}\to\it{\Omega}^1_{X/Y}\to 0}
\end{equation}
Tous ces homomorphismes sont d\'efinis par la condition d'\^etre de
nature locale (ce qui ram\`ene au cas affine) et de commuter avec
les op\'erateurs $\rd$. L'exactitude de~\eqref{eq:II.4.2} est
classique et triviale, et se transcrit dans le cas affine en la suite
exacte (correspondant \`a un homomorphisme $B\to C$ de
$A$-alg\`ebres):
\begin{equation*}
\label{eq:II.4.2bis}
\tag{4.2\textrm{bis}}
{\Omega^1_{B/A}\otimes_{B}C\to\Omega^1_{C/A}\to\Omega^1_{C/B}\to 0}
\end{equation*}

\begin{lemme}
\label{II.4.1}
Soit
\marginpar{36}
$f\colon X\to Y$ un morphisme de~$S$-pr\'esch\'emas. Si $f$
est non ramifi\'e (\resp \'etale) alors
$f^*\big(\it{\Omega}^1_{Y/S}\big)\to\it{\Omega}^1_{X/S}$
est surjectif (\resp un isomorphisme). La r\'eciproque est vraie
dans le cas \og non ramifi\'e\fg, si $f$ est suppos\'e localement de
type fini.
\end{lemme}

Le cas non ramifi\'e r\'esulte de la suite
exacte~\eqref{eq:II.4.2} et de (I~\Ref{I.3.1}), mais peut aussi se
voir directement comme le cas \'etale. Consid\'erons le diagramme
\begin{equation*}
\xymatrix@C=1.5cm{
X \ar[r]^-{\Delta_{X/Y}} & X\times_{Y}X \ar[d] \ar[r] &
X\times_{S}X \ar[d] \\ & Y \ar[r]^-{\Delta_{Y/S}} & Y\times_{S}Y
}
\end{equation*}
dans lequel $X\times_{Y}X$ s'identifie au produit fibr\'e de~$Y$ et
$X\times_{S}X$ sur~$Y\times_{S}Y$. Comme $f$ est non ramifi\'e,
$X\to X\times_{Y}X$ est une immersion ouverte, donc le faisceau
\og conormal\fg de l'immersion compos\'ee $\Delta_{X/S}$ de cette
derni\`ere avec $X\times_{Y}X\to X\times_{S}X$ est isomorphe \`a
l'image inverse sur~$X$ du faisceau conormal pour l'immersion
$X\times_{Y}X\to X\times_{S}X$. D'autre part, $X\to Y$ \'etant
\'etale donc plat, $X\times_{S}X\to Y\times_{S}Y$ est plat, donc le
faisceau conormal pour l'immersion $X\times_{Y}X\to X\times_{S}X$ est
isomorphe \`a l'image inverse du faisceau conormal pour l'immersion
$Y\to Y\times_{S}Y$, \ie l'image inverse de~$\it{\Omega}^1_{Y/S}$. La
conclusion en r\'esulte.

\begin{lemme}
\label{II.4.2}
Soit $X=Y[t_{1},\dots,t_{n}]$, $Y$ \'etant un
$S$-pr\'esch\'ema. Alors la suite d'homomorphismes canoniques
$$
0\to f^*\big(\it{\Omega}^1_{Y/S}\big)\to
\it{\Omega}^1_{X/S}\to\it{\Omega}^1_{X/Y}\to 0
$$
est exacte et $\it{\Omega}^1_{X/Y}$ est libre de base les
$\rd_{X/Y}t_{i}$.
\end{lemme}

La v\'erification (purement affine) est imm\'ediate. (N.B. on
conna\^it d\'ej\`a l'exactitude de~\eqref{eq:II.4.2}).

Combinant ces deux \'enonc\'es et d\'efinition~\Ref{II.1.1}, on
trouve

\begin{theoreme}
\label{II.4.3}
Soient $f\colon X\to Y$ un morphisme lisse de~$S$-pr\'esch\'emas,
alors:
\begin{enumerate}
\item[(i)] La suite d'homomorphismes canoniques
$$
0\to f^*\big(\it{\Omega}^1_{Y/S}\big)\to
\it{\Omega}^1_{X/S}\to\it{\Omega}^1_{X/Y}\to 0
$$
est exacte.
\item[(ii)] $\it{\Omega}^1_{X/Y}$ est localement libre, son rang $n$
en $x$ est \'egal \`a la dimension relative de~$f$ en~$x$.
\end{enumerate}
\end{theoreme}

\begin{corollaire}
\label{II.4.4}
L'homomorphisme
\marginpar{37}
$f^*\big(\it{\Omega}^1_{Y/S}\big)\to\it{\Omega}^1_{X/S}$ est
injectif, son image dans $\it{\Omega}^1_{X/S}$ est localement facteur
direct.
\end{corollaire}

Soit $u\colon F\to G$ un homomorphisme de Modules sur le
pr\'esch\'ema $X$, on dit qu'il est \emph{universellement
injectif}
\index{universellement injectif (morphisme de modules)|hyperpage}%
en $x\in X$, si l'homomorphisme $F_{x}\to G_{x}$ de
$\cal{O}_{x}$-modules est injectif, et reste tel par tensorisation
avec toute $\cal{O}_{x}$-alg\`ebre (ou, ce qui revient au m\^eme
d'ailleurs, avec tout $\cal{O}_{x}$-module). Il suffit par exemple
qu'il existe un voisinage ouvert $U$ de~$x$ tel que $u$ induise un
isomorphisme de~$F|U$ sur un facteur direct de~$G|U$, cette condition
est aussi n\'ecessaire lorsque $F$ et $G$ sont libres (et de type
fini) dans un voisinage de~$x$, de fa\c con pr\'ecise dans ce cas
les conditions suivantes sont \'equivalentes:
\begin{enumerate}
\item[(i)] $u$ est injectif en $x$ et $\Coker u$ libre en $x$;
\item[(ii)] Il existe un voisinage ouvert $U$ de~$x$ tel que $u$
induise un isomorphisme de~$F|U$ sur un facteur direct de~$G|U$;
\item[(iii)] $u$ est universellement injectif en $x$;
\item[(iv)] l'homomorphisme
\ifthenelse{\boolean{orig}}
{$\cal{F}_{x}\otimes \kres(x)\to\cal{G}_{x}\otimes \kres(x)$}
{$F_{x}\otimes \kres(x)\to G_{x}\otimes \kres(x)$}
sur les \og fibres\fg restreintes induit par $u$ est injectif;
\item[(v)] L'homomorphisme transpos\'e $\check{G}\to\check{F}$ est
surjectif au point $x$ (ou encore, ce qui revient au m\^eme, au
voisinage de~$x$).
\end{enumerate}
(D\'emonstration circulaire, (iv)$\To$(v) r\'esulte de
Nakayama, d'autre part (v)$\To$(i)
\ifthenelse{\boolean{orig}}
{puisque un}
{puisqu'un}
faisceau
quotient localement libre est n\'ecessairement facteur
direct). G\'eom\'etriquement, la situation envisag\'ee signifie
que $u$ correspond \`a un isomorphisme du fibr\'e vectoriel dont
le faisceau des sections
\ifthenelse{\boolean{orig}}{est~$\cal{F}$,}{est~$F$,}
sur un sous-fibr\'e du
fibr\'e vectoriel analogue d\'efini par $G$. Bien entendu, il ne
suffit pas pour cela que
\ifthenelse{\boolean{orig}}{$\cal{F}\to\cal{G}$}{$F\to G$}
soit injectif.

\begin{corollaire}
\label{II.4.5}
Soit $f\colon X\to Y$ un morphisme de~$S$-pr\'esch\'emas,
localement de type fini, $x\in X$, $y=f(x)$, $s$ la projection de~$x$
et $y$ sur~$S$. On suppose $Y$ lisse en $y$ sur~$S$. Conditions
\'equivalentes:
\begin{enumerate}
\item[(i)] $f$ est lisse en $x$.
\item[(ii)] $X$ est lisse sur
\ifthenelse{\boolean{orig}}{$Y$}{$S$}
en~$x$, et
$f^*\big(\it{\Omega}^1_{Y/S}\big)\to\it{\Omega}^1_{X/S}$ est
universellement injectif en~$x$, \ie c'est un homomorphisme injectif en~$x$ et son conoyau $\it{\Omega}^1_{X/Y}$ est libre
\ifthenelse{\boolean{orig}}{en~$x$).}{en~$x$.}
\end{enumerate}
\end{corollaire}

La n\'ecessit\'e r\'esulte de~\Ref{II.1.3} (iii) et
de~\Ref{II.4.3} (i)\,(ii), prouvons la suffisance. Comme les
$\rd g$ ($g\in\cal{O}_{x}$) engendrent le module
$\it{\Omega}^1_{X/Y}$ en $x$, on peut trouver des
\marginpar{38}
$g_{i}$ ($1\le i\le n$) tels que les images des $\rd g_{i}$ dans
$\big(\it{\Omega}^1_{X/Y}\big)_{x}$ forment une base de ce
module. Prenant $X$ assez petit, on peut supposer que les $g_{i}$
proviennent de sections de~$\cal{O}_{X}$, et d\'efinissent donc un
$Y$-morphisme $g\colon X\to Y'=Y[t_{1},\dots,t_{n}]$. Utilisant
l'hypoth\`ese et lemme~\Ref{II.4.2}, on voit facilement que
l'homomorphisme correspondant
$g^*\big(\it{\Omega}^1_{Y'/S}\big)\to\it{\Omega}^1_{X/S}$ est
bijectif en $x$, ce qui nous ram\`ene \`a prouver le

\begin{corollaire}
\label{II.4.6}
Soit $f\colon X\to Y$ un morphisme de~$S$-pr\'esch\'emas
lisses. Pour que $f$ soit \'etale en $x\in X$, il faut et il suffit
que $f^*\big(\it{\Omega}^1_{Y/S}\big)\to\it{\Omega}^1_{X/S}$ soit
un isomorphisme en $x$.
\end{corollaire}

On sait que c'est n\'ecessaire par~\Ref{II.4.1}, et cette condition
implique que $f$ est non ramifi\'e en $x$ par le m\^eme lemme. En
vertu de~\Ref{II.2.2}, on est ramen\'e au cas o\`u
$S=\Spec(k)$. Comme $Y$ est lisse sur~$k$, il est r\'egulier, donc a
fortiori normal, et en vertu de (I~\Ref{I.9.5} (ii)) on est ramen\'e
\`a prouver que $\cal{O}_{y}\to\cal{O}_{x}$ est injectif, ou encore
que $\cal{O}_{y}$ et $\cal{O}_{x}$ ont m\^eme dimension. Or ces
dimensions sont respectivement les rangs de~$\it{\Omega}^1_{Y/k}$ et
$\it{\Omega}^1_{X/k}$ en $y$ \resp $x$, donc \'egaux en vertu de
l'hypoth\`ese.
\begin{remarques}
\label{II.4.7}
$X$ et $Y$ \'etant suppos\'es lisses sur~$S$, le
crit\`ere~\Ref{II.4.5} (ii) de lissit\'e de~$f\colon X\to Y$ peut
encore s'\'enoncer en disant que pour tout $x\in X$, l'application
\emph{tangente} (relativement \`a la base $S$) de~$f$ en $x$,
\ie la transpos\'ee de l'homomorphisme des $\kres(x)$ espaces vectoriels
de dimension finie, fibres restreintes de
$f^*\big(\it{\Omega}^1_{Y/S}\big)$ et $\it{\Omega}^1_{X/S}$ en
$x$, est \emph{surjective}. C'est l\`a une hypoth\`ese bien
famili\`ere en particulier parmi les gens travaillant avec les
espaces analytiques. L'hypoth\`ese de non singularit\'e qu'ils
font d'ordinaire (qui signifie que $X$ et $Y$ sont \og lisses sur~$\CC$\fg, \cf \No \Ref{II.5}) ne semble due qu'\`a la peur qu'inspirent
encore \`a bien des g\'eom\`etres les points singuliers des
vari\'et\'es alg\'ebriques ou espaces analytiques.
\end{remarques}

Signalons le cas particulier suivant de~\Ref{II.4.6}:
\begin{corollaire}
\label{II.4.8}
Soient $X$ un $S$-pr\'esch\'ema, $g\colon X\to
S[t_{1},\dots,t_{n}]$ un $S$-morphisme, d\'efini par les sections
$g_{i}$ ($1\le i\le n$) de~$\cal{O}_{X}$, $x$ un point de~$X$ tel que
$X$ soit lisse sur~$S$ en $x$. Pour que $g$ soit \'etale en $x$, il
faut et il suffit que les $\rd g_{i}$
\marginpar{39}
($1\le i\le n$) forment une base de~$\it{\Omega}^1_{X/S}$ en $x$ (ou,
ce qui revient au m\^eme, que leurs images dans
$\it{\Omega}^1_{X/S}(x) =
\big(\it{\Omega}^1_{X/S}\big)_{x}\otimes_{\cal{O}_{x}}\kres(x)$ forment
une base de cet espace vectoriel sur~$\kres(x)$).
\end{corollaire}

Soient $X$ un pr\'esch\'ema, $Y$ un sous-pr\'esch\'ema
ferm\'e de~$X$ d\'efini par un faisceau coh\'erent $\cal{J}$
d'id\'eaux. Donc $\cal{J}/\cal{J}^2$ peut \^etre consid\'er\'e
comme un faisceau coh\'erent sur~$Y$ (le \emph{faisceau conormal}
\index{conormal (faisceau)|hyperpage}%
\index{faisceau conormal|hyperpage}%
de~$Y$ dans~$X$). Si maintenant $X$ est un $S$-pr\'esch\'ema, on a
une suite exacte canonique de faisceaux quasi-coh\'erents sur~$Y$
\begin{equation}
\label{eq:II.4.3}
{\cal{J}/\cal{J}^2\lto{\rd}
\it{\Omega}^1_{X/S}\otimes_{\cal{O}_{X}}\cal{O}_{Y}\to
\it{\Omega}^1_{Y/S}\to 0}
\end{equation}
dont la partie de droite n'est autre que~\eqref{eq:II.4.2} (avec le
r\^ole de~$X$ et $Y$ interchang\'es, compte tenu que
$\it{\Omega}^1_{Y/X}=0$), tandis que l'homomorphisme
$\cal{J}/\cal{J}^2\to\it{\Omega}^1_{X/S}\otimes_{\cal{O}_{X}}\cal{O}_{Y}$
est d\'eduit de l'homomorphisme (en g\'en\'eral non
lin\'eaire) $g\to\rd g$ par passage aux
quotients. L'exactitude de~\eqref{eq:II.4.3} est classique et
d'ailleurs triviale, et s'interpr\`ete dans le cas affine par la
suite exacte suivante (correspondante \`a un homomorphisme surjectif
$B\to C$ de~$A$-alg\`ebres, de noyau $J$):
\begin{equation*}
\label{eq:II.4.3bis}
\tag{4.3\textrm{bis}}
{J/J^2\to\Omega^1_{B/A}\otimes_{B}C\to\Omega^1_{C/A}\to 0\qquad
(C=B/J)}
\end{equation*}
(suite exacte qui avait d\'ej\`a \'et\'e utilis\'ee
implicitement dans la d\'emonstration de (I~\Ref{I.7.2})!).

\begin{proposition}
\label{II.4.9}
Soient $X$ un $S$-pr\'esch\'ema, $Y$ un sous-pr\'esch\'ema
ferm\'e de~$X$ d\'efini par un faisceau coh\'erent $\cal{J}$
d'id\'eaux sur~$X$, $x$ un point de~$X$, $g_i$ ($1 \leq i \leq n$)
des sections de~$\cal{O}_X$, d\'efinissant un $S$-morphisme
$$
g \colon X \to S[t_1,\dots,t_n] = X'
$$
enfin $p$ un entier, $0 \leq p \leq n$. On suppose $X$ \emph{lisse sur~$S$ en $x$}. Les conditions suivantes sont \'equivalentes:
\begin{enumerate}
\item[(i)] Il existe un voisinage ouvert $X_1$ de~$x$ dans $X$ tel que
$g|X_1$ soit \emph{\'etale} et que $Y_1=Y \cap X_1$ (trace de~$Y$
sur~$X_1$) soit \emph{l'image inverse} du sous-pr\'esch\'ema
ferm\'e $Y'=S[t_{p+1},\dots,t_n]$
\ifthenelse{\boolean{orig}}
{\emph{de}}
{de}
$X'=S[t_1,\dots,t_n]$ (\ie les $g_i$ ($1 \leq i \leq p$) engendrent
$\cal{J}|X_1$):
$$
\xymatrix{ Y_1 \ar[d] \ar[r] & X_1 \ar[d]^{\hbox{\'etale}} \\
Y'=S[t_{p+1},\dots, t_n] \ar[r] & X'=S[t_1,\dots,t_n]}
$$
\item[(ii)]
$Y$
\marginpar{40}
est \emph{lisse sur~$S$ en $x$}, les $g_i$ ($1 \leq i \leq p$)
d\'efinissent des \emph{\'el\'ements de} $\cal{J}_x$, les $d
g_i(x)$ ($1 \leq i \leq n$) forment une \emph{base de}
$\it{\Omega}^1_{X/S}(x)$ sur~$\kres(x)$, les $dg'_i(x)$ ($p+1 \leq i \leq
n$) forment une \emph{base de} $\it{\Omega}^1_{Y/S}(x)$ sur~$\kres(x)$
(o\`u les $g'_i$ d\'esignent les restrictions des $g_i$ \`a~$Y$;
les diff\'erentielles sont prises par rapport \`a~$S$).
\item[(iii)] Les $g_i$ ($1 \leq i \leq p$) d\'efinissent un
\emph{syst\`eme de g\'en\'erateurs} de~$\cal{J}_x$, et les
$dg_i(x)$ ($1 \leq i \leq n$) forment une \emph{base de}
$\it{\Omega}^1_{X/S}(x)$ sur~$\kres(x)$.
\item[(iv)] $Y$ est \emph{ lisse sur~$S$} en $x$, les $g_i$
\ifthenelse{\boolean{orig}}{\ignorespaces}{($1 \leq i \leq p$)}
forment un \emph{syst\`eme minimal de g\'en\'erateurs de}
$\cal{J}_x$, les $dg'_i(x)$ ($p+1 \leq i \leq n$) forment une
\emph{base de} $\it \Omega^1_{Y/S}(x)$ sur~$\kres(x)$.
\end{enumerate}

De plus, sous ces conditions, $\cal{J}/\cal{J}^2$ est un Module libre
sur~$Y$ en $x$, admettant comme \emph{base en $x$} les classes des
$g_i$ ($1 \leq i \leq p$), \emph{et l'homomorphisme} canonique
$\cal{J}/\cal{J}^2 \to \it{\Omega}^1_{X/S} \otimes \cal{O}_Y$ est
\emph{universellement injectif en $x$}.
\end{proposition}

\begin{remarquestar}
Cela implique que $p$ est bien d\'etermin\'e par les autres
conditions, soit comme \emph{rang} du Module libre
$\cal{J}/\cal{J}^2$ sur~$Y$ en $x$, ou encore le \emph{nombre minimum
de g\'en\'erateurs} de~$\cal{J}_x$ sur~$X$, ou enfin par le fait
que la dimension relative de~$Y$ rel. $S$ en $x$ est $n-p$.
\end{remarquestar}

\subsubsection*{D\'emonstration}
Supposons d'abord (i) v\'erifi\'e. Alors par (I~\Ref{I.4.6}
(iii)) $Y_1$ est \'etale sur~$Y'$, donc par d\'efinition il est
lisse sur~$S$ en $x$ (de dimension relative $n-p$), il en est donc de
m\^eme de~$Y$. Il r\'esulte alors de \eqref{II.4.8} que les $dg_i$
($1 \leq i \leq n$) forment une base de~$\it{\Omega}^1_{X/S}$ en $x$,
et que les $dg'_i$ ($p+1 \leq i \leq n$) une base de
$\it{\Omega}^1_{Y/S}$ en $x$, d'o\`u il r\'esulte par la suite
exacte \eqref{eq:II.4.3} que les $g_i$ ($1 \leq i \leq p$) sont
lin\'eairement ind\'ependants dans $\cal{J}/\cal{J}^2$
(consid\'er\'e comme Module sur~$Y$) en $x$; comme les $g_i$ ($1
\leq i \leq p$) engendrent $\cal{J}_x$, il s'ensuit que les $g_i \mod
\cal{J}_x^2$ forment une \emph{base de} $\cal{J}/\cal{J}^2$ en
$x$. Cela implique d'une part que les $g_i$ ($1 \leq i \leq p$)
forment un syst\`eme \emph{minimal} de g\'en\'erateurs de
$\cal{J}_x$, d'autre part que l'homomorphisme $\cal{J}/\cal{J}^2 \to
\it{\Omega}^1_{X/S} \otimes \cal{O}_Y$ de \eqref{eq:II.4.3} est
universellement injectif en $x$ (car applique une base d'un Module
libre en $x$ sur une partie d'une base d'un Module libre en $x$ --- N.B. il s'agit de~$Y$-Modules). Cela prouve que (i) implique (ii),
(iii), (iv), ainsi que les derni\`eres assertions de
proposition~\Ref{I.4.9}.

(iii) implique (i) en vertu de corollaire~\Ref{I.4.8}.

(ii)
\marginpar{41}
implique (i). En effet, la premi\`ere hypoth\`ese dans (ii)
signifie que (quitte \`a remplacer $X$ par un voisinage ouvert de
$x$ dans $X$) $g$ induit un morphisme $h\colon Y \to
Y'$. D'apr\`es~\Ref{I.4.8}, les deux autres hypoth\`eses de (ii)
signifient que $g$ est \'etale en $x$, et $h$ \'etale en $x$. Soit
alors $Y''$ l'image inverse de~$Y'$ par $g$. Donc $Y$ est un
sous-pr\'esch\'ema ferm\'e de~$Y''$, qui est \'etale sur~$Y'$
en $x$ par (I~\Ref{I.4.6}~(iii)) puisque $g$ est \'etale en
$x$. Donc le morphisme d'immersion $Y \to Y''$ est lui-m\^eme
\'etale (I~\Ref{I.4.8}) donc une immersion ouverte (I~\Ref{I.5.8} ou
I~\Ref{I.5.2}), donc rempla\c cant encore $X$ par un voisinage
ouvert convenable $X_1$ de~$x$, on obtient (i).

Ce qui pr\'ec\`ede \'etablit l'\'equivalence des conditions
\ifthenelse{\boolean{orig}}
{(i) (ii) (iii)}
{(i), (ii), (iii)}, et le fait qu'elles impliquent (iv), il reste \`a
prouver que (iv)$\To$(ii), ce qui est imm\'ediat (compte
tenu que $\it{\Omega}^1_{X/S}$ est libre sur~$X$ en $x$) une fois
qu'on sait que le fait que $Y$ est lisse sur~$S$ en $x$ implique que
$\cal{J}/\cal{J}^2$ est libre sur~$Y$ en $x$, et l'homomorphisme
$\cal{J}/\cal{J}^2 \to \it{\Omega}^1_{X/S}\otimes \cal{O}_Y$
universellement injectif en $x$. Ce dernier point est inclus dans le
\begin{theoreme}
\label{II.4.10}
Soient $X$ un $S$-pr\'esch\'ema lisse, $Y$ un
sous-pr\'esch\'ema ferm\'e de~$X$ d\'efini par un faisceau
coh\'erent $\cal{J}$ d'id\'eaux sur~$X$, $x$ un point de~$X$. Les
conditions suivantes sont \'equivalentes:
\begin{enumerate}
\item[(i)] $Y$ est \emph{lisse sur~$S$ en $x$}.
\item[(ii)] Il existe un voisinage ouvert $X_1$ de~$x$ dans $X$ et un
$S$-morphisme \emph{\'etale}
$$
g\colon X_1 \to X'=S[t_1,\dots,t_n]
$$
tel que $Y_1=Y \cap X_1$ (trace de~$Y$ sur~$X_1$) soit le
sous-pr\'esch\'ema de~$X_1$ \emph{image inverse} par $g$ du
sous-pr\'esch\'ema ferm\'e $Y'=S[t_{p+1},\dots,t_n]$ de
$X'=S[t_1,\dots,t_n]$, pour un $p$ convenable.
\item[(iii)] Il existe des \emph{g\'en\'erateurs $g_i$ ($1 \leq i
\leq p$) de~$\cal{J}_x$}, tels que les $dg_i$ forment une partie d'une
base de~$\it{\Omega}^1_{X/S}$ en $x$ (ou, ce qui revient au m\^eme,
tel que les $dg_i(x)$ dans $\it{\Omega}^1_{X/S}$ soient
lin\'eairement ind\'ependant sur~$\kres(x)$).
\item[(iv)] Le faisceau $\cal{J}/\cal{J}^2$ est libre sur~$Y$ en $x$,
et l'homomorphisme canonique
$$
d\colon \cal{J}/\cal{J}^2 \to \it{\Omega}^1_{X/S}\otimes \cal{O}_Y
$$
est universellement injectif en~$x$; ou encore: la suite
d'homomorphismes canoniques
$$
0 \to \cal{J}/\cal{J}^2 \to \it{\Omega}^1_{X/S}\otimes \cal{O}_Y \to
\it{\Omega}^1_{Y/S} \to 0
$$
est exacte en $x$, et $\it{\Omega}^1_{Y/S}$ est localement libre
en~$x$.
\end{enumerate}
\end{theoreme}

\subsubsection*{D\'emonstration}
On
\marginpar{42}
sait d\'ej\`a que (ii) implique (i), (iii), (iv) d'apr\`es ce qui
pr\'ec\`ede. Prouvons que (i)$\To$(ii) (ce qui
ach\`evera en m\^eme temps la d\'emonstration
\ifthenelse{\boolean{orig}}
{de~\Ref{I.4.9})}
{de~\Ref{I.4.9}).}
En vertu de th\'eor\`eme \Ref{II.4.3}~(ii), les deux derniers
termes dans la suite exacte \eqref{eq:II.4.3} sont des Modules libres
sur~$Y$. Donc, comme les images dans
$\it{\Omega}^1_{X/S}\otimes_{\cal{O}_X} \cal{O}_Y$ des $dg$ ($g \in
\cal{O}_X$) engendrent ce module en $x$, donc leurs images dans
$\it{\Omega}^1_{Y/S}$ engendrent ce dernier en~$x$, on peut trouver
des $g_i$ ($p+1 \leq i \leq n$) dans $\cal{O}_X$ tels que les $dg'_i$
forment une base de~$\it{\Omega}^1_{Y/S}$, puis (en vertu de
l'exactitude de \eqref{eq:II.4.3}) compl\'eter le syst\`eme des
$dg_i$ ($p+1 \leq i \leq n$) en une base du terme m\'edian par des
\'el\'ements de la forme $dg_i$ ($1 \leq i \leq n$) o\`u les
$g_i$ ($1 \leq i \leq p$) \emph{sont dans $\cal{J}_x$.} Les $g_i$
proviennent de sections de~$\cal{O}_X$ sur un voisinage de~$x$ dans
$X$, qu'on peut supposer \'egal \`a $X$. On est alors sous les
conditions de~\Ref{II.4.8}~(ii), et on a \'etabli que cela implique
la condition~\Ref{II.4.8}~(i), d'o\`u~\Ref{II.4.10}~(ii).

L'implication (iii)$\To$(ii) dans~\Ref{II.4.10} r\'esulte
aussit\^ot de l'implication (iii)$\To$(i)
dans~\Ref{II.4.8}. Donc
\ifthenelse{\boolean{orig}}
{(i) (ii) (iii)}
{(i), (ii), (iii)}
sont \'equivalents, et
impliquent (iv). Enfin, il est trivial que (iv) implique (iii), compte
tenu que des $g_i \in \cal{J}_x$ qui forment une base de~$\cal{J}_x
\mod \cal{J}_x^2$ engendrent $\cal{J}_x$ (Nakayama).

De plus, la d\'emonstration qui pr\'ec\`ede montre ceci:
\begin{corollaire}
\label{II.4.11}
Soient $X$ un $S$-pr\'esch\'ema, $Y$ un sous-pr\'esch\'ema
ferm\'e d\'efini par un faisceau coh\'erent $\cal{J}$
d'id\'eaux sur~$X$, $x$ un point de~$Y$. On suppose \emph{$X$ et $Y$
lisses sur~$S$ en $x$}. Soient $g_i$ des sections de~$\cal{J}$ ($1
\leq i \leq p$). Les conditions suivantes sont \'equivalentes:
\begin{enumerate}
\item[(i)] Les $g_i$ \emph{engendrent} $\cal{J}_x$ et les $dg_i(x)$
sont \emph{lin\'eairement ind\'ependants} dans
$\it{\Omega}^1_{X/S}(x)$
\ifthenelse{\boolean{orig}}{sur~$(x)$.}{sur~$\kres(x)$.}
\item[(ii)] Les $g_i \mod \cal{J}^2$ forment une base de
$\cal{J}/\cal{J}^2$ en $x$.
\item[(iii)] Les $g_i$ forment un syst\`eme minimal de
g\'en\'erateurs de~$\cal{J}_x$.
\item[(iv)] On peut trouver d'autres sections $g_i$ ($p+1 \leq i \leq
n$) de~$\cal{O}_X$ sur un voisinage $X_1$ \emph{de} $X$,
d\'efinissant avec les pr\'ec\'edents un morphisme
\emph{\'etale} $X_1 \to X'=S[t_1,\dots,t_n]$ tel que $Y_1=Y \cap
X_1$ soit \emph{l'image inverse} par $g$ du sous-pr\'esch\'ema
ferm\'e $Y'=S[t_{p+1},\dots,t_n]$ de~$X'=S[t_{1},\dots,t_n]$.
\end{enumerate}
\end{corollaire}

En
\marginpar{43}
particulier:

\begin{corollaire}
\label{II.4.12}
Soient $X$ un $S$-pr\'esch\'ema, $F$ une section de~$\cal{O}_X$,
$Y$ le sous-pr\'esch\'ema des z\'eros de~$F$ (d\'efini par
l'Id\'eal coh\'erent $F\cdot\cal{O}_X$), $x$ un point de~$Y$. On
suppose $X$ simple sur~$S$ en $x$. Pour que $Y$ soit lisse sur~$S$ en
$x$, il faut et il suffit que ou bien $F$ soit nul au voisinage de
$x$, ou bien que $dF(x) \neq 0$ (o\`u $dF(x)$ d\'esigne l'image de
$dF$ dans l'espace vectoriel $\it{\Omega}^1_{X/S}(x)$ sur~$\kres(x)$).
\end{corollaire}

C'est suffisant en vertu de~\Ref{II.4.10} crit\`ere~(iii). C'est
n\'ecessaire, car comme $\cal{J}$ est engendr\'e par un
\'el\'ement, il faut d'abord que $\cal{J}/\cal{J}^2$ au point $x$
soit libre de rang $\leq 1$. Si ce rang est $0$,
\ie $\cal{J}/\cal{J}^2=0$ en $x$, il s'ensuit que $\cal{J}=0$ en $x$
par Nakayama, \ie $F$ est nul au voisinage de~$x$. Si ce rang est
$1$, alors $F$ forme un syst\`eme minimal de g\'en\'erateurs de
$\cal{J}$ en $x$, et on conclut par \eqref{II.4.11}, \'equivalence de
(i) et~(iii)).
\begin{corollaire}
\label{II.4.13}
Soient $Y$ un $S$-pr\'esch\'ema localement de type fini, $S'$ un
$S$-pr\'e\-sch\'ema \emph{plat}, $Y'=Y \times_S S'$, $x'$ un point
de~$Y'$ et $x$ son image canonique dans~$Y$. Pour que $Y$ soit lisse
sur~$S$ en~$X$, il faut et il suffit que $Y'$ soit lisse sur~$S'$ en~$x'$. En particulier, si $S'\to S$ est plat et surjectif, $Y$ est
lisse sur~$S$
\ifthenelse{\boolean{orig}}{\emph{sss}}{si et seulement si\xspace} $Y'$ est lisse sur~$S'$.
\end{corollaire}

Il n'y a \`a prouver que la suffisance (la n\'ecessit\'e a
\'et\'e vue dans~\Ref{II.1.3}~(ii)). On peut supposer (rempla\c cant $Y$ par un voisinage convenable de~$x$, $Y'$ par l'image inverse
de ce dernier) que $Y$ est affine de type fini sur~$S$ affine, donc
$Y$ est isomorphe \`a un sous-pr\'esch\'ema ferm\'e d'un
sch\'ema $S[t_1,\dots,t_n]$. Par suite, $Y'$ s'identifie \`a un
sous-pr\'esch\'ema ferm\'e de~$X'=X \times_S S'$. Comme $X$ est
lisse sur~$S$, donc $X'$ lisse sur~$S'$, on peut appliquer les
crit\`eres de lissit\'e~\Ref{II.4.10}. Ici, le crit\`ere (iv)
donne le r\'esultat aussit\^ot.

\begin{remarques}
\label{II.4.14}
Le crit\`ere (iii) de th\'eor\`eme~\Ref{II.4.10} m\'erite
d'\^etre appel\'e \emph{crit\`ere jacobien de lissit\'e}.
\index{crit\`ere jacobien de lissit\'e|hyperpage}%
\index{jacobien (crit\`ere de lissit\'e)|hyperpage}%
\index{lissit\'e (crit\`ere jacobien de)|hyperpage}%
Il permet de reconna\^{\i}tre, th\'eoriquement, si un
$S$-pr\'esch\'ema donn\'e $Y$ est lisse sur~$S$ en un point $x$
de~$Y$, puisque il existe toujours un voisinage de~$Y$ isomorphe \`a
un sous-pr\'esch\'ema d'un $S$-pr\'esch\'ema lisse $X$, par
exemple $X=S[t_1,\dots,t_n]$. C'est d'ailleurs pour
$X=S[t_1,\dots,t_n]$, $S=\Spec(A)$, qu'on \'enonce d'habitude le
crit\`ere jacobien (bien
\marginpar{44}
entendu, dans le cas classique envisag\'e par Zariski, $A$ \'etait
un corps). On laisse au lecteur de donner l'\'enonc\'e relatif
\`a la donn\'ee d'un id\'eal $J$ de~$A[t_1,\dots,t_n]$ et d'un
id\'eal premier le contenant, auquel on est ainsi conduit. Notons
qu'il semble bien \`a l'heure actuelle (et surtout depuis que Nagata
est parvenu \`a g\'en\'eraliser par des m\'ethodes
non-diff\'erentielles le th\'eor\`eme de Zariski disant que
l'ensemble des points r\'eguliers d'un sch\'ema alg\'ebrique est
ouvert) que le crit\`ere jacobien n'a gu\`ere d'int\'er\^et
que sous la forme o\`u nous le donnons ici (\ie en utilisant
exclusivement des diff\'erentielles \emph{relatives} et non pas des
diff\'erentielles \emph{absolues}, \ie relatives \`a l'anneau de
constantes absolues $\ZZ$). Comme bien souvent, la consid\'eration
des diff\'erentielles est plus commode ici que celle des
d\'erivations. Notons enfin que si $Y$ est lisse sur~$S$ en $x$, de
dimension relative $n$, alors il existe un voisinage ouvert de~$x$
dans~$Y$ isomorphe \`a un sous-pr\'esch\'ema de
$X=S[t_1,\dots,t_n]$ avec $n=m+1$, comme il r\'esulte de la
d\'efinition et
\ifthenelse{\boolean{orig}}{de~I~\Ref{I.7.6}).}{de~(I~\Ref{I.7.6}).}

Soient $A$ un anneau noeth\'erien, $x_i$ ($1 \leq i \leq n$) des
\'el\'ements de~$A$, $J$ l'id\'eal engendr\'e par les
$x_i$. On dit que les $x_i$ forment un \emph{syst\`eme r\'egulier
de g\'en\'erateurs}
\index{regulier (systeme de generateurs)@r\'egulier (syst\`eme de g\'en\'erateurs)|hyperpage}%
\index{systeme regulier de generateurs)@syst\`eme r\'egulier de g\'en\'erateurs|hyperpage}%
de~$J$ si l'homomorphisme surjectif canonique
$$
(A/J)[t_1,\dots,t_n] \to \gr^J(A)
$$
d\'efini par les $x_i$ (o\`u le deuxi\`eme membre d\'esigne
l'anneau gradu\'e associ\'e \`a $A$ filtr\'e par les
puissances de~$J$) est un \emph{isomorphisme}. Cette condition
signifie aussi que
\begin{enumerate}
\item[(i)] L'homomorphisme surjectif canonique
$$
S_{A/J} (J/J^2) \to \gr^J(A)
$$
(o\`u, le premier membre d\'esigne l'alg\`ebre sym\'etrique du
$A/J$-module $J/J^2$) est un isomorphisme, et
\item[(ii)] $J/J^2$ est libre et admet pour base les classes des $x_i
\mod J^2$.
\end{enumerate}
Sous cette forme, on voit que si $J\neq A$, les $x_i$ forment un
\emph{syst\`eme minimal de g\'en\'erateurs} de~$J$, et que
\emph{tout autre syst\`eme minimal de g\'en\'erateurs} de~$J$
\emph{est un syst\`eme r\'egulier de g\'en\'erateurs}
(N.B. \og minimal\fg est pris au sens strict: nombre minimum
d'\'el\'ements, qui n'est \'equivalent au sens: minimal pour
l'inclusion, que si $A$ est local); d'autre part, si $J=A$, tout
syst\`eme de g\'en\'erateurs de~$J$ est r\'egulier.

La
\marginpar{45}
condition de r\'egularit\'e d'un syst\`eme de
g\'en\'erateurs d'un id\'eal est stable par localisation par un
ensemble multiplicativement stable quelconque, et d'autre part on voit
tout de suite que pour que $(x_i)$ soit un syst\`eme minimal de
g\'en\'erateurs de~$J$, il suffit d\'ej\`a que pour tout id\'eal
\emph{maximal $\goth m$ contenant $J$}, les $x_i$ d\'efinissent un
syst\`eme r\'egulier de g\'en\'erateurs de~$JA_{\goth m}$ dans
$A_{\goth m}$. Cela nous ram\`ene donc au cas o\`u $A$ est un
anneau local d'id\'eal maximal $\goth m$, et o\`u les $x_i$ sont
dans $\goth m$. \emph{Alors les $x_i$ forment un syst\`eme
r\'egulier de g\'en\'erateurs de~$J$ si et seulement si ils
forment une $A$-suite au sens de Serre}\footnote{Nous dirons
maintenant plut\^ot \og suite $A$-r\'eguli\`ere\fg, \cf EGA
$0_{\textup{IV}}$~15.1.7 et 15.1.11.}, \ie si pour tout $i$ tel que $1
\leq i \leq n$, $x_i$ est non-diviseur de~$0$ dans
\ifthenelse{\boolean{orig}}{$A/(x_1,\dots,x_i)A$.}
{$A/(x_1,\dots,x_{i-1})A$.}

Enfin, dans le cas o\`u $A$ est une alg\`ebre sur un anneau $B$,
et o\`u $A/J$ est isomorphe comme $B$-alg\`ebre \`a $B$ (de
sorte que $J$ est le noyau d'un homomorphisme de~$B$-alg\`ebres
$A\to B$), alors les $x_i$ forment un syst\`eme r\'egulier de
g\'en\'erateurs de~$J$ si et seulement si l'homomorphisme
canonique
$$
B[[ t_1,\dots,t_n]]\to \hat{A}
$$
d\'efini par les $x_i$ (o\`u le deuxi\`eme membre d\'esigne le
compl\'et\'e s\'epar\'e $\varprojlim A/J^{n+1}$ de~$A$ pour la
topologie d\'efinie par les puissances de~$J$) est un
\emph{isomorphisme} (il est en tout cas \emph{surjectif}).

Tous ces faits sont bien connus, et figurent sans doute dans le cours
de Serre d'alg\`ebre commutative r\'edig\'e par Gabriel, \`a
peu de choses pr\`es. (O\`u on trouve $N$ autres
caract\'erisations des $A$-suites, dans le cas o\`u $A$ est un
anneau local).

Soit $J$ un id\'eal dans un anneau noeth\'erien $A$. On dira que
$J$ est un \emph{id\'eal r\'egulier}
\index{ideal regulier@id\'eal r\'egulier|hyperpage}%
\index{regulier (ideal)@r\'egulier (id\'eal)|hyperpage}%
si pour tout id\'eal premier $\goth p$ de~$A$, $JA_{\goth p}$ admet
un syst\`eme r\'egulier de g\'en\'erateurs. Il suffit
\'evidemment de le v\'erifier pour $\goth p \supset J$, et on peut
de plus se borner \`a $\goth p$ maximal. Plus g\'en\'eralement,
soit $\cal{J}$ un Id\'eal sur un pr\'esch\'ema localement
noeth\'erien $X$, on dit que $\cal{J}$ est un \emph{Id\'eal
r\'egulier} si pour tout $x\in X$, $\cal{J}_x$ est un id\'eal de
$\cal{O}_x$ qui admet un syst\`eme r\'egulier de
g\'en\'erateurs. Cela \'equivaut \`a la conjonction des deux
conditions suivantes:
\begin{enumerate}
\item[(a)] L'homomorphisme canonique surjectif
$$
\ifthenelse{\boolean{orig}}
{S_{\cal{O}/\cal{J}}(\cal{J}/\cal{J}^2) \to \gr^{\cal{J}}(\cal{O}_X)}
{S_{\cal{O}_X/\cal{J}}(\cal{J}/\cal{J}^2) \to \gr^{\cal{J}}(\cal{O}_X)}
$$
est un isomorphisme et
\item[(b)] Le faisceau de~$\cal{O}_X/\cal{J}$-Modules
$\cal{J}/\cal{J}^2$ est localement libre.
\end{enumerate}

On
\marginpar{46}
dit alors aussi que le sous-pr\'esch\'ema $Y$ de~$X$ d\'efini
par $\cal{J}$ (donc tel que $\cal{O}_Y$ prolong\'e par $0$ soit
isomorphe \`a $\cal{O}_X/\cal{J}$) est \emph{r\'eguli\`erement
immerg\'e}
\index{regulierement immerge (preschema)@r\'eguli\`erement immerg\'e (pr\'esch\'ema)|hyperpage}%
dans~$X$, et on d\'efinit de m\^eme (de fa\c con
\'evidente) la notion de morphisme d'\emph{immersion
r\'eguli\`ere},
\index{immersion r\'eguli\`ere|hyperpage}%
\index{reguliere (immersion)@r\'eguli\`ere (immersion)|hyperpage}%
(\resp \emph{r\'eguli\`ere en un point $x$}), morphisme
d'immersion $Y\to X$ identifiant $Y$ (\resp un voisinage convenable de
$x$), \`a un sous-pr\'esch\'ema ferm\'e r\'eguli\`erement
immerg\'e dans un ouvert de~$X$. (Il ne faut pas dire:
sous-pr\'esch\'ema r\'egulier, car cela signifierait que les
anneaux locaux de~$Y$ sont r\'eguliers). Enfin, des sections $x_i$
de~$\cal{J}$ sont appel\'ees \emph{syst\`eme r\'egulier de
g\'en\'erateurs} si pour tout $x\in X$, les \'el\'ements
correspondants de~$\cal{O}_x$ forment un syst\`eme r\'egulier de
g\'en\'erateurs de~$\cal{J}_x$ (terminologie compatible avec celle
introduite pour des g\'en\'erateurs d'un id\'eal d'un
anneau). Cela signifie aussi que l'homomorphisme surjectif canonique
$$
\cal{O}_Y[t_1,\dots,t_n] \to \gr^{\cal{J}} (\cal{O}_X)
$$
d\'efini par les $x_i$ est un isomorphisme. Si on sait par avance
que l'Id\'eal $\cal{J}$ est r\'egulier, cela signifie aussi,
simplement, que en tout point $x$ \emph{de~$Y$}, les $x_i$
d\'efinissent une \emph{base} de~$\cal{J}/\cal{J}^2$ sur~$\cal{O}_{Y,x}$. (N.B. cette condition est vide si $Y$ est
vide). Ainsi, pour que $\cal{J}$ admette un syst\`eme r\'egulier
de g\'en\'erateurs, il faut et il suffit que $\cal{J}$ soit
r\'egulier, et le $\cal{O}_Y$-Module $\cal{J} /\cal{J}^2$ soit
globalement libre (et non-seulement localement libre), \ie que
l'homomorphisme canonique $S_{\cal{O}_Y}(\cal{J}/\cal{J}^2) \to
\gr^{\cal{J}} (\cal{O}_X)$ soit surjectif, et que le
$\cal{O}_Y$-Module $\cal{J}/\cal{J}^2$ soit globalement libre.

Un \emph{anneau augment\'e est dit r\'egulier}
\index{anneau augment\'e r\'egulier|hyperpage}%
\index{regulier (anneau augmente@r\'egulier (anneau augment\'e)|hyperpage}%
si l'id\'eal de l'augmentation est r\'egulier. Ainsi, si $A$ est
un anneau local, consid\'er\'e comme augment\'e dans son corps
r\'esiduel $k$, alors $A$ est un anneau local r\'egulier si et
seulement si c'est un anneau augment\'e r\'egulier.

(\`A vrai dire, il semble qu'il \'etait inutile de commencer par
faire le sorite pr\'eliminaire pour les anneaux, il y a
int\'er\^et \`a commencer avec les faisceaux tout de suite. Si
on veut quelque chose dans le cas noeth\'erien, c'est la
d\'efinition adopt\'ee ici --- a priori moins stricte que celle par
les $A$-suites de Serre --- qui semble pr\'ef\'erable pour les
besoins du calcul diff\'erentiel. Bien entendu, pour bien faire, il
faudrait d\'evelopper aussi au moins une partie de la th\'eorie
des morphismes lisses dans le cadre non-noeth\'erien\footnote{Comme
il est dit dans l'avant-propos, c'est chose faite maintenant, \cf EGA
IV 17, 18}, probablement en
\marginpar{47}
partant du crit\`ere jacobien, de fa\c con \`a obtenir si
possible toutes les propri\'et\'es formelles essentielles des
morphismes lisses et des morphismes \'etales \ie lisses et
quasi-finis; les r\'eciproques seules faisant appel \`a des
hypoth\`eses noeth\'eriennes).
\end{remarques}

Apr\`es ces longs pr\'eliminaires terminologiques, un petit
th\'eor\`eme:
\begin{theoreme}
\label{II.4.15}
Soient $X$ un $S$-pr\'esch\'ema localement de type fini, $Y$ un
sous-pr\'esch\'ema ferm\'e de~$X$ d\'efini par un faisceau
coh\'erent $\cal{J}$ d'id\'eaux sur~$X$, $x$ un point de~$X$. On
suppose maintenant \emph{$Y$ lisse sur~$S$ en $x$} (et rien sur~$X$). Alors les conditions suivantes sont \'equivalentes:
\begin{enumerate}
\item[(i)] $X$ est lisse sur~$S$ en $x$
\item[(ii)] L'immersion $i\colon Y \to X$ est r\'eguli\`ere en
$x$, \ie $\cal{J}_x$ est un id\'eal r\'egulier de~$\cal{O}_x$.
\end{enumerate}
\end{theoreme}

\begin{corollaire}
\label{II.4.16}
Supposons $Y$ \emph{lisse} sur~$S$. Pour que $X$ soit lisse sur~$S$
dans un voisinage de~$Y$ (\ie aux points de~$Y$) il faut et il suffit
que $Y$ soit r\'eguli\`erement plong\'e dans $X$, \ie que
l'immersion $i \colon Y \to X$ soit r\'eguli\`ere.
\end{corollaire}

\subsubsection*{D\'emonstration}
(i) implique (ii). On applique \Ref{II.4.10} crit\`ere~(ii), comme
$g\colon X_1 \to X$ est \emph{plat}, pour montrer que l'image inverse
par $g$ du sous-pr\'esch\'ema $Y'$ de~$X'$ est
r\'eguli\`erement plong\'e, on est ramen\'e \`a prouver que
$Y'=S[t_{p+1},\dots, t_n]$ est r\'eguli\`erement plong\'e dans
$S[t_1,\dots,t_n]$, ce qui est trivial (les $t_i$ ($1 \leq i \leq
p$) forment un syst\`eme r\'egulier de g\'en\'erateurs de
l'Id\'eal d\'efinissant $Y'$ dans $X'$).

(ii) implique (i). Soit $g_i$ ($1 \leq i \leq p$) un syst\`eme
r\'egulier de g\'en\'erateurs de~$\cal{J}_x$ et soient $g_i$
($p+1 \leq i \leq n$) des \'el\'ements de~$\cal{O}_{X,x}$ tels que
leurs images $g'_i$ dans $\cal{O}_{Y,x}$ d\'efinissent un morphisme
\emph{\'etale}
$$
Y_1 \to Y'=S[t_{p+1},\dots, t_n]
$$
d'un voisinages $Y_1$ de~$Y$ dans $Y'$. Les $g_i$ ($1 \leq i \leq n$)
proviennent de sections (de m\^eme nom) de~$\cal{O}_X$ sur un
voisinage $X_1$ de~$x$, et on peut supposer $X_1=X$, $Y_1=Y$. On
obtient ainsi un morphisme
$$
g\colon X \to X'=S[t_1,\dots,t_n]
$$
et tout revient \`a montrer que ce morphisme est \emph{\'etale}
en~$x$. Prenant $X_1$ assez petit,
\marginpar{48}
on peut supposer que les $g_i$ ($1 \leq i \leq p$) forment un
syst\`eme r\'egulier de g\'en\'erateurs de~$\cal{J}$ sur tout
$X$. En particulier, ils engendrent $\cal{J}$, donc le
sous-pr\'esch\'ema $Y$ de~$X$ s'identifie \`a l'image inverse
par $g$ du sous-pr\'esch\'ema $Y'$ de~$X'$. Soit $x'=g(x)$, alors
la fibre de~$X' \to X$ en $x'$ est
\ifthenelse{\boolean{orig}}{donc}{\ignorespaces}
identique \`a la fibre de~$Y \to Y'$ en $x$, donc est \'etale sur~$\kres(x')$, donc $g$ est \emph{non ramifi\'e} en $x$, reste \`a
prouver que $g$ est \emph{plat} en $x$. Or le gradu\'e associ\'e
\`a $\cal{O}_{X',x'}$ filtr\'e par les puissances de
$\cal{J}'_{x}$ est \emph{libre} sur~$\cal{O}_{Y',x'}$ en tous
degr\'es, d'autre part le gradu\'e associ\'e \`a
$\cal{O}_{X,x}$ filtr\'e par les puissances de
$\cal{J}_x=\cal{J}'_x\cal{O}_{X,x}$ est isomorphe (sous
l'homomorphisme canonique) au produit tensoriel du pr\'ec\'edent
par $\cal{O}_{Y,x}$ (puisque l'un et l'autre anneau sont des anneaux
de polyn\^omes \`a $n-p$ ind\'etermin\'ees, \`a anneau de
constantes $\cal{O}_{Y',x'}$, \resp $\cal{O}_{Y,x}$), enfin sur~$\cal{O}_{X',x'}/\cal{J}'_{x'}= \cal{O}_{Y',x'}$
$\cal{O}_{X,x}/\cal{J}_x =\cal{O}_{Y,x}$ est plat.

D'apr\`es un crit\`ere g\'en\'eral de platitude (valable pour
un homomorphisme local d'anneaux locaux noeth\'eriens $A' \to A$,
$A'$ \'etant muni d'un id\'eal $J' \neq A'$ tel que le gradu\'e
associ\'e soit libre sur~$A'/J'$ en toute dimension) il s'ensuit que
$X$ est plat sur~$X'$ en $x$, cqfd.

\begin{corollaire}
\label{II.4.17}
Soient $X$ un pr\'esch\'ema localement de type fini sur~$Y$, $i$
une section de~$X$ sur~$Y$, $y$ un point de~$Y$, $x=i(y)$, $\cal{J}$
le faisceau d'id\'eaux sur~$X$ d\'efini par le
sous-pr\'esch\'ema $i(Y)$ (que nous supposons ferm\'e pour
simplifier l'\'enonc\'e, condition v\'erifi\'ee si $X$ est un
sch\'ema).

Les conditions suivantes sont \'equivalentes:
\begin{enumerate}
\item[(i)] $X$ est lisse sur~$Y$ en $x$
\item[(ii)] $i$ est une immersion r\'eguli\`ere en $y$
\item[(iii)] La $\cal{O}_y$-alg\`ebre compl\'et\'ee de
$\cal{O}_x$ pour la topologie d\'efinie par les puissances de
$\cal{J}_x$ est isomorphe \`a une alg\`ebre de s\'eries
formelles $\cal{O}_y[[ t_1,\dots,t_n]]$.
\item[(iii~bis)] Il existe un voisinage ouvert $U$ de~$y$ tel que le
faisceau d'alg\`ebres $\varprojlim i^*(\cal{O}_X/\cal{J}^{n+1})$ sur~$\cal{O}_Y$ soit isomorphe \`a un faisceau de la forme
$\cal{O}_Y[[ t_1,\dots, t_n]]$ au-dessus de~$U$.
\item[(iv)] Il existe un voisinage ouvert $U$ de~$y$, et un voisinage
ouvert $V$ de~$x$, et enfin un $Y$-morphisme $g \colon V \to
U[t_1,\dots, t_n]$, tel que $g$ soit \'etale, que
\ifthenelse{\boolean{orig}}{$s$}{$i$}
induise une section de~$V$ sur~$U$, transform\'ee par $g$ en la
section nulle de~$U[t_1,\dots, t_n]$ sur~$U$.
\end{enumerate}
\end{corollaire}

L'\'equivalence
\marginpar{49}
de (i) et (ii) est un cas particulier de
th\'eor\`eme~\Ref{II.4.15}, en faisant $Y=S$, l'\'equivalence de
(ii) et (iii) (et moralement de (ii) et (iii~bis)) a \'et\'e
signal\'e avec les \og rappels\fg. quant \`a l'\'equivalence de (i)
et (iv), elle se d\'eduit facilement de
th\'eor\`eme~\Ref{II.4.10} (\'equivalence des conditions (i) et
(ii) dudit).

\begin{corollaire}
\label{II.4.18}
Soit $X$ un pr\'esch\'ema lisse au-dessus de~$S$. Alors le
morphisme diagonal
$$
\Delta_{X/S}\colon X\to X\times_S X
$$
est une
\ifthenelse{\boolean{orig}}
{\emph{immersion r\'eguli\`ere}.}
{\emph{immersion r\'eguli\`ere},}
ou comme on dit encore, $X$ est \og \emph{diff\'erentiablement
lisse}\fg sur~$S$.
\index{diff\'erentiablement lisse|hyperpage}%
\end{corollaire}


En effet, c'est un cas particulier de corollaire~\Ref{II.4.16}, puisque
$X$ et $X\times_S X$ sont tous deux lisses sur~$S$.

\setcounter{subsection}{17}

\begin{remarques}
\label{rem:II.4.18}
Rappelons (I~\Ref{I.1}) que si $X$ est un pr\'esch\'ema au-dessus
de~$S$, on introduit les faisceaux quasi-coh\'erents
\ifthenelse{\boolean{orig}}
{d'Alg\`ebres}
{d'alg\`ebres}
$\cal{P}_{X/S}^n=\cal{O}_{X\times_S X}/\cal{I}_X^{n+1}$ sur~$X$,
(o\`u $\cal{I}_X$ d\'esigne le faisceau
\ifthenelse{\boolean{orig}}
{d'Id\'eaux}
{d'id\'eaux}
qui d\'efinit la diagonale dans $X\times_S X$), consid\'er\'e comme
faisceau de~$\cal{O}_X$-alg\`ebres gr\^ace \`a la premi\`ere
projection $\pr_1\colon X\times_S X\to X$. Les
$\cal{P}_{X/S}^n$ forment un syst\`eme projectif d'Alg\`ebres sur~$X$, dont la limite projective est not\'ee $\cal{P}_{X/S}^{\infty}$
et n'est autre que le faisceau structural du compl\'et\'e formel
de~$X\times_S X$ le long de la diagonale (en supposant maintenant $X$
localement de type fini sur~$S$, donc les $\cal{P}_{X/S}^n$
coh\'erents). Dire que $X$ est diff\'erentiablement lisse sur~$S$
(\ie que le morphisme diagonal $\Delta_{X/S}$ est une immersion
r\'eguli\`ere) signifie aussi que $\cal{P}_{X/S}^{\infty}$ est
r\'egulier, en tant que faisceau d'alg\`ebres augment\'e vers
$\cal{O}_X$, \ie que $\it{\Omega}_{X/S}^1$ est localement libre et
l'homomorphisme surjectif canonique
$$
S_{\cal{O}_X}(\it{\Omega}_{X/S}^1)\to
\mathrm{gr}_*(\cal{P}_{X/S}^{\infty})
$$
est un isomorphisme, ou enfin que tout point de~$X$
\ifthenelse{\boolean{orig}}{\`a}{a}
un voisinage ouvert sur lequel le faisceau
\ifthenelse{\boolean{orig}}
{d'Alg\`ebres}
{d'alg\`ebres}
augment\'ees $\cal{P}_{X/S}^{\infty}$ soit isomorphe \`a un
faisceau $\cal{O}_X[[t_1,\dots,t_n]]$.

Soit $s$ une section de~$X$ sur~$S$, $\cal{J}$ le faisceau
d'id\'eaux sur~$X$ qu'elle d\'efinit (supposant
\marginpar{50}
pour simplifier que $s(S)$ est ferm\'e), on a alors des
isomorphismes canoniques de~$\cal{O}_X$-alg\`ebres augment\'ees:
\begin{equation}
\label{eq:II.4.4}
s^*(\cal{P}_{X/S}^n)=\cal{O}_X/\cal{J}^{n+1}\quoi,\quad
s^*(\cal{P}_{X/S}^{\infty})= \varprojlim_{n}
\cal{O}_X/\cal{J}^{n+1}
\end{equation}
Ces isomorphismes sont fonctoriels dans un sens \'evident par
changement de base, et (compte tenu de ce fait) redonnent une
caract\'erisation des faisceaux d'alg\`ebres $\cal{P}_{X/S}^n$ sur~$S$. Si par exemple $S=\Spec(k)$, $k$ un corps, alors la donn\'ee
d'une section~$s$ de~$X$ sur~$S$ \'equivaut \`a la donn\'ee d'un
point $x$ de~$X$ rationnel sur~$k$, et les formules
pr\'ec\'edentes signifient que l'on a un isomorphisme de
$k$-alg\`ebres
\begin{equation}
\label{eq:II.4.5}
\cal{P}_{X/S}^n(x)=\cal{O}_x/\goth{m}_x^{n+1}
\end{equation}
ce qui justifie le nom: \og \emph{faisceau des parties principales
d'ordre~$n$ sur~$X$ rel. \`a~$S$}\fg
\index{parties principales d'ordre~$n$ (faisceau des, syst\`eme des)|hyperpage}%
\label{indnot:ad1}\oldindexnot{$\cal{P}_{X/Y}^n$|hyperpage}%
donn\'e \`a~$\cal{P}_{X/S}^n$. On voit de plus sur
\eqref{eq:II.4.4} que \emph{si $X$ est diff\'erentiablement lisse
sur~$S$ en tout point de~$s(S)$, alors $X$ est lisse sur~$S$ en tout
point de~$s(S)$}, (corollaire~\Ref{II.4.17}) \emph{la r\'eciproque
\'etant \'egalement vraie} (corollaire~\Ref{II.4.18}). Compte tenu
de~\Ref{II.4.13}, on en conclut facilement que si $X$ est un
$S$-pr\'esch\'ema localement de type fini, \emph{$X$ est lisse sur~$S$ si et seulement si il est plat sur~$S$ et diff\'erentiellement
lisse sur~$S$.} (N.B. l'hypoth\`ese de platitude est essentielle,
comme on voit en prenant pour $X$ un sous-pr\'esch\'ema ferm\'e
de~$S$).

Notons encore, \`a titre de rappel, qu'on obtient une
\emph{deuxi\`eme structure d'alg\`ebre} sur~$\cal{P}_{X/S}^n$
gr\^ace \`a la projection $\pr_2\colon X\times_S X\to X$,
se d\'eduisant d'ailleurs de la pr\'ec\'edente \`a l'aide de
l'\emph{involution canonique} du faisceau d'anneaux $\cal{P}_{X/S}^n$,
induit par l'automorphisme de sym\'etrie de~$X\times_S X$. On note
par $d_{X/S}^n$
\label{indnot:bb}\oldindexnot{$d_{X/S}^n$|hyperpage}%
ou simplement $d^n$, l'homomorphisme de faisceaux d'anneaux
\begin{equation}
\label{eq:II.4.6}
d_{X/S}^n\colon \cal{O}_X\to\cal{P}_{X/S}^n
\end{equation}
qui correspond \`a cette deuxi\`eme structure
d'Alg\`ebre. Compte tenu de l'isomorphisme~\eqref{eq:II.4.4}, cet
homomorphisme transforme une section $f$ de~$\cal{O}_X$ en une section
$d^n(f)$ de~$\cal{P}_{X/S}^n$ dont l'image inverse par une section $s$
de~$X$ sur~$S$ s'identifie \`a l'image canonique de~$f$ dans
$\Gamma(X,\cal{O}_X/\cal{J}^{n+1})$. Cela justifie le nom de
\og \emph{syst\`eme des parties principales d'ordre $n$ de~$f$}\fg
\index{systeme des parties principales d'ordre n d'une section@syst\`eme des parties principales d'ordre~$n$ d'une section|hyperpage}%
donn\'e \`a $d^nf$, notamment dans le cas o\`u $S=\Spec(k)$,
\marginpar{51}
envisag\'e dans la formule~\eqref{eq:II.4.5}.

Pour finir, notons que l'homomorphisme \eqref{eq:II.4.6} peut \^etre
consid\'er\'e comme l'\emph{op\'e\-ra\-teur diff\'erentiel
d'ordre}~$\leq n$\kern1pt\footnote{Pour tout ce qui concerne le pr\'esent
alin\'ea, on pourra consulter EGA~IV~16.8 \`a 16.12.}
(relativement au pr\'esch\'ema des constantes $S$)
\emph{universel}
\index{op\'erateur diff\'erentiel d'ordre~$\leq n$ universel|hyperpage}%
sur~$\cal{O}_X$, en convenant d'appeler op\'erateur diff\'erentiel
d'ordre~$\leq n$
\index{op\'erateur diff\'erentiel d'ordre~$\leq n$|hyperpage}%
de~$\cal{O}_X$ dans un Module $F$, un homomorphisme de faisceaux $D$
qui se factorise en
$$
D\colon \cal{O}_X\lto{d^n} \cal{P}_{X/S}^n\lto{u} F
$$
o\`u $u$ est un homomorphisme \emph{de~$\cal{O}_X$-Modules},
d'ailleurs uniquement d\'etermin\'e par~$D$. Cette d\'efinition
concorde avec la d\'efinition r\'ecurrente intuitive ($D$ est un
op\'erateur diff\'erentiel d'ordre $\leq n$ si pour toute section
$g$ de~$\cal{O}_X$ sur un ouvert $U$ de~$X$, $f\mto D(fg)-D(f)$ est
un op\'erateur diff\'erentiel d'ordre $\leq n-1$ sur~$U$). Il
s'ensuit que \emph{si $X$ est diff\'erentiablement lisse sur~$S$, le
faisceau d'anneaux des op\'erateurs diff\'erentiels de tous ordres
a la structure simple bien connue} en calcul diff\'erentiel sur les
vari\'et\'es diff\'erentiables, et en particulier admet
localement une base de~$\cal{O}_X$-Module
\ifthenelse{\boolean{orig}}{ferm\'e }{form\'e }
des \emph{puissances divis\'ees} en des op\'erateurs permutables
$\delta/\delta x_i$ $(1\leq i\leq n)$. Si~$S$ est un faisceau de
$\QQ$-alg\`ebres ($\QQ=$ corps des rationnels) il suffit de prendre
les polyn\^omes ordinaires en les $\delta/\delta x_i$. Dans ce cas
d'ailleurs, et tr\`es exceptionnellement, pour que $X$ soit
diff\'erentiablement lisse sur~$S$, il suffit d\'ej\`a que
$\it{\Omega}_{X/S}^1$ soit localement libre.
\end{remarques}

\begin{remarque}
\label{II.4.19}
\index{Cohen-Macaulay (sch\'ema de)|hyperpage}%
La terminologie \og immersion r\'eguli\`ere\fg, \og id\'eal
r\'egulier\fg, etc. introduite dans ce num\'ero a rencontr\'e
une opposition assez vive et g\'en\'erale (Chevalley, Serre). On a
propos\'e \og id\'eal de Cohen-Macaulay\fg ou \og id\'eal de
Macaulay\fg ou \og id\'eal macaulayen\fg, ce qui moralement obligerait
\`a adopter aussi \og immersion de Cohen-Macaulay\fg ou \og immersion de
Macaulay\fg. Cette terminologie cependant conflicte avec une autre
d\'ej\`a employ\'ee dans de futures r\'edactions du
multiplodoque, o\`u un morphisme (de type fini) est dit
\og Cohen-Macaulay\fg en un point s'il est plat en ce point, et si la
fibre passant par ce point y a un anneau local qui soit un anneau de
Cohen-Macaulay. En attendant de trouver une solution satisfaisante,
nous garderons sous toutes r\'eserves la terminologie introduite
dans ce num\'ero\footnote{C'est celle adopt\'ee dans
EGA~$\mathrm{0}_\mathrm{IV}$~15.1.7.}.
\end{remarque}

\section{Cas d'un corps de base}
\label{II.5}
\marginpar{52}

\begin{proposition}
\label{II.5.1}
Soient $k$ un corps, $X$ un pr\'esch\'ema de type fini sur~$k$,
$x$ un point de~$X$ et $n$ la dimension de~$X$ en~$x$,
\ifthenelse{\boolean{orig}}
{$f\colon X\to\Spec k[t_1,\dots,tn]=Y$}
{$f\colon X\to\Spec k[t_1,\dots,t_n]=Y$}
un morphisme, d\'efini par des \'el\'ements
\ifthenelse{\boolean{orig}}
{$f_i\in\Gamma(X:\cal{O}_X)$.}
{$f_i\in\Gamma(X,\cal{O}_X)$.}
Les conditions suivantes sont \'equivalentes (et entra\^inent que
$X$ est lisse sur~$k$ en $x$, et a fortiori r\'egulier en $x$
d'apr\`es~\Ref{II.3.1}):
\begin{enumerate}
\item[(i)] $f$ est \'etale en $x$.
\item[(ii)] Les $df_i$ forment une base de~$\it{\Omega}_{X/k}^1$ en $x$.
\item[(iii)] Les $df_i$ engendrent $\it{\Omega}_{X/k}^1$ en $x$.
\end{enumerate}
\end{proposition}

Comme (i) implique que $X$ est lisse sur~$k$ en $x$, l'implication
(i)$\To$(ii) est un cas particulier de~\Ref{II.4.8},
(ii)$\To$(iii) est trivial, reste \`a prouver
(iii)$\To$(i). Or si (iii) est v\'erifi\'e, $f$ est net
en $x$ en vertu de lemme~\Ref{II.4.1}, donc (rempla\c cant $x$ par
un voisinage ouvert) quasi-fini, donc dominant par raison de
dimensions. Comme $Y$ est r\'egulier, il s'ensuit que $f$ est
\'etale par (I~\Ref{I.9.5}~(ii)) ou (I~\Ref{I.9.11}).

\begin{corollaire}
\label{II.5.2}
Sous les conditions pr\'eliminaires de~\Ref{II.5.1}, supposons que
\ifthenelse{\boolean{orig}}
{$\kappa(x)$}
{$\kres(x)$}
soit une extension \emph{finie s\'eparable} de~$k$, et
que les $f_i$ ($1\leq i\leq n$) d\'efinissent des \'el\'ements
de~$\goth{m}_x$. Alors les conditions pr\'ec\'edentes
\'equivalent \`a
\begin{enumerate}
\item[(iv)] Les $f_i$ forment un syst\`eme de g\'en\'erateurs de
$\goth{m}_x$ (ou encore: les $f_i$ mod $\goth{m}_x^2$ forment une base
de~$\goth{m}_x/\goth{m}_x^2$
\ifthenelse{\boolean{orig}}
{sur~$\kappa(x)$}
{sur~$\kres(x)$}).
\end{enumerate}
\end{corollaire}

En effet, (iv)$\To$(iii) en vertu de la suite exacte
\begin{equation}
\label{eq:II.5.1}
\ifthenelse{\boolean{orig}}
{\goth{m}_x/\goth{m}_x^2\to\it{\Omega}_{\cal{O}_x/k}^1\to
\it{\Omega}_{\kappa(x)/k}^1\to 0}
{\goth{m}_x/\goth{m}_x^2\to\Omega_{\cal{O}_x/k}^1\to
\Omega_{\kres(x)/k}^1\to 0}
\end{equation}
\ifthenelse{\boolean{orig}}
{et compte tenu de~$\Omega_{\kappa(x)/k}^1=0$ puisque $\kappa(x)$ est
\'etale sur~$k$. D'autre part (ii) implique (iv), car comme $X$ et
$\Spec(\kappa(x))$ sont lisses sur~$k$ en $x$, on peut dans la suite
exacte pr\'ec\'edente mettre un $0$ sur la gauche en vertu de
\Ref{II.4.10}~(iv).}
{et compte tenu de~$\Omega_{\kres(x)/k}^1=0$ puisque $\kres(x)$ est
\'etale sur~$k$. D'autre part (ii) implique (iv), car comme $X$ et
$\Spec(\kres(x))$ sont lisses sur~$k$ en $x$, on peut dans la suite
exacte pr\'ec\'edente mettre un $0$ sur la gauche en vertu de
\Ref{II.4.10}~(iv).}

\begin{corollaire}
\label{II.5.3}
Soit $x$ un point de~$X$ (de type fini sur~$k$). Si $X$ est lisse sur~$k$ en $x$, alors $\cal{O}_x$ est r\'egulier, la r\'eciproque
\'etant vraie si
\ifthenelse{\boolean{orig}}
{$\kappa(x)$}
{$\kres(x)$}
est une extension finie s\'eparable
de~$k$.
\end{corollaire}

En
\marginpar{53}
effet, la r\'eciproque r\'esulte de~\Ref{II.5.2}, en prenant un
syst\`eme r\'egulier $(f_i)$ de g\'en\'erateurs de
$\goth{m}_x$. (N.B. au lieu de~\Ref{II.5.2}, on peut aussi invoquer le
th\'eor\`eme~\Ref{II.4.15}). On conclut:

\begin{proposition}
\label{II.5.4}
Soit $X$ un pr\'esch\'ema de type fini sur~$k$. Si $X$ est lisse
sur~$k$, il est r\'egulier, la r\'eciproque \'etant vraie si $k$
est parfait.
\end{proposition}

Pour la r\'eciproque, on note qu'en vertu de~\Ref{II.5.3}, $X$ est
lisse sur~$k$ en tout point ferm\'e, donc partout puisque l'ensemble
des points o\`u il est lisse est ouvert.

\begin{theoreme}
\label{II.5.5}
Soient $X$ un pr\'esch\'ema de type fini sur~$k$, $x$ un point de
$X$, $n$ la dimension de~$X$ en $x$, $k'$ une extension parfaite de
$k$. Les conditions suivantes sont \'equivalentes:
\begin{enumerate}
\item[(i)] $X$ est lisse sur~$k$ en~$x$.
\ifthenelse{\boolean{orig}}
{\item[(ii)] $\Omega_{X/k}^1$ est libre de rang $n$ en~$x$.
\item[(ii~bis)] $\Omega_{X/k}^1$ est engendr\'e par $n$}
{\item[(ii)] $\it{\Omega}_{X/k}^1$ est libre de rang $n$ en~$x$.
\item[(ii~bis)] $\it{\Omega}_{X/k}^1$ est engendr\'e par $n$}
\'el\'ements en~$x$.
\item[(iii)] $X$ est diff\'erentiablement lisse sur~$k$ en $x$.
\item[(iv)] Il existe un voisinage ouvert $U$ de~$x$ tel que
$U\otimes_k k'$ soit r\'egulier (\ie les anneaux locaux de ses
points sont r\'eguliers).
\end{enumerate}
\end{theoreme}

On a (i)$\To$(ii) par \Ref{II.4.3}~(ii), (ii)$\To$(ii~bis) trivialement et (ii~bis)$\To$(i) gr\^ace
\`a~\Ref{II.5.1}. Comme $X$ est plat sur~$k$, on a
(i)$\Leftrightarrow$(iii) en vertu de~\Ref{II.4.18}. On a
(i)$\To$(iv) puisque lisse est invariant par extension de la
base et implique r\'egulier, et (iv)$\To$(i) car par
Proposition~\Ref{II.5.4}, on voit que $U\otimes_k k'$ est simple sur~$k'$, donc $U$ est simple sur~$k$ par~\Ref{II.4.13}.

Prenant pour $x$ le point g\'en\'erique de~$X$ suppos\'e
irr\'eductible, on trouve:

\begin{corollaire}
\label{II.5.6}
Soit $K$ un anneau d'Artin local localis\'e d'une alg\`ebre de
type fini sur le corps~$k$ (par exemple, $K$ peut \^etre une
extension de type fini de~$k$), soit~$n$ le degr\'e de transcendance
sur~$K$ de son corps r\'esiduel. Conditions \'equivalentes:
\begin{enumerate}
\item[(i)] $K$ est une extension finie s\'eparable d'une extension
transcendante pure $k(t_1,\dots,t_n)$ de~$k$.
\item[(ii)]
\marginpar{54}
\ifthenelse{\boolean{orig}}{$\Omega^1_{X/k}$}{$\it{\Omega}^1_{X/k}$}
est un $K$-module libre de rang~$n$.
\item[(ii~bis)]
\ifthenelse{\boolean{orig}}{$\Omega^1_{X/k}$}{$\it{\Omega}^1_{X/k}$}
est un $K$-module admettant $n$
g\'en\'erateurs.
\item[(iii)] Le compl\'et\'e $O'$ de~$K\otimes_k K$ pour la
topologie d\'efinie par les puissances de l'id\'eal d'augmentation
$K\otimes_k K\to K$ est une
\ifthenelse{\boolean{orig}}{$K$-alg\`ebr\'e}{$K$-alg\`ebre}
augment\'ee \og r\'eguli\`ere\fg, \ie isomorphe \`a une
alg\`ebre de s\'eries formelles en~$K$ (N.B. Si~$K$ est un corps,
cela \'equivaut \`a dire que~$O'$ est un anneau local
r\'egulier).
\item[(iv)] $K$ est une extension s\'eparable de~$k$.
\end{enumerate}
\end{corollaire}

En effet, on peut toujours consid\'erer $K$ comme l'anneau local du
point g\'en\'erique d'un sch\'ema de type fini
irr\'eductible~$X$ sur~$k$, et les conditions envisag\'ees sont
les conditions de m\^eme nom dans~\Ref{II.5.5}, en prenant dans~(iv)
pour~$k'$ une extension alg\'ebriquement close de~$k$. Seule
l'implication $K$ s\'eparable sur~$k$ $\To$ $X$ lisse sur~$k$ en~$x$, demande une d\'emonstration. Or on est aussit\^ot
ramen\'e gr\^ace \`a~\Ref{II.4.13} au cas o\`u le corps de
base est~$k'$, donc alg\'ebriquement clos, donc o\`u il existe un
point $a$ de~$X$ rationnel sur~$k$. Mais alors $X$ est lisse sur~$k$
en~$a$ d'apr\`es~\Ref{II.5.4}, a fortiori il est lisse sur~$k$ en le
point g\'en\'erique~$x$, cqfd\footnote{\Cf Errata \`a la fin du pr\'esent Exp\ptbl II
(p\ptbl \pageref{II.fin.errata})}.

On remarquera que dans le cas o\`u~$K$ est une extension de type
fini de~$k$, l'\'equivalence de
\ifthenelse{\boolean{orig}}
{(i) (ii) (ii~bis) (iv)}
{(i), (ii), (ii~bis), (iv)}
est bien
connue, mais que nous ne nous sommes servis d'aucune de
\ifthenelse{\boolean{orig}}{ses}{ces}
\ifthenelse{\boolean{orig}}{\'equivalence}{\'equivalences}
d\'ej\`a connues. Bien entendu, la
proposition~\Ref{II.5.1} contient comme cas particulier qu'une suite
d'\'el\'ements~$x_i$ ($1\leq i\leq n$) est un \og base de
transcendance s\'eparante\fg de~$K$ sur~$k$ si et seulement si
les~$dx_i$ forment une base du $K$-module $\Omega^1_{K/k}$, fait bien
connu par ailleurs.

\begin{corollaire}
\label{II.5.7} Soit $X$ un pr\'esch\'ema de type fini sur un
corps~$k$. Pour que $X$ soit lisse sur~$k$, il faut et il suffit que
\ifthenelse{\boolean{orig}}{$\Omega^1_{X/k}$}{$\it{\Omega}^1_{X/k}$}
soit localement libre, et que les anneaux locaux des
points g\'en\'eriques des composantes irr\'eductibles de~$X$
soient des extensions s\'eparables de~$k$ (cette derni\`ere
condition \'etant automatiquement v\'erifi\'ee si $k$ est
parfait et $X$ r\'eduit).
\end{corollaire}

On peut supposer $X$ connexe, soit $n$ le rang de
\ifthenelse{\boolean{orig}}{$\Omega^1_{X/k}$}{$\it{\Omega}^1_{X/k}$}
suppos\'e localement libre. D'apr\`es l'hypoth\`ese
et~\Ref{II.5.6}, c'est aussi le degr\'e de transcendance des
extensions
\marginpar{55}
de~$k$ d\'efinies par les anneaux locaux des points
g\'en\'eriques de~$X$, donc toutes les composantes
irr\'eductibles de~$X$ sont de dimension~$n$. On conclut alors
gr\^ace \`a~\Ref{II.5.5}.

On fera attention que si $K$ est une extension finie (non
n\'ecessairement s\'eparable) de~$k$, alors $\Omega^1_{K/k}$ est
un $k$-module libre, donc posant $X=\Spec(K)$, $\Omega^1_{X/k}$ est un
faisceau localement libre, et $X$ est r\'eduit, sans que $X$ soit
n\'ecessairement lisse sur~$k$. \'Etendant alors les scalaires \`a
la cl\^oture alg\'ebrique de~$k$, on trouve un exemple analogue,
o\`u $k$ est alg\'ebriquement clos, mais $X$ en revanche
n'\'etant pas r\'eduit.

\begin{corollaire}\label{II.5.8}
Soient $X$ un pr\'esch\'ema de type fini sur le corps~$k$, $x$ un
point de~$X$, $n$ la dimension de~$X$ en~$x$, $p$ la dimension
\ifthenelse{\boolean{orig}}
{de~$\cal{O}_x$}
{de~$\cal{O}_x$,}
\ie la codimension dans $X$ de l'adh\'erence $Y$ de
$x$ dans~$X$; donc $n-p$ le degr\'e de transcendance
\ifthenelse{\boolean{orig}}
{de~$\kappa(x)$}
{de~$\kres(x)$}
sur~$k$. Soient $f_i$ ($1\leq i\leq n$) des \'el\'ements
de~$\cal{O}_x$, tels que $f_i\in \goth{m}_x$ pour $1\leq i\leq p$. Les
conditions suivantes sont \'equivalentes

\begin{enumerate}
\item[(i)] le germe de morphisme en~$x$
$$
X\to \Spec\bigl(k[t_1,\dots,t_n]\bigr)
$$
d\'efini par les $f_i$ est \'etale en~$x$.
\item[(ii)] Les $f_i$ ($1\leq i\leq p$)
\ifthenelse{\boolean{orig}}
{engendrent~$\goth{m}_x$}
{engendrent~$\goth{m}_x$,}
\ie forment un syst\`eme r\'egulier de param\`etres
de~$\cal{O}_x$, et les classes dans
\ifthenelse{\boolean{orig}}
{$\kappa(x)$ des~$f_j$ ($p+1\leq j\leq n$) forment une base de transcendance s\'eparante, (\ie les~$d\overline{f}_j$ ($p+1\leq j\leq n$) forment une base
de~$\Omega^1_{\kappa(x)/k}$, ou encore engendrent~$\Omega^1_{\kappa(x)/k}$).}
{$\kres(x)$ des~$f_j$ ($p+1\leq j\leq n$) forment une base de transcendance s\'eparante, (\ie les~$d\overline{f}_j$ ($p+1\leq j\leq n$) forment une base
de~$\Omega^1_{\kres(x)/k}$, ou encore engendrent~$\Omega^1_{\kres(x)/k}$).}
\end{enumerate}
\end{corollaire}

\ifthenelse{\boolean{orig}}
{Supposons~(i) v\'erifi\'e. Il en r\'esulte que les $df_i(x)$
forment une base de~$\it{\Omega}^1_{X/k}(x)$ \eqref{II.4.8} donc leurs
images $d\overline{f}_i(x)$ dans $\Omega^1_{\kappa(x)/k}$ engendrent
cet espace vectoriel sur~$k$. Comme les $\overline{f}_i$ pour $1\leq
i\leq p$ sont nuls, il s'ensuit qu'il suffit de prendre les
$d\overline{f}_i(x)$ avec $p+1\leq i\leq n$. Comme le degr\'e de
transcendance de~$\kappa(x)$ sur~$k$ est~$n-p$, il r\'esulte alors
du corollaire~\Ref{II.5.6} crit\`ere~(iii) (appliqu\'e
\`a~$K=\kappa(x)$) que $Y$ est lisse sur~$k$ en son point
g\'en\'erique~$x$, et que les~$d\overline{f}_i(x)$ ($p+1\leq i\leq
n$) forment une \emph{base} de~$\Omega^1_{\kappa(x)/k}$
sur~$\kappa(x)$. Par suite, la condition~(ii) de~\Ref{II.4.9} est
v\'erifi\'ee, donc aussi la condition~(iii) et en particulier
les~$f_i$ ($1\leq i\leq p$) forment un syst\`eme de
g\'en\'erateurs de~$\goth{m}_x$. Comme}
{Supposons~(i) v\'erifi\'e. Il en r\'esulte que les $df_i(x)$
forment une base de~$\it{\Omega}^1_{X/k}(x)$ \eqref{II.4.8} donc leurs
images $d\overline{f}_i(x)$ dans $\Omega^1_{\kres(x)/k}$ engendrent
cet espace vectoriel sur~$k$. Comme les $\overline{f}_i$ pour $1\leq
i\leq p$ sont nuls, il s'ensuit qu'il suffit de prendre les
$d\overline{f}_i(x)$ avec $p+1\leq i\leq n$. Comme le degr\'e de
transcendance de~$\kres(x)$ sur~$k$ est~$n-p$, il r\'esulte alors
du corollaire~\Ref{II.5.6} crit\`ere~(iii) (appliqu\'e
\`a~$K=\kres(x)$) que $Y$ est lisse sur~$k$ en son point
g\'en\'erique~$x$, et que les~$d\overline{f}_i(x)$ ($p+1\leq i\leq
n$) forment une \emph{base} de~$\Omega^1_{\kres(x)/k}$
sur~$\kres(x)$. Par suite, la condition~(ii) de~\Ref{II.4.9} est
v\'erifi\'ee, donc aussi la condition~(iii) et en particulier
les~$f_i$ ($1\leq i\leq p$) forment un syst\`eme de
g\'en\'erateurs de~$\goth{m}_x$. Comme}
\marginpar{56}
$\cal{O}_x$ est de dimension~$p$, ils forment donc un syst\`eme
r\'egulier de param\`etres en~$x$. Cela prouve~(ii).

Supposons~(ii) v\'erifi\'e. En vertu de la suite
exacte~\eqref{II.5.1}, il s'ensuit que les~$df_i(x)$
\ifthenelse{\boolean{orig}}{engendrent~$\Omega^1_{X/k}$,}
{engendrent~$\it{\Omega}^1_{X/k}$,}
d'o\`u~(i) gr\^ace \`a~prop\ptbl \Ref{II.5.1}.

\begin{corollaire}
\label{II.5.9} Soient~$X$ un pr\'esch\'ema de type fini sur le
corps~$k$, $x$ un point de~$X$, $n$ la dimension de~$X$ en~$x$, $p$ la
dimension
\ifthenelse{\boolean{orig}}
{de~$\cal{O}_x$}
{de~$\cal{O}_x$,}
\ie la codimension de l'adh\'erence~$Y$
de~$x$ dans~$X$, donc~$n-p$ le degr\'e de transcendance
\ifthenelse{\boolean{orig}}
{de~$\kappa(x)$}
{de~$\kres(x)$}
sur~$k$. Conditions \'equivalentes:
\begin{enumerate}
\item[(i)] $\cal{O}_x$ est r\'egulier
\ifthenelse{\boolean{orig}}
{et~$\kappa(x)$}
{et~$\kres(x)$}
 est une
extension s\'eparable de~$k$.
\item[(ii)] $X$ est lisse sur~$k$ en~$x$, et l'homomorphisme canonique
$$
\ifthenelse{\boolean{orig}}
{\goth{m}_x/\goth{m}_x^2\to
\Omega^1_{\cal{O}_x/k}\otimes_{\cal{O}_x}
\kappa(x)=\it{\Omega}^1_{X/k}(x)}
{\goth{m}_x/\goth{m}_x^2\to
\Omega^1_{\cal{O}_x/k}\otimes_{\cal{O}_x}
\kres(x)=\it{\Omega}^1_{X/k}(x)}
$$
est injectif.
\item[(iii)] Il y a des~$f_i\in \cal{O}_x$ ($1\leq i\leq n$)
avec~$f_i\in \goth{m}_x$ pour~$1\leq i\leq p$, tels que le germe de
morphisme en~$x$ de~$X$ dans~$\Spec\bigl(k[t_1,\dots,t_n]\bigr)$
d\'efini par les~$f_i$ soit \'etale en~$x$ (\ie par~\Ref{II.5.1}
tels que les~$df_i(x)$ engendrent~$\it{\Omega}^1_{X/k}(x)$).
\item[(iv)] Il y a des~$f_i\in \cal{O}_x$ ($1\leq i\leq n$) tels que
les~$f_i$ ($1\leq i\leq p$) engendrent~$\goth{m}_x$ et que
les~$df_j(x)$ ($p+1\leq j\leq n$) 
\ifthenelse{\boolean{orig}}
{engendrent~$\Omega^1_{\kappa(x)/k}$ sur~$\kappa(x)$.}
{engendrent~$\Omega^1_{\kres(x)/k}$ sur~$\kres(x)$.}
\end{enumerate}
\end{corollaire}

L'\'equivalence de~(iii) et~(iv) r\'esulte du
corollaire~\Ref{II.5.8}, ces conditions \'equivalent aussi
d'apr\`es~\Ref{II.4.9} au fait que~$X$ est lisse sur~$k$ en~$x$ et
que la condition~(ii)
\ifthenelse{\boolean{orig}}{de~\Ref{II.4.10},}{de~\Ref{II.4.10}}
est v\'erifi\'ee. Donc elles \'equivalent au fait que~$X$ est
lisse sur~$k$ en~$x$ et que la condition~(iv) de~\Ref{II.4.10} est
v\'erifi\'ee, donc \`a~\Ref{II.5.9}~(ii). Ou au fait que~$X$ est
lisse sur~$k$ en~$x$ et que la condition~(i) de~\Ref{II.4.10} est
v\'erifi\'ee, qui ici signifie simplement
\ifthenelse{\boolean{orig}}
{que~$\kappa(x)$}
{que~$\kres(x)$}
est
s\'eparable sur~$k$. Cela implique~\Ref{II.5.9}~(i), il reste \`a
prouver que~\Ref{II.5.9}~(i) l'implique, \ie \`a prouver le

\begin{corollaire}
\label{II.5.10}
Soit~$x$ un point
\ifthenelse{\boolean{orig}}{d'une}{d'un}
pr\'esch\'ema de type fini sur le corps~$k$, tel
\ifthenelse{\boolean{orig}}
{que~$\kappa(x)$}
{que~$\kres(x)$}
soit s\'eparable sur~$k$. Pour que~$X$ soit lisse sur~$k$ en~$x$,
il faut et il suffit qu'il soit r\'egulier en~$x$ (\ie que l'anneau
local~$\cal{O}_x$ soit r\'egulier).
\end{corollaire}

En effet, s'il en est ainsi, on peut \'evidemment trouver
des~$f_i\in \cal{O}_x$ ($1\leq i\leq n$) satisfaisant la
condition~\Ref{II.5.9}~(iv).

\subsection*{Errata}
\label{II.fin.errata}
Dans
\marginpar{57}
le pr\'esent num\'ero, d\'emonstration de~\Ref{II.5.6}, on
a utilis\'e le fait qu'un sch\'ema de type fini r\'eduit non
vide sur un corps alg\'ebriquement clos admet au moins un point
r\'egulier (donc lisse), fait qui se d\'emontre d'habitude par
voie diff\'erentielle (via le th\'eor\`eme de Zariski que
l'ensemble des points r\'eguliers de~$X$ est ouvert). Si on veut
\'eviter un cercle vicieux, il faut d\'emontrer que si~$K/k$ est
une extension s\'eparable de type fini, et si les~$f_i\in K$ sont
tels que~$d_{K/k}f_i$ forment une base de~$\Omega^1_{K/k}$, ($1\leq
i\leq n$), alors~$n$ est le degr\'e de transcendance de~$K$
\ifthenelse{\boolean{orig}}
{sur~$k$}
{sur~$k$,}
\ie les~$f_i$ sont alg\'ebriquement ind\'ependants. La
d\'emonstration de ce fait \`a l'aide du crit\`ere de
\ifthenelse{\boolean{orig}}
{Mac-Lane}
{Mac\kern2pt Lane}
est bien connue, \cf Bourbaki, Alg\`ebre, Chap\ptbl V par\ptbl 9 th\ptbl 2: on
prend un polyn\^ome~$g\in k[t_1,\dots,t_n]$ de degr\'e minimal
tel que~$g(f_1,\dots,f_n)=0$, on a donc
$$
\sum\frac{dg}{d t_i}(f_1,\dots,f_n) df_i=0
$$
d'o\`u (puisque les~$df_i$ forment une base de~$\Omega^1_{K/k}$) le
fait que les~$dg/d t_i$ annulent~$(f_1,\dots,f_n)$, donc sont nuls
d'apr\`es le caract\`ere minimal de~$g$, donc si~$k$ est de
caract\'eristique~$0$ on a~$g=0$, et si~$k$ est de
caract\'eristique~$p\neq 0$ on a~$g=h(t_1^p,\dots,t_n^p)$.
Utilisant le crit\`ere de
\ifthenelse{\boolean{orig}}
{Mac-Lane}
{Mac\kern2pt Lane},
on voit que le
polyn\^ome~$h\in k[t_1,\dots,t_n]$ annule aussi~$(f_1,\dots,f_n)$,
d'o\`u encore~$g=0$ par le caract\`ere minimal de~$g$.

\chapter{Morphismes lisses: propri\'et\'es~de~prolongement}
\label{III}
\marginpar{58}

\section{Homomorphismes formellement lisses}
\label{III.1}

Dans~\Ref{II}, nous nous sommes born\'es \`a la consid\'eration
d'homomorphismes de type fini, et par cons\'equent, dans les
homomorphismes locaux~$A\to B$ d'anneaux locaux, au cas o\`u~$B$ est
isomorphe \`a une alg\`ebre localis\'ee d'une $A$-alg\`ebre de
type fini. Ce cas est insuffisant pour diverses applications,
notamment en g\'eom\'etrie formelle ou en g\'eom\'etrie
analytique. Par exemple, l'anneau de s\'eries
formelles~$B=A[[t_1,\dots,t_n]]$ a (du point de vue de la
g\'eom\'etrie formelle) les propri\'et\'es d'une alg\`ebre
lisse sur~$A$. En g\'eom\'etrie analytique, il en est de m\^eme
de l'anneau local d'un point~$(y,z)$ d'un produit $Y\times\CC^n$
consid\'er\'e comme alg\`ebre sur l'anneau local de~$y$;
d'ailleurs, la compl\'et\'ee de cette alg\`ebre est isomorphe
\`a l'alg\`ebre des s\'eries formelles en~$n$
ind\'etermin\'ees sur le compl\'et\'e de l'anneau de
base~$\cal{O}_x$. C'est ce qui conduit \`a poser la d\'efinition
qui suit.

\begin{definition}
\label{III.1.1}
Soit~$u\colon A\to B$ un homomorphisme local d'anneaux locaux
(noeth\'eriens, on le rappelle). On 
\ifthenelse{\boolean{orig}}
{suppose~$\kappa(B)$ fini sur~$\kappa(A)$}
{suppose~$\kres(B)$ fini sur~$\kres(A)$}.
On dit que~$u$ est un \emph{homomorphisme
formellement lisse},
\index{formellement lisse (homomorphisme, alg\`ebre)|hyperpage}%
\index{lisse (formellement)|hyperpage}%
ou que l'alg\`ebre~$B$ est \emph{formellement lisse sur}~$A$, s'il
existe une $\overline{A}$-alg\`ebre locale
\ifthenelse{\boolean{orig}}{finie,}{finie~$A'$,}
libre sur~$\overline{A}$, telle que les composants locaux de l'anneau
semi-local~$\overline{B}\otimes_{\overline{A}} A'=B'$ soient
$A'$-isomorphes \`a des alg\`ebres des s\'eries formelles
sur~$A'$\kern1pt\footnote{Pour une d\'efinition plus g\'en\'erale et
plus conceptuelle, motiv\'ee par~\Ref{III.2.1} ci-dessous,
cf.\ EGA~$\mathrm{0}_{\mathrm{IV}}$~19.3.1.}.
\end{definition}
(On d\'enote par~$\overline{A}$, $\overline{B}$ les anneaux
compl\'et\'es de~$A$, $B$). Comme~$B'$ est fini et libre
sur~$\overline{B}$, c'est bien un anneau semi-local, compos\'e
direct d'anneaux locaux complets, dont chacun est encore un module
libre sur~$\overline{B}$, donc a m\^eme dimension que~$\overline{B}$
donc que~$B$. Il s'ensuit que le nombre de variables~$t_i$ dans les
anneaux de s\'eries formelles envisag\'es dans~\Ref{III.1.1} est
\'egal \`a $\dim \overline{B}-\dim \overline{A}=\dim B-\dim A$, et
en particulier ind\'ependant du composant local
\marginpar{59}
choisi. On voit tout de suite que c'est aussi la dimension de l'anneau
$B\otimes k=B/\goth{m} B$, o\`u~$k=A/\goth{m}$ est le corps
r\'esiduel de~$A$; on l'appellera la \emph{dimension relative
de}~$B$
\index{dimension relative|hyperpage}%
\emph{par rapport \`a}~$A$.

\begin{remarques}
\label{III.1.2}
Il est \'evident que la d\'efinition~\Ref{III.1.1} ne d\'epend
que de l'homomorphisme sur les compl\'et\'es $\overline{A} \to
\overline{B}$ d\'eduit de~$A\to B$, ce qui justifie dans une
certaine mesure la terminologie. Nous nous repentons ici de la
d\'efinition~I~\Ref{I.3.2}~b) et I~\Ref{I.4.1}~b), qui risque
d'induire en erreur, et pr\'ef\'erons dire \og formellement non
ramifi\'e\fg
\index{formellement non ramifi\'e (\resp net)|hyperpage}%
et \og formellement \'etale\fg
\index{formellement etale@formellement \'etale|hyperpage}%
dans les cas envisag\'es dans ces d\'efinitions, r\'eservant la
terminologie \og non ramifi\'ee\fg et \og \'etale\fg au cas o\`u $B$
est localis\'ee d'une $A$-alg\`ebre de type fini\footnote{Ou
mieux, \og essentiellement non ramifi\'e\fg \resp \og essentiellement
\'etale\fg, comparer EGA IV 18.6.1.}.
\index{essentiellement etale@essentiellement \'etale|hyperpage}%
\index{essentiellement non ramifi\'e (\resp net)|hyperpage}%
Le lecteur v\'erifiera aussit\^ot que \og formellement \'etale\fg
\'equivaut \`a \og formellement lisse et quasi-fini\fg. Enfin,
signalons qu'il y a une d\'efinition raisonnable de \og formellement
lisse\fg sans aucune hypoth\`ese pr\'ealable sur l'extension
r\'esiduelle
\ifthenelse{\boolean{orig}}
{$\kappa(B)/\kappa(A)$}
{$\kres(B)/\kres(A)$}
(suppos\'ee ici finie),
englobant entre autres les homomorphismes locaux $A\to B$ tels que $B$
soit \emph{plat} sur~$A$ et $B/\goth{m}B$ une \emph{extension
s\'eparable} de~$A/\goth{m} =k$ (pas n\'ecessairement de type
fini); par exemple, un $p$-anneau de Cohen est formellement
\ifthenelse{\boolean{orig}}
{simple}
{lisse}
sur
l'anneau des entiers $p$-adiques. C'est la propri\'et\'e de
rel\`evement des homomorphismes (comparer~\Ref{III.2.1}) qui doit
devenir d\'efinition dans ce cas g\'en\'eral. Pour les
applications que nous avons en vue, le cas trait\'e dans la
d\'efinition~\Ref{III.1.1} nous suffira; par la suite, dans
\og formellement lisse\fg nous sous-entendrons \og \`a extension
r\'esiduelle finie\fg.
\end{remarques}

\begin{lemme}
\label{III.1.3}
Si $B$ est formellement lisse sur~$A$, $B$ est plat sur~$A$.
\end{lemme}

\ifthenelse{\boolean{orig}}
{\itshape}{}
Comme la platitude est invariante par compl\'etion, on peut supposer
$A$ et $B$ complets. Comme la platitude est invariante par extension
plate locale (donc fid\`element plate) de l'anneau de base, on est
ramen\'e en vertu de d\'efinition~\Ref{III.1.1} au cas o\`u $B$
est une alg\`ebre de s\'eries formelles sur~$A$. Mais alors en
tant que $A$-module, $B$ est isomorphe \`a un produit de~$A$-modules
isomorphes \`a $A$, donc (l'anneau de base $A$ \'etant
noeth\'erien) est $A$-plat comme produit de~$A$-modules plats.
\ifthenelse{\boolean{orig}}
{\upshape}{}

\medskip
Mettons-nous sous les conditions de~\Ref{III.1.1}. Comme les
extensions
\marginpar{60}
r\'esiduelles des composants locaux de~$B'$ sur~$A'$ sont triviales,
il s'ensuit que $L\otimes_k k'$ est une $k'$-alg\`ebre artinienne
dont les composants locaux ont des extensions r\'esiduelles
triviales (o\`u $L,\, k,\, k'$ sont les corps r\'esiduels de~$A,\,
B,\, A'$). Cette condition n\'ecessaire pour que l'extension finie
libre $A'$ satisfasse la condition \'enonc\'ee dans~\Ref{III.1.1}
est aussi suffisante, comme il r\'esulte aussit\^ot de
\Ref{III.1.4}~(i) et~\Ref{III.1.5} ci-dessous.

\begin{proposition}
\label{III.1.4}
Soit $A\to B$ un homomorphisme local d'anneaux locaux, \`a extension
r\'esiduelle finie
\ifthenelse{\boolean{orig}}{\ignorespaces}{et}
soit $A'$ une $A$-alg\`ebre finie locale sur~$A$, de sorte que
$B'=B\otimes_A A'$ est finie sur~$B$, donc un anneau semi-local
\'egalement noeth\'erien. \textup{(i)} Si $B$ est formellement lisse
sur~$A$, alors les localis\'es de B en ses id\'eaux maximaux sont
formellement lisses sur~$A'$. \textup{(ii)} La r\'eciproque est vraie
si $A'$ est libre sur~$A$.
\end{proposition}

On est aussit\^ot ramen\'e au cas o\`u $A$, $B$ sont complets.

(i) Soit $A''$ une extension finie libre locale de~$A$ telle que les
composants locaux de~$B''=B\otimes_A A''$ soient des alg\`ebres de
s\'eries formelles sur~$A''$. Faisant l'extension des scalaires
$A''\to A''\otimes_A A' \to A'''$, o\`u $A'''$ est un composant
local de~$A''\otimes_A A'$, on voit que les composants locaux de
$B''\otimes_{A''} A'''=B\otimes_A A'''$ sont des alg\`ebres de
s\'eries formelles sur~$A'''$. Or on a aussi
\ifthenelse{\boolean{orig}}
{$B\otimes_AA'''=(B\otimes_A A')\otimes_{A'} A'''=B\otimes_{A'} A'''$}
{$B\otimes_AA'''=(B\otimes_A A')\otimes_{A'} A'''=B'\otimes_{A'} A'''$},
d'autre
part comme $A''$ est libre sur~$A$, $A''\otimes_A A'$ est libre sur~$A'$ et il en est de m\^eme par suite de~$A'''$ qui en est facteur
direct, ce qui prouve que $B'$ est formellement lisse sur~$A'$.

(ii) Soit $A''$ une $A'$-alg\`ebre finie libre locale telle que les
composants locaux de~$B'\otimes_{A'} A'' = B\otimes_A A''$ soient des
alg\`ebres de s\'eries formelles sur~$A''$. Comme $A'$ est libre
sur~$A$, $A''$ l'est aussi, donc $B$ est formellement \emph{lisse}
sur~$A$.

\begin{proposition}
\label{III.1.5}
Soit $A\to B$ un homomorphisme local d'anneaux locaux, \`a extension
r\'esiduelle \emph{triviale}. Pour que $B$ soit formellement lisse
sur~$A$, il faut et il suffit que $\overline{B}$ soit isomorphe \`a
une alg\`ebre de s\'eries formelles sur~$\overline{A}$.
\end{proposition}

Il n'y a \`a prouver que la n\'ecessit\'e, et on peut supposer
$A$ et $B$ complets. Soit $\goth{m}$ ($\goth{n}$) l'id\'eal maximal
de~$A$ ($B$) et soient $t_1,\dots, t_n$ des \'el\'ements de
$\goth{n}$ qui d\'efinissent une base de l'espace vectoriel
$$
(\goth{n}/\goth{n}^2)/ \Im (\goth{m}/\goth{m}^2)=
\goth{n}/(\goth{n}^2 +\goth{m}B)
$$
Ces \'el\'ements d\'efinissent donc un homomorphisme de
$A$-alg\`ebres locales
$$
B_1= A[[t_1,\dots,t_n]] \to B
$$
prouvons
\marginpar{61}
que c'est un isomorphisme. Il suffit de prouver que pour toute
puissance $\goth{m}^q$ de~$\goth{m}$, on obtient un isomorphisme en
r\'eduisant mod~$\goth{m}^q$ (puisque $B_1$ et $B$ sont les
limites projectives des anneaux correspondants r\'eduits
mod~$\goth{m}^q$, $q$ variable). Comme $B$ et $B_1$ sont des
$A$-modules plats donc les gradu\'es associ\'es \`a la
filtration $\goth{m}$-adique s'obtiennent en tensorisant par $\gr
(A)$, sur~$k=A/\goth{m}$, les anneaux $B_1/\goth{m}B_1$
\resp $B/\goth{m}B$, on est ramen\'e \`a montrer que
$B_1/\goth{m}B_1\to B/\goth{m}B$ est un isomorphisme. Compte tenu
de~\Ref{III.1.3}, on est ainsi ramen\'e au cas o\`u $A$ est un
\emph{corps} $k$. D'autre part, si $A'$ est $A$-alg\`ebre finie
libre locale telle que $B\otimes_A A'$ soit une alg\`ebre de
s\'eries formelles sur~$A'$ (N.B. cette alg\`ebre est locale,
puisque l'extension r\'esiduelle de~$B$ sur~$A$ est triviale), pour
prouver que $B_1\to B$ est un isomorphisme, il suffit de prouver que
$B_1\otimes_A A'\to B\otimes_A A'$ l'est. Cela nous ram\`ene au cas
o\`u $B$ est d\'ej\`a une alg\`ebre de s\'eries formelles
(il fallait commencer par cette r\'eduction, avant de se ramener au
cas d'un corps de base). Mais alors $B$ est un anneau local
r\'egulier \`a corps de repr\'esentants $k$, et il est bien
connu (et imm\'ediat par consid\'eration des gradu\'es
associ\'es \`a la filtration $\goth{n}_1 $-adique et
$\goth{n}$-adique sur~$B_1$ et $B$) que $B_1 \to B$ est un
isomorphisme, ce qui ach\`eve la d\'emonstration.

\begin{corollaire}
\label{III.1.6}
Si $B$ est formellement lisse sur~$A$, alors il existe une
$A$-alg\`ebre finie locale $A'$ telle que les composants locaux de
$\overline{B}\otimes_{\overline{A}}
\overline{A'}=\overline{(B\otimes_A A')}$ soient isomorphes \`a des
alg\`ebres de s\'eries formelles sur~$\overline{A'}$.
\end{corollaire}
En effet, si $L/k$ est l'extension r\'esiduelle de~$B/A$, on
consid\`ere une extension $k'/k$, telle que les extensions
r\'esiduelles dans la $k'$-alg\`ebre $L\otimes_k k'$ soient
triviales. Soit $A'$ une alg\`ebre finie libre sur~$A$ telle que
$A'/\goth{m}A'=k'$ (on sait qu'il en existe, par exemple en se
ramenant de proche en proche au cas o\`u $k'/k$ est monog\`ene, et
alors on rel\`eve dans $A$ les coefficients du polyn\^ome minimal
d'un g\'en\'erateur de~$k'$ sur~$k$). Elle est locale. Alors
$B\otimes_A A'$ a en ses id\'eaux maximaux des extensions
r\'esiduelles triviales au-dessus de celle $k'$ de~$A'$, et on
ach\`eve \`a l'aide de~\Ref{III.1.5}.

\begin{corollaire}
\label{III.1.7}
Soit $A\to B$ un homomorphisme local d'anneaux locaux. Pour que $B$
soit formellement lisse sur~$A$, il faut et il suffit que $B$ soit
plat sur~$A$ et que $B/\goth{m}B$ soit formellement lisse sur~$k=A/\goth{m}$.
\end{corollaire}
Faisant une extension finie libre locale $A'$ convenable de~$A$ et
utilisant \Ref{III.1.4}~(ii), on est ramen\'e au cas o\`u
l'extension r\'esiduelle de~$B/A$ est triviale.
\marginpar{62}
On sait d'ailleurs par \Ref{III.1.4}~(i) et~\Ref{III.1.3} que les
conditions \'enonc\'ees sont n\'ecessaires. Pour la suffisance,
il suffit de remarquer que la d\'emonstration de~\Ref{III.1.5}
prouve, sous les hypoth\`eses faites ici, que $B$ est une
alg\`ebre de s\'eries formelles sur~$A$ (supposant $A$ et $B$
complets, ce qui est loisible).

\begin{remarque}
\label{III.1.8}
Il ne serait pas difficile de d\'evelopper, pour les homomorphismes
formellement lisses, l'analogue de toutes les propri\'et\'es des
morphismes lisses, \'etudi\'ees dans~\Ref{II}. Pour les
propri\'et\'es diff\'erentielles, cela demande cependant une
modification de la d\'efinition habituelle des diff\'erentielles
de K\"ahler (\cf I~\Ref{I.1}), les produits tensoriels
compl\'et\'es rempla\c cant les produits tensoriels
ordinaires. Nous nous contentons d'\'evoquer ici ces ab\^imes,
ce qui pr\'ec\`ede \'etant suffisant pour notre propos.

Il reste \`a faire le lien entre la lissit\'e formelle, et la
notion de lissit\'e d\'evelopp\'ee dans~\Ref{II} (et dont nous
n'avons encore fait nul usage):
\end{remarque}

\begin{proposition}
\label{III.1.9}
Soit $A\to B$ un homomorphisme local, $B$ \'etant localis\'ee
d'une $A$-alg\`ebre de type fini. Pour que $B$ soit lisse sur~$A$,
il faut et il suffit qu'il soit formellement lisse sur~$A$.
\end{proposition}

Utilisant~\Ref{III.1.7} et II~\Ref{II.2.1}, on est ramen\'e au cas
o\`u $A$ est un corps.

Utilisant~\Ref{III.1.4}~(ii) et II~\Ref{II.4.13} une extension convenable
$k'$ de~$k$ nous ram\`ene au cas o\`u l'extension r\'esiduelle
pour $B/k$ est triviale. En vertu de~\Ref{III.1.5}
(\resp II~\Ref{II.5.2}) $B$ est alors lisse sur~$k$ (\resp formellement
\emph{lisse sur} $k$) si et seulement si $B$ est un anneau local
r\'egulier (\resp son compl\'et\'e est une alg\`ebre de
s\'eries formelles sur~$k$). Or il est bien connu que ces deux
conditions sont \'equivalentes (l'extension r\'esiduelle \'etant
triviale).

\section[Propri\'et\'e de rel\`evement]{Propri\'et\'e de rel\`evement caract\'eristique des
ho\-mo\-mor\-phis\-mes formellement lisses}
\label{III.2}

\begin{theoreme}
\label{III.2.1}
Soit $A \to B$ un homomorphisme local d'anneaux locaux d\'efinissant
une extension r\'esiduelle finie. Les conditions suivantes sont
\'equivalentes:
\begin{enumerate}
\item[(i)]
\marginpar{63}
$B$ est formellement lisse sur~$A$
\item[(ii)] Pour tout homomorphisme local $A \to C$, o\`u $C$ est un
anneau local \emph{complet}, tout id\'eal $J$ de~$C$ contenu dans
\ifthenelse{\boolean{orig}}
{$r(C)$}
{$\goth{r}(C)$},
et tout $A$-homomorphisme local $B \to C/J$, il existe un
$A$-homomorphisme (n\'ecessairement local) $B \to C$ qui le
rel\`eve.
\item[(iii)] Pour toute $A$-alg\`ebre $C$ (pas n\'ecessairement un
anneau noeth\'erien) tout id\'eal nilpotent $J$ de~$C$, et tout
$A$-homomorphisme $B \to C$ continu (\ie s'annulant sur une puissance
\ifthenelse{\boolean{orig}}
{de~$r(B)$}
{de~$\goth{r}(B)$}), il existe un $A$-homomorphisme $B \to C$
(n\'ecessairement continu lui aussi) qui le rel\`eve.
\item[(iv)] M\^eme \'enonc\'e que \textup{(ii)} et \textup{(iii)}, mais $C$
\'etant un anneau local artinien, fini au-dessus de~$A$.
\item[(iv~bis)] Comme \textup{(iv)}, mais $J$ \'etant de plus de carr\'e
nul.
\end{enumerate}
\end{theoreme}

\begin{remarquestar}
Pour la suite de cet expos\'e, nous nous servirons seulement de
l'implication (iv)$\To$(i) ou (iv~bis)$\To$(i);
l'implication directe (i)$\To$(ii) sera prouv\'ee par une
autre m\'ethode au \no suivant lorsque $B$ est localis\'ee d'une
alg\`ebre de type fini sur~$A$. Rappelons que dans la \og bonne\fg
th\'eorie des th\'eor\`emes de Cohen\footnote{\Cf EGA
$0_{\textup{IV}}$ 19.3, 19.8}, la propri\'et\'e (ii) ou (iii)
devient la d\'efinition des homomorphismes formellement lisses,
alors que \Ref{III.1.1} devient une propri\'et\'e caract\'eristique
valable seulement dans le cas d'une extension r\'esiduelle finie. On
fera attention que des propri\'et\'es (ii) et (iii) aucune n'est
plus g\'en\'erale que l'autre. On pourrait donner une
propri\'et\'e (\'equivalente) qui les coiffe toutes deux, en
introduisant un anneau lin\'eairement topologis\'e $C$,
\emph{s\'epar\'e} et \emph{complet}, un id\'eal \emph{ferm\'e
topologiquement nilpotent} de~$C$, et un homomorphisme continu $A \to
C$ (faisant donc de~$C$ une $A$-alg\`ebre topologique); nous
laisserons cette modification au lecteur.
\end{remarquestar}

\subsubsection*{D\'emonstration de~\Ref{III.2.1}}
Nous prouverons (i)$\To$(iii)$\To$(ii), puis (iv)$\To$(i). Comme (ii)$\To$(iv) est trivial, et
l'\'equivalence de (iv) et (iv~bis) se voit par une r\'ecurrence
imm\'ediate sur l'entier $n$ tel que $J^n = 0$, cela ach\`evera la
d\'emonstration.

(i)$\To$(iii). Une r\'ecurrence imm\'ediate nous
ram\`ene au cas o\`u $J^2 = 0$. Comme $C$ est fini sur~$A$, il
existe une puissance $\goth{m}^q$ de l'id\'eal maximal de~$A$ qui annule
$C$. Divisant par $\goth{m}^q$, et notant que $B/\goth{m}^qB$ est
encore formellement lisse
\marginpar{64}
sur~$A/\goth{m}^q$ par \Ref{III.1.4}~(i), on peut supposer $A$
artinien. Comme $B$ est plat sur~$A$ par~\Ref{III.1.3}, $B$ \emph{est
libre sur} $A$ puisque $A$ est artinien. Donc il existe un
\emph{homomorphisme de} $A$-\emph{modules}
$$
w \colon B \to C
$$
qui rel\`eve l'homomorphisme donn\'e $u \colon B \to C/J$. Posons
$$
f(x,y) = w(xy) - w(x)w(y) \qquad (x,y \in B)
$$
alors $f(x,y) \in J$, et $f$ est donc une application
$A$-bilin\'eaire de~$B \times B$ dans $J$. Pour qu'il existe un
rel\`evement $v \colon B \to C$ de~$u$ qui soit un homomorphisme
d'alg\`ebres, il faut et il suffit qu'il existe une application
$A$-lin\'eaire $g \colon B \to J$ telle que $v = w + g$ soit un
homomorphisme d'alg\`ebres, ce qui s'\'ecrit
\begin{align*}
g(1) &= 1 - w(1) \\
g(x,y) - u(x)g(y) -u(y)g(x) &= -f(x,y) \qquad (\text{pour } x,y \in
B)
\end{align*}
C'est l\`a un syst\`eme d'\'equations \emph{lin\'eaires} dans
$\Hom_A(B,J)$, \`a seconds membres dans $J$, donc il a une solution
si et seulement si le syst\`eme correspondant dans $\Hom_A(B,J)
\otimes_A A'$, \`a seconds membres dans $J' = J \otimes_A A'$, a une
solution, --- $A'$ d\'esignant une alg\`ebre fid\`element plate
sur~$A$. Or soit $A'$ une alg\`ebre finie et libre sur~$A$, locale,
telle que $B' = B \otimes_A A'$ soit une alg\`ebre de s\'eries
formelles sur~$A'$ (N.B. on peut dans notre d\'emonstration supposer
$A$ et $B$ complets, comme on constate aussit\^ot). Comme $A'$ est
libre de type fini sur~$A$, on a
$$
\Hom_A(B,J) \otimes_A A' = \Hom_{A'}(B',J')
$$
et on constate que le syst\`eme d'\'equations obtenu dans
$\Hom_{A'}(B',J')$ est celui qui d\'etermine les homomorphismes de
$A'$-alg\`ebres $B' \to C' = C \otimes_A A'$ qui rel\`event
l'homomorphisme $u' \colon B' \to C'/J'$ d\'eduit de~$u$ par
extension des scalaires en \og corrigeant\fg par un homomorphisme de
$A'$-modules $g' \colon B' \to J'$ l'homomorphisme de~$A'$-modules $w'
\colon B' \to C'$ d\'eduit de~$w$ par extension des
scalaires. (Noter que $B$ engendre $B'$ comme $A'$-module). On est
ainsi ramen\'e \`a prouver (iii) lorsque $B$ est une
\emph{alg\`ebre de s\'eries formelles} sur~$A$, $B = A[[t_1,
\dots, t_n]]$. Relevons alors de fa\c con quelconque les images dans
$C/J$ des $t_i$ en des \'el\'ements $z_i$ de~$C$. Comme les $z_i$~mod~$J$ sont nilpotents ($u \colon B \to C/J$ \'etant continu) il en
est de m\^eme des $z_i$ (puisque $J$ est nilpotent), donc les $z_i$
d\'efinissent un homomorphisme continu de~$A$-alg\`ebres
topologiques de~$B$ dans $C$ discret, relevant \'evidemment $u$,
cqfd.

(iii)$\To$(ii).
\marginpar{65}
Soit $\goth{n}$ l'id\'eal maximal de~$C$, et
pour tout entier $q > 0$, soit
$$
C_q = C/\goth{n}^q \text{\quoi, \quad} J_q = (J +
\goth{n}^q)/\goth{n}^q
$$
donc $C_q/J_q$ s'identifie \`a une alg\`ebre quotient de~$C/J$,
d'autre part l'homomorphisme compos\'e $u_q \colon B \to C/J \to
C_q/J_q$ est continu de~$B$ dans $C_q/J_q$ discret, et $J_q$ est un
id\'eal nilpotent dans $C_q$. On construit alors de proche en proche
des $A$-homomorphismes
$$
v_q \colon B \to C_q
$$
tels que (a) $v_q$ rel\`eve $u_q$ et (b) $v_q$ rel\`eve
$v_{q-1}$. La possibilit\'e de la r\'ecurrence se v\'erifie
ais\'ement, car comme
$$
u_q \colon B \to C/(J + \goth{n}^q) \text{\quad et \quad} v_{q-1}
\colon B \to C/\goth{n}^{q-1}
$$
d\'efinissent le m\^eme homomorphisme
$$
B \to C/((J + \goth{n}^q) + \goth{n}^{q-1}) = C/(J + \goth{n}^{q-1}) =
C_{q-1}/J_{q-1}
$$
\`a savoir $u_{q-1}$, ils d\'efinissent un homomorphisme
$$
B \to C/J_q' \text{\quad o\`u \quad} J_q' = (J + \goth{n}^q)\cap
\goth{n}^{q-1} \supset \goth{n}^q
$$
(dont ils proviennent l'un et l'autre par r\'eduction). On est donc
ramen\'e \`a relever un homomorphisme $B \to C/J_q'$ de~$B$ dans
un quotient de~$C_q$ par un id\'eal $J_q'/\goth{n}^q$ contenu dans
$J_q$, donc nilpotent, et cela est possible d'apr\`es
l'hypoth\`ese (iii).

Ceci fait, les $v_q$ d\'efinissent un homomorphisme de~$B$ dans la
limite projective $C$ des $C_q$. Comme $J$ est ferm\'e, $J$ est la
limite projective des $J_q$, donc $v$ rel\`eve $u$, cqfd.

(iv)$\To$(i). On constate d'abord aussit\^ot que si (iv)
est v\'erifi\'e, (iv) reste v\'erifi\'e pour las composants
locaux de~$B \otimes_A A'$ sur~$A'$, si $A'$ est une alg\`ebre
locale finie sur~$A$. Prenant $A'$ libre sur~$A$ et telle que les
extensions r\'esiduelles de~$B'$ au-dessus de~$A'$ soient triviales,
on est ramen\'e gr\^ace \`a~\Ref{III.1.4} (ii) au cas o\`u
l'extension r\'esiduelle de~$B$ sur~$A$ est triviale. Nous allons
prouver alors le r\'esultat un peu pr\'ecis:

\begin{corollaire}
\label{III.2.2}
Sous les conditions de~\Ref{III.2.1}, supposons de plus l'extension
r\'esiduelle de~$B$ au-dessus de~$A$ triviale. Alors les conditions
\'equivalentes
\marginpar{66}
de~\Ref{III.2.1} \'equivalent aussi aux deux conditions suivantes
(en supposant dans \textup{(v)} $A$ et $B$ complets):
\begin{enumerate}
\item[(iv~ter)] Comme \textup{(iv)}, mais l'anneau local artinien $C$ fini sur~$A$ \'etant restreint \`a avoir une extension r\'esiduelle
triviale (et de plus, si on y tient, l'id\'eal $J$ \'etant de
carr\'e nul).
\item[(v)] Il existe un $A$-homomorphisme local (o\`u $n =
\dim{\goth{n}/(\goth{n}^2+\goth{m}B})$)
$$
u \colon B \to B_1 = A[[t_1, \dots, t_n]]
$$
induisant un \emph{isomorphisme}
$$
\goth{n}/(\goth{n}^2+\goth{m}B) \isomto
\goth{n}_1/(\goth{n}_1^2+\goth{m}B_1)
$$
o\`u $\goth{n}$ ($\goth{n}_1$) est l'id\'eal maximal de~$B$
($B_1$), $\goth{m}$ celui de~$A$.
\end{enumerate}
\end{corollaire}

\subsubsection*{D\'emonstration} Comme (iv~bis) implique \'evidemment (iv~ter) --- en faisant abstraction du canular de l'id\'eal de carr\'e nul ---,
il suffira de prouver (iv~ter)$\To$(v)$\To$(i).

(iv~ter)$\To$(v). Choisissons une base $a_1,\dots,a_n$ de
$\goth{n}/(\goth{n}^2+\goth{m}B)$, ce qui d\'efinit donc un
homomorphisme local de~$A$-alg\`ebres
$$
B \to B_1/(\goth{n}_1^2+\goth{m}B_1) = k[t_1, \dots, t_n]/(t_1, \dots,
t_n)^2
$$
qu'on peut relever de proche en proche, en vertu de (iv~ter) en des
homomorphismes de~$A$-alg\`ebres de~$B$ dans $B_1/\goth{n}_1^2$,
$B_1/\goth{n}_1^3$ etc, d'o\`u en passant \`a la limite projective
l'homomorphisme $B \to B_1$ ayant la propri\'et\'e voulue.

(v)$\To$(i). Comme dans le diagramme commutatif
$$
\ifthenelse{\boolean{orig}}
{\xymatrix{ \goth{m}/\goth{m}^2 \ar[r] & \goth{n}/\goth{n}^2 \ar[r]
\ar[d] & \goth{n}/(\goth{n}^2+\goth{m}B) \ar[r] \ar[d] & 0 \\
\goth{m}/\goth{m}^2 \ar[r] & \goth{n}_1/\goth{n}_1^2 \ar[r] &
\goth{n}_1/(\goth{n}_1^2+\goth{m}B_1) \ar[r] & 0}}
{\xymatrix{ \goth{m}/\goth{m}^2 \ar[r] \ar[d] & \goth{n}/\goth{n}^2
\ar[r] \ar[d] & \goth{n}/(\goth{n}^2+\goth{m}B) \ar[r] \ar[d] & 0 \\
\goth{m}/\goth{m}^2 \ar[r] & \goth{n}_1/\goth{n}_1^2 \ar[r] &
\goth{n}_1/(\goth{n}_1^2+\goth{m}B_1) \ar[r] & 0}}
$$
les deux lignes sont exactes, et les fl\`eches verticales
extr\^emes surjectives, la fl\`eche m\'ediane est surjective et
il s'ensuit ($B$ \'etant complet) que $B \to B_1$ est
\emph{surjectif}. Soient $x_i$ ($1 \leq i \leq n$) des
\'el\'ements de~$B$ qui rel\`event les $t_i$. Ils
d\'efinissent donc un homomorphisme de~$A$-alg\`ebre $B_1 \to B$,
qui sera surjectif pour la m\^eme raison que $u$, et dont le
compos\'e avec $u$ est l'identit\'e par construction. Donc $B_1
\to B$ est aussi injectif, et est par suite un isomorphisme. On trouve
donc

\begin{corollaire}
\label{III.2.3}
Sous
\marginpar{67}
les conditions de \Ref{III.2.2}~\textup{(v)}, $u$ est n\'ecessairement
un isomorphisme.
\end{corollaire}

Cela ach\`eve de prouver que $B$ est formellement lisse sur~$A$. On
a d'ailleurs en m\^eme temps retrouv\'e~\Ref{III.1.5} (mais il n'y
a gu\`ere de m\'erite \`a \c ca).

\section[Prolongement infinit\'esimal local des morphismes]{Prolongement infinit\'esimal local des morphismes dans un
$S$-sch\'ema lisse}
\label{III.3}

\begin{theoreme}
\label{III.3.1}
Soit $f \colon X \to Y$ un morphisme localement de type fini.
Conditions \'equivalentes:
\begin{enumerate}
\item[(i)] $f$ est lisse.
\item[(ii)] Pour tout pr\'esch\'ema $Y'$ sur~$Y$, tout
sous-pr\'esch\'ema ferm\'e $Y'_0$ de~$Y'$ ayant m\^eme espace
sous-jacent que $Y'$, tout $Y$-morphisme $g_0 \colon Y'_0 \to X$ et
tout $z \in Y'_0$, il existe un voisinage ouvert $U$ de~$z$ dans $Y'$
et un prolongement $g$ de~$g_0|Y'_0 \cap U$ en un $Y$-morphisme $U \to
X$.
\item[(ii~bis)] Pour $Y'$, $Y'_0$ et $z$ comme dans \textup{(ii)}, posant
$X' = X \times_Y Y'$, $X'_0 = X \times_Y Y'_0$, toute section de
$X'_0$ sur~$Y'_0$ se prolonge en une section de~$X'$ au-dessus d'un
voisinage ouvert $U$ de~$z$.
\item[(iii)] Pour tout $Y$-sch\'ema $Y'$, spectre d'un anneau artinien
local fini sur quelque $\mathcal{O}_y$ ($y \in Y$), tout
sous-pr\'esch\'ema ferm\'e non vide~$Y'_0$ de~$Y'$, et tout
$Y$-morphisme $g_0 \colon Y'_0 \to X$, il existe un $Y$-morphisme $g
\colon Y' \to X$ qui prolonge $g_0$.
\item[(iii~bis)] Pour tout $Y'$, $Y'_0$ comme dans \textup{(iii)},
posant $X' = X \times_Y Y'$, $X'_0 = X \times_Y Y'_0$, toute section
de~$X'_0$ sur~$Y'_0$ se prolonge en une section de~$X'$ sur~$Y'$.
\end{enumerate}
\end{theoreme}

\subsubsection*{D\'emonstration} L'\'equivalence de (ii) et (ii~bis) d'une part, de
(iii) et (iii~bis) d'autre part, est triviale, ainsi que l'implication
(ii)$\To$(iii). Il reste donc \`a prouver (i)$\To$(ii) et (iii)$\To$(i).

(i)$\To$(ii). Soit $x = g_0(z)$. Rempla\c cant $X$ par un
voisinage ouvert convenable de~$x$, et $Y'$ par le pr\'esch\'ema
induit sur l'ouvert image r\'eciproque de ce dernier par~$g_0$, on
peut supposer que $X$ est \'etale au-dessus de~$Y[t_1, \dots, t_n]$.
Consid\'erons le $Y$-morphisme compos\'e $Y'_0 \to X \to Y[t_1,
\dots, t_n]$, il est d\'efini par $n$ sections du faisceau
$\mathcal{O}_{Y'_0}$, qui peuvent donc se prolonger au voisinage de
$z$ en des sections de~$\mathcal{O}_{Y'}$, donc on peut supposer que
\marginpar{68}
le morphisme en question a \'et\'e prolong\'e en un
$Y$-morphisme $Y' \to Y'_0$. En vertu de (I~\Ref{I.5.6}) il existe
alors un unique $Y$-morphisme $g \colon Y' \to X$ qui rel\`eve le
pr\'ec\'edent, et prolonge en m\^eme temps $g_0$, cqfd.

(iii)$\To$(i). Comme l'ensemble des points o\`u $f$ est
lisse est ouvert, il suffit de prouver qu'il contient tout $x \in X$
qui est \emph{ferm\'e} dans sa fibre. Soit $y = f(x)$, alors
$\mathcal{O}_x$ est une alg\`ebre sur~$\mathcal{O}_y$, localis\'ee
d'une alg\`ebre de type fini, \`a extension r\'esiduelle finie.
D'autre part, l'hypoth\`ese (iii) implique que tout homomorphisme de
$\mathcal{O}_x$ dans une alg\`ebre $A/J$, o\`u $A$ est une
alg\`ebre artinienne locale finie sur~$\mathcal{O}_y$, et $J$ un
id\'eal contenu dans son radical, se rel\`eve en un homomorphisme
de~$\mathcal{O}_x$ dans l'alg\`ebre $A$ (compte tenu qu'un morphisme
d'un $\Spec(B)$, $B$ anneau local, dans $X$ est d\'etermin\'e
biunivoquement par un homomorphisme local d'un $\mathcal{O}_x$ ($x \in
X$) dans $B$). Par~\Ref{III.2.1} il s'ensuit que $\mathcal{O}_x$ est
formellement lisse sur~$\mathcal{O}_y$, donc lisse sur~$\mathcal{O}_y$
en vertu de~\Ref{III.1.9}.

\begin{corollaire}
\label{III.3.2}
Soit $f\colon X\to Y$ comme dans~\Ref{III.3.1}. Conditions
\'equivalentes:
\begin{enumerate}
\item[(i)] $f$ est \'etale
\item[(ii)] On a la condition~\textup{(ii)} de~\Ref{III.3.1} avec
{\emph{unicit\'e}} du prolongement $g$ de~$g_0$ \`a $U$.
\item[(iii)] On a la condition~\textup{(iii)} de~\Ref{III.3.1} avec
{\emph{unicit\'e}} de~$g$.
\end{enumerate}
\end{corollaire}

Il suffit de noter, dans la d\'emonstration de (i)$\To$(ii)
ci-dessus, qu'on ne peut avoir unicit\'e (quand $Y'_0$ n'est pas
identique \`a~$Y'$ au voisinage de~$z$) que si $n=0$ (condition
qu'on sait \^etre suffisante).

\begin{corollaire}
\label{III.3.3}
Soit $X$ un pr\'esch\'ema localement de type fini sur un anneau
local {\emph{complet}} $A$, $y$ le point ferm\'e de~$Y=\Spec(A)$ et
$x$ un point de~$f^{-1}(y)$ {\emph{rationnel}}
\ifthenelse{\boolean{orig}}
{sur~$\kappa(y)$}
{sur~$\kres(y)$}. Si $X$
est {\emph{lisse sur~$A$ en~$x$}}, alors il existe une section $s$ de
$X$ sur~$Y$ \og passant par~$x$\fg \ie telle que $s(y)=x$.
\end{corollaire}

En particulier, si $X$ est lisse sur~$A$, alors l'application
naturelle
$$
\Gamma(X/Y)\to\Gamma(X\otimes_A k/k)
$$
des sections de~$X$ sur~$Y$ dans l'ensemble des points de la fibre
$f^{-1}(y)=X\otimes_A k$ rationnels sur~$k$, est surjective. Ce fait
\'etait surtout bien connu et utilis\'e lorsque $A$ est un anneau
de valuation discr\`ete et $X$
\marginpar{69}
est propre sur~$A$ (en fait, projectif sur~$A$), auquel cas les
sections de~$X$ sur~$Y$ (\ie les \og points de~$X$ \`a valeurs
dans~$A$\fg) s'identifient aussi aux sections rationnelles, \ie aux
points de~$X\otimes_A K=X_K$ (qui est un sch\'ema propre et simple
sur~$K$) \`a valeurs dans $K=$ corps des fractions
\ifthenelse{\boolean{orig}}
{de~$A$}
{de~$A$,}
\ie aux
points de~$X$ rationnels sur~$K$.

\section{Prolongement infinit\'esimal local des $S$-sch\'emas lisses}
\label{III.4}

\begin{theoreme}
\label{III.4.1}
Soient $Y$ un pr\'esch\'ema localement noeth\'erien, $Y_0$ un
sous-pr\'e\-sch\'ema ferm\'e ayant m\^eme espace sous-jacent,
$X_0$ un $Y_0$-pr\'esch\'ema lisse, $x$ un point de~$X_0$. Alors
il existe un voisinage ouvert $U_0$ de~$x$, un pr\'esch\'ema $X$
lisse sur~$Y$, et un $Y_0$-isomorphisme:
$$
h\colon U_0\isomto X\times_Y Y_0.
$$
De plus, si $(U'_0,X',h')$ est une autre solution de ce probl\`eme,
alors \og elle est isomorphe \`a la premi\`ere au voisinage
de~$x$\fg.
\end{theoreme}

On laisse au lecteur de pr\'eciser ce qu'on veut dire par l\`a. On
peut noter que pour $U_0$ donn\'e, une solution du
\ifthenelse{\boolean{orig}}{Pb}{probl\`eme} pos\'e
revient \`a la donn\'ee, sur~$U_0$, d'un faisceau d'alg\`ebres
$\cal{B}$ sur~$f_0^{-1}(\cal{O}_Y)$ (o\`u $f_0$ est l'application
continue sous-jacente au morphisme structural $U_0\to Y_0$), muni d'un
homomorphisme d'anneaux $\cal{B}\to \cal{O}_{U_0}$ compatible avec
l'homomorphisme $f^{-1}(\cal{O}_Y)\to f^{-1}(\cal{O}_{Y_0})$, tels que

(a) cet homomorphisme induit un isomorphisme
$$
\cal{B}\otimes_{f^{-1}(\cal{O}_Y)} f^{-1}(\cal{O}_{Y_0})\isomto
\cal{O}_{Y_0};
$$

(b) $U_0$ muni de~$\cal{B}$ devient un $Y$-pr\'esch\'ema lisse.

De cette fa\c con, le sens pr\'ecis de l'assertion d'unicit\'e
locale devient particuli\`erement \'evident.

\subsubsection*{D\'emonstration} On peut supposer d\'ej\`a que $X_0$ est
\'etale \hbox{sur un $Y_0[t_1,\dots,t_n]=Y'_0$}. Or ce dernier peut
\^etre consid\'er\'e comme un sous-pr\'esch\'ema ferm\'e
de~$Y'{=} Y[t_1,\dots,t_n]$, ayant m\^eme espace sous-jacent. Par
(I~\Ref{I.8.3}), il existe un $X$ \'etale sur~$Y'$, et un
$Y'_0$-isomorphisme $X\times_{Y'}Y'_0\isomto X'$. On a gagn\'e pour
l'existence. Pour l'unicit\'e, on utilise la
propri\'et\'e~\Ref{III.3.1}~(ii)
\marginpar{70}
des morphismes lisses, en tenant compte du

\begin{lemme}
\label{III.4.2}
Soient $Y$ un pr\'esch\'ema, $Y_0$ un sous-pr\'esch\'ema
ferm\'e d\'efini par un faisceau d'id\'eaux $\cal{J}$ localement
nilpotent, $X$ et $X'$ des $Y$-pr\'esch\'emas, $u\colon X\to X'$
un $Y$-morphisme. On suppose $X$ plat sur~$Y$. Pour que $u$ soit un
isomorphisme, il faut et suffit que $u_0\colon X\times_Y Y_0\to
X'\times_Y Y_0$ le soit.
\end{lemme}

D\'emonstration facile, en passant au cas affine et regardant les
gradu\'es associ\'es. On notera d'ailleurs que l'\'enonc\'e
analogue, obtenu en rempla\c cant \og isomorphisme\fg par \og immersion
ferm\'ee\fg, est \'egalement valable, et sans hypoth\`ese de
platitude.

\begin{remarque}
\label{III.4.3}
Il est essentiel de noter que le prolongement local $X$ obtenu
dans~\Ref{III.4.1} {\emph{n'est pas canonique}}, en d'autres termes
l'isomorphisme local entre deux solutions n'est pas unique, \ie il
existe en g\'en\'eral des $Y$-automorphismes non triviaux de~$X$
induisant l'identit\'e sur le sous-pr\'esch\'ema ferm\'e
$X_0=X\times_Y Y_0$. C'est pour cela qu'il faut s'attendre, pour la
construction de prolongements infinit\'esimaux {\emph{globaux}} de
pr\'esch\'emas
\ifthenelse{\boolean{orig}}
{simples}
{lisses}, \`a l'existence d'une obstruction de
nature cohomologique, qui sera pr\'ecis\'ee plus bas (\no \Ref{III.6}).
\end{remarque}

\section{Prolongement infinit\'esimal global des morphismes}
\label{III.5}

Soient $T$ un espace topologique, $\cal{G}$ un faisceau de groupes
sur~$X$,
\ifthenelse{\boolean{orig}}{$P$}{$\cal{P}$}
un faisceau d'ensembles sur~$T$ sur lequel $\cal{G}$ op\`ere (\`a
droite, pour fixer les id\'ees). On dit que $\cal{P}$ est
\emph{formellement principal homog\`ene}
\index{formellement principal homog\`ene (faisceau)|hyperpage}%
sous~$\cal{G}$, si l'homomorphisme bien connu
$$
\cal{G}\times \cal{P}\to \cal{P}\times \cal{P}
$$
de faisceaux d'ensembles, d\'eduit des op\'erations de~$\cal{G}$
sur~$\cal{P}$, est un {\emph{isomorphisme}}. Il revient au m\^eme de
dire que pour tout $x\in T$, $\cal{P}_x$ est {\emph{vide ou un espace
principal homog\`ene}} sous le groupe ordinaire $\cal{G}_x$, ou
aussi que pour tout ouvert $U$ de~$T$, $\cal{P}(U)$ est vide ou un
espace principal homog\`ene sous le groupe ordinaire
$\cal{G}(U)$. On dit que
\ifthenelse{\boolean{orig}}{$P$}{$\cal{P}$}
est un {\emph{faisceau principal homog\`ene}}
\index{principal homog\`ene (faisceau)|hyperpage}%
\ifthenelse{\boolean{orig}}{sous~$G$}{sous~$\cal{G}$}
s'il l'est formellement et si en plus les $\cal{P}_x$ sont non vides
(en d'autres termes, si {\emph{tous}} les $\cal{P}_x$ sont des espaces
principaux homog\`enes,
\marginpar{71}
donc non vides, sous les $\cal{G}_x$)\footnote{Il semble
pr\'ef\'erable d'adopter le terme plus court et plus parlant de
\ifthenelse{\boolean{orig}}{\og {\emph{torseur sous}}~$G$\fg,}
{\og {\emph{torseur sous}}~$\cal{G}$\fg,}
introduit dans la th\`ese de J\ptbl \textsc{Giraud}.}. Rappelons que l'ensemble
des classes (\`a isomorphisme pr\`es) de faisceaux principaux
homog\`enes sous~$\cal{G}$ s'identifie \`a l'ensemble de
cohomologie $\H^1(T,\cal{G})$, qui est aussi le groupe de cohomologie
usuel de~$T$ \`a coefficients dans $\cal{G}$ lorsque $\cal{G}$ est
commutatif. On a ainsi, pour tout $\cal{P}$ principal homog\`ene,
une classe caract\'eristique $c(\cal{P})\in\H^1(T,\cal{G})$, dont la
trivialit\'e est n\'ecessaire et suffisante pour que $\cal{P}$
soit trivial (\ie isomorphe \`a~$\cal{G}$, sur lequel $\cal{G}$
op\`ere par translations \`a droite), ou encore pour que $\cal{P}$
ait une section.

\begin{proposition}
\label{III.5.1}
Soient $S$ un pr\'esch\'ema, $X$ et $Y$ des pr\'esch\'emas
sur~$S$, $Y_0$ un sous-pr\'esch\'ema ferm\'e de~$Y$ d\'efini
par un Id\'eal $\cal{J}$ sur~$Y$ {\emph{de carr\'e nul}}. Soit
$g_0$ un $S$-morphisme de~$Y_0$ dans~$X$, et $\cal{P}(g_0)$ le
faisceau sur~$Y$ dont les sections sur un ouvert~$U$ sont les
prolongements $g\colon U\to X$ de~$g_0|U\cap Y_0$ en un $S$-morphisme
$g$. Alors $\cal{P}(g_0)$ est un faisceau {\emph{formellement
principal homog\`ene}} (de fa\c con naturelle) sous le faisceau
en groupes commutatif
$$
\cal{G}=\SheafHom_{\cal{O}_{Y_0}}(g_0^*(\it{\Omega}^1_{X/S}),\cal{J})
$$
\end{proposition}

Posons $\cal{P}=\cal{P}(g_0)$. Nous devons d\'efinir pour tout
ouvert $U$ de~$Y$ une application
$$
\cal{P}(U)\times \cal{G}(U)\to \cal{P}(U)
$$
de fa\c con que (a) pour $g\in \cal{P}(U)$ fix\'e, l'application
$s\mto gs$ de~$\cal{G}(U)$ dans~$\cal{P}(U)$ est bijective (b)
$\cal{P}(U)$ devient un ensemble \`a groupe d'op\'erateurs
$\cal{G}(U)$ (c) les applications pr\'ec\'edentes sont compatibles
avec les op\'erateurs de restriction pour un ouvert $V\subset U$. La
v\'erification de (c) sera triviale, on peut donc pour simplifier
supposer $U=Y$. La v\'erification de (b) (qui est, si on veut, de
nature locale) sera laiss\'ee au lecteur, nous nous bornerons donc,
pour un $g\in \cal{P}(Y)$ donn\'e, de d\'efinir une bijection
naturelle de~$\cal{G}(Y)$ sur~$\cal{P}(Y)$. Donc on suppose donn\'e
d\'ej\`a un $S$-morphisme $g\colon Y\to X$, et on cherche une
bijection canonique
\begin{equation*}
\label{eq:III.5.1.*}
\tag{$*$}
\Hom_{\cal{O}_{Y_0}}(g_0^*(\it{\Omega}^1_{X/S}),\cal{J})\isomto
\cal{P}(Y)
\end{equation*}
o\`u $\cal{P}(Y)$ est l'ensemble des $S$-morphismes $g'$ de~$Y$
dans~$X$ induisant le m\^eme morphisme $g_0\colon Y_0\to X$ que~$g$.
La donn\'ee d'un tel $g'$ est \'equivalente \`a la donn\'ee
d'un $S$-morphisme $h\colon Y\to X\times_S X$ tel que
${\pr}_1\circ h=g$, et $h\circ i=(g_0,g_0)$, o\`u
\marginpar{72}
\ifthenelse{\boolean{orig}}
{${\pr}_1\colon X\times_S X$}
{${\pr}_1\colon X\times_S X \to X$}
est la premi\`ere projection, $i\colon Y_0\to Y$ l'immersion
canonique, et $(g_0,g_0)\colon Y_0\to X\times_S X$ est le morphisme
$\Delta_{X/S}g_0$ de
\ifthenelse{\boolean{orig}}{composante}{composantes}
$g_0,g_0$:
$$
\xymatrix@C=2.2cm{ X\times_S X \ar[d]_{{\pr}_1} & Y_0
\ar[l]_-{h_0=(g_0,g_0)} \ar[d]^i \\ X & Y \ar[l]_-g
\ar@{-->}[lu]_-h }
$$
Comme $h_0$ se factorise par l'immersion diagonale $\Delta_{X/S}$ et
que $Y$ est dans le voisinage infinit\'esimal d'ordre $1$ de~$Y_0$
(\ie $\cal{J}^2=0$) les $h$ cherch\'es se factorisent
n\'ecessairement (de fa\c con unique) par le voisinage
infinit\'esimal du premier ordre de la diagonale, lequel s'identifie
(en tant que $X$-pr\'esch\'ema gr\^ace \`a ${\pr}_1$)
au spectre $X'$ du faisceau d'alg\`ebres
$\cal{O}_X+\it{\Omega}^1_{X/S}$ (o\`u le deuxi\`eme terme est
consid\'er\'e comme un Id\'eal de carr\'e nul), le morphisme
diagonal $X\to X'$ correspondant \`a l'augmentation canonique de ce
faisceau d'alg\`ebres. Posons $Y'=X'\times_X Y$, $Y'_0=Y'\times_Y
Y_0=X'\times_X Y_0$, de sorte que les $h$ cherch\'es sont en
correspondance biunivoque avec les sections $u$ de~$Y'$ sur~$Y$ qui
prolongent une section donn\'ee $u_0$ de~$Y'_0$ sur~$Y_0$. On peut
d'ailleurs identifier $Y'$ au spectre du faisceau d'alg\`ebres sur~$Y$ $\cal{A}=g^*(\cal{O}_X+\it{\Omega}^1_{X/S})=\cal{O}_Y+
g^*(\it{\Omega}^1_{X/S})$, et $Y'_0$ au faisceau d'alg\`ebres
$\cal{A}_0=\cal{A}\otimes_{\cal{O}_Y}\cal{O}_{Y_0}=\cal{O}_{Y_0}+
g^*_0(\it{\Omega}^1_{X/S})$, alors la section $u_0$ est celle
d\'efinie par l'augmentation canonique de~$A_0$
dans~$\cal{O}_{Y_0}$. Donc $\cal{P}(Y)$ s'identifie \`a l'ensemble
des homomorphismes d'alg\`ebres $\cal{A}\to \cal{O}_Y$ qui induisent
l'augmentation canonique $\cal{A}_0\to \cal{O}_{Y_0}$. Or les
homomorphismes d'alg\`ebres $\cal{A}\to \cal{O}_Y$
\ifthenelse{\boolean{orig}}{correspondant}{correspondent}
biunivoquement aux homomorphismes de Modules $\cal{M}\to \cal{A}$
(posant pour simplifier $\cal{M}=g^*(\it{\Omega}^1_{X/S})$), et on
s'int\'eresse \`a ceux qui induisent l'homomorphisme {\emph{nul}}
$\cal{M}_0\to \cal{O}_{Y_0}$ (o\`u
$\cal{M}_0=\cal{M}\otimes_{\cal{O}_Y}\cal{O}_{Y_0}$) \ie qui
appliquent $\cal{M}$ dans l'id\'eal $\cal{J}$ de l'augmentation. On
trouve donc l'ensemble
$\Hom_{\cal{O}_Y}(\cal{M},\cal{J})=\Hom_{\cal{O}_{Y_0}}(\cal{M}_0,\cal{J})$
(puisque $\cal{J}$ est annul\'e par $\cal{J}$). C'est la bijection
canonique cherch\'ee \eqref{eq:III.5.1.*}.

Tenant compte de l'implication (i)$\To$(iii)
dans~\Ref{III.3.1}, on trouve:

\begin{corollaire}
\label{III.5.2}
\marginpar{73}
Si $X$ est lisse sur~$S$ (du moins aux points de~$g_0(Y_0)$)
alors~$\cal{P}$ est m\^eme un \emph{faisceau principal homog\`ene}
sous le faisceau en groupes commutatifs~$\cal{G}$, qui en l'occurrence
peut aussi s'\'ecrire
\begin{equation*}
\cal{G}=g_0^*(\goth{g}_{X/S})\otimes_{\cal{O}_{Y_0}}\cal{J}
\end{equation*}
o\`u $\goth{g}_{X/S}$
\label{indnot:cb}\oldindexnot{$\goth{g}_{X/S}$|hyperpage}%
est le faisceau sur~$X$ dual de~$\it{\Omega}_{X/S}^1$, \ie le
\emph{faisceau tangent}
\index{faisceau tangent|hyperpage}%
\index{tangent (faisceau)|hyperpage}%
(ou \emph{faisceau des d\'erivations})
\index{derivations (faisceau des)@d\'erivations (faisceau des)|hyperpage}%
\index{faisceau des derivations@faisceau des d\'erivations|hyperpage}%
de~$X$ par rapport \`a~$S$. (Cette derni\`ere formule provient du
fait que~$\it{\Omega}_{X/S}^1$ est alors libre de type fini).
\end{corollaire}

En particulier, \`a ce faisceau principal homog\`ene correspond
une classe de cohomologie dans $\H^1(Y_0,\cal{G})$, dont l'annulation
est n\'ecessaire et suffisante pour l'existence d'un
$S$-morphisme~$g$ prolongeant~$g_0$. Et s'il existe un tel
prolongement, l'ensemble de tous les prolongements possibles est un
espace homog\`ene sous le groupe~$\H^0(Y_0,\cal{G})$.

Dans l'application des m\'ethodes de la g\'eom\'etrie formelle,
la situation est le plus souvent la suivante: on donne deux
$S$-pr\'esch\'emas $X$ et~$Y$, un Id\'eal coh\'erent~$\cal{I}$
sur~$S$, on d\'esigne par~$S_n$ le sous-pr\'esch\'ema ferm\'e
de~$S$ d\'efini par~$\cal{I}^{n+1}$, et on pose
\begin{equation*}
X_n=X\times_S S_n \quoi,\quad Y_n=Y\times_S S_n.
\end{equation*}
On suppose qu'on a un $S_n$-morphisme:
\begin{equation*}
g_n\colon Y_n\to X_n
\end{equation*}
(ou, ce qui revient au m\^eme, un~$S$-morphisme $Y_n\to X$ ou encore
un $S_{n+1}$-morphisme $Y_n\to X_{n+1}$,
\ifthenelse{\boolean{orig}}
{puisque un}
{puisqu'un}
tel morphisme
induit n\'ecessairement $Y_n\to X_n$), et on cherche \`a le
prolonger en un $S_{n+1}$-morphisme
\begin{equation*}
g_{n+1}\colon Y_{n+1}\to X_{n+1}
\end{equation*}
(Si on peut continuer ind\'efiniment, on obtient donc un morphisme
$\widehat{Y}\to\widehat{X}$ pour les pr\'esch\'emas formels
compl\'et\'es de~$Y$ et~$X$ pour les Id\'eaux
\ifthenelse{\boolean{orig}}{$\cal{IO}_X$ et~$\cal{IO}_Y$).}
{$\cal{IO}_Y$ et~$\cal{IO}_X$).}
On peut appliquer~\Ref{III.5.1} en y rempla\c cant $(S,X,Y,Y_0,g_0)$
par $(S_{n+1},X_{n+1},Y_{n+1},Y_n,g_n)$, le faisceau $\cal{G}$ devient
ici le faisceau des homomorphismes de Modules
de~$g_n^*(\it{\Omega}_{X_{n+1}/S_{n+1}}^1)$
dans~$\cal{J}=\cal{I}^{n+1}\cal{O}_Y/\cal{I}^{n+2}\cal{O}_Y$.
Comme~$\cal{J}$ est annul\'e par $\cal{IO}_Y$, on peut remplacer
alors~$g_n^*(\it{\Omega}_{X_{n+1}/S_{n+1}}^1)$
\marginpar{74}
par le faisceau qu'il induit sur~$Y_0$,
savoir~$h_0^*(\it{\Omega}_{X/S}^1)$, o\`u~$h_0$ est le compos\'e
$Y_0\to Y_n\to X_{n+1}$, ou encore le compos\'e $Y_0 \to X_0\to
X_{n+1}$, o\`u $g_0\colon Y_0\to X_0$ est induit par $g_n$. Comme
l'image inverse de~$\it{\Omega}_{X_{n+1}/S_{n+1}}^1$
sur~$X_0=X_{n+1}\times_{S_{n+1}}S_0$
est~$\it{\Omega}_{X_{0}/S_{0}}^1$, on voit qu'on a aussi
\begin{equation*}
\cal{G}=\SheafHom_{\cal{O}_{Y_0}} (g_0^*(\it{\Omega}_{X_0/S_0}^1),
\cal{I}^{n+1}\cal{O}_Y/\cal{I}^{n+2}\cal{O}_Y)
\end{equation*}
Donc on obtient le

\begin{corollaire}
\label{III.5.3}
Soient $S,X,Y,\cal{I},g_n$ comme ci-dessus, soit $\cal{P}(g_n)$ le
faisceau sur~$Y$ dont les sections sur un ouvert $U$ sont les
prolongements $g_{n+1}$ de~$g_n$ en un $S_{n+1}$-morphisme $Y_{n+1}\to
X_{n+1}$. Alors
\ifthenelse{\boolean{orig}}{$\cal{P}$}{$\cal{P}(g_{n})$}
est un faisceau formellement principal homog\`ene sous le faisceau
en groupes
\begin{equation*}
\cal{G}=\SheafHom_{\cal{O}_{Y_0}} (g_0^*(\it{\Omega}_{X_0/S_0}^1),
\gr_{\cal{IO}_Y}^{n+1}(\cal{O}_Y))
\end{equation*}
\end{corollaire}

En particulier:

\begin{corollaire}
\label{III.5.4}
Si de plus $X$ est lisse sur~$S$ (du moins aux points de~$g_0(Y_0)$)
alors
\ifthenelse{\boolean{orig}}{$\cal{P}$}{$\cal{P}(g_{n})$}
est m\^eme un faisceau principal homog\`ene.
\ifthenelse{\boolean{orig}}{}
{\enlargethispage{.5cm}}%
En particulier, il d\'efinit une classe d'obstruction dans
$\H^1(Y_0,\cal{G})$, dont l'annulation est n\'ecessaire et
suffisante pour l'existence d'un prolongement global~$g_{n+1}$
de~$g_n$. Et s'il existe un tel prolongement, l'ensemble de tous les
prolongements globaux est un espace principal homog\`ene sous
$\H^0(Y_0,\cal{G})$. Enfin, dans le cas envisag\'e, le
faisceau~$\cal{G}$ peut aussi s'\'ecrire
\begin{equation*}
\cal{G}=g_0^*(\goth{g}_{X_0/S_0})\otimes_{\cal{O}_{Y_0}}
\gr_{\cal{IO}_Y}^{n+1}(\cal{O}_Y)
\end{equation*}
\end{corollaire}

Proc\'edant de proche en proche, on voit donc que si tous les
$\H^1(Y_0,\cal{G}_n)$ son nuls (o\`u
$\cal{G}_n=g_0^*(\goth{g}_{X_0/S_0})\otimes
\gr_{\cal{IO}_Y}^{n}(\cal{O}_Y)$), alors partant avec un $g_k$
quelconque, on peut le prolonger successivement en $g_{k+1},\dots$ En
particulier, si $\cal{I}$ est nilpotent, on pourra trouver un
prolongement $g$ de~$g_k$ \`a $Y$. La condition de nullit\'e des
$\H^1$ est v\'erifi\'ee en particulier si $Y_0$ est affine. On
trouve donc:

\begin{corollaire}
\label{III.5.5}
Dans
\marginpar{75}
l'\'enonc\'e du th\'eor\`eme~\Ref{III.3.1}, on obtient
une condition n\'ecessaire et suffisante \'equivalente aux autres
en supposant que le~$Y^\prime$ qui intervient dans \textup{(ii)} (ou \textup{(ii~bis)})
est affine, et en exigeant l'existence d'un \emph{prolongement global}
$g$ de~$g_0$ \`a tout~$Y^\prime$.
\end{corollaire}

On notera que la d\'emonstration de~\Ref{III.3.1} ne pouvait pas
nous donner ce r\'esultat directement.

Un cas important est celui o\`u $Y$ est \emph{plat} sur~$S$, alors
on a donc
\begin{equation*}
\gr^n(\cal{O}_Y)=\gr^n(\cal{O}_S)\otimes_{\cal{O}_{S_0}}\cal{O}_{Y_0}
\end{equation*}
et lorsque de plus les $\gr^n(\cal{O}_S)$ sont \emph{localement libres
sur}~$S$, on trouve
\begin{equation*}
\cal{G}_n=
\SheafHom_{\cal{O}_{Y_0}}(g_0^*(\it{\Omega}_{X_0/S_0}^1),\cal{O}_{Y_0})
\otimes_{\cal{O}_{S_0}}\gr^n(\cal{O}_S),
\end{equation*}
ou encore, si $\it{\Omega}_{X_0/S_0}^1$ est lui aussi localement libre
(par exemple $X$ lisse sur~$S$)
\begin{equation*}
\cal{G}_n=g_0^*(\goth{g}_{X_0/S_0})
\otimes_{\cal{O}_{S_0}}\gr^{n}(\cal{O}_S)
\end{equation*}
Si par exemple $S$ est affine d'anneau affine~$A$, $\cal{I}$ \'etant
d\'efini par un id\'eal $I$ de~$A$, on trouve
\begin{equation*}
\H^i(Y_0,\cal{G}_n)=\H^i(Y_0,\cal{G}_0) \otimes_A\gr_I^n(A)
\end{equation*}
pour tout $i$ (en effet, la question est locale sur~$S_0$, et on est
ramen\'e au cas o\`u on tensorise par un Module libre).
\emph{Dans ce cas, la nullit\'e de~$\H^1(Y_0,\cal{G}_0)$ implique
que toutes les obstructions aux prolongements successifs de~$g_n$ sont
nulles}. On obtient donc:

\begin{corollaire}
\label{III.5.6}
Soient $(S,X,Y,\cal{I},g_n)$ comme plus haut, supposons de plus $X$
lisse sur~$S$ et $Y$ plat sur~$S$, enfin $S$ affine, et les
$\gr^n(\cal{O}_S)=\cal{I}^n/\cal{I}^{n+1}$ localement libres. Alors
l'obstruction \`a construire~$g_{n+1}$ se trouve dans
$\H^1(Y_0,\cal{G}_0)\otimes_A\gr_I^{n+1}(A)$ (o\`u~$A$ est l'anneau
de~$S$, $I$ l'id\'eal de~$A$ d\'efinissant $\cal{I}$), en posant
\begin{equation*}
\cal{G}_0=g_0^*(\goth{g}_{X_0/S_0})
\end{equation*}
Si $\H^1(Y_0,\cal{G}_0)=0$, alors $g_n$ peut se prolonger en un
$\widehat{S}$-morphisme $\widehat{g}\colon\widehat{Y}\to\widehat{X}$.
\end{corollaire}

Bien entendu, ce r\'esultat resterait valable tel quel, si on
partait, au lieu de~$S$-pr\'esch\'emas ordinaires $X$ et~$Y$, de
$\widehat{S}$-pr\'esch\'emas formels $\widehat{\cal{I}}$-adiques
\marginpar{76}
$\goth{X}$ et~$\goth{Y}$. Il permet de prouver par exemple que
certains sch\'emas formels propres sur un anneau local complet (par
exemple) sont en fait alg\'ebriques. En effet, proc\'edant comme
dans le lemme~\Ref{III.4.2}, on trouve:

\begin{corollaire}
\label{III.5.7} Sous les conditions de~\Ref{III.5.6}, si $g_0$ est
un isomorphisme il en est de m\^eme de~$\widehat{g}$.
\end{corollaire}

\noindent
(N.B. le m\^eme r\'esultat vaut pour les immersions ferm\'ees).

On obtient ainsi:

\begin{proposition}
\label{III.5.8} Soient $A$ un anneau local complet d'id\'eal
maximal $\goth{m}$,
\ifthenelse{\boolean{orig}}{\ignorespaces}{de}
corps r\'esiduel~$k$, soient $\goth{X}$ et $\goth{Y}$ deux
pr\'esch\'emas formels $\goth{m}$-adiques sur~$A$, plats sur~$A$
(\ie pour tout~$n$, $X_n$ et $Y_n$ sont plats
sur~$A_n=A/\goth{m}^{n+1}$), on suppose $X_0=\goth{X}\otimes_Ak$ lisse
sur~$k$ et $\H^1(X_0,\goth{g}_{X_0/k})=0$. Alors tout $k$-isomorphisme
de~$Y_0$ sur~$X_0$ se prolonge en un $A$-isomorphisme de~$\goth{Y}$
sur~$\goth{X}$; ce prolongement est unique si de plus
$\H^0(X_0,\goth{g}_{X_0/k})=0$.
\end{proposition}

Cela donne en particulier un r\'esultat d'\emph{unicit\'e} de
pr\'esch\'ema formel lisse sur~$A$ se r\'eduisant suivant un
pr\'esch\'ema $X_0$ donn\'e (moyennant
$\H^1(X_0,\goth{g}_{X_0/k})=0$). De plus, si~$\goth{X}$ et~$\goth{Y}$
proviennent de sch\'emas ordinaires propres sur~$A$, soient $X$
et~$Y$, alors on sait d'apr\`es le th\'eor\`eme d'existence de
faisceaux en G\'eom\'etrie formelle (\cf expos\'e au
S\'eminaire Bourbaki, \No 182\footnote{\Cf EGA III 5.4.1 pour la
d\'emonstration.}) qu'il y a correspondance biunivoque entre les
$A$-isomorphismes $Y\isomto X$ et le $A$-isomorphismes des
compl\'et\'es formels, donc

\begin{corollaire}
\label{III.5.9} L'\'enonc\'e pr\'ec\'edent~\Ref{III.5.8} reste
valable en y rempla\c cant $\goth{X}$ et $\goth{Y}$ par des
$A$-sch\'emas ordinaires $X$ et~$Y$, \emph{propres} sur~$A$.
\end{corollaire}

Enfin, lorsque $\goth{X}$ est un sch\'ema formel propre sur~$A$, et
que $\goth{Y}$ est de la forme~$\widehat{Y}$ o\`u~$Y$ est un
sch\'ema ordinaire propre sur~$A$, alors la
proposition~\Ref{III.5.8} donne des conditions suffisantes pour qu'on
puisse trouver un isomorphisme de~$\goth{X}$ avec~$\widehat{Y}$, donc
pour que le sch\'ema formel $\goth{X}$ soit en fait
\og alg\'ebrique\fg (\ie isomorphe \`a un~$\widehat{X}$, $X$ un
sch\'ema ordinaire propre sur~$A$, lequel sera alors canoniquement
d\'etermin\'e). C'est ce qui a lieu notamment si $X_0=\PP_k^r$ (ou
plus g\'en\'eralement, si $X_0$ est un sch\'ema de
\ifthenelse{\boolean{orig}}{S\'ev\'eri--Brauer}{Severi-Brauer},
\ie devient
\marginpar{77}
isomorphe \`a l'espace projectif type sur la cl\^oture
alg\'ebrique de~$k$): tout sch\'ema formel propre et plat sur~$A$,
de fibre~$\PP_k^r$, est alg\'ebrisable, et de fa\c con plus
pr\'ecise est isomorphe au compl\'et\'e formel $\goth{m}$-adique
de~$\PP_A^r$. En particulier (gr\^ace au \og th\'eor\`eme
d'existence\fg) tout sch\'ema ordinaire propre sur~$A$, de
fibre~$\PP_k^r$, est isomorphe \`a~$\PP_A^r$ ($A$ \'etant un
anneau local complet). Utilisant la th\'eorie de la descente, on
peut prouver que si~$A$ n'est pas complet, $X$ devient isomorphe
\`a~$\PP^r$ en faisant une extension $A\to A^\prime$ finie et
\'etale de la base (et sous cette forme, le r\'esultat reste
valable pour une fibre qui est un sch\'ema de
\ifthenelse{\boolean{orig}}{S\'ev\'eri--Brauer}{Severi-Brauer}).

\section{Prolongement infinit\'esimal global des $S$-sch\'emas lisses}
\label{III.6}

Sous les conditions du th\'eor\`eme~\Ref{III.4.1}, on se propose
de chercher s'il existe un pr\'esch\'ema~$X$ lisse sur~$Y$ tel que
$X\times_YY_0$ soit $Y_0$-isomorphe \`a~$X_0$, sachant qu'un tel
sch\'ema \og existe localement sur~$X_0$\fg. Reprenant la m\'ethode
de construction de proche en proche, on est conduit \`a
remplacer~$Y$ par la lettre~$S$, \`a supposer qu'on se donne un
sous-pr\'esch\'ema ferm\'e $S_0$ de~$S$ d\'efini par un
faisceau d'id\'eaux~$\cal{I}$, (qu'il n'est plus n\'ecessaire de
supposer localement nilpotent), \`a introduire les
sous-pr\'esch\'emas ferm\'es~$S_n$ de~$S$ d\'efinis par
les~$\cal{I}^{n+1}$, et \`a supposer qu'on s'est donn\'e un
sous-pr\'esch\'ema~$X_{n}$ lisse sur~$S_n$. On se propose de
trouver un $S_{n+1}$-pr\'esch\'ema~$X_{n+1}$ \og qui se r\'eduit
suivant~$X_{n}$\fg, \ie muni d'un isomorphisme
\begin{equation*}
X_{n+1}\times_{S_{n+1}}S_n\isomfrom X_n
\end{equation*}
qui soit \emph{lisse} sur~$S_{n+1}$ (ou, ce qui revient au m\^eme
par II~\Ref{II.2.1}, \emph{plat} sur~$S_{n+1}$). Comme nous l'avons
signal\'e dans \No \Ref{III.4}, une telle donn\'ee revient \`a la
donn\'ee d'un faisceau d'alg\`ebres~$\cal{B}$ sur~$f^{-1}(\cal{O}_{S_{n+1}})$ (o\`u $f$ est l'application continue
sous-jacente au morphisme structural $X_n\to S_n$), muni d'une
augmentation $\cal{B}\to\cal{O}_{X_n}$ compatible avec l'augmentation
$f^{-1}(\cal{O}_{S_{n+1}})\to f^{-1}(\cal{O}_{S_n})$, et satisfaisant
\`a deux conditions (a) et (b) que nous ne r\'e\'ecrirons pas,
nous bornant \`a noter qu'elles sont de \emph{nature locale} sur
l'espace topologique sous-jacent \`a~$X_n$. On sait
d'apr\`es~\Ref{III.4.1} qu'une solution existe localement. Elle est
de plus unique \`a isomorphisme (non unique) pr\`es, du moins
localement. Commen\c cons par pr\'eciser ce point:

\begin{proposition}
\label{III.6.1}
Soit
\marginpar{78}
$X_{n+1}$ sur~$S_{n+1}$ se r\'eduisant suivant $X_n$
sur~$S_{n}$. Alors le faisceau (sur l'espace topologique sous-jacent
\`a~$X_{n}$, ou encore \`a~$X_0$) des $S_{n+1}$-automorphismes
de~$X_{n+1}$ qui induisent l'identit\'e sur~$X_n$ est canoniquement
isomorphe~\`a
\begin{equation*}
\cal{G}=\goth{g}_{X_0/S_0}\otimes_{\cal{O}_{S_0}}
\gr_{\cal{I}}^{n+1}(\cal{O}_S)
\end{equation*}
(en tant que faisceau de groupes).
\end{proposition}

En effet, par~\Ref{III.5.4} et~\Ref{III.4.2} ce faisceau est un
faisceau principal homog\`ene sous~$\cal{G}$. Comme il est muni
d'une section privil\'egi\'ee (l'automorphisme identique
de~$X_{n+1}$), il s'identifie donc comme faisceau d'ensembles
\`a~$\cal{G}$. Il faut v\'erifier que cette identification est
compatible avec les structures de groupe. C'est facile, et d'ailleurs
un cas particulier d'un r\'esultat plus g\'en\'eral sur la
compatibilit\'e des structures de fibr\'es principaux,
dans~\Ref{III.5.1} et~\Ref{III.5.3}, avec la composition des
morphismes (r\'esultat que nous n'\'enon\c cons pas ici, mais
qui se doit de figurer dans le hyperplodoque).

En particulier, le faisceau sur~$X_0$ des germes d'automorphismes
de~$X_{n+1}$ (avec les structures explicit\'ees) est
\emph{commutatif}. Il s'ensuit que si~$X_{n+1}^\prime$ est une autre
solution du probl\`eme, isomorphe \`a~$X_{n+1}$ au-dessus de
l'ouvert~$U$ de~$X_0$, alors l'isomorphisme de~$\SheafAut(X_{n+1})|U$
sur~$\SheafAut(X_{n+1}^\prime)|U$ d\'eduit par transport de
structure d'un isomorphisme $X_{n+1}|U\isomto X_{n+1}^\prime|U$,
\emph{ne d\'epend pas} du choix de ce dernier. (Ce n'est d'ailleurs
autre que l'isomorphisme identique de~$\cal{G}$, lorsque on identifie
l'un et l'autre faisceau d'automorphismes \`a~$\cal{G}$ gr\^ace
\`a~\Ref{III.6.1}).

On d\'eduit de~\Ref{III.6.1}:

\begin{corollaire}
\label{III.6.2}
Soient~$X_{n+1}$, $X_{n+1}'$ lisses sur~$S_{n+1}$ et \og se
r\'eduisant suivant~$X_n$\fg. Alors le faisceau (sur l'espace
sous-jacent \`a~$X_0$) des~$S_{n+1}$-isomorphismes de~$X_{n+1}$
sur~$X_{n+1}'$ induisant l'identit\'e sur~$X_n$, est de fa\c con
naturelle un faisceau principal homog\`ene sous~$\cal{G}$.
\end{corollaire}

Cela exprime en effet que~$X_{n+1}$ et~$X_{n+1}'$ sont isomorphes
localement, et que le faisceau des germes d'automorphismes du premier
est~$\cal{G}$.

Notons maintenant qu'en vertu de~\Ref{III.4.1}, on peut toujours
trouver un recouvrement~$(U_i)$
\marginpar{79}
de~$X_n$ par des ouverts (qu'on peut supposer affines), et pour
tout~$i$ un sch\'ema lisse~$X^i$ sur~$S_{n+1}$ se r\'eduisant
\ifthenelse{\boolean{orig}}{suivant~$U_i=U_{in}$.}{suivant~$U_i$.}
Supposons pour simplifier~$X_n$ \emph{s\'epar\'e}, donc
les~$U_{ij}=U_i\cap U_j$ sont encore des ouverts \emph{affines}
de~$X_n$. Comme le~$\H^1$ d'un tel ouvert, \`a valeurs dans le
faisceau quasi-coh\'erent~$\cal{G}$, est nul, on en d\'eduit par
\ifthenelse{\boolean{orig}}{\ignorespaces}{le}
corollaire~\Ref{III.6.2} que~$X^i|U_{ij}$ est isomorphe
\`a~$X^j|U_{ij}$; soit
$$
f_{ji}:X^i|U_{ij}\isomto X^j|U_{ij}
$$
un tel isomorphisme. Il est d\'etermin\'e \`a une section
pr\`es de~$\cal{G}$ sur~$U_{ij}$. Posons, pour tout
\ifthenelse{\boolean{orig}}{triple}{triplet}
d'indices:
$$
f_{ji}^{(k)}=f_{ji}|U_{ijk}\quad\text{o\`u }U_{ijk}=U_i\cap U_j\cap
U_k.
$$
Si on avait
\begin{equation*}
\label{eq:III.6.1}
\tag{1} {f_{kj}^{(i)}f_{ji}^{(k)}=f_{ki}^{(j)},}
\end{equation*}
il s'ensuivrait que les~$X^i$ se \og recollent\fg par les~$f_{ji}$, donc
qu'ils d\'efinissent une solution~$X=X_{n+1}$ du probl\`eme
cherch\'e. Une telle solution existe plus g\'en\'eralement, si
on peut modifier les~$f_{ji}$ en des~$f_{ji}'$:
\begin{equation*}
\label{eq:III.6.2}
\tag{2} {f_{ji}'=f_{ji}g_{ji}\quad (g_{ji}\in\Gamma(U_{ij},\cal{G}))}
\end{equation*}
de telle fa\c con que les~$f_{ji}'$ satisfassent la condition de
transitivit\'e ci-dessus. Cette condition suffisante pour
l'existence d'une solution est aussi n\'ecessaire, comme on voit en
se rappelant qu'une telle solution~$X$ doit, sur chaque~$U_i$,
\^etre isomorphe \`a~$X^i$, ce qui permet donc de choisir des
isomorphismes
$$
f_i:X|U_i\isomto X^i
$$
et de d\'efinir des
$$
f_{ji}'=(f_j|U_{ij})(f_i|U_{ij})^{-1}:X^i|U_{ij}\isomto X^j|U_{ij}
$$
satisfaisant la condition de recollement.

Or posons
\begin{equation*}
\label{eq:III.6.3}
\tag{3}
\ifthenelse{\boolean{orig}}
{{f_{ijk}=(f_{ki}^{(j)})^{-1}f_{kj}^if_{ji}^k,}}
{{f_{ijk}=(f_{ki}^{(j)})^{-1}f_{kj}^{(i)}f_{ji}^{(k)},}}
\end{equation*}
c'est un automorphisme de~$X^i|U_{ijk}$, que nous identifierons \`a
une section de~$\cal{G}$ gr\^ace \`a~\Ref{III.6.1}. On constate
que c'est un $2$-\emph{cocycle} $f$ du recouvrement ouvert
\marginpar{80}
$\cal{U}=(U_i)$, \`a coefficients dans~$\cal{G}$, par un petit
calcul formel laiss\'e au lecteur. Le m\^eme calcul montre que
moyennant~\eqref{eq:III.6.2}, la condition de
recollement~\eqref{eq:III.6.1} \emph{pour les} $f_{ij}'$ \'equivaut
\`a la formule
\begin{equation*}
\label{eq:III.6.4}
\tag{4} {f=dg,}
\end{equation*}
o\`u $g=(g_{ij})$ est consid\'er\'e comme une $1$-cocha\^ine
de~$\cal{U}$ \`a coefficients dans~$\cal{G}$. Donc \emph{la
condition n\'ecessaire et suffisante pour l'existence d'une solution
du probl\`eme est que la classe de cohomologie dans
$\H^2(\cal{U},\cal{G})$ d\'efinie par le cocycle~\eqref{eq:III.6.3}}
soit nulle. Notons d'ailleurs que puisque~$\cal{U}=(U_i)$ est un
recouvrement affine de~$X_0$ qui est un \emph{sch\'ema},
$\H^2(\cal{U},\cal{G})$ s'identifie \`a $\H^2(X_0,\cal{G})$. Il est
imm\'ediat d'ailleurs que la classe de cohomologie ainsi obtenue
dans $\H^2(X_0,\cal{G})$ ne d\'epend pas du recouvrement affine
consid\'er\'e. On l'appellera la \emph{classe d'obstruction au
prolongement de~$X_n$ en un sch\'ema~$X_{n+1}$ lisse sur~$S_{n+1}$}.

Supposons cette obstruction nulle. Alors un raisonnement esquiss\'e
plus haut montre que toute solution~$X=X_{n+1}$ est isomorphe \`a
une solution obtenue par recollement \`a partir
d'isomorphismes~$f_{ji}'$, qu'on peut supposer sous la
forme~\eqref{eq:III.6.2}, la condition de recollement n'\'etant
autre que~\eqref{eq:III.6.3}. L'ensemble des $g$ admissibles est donc
un espace principal homog\`ene sous le groupe
$\mathrm{Z}^1(\cal{U},\cal{G})$ des $1$-cocycles de~$\cal{U}$ \`a
coefficients dans~$\cal{G}$. De plus, on constate tout de suite que
\emph{deux cocha\^ines $g$ et $g'$} (telles que $dg=dg'=f$)
\emph{d\'efinissent des solutions isomorphes si et seulement si le
cocycle $g-g'$ est de la forme $dh$}, o\`u
$h=(h_i)\in\mathrm{C}^0(\cal{U},\cal{G})$. On trouve donc:

\begin{theoreme}
\label{III.6.3}
\index{obstruction au prolongement d'un sch\'ema en un sch\'ema lisse|hyperpage}%
Soient $(S,\cal{I},X_n)$ comme ci-dessus, $X_n$ \'etant suppos\'e
s\'epar\'e\footnote{Cette condition est en fait inutile, et on
peut \'eviter les calculs de cocycles plus haut. \Cf le livre de J.
\textsc{Giraud}, Cohomologie Non Ab\'elienne (\`a para\^itre dans
Springer Verlag 1971). Comparer remarques~\Ref{III.6.4}.}. Alors on
peut d\'efinir canoniquement une classe d'obstruction dans
$\H^2(X_0,\cal{G})$ (o\`u~$\cal{G}$ est d\'efini
dans~\Ref{III.6.1}), dont l'annulation est n\'ecessaire et
suffisante pour l'existence d'un sch\'ema~$X_{n+1}$ lisse
sur~$S_{n+1}$ se r\'eduisant suivant~$X_n$. Si cette obstruction
est nulle, alors l'ensemble des classes, \`a isomorphisme pr\`es
(induisant l'identit\'e sur~$X_n$)
de~$S_{n+1}$-pr\'esch\'emas~$X_{n+1}$ se r\'eduisant
suivant~$X_n$, est de fa\c con naturelle un espace principal
homog\`ene sous $\H^1(X_0,\cal{G})$.
\end{theoreme}

\begin{remarques}
\label{III.6.4}
\`A partir de~\Ref{III.6.1}, les raisonnements faits ici sont purement
formels, et se transcrivent avantageusement dans le cadre des
cat\'egories locales,
\marginpar{81}
ou m\^eme des cat\'egories fibr\'ees g\'en\'erales. La
classe d'obstruction \`a l'existence d'un objet \og global\fg d'une
cat\'egorie (dont on peut trouver un objet \og localement\fg, et dont
deux objets sont toujours \og isomorphes localement\fg, le groupe des
automorphismes de tout objet \'etant commutatif) ainsi obtenu dans
le contexte g\'en\'eral, contient comme cas particulier le
\og deuxi\`eme homomorphisme bord\fg dans une suite exacte de faisceaux
de groupes non n\'ecessairement commutatifs, (\'etudi\'e par
exemple par Grothendieck dans Kansas ou
\ifthenelse{\boolean{orig}}
{Tohoku}
{T\^ohoku}).
Le calcul b\'eb\^ete par cocycles fait ici doit donc \^etre regard\'e
comme un pis-aller, d\^u \`a la non existence d'un texte de
r\'ef\'erence satisfaisant.
\end{remarques}

\subsection{}
\label{III.6.5}
On notera que dans~\Ref{III.6.3}, il n'y a pas en g\'en\'eral
d'\'el\'ement privil\'egi\'e dans l'espace principal
homog\`ene envisag\'e sous $\H^1(X_0,\cal{G})$. Cela se traduit
notamment par le fait que l'on obtient (localisant sur~$S$) un
faisceau principal homog\`ene sur~$S_0$, de groupe structural
$\R^1f_*(\cal{G})$, qui n'est pas n\'ecessairement trivial, \ie qui
d\'efinit une classe de cohomologie dans
$\H^1(S_0,\R^1f_*(\cal{G}))$ qui n'est pas n\'ecessairement nulle.
(Lorsque l'on suppose que la classe $d\in\H^2(X_0,\cal{G})$ n'est pas
nulle, mais nulle \og localement au-dessus de~$S$\fg, \ie d\'efinit
une section nulle de~$\R^2f_*(\cal{G})$, \ie un \'el\'ement nul
dans $\H^0(S_0,\R^2f_*(\cal{G}))$).

\subsection{}
\label{III.6.6}
On ne sait pour l'instant \`a peu pr\`es rien sur le m\'ecanisme
alg\'ebrique g\'en\'eral des classes de cohomologie introduites
dans ce num\'ero et ses relations avec le num\'ero
pr\'ec\'edent, et on ne sait rien en dire de pr\'ecis dans les
cas particuliers les plus simples, tel le cas des sch\'emas
ab\'eliens sur des anneaux artiniens\footnote{\label{III.6.6.p24}%
\ifthenelse{\boolean{orig}}
{On sait maintenant que cette obstruction est toujours nulle dans ce
cas.}
{On sait maintenant que cette obstruction est toujours nulle dans ce
cas \lcrochetbf ajout\'e en 2003 (MR): \cf F\ptbl Oort, Finite group schemes, local moduli for abelian varieties and
lifting problems, Algebraic Geometry Oslo 1970, Wolters-Noordhoff, 1972, p\ptbl223--254\rcrochetbf.}}. On esp\`ere qu'il se trouvera des gens pour chiader la
question, qui semble particuli\`erement int\'eressante. Elle est
intimement li\'ee en particulier \`a la \og th\'eorie des
modules\fg des structures alg\'ebriques.

\begin{corollaire}
\label{III.6.7}
Supposons que $\H^2(X_0,\cal{G})=0$, alors un~$X_{n+1}$ existe, et il
est unique \`a isomorphisme pr\`es si de plus
$\H^1(X_0,\cal{G})=0$.
\end{corollaire}

En particulier, on en conclut, en proc\'edant de proche en proche
(et remarquant qu'un sch\'ema affine est acyclique pour un faisceau
quasi-coh\'erent):

\begin{corollaire}
\label{III.6.8}
Sous
\marginpar{82}
les conditions du th\'eor\`eme~\Ref{III.4.1}, si~$X_0$ est
affine, alors il existe un~$X$ lisse sur~$Y$ se r\'eduisant
suivant~$X_0$, et cet~$X$ est unique \`a isomorphisme (non unique)
pr\`es.
\end{corollaire}

On notera que la d\'emonstration directe du th\ptbl \Ref{III.4.1} ne
pouvait nous donner ce r\'esultat.

\begin{corollaire}
\label{III.6.9}
Sous les conditions de~\Ref{III.6.3}, supposons~$S$ affine
d'anneau~$A$, $\cal{I}$ d\'efini par un id\'eal~$I$ de~$A$, enfin
les $\gr^n_{\cal{I}}(\cal{O}_S)=\cal{I}^n/\cal{I}^{n+1}$ localement
libres. Alors $\H^i(X_0,\cal{G})$ s'identifie \`a
$\H^i(X_0,\cal{G}_0)\otimes_A\gr^{n+1}_I(A)$, o\`u
$$
\cal{G}_0=\goth{g}_{X_0/S_0},
$$
donc la classe d'obstruction au prolongement de~$X_n$ se trouve dans
$\H^2(X_0,\cal{G}_0)\otimes_A\gr^{n+1}_I(A)$, et si elle est nulle,
l'ensemble des classes (\`a isomorphisme pr\`es) de solutions est
un espace principal homog\`ene sous
$\H^1(X_0,\cal{G}_0)\otimes_A\gr^{n+1}_I(A)$.
\end{corollaire}

En particulier:

\begin{corollaire}
\label{III.6.10}
Sous les conditions de~\Ref{III.6.9} supposons
$$
\H^2(X_0,\goth{g}_{X_0/S_0})=0,
$$
alors il existe un sch\'ema formel $\hat{I}$-adique
\ifthenelse{\boolean{orig}}{\ignorespaces}{$\goth{X}$}
sur le compl\'et\'e formel $\cal{I}$-adique~$\hat{S}$ de~$S$, qui
soit \og lisse sur~$S$\fg (\ie les~$\goth{X}_p$ sont lisses sur
les~$S_p$) et qui se r\'eduise suivant~$X_n$, \ie muni d'un
isomorphisme
$$
\goth{X}\times_S S_n\isomfrom X_n.
$$
Si de plus $\H^1(X_0,\goth{g}_{X_0/S_0})=0$, alors un tel~$\goth{X}$
est unique \`a isomorphisme pr\`es.
\end{corollaire}

En effet, on construit de proche en proche $X_{n+1}$,~$X_{n+2}$, etc.,
d'o\`u~$\goth{X}$ en passant \`a la limite inductive des~$X_i$.
L'assertion d'unicit\'e figure d\'ej\`a au \no
pr\'ec\'edent.

\section[Application \`a la construction de sch\'emas formels...]{Application \`a la construction de sch\'emas formels
et de sch\'emas ordinaires lisses sur un anneau local complet~$A$}
\label{III.7}
Les r\'esultats du \no pr\'ec\'edent permettent parfois de
prouver l'existence
\marginpar{83}
d'un sch\'ema formel $\goth{m}$-adique sur un tel anneau, se
r\'eduisant suivant un sch\'ema lisse~$X_0$ sur~$k$
donn\'e. Distinguons deux cas:

a) \emph{$A$ est \og d'\'egales caract\'eristiques\fg} (c'est le cas
en particulier si~$k$ est de caract\'eristique~$0$). Alors on sait
qu'il existe un \emph{sous-corps de repr\'esentants de}~$A$, \ie
un sous-corps~$k'$ tel que~$A\to k$ induise un isomorphisme~$k'\isomto
k$. \emph{Alors il existe m\^eme un sch\'ema ordinaire lisse
sur~$A$ se r\'eduisant suivant~$X_0$}, savoir~$X=X_0\otimes_kA$, $A$~\'etant consid\'er\'e comme une alg\`ebre sur~$k$ gr\^ace
\`a l'homomorphisme~$k\to k'\to A$ d\'efini par~$k'$. Il faut
cependant noter que cette construction n'est pas \og naturelle\fg; il est
facile de se convaincre (d\'ej\`a dans le cas
o\`u~$A=k[t]/(t^2)$, alg\`ebre des nombres duaux) qu'un autre
homomorphisme de rel\`evement~$k\to A$ (d\'efini en l'occurrence
par une d\'erivation absolue de~$k$ dans lui-m\^eme) d\'efinit
un~$X'$ sur~$A$ qui en g\'en\'eral \emph{n'est pas isomorphe
\`a~$X$} (si~$\H^1(X_0,\goth{g}_{X_0/k})\neq0$). Il serait
d'ailleurs int\'eressant d'\'etudier (pour~$k$ de
caract\'eristique~$0$, ou non parfait et de
\ifthenelse{\boolean{orig}}
{car.}
{caract\'eristique}~$p>0$) quels sont
les~$X$ lisses sur~$A$ qu'on obtient ainsi, \`a quelle condition
deux homomorphismes~$k\to A$ d\'efinissent des $A$-sch\'emas
isomorphes. N\'eanmoins, l'existence de~$k'$ suffit \`a
entra\^iner que la premi\`ere obstruction au rel\`evement
de~$X_0$, qui est
dans~$\H^2(X_0,\goth{g}_{X_0/k})\otimes_k\goth{m}/\goth{m}^2$, est
n\'ecessairement nulle. Bien entendu, quand on a relev\'e
alors~$X_0$ en~$X_1$ lisse sur~$A/\goth{m}^2$, la nouvelle obstruction
\`a la construction de~$X_2$ ne sera en g\'en\'eral pas nulle:
elle sera fonction d'un \'el\'ement variable dans un certain
espace principal homog\`ene
sous~$\H^1(X_0,\cal{G}_0)\otimes\goth{m}/\goth{m}^2$ et se trouve
dans~$\H^2(X_0,\cal{G}_0)\otimes\goth{m}^2/\goth{m}^3$: il
conviendrait d'\'etudier la situation de fa\c con
d\'etaill\'ee\footnote{Elle est sans doute d\'ecrite par
l'op\'eration crochet de Kodaira-Spencer (\cf S\'eminaire Cartan,
1960/61, Exp\ptbl 4).}.

b) \emph{$A$ est d'in\'egales caract\'eristiques.} Dans ce cas, on
ignore tout, (sauf si par chance $\H^2(X_0,\goth{g}_{X_0/k})=0$,
auquel cas on peut construire un sch\'ema formel $\goth{m}$-adique
simple sur~$A$ se r\'eduisant suivant~$k$). M\^eme
si~$A=\ZZ/p^2\ZZ$ et si~$X_0$ est un sch\'ema \og ab\'elien\fg de
dimension~$2$, on ne sait pas si on peut le relever en un~$X=X_1$ lisse
sur~$A$\kern1pt\footnote{C'est maintenant prouv\'e, \cf \
note~\Ref{III.6.6.p24} page~\pageref{III.6.6.p24}.}, d'autre part, on
n'a pas d'exemple d'un~$X_0$ dont on ait prouv\'e qu'il ne provient
pas d'un sch\'ema ordinaire~$X$ lisse sur~$A$. (J'ai l'impression
que cela doit exister, avec~$X_0$ une surface
projective)\footnote{\label{nIII.3}Un tel exemple a \'et\'e depuis
construit par
\ifthenelse{\boolean{orig}}
{J.P\ptbl Serre (Proc.\ Nat.\ Acad.\ Sc.\ USA, vol\ptbl 47,
\No 1, pp~108--109, 1961)}
{J-P\ptbl Serre (Proc.\ Nat.\ Acad.\ Sc.\ USA  \textbf{47} (1961),
\No 1, p\ptbl108--109)}, du moins pour certaines
dimensions. D\ptbl Mumford a donn\'e un exemple (non publi\'e) avec
une \emph{surface} alg\'ebrique.}. Signalons simplement que
d'apr\`es le th\'eor\`eme de Cohen, il existe un $p$-anneau de
Cohen~$B$ de corps r\'esiduel~$k$ et un homomorphisme~$B\to A$
induisant l'isomorphisme identique sur les corps r\'esiduels; par
suite, le r\'esultat \og le plus fort\fg de rel\`evement serait
obtenu en prenant pour~$A$
\marginpar{84}
un $p$-anneau de Cohen: s'il existe une solution (ordinaire ou
formelle) au-dessus d'un tel anneau, il en existe une au-dessus de
tout anneau local complet de corps r\'esiduel~$k$. En particulier,
comme pour un $p$-anneau de Cohen~$\goth{m}/\goth{m}^2$ s'identifie
canoniquement \`a~$k$, on voit que \emph{pour tout sch\'ema
lisse~$X_0$ sur un corps~$k$ de caract\'eristique~$p>0$, il existe
une classe de cohomologie dans~$\H^2(X_0,\goth{g}_{X_0/k})$}
premi\`ere obstruction au rel\`evement de~$X_0$ en un sch\'ema
lisse sur un $p$-anneau de Cohen; on ignore si elle peut \^etre non
nulle\footnote{Elle peut \^etre non nulle, comme signal\'e en
note~\Ref{nIII.3} \`a la page~\pageref{nIII.3}.}.

M\^eme si on arrive de proche en proche \`a construire les~$X_n$
se r\'eduisant suivant~$X_0$, cela ne donne en g\'en\'eral qu'un
sch\'ema \emph{formel}~$\goth{X}$ lisse sur~$A$, se r\'eduisant
suivant~$X_0$. Lorsque~$X_0$ est propre sur~$A$, il reste la question
si~$\goth{X}$ est en fait alg\'ebrisable, pour pouvoir obtenir un
sch\'ema \emph{ordinaire} propre sur~$A$ et simple sur~$A$, se
r\'eduisant suivant~$X_0$. Le seul crit\`ere connu (signal\'e
dans le S\'eminaire Bourbaki, et qui figure dans les \'El\'ements,
Chap\ptbl III, 4.7.1) est le suivant: si~$\goth{X}$ est propre sur~$A$, et
si~$\cal{L}$ est un faisceau inversible sur~$\goth{X}$ tel que le
faisceau induit~$\cal{L}_0$ sur~$X_0$ soit ample (\ie une puissance
tensorielle convenable~$\cal{L}_0^{\otimes n}$, $n>0$, provient d'une
immersion projective de~$X_0$) alors il existe un sch\'ema~$X$
projectif sur~$A$, et un faisceau inversible ample sur~$X$, tels
que~$(\goth{X},\cal{L})$ s'en d\'eduise par compl\'etion
$\goth{m}$-adique. Cela nous am\`ene donc, \'etant donn\'e un
faisceau localement libre~$\cal{E}_0$ sur~$X_0$ (que nous choisirons
inversible ample pour notre propos), de le prolonger en un faisceau
localement libre~$\cal{E}$ sur~$\goth{X}$. Pour ceci, on est
ramen\'e \`a construire de proche en proche des faisceaux
localement libres~$\cal{E}_n$ sur les~$X_n$. La discussion est toute
analogue \`a celle du \No \Ref{III.6}, (\cf remarque~\Ref{III.6.4}), le r\^ole
essentiel \'etant jou\'e par le \emph{faisceau des automorphismes}
d'un~$\cal{E}_{n+1}$ qui induisent l'identit\'e sur~$\cal{E}_n$: on
montre aussit\^ot que ce faisceau s'identifie \`a
\begin{equation*}
\label{eq:III.7.*}\tag{$*$}{}
\cal{G}= \SheafHom_{\cal{O}_{X_0}}(\cal{E}_0,
\cal{E}_0\otimes\gr_{\cal{I}}^{n+1}(\cal{O}_X))=
\SheafHom_{\cal{O}_{X_0}}(\cal{E}_0,
\cal{E}_0)\otimes\gr_{\cal{I}}^{n+1}(\cal{O}_X)
\end{equation*}
qui est encore un faisceau en groupes commutatifs. On trouve:

\begin{proposition}\label{III.7.1}
Soit~$S$ un pr\'esch\'ema muni d'un faisceau quasi-coh\'erent
d'id\'eaux~$\cal{I}$, $X$ un pr\'esch\'ema sur~$S$, $S_n$ le
sous-pr\'esch\'ema de~$S$ d\'efini par~$\cal{I}^{n+1}$,
et~$X_n=X\times_SS_n$ (pour tout entier~$n$). Soit~$\cal{E}_n$ un
faisceau localement libre
\marginpar{85}
sur~$X_n$, on se propose de le prolonger en un faisceau localement
libre~$\cal{E}_{n+1}$ sur~$X_{n+1}$. Alors~$\cal{E}_n$ d\'efinit une
classe d'obstruction canonique dans~$\H^2(X_0,\cal{G})$,
o\`u~$\cal{G}$ est le faisceau quasi-coh\'erent donn\'e par la
formule ci-dessus, classe dont l'annulation est n\'ecessaire et
suffisante pour l'existence d'un~$\cal{E}_{n+1}$
prolongeant~$\cal{E}_n$. Si cette classe est nulle, alors l'ensemble
des classes, \`a isomorphisme pr\`es (induisant l'identit\'e
sur~$\cal{E}_n$) de solutions~$\cal{E}_{n+1}$, est un espace principal
homog\`ene sous~$\H^1(X_0,\cal{G})$.
\end{proposition}

Cette proposition donne lieu aux corollaires habituels. Signalons
seulement que si~$X$ est \emph{plat} sur~$S$, alors on peut \'ecrire
$$
\cal{G}= \SheafHom_{\cal{O}_{X_0}}(\cal{E}_0, \cal{E}_0)
\otimes_{\cal{O}_{S_0}} \gr_{\cal{I}}^{n+1}(\cal{O}_S)
$$
d'o\`u, si~$S$ est affine d'anneau~$A$ et si
les~$\cal{I}^n/\cal{I}^{n+1}$ sont localement libres, la condition
suffisante
$$
\H^2(X_0,\cal{G}_0)=0 \qquad\text{avec}\qquad
\cal{G}_0=\SheafHom_{\cal{O}_{X_0}}(\cal{E}_0, \cal{E}_0)
$$
pour l'existence d'un~$\cal{E}_{n+1}$, donc de proche en proche pour
l'existence de prolongements successifs~$\cal{E}_m$ ($m=n,n+1,$
\ifthenelse{\boolean{orig}}{etc...}{etc.}).

Revenant \`a la situation de d\'epart, on trouve donc:

\begin{proposition}\label{III.7.2}
Soient~$A$ un anneau local complet, $\goth{X}$ un sch\'ema formel
propre et plat sur~$A$, tel que~$X_0$ soit projectif et
que~$\H^2(X_0,\cal{O}_{X_0})=0$. Alors il existe un sch\'ema~$X$
projectif sur~$A$ dont le compl\'et\'e formel $\goth{m}$-adique
est isomorphe \`a~$\goth{X}$.
\end{proposition}

Conjuguant avec~\Ref{III.6.10}, on trouve:

\begin{theoreme}\label{III.7.3}
Soient~$A$ un anneau local complet de corps r\'esiduel~$k$, $X_0$ un
sch\'ema projectif et lisse sur~$k$, tel que
$$
\H^2(X_0,\goth{g}_{X_0/k})=\H^2(X_0,\cal{O}_{X_0})=0
$$
Alors il existe un sch\'ema lisse et projectif~$X$ sur~$A$, se
r\'eduisant suivant~$X_0$.
\end{theoreme}

Plus g\'en\'eralement, si on se donne un~$X_n$ lisse
sur~$A_n=A/\goth{m}^{n+1}$ se r\'eduisant suivant~$X_0$, alors il
existe un~$X$ lisse et propre sur~$A$ et un
isomorphisme~$X\otimes_AA_n=X_n$.

\begin{corollaire}\label{III.7.4}
Toute courbe lisse et propre sur~$k$ provient par r\'eduction d'une
courbe lisse et propre sur~$A$.
\end{corollaire}

C'est ce r\'esultat qui sera l'outil essentiel (avec le
th\'eor\`eme d'existence de faisceaux en G\'eom\'etrie
Formelle) pour \'etudier le groupe fondamental de~$X_0$ par voie
transcendante.
\marginpar{86}

\chapter{Morphismes plats}\label{IV}
\marginpar{87}
Nous donnons ici surtout les propri\'et\'es de platitude qui nous
ont servi dans les expos\'es pr\'ec\'edents. Une \'etude plus
d\'etaill\'ee se trouvera au Chapitre~IV des \og \'El\'ements de
G\'eom\'etrie Alg\'ebrique\fg en pr\'eparation\footnote{\Cf EGA IV 11 et 12.}, o\`u on \'etudie de fa\c con syst\'ematique
la situation suivante: $X$ \'etant localement de type fini sur~$Y$
localement noeth\'erien, et~$\cal{F}$ coh\'erent sur~$X$ et
$Y$-plat, donner des relations entre les propri\'et\'es de~$Y$,
celles de~$\cal{F}$, et celles des faisceaux coh\'erents induits
par~$\cal{F}$ sur les fibres de~$X\to Y$ (du point de vue notamment de
la dimension, de la dimension cohomologique, de la profondeur
\ifthenelse{\boolean{orig}}{etc...}{etc.}). On a notamment une fa\c con syst\'ematique d'obtenir
des th\'eor\`emes du type \emph{Seidenberg} ou \emph{Bertini}
(pour les sections hyperplanes). Le r\'esultat essentiel pour
l'application des m\'ethodes de platitude dans ce contexte est le
suivant (qui sera d\'emontr\'e plus bas): Si~$Y$ est int\`egre,
$X$ de type fini sur~$Y$, $\cal{F}$ coh\'erent sur~$X$, il existe un
ouvert non vide~$U$ de~$Y$ tel que~$\cal{F}$ soit $Y$-plat aux points
de~$X$ au-dessus de~$U$. Une deuxi\`eme fa\c con, sans doute
encore plus importante, dont la platitude s'introduit en
G\'eom\'etrie Alg\'ebrique, est la \emph{th\'eorie de
descente}: voir par exemple les deux expos\'es de Grothendieck sur
le sujet au S\'eminaire Bourbaki\footnote{et, pour un expos\'e
plus d\'etaill\'e, les Expos\'es~\Ref{VIII} et~\Ref{IX} plus
bas.}. La platitude semble ainsi une des notions techniques centrales
en G\'eom\'etrie Alg\'ebrique.

Rappelons que la notion de platitude et fid\`ele platitude a
\'et\'e introduite par Serre dans GAGA. Un expos\'e des
\no \Ref{IV.1}~et~\Ref{IV.2} suivants se trouvera aussi dans Alg.\ Comm.\ de Bourbaki
(qui bien entendu, comme le titre du livre l'indique, ne se borne pas
au cas d'anneaux de base commutatifs)\footnote{N. Bourbaki,
Alg\`ebre Commutative, Chap\ptbl I (Modules Plats), Act.\ Sc.\ Ind 1290,
Paris, Hermann (1961).}.

Contrairement aux expos\'es pr\'ec\'edents, nous ne supposons
pas que les anneaux envisag\'es sont n\'ecessairement
noeth\'eriens.

\section{Sorites sur les modules plats}
\label{IV.1}

Un module $M$ sur l'anneau~$A$ est dit \emph{plat}
\index{plat, fid\`element plat (module)|hyperpage}%
(ou $A$-plat si on veut pr\'eciser~$A$)
\marginpar{88}
si le foncteur
$$
T_M \colon N \mto M \otimes_{A}N
$$
(qui est en tous cas exact \`a droite) est
\ifthenelse{\boolean{orig}}
{\emph{exact}}
{\emph{exact},}
\ie transforme monomorphismes en monomorphismes. Il revient au
m\^eme de dire que le premier foncteur d\'eriv\'e \`a droite,
ou tous les foncteurs d\'eriv\'es \`a droite, sont nuls,
\ie que l'on a%
\ifthenelse{\boolean{orig}}
{}
{\enlargethispage{.5cm}}
$$
\Tor^{A}_{1}(M,N)=0 \qquad \text{pour tout $N$,}
$$
\resp qu'on a
$$
\Tor^{A}_{i}(M,N)=0 \qquad \text{pour $i>0$, tout $N$.}
$$
Comme les $\Tor_{i}$ commutent aux limites inductives, il suffit
d'ailleurs de v\'erifier ces conditions pour $N$ de type fini, et
m\^eme (prenant alors une suite de composition de~$N$ \`a quotients
monog\`enes) que l'on ait
$$
\Tor^{A}_{1}(M,N)=0
$$
si $N$ monog\`ene, \ie de la forme $A/I$, $I$ \'etant un
id\'eal de~$A$. Notons d'ailleurs que
$$
\Tor^{A}_{1}(M,A/I)=0 \quad\Longleftrightarrow\quad I\otimes_{A}M \rightarrow
M=A\otimes_{A}M\text{ est \emph{injectif}},
$$
comme on voit sur la suite exacte des $\Tor$,
compte tenu de~$\Tor^{A}_{1}(M,A)=0$. Donc $M$ plat \'equivaut \`a
dire que pour tout id\'eal~$I$, l'homomorphisme naturel
$$
I\otimes_{A}M \to IM
$$
est un isomorphisme. Il suffit d'ailleurs de le v\'erifier pour $I$
de type fini, a fortiori il suffit de v\'erifier que le foncteur
$M\otimes$ est exact sur les modules \emph{de type fini}.

Comme chaque fois qu'on a un foncteur exact $T$, si on identifie, pour
un sous-objet $N'$ de~$N$, $T(N')$ \`a un sous-objet de~$T(N)$, on a
$$
T(N' \cap N'')=T(N') \cap T(N'')
$$
$$
T(N' + N'')=T(N') + T(N'')
$$
pour deux sous-objets $N'$, $N''$ de~$N$.

Une somme directe de modules plats, un facteur direct d'un module
plat, est plat. En particulier, $A$ \'etant plat, un module
\emph{libre}, donc aussi un module \emph{projectif}, est plat. Le
produit tensoriel de deux modules plats est plat, et si $M$ est plat
sur~$A$, alors $M\otimes_{A}B$ est plat sur~$B$ pour tout changement
de base
\ifthenelse{\boolean{orig}}{$A \ B$}{$A \to B$}
(\`a cause de l'associativit\'e du produit tensoriel et du fait
qu'un compos\'e de foncteurs
\marginpar{89}
exacts est exact). Si $M$ est plat sur~$B$, $B$ plat sur~$A$, alors
$M$ est plat sur~$A$ (m\^eme raison).

La suite exacte des $\Tor$, plus la \og commutativit\'e\fg du $\Tor$,
donne:

\begin{proposition}
\label{IV.1.1}
Soit $0 \to M' \to M \to M'' \to 0$ une suite exacte de~$A$-modules,
$M''$ \'etant plat. Alors
\begin{enumerate}
\item[(i)] cette suite reste exacte par tensorisation par n'importe
quel $A$-module $N$
\item[(ii)] pour que $M$ soit plat, il faut et il suffit que $M'$ le
soit.
\end{enumerate}
\end{proposition}

On peut donc dire que du point de vue du comportement par produits
tensoriels, les modules plats sont \og aussi bons\fg que les modules
libres ou projectifs (et la suite exacte de~\Ref{IV.1.1} en
particulier est \og aussi bonne\fg que si elle splittait).

Soit $S$ une partie multiplicativement stable de~$A$, alors $S^{-1}A$
est plat sur~$A$, car $S^{-1}A \otimes N=S^{-1}N$ est un foncteur
exact en~$N$. Si $M$ est $A$-plat, alors $S^{-1}M=S^{-1}A \otimes M$
est $S^{-1}A$-plat, la r\'eciproque \'etant vraie si $M \to
S^{-1}M$ est un
\ifthenelse{\boolean{orig}}
{isomorphisme}
{isomorphisme,}
\ie \ifthenelse{\boolean{orig}}{Si}{si}
les $s \in S$ sont bijectifs dans $M$ (\`a cause de la
transitivit\'e de la platitude, $S^{-1}A$ \'etant plat sur~$A$).
Plus g\'en\'eralement, le cas d'un morphisme de pr\'esch\'emas
$X \to Y$ et d'un faisceau quasi-coh\'erent $\cal{F}$ sur~$X$ dont
on veut \'etudier la platitude par rapport \`a~$Y$ conduit \`a
la situation avec deux anneaux:

\begin{proposition}
\label{IV.1.2}
Soient $A \to B$ un homomorphisme d'anneaux, $M$ un $B$-module, $T$
une partie multiplicativement stable de~$B$.
\begin{enumerate}
\item[(i)] Si $M$ est $A$-plat, alors $T^{-1}M$ est $A$-plat (donc
aussi $S^{-1}A$-plat pour toute partie multiplicativement stable $S$
de~$A$ s'envoyant dans $T$).
\item[(ii)] Inversement, si $M_{\goth{n}}$ est plat sur~$A_{\goth{n}}$
pour tout id\'eal maximal $\goth{n}$ de~$B$, $M_{\goth{n}}$ est plat
sur~$A$ (ou, ce qui revient au m\^eme, sur~$A_{\goth{m}}$, o\`u
$\goth{m}$ est l'id\'eal premier de~$A$ image inverse de~$\goth{n}$)
alors $M$ est $A$-plat.
\end{enumerate}
\end{proposition}

On a en effet la formule, fonctorielle par rapport au $A$-module~$N$:
$$
T^{-1}M \otimes_{A} N=T^{-1}(M \otimes_{A} N)
$$
car les deux membres sont fonctoriellement isomorphes \`a $T^{-1}B
\otimes_{B} M \otimes_{B} N_{(B)}$ (avec $N_{(B)}=N \otimes_{A} B$) en
vertu des formules d'associativit\'e de~$\otimes$. Il s'ensuit
aussit\^ot que si $M \otimes_{A} N$ est exact en~$N$, il en est de
m\^eme de~$T^{-1}M \otimes_{A} N$ (comme compos\'e de deux
foncteurs exacts), d'o\`u~(i). Et il s'ensuit de m\^eme (ii), car
pour v\'erifier l'exactitude d'une suite de~$B$-modules, il suffit
de v\'erifier l'exactitude des localis\'es en tous les id\'eaux
maximaux de~$B$.

\begin{proposition}
\label{IV.1.3}
\begin{enumerate}
\item[(i)] Soit
\marginpar{90}
$M$ un $A$-module plat. Si $x \in A$ est non diviseur
de~$0$ dans $A$, il est non-diviseur de~$0$ dans $M$. En particulier,
si $A$ est int\`egre, $M$ est sans torsion.
\item[(ii)] Supposons que $A$ soit int\`egre et que pour tout
id\'eal maximal $\goth{m}$ de~$A$, $A_{\goth{m}}$ soit principal
(par exemple $A$ anneau de Dedekind, ou m\^eme principal). Pour que
le $A$-module $M$ soit plat, il faut et
\ifthenelse{\boolean{orig}}{\ignorespaces}{il}
suffit qu'il soit sans torsion.
\end{enumerate}
\end{proposition}

On obtient (i) en notant que l'homoth\'etie $x$ dans $M$ s'obtient
en tensorisant par~$M$ l'homoth\'etie $x$ dans $A$. Pour (ii), on
peut supposer d\'ej\`a $A$ principal gr\^ace \`a~\Ref{IV.1.2}
(ii); il faut montrer que si $M$ est sans torsion, alors pour tout
id\'eal $I$ de~$A$, l'injection $I \to A$ tensoris\'ee par $M$ est
une injection, ce qui signifie que le g\'en\'erateur $x$ de~$I$
est non diviseur de~$0$ dans~$M$, O.K.

\section{Modules fid\`element plats}
\label{IV.2}

Un foncteur $F$ d'une cat\'egorie dans une autre est dit
\emph{fid\`ele} si pour tout $X$, $Y$, l'application $\Hom (X,Y) \to
\Hom (F(X),F(Y))$ est injective. S'il s'agit d'un foncteur additif de
cat\'egories additives, il revient au m\^eme de dire que $F(u)=0$
implique $u=0$, et cela implique que $F(X)=0$ implique $X=0$. Pour que
$F$ soit \emph{fid\`ele et exact}, il faut et il suffit que la
condition suivante soit v\'erifi\'ee: pour toute suite $M' \to M
\to M''$ de morphismes dans $\cal{C}$, la suite transform\'ee
\ifthenelse{\boolean{orig}}{$F(M)\to F(M') \to F(M'')$}
{$F(M') \to F(M) \to F(M'')$}
est exacte \emph{si} et \emph{seulement si} la
pr\'ec\'edente l'est. Ou encore: $F$ est exact, et $F(X)=0$
implique $X=0$. (N.B. pour pouvoir parler d'exactitude, il faut
supposer les cat\'egories en jeu \emph{ab\'eliennes}). Supposons
qu'on ait une famille $(M_{i})$ d'objets non nuls de~$\cal{C}$ tels
que tout objet non nul de~$\cal{C}$ ait un sous-objet admettant un
quotient isomorphe \`a un $M_{i}$. Alors $F$ fid\`ele et exact
\'equivaut \`a: $F$ exact, et $F(M_{i}) \neq 0$ pour tout~$i$. Si
$\cal{C}$ est la cat\'egorie des modules sur un anneau~$A$, on peut
prendre par exemple pour $(M_{i})$ la famille des $A/\goth{m}$,
$\goth{m}$ parcourant les id\'eaux maximaux de~$A$. (En effet, tout
module non nul admet un sous-module non nul monog\`ene, donc
isomorphe \`a un $A/I$, $I$ id\'eal~$\neq A$, lequel par Krull
admet un quotient $A/\goth{m}$). De ces sorites, on d\'eduit en
particulier:

\begin{proposition}
\label{IV.2.1}
Soit $M$ un $A$-module. Conditions \'equivalentes:
\begin{enumerate}
\item[(i)] Le foncteur $M \otimes_{A}$ est fid\`ele et exact.
\item[(i~bis)]
\marginpar{91}
$M$ est plat, et $M \otimes_{A} N=0$ implique $N=0$
\item[(i~ter)] $M$ est plat, et $M \otimes A/\goth{m} \neq 0$ pour
tout id\'eal maximal $\goth{m}$ de~$A$.
\item[(ii)] Pour toute suite d'homomorphismes $N' \to N \to N''$, la
suite tensoris\'ee par $M$ est exacte si et seulement si la suite
initiale l'est.
\end{enumerate}
\end{proposition}

On dit alors que $M$ est un $A$-module \emph{fid\`element plat}.
\index{fidelement plat (module)@fid\`element plat (module)|hyperpage}%
\index{plat, fid\`element plat (module)|hyperpage}%
En particulier, si $M$ est fid\`element plat, alors $N \to N'$ est
un monomorphisme (\'epimorphisme, isomorphisme) si et seulement si
l'homomorphisme tensoris\'e par~$M$ l'est. Un module fid\`element
plat est \emph{fid\`ele}, puisque l'homoth\'etie $f$ dans~$M$
s'obtient en tensorisant par~$M$ l'homoth\'etie $f$ dans~$A$.

On voit comme
\ifthenelse{\boolean{orig}}{dans~\Ref{IV.1}}{dans~\Ref{IV.1}}
les propri\'et\'es de
transitivit\'e habituelles: le produit tensoriel de deux modules
fid\`element plats est fid\`element plat, si $M$ est
fid\`element plat sur~$A$, $M \otimes_{A}B$ est fid\`element plat
sur~$B$ pour toute extension de la base $A \to B$, si $B$ est une
$A$-alg\`ebre qui est fid\`element plate sur~$A$ et si $M$ est un
$B$-module fid\`element plat, c'est un $A$-module fid\`element
plat.

\begin{corollaire}
\label{IV.2.2}
Soit $A \to B$ un homomorphisme local d'anneaux locaux, $M$ un
$B$-module de type fini. Pour que $M$ soit fid\`element plat
sur~$A$, il faut et il suffit qu'il soit plat sur~$A$ et non nul.
\end{corollaire}

R\'esulte du crit\`ere (i~ter) et de Nakayama. En particulier,
\emph{pour que $B$ soit $A$-plat, il faut et il suffit qu'il soit
fid\`element $A$-plat}.

\begin{proposition}
\label{IV.2.3}
Soit $A \to B$ un homomorphisme d'anneaux, $M$ un $B$-module qui est
fid\`element plat sur~$A$. Pour tout id\'eal premier $\goth{p}$
de~$A$, il existe un id\'eal premier~$\goth{q}$ de~$B$ qui
l'induise.
\end{proposition}

Divisant par $\goth{p}$, on est ramen\'e au cas o\`u
$\goth{p}=0$. Localisant en l'id\'eal premier $0$, on est ramen\'e
au cas o\`u $A$ est un corps. Mais $M$ \'etant fid\`element plat
sur~$A$ est non nul, a fortiori $B \neq 0$, donc $B$ a un id\'eal
premier, qui ne peut qu'induire l'unique id\'eal premier de~$A$!
G\'eom\'etriquement, on peut dire que l'existence d'un faisceau
quasi-coh\'erent $\cal{F}$ sur~$X=\Spec (B)$ qui soit
\og fid\`element plat\fg relativement \`a~$A$, implique que $X \to
Y=\Spec (A)$ est \emph{surjectif}.

\begin{corollaire}
\label{IV.2.4}
Supposons $M$ plat sur~$A$, de type fini sur~$B$ et
\ifthenelse{\boolean{orig}}{supp. $M=\Spec(B)$}{$\Supp M=\Spec(B)$}
(\ie $M_{\goth{q}}\neq 0$ pour tout id\'eal premier $\goth{q}$
de~$B$). Alors les id\'eaux premiers $\goth{q}$ de~$B$ contenant
$\goth{p}B$ minimaux induisent~$\goth{p}$.
\end{corollaire}

On
\marginpar{92}
est encore ramen\'e au cas $\goth{p}=0$ (car les hypoth\`eses se
conservent toutes en divisant), donc $A$ int\`egre. On est
ramen\'e \`a l'\'enonc\'e suivant:

\begin{corollaire}
\label{IV.2.5}
$M$ \'etant comme dessus, tout id\'eal premier minimal $\goth{q}$
de~$B$ induit un id\'eal premier $\goth{p}$ de~$A$ qui est minimal.
\end{corollaire}

En effet, localisant en $\goth{p}$ et $\goth{q}$, on est ramen\'e
\`a prouver que si $A$ et $B$ sont locaux et l'homomorphisme $A \to
B$ local, $M$ un $B$-module non nul plat sur~$A$, et si $B$ est de
dimension $0$, alors $A$ est de dimension $0$. Par~\Ref{IV.2.2}
et~\Ref{IV.2.3}, on conclut que tout id\'eal premier de~$A$ est
induit par un id\'eal premier de~$B$, donc par l'id\'eal maximal
\ifthenelse{\boolean{orig}}{de~$A$,}{de~$B$,}
donc est l'id\'eal maximal, cqfd. G\'eom\'etriquement,
\Ref{IV.2.5} signifie que toute composante irr\'eductible de
\ifthenelse{\boolean{orig}}{$\overline{X}=\Spec(B)$} {$X=\Spec (B)$}
domine quelque composante irr\'eductible de~$Y=\Spec (A)$ (moyennant
l'existence d'un faisceau quasi-coh\'erent de type fini sur~$X$, de
support~$X$, et plat par rapport \`a~$Y$).

On notera qu'on n'a pas eu \`a supposer dans~\Ref{IV.2.4} $M$
fid\`element plat sur~$A$, mais rien ne garantit alors l'existence
d'un id\'eal premier contenant $\goth{p}B$ (donc d'un minimal parmi
de tels).

\begin{proposition}
\label{IV.2.6}
Soit $i \colon A \to B$ un homomorphisme d'anneaux. Conditions
\'equivalentes:
\begin{enumerate}\item[(i)] $B$ est un $A$-module fid\`element plat.
\item[(ii)] $B$ est plat sur~$A$, et $\Spec (B) \to \Spec (A)$ est
surjectif
\item[(ii~bis)] $B$ est plat sur~$A$, et tout id\'eal maximal est
induit par un id\'eal de~$B$.
\item[(iii)] $i$ est injectif et $\Coker i$ est un $A$-module plat.
\item[(iv)] Le foncteur $M_{(B)}=M \otimes_{A} B$ en le $A$-module $M$
est exact, et l'homomorphisme fonctoriel canonique $M \to M_{(B)}$ est
injectif.
\item[(iv~bis)] Pour tout id\'eal $I$ de~$A$, $I \otimes_{A} B \to
IB$ est un isomorphisme, et l'image inverse de~$IB$ dans $A$ est
\'egale \`a~$I$.
\end{enumerate}
\end{proposition}

On a (i)$\To$(ii) par~\Ref{IV.2.3},
(ii)$\To$(ii~bis) est trivial, (ii~bis)$\To$(i)
par le crit\`ere (i~ter) de~\Ref{IV.2.1}. On a
(iii)$\To$(iv) par~\Ref{IV.1.1}, (iv)$\To$(iv~bis)
trivialement (en faisant $M=A/I$ dans la deuxi\`eme
condition~(iv~bis)), et (iv~bis)$\To$(i) en vertu du
crit\`ere de platitude par id\'eaux vu au d\'ebut de~\Ref{IV.1}
et du crit\`ere \Ref{IV.2.1}~(i~ter). Enfin,
(iv)$\To$(iii) par une r\'eciproque facile
de~\Ref{IV.1.1}, et (i)$\To$(iv) car
\marginpar{93}
si $N$ est le noyau de~$M \to M \otimes_{A}B =T(M)$, alors ($T$~\'etant exact) $N \to T(N)$ est nul d'o\`u $T(N)=N \otimes_{A}B
=0$, d'o\`u $N=0$, cqfd.

\section{Relations avec la compl\'etion}
\label{IV.3}
Soient $A$ un anneau noeth\'erien, $I$ un id\'eal dans $A$,
$\widehat{A}$ le s\'epar\'e compl\'et\'e de~$A$ pour la
topologie $I$-pr\'eadique, et pour tout $A$-module $M$, soit
$\widehat{M}$ son compl\'et\'e pour la topologie
$I$-pr\'eadique. C'est un $\widehat{A}$-module, d'o\`u un
homomorphisme canonique
$$
M \otimes_A \widehat{A} \to \widehat{M}.
$$
Lorsque $M$ parcourt les modules \emph{de type fini}, le foncteur $M
\mto \widehat{M}$ est exact, comme il r\'esulte facilement du
\emph{th\'eor\`eme de Krull: Si $N \subset M$, la topologie de~$N$
est celle induite par la topologie de~$M$}. Comme $M \otimes_A
\widehat{A}$ est exact \`a droite, on en conclut ais\'ement (en
r\'esolvant $M$ par $L \to L' \to M$, avec $L$ et $L'$ libres de
type fini) que l'homomorphisme fonctoriel plus haut est un
\emph{isomorphisme} ($\widehat{M}$ \'etant aussi exact \`a droite)
et par cons\'equent que $M \otimes_A \widehat{A}$ est aussi un
foncteur \emph{exact} en $M$. Par suite:
\begin{proposition}
\label{IV.3.1}
Soient $A$ un anneau noeth\'erien, $I$ un id\'eal de~$A$, alors le
compl\'et\'e s\'epar\'e $\widehat{A}$ de~$A$ (pour la
topologie $I$-pr\'eadique) est \emph{plat} sur~$A$.
\end{proposition}

\begin{corollaire}
\label{IV.3.2}
Pour que $\widehat{A}$ soit fid\`element plat sur~$A$, il faut et il
suffit que $I$ soit contenu dans le radical de~$A$.
\end{corollaire}

En effet, il suffit d'appliquer le crit\`ere \Ref{IV.2.1} (i~ter).

Ces r\'esultats r\'esument tout ce qu'on sait dire, du point de
vue de l'alg\`ebre lin\'eaire, sur les relations entre $A$ et
$\widehat{A}$. Le corollaire \Ref{IV.3.2} est surtout utilis\'e lorsque $A$
est un anneau local noeth\'erien et que $I$ est contenu dans
l'id\'eal maximal $\goth{m}$ (et le plus souvent, lui est \'egal).

\section{Relations avec les modules libres}
\label{IV.4}

\begin{proposition}
\label{IV.4.1}
Soient $A$ un anneau, $I$ un id\'eal de~$A$, $M$ un
$A$-module. Supposons qu'on soit sous l'une ou l'autre des
hypoth\`eses suivantes:
\begin{enumerate}
\item[(a)] $I$ est nilpotent
\item[(b)] $A$ est noeth\'erien, $I$ est dans le radical de~$A$, et
$M$ est de type fini.
\end{enumerate}
Pour que $M$ soit libre sur~$A$, il faut et il suffit que $M \otimes
A/I$ soit libre sur~$A/I$, et que $\Tor_1^A(M,A/I) = 0$.
\end{proposition}

C'est
\marginpar{94}
n\'ecessaire, prouvons la suffisance. Soit $(e_i)$ une famille
d'\'el\'ements de~$M$ dont l'image dans $M \otimes A/I = M/IM$ y
d\'efinit une base sur~$A/I$ (c'est une famille finie dans le cas
(b)). Soit $L$ le $A$-module libre construit sur le m\^eme ensemble
d'indices, on a donc un homomorphisme $L \to M$ tel que la
tensorisation $T$ par $A/I$ induit un isomorphisme $T(L)
\isomto T(M)$. Si $Q$ est le conoyau de~$L \to M$, on a
donc $T(Q) = 0$, d'o\`u $Q=0$ en vertu de Nakayama (valable sous
l'une ou l'autre condition (a) ou (b)). Donc $L \to M$ est surjectif,
soit $R$ son noyau, on a donc une suite exacte
$$
0 \to R \to L \to M \to 0
$$
d'o\`u, comme $\Tor_1^A(M,A/I)=0$, une suite exacte $0 \to T(R) \to
T(L) \to T(M) \to 0$, d'o\`u $T(R)=0$, d'o\`u encore
\ifthenelse{\boolean{orig}}{$T(R)=0$}{$R=0$}
en vertu de Nakayama (tenant compte que dans le cas (b), $R$ est de
type fini puisque $A$ \'etait suppos\'e noeth\'erien).

\ifthenelse{\boolean{orig}}
{\begin{corollairestar}}
{\begin{corollaire}\label{IV.4.2}}
On peut remplacer la condition $\Tor_1^A(M,A/I)=0$ par:
l'homomorphisme canonique surjectif
\begin{equation*}
\label{eq:IV.2.*}
\tag{$*$} { \gr_I^0(M) \otimes_{A/I} \gr_I(A) \to \gr_I(M) }
\end{equation*}
est un isomorphisme.
\ifthenelse{\boolean{orig}}
{\end{corollairestar}}
{\end{corollaire}}

En effet, si $M$ est libre, cela est certainement v\'erifi\'e. Il
faut donc prouver que si $M \otimes A/I$ est libre sur~$A/I$ et la
condition sur les $\gr$ v\'erifi\'ee, alors $M$ est libre. On
reprend la d\'emonstration ci-dessus en construisant $L \to M$, il
r\'esulte de l'hypoth\`ese que cet homomorphisme induit un
isomorphisme pour les gradu\'es associ\'es, donc son noyau est
contenu dans l'intersection des $I^nL$, donc est nul (comme il est
trivial dans (a), et bien connu dans (b)). Cqfd.

\ifthenelse{\boolean{orig}}{\refstepcounter{subsection}}{}

\begin{corollaire}
\label{IV.4.3}
Supposons que $A/I$ soit un corps. Alors les conditions suivantes sur~$M$ sont \'equivalentes:
\begin{enumerate}
\item[(i)] $M$ est libre
\item[(ii)] $M$ est projectif
\item[(iii)] $M$ est plat
\item[(iv)] $\Tor_1^A(M, A/I) = 0$
\item[(v)] L'homomorphisme canonique \eqref{eq:IV.2.*} est bijectif.
\end{enumerate}
\end{corollaire}

En effet, dans le cas envisag\'e, $M \otimes A/I$ est
automatiquement libre.

Le r\'esultat pr\'ec\'edent est valable dans les deux cas
suivants:
\begin{enumerate}
\item[(a)] $M$ est un module \emph{quelconque} sur un anneau local $A$
dont l'id\'eal maximal $I$ est \emph{nilpotent} (par exemple un
anneau local artinien).
\item[(b)] $M$ est un module \emph{de type fini} sur un anneau
\emph{local noeth\'erien}.
\end{enumerate}

Rappelons
\marginpar{95}
pour m\'emoire:
\begin{corollaire}
\label{IV.4.4}
Supposons que $A$ soit un anneau \emph{local noeth\'erien
int\`egre} d'id\'eal maximal $\goth{m} = I$, de corps r\'esiduel
$k = A/I$, de corps des fractions $K$. Soit $M$ un module de type fini
sur~$A$. Alors les conditions \'equivalentes \textup{(i)} \`a \textup{(v)}
pr\'ec\'edentes \'equivalent aussi \`a
\ifthenelse{\boolean{orig}}{}
{\enlargethispage{.5cm}}%
\begin{enumerate}
\item[(vi)] $M \otimes_A K$ et $M \otimes_A k$ sont des espaces
vectoriels de m\^eme dimension (\ie le rang de~$M$ sur~$A$ est
\'egal au nombre minimum de g\'en\'erateurs du $A$-module $M$).
\end{enumerate}
\end{corollaire}

D\'emonstration imm\'ediate; on laisse au lecteur le soin de
g\'en\'eraliser au cas o\`u~$A$ est seulement suppos\'e sans
\'el\'ements nilpotents: il faut alors exiger que les rangs de~$M$
pour les id\'eaux premiers minimaux de~$A$ soient \'egaux \`a la
dimension de l'espace vectoriel $M \otimes_A k$.

\section{Crit\`eres locaux de platitude}
\label{IV.5}

\begin{proposition}
\label{IV.5.1}
Soit $A$ un anneau muni d'un id\'eal $I$, $M$ un
$A$-module. Supposons
$$
\Tor_1^A(M, A/I^n) = 0 \qquad \text{pour $n>0$}
$$
alors l'homomorphisme canonique surjectif
\begin{equation*}
\label{eq:IV.5.1.*}
\tag{$*$} {\gr_I^0(M) \otimes_{A/I} \gr_I(A) \to \gr_I(M)}
\end{equation*}
est un isomorphisme. La r\'eciproque est vraie si $I$ est nilpotent.
\end{proposition}

L'hypoth\`ese signifie que les homomorphismes
$$
I^n \otimes_A M \to I^n M
$$
sont des isomorphismes, d'o\`u aussit\^ot le fait que les
homomorphismes
$$
I^n / I^{n+1} \otimes_A M \to I^n M / I^{n+1} M
$$
sont des isomorphismes. R\'eciproquement, supposons cette condition
v\'erifi\'ee et $I$ nilpotent, prouvons $\Tor_1^A(M,A/I^n) =0$
pour tout $n$. C'est vrai pour $n$ grand, proc\'edons par
r\'ecurrence descendante sur~$n$, en le supposant prouv\'e pour
$n+1$. On a un diagramme commutatif
$$
\xymatrix{ & M \otimes I^{n+1} \ar[r] \ar[d] & M \otimes I^n \ar[r]
\ar[d] & M \otimes(I^n/I^{n+1}) \ar[r] \ar[d] & 0 \\ 0 \ar[r] &
MI^{n+1} \ar[r] & MI^n \ar[r] & MI^n/MI^{n+1} \ar[r] & 0 }
$$
o\`u
\marginpar{96}
les lignes sont exactes. Par hypoth\`ese, la derni\`ere fl\`eche
verticale est un isomorphisme, et l'hypoth\`ese de r\'ecurrence
signifie aussi que la premi\`ere fl\`eche verticale l'est. Il en
est donc de m\^eme de la fl\`eche verticale m\'ediane, ce qui
ach\`eve la d\'emonstration.

La proposition suivante a \'et\'e d\'egag\'ee au moment du
S\'eminaire par Serre; elle permet des simplifications
substantielles dans le pr\'esent num\'ero.
\begin{proposition}
\label{IV.5.2}
Soient $A \to B$ un homomorphisme d'anneaux, $M$ un $A$-module. Les
conditions suivantes sont \'equivalentes:
\begin{enumerate}
\item[(i)] Pour tout $B$-module $N$, on a $\Tor_1^A(M,N) = 0$;
\item[(ii)] $\Tor_1^A(M,B) = 0$, et $M_{(B)} = M \otimes_A B$ est
$B$-plat.
\end{enumerate}
\end{proposition}

On a un isomorphisme fonctoriel
$$
M \otimes_A N = (M \otimes_A B) \otimes_B N
$$
qui exprime le premier membre, consid\'er\'e comme foncteur en
$M$, comme un compos\'e de deux foncteurs $M \mto M \otimes_A B$
et $P \mto P \otimes_B N$. Comme le premier transforme modules
libres sur~$A$ en modules libres sur~$B$, donc projectifs en
projectifs, on a la suite spectrale des foncteurs compos\'es
$$
\Tor_n^A(M,N) \From \Tor_p^B(\Tor_q^A(M,B),N)
$$
d'o\`u une suite exacte pour les termes de bas degr\'e
$$
0 \from \Tor_1^B(M \otimes_A B,N) \from \Tor_1^A(M,N)
\from \Tor_1^A(M,B) \otimes_A N
$$
Si (i) est v\'erifi\'e, alors on conclut de cette suite exacte
$\Tor_1^B(M \otimes_A B,N) = 0$ pour tout $N$, \ie $M \otimes_A B$
est $B$-plat, d'o\`u (ii). Si inversement (ii) est v\'erifi\'e,
alors dans la suite exacte les termes entourant $\Tor_1^A(M,N)$ sont
nuls, donc on a (i).

\begin{corollaire}
\label{IV.5.3}
Supposons que $B = A/I$, alors les conditions pr\'ec\'edentes
\'equivalent \`a la suivante:
\begin{enumerate}
\item[(iii)] $\Tor_1^A(M,N) = 0$ pour tout $A$-module $N$ annul\'e
par une puissance de~$I$.
\end{enumerate}
\end{corollaire}

En effet, (i) signifie qu'il en est ainsi si $N$ est annul\'e par
$I$. On en d\'eduit (iii) en appliquant l'hypoth\`ese aux
$I^nN/I^{n+1}N$.

\begin{corollaire}
\label{IV.5.4}
Sous les conditions de~\Ref{IV.5.3}, les conditions envisag\'ees
impliquent que
\marginpar{97}
l'ho\-mo\-mor\-phisme fonctoriel
\begin{equation*}
\label{eq:IV.5.*}
\tag{$*$} {\gr_I^0(M) \otimes_{A/I} \gr_I(A) \to \gr_I(M)}
\end{equation*}
est un isomorphisme, et que $M \otimes_A A/I$ est plat sur~$A/I$.
\end{corollaire}

Il suffit d'appliquer (iii) et \eqref{IV.5.1}. Utilisant la
r\'eciproque de~\Ref{IV.5.1} dans le cas $I$ nilpotent, on trouve:
\begin{corollaire}
\label{IV.5.5}
Soient $A$ un anneau muni d'un id\'eal nilpotent $I$, $M$ un
$A$-module. Les conditions suivantes sont \'equivalentes:
\begin{enumerate}
\item[(i)] $M$ est $A$-plat
\item[(ii)] $M \otimes_A A/I$ est $A/I$-plat, et $\Tor_1^A(M,A/I) = 0$
\item[(iii)] $M \otimes_A A/I$ est $A/I$-plat, et l'homomorphisme
canonique \eqref{eq:IV.5.*} sur les gradu\'es est un isomorphisme.
\end{enumerate}
\end{corollaire}

En effet, ce sont respectivement les conditions (iii) et (ii)
pr\'ec\'edentes, et celles de corollaire~\Ref{IV.5.4}.

Ne supposons plus $I$ nilpotent, alors on aura seulement a priori
dans~\Ref{IV.5.5} les implications (i)$\To$(ii)$\To$(iii). D'autre part, comme la condition (iii) reste
stable en divisant par une puissance de~$I$, on voit en vertu
de~\Ref{IV.5.5} qu'elle implique
\begin{enumerate}
\item[(iv)] \emph{pour tout entier $n$, $M \otimes A/I^n$ est plat sur~$A/I^n$.}
\end{enumerate}
On se propose de donner des conditions moyennant lesquelles on peut
en conclure~(i), \ie que $M$ est $A$-plat. Je dis qu'il suffit pour
ceci que $A$ soit noeth\'erien et que~$M$ satisfasse la condition de
finitude suivante: \emph{pour tout module de type fini $N$ sur~$A$, $M
\otimes_A N$ est s\'epar\'e pour la topologie
$I$-pr\'eadique}. (Il suffirait de le v\'erifier si~$N$ est un
id\'eal de type fini dans $A$). En effet, prouvons que sous ces
conditions, si $N' \to N$ est un monomorphisme de modules de type
fini, $M \otimes_A N' \to M \otimes_A N$ est un monomorphisme. Il
suffit en effet de montrer que le noyau est contenu dans les $I^n(M
\otimes_A N') = \Im(M \otimes_A I^n N' \to M \otimes_A N')$, ou encore
dans les $\Im(M \otimes_A V'_n \to M \otimes_A N') = \Ker (M \otimes_A
N' \to M \otimes_A (N'/V'_n))$, o\`u $V'_n$ parcourt un syst\`eme
fondamental d\'enombrable de voisinages de 0 dans $N'$ (muni de sa
topologie $I$-adique). D'apr\`es le th\'eor\`eme de Krull, la
topologie $I$-adique de~$N'$ est induite par celle de~$N$, on peut
donc prendre $V'_n = N' \cap I^n N$. Consid\'erons alors le
diagramme commutatif
$$
\xymatrix{ M \otimes_A N' \ar[r] \ar[d] & M \otimes_A (N'/V'_n) \ar[d]
\\ M \otimes_A N \ar[r] & M \otimes_A (N / I^n N) }
$$
\marginpar{98}%
Comme $N'/V'_n$ et $N/I^n N$ sont annul\'es par $I^n$, le
deuxi\`eme homomorphisme vertical s'identifie \`a celui d\'eduit
de l'homomorphisme \emph{injectif} $N'/V'_n \to N/I^n N$ en
tensorisant sur~$A/I^n$ avec le $(A/I^n)$-module \emph{plat} $M
\otimes_A A/I^n$, il est donc \emph{injectif}. Par suite, le noyau de
$M \otimes_A N' \to M \otimes_A N$ est contenu dans celui de~$M
\otimes_A N' \to M \otimes_A (N'/V'_n),$ ce qu'on voulait.

La condition de \og finitude\fg envisag\'ee sur~$M$ est
v\'erifi\'ee en particulier si $M$ est un module de type fini sur
une $A$-alg\`ebre
\ifthenelse{\boolean{orig}}
{noeth\'erienne,}
{noeth\'erienne}
$B$ telle que $IB$ soit contenu
dans le radical de~$B$: en effet, alors $M \otimes_A N$ est un module
de type fini sur~$B$ pour tout module de type fini $N$ sur~$A$, donc
s\'epar\'e par Krull pour la topologie $I$-adique $=$ sa topologie
$(IB)$-adique. On trouve ainsi:

\begin{theoreme}
\label{IV.5.6}
Soient $A\to B$ un homomorphisme d'anneaux noeth\'eriens, $I$ un
id\'eal de~$A$ tel que $IB$ soit contenu dans le radical de~$B$, $M$
un $B$-module de type fini. Les conditions suivantes sont
\'equivalentes:
\begin{enumerate}
\item[(i)] $M$ est $A$-plat
\item[(ii)] $M\otimes_A A/I$ est $A/I$-plat, et $\Tor_1^A(M,A/I)=0$
\item[(iii)] $M\otimes_A A/I$ est $A/I$-plat, et l'homomorphisme
canonique
$$
\gr_I^0(M)\otimes_{A/I}\gr_I(A)\to\gr_I(M)
$$
est un isomorphisme.
\item[(iv)] Pour tout entier $n$, $M\otimes_A A/I^n$ est plat
sur~$A/I^n$.
\end{enumerate}
\end{theoreme}

Ce r\'esultat s'applique surtout lorsque $A$, $B$ sont des anneaux
\emph{locaux} noeth\'eriens, $A\to B$ un homomorphisme local, et~$I$
un id\'eal de~$A$ contenu dans son id\'eal maximal (et on peut
r\'eduire aussit\^ot~\Ref{IV.5.6} \`a ce cas). Un cas
int\'eressant est celui o\`u $A/I$ est un corps, \ie $I$ maximal,
auquel cas la condition que $M\otimes_A(A/I)$ est plat sur~$A/I$
devient inutile; de plus, comme alors les $A/I^n$ sont des anneaux
locaux artiniens, la condition (iv) signifie que les
$M\otimes_A(A/I^n)$ sont \emph{libres} sur les~$A/I^n$.

\begin{corollaire}
\label{IV.5.7}
Soit
\marginpar{99}
$A\to B$ un homomorphisme local d'anneaux locaux noeth\'eriens,
$u\colon M'\to M$ un homomorphisme de~$B$-modules de type fini,
supposons $M$ plat sur~$A$. Alors les conditions suivantes sont
\'equivalentes:
\begin{enumerate}
\item[(i)] $u$ est injectif, et $\Coker u$ est plat sur~$A$.
\item[(ii)] $u\otimes_A k\colon M'\otimes_A k\to M\otimes_A k$ est
injectif
\end{enumerate}
(o\`u $k$ d\'esigne le corps r\'esiduel de~$A$).
\end{corollaire}

(i)$\To$(ii) en vertu de~\Ref{IV.1.1}, prouvons la
r\'eciproque. Tout d'abord $u$ est injectif, car il suffit de le
v\'erifier sur les gradues associ\'es, o\`u cela r\'esulte
d'un carr\'e commutatif que le lecteur \'ecrira. Soit $M''$ son
$\Coker$, on a donc une suite exacte
$$
0\to M'\to M\to M''\to 0
$$
d'o\`u par la suite exacte des $\Tor$, compte tenu de
l'hypoth\`ese~(ii) et de~$\Tor_1^A(M,k)=0$, la relation
$\Tor_1^A(M'',k)=0$, donc $M''$ est plat sur~$A$ par le
th\'eor\`eme~\Ref{IV.5.6}.

\begin{corollaire}
\label{IV.5.8}
Sous les conditions de~\Ref{IV.5.6}, soit $J$ un id\'eal de~$B$
contenant $IB$ et contenu dans le radical. Soient $\hat A$ le
compl\'et\'e $I$-adiques de~$A$ et~$\hat B$ et~$\hat M$ les
compl\'et\'es $J$-adiques de~$B$ et~$M$. Pour que $M$ soit
$A$-plat, il faut et il suffit que $M$ soit $\hat A$-plat.
\end{corollaire}

(N.B. la suffisance r\'esulterait d\'ej\`a facilement
de~\Ref{IV.3.2}). On utilise le crit\`ere (iii) de~\Ref{IV.5.6}
dans la situation $(A,B,I,M)$ et dans la situation $(\hat A,\hat
B,I\hat A, \hat M)$. On constate que les conditions obtenues pour
l'un et l'autre cas sont \'equivalentes, gr\^ace
\`a~\Ref{IV.3.2}.

\begin{corollaire}
\label{IV.5.9}
Soient $A\to B\to C$ des homomorphismes locaux d'anneaux locaux
noe\-th\'e\-riens, $M$ un $C$-module de type fini (N.B. $C$
n'intervient que pour pouvoir mettre une condition de finitude sur~$M$). On suppose $B$ plat sur~$A$. Soit $k$ le corps r\'esiduel
de~$A$. Conditions \'equivalentes:
\begin{enumerate}
\item[(i)] $M$ est plat sur~$B$.
\item[(ii)] $M$ est plat sur~$A$, et $M\otimes_A k$ est plat sur~$B\otimes_A k$.
\end{enumerate}
\end{corollaire}

L'implication (i)$\To$(ii) est triviale, prouvons
(ii)$\To$(i). On applique le crit\`ere~(iii)
de~\Ref{IV.5.6} \`a $(B,C,\goth{m}B=I,M)$, comme
$M\otimes_B(B/I)=M\otimes_B(B\otimes_A k) = M\otimes_A k$, la
premi\`ere condition de ce crit\`ere signifie pr\'ecis\'ement
que $M\otimes_A k$ est plat sur~$B\otimes_A k$,
\marginpar{100}
parfait. La deuxi\`eme condition du crit\`ere est
v\'erifi\'ee parce que $M$ est plat sur~$A$ et~$B$ plat sur~$A$,
par une formule d'associativit\'e du produit tensoriel. --- Bien
entendu, se r\'ef\'erant \`a~\Ref{IV.5.5} au lieu
de~\Ref{IV.5.6}, on obtient un \'enonc\'e analogue sans condition
noeth\'erienne et de finitude, quand on suppose en revanche que
l'id\'eal~$\goth{m}$ de~$A$ est nilpotent. (Le fait que $\goth{m}$
ait \'et\'e pris maximal n'est d'ailleurs pas intervenu; mais
c'est en un sens le cas \og $\goth{m}$ maximal\fg qui est \og le meilleur
possible\fg).

\section{Morphismes plats et ensembles ouverts}
\label{IV.6}

Rappelons d'abord quelques r\'esultats sur les ensembles
constructibles, qui sont d'ailleurs d\'emontr\'es dans des notes
en circulation du S\'eminaire Dieudonn\'e-Rosenlicht sur les
Sch\'emas\footnote{\Cf EGA~0$_\mathrm{III}$~9, EGA~IV~1.8 et~1.10.}.

Soit $X$ un espace topologique. On dit avec Chevalley qu'une partie
de~$X$ est \emph{constructible}
\index{constructible (partie)|hyperpage}%
si elle est r\'eunion finie de parties localement ferm\'ees.

\begin{lemme}
\label{IV.6.1}
Soit $X$ un espace topologique noeth\'erien, soit $Z$ une partie
de~$X$. Pour que~$Z$ soit constructible, il faut et il suffit que
pour toute partie ferm\'ee irr\'eductible $Y$ de~$X$, $Z\cap Y$
est non dense dans $Y$ ou contient une partie ouverte non vide de
l'espace~$Y$.
\end{lemme}

On en d\'eduit, utilisant un lemme bien connu d'Alg\`ebre
Commutative:

\begin{lemme}[Chevalley]
\label{IV.6.2}
Soit $f\colon X\to Y$ un morphisme de type fini de
pr\'esch\'emas, avec $Y$ noeth\'erien. Alors $f(X)$ est
constructible.
\end{lemme}

\begin{lemme}
\label{IV.6.3}
Soient $X$ un espace topologique noeth\'erien dont toute partie
ferm\'ee irr\'eductible admet un point g\'en\'erique, $U$ une
partie constructible de~$X$, $x\in X$. Pour que $U$ soit un voisinage
de~$x$, il faut et il suffit que toute g\'en\'erisation $y$ de~$x$
(\ie tout $y\in X$ tel que $x\in \bar{y}$) soit dans~$U$.
\end{lemme}

En particulier

\begin{corollaire}
\label{IV.6.4}
Soit $X$ un espace topologique noeth\'erien dont toute partie
ferm\'ee irr\'eductible admet un point g\'en\'erique, $U$ une
partie de~$X$. Pour que $U$ soit ouverte, il faut et il suffit
qu'elle satisfasse les deux conditions suivantes:
\begin{enumerate}
\item[(a)] $U$ contient toute g\'en\'erisation de chacun de ses
points
\item[(b)] si $x\in U$, alors
\marginpar{101}
$U\cap\bar{x}$ contient une partie ouverte non vide de
l'espace~$\bar{x}$.
\end{enumerate}
\end{corollaire}

En effet, $U$ est n\'ecessairement constructible gr\^ace
\`a~\Ref{IV.6.1}, et on applique le crit\`ere~\Ref{IV.6.2} qui
prouve que $U$ est un voisinage de chacun de ses points.

\begin{corollaire}
\label{IV.6.5}
Soit $f\colon X\to Y$ un morphisme de type fini de pr\'esch\'emas,
avec $Y$ localement noeth\'erien, $x$ un point de~$X$, $y=f(x)$.
Pour que $f$ transforme toute voisinage de~$x$ en un voisinage de~$y$,
il faut et il suffit que pour toute g\'en\'erisation $y'$ de~$y$,
il existe une g\'en\'erisation $x'$ de~$x$ telle que $f(x')=y'$.
\end{corollaire}

On peut \'evidemment supposer que $X$ et~$Y$ sont affines, donc
noeth\'eriens. La condition est suffisante, car il suffit de
prouver que $f(X)$ est un voisinage de~$y$, or $f(X)$ est
constructible par~\Ref{IV.6.1}, et il suffit d'appliquer le
crit\`ere~\Ref{IV.6.3}. La condition est n\'ecessaire, car soit
$Y'=\overline{y'}$, et soit $F$ la r\'eunion des composantes
irr\'eductibles de~$f^{-1}(Y')$ qui ne contiennent pas $x$. Alors
$X-F$ est un voisinage ouvert de~$x$, donc son image est un voisinage
de~$y$, et a fortiori contient $y'$, donc il existe $x'_1\in X-F$ tel
que $f(x'_1)=y'$. Consid\'erons une composante irr\'eductible
de~$f^{-1}(Y')$ contenant $x'_1$, elle contient n\'ecessairement $x$
(car autrement elle serait contenue dans $F$), soit $x'$ son point
g\'en\'erique. C'est une g\'en\'erisation de~$x$, et $f(x')$
est une g\'en\'erisation de~$f(x'_1)=y'$ contenue dans~$Y'$, donc
est \'egal \`a $y'$, cqfd.

\begin{theoreme}
\label{IV.6.6}
Soient $f\colon X\to Y$ un morphisme localement de type fini, avec~$Y$
localement noeth\'erien, $F$ un faisceau coh\'erent sur~$X$ de
support~$X$, plat par rapport \`a~$Y$. Alors $f$ est un morphisme
ouvert (\iev transforme ouverts en ouverts).
\end{theoreme}

Il suffit de prouver le crit\`ere~\Ref{IV.6.5} pour tout point $x\in
X$. Or les g\'en\'erisations $x'$ de~$x$ correspondent aux
id\'eaux premiers de~$\cal{O}_x$, celles~$y'$ de~$y$ correspondent
aux id\'eaux premiers de~$\cal{O}_y$, et il faut donc v\'erifier
que tout id\'eal premier de~$\cal{O}_y$ est induit par une id\'eal
premier de~$\cal{O}_x$. Or $F_x$ est un $\cal{O}_x$-module non nul et
$\cal{O}_y$ plat, donc fid\`element plat sur~$\cal{O}_y$
par~\Ref{IV.2.2}. On peut donc appliquer~\Ref{IV.2.3}, ce qui
ach\`eve la d\'emonstration.

\begin{remarquesstar}
Comme la platitude se conserve par extension de la base, on voit que
sous les conditions de~\Ref{IV.6.5} $f$ est m\^eme
\emph{universellement ouvert}. J'ignore cependant, lorsque $Y$ est
int\`egre et $X$ de type fini sur~$Y$, si $f$ induit sur toute
composante $X_i$ de~$X$ un morphisme ouvert, ou m\^eme seulement
\'equidimensionnel\footnote{La r\'eponse \`a la deuxi\`eme question est
affirmative, celle \`a la premi\`ere n\'egative m\^eme si $f$
est \'etale; \cf EGA~IV~12.1.1.5 et EGA~Err$_{\rm IV}$~33.},
\ie
\marginpar{102}
dont toutes les composantes des fibres ont m\^eme dimension (on sait
seulement que $X_i$ \emph{domine} $Y$). La question est li\'ee
\`a la suivante: soit $A\to B$ un homomorphisme local d'anneaux
locaux noeth\'eriens, tel que $B$ soit plat sur~$A$ et~$\goth{m}B$
soit un id\'eal de d\'efinition de~$B$, (ce qui implique
d'ailleurs $\dim B=\dim A$). Est-il vrai que pour tout id\'eal
premier minimal $\goth{p}_i$ de~$B$, on a $\dim B/\goth{p}_i=\dim B$?
Signalons seulement que la r\'eponse \`a la premi\`ere question
est n\'egative quand on remplace l'hypoth\`ese de platitude
de~\Ref{IV.6.5} par la seule hypoth\`ese que $f$ soit
universellement ouvert.
\end{remarquesstar}

\begin{lemme}
\label{IV.6.7}
Soient $A$ un anneau int\`egre noeth\'erien, $B$ une
$A$-alg\`ebre de type fini, $M$ un $B$-module de type fini. Alors
il existe un \'el\'ement non nul $f$ de~$A$ tel que $M_f$ soit un
module libre (a fortiori plat) sur~$A_f$.
\end{lemme}

Soit $K$ le corps des fractions de~$A$, alors $B\otimes_A K$ est une
alg\`ebre de type fini sur~$K$, et $M\otimes_A K$ un module de type
fini sur cette derni\`ere. Soit $n$ la dimension du support de ce
module, nous raisonnerons par r\'ecurrence sur~$n$.
\ifthenelse{\boolean{orig}}
{Si $n<0$}
{Si $n<0$,}
\ie si $M\otimes_A K=0$, alors prenant un nombre fini de
g\'en\'erateurs de~$M$ sur~$B$, on voit qu'il existe un $f\in A$
qui annule ces g\'en\'erateurs, donc $M$, d'o\`u $M_f=0$ et on a
gagn\'e. Supposons $n\ge 0$. On sait que le $B$-module~$M$ admet
une suite de composition dont les quotients successifs sont isomorphes
\`a des modules $B/\goth{p}_i$, les $\goth{p}_i$ \'etant des
id\'eaux premiers de~$B$. Comme une extension de modules libres est
libre, on est ramen\'e au cas o\`u $M$ lui-m\^eme est de la
forme $B/\goth{p}$, ou encore identique \`a $B$, $B$ \'etant une
$A$-alg\`ebre \emph{int\`egre}. Appliquant le lemme de
normalisation de Noether \`a la $K$-alg\`ebre $B\otimes_A K$, on
voit facilement qu'il existe un \'el\'ement $f$ non nul de~$A$ tel
que $B_f$ soit entier sur le sous-anneau $A_f[t_1,\dots,t_n]$, o\`u
les $t_i$ sont des ind\'etermin\'ees. Donc on peut d\'ej\`a
supposer $B$ est entier sur~$C=A[t_1,\dots,t_n]$, c'est donc un
$C$-module de type fini sans torsion. Soit $m$ son rang, il existe
donc une suite exacte de~$C$-modules:
$$
0\to C^m\to B\to M'\to 0
$$
o\`u $M'$ est un $C$-module de torsion. Il s'ensuit que la
dimension de Krull du~$C\otimes_A K$-module $M'\otimes_A K$ est
strictement inf\'erieure \`a celle~$n$ de~$C\otimes_A K$.
D'apr\`es l'hypoth\`ese de r\'ecurrence, il s'ensuit que, \`a
condition de localiser par rapport \`a un $f$ non nul convenable
de~$A$, on peut supposer que $M'$ est un $A$-module
\marginpar{103}
libre. D'autre part, $C^m$ est un $A$-module libre. Donc $B$ est
alors un $A$-module libre, on a fini.

\begin{lemme}
\label{IV.6.8}
Soient $A$ un anneau noeth\'erien, $B$ une alg\`ebre de type fini
sur~$A$, $M$~un $B$-module de type fini, $\goth{p}$ un id\'eal
premier de~$B$, $\goth{q}$ l'id\'eal premier qu'il induit sur~$A$. On suppose $M_\goth{p}$ plat sur~$A_{\goth{q}}$ (ou sur~$A$,
c'est pareil). Alors il existe un $g\in B-\goth{p}$ tel que
\begin{enumerate}
\item[(a)] $(M/\goth{q}M)_g$ est plat sur~$A/\goth{q}$.
\item[(b)] $\Tor_1^A(M,A/\goth{q})_g=0$.
\end{enumerate}
\end{lemme}
En effet, appliquant~\Ref{IV.6.7} \`a $(A/\goth(q), B/\goth{q}B,
M/\goth{q}M)$ on voit d'abord qu'il existe un~$f$ dans $A-\goth{q}$
tel que $(M/\goth{q}M)_f$ soit plat sur~$A/\goth{q}$. D'autre part,
comme $M_\goth{p}$ est plat sur~$A$, on a
$\Tor_1^A(M,A/\goth{q})_\goth{p}=\Tor_1^A(M_\goth{p},A/\goth{q})=0$,
donc comme $\Tor_1^A(M,A/\goth{q})$ est un $B$-module de type fini, il
existe un $g\in B-\goth{p}$ tel qu'on ait (b). On peut alors
(rempla\c cant $g$ par $gf$) supposer qu'on a en m\^eme temps (a),
ce qui prouve le corollaire.

\begin{corollaire}
\label{IV.6.9}
Avec les notations de~\Ref{IV.6.8}, pour tout id\'eal premier
$\goth{p}'$ de~$B$ contenant~$\goth{p}$ et ne contenant pas~$g$,
$M_{\goth{p}'}$ est plat sur~$A$ (ou, ce qui revient au m\^eme, sur~$A_{\goth{q}'}$, o\`u $\goth{q}'$ est l'id\'eal premier de~$A$
induit par~$\goth{p}')$.
\end{corollaire}
Il suffit d'appliquer le crit\`ere~\Ref{IV.5.6} (ii) au syst\`eme
$(A, B_{\goth{q}'}, \goth{q}, M_{\goth{q'}})$, en utilisant la
localisation des~$\Tor$.
\begin{theoreme}
\label{IV.6.10}
Soit $f\colon X\to Y$ un morphisme de type fini, avec $Y$ localement
noeth\'erien, et soit $F$ un faisceau coh\'erent sur~$X$. Soit $U$
l'ensemble des points $x\in X$ tels que $F_x$ soit plat sur~$\cal{O}_{f(x)}$. Alors $U$ est un ensemble \emph{ouvert}.
\end{theoreme}

\subsubsection*{D\'emonstration} On peut supposer $X$ et $Y$ affines, d'anneaux $B$
et $A$, donc $F$ d\'efini par un $B$-module $M$ de type fini. On
applique le crit\`ere~\Ref{IV.6.4}. La condition (a) est
v\'erifi\'ee trivialement par \Ref{IV.1.2}~(i), reste \`a
v\'erifier la condition (b) de~\Ref{IV.6.4}. C'est ce qui a
\'et\'e fait dans le lemme~\Ref{IV.6.8} et
corollaire~\Ref{IV.6.9}.

Dans beaucoup de questions, la forme plus faible suivante du
th\'eor\`eme~\Ref{IV.6.10} est suffisante (qui r\'esulte
d\'ej\`a du lemme~\Ref{IV.6.7}, et ne n\'ecessite donc ni la
technique
\marginpar{104}
des constructibles, ni le th\'eor\`eme~\Ref{IV.5.6}):
\begin{corollaire}
\label{IV.6.11}
Sous les conditions de~\Ref{IV.6.10}, si on suppose $Y$ int\`egre,
alors il existe un ouvert non vide~$V$ dans $Y$ tel que $F$ soit plat
relativement \`a $Y$ en tous les points de~$f^{-1}(V)$.
\end{corollaire}
En effet, l'ensemble ouvert $U$ contient la fibre du point
g\'en\'erique de~$Y$ (puisque l'anneau local de ce point est un
corps), donc il contient un ouvert de la forme $f^{-1}(V)$, $X$
\'etant de type fini sur~$Y$. De~\Ref{IV.6.11}, on conclut aussi
facilement le r\'esultat suivant, o\`u~$Y$ est suppos\'e
noeth\'erien (mais pas n\'ecessairement int\`egre): il existe
une partition de~$Y$ en des parties localement ferm\'ees $Y_i$
telles que (munissant $Y_i$ de la structure r\'eduite induite) $F$
induise sur chaque $X_i=X\times_Y Y_i$ un faisceau plat par rapport
\`a~$Y_i$.

\chapter{Le groupe fondamental: g\'en\'eralit\'es}
\label{V}
\marginpar{105}

\setcounter{section}{-1}
\section{Introduction}

Le pr\'esent S\'eminaire est la suite du S\'eminaire 1960. Nous
r\'ef\'erons \`a ce dernier par des sigles tels que I~\Ref{I.9.7}
qui signifie: S\'eminaire de G\'eom\'etrie Alg\'ebrique,
expos\'e~I, \No9.7. Les num\'eros des expos\'es de 1961 suivront
ceux de 1960. Nous r\'ef\'erons aux \'El\'ements de
G\'eom\'etrie Alg\'ebrique de Dieudonn\'e-Grothendieck par des
sigles tels que (EGA~I 8.7.3).

Le pr\'esent expos\'e r\'esume (avec de l\'egers
compl\'ements) les derniers expos\'es de 1960, qui n'avaient pas
\'et\'e r\'edig\'es.

Comme en 1961, nous nous limiterons en r\`egle g\'en\'erale
\`a des pr\'esch\'emas localement noeth\'eriens, bien que
souvent cette restriction soit inessentielle. Nous admettrons dans
l'expos\'e \Ref{VI} la th\'eorie de la descente fid\`element plate,
r\'esum\'ee dans S\'eminaire Bourbaki \No190. S'il y a lieu,
nous en donnerons un expos\'e plus d\'etaill\'e dans un
expos\'e ult\'erieur\footnote{\Cf Exp\ptbl \Ref{VI} et Exp\ptbl \Ref{VIII}}, une fois
que le lecteur aura eu l'occasion de se convaincre de l'utilit\'e de
cette technique, pour la th\'eorie du groupe fondamental.

\section{Pr\'esch\'ema \`a groupe fini d'op\'erateurs,
pr\'esch\'ema quotient}
\label{V.1}

Soient $X$ un pr\'esch\'ema, $G$ un groupe fini op\'erant sur~$X$ par automorphismes, \`a droite pour fixer les id\'ees. Si $X$
est affine d'anneau $A$, $G$ op\`ere donc par automorphismes \`a
gauche sur~$A$.

Pour tout pr\'esch\'ema $Z$, $G$ op\`ere \`a gauche sur
l'ensemble $\Hom(X,Z)$,
on peut donc consid\'erer l'ensemble
$$
\Hom(X,Z)^G
$$
des morphismes invariants par $G$. Il d\'epend fonctoriellement de
$Z$,
on peut se demander si ce foncteur est \og repr\'esentable\fg,
\ie isomorphe \`a un foncteur $Z\mto\Hom(Y,Z)$.
\marginpar{106}
Cela signifie qu'on peut trouver un pr\'esch\'ema $Y$, et un
morphisme invariant par $G$
$$
p\colon X\to Y
$$
tel que pour tout $Z$, l'application correspondante $g\mto gp$
$$
\Hom(Y,Z)\to \Hom(X,Z)^G
$$
soit bijective. On dit alors que $(Y,p)$ est un
\emph{pr\'esch\'ema quotient}
\index{quotient (pr\'esch\'ema)|hyperpage}%
de~$X$
\ifthenelse{\boolean{orig}}{par~$g$}{par~$G$}
(il est d\'etermin\'e \`a isomorphisme unique pr\`es).

\begin{proposition}
\label{V.1.1}
Soient $A$ un anneau sur lequel le groupe fini $G$ op\`ere \`a
gauche, $B=A^G$ le sous-anneau des invariants de~$A$, $X=\Spec(A)$ et
$Y=\Spec(B)$, $p\colon X\to Y$ le morphisme canonique (\'evidemment
invariant par $G$). Alors
\begin{enumerate}
\item[(i)] $A$ est entier sur~$B$ \ie $p$ est un morphisme
\emph{entier}.
\item[(ii)] Le morphisme $p$ est surjectif, ses fibres sont les
trajectoires de~$G$, la topologie de~$Y$ est quotient de celle de~$X$.
\item[(iii)] Soit $x\in X$, $y=p(x)$, $G_x$ le stabilisateur de~$x$, alors
\ifthenelse{\boolean{orig}}
{$\kappa(x)$ est une extension alg\'ebrique quasi-galoisienne de
$\kappa(y)$ et l'application canonique de~$G_x$ dans le groupe
$\Gal(\kappa(x)/\kappa(y))$ des $\kappa(y)$-automorphismes de
$\kappa(x)$}
{$\kres(x)$ est une extension alg\'ebrique quasi-galoisienne de
$\kres(y)$ et l'application canonique de~$G_x$ dans le groupe
$\Gal(\kres(x)/\kres(y))$ des $\kres(y)$-automorphismes de
$\kres(x)$}
est
\ifthenelse{\boolean{orig}}{surjectif.}{surjective.}
\item[(iv)] $(Y,p)$ est un pr\'esch\'ema quotient de~$X$ par $G$.
\end{enumerate}
\end{proposition}

Les \'enonc\'es
\ifthenelse{\boolean{orig}}
{(i) (ii) (iii)}
{(i), (ii), (iii)}
sont bien connus en alg\`ebre
commutative\footnote{\Cf N. Bourbaki, Alg. Comm. Chap\ptbl 5, \S 1 et \S
2, th\ptbl 2.} et sont mis seulement pour m\'emoire, sauf l'assertion
sur la topologie, qui provient du fait g\'en\'eral suivant,
cons\'equence facile du th\'eor\`eme de Cohen-Seidenberg: un
morphisme entier est ferm\'e (\ie transforme ferm\'es en
ferm\'es). Notons tout de suite:

\begin{corollaire}
\label{V.1.2}
Sous les conditions pr\'ec\'edentes, l'homomorphisme naturel
$$\cal{O}_Y{\to} p_*(\cal{O}_X)^G$$ est un isomorphisme.
\end{corollaire}

Cela r\'esulte aussit\^ot de la formule
$$
(S^{-1}A)^G = S^{-1}(A^G)
$$
valable pour toute partie multiplicativement stable $S$ de~$B=A^G$
(formule qui se module, et s'\'enonce plus g\'en\'eralement pour
un changement de base $A\to A'$ qui est \emph{plat}),
\ifthenelse{\boolean{orig}}{appliqu\'e}{appliqu\'ee}
au cas o\`u $S$ est engendr\'e par un \'el\'ement $f$ de~$B$.

L'assertion (ii) et cor\ptbl \Ref{V.1.2} impliquent facilement (iv); plus
g\'en\'eralement, on aura ceci:

\begin{proposition}
\label{V.1.3}
Soient
\marginpar{107}
$X$ un pr\'esch\'ema \`a groupe d'automorphismes finis
$G$, $p\colon{}X\to{}Y$ un morphisme affine invariant tel que
$\cal{O}_Y\isomto p_*(\cal{O}_X)^G$. Alors les conclusions
\ifthenelse{\boolean{orig}}
{\textup{(i) (ii) (iii) (iv)}}
{\textup{(i), (ii), (iii), (iv)}}
de~\Ref{V.1.1} sont encore valables.
\end{proposition}

En effet, pour
\ifthenelse{\boolean{orig}}
{(i) (ii) (iii)}
{(i), (ii), (iii),}
on peut supposer $Y$ donc $X$ affine, et
si $B,A$ sont leurs anneaux, l'hypoth\`ese implique $B=A^G$, il
suffit d'appliquer~\Ref{V.1.1}. Pour (iv), on utilise (ii) et
$\cal{O}_Y=p_*(\cal{O}_X)^G$.

\begin{corollaire}
\label{V.1.4}
Sous les conditions de~\Ref{V.1.3}, pour tout ouvert $U$ de~$Y$, $U$
est un quotient de~$X|U=p^{-1}(U)$ par $G$.
\end{corollaire}

En effet, $p^{-1}(U)\to{}U$ induit par $p$ satisfait aux m\^emes
hypoth\`eses que~$p$.

Si maintenant $X$ est un $Z$-pr\'esch\'ema et les op\'erations
de~$G$ sont des $Z$-automorphismes, alors par (iv) $Y$ est un
$Z$-pr\'esch\'ema. Ceci dit:

\begin{corollaire}
\label{V.1.5}
Pour que $X$ soit affine \resp s\'epar\'e sur~$Z$, il faut et il
suffit que~$Y$ le soit. Si $X$ est de type fini sur~$Z$, il est
\emph{fini} sur~$Y$; si de plus $Z$ est localement noeth\'erien, $Y$
est de type fini sur~$Z$.
\end{corollaire}

Comme $X$ est affine et a fortiori s\'epar\'e sur~$Y$, si $Y$ est
affine \resp s\'epar\'e sur~$Z$, $X$~l'est
aussi. R\'eciproquement, supposons $X$ affine sur~$Z$, prouvons que
$Y$ l'est: on peut gr\^ace \`a~\Ref{V.1.4} supposer $Z$ affine, et
on est ramen\'e \`a prouver que si $X$ est affine, $Y$~l'est, ce
qui r\'esulte de la d\'etermination explicite de~$Y$ comme
$\Spec(A)^G$ faite dans~\Ref{V.1.1}. De m\^eme, comme
$p\colon{}X\to{}Y$ est entier donc universellement ferm\'e, et
surjectif, il s'ensuit que si $X$ est s\'epar\'e sur~$Z$, $Y$
l'est aussi (lemme \`a d\'egager!); en effet, dans le diagramme
$$
\xymatrix@C=1.5cm{ X \times_ZX\ar[r]^-{p\times_Zp} & Y\times_ZY \\ X
\ar[u]^{\Delta_{X/Z}}\ar[r]^{p} & Y \ar[u]_{\Delta_{Y/Z}} }
$$
le morphisme $X\times_ZX\to Y\times_ZY$ est ferm\'e, donc transforme
la diagonale (ferm\'ee) de~$X\times_ZX$ en une partie ferm\'ee de
$Y\times_ZY$, qui n'est d'ailleurs autre que la diagonale de ce
dernier puisque $p$ est surjectif.
--- Si $X$ est de type fini sur~$Z$, il l'est a fortiori sur~$Y$ donc
il est fini sur~$Y$ (puisqu'il est d\'ej\`a entier sur~$Y$). Supposons de plus $Z$ localement noeth\'erien, prouvons que
$Y$ est de type fini sur~$Z$. On peut gr\^ace \`a~\Ref{V.1.4}
supposer $Z$ affine. Comme l'espace topologique $X$ est quasi-compact
et que $p\colon X\to Y$ est surjectif, $Y$ est
\ifthenelse{\boolean{orig}}{egalement}{\'egalement}
quasi-compact donc r\'eunion finie d'ouverts
\marginpar{108}
affines, et par~\Ref{V.1.4} on est ramen\'e au cas o\`u $Y$ est
affine, donc $X$ affine. Mais alors l'anneau $A$ de~$X$ est une
alg\`ebre de type fini sur l'anneau $C$ de~$Z$ qui est
noeth\'erien, et il est connu que $B=A^G$ est alors \'egalement
une alg\`ebre de type fini sur~$C$ (car $A$ sera enti\`ere, donc
finie sur une sous-alg\`ebre $B'$ de~$B$ de type fini sur~$C$, donc
comme $B'$ est noeth\'erien, $B$ est \'egalement
\ifthenelse{\boolean{orig}}{fini}{finie}
sur~$B'$, donc de type fini sur~$C$).

\begin{corollaire}
\label{V.1.6}
Pour que $X$ soit affine \resp un sch\'ema, il \fets que $Y$ le
soit.
\end{corollaire}

\begin{definition}
\label{V.1.7}
Soit
\marginpar{109}
$X$ un pr\'esch\'ema o\`u un groupe fini $G$ op\`ere
\`a droite. On dit que $G$ \emph{op\`ere de fa\c con admissible}
s'il existe un morphisme $p\colon X\to Y$ ayant les propri\'et\'es
de~\Ref{V.1.3} (ce qui implique que $X/G$ existe et est isomorphe
\`a~$Y$).
\end{definition}

\begin{proposition}
\label{V.1.8}
Soit $X$ un pr\'esch\'ema o\`u le groupe fini $G$ op\`ere
\`a droite. Pour que $G$ op\`ere de fa\c con admissible, il faut
et il suffit que $X$ soit r\'eunion d'ouverts affines invariants par
$G$, ou encore que toute trajectoire de~$G$ dans $X$ soit contenue
dans un ouvert affine.
\end{proposition}

Cette derni\`ere condition est \'evidemment impliqu\'ee par la
premi\`ere, et \`a son tour elle l'implique; car soit $T$ une
trajectoire de~$G$, $U$ un ouvert affine la contenant, l'intersection
des transform\'es de~$U$ par les
\ifthenelse{\boolean{orig}}{$g$ $G$}{$g$ dans $G$}
est alors un ouvert $U'$ stable par~$G$, contenant $T$ et contenu dans
l'ouvert affine $U$. Comme dans $U$, toute partie finie a un
syst\`eme fondamental de voisinages ouverts affines, il existe un
voisinage ouvert affine $V$ de~$T$ contenu dans $U'$. Ces
transform\'es par les
\ifthenelse{\boolean{orig}}{$g$ $G$}{$g$ dans $G$}
sont donc affines et contenus dans $U'$ qui est \emph{s\'epar\'e},
donc leur intersection $U''$ est un ouvert affine qui est invariant
par $G$ et contient $T$.
--- Ceci pos\'e, la condition envisag\'ee dans~\Ref{V.1.8} est
\emph{n\'ecessaire}, car on prendra les images inverses $X_i$
d'ouverts affines $Y_i$ recouvrant~$Y$. Elle est suffisante, car on
peut alors par~\Ref{V.1.1} construire les quotients $Y_i=X_i/G$; dans
chaque $Y_i$ l'image de~$X_i\cap X_j$ est un ouvert $Y_{ij}$
s'identifiant \`a $X_{ij}/G$ par~\Ref{V.1.4}, en particulier on en
d\'eduit des isomorphismes $Y_{ij}\isomto Y_{ji}$ permettant de
recoller les $Y_i$ pour construire $Y$. Serre pr\'ef\`ere
construire directement l'espace topologique quotient $Y$ de~$X$ par
$G$, mettre dessus le faisceau $p_*(\cal{O}_X)^G$ et v\'erifier que
$Y$ devient un pr\'esch\'ema et qu'on est alors sous les
conditions de~\Ref{V.1.3}.


\setcounter{subsection}{6}

\begin{corollaire}
\label{cor:V.1.7}
Si $G$ op\'erant sur~$X$ est admissible, il en est de m\^eme pour
tout sous-groupe $H$ de~$G$ (donc $X/H$ existe).
\end{corollaire}

Cela peut aussi se v\'erifier directement sur la
situation~\Ref{V.1.3}, en notant qu'on peut toujours supposer $X$
affine sur un $Z$ et les $s\in G$ op\`erent par $Z$-automorphismes
(on prend par exemple $Z=Y$); on a en effet:

\begin{corollaire}
\label{cor:V.1.8}
Supposons $X$ affine sur~$Z$, et les op\'erations de~$G$ des
$Z$-auto\-mor\-phismes. Alors $G$ op\`ere sur~$X$ de fa\c con
admissible. Si $X$ est d\'efini par un faisceau quasi-coh\'erent
$\cal{A}$ d'alg\`ebres, $Y$ est d\'efini par le faisceau
$\cal{A}^G$ des invariants de~$\cal{A}$ par~$G$.
\end{corollaire}

\begin{proposition}
\label{V.1.9}
Supposons que $G$ op\`ere de fa\c con admissible sur~$X$, et que
$X/G=Y$ soit un pr\'esch\'ema sur~$Z$. Consid\'erons un
morphisme de changement de base $Z'\to Z$, posons $X'=X\times_Z Z',
Y'=Y\times_Z Z'$, de sorte que $G$ op\`ere encore par transport de
structure sur~$X'$, le morphisme $p'\colon X'\to Y'$ \'etant
invariant. Si $Z'$ est \emph{plat} sur~$Z$, alors $p'$ satisfait
encore les hypoth\`eses
\ifthenelse{\boolean{orig}}
{de~\Ref{V.1.3}}
{de~\Ref{V.1.3},}
\ie $\cal{O}_Y'\to
p_*'(\cal{O}_{X'})^G$ est un isomorphisme ($p'$ \'etant de toutes
fa\c cons affine). Donc $G$ op\`ere de fa\c con admissible sur~$X'$, et $(X/G)\times_Z Z'\thickapprox (X\times_Z Z')/G$.
\end{proposition}

On peut \'evidemment supposer $Z=Y$, on est ramen\'e au cas o\`u
de plus $Y$ et $Y'$ sont affines. Il faut montrer que si $B$ est le
sous-anneau des invariants de~$G$ op\'erant dans $A$, et si $B'$ est
une alg\`ebre sur~$B$ plate sur~$B$, alors $B'$ est la
sous-alg\`ebre des invariants de~$A'=A\otimes_B B'$. C'est
imm\'ediat, car la suite exacte
\ifthenelse{\boolean{orig}}{}
{\enlargethispage{.5cm}}%
$$
0 \to B \lto{i} A \lto{j} A^{(G)}
$$
(o\`u le dernier terme signifie une puissance de~$A$, et o\`u
$j(x)$ est le syst\`eme des $s\cdot x-x$, $s\in G$) reste exacte par
tensorisation par
\ifthenelse{\boolean{orig}}
{$A'$}
{$B'$}.

On fera attention que l'hypoth\`ese de platitude \'etait
essentielle pour la validit\'e du r\'esultat; en particulier, si
$Y'$ est un sous-pr\'esch\'ema ferm\'e de
\ifthenelse{\boolean{orig}}
{$X$}
{$Y$}
(par exemple m\^eme un point ferm\'e de
\ifthenelse{\boolean{orig}}
{$X$}
{$Y$}),
$X'$ son image inverse dans $X$,
alors $Y'$ \emph{ne s'identifie pas} en g\'en\'eral \`a
$X'/G$. Nous verrons qu'il en est n\'eanmoins ainsi si $X$ est
\'etale sur~$Y$.

Pour finir, donnons un formalisme aussi commode que trivial. Soit $Y$
un pr\'esch\'ema. Comme dans la cat\'egorie des
pr\'esch\'emas, les sommes directes existent, on peut pour tout
ensemble $E$ consid\'erer le pr\'esch\'ema somme d'une famille
$(Y_i)_{i\in E}$ de pr\'esch\'emas tous identiques \`a $Y$, ce
pr\'esch\'ema sera not\'e $Y\times E$.
\marginpar{110}
Il est caract\'eris\'e par la formule
\begin{equation*}
\label{eq:V.1.*}
\tag{$*$} \Hom(Y\times E,Z) = \Hom(E,\Hom(Y,Z))
\end{equation*}
o\`u le deuxi\`eme $\Hom$ d\'esigne \'evidemment l'ensemble
des applications de l'ensemble $E$ dans l'ensemble $\Hom(Y,Z)$. On a
un morphisme canonique
$$
Y\times E\to Y
$$
faisant de~$Y\times E$ un pr\'esch\'ema sur~$Y$. Comme les
produits fibr\'es commutent aux sommes directes (dans la
cat\'egorie des pr\'esch\'emas) on aura, si $Y$ est un
pr\'esch\'ema sur un autre $Z$, pour un changement de base $Z'\to
Z$:
$$
(Y\times E)\times_Z Z'=(Y\times_Z Z')\times E
$$
(formule surtout utile si $Z=Y$). D'autre part, on conclut
trivialement de la d\'efinition
$$
(Y\times E)\times F=Y\times(E\times F)=(Y\times E)\times_Y(Y\times F)
$$
(la derni\`ere formule cependant r\'esultant de la
commutativit\'e signal\'ee plus haut).

Pour $Y$ fix\'e, on peut regarder $Y\times E$ comme un foncteur en
$E$, \`a valeurs dans les pr\'esch\'emas sur~$Y$, foncteur qui
commute aux produits finis d'apr\`es la formule pr\'ec\'edente,
(ce qui permet par exemple \`a tout groupe ordinaire $G$ de faire
correspondre un sch\'ema en groupes $Y\times G$ sur~$Y$, qui sera
fini sur~$Y$ si $Y$ l'est, \ifthenelse{\boolean{orig}}{etc...}{etc.}). Plus g\'en\'eralement, ce
foncteur est \og exact \`a gauche\fg, mais nous n'aurons pas \`a nous
en servir ici. Ce foncteur commute aussi trivialement aux sommes
directes, et il est aussi \og exact \`a droite\fg, comme on voit
aussit\^ot sur la formule de d\'efinition \eqref{eq:V.1.*}. En
particulier, si le groupe fini $G$ op\`ere \`a droite dans
l'ensemble $E$, alors il op\`ere \`a droite dans $Y\times E$, et
on a
$$
(Y\times E)/G = Y\times (E/G)
$$
o\`u en fait le quotient du premier membre satisfait aux conditions
de~\Ref{V.1.3} (c'est imm\'ediat).

\section{Groupes de d\'ecomposition et d'inertie. Cas \'etale}
\label{V.2}
Soit $G$
\ifthenelse{\boolean{orig}}{\ignorespaces}{un}
groupe fini op\'erant \`a droite sur le pr\'esch\'ema $X$. Si
$x\in X$, on appelle
\marginpar{111}
\emph{groupe de d\'ecomposition de~$x$}
\index{decomposition (groupe de)@d\'ecomposition (groupe de)|hyperpage}%
\index{groupe de d\'ecomposition|hyperpage}%
le stabilisateur $G_d(x)$ de~$x$. Ce groupe op\`ere canoniquement
(\`a gauche) sur le corps r\'esiduel
\ifthenelse{\boolean{orig}}
{$\kappa(x)$}
{$\kres(x)$}, et l'ensemble
des \'el\'ements de~$G_d(x)$ qui op\`erent trivialement est
appel\'e \emph{groupe d'inertie}
\index{groupe d'inertie|hyperpage}%
\index{inertie (groupe d')|hyperpage}%
de~$x$, not\'e~$G_i(x)$.

Supposons que $G$ op\`ere sur~$X$ de fa\c con admissible et que
$Y$ soit un pr\'esch\'ema sur un pr\'esch\'ema
$Z$. Fixons-nous un $z\in Z$, et une extension alg\'ebriquement
close $\Omega$
\ifthenelse{\boolean{orig}}
{de~$\kappa(z)$ ayant un degr\'e de transcendance
sup\'erieur \`a celui des $\kappa(x)/\kappa(z)$}
{de~$\kres(z)$ ayant un degr\'e de transcendance
sup\'erieur \`a celui des $\kres(x)/\kres(z)$}, o\`u $x$ est un
point de~$X$ au-dessus de~$z$. On peut regarder $\Spec(\Omega)$ comme
un $Z$-sch\'ema, et les points de~$X$ \`a valeurs dans $\Omega$
correspondent aux homomorphismes
\ifthenelse{\boolean{orig}}
{de~$\kappa(z)$-alg\`ebres $\kappa(x)\to\Omega$}
{de~$\kres(z)$-alg\`ebres $\kres(x)\to\Omega$},
o\`u $x$ est un point de~$X$ au-dessus de~$z$;
comme $\Omega$ a \'et\'e prise assez grande, tout point $x$ de~$X$
au-dessus de~$z$ est la localit\'e d'un point de~$X$ \`a valeurs
dans $\Omega$. Si $X(\Omega)$ et $Y(\Omega)$ d\'esignent
respectivement l'ensemble des points de~$X$ et $Y$ \`a valeurs dans
$\Omega$, on a une application naturelle
$$
X(\Omega)\to Y(\Omega),
$$
d'autre part, $G$ op\`ere sur~$X(\Omega)$ et l'application
pr\'ec\'edente est invariante par $G$. Ceci pos\'e, les
conclusions (ii) et (iii) de~\Ref{V.1.3} s'interpr\`etent aussi
ainsi: \emph{l'application pr\'ec\'edente est surjective et
identifie $Y(\Omega)$ au quotient $X(\Omega)/G$. De plus, si $x$ est
la localit\'e de~$a\in X(\Omega)$, alors le stabilisateur de~$a$
dans $G$ n'est autre que le groupe d'inertie $G_i(x)$}. Tout ceci est
d'ailleurs vrai sans supposer $\Omega$ \og assez grand\fg, cette
derni\`ere hypoth\`ese sert uniquement \`a assurer qu'on peut
caract\'eriser le groupe d'inertie de tout \'el\'ement de~$X$
au-dessus de~$z$ comme un stabilisateur \og g\'eom\'etrique\fg. On en
conclut par exemple aussit\^ot:
\begin{proposition}
\label{V.2.1}
Faisons une extension de la base $Z'\to Z$, d'o\`u
$X'=X\times_ZZ'$. Soit $x'$ un point de~$X'$, $x$ son image dans $X$,
alors on a $G_i(x)=G_i(x')$.
\end{proposition}

Il suffit, dans les consid\'erations ci-dessus, de prendre pour
$\Omega$ une extension assez grande
\ifthenelse{\boolean{orig}}
{de~$\kappa(z')$}
{de~$\kres(z')$}
(o\`u $z,z'$ sont les images de~$x,x'$ dans $Z,Z'$).

\begin{proposition}
\label{V.2.2}
Sous les conditions de~\Ref{V.1.3}, supposons $Y$ localement
noeth\'erien, $X$ fini sur~$Y$. Soit $H$ un sous-groupe de~$G$,
consid\'erons $X'=X/H$ (\cf \Ref{cor:V.1.7}), soit
\marginpar{112}
$x\in X$, $x'$ son image dans $X'$ et $y$ son image dans $Y$.
\begin{enumerate}
\item[(i)] Si $H\supset G_d(x)$, alors l'homomorphisme
$\cal{O}_y\to\cal{O}_{x'}$ induit un isomorphisme sur les
compl\'et\'es.
\item[(ii)] Si $H\supset G_i(x)$, alors l'homomorphisme $\cal{O}_y\to
\cal{O}_{x'}$ est \'etale \ie $X'$ est \'etale sur~$Y$ en $x'$.
\end{enumerate}
\end{proposition}
Soit $Y_1=\Spec(\widehat{\cal{O}_y})$, faisons le changement de base
$Y_1\to Y$, on trouve un $X_1=X\times_Y Y_1$ fini sur~$Y_1$, sur
lequel $G$ op\`ere, le quotient \'etant $Y_1$
par~\Ref{V.1.9}. Soit~$y_1$ l'unique point de~$Y_1$ au-dessus de~$y$,
comme
\ifthenelse{\boolean{orig}}
{$\kappa(y)=\kappa(y_1)$}
{$\kres(y)=\kres(y_1)$},
il s'ensuit que la fibre de~$X$ en~$y$
est isomorphe \`a celle de~$X_1$ en $y_1$, d'o\`u un unique point
$x_1$ de~$X_1$ au-dessus de~$x$. D'ailleurs par~\Ref{V.1.9} on aura
$X_1/H=X'_1=X'\times_Y Y_1$, soit $x'_1$ l'image de~$x_1$ dans $X'_1$,
il est au-dessus de~$x'$, et on v\'erifie facilement ($X'$ \'etant
de type fini sur~$Y$) que l'homomorphisme
$\cal{O}_{x'}\to\cal{O}_{x'_1}$ induit un isomorphisme sur les
compl\'et\'es. Donc on est ramen\'e au cas o\`u $Y$ est le
spectre d'un anneau local complet, soit $B$, donc~$X$ le spectre d'un
anneau fini $A$ sur~$B$, compos\'e d'un nombre fini d'anneaux locaux~$A_x$ \ifthenelse{\boolean{orig}}{correspondants}{correspondant}
aux points $x_i$ de~$X$ sur~$Y$. Si $A_0$ correspond \`a $x=x_0$,
alors $A$ s'identifie \`a l'anneau $\Hom_{G_d}(G,A_0)$ des fonctions
$f\colon G\to A_0$ telles que $f(st)=sf(t)$ pour $s\in G_d$, les
op\'erations de~$G$ sur ces fonctions \'etant d\'efinies par
$(uf)(t)=f(tu)$. On voit donc que si $H$ est un sous-groupe quelconque
de~$G$, alors $A^H$ est l'anneau des fonctions $f\colon G\to A_0$
telles que
$$
f(stu)=sf(t)\qquad\text{pour}\ s\in G_d,\ u\in H
$$
donc c'est un anneau semi-local dont les composants locaux
correspondent aux doubles classes $G_daH$ dans $G$, \`a la double
classe d\'efinie par $a\in G$ correspondant (gr\^ace \`a
l'application $f\mto f(a)$) le sous-anneau $A_0^{H(a)}$ de~$A_0$,
o\`u $H(a)=G_d\cap aHa^{-1}$. D'ailleurs, le composant local de
$A^H$ correspondant \`a l'image $x'$ de~$x$ est aussi celui
correspondant \`a la double classe $G_dH$ de l'\'el\'ement
neutre, son composant local est donc $A_0^{G_d\cap H}$. Si donc
\ifthenelse{\boolean{orig}}{$G_d\in H$,}{$G_d\subset H$,}
on trouve $A_0^{G_d}=A^G=B$, ce qui prouve (i). Pour
prouver (ii), on peut, en passant \`a une extension finie plate
convenable de~$A$, et utilisant~\Ref{V.2.1}, se ramener au cas o\`u
l'extension r\'esiduelle
\ifthenelse{\boolean{orig}}
{$\kappa(x)/\kappa(y)$}
{$\kres(x)/\kres(y)$}
est triviale. Mais
alors $G_i(x)=G_d(x)$, et on est ramen\'e au cas pr\'ec\'edent.

\begin{corollaire}
\label{V.2.3}
Sous
\marginpar{113}
les conditions de~\Ref{V.2.2}, supposons $G_i(x)=(e)$, alors X
est \'etale sur~$Y$ en $x$. Donc si $G_i(x)=(e)$ pour tout $x \in
X$, alors $X \to Y$ est un morphisme \'etale.
\end{corollaire}

Il y a une r\'eciproque partielle:

\begin{corollaire}
\label{V.2.4} Supposons $X$ connexe et le groupe $G$ fid\`ele
sur~$X$. Pour que $p \colon X \to Y=X/G$ soit \'etale, il faut et il
suffit que les groupes d'inertie des points de~$X$ soient r\'eduits
\`a l'\'el\'ement neutre. S'il en est ainsi, $G$ s'identifie au
groupe de tous les $Y$-automorphismes du $Y$-sch\'ema~$X$.
\end{corollaire}

Compte tenu de~\Ref{V.2.3}, on peut supposer $X$ \'etale sur~$Y$. Mais si un $s \in G$ est dans un $G_i(x)$, il r\'esulte alors
de I~\Ref{I.5.4} que $s$ op\`ere trivialement sur~$G$, donc est
l'\'el\'ement unit\'e puisque $G$ est fid\`ele, ce qui prouve
la premi\`ere assertion. Soit $u$ un $Y$-automorphisme de~$X$, soit
$x \in X$. D'apr\`es la proposition~\Ref{V.1.3}, il existe un $s \in
G$ tel que $s(x)=u(x)$, et induisant le m\^eme homomorphisme
r\'esiduel
\ifthenelse{\boolean{orig}}
{$\kappa(x)\to\kappa(x')$}
{$\kres(x)/\kres(x')$}
que $u$. Par \loccit on a donc
$s=u$, ce qui ach\`eve la d\'emonstration.

\begin{remarque}
\label{V.2.5}
L'hypoth\`ese que $G$ op\`ere fid\`element n'est \'evidemment
par surabondante dans le corollaire~\Ref{V.2.4}. Il en est de m\^eme
de l'hypoth\`ese que $X$ est connexe, comme on voit par exemple en
prenant $X=Y \times E$, $E$ \'etant un ensemble fini, et $G$ le
groupe des permutations de~$E$: $G$ op\`ere avec force inertie,
n\'eanmoins $(Y \times E)/G=Y \times (E/G)=Y$, et $X$ est \'etale
sur~$Y$. Prenant pour $G$ un groupe strictement plus petit que le
groupe sym\'etrique de~$E$, mais op\'erant transitivement sur~$E$,
on voit qu'il y aura aussi des $Y$-automorphismes de~$X$ ne provenant
pas de~$G$.
\end{remarque}

L'exemple type d'un groupe $G$ op\'erant sans inertie est celui de
$Y \times G$, sur lequel on fait op\'erer $G$ gr\^ace \`a ses
op\'erations sur le facteur $G$ par translations \`a droite: un
$Y$-pr\'esch\'ema $X$ \`a groupe d'op\'erateurs \`a droite
$G$ est dit \emph{trivial} s'il est isomorphe \`a $Y \times G$.

Pour faire le lien entre le pr\'esch\'emas \`a groupes finis
d'op\'erateurs et la notion de fibr\'e principal dans une
cat\'egorie (lien dont nous n'aurons pas besoin d'ailleurs pour la
suite du s\'eminaire, mais important dans d'autres contextes) les
consid\'erations suivantes sont utiles. Nous fixons un
pr\'esch\'ema de base $Y$, et nous pla\c cons dans le
cat\'egorie des $Y$-pr\'esch\'emas. Si $G$ est un groupe
\marginpar{114}
fini, nous poserons pour abr\'eger $G_Y=Y \times G$, c'est donc un
sch\'ema en groupes finis sur~$Y$ (\cf \numero 1), et si $X$ est
un $Y$-pr\'esch\'ema, on a
$$
X \times_Y G_Y = X \times G
$$
(m\^eme r\'ef\'erence). La donn\'ee d'un $Y$-morphisme $X
\times_Y G_Y \to X$ \'equivaut donc \`a la donn\'ee d'un
$Y$-morphisme $X \times G \to X$, \ie \`a la donn\'ee pour tout
$g \in G$ d'un $Y$-morphisme $T_g \colon X \to X$. On constate
aussit\^ot que pour que la donn\'ee des $T_g$ d\'efinisse sur~$X$ une structure de pr\'esch\'ema \`a groupe d'op\'erateurs
\`a droite $G$, (\ie $T_{gg'}=T_{g'}T_g$, $T_e=\id_x$) il faut et il
suffit que le $Y$-morphisme correspondant $X \times_Y G_Y \to X$
d\'efinisse sur~$X$ une structure de~$Y$-pr\'esch\'ema \`a
$Y$-sch\'ema en groupes d'op\'erateurs, (au sens g\'en\'eral
des objets \`a $\cal{C}$-groupe d'op\'erateurs dans une
cat\'egorie $\cal{C}$). Supposons qu'il en soit ainsi. Rappelons que
$X$ est dit \emph{formellement principal homog\`ene}
\index{formellement principal homog\`ene (pr\'esch\'ema)|hyperpage}%
sous~$G_Y$\kern1pt\footnote{on dit plut\^ot maintenant: $X$ est un
pseudo-torseur sous $G_Y$.} si le morphisme canonique
$$
X \times_Y G_Y \to X \times_Y X
$$
dont les composantes sont respectivement $\pr_1$ et le morphisme de
multiplication $\pi \colon X \times_Y G_Y \to X$, est un
isomorphisme. En l'occurrence, identifiant le premier membre \`a $X
\times G$, le morphisme consid\'er\'e est celui qui, \`a tout $g
\in G$, associe le morphisme
$$
(\id_X,T_g)=(\id_X \times_Y T_g)\Delta_{X/Y} \colon X \to X \times_Y X
$$
et par suite, dire que $X$ est formellement principal homog\`ene
sous $G_Y$ signifie aussi que $X \times_Y X$ est isomorphe \`a la
somme directe des transform\'ees de la diagonale par les
\'el\'ements $(e,g)$ de~$G \times G$ (op\'erant sur~$X \times_Y
X$ de fa\c con \'evidente) o\`u $e$ d\'esigne l'\'el\'ement
unit\'e de~$G$. Si on ne veut pas distinguer la gauche et la droite
et donner une formule qui reste applicable au produit de plus de deux
facteurs identiques \`a $X$, on peut formuler la condition en disant
que le morphisme canonique
$$
X \times_G (G \times G) \to X \times_Y X
$$
obtenu en attachant au couple $(g,g')$ le morphisme
$$
(T_g,T_{g'})=(T_g \times_Y T_{g'})\Delta_{X/Y} \colon X \to X \times_Y
X
$$
et en faisant op\'erer $G$ \`a gauche sur~$G \times G$ par
l'homomorphisme diagonal;
$$
s(g,g')=(sg,sg'),
$$
est un \emph{isomorphisme}.

La
\marginpar{115}
notion d'\emph{espace principal homog\`ene}
\index{principal homog\`ene (pr\'esch\'ema, fibr\'e)|hyperpage}%
est d\'eduite de celle de l'espace formellement principal
homog\`ene en ajoutant un axiome suppl\'ementaire, assurant que le
\og quotient\fg de~$X$ par $G_Y$ existe et est pr\'ecis\'ement
l'objet unit\'e \`a droite de la cat\'egorie, ici $Y$. Cet
axiome peut varier suivant le contexte, et s'explicite souvent le plus
commod\'ement (dans le yoga de la \og descente\fg) en exigeant que
l'objet \`a op\'erateurs devienne \og trivial\fg \ie isomorphe au
produit $X \times_Y G_Y$ (en l'occurrence $X \times G$) par un
changement de base convenable, de type pr\'ecis\'e (de telle fa\c con, en pratique, \`a permettre la technique de descente;
\cf Grothendieck, Technique de descente et th\'eor\`emes
d'existence en G\'eom\'etrie Alg\'ebrique, S\'em. Bourbaki
\No 190, pages 26 \`a 28)\footnote{\Cf Exp\ptbl \Ref{VIII} pour la
th\'eorie de la descente plate.}. Dans cet ordre d'id\'ees,
signalons ici la caract\'erisation des fibr\'es principaux
homog\`enes de groupe $G$ (au sens de \loccit):

\begin{proposition}
\label{V.2.6}
Soient $Y$ un pr\'esch\'ema localement noeth\'erien, $X$ un
$Y$-pr\'esch\'ema \`a groupe fini $G$ d'op\'erateurs
op\'erant \`a droite. Les conditions suivantes sont
\'equivalentes:
\begin{enumerate}
\item[(i)] $X$ est fini sur~$Y$, $Y=X/G$, les groupes d'inertie des
points de~$X$ sont r\'eduits \`a l'unit\'e.
\item[(ii)] Il existe un changement de base fid\`element plat et
quasi-compact $Y_1 \to Y$ tel que $X_1=X \times_Y Y_1$ soit un
$Y_1$-pr\'esch\'ema \`a op\'erateurs trivial, \ie isomorphe
\`a $Y_1 \times G$.
\item[(ii~bis)] Comme (ii), mais $Y_1 \to Y$ \'etant fini,
\'etale, surjectif.
\item[(iii)] $X$ est formellement principal homog\`ene sous $G_Y$,
et fid\`element plat et quasi-compact sur~$Y$.
\end{enumerate}
\end{proposition}
\textit{D\'emonstration} (i)$\To$(ii~bis) On prendra
$Y_1=X$, notant que $X \to Y$ est bien fini, \'etale par~\Ref{V.2.3}
et surjectif. Montrons que $X_1$ est alors trivial sur~$Y_1$, ce que
r\'esultera du

\begin{corollaire}
\label{V.2.7} Si \textup{(i)} est v\'erifi\'e et si $X$ admet une
section sur~$Y$, alors $X$ est un espace \`a op\'erateurs trivial.
\end{corollaire}

En effet, cette section permet de d\'efinir un $G$-morphisme $X
\times G \to X$, surjectif puisque $G$ est transitif sur les fibres de
$X$, injectif puisque $G$ op\`ere sans inertie; enfin, c'est un
isomorphisme local en vertu de I~\Ref{I.5.3} puisque $X$ est \'etale
sur~$Y$. Donc c'est un isomorphisme.

(ii~bis)
\marginpar{116}
implique trivialement (ii), qui implique (i) car les ingr\'edients
de (i) sont \og invariants\fg par extension fid\`element plate
quasi-compacte de la base (\cf S\'eminaire Bourbaki cit\'e plus
haut pour \og fini\fg; pour les groupes d'inertie, on
applique~\Ref{V.2.1}, et pour $Y=X/G$, une r\'eciproque
de~\Ref{V.1.9} dans le cas d'un changement de base
\emph{fid\`element plat}, que nous avions oubli\'e d'expliciter).

Nous avons prouv\'e (i)$\To$(iii) en passant en prouvant
(i)$\To$(ii~bis). Enfin (iii)$\To$(ii), car la
premi\`ere hypoth\`ese dans (iii) signifie pr\'ecis\'ement que
$X$ devient trivial en faisant le changement de base $Y_1=X$; d'o\`u
(ii) puisque $X$ est fid\`element plat et quasi-compact sur~$Y$.

\begin{definition}
\label{V.2.8} Un $Y$-pr\'esch\'ema $X$ \`a groupe
d'op\'erateurs \`a droite $G$ satisfaisant les conditions
\'equivalentes de~\Ref{V.2.6} est appel\'e un \emph{rev\^etement
principal de~$Y$, de groupe de Galois~$G$}.
\index{groupe de Galois d'un rev\^etement principal|hyperpage}%
\index{principal (revetement)@principal (rev\^etement)|hyperpage}%
\index{revetement principal@rev\^etement principal|hyperpage}%
\end{definition}

\section{Automorphismes et morphismes de rev\^etements \'etales}
\label{V.3}

\begin{proposition}
\label{V.3.1} Soit $X$ \'etale s\'epar\'e de type fini sur~$Y$ localement noeth\'erien, soit $G$ un groupe fini op\'erant sur~$X$ par $Y$-automorphismes. Alors $G$ op\`ere de fa\c con
admissible et le pr\'esch\'ema quotient $X/G$ est \'etale sur~$Y$.
\end{proposition}

On ne suppose pas $X$ fini sur~$Y$, cependant $X$ est quasi-projectif
sur~$Y$ d'o\`u l'existence de~$X/G$ gr\^ace
\`a~\Ref{V.1.8}. Prouvons d'abord

\begin{corollaire}
\label{V.3.2} Le morphisme $X \to X/G$ est \'etale.
\end{corollaire}

On peut supposer \'evidemment $G$ transitif sur l'ensemble des
composantes connexes de~$X$, puis par consid\'eration du
stabilisateur d'une composante connexe, que $X$ est m\^eme
connexe. Enfin, on peut supposer que $G$ op\`ere fid\`element.
Mais alors on voit comme dans~\Ref{V.2.4} que $G$ op\`ere sans
inertie, donc par~\Ref{V.2.3} il s'ensuit que $X \to X/G$ est
\'etale. On conclut gr\^ace au

\begin{lemme}[remords \`a l'expos\'e I]
\label{V.3.3} Soient $X \to X' \to
Y$ des morphismes de type fini, $x$ un point de~$X$, $x'$ et $y$ ses
images. On suppose $Y$ localement noeth\'erien. Si deux des
morphismes envisag\'es sont \'etales aux points marqu\'es, il en est
de m\^eme du troisi\`eme.
\end{lemme}

Il reste seulement \`a regarder le cas o\`u $X \to X'$ et $X \to
Y$ sont \'etales en $x$ et prouver que $X' \to Y$ l'est en $x'$ (ce
qui est le cas dont nous avons besoin
\marginpar{117}
pour~\Ref{V.3.1}). Faisant une extension plate convenable de la base
$Y$, on est ramen\'e au cas o\`u l'extension r\'esiduelle
\ifthenelse{\boolean{orig}}
{$\kappa(x)/\kappa(y)$}
{$\kres(x)/\kres(y)$}
est triviale. Consid\'erons les homomorphismes
$\cal{O}_y \to \cal{O}_{x'} \to \cal{O}_x$ et les homomorphismes
d\'eduits par passage aux compl\'et\'es, l'hypoth\`ese
signifie que $\hat{\cal{O}_y} \to \hat{\cal{O}_x}$ et
$\hat{\cal{O}_{x'}} \to \hat{\cal{O}_x}$ sont des isomorphismes,
d'o\`u aussit\^ot que $\hat{\cal{O}_y} \to \hat{\cal{O}_{x'}}$ en
est un, ce qui prouve le lemme.

\begin{corollaire}
\label{V.3.4} Si $X$ est fini et \'etale sur~$Y$, alors $X/G$
est fini et \'etale sur~$Y$.
\end{corollaire}

\begin{proposition}
\label{V.3.5} Soient $X$, $X'$ deux rev\^etements \'etales de
$Y$. Alors tout $Y$-morphisme $f \colon X \to X'$ se factorise en le
produit d'un morphisme \'etale surjectif $X \to X''$ et de
l'immersion canonique $X'' \to X'$ d'une partie $X''$ de~$X'$ \`a la
fois ouverte et ferm\'ee.
\end{proposition}

On sait (I~\Ref{I.4.8}) que $f$ est \'etale, donc un morphisme
ouvert, d'autre part $X$ \'etant fini sur~$Y$, $f$ est ferm\'e,
donc $f(X)=X''$ est une partie \`a la fois ouverte et ferm\'ee de~$X'$. On a fini (N.B. il suffisait que $X'$, au lieu d'un
rev\^etement \'etale, soit non ramifi\'e sur~$Y$).

\begin{corollaire}
\label{V.3.6} Avec les notations pr\'ec\'edentes, $X \to X'$
est un \'epimorphisme strict dans la cat\'egorie des
pr\'esch\'emas, et $X' \to X''$ est un monomorphisme (et m\^eme
un monomorphisme strict) dans la cat\'egorie des pr\'esch\'emas.
\end{corollaire}

La premi\`ere assertion signifie par d\'efinition que la suite de
morphismes
$$
\xymatrix@C=.7cm{X \times_{X''} X \ar@<2pt>[r]^-{\pr_1}
\ar@<-2pt>[r]_-{\pr_2} & X \ar[r] & X''}
$$
est exacte, et cela r\'esulte du fait que $X \to X''$ est fini et
fid\`element plat, comme on voit facilement (\cf Grothendieck, \loccit). L'assertion duale pour $X'' \to X'$ est encore plus triviale.

Le corollaire~\Ref{V.3.6} nous sera utile pour la th\'eorie du
groupe fondamental au \No suivant; il est possible (pour ceux
qui n'aiment pas la notion d'\'epimorphisme strict) de remplacer le
corollaire~\Ref{V.3.6} par telle variante que le lecteur arrangera
\`a son go\^ut personnel. Profitons seulement de l'occasion pour
signaler qu'une factorisation
\marginpar{118}
$f=f'f''$, avec $f''$ un \'epimorphisme strict et $f'$ un
monomorphisme, est n\'ecessairement unique \`a isomorphisme unique
pr\`es (dans toute cat\'egorie); cependant, il peut exister en
m\^eme temps une factorisation $f=f_1f_2$ ayant les
propri\'et\'es duales: $f_2$ est un \'epimorphisme, $f_1$ un
monomorphisme strict, (\'egalement unique \`a isomorphisme unique
pr\`es), qui ne soit pas isomorphe \`a la pr\'ec\'edente: il
suffit de prendre par exemple la cat\'egorie des espaces vectoriels
topologiques (s\'epar\'es, si on y tient), et pour $u \colon X \to
X'$ un morphisme tel que $u(X)$ ne soit pas ferm\'e.

\begin{proposition}
\label{V.3.7} Soient $Y$ un pr\'esch\'ema \emph{connexe} localement
noeth\'erien, $y$ un point de~$Y$, $\mathit{\Omega}$ une extension
alg\'ebriquement close
\ifthenelse{\boolean{orig}}
{de~$\kappa(y)$}
{de~$\kres(y)$}. Pour tout $X$ sur~$Y$, on
d\'esigne par $X(\Omega)$ l'ensemble des points de~$X$
\`a valeurs dans $\Omega$. Soient $X$, $X'$ des rev\^etements
\'etales de~$Y$ et $u \colon X \to X'$ un $Y$-morphisme tel que
l'application correspondante $X(\Omega) \to
X'(\Omega)$ soit un isomorphisme. Alors $u$ est un
isomorphisme.
\end{proposition}

On est imm\'ediatement ramen\'e au cas o\`u $X'$ est
connexe. Comme $X \to X'$ est fini et \'etale, on sait que le nombre
g\'eom\'etrique de points dans une fibre de~$X \to X'$ est
constant, et \'egal \`a 1 si et seulement si le morphisme
consid\'er\'e est un isomorphisme. Or ce nombre est aussi le
nombre d'\'el\'ements dans une fibre de~$X(\Omega) \to
X'(\Omega)$, d'o\`u la conclusion.
\ifthenelse{\boolean{orig}}{}
{\enlargethispage{.5cm}}%

\section{Conditions axiomatiques d'une th\'eorie de Galois}
\label{V.4}
Soit $\cal{C}$ une cat\'egorie, $F$ un foncteur covariant de
$\cal{C}$ dans la cat\'egorie des ensembles finis. Supposons les
conditions suivantes satisfaites:

\begin{enumerateb}

\item[(G 1)] $\cal{C}$ a un objet final\footnote{Rappelons qu'un objet
$e$ de~$\cal{C}$ est appel\'e \emph{objet final} si pour tout $X$ dans
$\cal{C}$, $\Hom(X,e)$ a exactement un \'el\'ement. On d\'efinit
de fa\c con duale un \emph{objet initial} de~$\cal{C}$.} et le produit
fibr\'e de deux objets au-dessus d'un troisi\`eme dans $\cal{C}$
existe (cet axiome peut aussi s'\'enoncer en disant que dans
$\cal{C}$ les limites projectives finies existent).

\item[(G 2)] Les sommes finies dans $\cal{C}$ existent (donc aussi un
objet initial $\emptyset_{\cal{C}}$, jouant le r\^ole de l'ensemble
vide), ainsi que le quotient d'un objet de~$\cal{C}$ par un groupe
fini d'automorphismes.

\item[(G 3)] Soit $u \colon X \to Y$ un morphisme dans $\cal{C}$,
alors $u$ se factorise en un produit $X\lto{u'} Y'\lto{u''} Y$, avec $u'$ un \'epimorphisme \emph{strict} et
$u''$ un monomorphisme, qui est un isomorphisme sur un sommande direct
de~$Y$.

\item[(G 4)]
\marginpar{119}
Le foncteur $F$ est exact \`a gauche (\ie transforme unit\'e
\`a droite en unit\'e \`a droite, et commute aux produits
\ifthenelse{\boolean{orig}}{fibr\'ee}{fibr\'es}).

\item[(G 5)] $F$ commute aux sommes directes finies, transforme
\'epimorphismes stricts en \'epi\-mor\-phismes, et commute au
passage au quotient par un groupe fini d'automorphismes.

\item[(G 6)] Soit $u \colon X \to Y$ un morphisme dans $\cal{C}$ tel
que $F(u)$ soit un isomorphisme, alors $u$ est un isomorphisme.
\end{enumerateb}

Notre objet est construire un groupe topologique~$\pi$, limite
projective de groupes finis, et une \'equivalence de la
cat\'egorie $\cal{C}$ avec la cat\'egorie $\cal{C}(\pi)$
\label{indnot:eb}\oldindexnot{$\cal{C}(\pi)$|hyperpage}%
\emph{des ensembles finis o\`u} $\pi$ \emph{op\`ere
contin\^ument} (\ie de telle fa\c con que le stabilisateur d'un
point soit un sous-groupe ouvert, ou encore qu'il existe un groupe
quotient discret qui op\`ere d\'ej\`a sur l'ensemble
envisag\'e), l'\'equivalence construite transformant le foncteur
donn\'e $F$ en le foncteur d'inclusion \'evident de~$\cal{C}(\pi)$
dans la cat\'egorie des ensembles finis. On notera tout de suite que
la cat\'egorie $\cal{C}(\pi)$ construite \`a l'aide d'un groupe
topologique~$\pi$, et le foncteur d'inclusion pr\'ec\'edent,
satisfont bien aux conditions (G 1) \`a (G 6).

Nous proc\'edons en plusieurs \'etapes.
\begin{enumerateb}
\item[a)] \emph{Soit} $u \colon X \to Y$ \emph{dans}
$\cal{C}$. \emph{Pour que} $u$ \emph{soit un monomorphisme, il faut et
il suffit} $F(u)$ \emph{le soit}. (Utilise (G 1), (G 4), (G 6)).

En effet, dire que $u$ est un monomorphisme signifie que la projection
$\pr_1 \colon X \times_Y X \to X$ est un isomorphisme.

\item[b)] \emph{Tout objet} $X$ \emph{de} $\cal{C}$ \emph{est
artinien}.

En effet, si $X' \to X'' \to X$ sont des monomorphismes tels que
$F(X')$ et $F(X'')$ aient m\^eme image dans $F(X)$, alors par a)
$F(X') \to F(X'')$ est un isomorphisme, donc $X' \to X''$ est un
isomorphisme par (G 6).

\item[c)] \emph{Le foncteur} $F$ \emph{est strictement
pro-repr\'esentable}. (\cf Grothendieck, Technique de descente et
th\'eor\`emes d'existence en G\'eom\'etrie Alg\'ebrique, II,
S\'eminaire Bourbaki 195, f\'evrier 1960).

En effet, d'apr\`es \loccit prop\ptbl \Ref{V.3.1}, cela r\'esulte de
b) et (G 4). On peut donc trouver un syst\`eme projectif sur~$I$
ordonn\'e filtrant:
$$
P=(P_i)_{i \in I}
$$
\marginpar{120}%
dans $\cal{C}$, consid\'er\'e comme pro-objet de~$\cal{C}$, et un
isomorphisme fonctoriel
\begin{equation*}
\label{eq:V.4.*} \tag{$*$} F(X)=\Hom_{\Pro(\cal{C})}(P,X)=\varinjlim_{i} \Hom_{\cal{C}}(P_i,X)
\end{equation*}
Cet isomorphisme est r\'ealis\'e par un \'el\'ement
$$
\varphi \in \varprojlim_{i} F(P_i)=F(P)
$$
On peut supposer de plus que les homomorphismes de transition
$\varphi_{ji} \colon P_i \to P_j$ ($i \geq j$) sont des
\emph{\'epimorphismes}, et que \emph{tout \'epimorphisme} $P_i \to
P'$ soit \'equivalent \`a un \'epimorphisme $P_i \to P_j$ pour
$j \leq i$ convenable (ce qui d\'etermine le syst\`eme projectif~$P$ de fa\c con essentiellement unique).

Un objet $X \in \cal{C}$ est dit \emph{connexe}
\index{connexe (objet d'une cat\'egorie)|hyperpage}%
s'il n'est pas isomorphe \`a la somme de deux objets de~$\cal{C}$
non isomorphes \`a l'objet initial $\emptyset_{\cal{C}}$.

\item[d)] \emph{Les} $P_i$ \emph{sont connexes et non isomorphes
\`a} $\emptyset_{\cal{C}}$.

Si $X$ est une unit\'e \`a gauche, on a $F(X)=\emptyset$ par (G 5)
appliqu\'e \`a la somme directe d'une famille vide, et
r\'eciproquement par (G 6). Donc si $X'$ est un objet de~$\cal{C}$
qui n'est pas unit\'e \`a gauche, \ie tel que $F(X') \neq
\emptyset$, il n'existe aucun morphisme de~$X'$ dans $X$. Donc si un
$P_i$ est unit\'e \`a gauche, alors $i$ est un plus grand
\'el\'ement de l'ensemble d'indices ordonn\'e filtrant $I$, et
la formule \eqref{eq:V.4.*} signifierait $F(X)=\Hom(P_i,X)={}$\kern2pt ensemble
r\'eduit \`a un \'el\'ement pour tout $X$, ce qui est absurde
puisque $F(\emptyset_{\cal{C}})=\emptyset$. Donc les~$P_i$ sont non
isomorphes \`a $\emptyset_{\cal{C}}$.

Supposons qu'on ait $P_i=A \amalg B$, d'o\`u par (G 5)
$F(P_i)=F(A) \amalg F(B)$, en particulier l'\'el\'ement $a_i$
de~$F(P_i)$, correspondant par \eqref{eq:V.4.*} \`a
l'homomorphisme identique $P_i \to P_i$, est dans $F(A) \amalg
F(B)$, par exemple dans $F(A)$. Cela signifie qu'il existe un $j
\geq i$ tel que $\varphi_{ij} \colon P_j \to P_i$ se factorise en
$P_j \to A \to P_i=A \amalg B$, o\`u la deuxi\`eme fl\`eche
est le morphisme canonique. Donc $F(P_j) \to F(P_i)$ se factorise
en $F(P_j) \to F(A) \to F(P_i)=F(A) \amalg F(B)$, et comme $F(P_j)
\to F(P_i)$ est surjectif par (G 5), il s'ensuit que
$F(B)=\emptyset$, donc $B$ est isomorphe \`a $\emptyset_{\cal{C}}$.

\item[e)]
\marginpar{121}
\emph{Tout morphisme} $u \colon X \to Y$ \emph{dans} $\cal{C}$,
\emph{avec} $X$ \emph{non isomorphe \`a} $\emptyset_{\cal{C}}$
\emph{et} $Y$ \emph{connexe, est un \'epimorphisme strict. Tout
endomorphisme d'un objet connexe est un automorphisme}.

Consid\'erons la factorisation (G 3) de~$u$, comme $X \neq
\emptyset_{\cal{C}}$ il r\'esulte de (G 6) que $F(X) \neq \emptyset$ donc
$F(Y') \neq \emptyset$ donc $Y' \neq \emptyset_{\cal{C}}$, donc $Y$
\'etant connexe, $Y'$ s'identifie \`a~$Y$, donc $u$ est un
\'epimorphisme strict. Supposons que $u$ soit un endomorphisme de
l'objet connexe $X$, prouvons que c'est un automorphisme. En effet, on
peut supposer~$X$ non isomorphe \`a $\emptyset_{\cal{C}}$, donc $u$ est
un \'epimorphisme strict par ce qui pr\'ec\`ede, donc $F(u)$ est
un \'epimorphisme par (G 5), et comme $F(X)$ est un ensemble fini,
$F(u)$ est bijectif. Donc $u$ est un automorphisme par (G 6).

En particulier, \emph{tout endomorphisme d'un} $P_i$ \emph{est un
automorphisme}.

\item[f)] \emph{Les conditions suivantes sur un} $P_i$ \emph{sont
\'equivalentes}: (i) L'application injective naturelle
$\Hom(P_i,P_i) \to \Hom(P,P_i) \simeq F(P_i)$ est aussi surjective,
\ie pour tout $u \colon P \to P_i$ il existe un $v \colon P_i \to
P_i$ tel que $u=v\varphi_i$ (o\`u $\varphi_i$ est l'homomorphisme
canonique $P \to P_i$). (ii) Le groupe $\Aut(P_i)$ op\`ere de fa\c con transitive sur~$F(P_i)$. (iii)~Le groupe $\Aut(P_i)$ op\`ere de
fa\c con simplement transitive sur~$F(P_i)$.

En effet, identifiant $\Hom(P,P_i)$ \`a $F(P_i)$, l'application
envisag\'ee dans (i) n'est autre que $v \mto
F(v)(\varphi_i)$. L'\'equivalence des trois conditions provient
alors du fait que $\Hom(P_i,P_i)=\Aut(P_i)$ et que l'application
pr\'ec\'edente est d\'ej\`a injective.

Un $P_i$ satisfaisant les conditions \'equivalentes
\ifthenelse{\boolean{orig}}
{(i) (ii) (iii)}
{(i), (ii), (iii)}
de~f) est appel\'e \emph{galoisien}.
\index{galoisien (objet)|hyperpage}%

\item[g)] \emph{Pour tout} $X$ \emph{dans} $\cal{C}$, \emph{il existe
un} $P_i$ galoisien \emph{tel que tout} $u \in \Hom(P,X)$ \emph{se
factorise en} $P \lto{\varphi_i} P_i \longrightarrow X$.

Soit $J=\Hom(P,X)=F(X)$, c'est un ensemble fini, donc il existe un
$P_j$ tel que tout $u \colon P \to X$ se factorise en $P\lto{j}P_j\longrightarrow X$, ou encore tel que le morphisme naturel
$$
P \to X^J \qquad (J=\Hom(P,X))
$$
se factorise en
$$
P \lto{\varphi_j} P_j \to X^J
$$
En
\marginpar{122}
vertu de (G~3), le morphisme $P_j \to X^J$ se factorise en le produit
d'un monomorphisme par un \'epimorphisme strict, que l'on peut
prendre de la forme $\varphi_{ij} \colon P_j \to P_i$. On est donc
ramen\'e \`a prouver que $P_i$ est galoisien. Soit $k$ un indice
$\geq j$ tel que tout morphisme $P \to P_i$ se factorise par
$P \lto{\varphi_k} P_k \longrightarrow P_i$. Notons que le
morphisme naturel $P_k \to X^J$ se factorise encore en le compos\'e
$$
P_k \lto{\varphi_{ik}} P_i \lto{U} X^J
$$
o\`u la premi\`ere \ifthenelse{\boolean{orig}}{fl\^eche}{fl\`eche} est un \'epimorphisme strict
par~e), et la deuxi\`eme un monomorphisme. On veut prouver que pour
un morphisme donn\'e $\psi \colon P_k \to P_i$, il existe un
endomorphisme $v$ de~$P_i$ tel que $\psi = v \varphi_{ik}$. Mais pour
tout $u \in \Hom(P_i,X)$, consid\'erons $u \psi \in \Hom(P_k,X)$,
il est donc de la forme $u'\varphi_{ik}$ avec un $u' \in \Hom(P_i,X)$
bien d\'etermin\'e. L'application $u \mto u'$ de~$J$ dans $J$
ainsi d\'efinie par $\psi$ est d'ailleurs injective car $\psi$ est
un \'epimorphisme en vertu de e); elle est donc bijective puisque
l'ensemble~$J$ est fini. L'application bijective $u \mto u'$ de~$J$
dans $J$ d\'efinit donc un isomorphisme $\alpha \colon X^J \isomto
X^J$ rendant commutatif le diagramme
$$
\xymatrix{P_k \ar[r]^-{\varphi_{ik}} \ar@{=}[d] & P_i \ar[r]^-U & X^J
\ar[d]^-{\alpha}_-{\simeq} \\ P_k \ar[r]^{\psi} & P_i \ar[r]^-U & X^J}
$$
D'apr\`es les propri\'et\'es d'unicit\'e de la factorisation
d'un morphisme en produit d'un monomorphisme par un \'epimorphisme
strict, il s'ensuit (puisque $\psi$ lui aussi est un \'epimorphisme
strict par e)) que l'on peut trouver un morphisme $v \colon P_i \to
P_i$ qui laisse le diagramme commutatif, cqfd.

On en conclut en particulier que \emph{les} $P_i$ \emph{galoisiens
forment un syst\`eme cofinal dans le syst\`eme des} $(P_j)$. On
aura donc, puisque pour un objet galoisien $P_i$ on a
$$
\Hom(P,P_i)=\Hom(P_i,P_i)=\Aut(P_i),
$$
par passage \`a la limite:
$$
\Hom(P,P) = \varprojlim_{i} \Hom(P,P_i) = \varprojlim_{i} \Hom(P_i,P_i) =
\varprojlim_{i} \Aut(P_i)
$$
o\`u
\ifthenelse{\boolean{orig}}
{\ignorespaces}
{\enlargethispage{.5cm}}
\marginpar{123}
la limite projective est prise sur les $P_i$ galoisiens. D'ailleurs,
moyennant l'identification $\Hom(P,P_i)=F(P_i)$, et compte tenu que
$F$ transforme \'epimorphismes en \'epimorphismes, on voit que les
homomorphismes de transition dans le syst\`eme projectif
pr\'ec\'edent sont surjectifs. On conclut de tout ceci:

\item[h)] \emph{On a}
$$
\Hom(P,P) = \Aut(P) = \varprojlim_{i} F(P_i) = \varprojlim_{i} \Aut(P_i),
$$
\emph{o\`u la limite projective est prise sur les} $P_i$
\emph{galoisiens}.

En particulier, $\Aut(P)$ appara\^it comme limite projective d'un
syst\`eme projectif de groupes finis (les homomorphismes de
transition \'etant surjectifs), on le munira de la topologie limite
projective des topologies discr\`etes. \emph{On d\'esignera
par}~$\pi$
\label{indnot:ec}\oldindexnot{$\pi$|hyperpage}%
\emph{et on appellera groupe fondamental}
\index{groupe fondamental d'une cat\'egorie galoisienne|hyperpage}%
(de~$\cal{C}$ muni de~$F$) \emph{le groupe oppos\'e
\`a}~$\Aut(P)$. Ce groupe op\`ere donc \emph{\`a droite}
sur~$P$, c'est la limite projective de groupes finis $\pi_i$
op\'erant \`a droite sur les $P_i$ galoisiens, $\pi_i$ \'etant
le groupe oppos\'e \`a $\Aut(P_i)$.

Compte tenu de l'isomorphisme fonctoriel
$$
F(X)=\Hom(P,X)
$$
et de la d\'efinition de~$\pi$, on voit donc que $\pi$ op\`ere
\emph{\`a gauche} sur~$F(X)$, et d'ailleurs de fa\c con continue
d'apr\`es g) (car avec les notations de g), c'est en fait $\pi_i$
qui op\`ere sur~$F(X)$). Il est trivial que pour tout morphisme $u
\colon X \to Y$ dans $\cal{C}$, le morphisme $F(u) \colon F(X) \to
F(Y)$ est compatible avec les op\'erations de~$\pi$. \emph{On peut
donc consid\'erer par la suite} $F$ \emph{comme un foncteur
covariant}
$$
F \colon \cal{C} \to \cal{C}(\pi)
$$
\emph{o\`u} $\cal{C}(\pi)$ \emph{est la cat\'egorie des ensembles
finis o\`u} $\pi$ \emph{op\`ere \`a gauche contin\^ument}.

Nous allons maintenant d\'efinir un foncteur en sens inverse:
$$
G \colon \cal{C} \from \cal{C}(\pi)
$$
par la formule
$$
G(E)=P \times_{\pi} E,
$$
o\`u $P \times_{\pi} E$ est d\'efini comme solution du
probl\`eme universel r\'esum\'e par
$$
\Hom_{\cal{C}}(P \times_{\pi} E,X) \isomto \Hom_{\pi}(E,\Hom(P,X))
$$
(o\`u
\marginpar{124}
dans le deuxi\`eme membre $\Hom(P,X)=F(X)$ est consid\'er\'e
comme un ensemble o\`u $\pi$ op\`ere \`a gauche). Il faut prouver
l'existence de l'objet $P \times_{\pi} E$.

\item[i)] \emph{Soit} $Q$ \emph{un objet de} $\cal{C}$ \emph{o\`u un
groupe fini} $G$ \emph{op\`ere \`a droite, et} $E$ \emph{un
ensemble fini o\`u} $G$ \emph{op\`ere \`a gauche. Alors} $G
\times_G E$ \emph{existe, et l'application canonique}
$$
F(Q) \times_G E \to F(Q \times_G E)
$$
\emph{est un isomorphisme}.

Comme les sommes directes finies existent dans $\cal{C}$ par (G 2), et
que $F$ y commute par (G 5), on est ramen\'e aussit\^ot au cas
o\`u $G$ op\`ere transitivement sur~$E$, car si les~$E_j$ sont les
trajectoires de~$G$ dans $E$, on aura
$$
Q \times_G E = \underset {j}{\amalg}\ Q \times_G E_j.
$$
Soit alors $a \in E$, soit $H$ son stabilisateur, on voit aussit\^ot
sur la d\'efinition que $Q \times_G E$ s'identifie \`a $Q/H$.
D'o\`u l'existence, gr\^ace \`a (G~2), et la propri\'et\'e
de commutation pour~$F$ gr\^ace \`a (G~5).

\item[j)] \emph{Soit} $E$ \emph{un objet de} $\cal{C}(\pi)$, \emph{et
soit} $P_i$ \emph{galoisien tel que} $\pi_i$ \emph{op\`ere
d\'ej\`a sur} $E$. \emph{Alors} $P_i \times_{\pi_i} E$
\emph{existe et on a un isomorphisme canonique}
$$
E \isomto F(P_i \times_{\pi_i} E)
$$
\emph{Si} $j \geq i$ \emph{est tel que} $P_j$ \emph{soit galoisien,
alors l'homomorphisme canonique} $P_j \times_{\pi_j} E \to\allowbreak P_i
\times_{\pi_i} E$ \emph{est un isomorphisme}.

La premi\`ere assertion est un cas particulier de i), compte tenu
que $\pi_i$ op\`ere de fa\c con simplement transitive sur~$F(P_i)$
qui est muni d'un point marqu\'e $\varphi_i$, d'o\`u un isomorphisme
$F(P_i) \times_{\pi_i} E \simeq E$. Pour la deuxi\`eme assertion on
utilise par exemple (G 6).

Soit, pout tout $j$, $\cal{C}_j$ la sous-cat\'egorie pleine de
$\cal{C}$ form\'ee des $X$ tels que $\Hom(P_j,X) \to \Hom(P,X)
\simeq F(X)$ soit bijectif. On sait par g) que $\cal{C}$ \emph{est
r\'eunion filtrante des} $\cal{C}_j$. On a donc pour $X \in
\cal{C}_j$:
\begin{align*}
\Hom_{\pi}(E,\Hom(P,X)) & \simeq \Hom_{\pi}(E,\Hom(P_j,X)) \simeq
\Hom_{\pi_j}(E,\Hom(P_j,X))\\ & \simeq \Hom(P_j \times_{\pi_j} E,X)
\end{align*}
et
\marginpar{125}
compte tenu de la derni\`ere assertion dans j) on trouve un
isomorphisme fonctoriel en l'objet $X$ de~$\cal{C}_j$:
$$
\Hom_{\pi}(E,\Hom(P,X)) \simeq \Hom(P_i \times_{\pi_i} E,X)
$$
Comme cela est vrai pour tout $j$ et que ces isomorphismes fonctoriels,
pour $j$ variable, s'induisent mutuellement, on conclut:

\item[k)] \emph{Sous les conditions de} j), \emph{le morphisme
compos\'e des morphismes canoniques}
$$
E \to \Hom(P_i,P_i \times_{\pi_i} E) \to \Hom(P,P_i \times_{\pi_i} E)
$$
\emph{fait de} $P_i \times_{\pi_i} E$ \emph{une solution du
probl\`eme universel d\'efinissant} $P \times_{\pi} E$,
\ie \emph{ce dernier existe et on a un isomorphisme}
$$
P \times_{\pi} E \isomto P_i \times_{\pi_i} E
$$

Cela ach\`eve la construction du foncteur $G(E)$. On a d'autre part
un homomorphisme fonctoriel
$$
\alpha \colon \id_{\cal{C}(\pi)} \to FG
$$
\ie un homomorphisme fonctoriel en l'objet $E$ de~$\cal{C}(\pi)$:
$$
\alpha(E) \colon E \to FG(E) = F(P \times_{\pi} E)
$$
savoir le compos\'e des morphismes canoniques
$$
E \to F(P) \times_{\pi} E \to F(P \times_{\pi} E)
$$
(o\`u le premier provient du point marqu\'e $\varphi \in F(P)$).
Conjuguant j) et k), on trouve:

\item[l)] \emph{L'homomorphisme} $\alpha$ \emph{est un isomorphisme}

On d\'efinit de m\^eme un homomorphisme fonctoriel
$$
\beta \colon GF \to \id_{\cal{C}}
$$
\ie un homomorphisme fonctoriel en l'objet $X$ de~$\cal{C}$:
$$
\beta(X) \colon P \times_{\pi} F(X) \to X
$$
comme
\marginpar{126}
associ\'e au $\pi$-homomorphisme
$$
F(X) \to \Hom(P,X)
$$
inverse de l'isomorphisme canonique $\Hom(P,X) \isomto F(X)$.

\item[m)] \emph{Les compos\'es}
\begin{gather*}
F(X) \lto{\alpha(F(X))} FGF(X) \lto{F(\beta(X))} F(X)\\
G(E) \lto{G(\alpha(E))} GFG(E) \lto{\beta(G(E))} G(E)
\end{gather*}
\emph{sont les isomorphismes identiques}.

\^Ane qui trotte.

Compte tenu de l) il s'ensuit:

\item[n)] \emph{L'homomorphisme} $\beta$ \emph{est un isomorphisme}

Nous avons ainsi obtenu le r\'esultat promis:
\end{enumerateb}

\begin{theoreme}
\label{V.4.1}
Soit $\cal{C}$ une cat\'egorie satisfaisant les conditions \emph{(G
1)}, \emph{(G 2)}, \emph{(G~3)} du d\'ebut du num\'ero, et $F$ un
foncteur covariant de~$\cal{C}$ dans la cat\'egorie des ensembles
finis, satisfaisant les conditions \emph{(G 4)}, \emph{(G 5)} et
\emph{(G 6)}. Alors les constructions canoniques pr\'ec\'edentes
d\'efinissent des \'equivalences de cat\'egories $F\colon
\cal{C}\to \cal{C}(\pi)$ et $G\colon \cal{C}(\pi)\to \cal{C}$
quasi-inverses l'une de l'autre. De fa\c con pr\'ecise, il existe
un pro-objet~$P$ de~$\cal{C}$, et un isomorphisme fonctoriel $F(X)
\isomfrom \Hom(P, X)$, $\pi$ est le groupe oppos\'e au groupe des
automorphismes de~$P$, topologis\'e de fa\c con convenable, de
fa\c con que~$\pi$ op\`ere de fa\c con continue sur les
ensembles $\Hom(P, X)\simeq F(X)$. Enfin $G$ est donn\'e par
$G(E)\simeq P\times_{\pi} E$.
\end{theoreme}

\begin{remarques}
\label{V.4.2}
L'\'enonc\'e des conditions (G 1) \`a (G 6) devient plus simple
et plus sympathique si on remplace (G 2) et (G 5) respectivement par:
\begin{enumerateb}
\item[(G' 2)] Les limites inductives finies dans
\ifthenelse{\boolean{orig}}{$C$}{$\cal{C}$}
existent.
\item[(G' 5)] Le foncteur $F$ est exact \`a droite (\ie commute aux
limites inductives finies).
\end{enumerateb}
Ces conditions sont en apparence plus fortes que (G~2) et (G~5), mais
il r\'esulte aussit\^ot du th\'eor\`eme de
structure~\Ref{V.4.1} qu'elles sont entra\^in\'ees par (G~1)
\`a (G~6). On notera cependant que dans les cas qui nous
int\'eresseront, la v\'erification de (G~2) et (G~5) semble
effectivement plus simple que celle de (G'~2) et (G'~5). J'ignore si
dans la condition (G~3), le fait que $u''$ soit un isomorphisme sur un
sommande direct de~$Y$ ne pourrait \^etre omis.
\end{remarques}

\section{Cat\'egories galoisiennes}
\label{V.5}
\marginpar{127}

\begin{definition}
\label{V.5.1}
On appelle cat\'egorie galoisienne
\index{galoisienne (cat\'egorie)|hyperpage}%
une cat\'egorie $\cal{C}$ \'equivalente \`a une cat\'egorie
$\cal{C}(\pi)$, o\`u $\pi$ est un groupe compact, limite projective
de groupes finis (\ie totalement disconnexe).

Pour la d\'efinition de~$\cal{C}(\pi)$, \cf d\'ebut du \No \Ref{V.4}. En
vertu du th\ptbl \Ref{V.4.1}, $\cal{C}$ est galoisienne si et seulement si
elle satisfait les conditions (G 1) \`a (G 3), et s'il existe un
foncteur $F$ de
\ifthenelse{\boolean{orig}}{$C$}{$\cal{C}$}
dans la cat\'egorie des ensembles finis satisfaisant les conditions
(G~4) \`a (G~6) (\ie qui est \emph{exact} et \emph{conservatif},
dans une terminologie g\'en\'erale). Un tel foncteur sera
appel\'e \emph{foncteur fondamental}
\index{foncteur fondamental|hyperpage}%
de la cat\'egorie galoisienne $\cal{C}$\kern1pt\footnote{Il semble pr\'ef\'erable d'adopter le terme plus parlant
de \og foncteur fibre\fg.};
\index{foncteur fibre|hyperpage}%
il est pro-repr\'esentable par un pro-objet que nous noterons~$P_F$;
un pro-objet $P$ tel que le foncteur $F$ associ\'e soit fondamental
est appel\'e \emph{pro-objet fondamental}.
\index{pro-objet fondamental|hyperpage}%
De cette fa\c con, la cat\'egorie des foncteurs fondamentaux sur~$\cal{C}$ est anti-\'equivalente \`a la cat\'egorie des
pro-objets fondamentaux; si $F$ et $P$ se correspondent, le groupe
$\Aut F$ est donc isomorphe \`a l'oppos\'e du groupe $\Aut P$,
donc le groupe not\'e $\pi$ dans le num\'ero pr\'ec\'edent
n'est autre que $\Aut P$. Rappelons qu'au num\'ero pr\'ec\'edent
nous avons construit, \`a partir d'un foncteur fondamental
\emph{donn\'e} $F$, une \'equivalence de~$\cal{C}$ avec
$\cal{C}(\pi)$ (o\`u $\pi=\Aut(F)$) qui transforme le foncteur $F$
en le foncteur canonique de~$\cal{C}(\pi)_x$ dans la cat\'egorie des
ensembles finis. Dans ce cas type $\cal{C}=\cal{C}(\pi)$, $F=$
foncteur canonique, le pro-objet fondamental associ\'e \`a $F$
n'est autre que le syst\`eme projectif des quotients discrets
$\pi_i$ de~$\pi$.
\end{definition}

Il peut \^etre utile d'expliciter la cat\'egorie des pro-objets de
$\cal{C}(\pi)$. On trouve:

\begin{proposition}
\label{V.5.2}
La cat\'egorie $\Prodash\cal{C}(\pi)$
\label{indnot:ed}\oldindexnot{$\Prodash\cal{C}(\pi)$|hyperpage}%
est canoniquement \'equivalente \`a la cat\'egorie
$\cal{C}'(\pi)$ des espaces, \`a groupe topologique $\pi$
d'op\'erateurs, qui sont compacts et totalement disconnexes.
\end{proposition}

Comme cette derni\`ere contient $\cal{C}(\pi)$ comme
sous-cat\'egorie pleine (correspondant aux espaces \`a
op\'erateurs compacts discrets) et que les limites projectives y
existent, on a en tous cas un foncteur canonique $g\colon
\Prodash\cal{C}(\pi) \to \cal{C}'(\pi)$, au syst\`eme projectif
$Q=(Q_i)$ correspondant
l'objet $X=\varprojlim_{i} Q_i$ de~$\cal{C}'(\pi)$. Pour d\'efinir
un foncteur en sens inverse, il suffit de d\'efinir un foncteur
contravariant de~$\cal{C}'(\pi)$ dans la cat\'egorie des foncteurs
$\cal{C}\to \Ens$ exacts \`a gauche, et on prendra pour tout
$X\in\cal{C}'(\pi)$ le foncteur $h(X)(E)= \Hom(X, E)$
\marginpar{128}
(le $\Hom$ \'etant pris dans $\cal{C}'(\pi)$). Il est imm\'ediat
par d\'efinition que les foncteurs $h$ et $g$ sont adjoints l'un de
l'autre, et que $hg$ est canoniquement isomorphe au foncteur identique
de~$\Prodash\cal{C}(\pi)$. Il reste \`a prouver (pour \'etablir
que $g$ et $h$ sont quasi-inverses l'un de l'autre) que tout objet de
$\cal{C}'(\pi)$ est isomorphe \`a un objet de la forme $g(Q)$,
o\`u $Q\in\Prodash{\cal C}(\pi)$, en d'autres termes: \emph{tout
espace $X$ \`a groupe topologique $\pi$ d'op\'erateurs, qui est
compact et totalement disconnexe, est isomorphe \`a une limite
projective d'espaces \`a op\'erateurs finis discrets}. Comme $X$
est limite projective de ses quotients finis discrets (en tant
qu'espace topologique sans op\'erateurs), on est ramen\'e \`a
montrer que dans l'ensemble de ces quotients, il y a un syst\`eme
cofinal qui est invariant par $\pi$. Il suffit pour cela de montrer
que pour un tel quotient $X'$, l'ensemble des transform\'es de ce
quotient par les op\'erations de~$\pi$ est fini, (on prendra alors
le sup desdits transform\'es, qui sera un quotient invariant
majorant $X'$). Ou encore, qu'il y a un sous-groupe invariant ouvert
$\pi'$ de~$\pi$ tel que les \'el\'ements de ce sous-groupe
invarient $X'$. Or $X'$ correspond \`a une partition finie de~$X$ en
ensembles ouverts $X_i$. Par raison de continuit\'e et de
compacit\'e de~$\pi$, il existe un voisinage $V$ de
l'\'el\'ement neutre de~$\pi$ tel que $s\in V$ implique
$s\cdot X_i\subset X_i$ pour tout~$i$, donc $s$ invarie~$X'$. Or on sait
que les sous-groupes invariants ouverts de~$\pi$ forment un
syst\`eme fondamental de voisinages de l'\'el\'ement neutre, ce
qui ach\`eve la d\'emonstration.

Remarquons qu'on voit encore plus simplement que la cat\'egorie
$\Ind \cal{C}(\pi)$ est canoniquement \'equivalente \`a la
cat\'egorie des ensembles o\`u $\pi$ op\`ere
contin\^ument. Nous n'en aurons pas besoin ici.

\begin{proposition}
\label{V.5.3}
Soient $\cal{C}$ une cat\'egorie galoisienne, $F$ un foncteur
fondamental
\ifthenelse{\boolean{orig}}{sur~$C$,}{sur~$\cal{C}$,}
$P=(P_i)$ le pro-objet associ\'e, normalis\'e de la fa\c con
habituelle. Soit
\ifthenelse{\boolean{orig}}{$X\in C$;}{$X\in\cal{C}$;}
pour que $X$ soit connexe, il faut et il suffit que $\pi$ op\`ere
transitivement sur~$E=F(X)$.
\end{proposition}

On est ramen\'e au cas type $\cal{C}=\cal{C}(\pi)$, $F=$ foncteur
canonique, o\`u c'est trivial.

\begin{corollaire}
\label{V.5.4}
Conditions \'equivalentes sur~$X$: \textup{(i)} $X$ est connexe et
$\not\simeq\emptyset_{\cal{C}}$ \textup{(ii)} le groupe $\pi$ est
transitif sur~$E=F(X)$, et $F(X)\ne\emptyset$ \textup{(iii)} $X$ est
isomorphe \`a un $P_i$.
\end{corollaire}

L'\'equivalence
\marginpar{129}
de (i) et (iii) r\'esulte aussi d\'ej\`a facilement de
\No \Ref{V.4},~e).

\begin{proposition}
\label{V.5.5}
Soit $Q=(Q_i)_{i\in I}$ un pro-objet
\ifthenelse{\boolean{orig}}{de~$C$,}{de~$\cal{C}$,}
normalis\'e de la fa\c con habituelle, et soit $G$ le foncteur
correspondant $G(X)=\Hom(Q, X)$ de~$\cal{C}$
\ifthenelse{\boolean{orig}}{(dans~$\Ens$.}{(dans~$\Ens$).}
Les conditions suivantes sont \'equivalentes:
\begin{enumerate}
\item[(i)] $G$ commute aux sommes directes finies
\item[(ii)] $G$ commute \`a la somme de deux objets
\item[(iii)] Les $Q_i$ sont connexes et
\ifthenelse{\boolean{orig}}{$\not\simeq\emptyset_C$}
{$\not\simeq\emptyset_\cal{C}$}
\item[(iv)] $Q$ est isomorphe \`a $\pi/H$, o\`u $H$ est un
sous-groupe ferm\'e de~$\pi$.
\item[(v)] Le foncteur $G$ est isomorphe au foncteur $E\mto
E^H$ (ensemble des invariants par $H$) d\'efini par un sous-groupe
ferm\'e $H$ de~$\pi$.
\end{enumerate}
\end{proposition}

N.B. dans l'\'enonc\'e de (iv) et (v), on suppose choisi un
foncteur fondamental, permettant d'identifier $\cal{C}$ \`a la
cat\'egorie~$\cal{C}(\pi)$.

\subsubsection*{D\'emonstration} On peut supposer
$\cal{C}=\cal{C}(\pi)$. L'implication (i)$\To$(ii) est
triviale, (ii)$\To$(iii) se prouve comme la
propri\'et\'e d) du \No \Ref{V.4}. Prouvons (iii)$\To$(iv). En
effet, on sait que $\varprojlim_i Q_i$ est non vide comme limite
projective d'ensembles finis non vides. Soit $a$ un point de
$\varprojlim_i Q_i$, il d\'efinit un homomorphisme d'espaces \`a
op\'erateurs
$$
\pi\to Q
$$
qui est \emph{surjectif}, car pour tout $i$ le compos\'e $\pi\to
Q\to Q_i$ l'est, puisque $\pi$ est transitif dans $Q_i$ en vertu
de~\Ref{V.5.3}. Si donc $H$ est le sous-groupe de stabilisateur de
$a$, on obtient un isomorphisme $\pi/H\isomto Q$. Les implications
(iv)$\To$(v) et (v)$\To$(i) sont \`a nouveau
triviales.

\begin{proposition}
\label{V.5.6}
Soient $\cal{C}$ une cat\'egorie galoisienne, $P$ un pro-objet
fondamental de~$\cal{C}$ et $F$ le foncteur fondamental
associ\'e. Soit $P'=(P'_i)_{i\in I}$ un pro-objet de~$\cal{C}$, mis
sous forme normale, et $F'$ le foncteur associ\'e $F'(X)=\Hom(P',X)$
de~$\cal{C}$ dans~$\Ens$. Conditions \'equivalentes:
\begin{enumerate}
\item[(i)] $P'\simeq P$, ou encore $F'\simeq F$.
\item[(ii)] $P'$ est fondamental, ou encore $F'$ est
fondamental.
\item[(iii)] $F'$ transforme somme de deux objets en somme, et
$X\not\simeq\emptyset_{\cal{C}}$ implique $F(X)\ne\emptyset$.
\item[(iv)]
\marginpar{130}
Les objets de~$\cal{C}$ connexes et $\ne\emptyset_{\cal{C}}$ sont
exactement les objets isomorphes \`a un~$P'_i$.
\end{enumerate}
\end{proposition}

On a trivialement (i)$\To$(iii) et (i)$\To$(ii),
de plus (ii)$\To$(iv) en vertu de~\Ref{V.5.4} (appliqu\'e
\`a $P'$ au lieu de~$P$). De plus (iii) ou (iv) implique en vertu
de~\Ref{V.5.5} que $P'$ est de la forme $\pi/H$ o\`u $H$ est un
sous-groupe ferm\'e de~$\pi$. Dans le cas (iii), il existe pour tout
sous-groupe invariant ouvert $\pi'$ de~$\pi$ un $\pi$-homomorphisme
$P'=\pi/H\to \pi/\pi'$, d'o\`u $H\subset \pi'$, d'o\`u $H=(0)$ et
par suite (i), cqfd.

\begin{corollaire}
\label{V.5.7}
Soit $\cal{C}$ une cat\'egorie galoisienne. Les pro-objets
fondamentaux sont isomorphes, les foncteurs fondamentaux sont
isomorphes.
\end{corollaire}

En d'autres termes, \emph{la cat\'egorie des foncteurs fondamentaux
est un groupo\"ide connexe~$\Gamma$},
\label{indnot:ef}\oldindexnot{$\Gamma$|hyperpage}%
qu'on peut appeler le \emph{groupo\"ide fondamental}
\index{groupoide fondamental@groupo\"ide fondamental|hyperpage}%
de la cat\'egorie galoisienne~$\cal{C}$. Si $\cal{C}=\cal{C}(\pi)$,
le groupe des automorphismes d'un objet du groupo\"ide fondamental
est isomorphe \`a $\pi$, cet isomorphisme \'etant bien
d\'etermin\'e \`a automorphisme int\'erieur pr\`es. (N.B. on
appelle \emph{groupo\"ide} une cat\'egorie o\`u tous les
morphismes sont des isomorphismes, groupo\"ide \emph{connexe}
groupo\"ide dont tous les objets sont isomorphes). Les pro-objets
fondamentaux de~$\cal{C}$ forment un groupo\"ide connexe
\'equivalent \`a l'\emph{oppos\'e} du groupo\"ide
fondamental. Si $F, F'$ sont deux foncteurs fondamentaux, associ\'es
\`a des pro-objets fondamentaux $P, P'$, alors $\Hom(F, F')=\Isom(F,
F')$ est parfois not\'e $\pi_{F', F}$ et joue le r\^ole d'un
\og ensemble de \emph{classes de chemins} de~$F$ \`a~$F'$\fg. En
particulier $\pi_{F, F}=\pi_F$ n'est autre que le \emph{groupe
fondamental de~$\cal{C}$} en $F$ construit dans le num\'ero
pr\'ec\'edent. Quant au pro-objet $P$ associ\'e \`a~$F$, il
joue le r\^ole d'un \emph{rev\^etement universel en~$F$} de
l'objet final $e_{\cal{C}}$ de~$\cal{C}$.

Il peut \^etre commode d'avoir une description de~$\cal{C}$ (\`a
\'equivalence pr\`es) en termes de son groupo\"ide fondamental
$\Gamma$, sans passer par un choix d'un objet particulier $F$ de ce
dernier. Or \`a tout objet $X$ de~$\cal{C}$ est associ\'e le
foncteur $E_X$ sur le groupo\"ide fondamental, d\'efini par
$$
E_X(F)=F(X),
$$
\`a valeurs dans~$\Ens$. (Un tel foncteur est connu en topologie
sous le nom de \og syst\`eme local\fg sur le groupo\"ide)
$F(X)=E_X(F)$ peut \^etre appel\'e la \emph{fibre} de~$X$ en $F$,
et le foncteur $E_X$ le foncteur-fibre associ\'e \`a~$X$. Le
foncteur $E_X$ a la propri\'et\'e
\marginpar{131}
suivante:
\ifthenelse{\boolean{orig}}
{\emph{\underline{pour tout} $F$, $E_X(F)$
\underline{est un ensemble fini o\`u le groupe topologique}
$\pi_F=\Aut(F)$ \underline{op\`ere contin\^ument}}}
{\emph{pour tout $F$, $E_X(F)$ est un ensemble fini
o\`u le groupe topologique $\pi_F=\Aut(F)$ op\`ere contin\^ument}}.
Pour un
foncteur covariant donn\'e $\xi$ du groupo\"ide fondamental
dans~$\Ens$, la condition pr\'ec\'edente \'equivaut d'ailleurs
\`a la m\^eme condition pour \emph{un} $F$ fix\'e
quelconque. Ceci pos\'e:

\begin{proposition}
\label{V.5.8}
Le foncteur $X\mto E_X$ est une \'equivalence de la cat\'egorie
$\cal{C}$ avec la cat\'egorie des foncteurs covariants du
groupo\"ide fondamental $\Gamma$ de~$\cal{C}$ dans~$\Ens$, qui
satisfont la condition
\ifthenelse{\boolean{orig}}
{soulign\'ee plus haut.}
{mise en \'evidence ci-dessus.}
\end{proposition}

En effet, soit $F_0$ un objet du groupo\"ide fondamental, et soit
$\pi_0=\Aut(F_0)$, alors le foncteur $\xi\mto \xi(F_0)$ est une
\'equivalence de la deuxi\`eme cat\'egorie envisag\'ee
dans~\Ref{V.5.8} avec la cat\'egorie $\cal{C}(\pi_0)$, comme on
constate aussit\^ot. D'autre part, le compos\'e de ce dernier et
du foncteur $X\mto E_X$ est l'\'equivalence naturelle
$\cal{C}\to\cal{C}(\pi_0)$. Il en r\'esulte que le foncteur
$X\mto E_X$ lui-m\^eme est une \'equivalence.

\begin{corollaire}
\label{V.5.9}
La cat\'egorie $\Prodash\cal{C}$
est \'equivalente canoniquement \`a la cat\'egorie des foncteurs
covariants $\xi$ du groupo\"ide fondamental $\Gamma$ dans la
cat\'egorie des espaces topologiques, satisfaisant la condition:
pour tout objet $F$ de~$\Gamma$, $\xi(F)$ est un espace compact
totalement disconnexe \`a groupe topologique $\pi_F$
d'op\'erateurs.
\end{corollaire}

Ici encore, on peut v\'erifier cette condition sur~$\xi$, il suffit
de la v\'erifier pour \emph{un}~$F$. La d\'emonstration est la
m\^eme que pour~\Ref{V.5.8}.

\begin{remarque}
\label{V.5.10}
Soit $(F_s)_{s\in S}$ une famille
\ifthenelse{\boolean{orig}}{d'objet}{d'objets}
du groupo\"ide fondamental~$\Gamma$. Posons pour $s, s'\in S$:
$$
\Hom(s, s')=\Hom(F_s, F_{s'})
$$
de sorte que $S$ devient lui-m\^eme un groupo\"ide connexe et
l'application $s\mto F_s$ un foncteur pleinement fid\`ele de~$S$
dans $\Gamma$, soit~$f$. Consid\'erant alors le foncteur $X\mto
E_X\circ f$ de~$\cal{C}$ dans la cat\'egorie de foncteurs
$\Hom(S,\Ens)$, on obtient une variante de~\Ref{V.5.8}
(et~\Ref{V.5.9}) avec $\Gamma$ remplac\'e par~$S$. L'\'enonc\'e
ainsi obtenu se r\'eduit au th\'eor\`eme~\Ref{V.4.1} lorsque $S$
est r\'eduit \`a un point, et n'est autre que~\Ref{V.5.8}
lui-m\^eme si $S$ est l'ensemble des objets de~$\Gamma$.
\end{remarque}

Nous allons utiliser~\Ref{V.5.9} pour d\'efinir un pro-objet
canonique de~$\cal{C}$. Pour ceci, nous consid\'erons le foncteur
de~$\Gamma$ dans la cat\'egorie des espaces topologiques
\marginpar{132}
(et m\^eme des groupes topologiques):
$$
f\colon F\mto \Aut(F)=\pi_F.
$$
Ce foncteur satisfait la condition envisag\'ee dans~\Ref{V.5.8},
l'espace \`a op\'erateurs $f(F)$ sous $\pi_F$ n'est autre que
$\pi_F$, consid\'er\'e comme espace \`a op\'erateurs sous
lui-m\^eme par les automorphismes int\'erieurs. Donc le foncteur
$f$ correspond \`a un pro-objet de~$\cal{C}$ d\'etermin\'e \`a
isomorphisme unique pr\`es, qui est m\^eme un pro-groupe de
$\cal{C}$ et qui est appel\'e le \emph{pro-groupe fondamental
de~${\cal C}$},
\index{pro-groupe fondamental d'une cat\'egorie galoisienne|hyperpage}%
jouant le r\^ole d'un syst\`eme local de groupes
fondamentaux. C'est donc un pro-groupe $\Pi$ de~$\cal{C}$ d\'efini
par la condition qu'on ait un isomorphisme fonctoriel en $F$
$$
F(\Pi)\simeq\pi_F.
$$

Si $X$ est un pro-objet quelconque de~$\cal{C}$, on a un morphisme
canonique
$$
\Pi\times X \to X
$$
qui fait de~$X$ un objet \`a groupe d'op\'erateurs \`a gauche
$G$ dans $\Prodash\cal{C}$. Il suffit pour ceci de noter que pour $F$
variable, on a une application canonique
$$
\Pi(F)\times X(F) \to X(F)
$$
\ie $$
\Aut(F)\times E_X(F) \to E_X(F), \quad \text{ou} \quad \pi_F\times
F(X)\to F(X)
$$
qui est fonctorielle en~$F$. Elle est aussi fonctorielle en $X$, donc
pour tout morphisme $X\to Y$ de pro-objets, le diagramme
$$
\xymatrix{
\Pi \times X \ar[r]\ar[d]& X\ar[d] \\ \Pi \times Y \ar[r]& Y
}
$$
correspondant est commutatif.

\begin{remarque}
\label{V.5.11}
On se gardera de confondre un pro-objet fondamental $P$ (qui n'est pas
muni d'une structure de groupe, et est connexe) avec le pro-groupe
fondamental (qui est un pro-\emph{groupe}, et en g\'en\'eral non
connexe). De fa\c con pr\'ecise, $G$ est connexe si
\ifthenelse{\boolean{orig}}{\ignorespaces}{et}
seulement si $\pi_F$ op\'erant sur lui-m\^eme par automorphismes
int\'erieurs est transitif, \ie si $\pi$ est r\'eduit \`a
l'\'el\'ement neutre, ou encore $\cal{C}$ \'equivalente \`a la
cat\'egorie des ensembles finis. Une autre diff\'erence
\marginpar{133}
essentielle est que $G$ est d\'etermin\'e \`a isomorphisme
unique pr\`es, et $P$ n'est d\'etermin\'e qu'\`a isomorphisme
non unique pr\`es.

Soit $E$ un ensemble fini et consid\'erons le foncteur constant sur
le groupo\"ide~$\Gamma$, de valeur $E$: il d\'efinit en vertu
de~\Ref{V.5.8} un objet de~$\cal{C}$, not\'e $E_{\cal{C}}$, et qui
peut aussi s'interpr\'eter comme la somme de~$E$ exemplaires de
l'objet final $e_{\cal{C}}$ de~$\cal{C}$. On peut consid\'erer
$E_{\cal{C}}$ comme un foncteur en $E$, de la cat\'egorie des
ensembles finis dans la cat\'egorie $\cal{C}$, et ce foncteur est
\emph{exact}, donc transforme groupes finis en $\cal{C}$-groupes,
\ifthenelse{\boolean{orig}}{etc...}{etc.} Si donc $X$ est un objet de~$\cal{C}$ sur lequel le groupe fini
$G$ op\`ere \`a droite, on voit qu'on peut consid\'erer $X$
comme un objet de~$\cal{C}$ ayant un $\cal{C}$-groupe d'op\'erateurs
\`a droite~$G_{\cal{C}}$. On dira donc par extension d'une
terminologie g\'en\'erale relative \`a des objets \`a
$\cal{C}$-groupes d'op\'erateurs, que $X$ est \emph{formellement
principal homog\`ene}
\index{formellement principal homog\`ene (objet)|hyperpage}%
sous $G$ si $X$ est formellement principal homog\`ene sous
$G_{\cal{C}}$, \ie si le morphisme canonique
$$
X\times G_{\cal{C}} \to X\times X
$$
d\'eduit de l'op\'eration de~$G_{\cal{C}}$ sur~$X$ \`a droite,
est un isomorphisme. On dit que $X$ est \emph{principal homog\`ene}
\index{principal homog\`ene (objet)|hyperpage}%
sous $G$ s'il l'est sous $G_{\cal{C}}$, \ie s'il l'est formellement,
et si de plus le quotient $X/G=X/G_{\cal{C}}$ est~$e_{\cal{C}}$. Si on
se fixe un foncteur fondamental, d'o\`u une \'equivalence de
$\cal{C}$ avec une cat\'egorie $\cal{C}(\pi)$, $X$ correspond \`a
un ensemble sur lequel $\pi$ op\`ere \`a gauche contin\^ument,
soit $E=F(X)$. Faire op\'erer $G$ sur~$X$ \`a droite revient alors
\`a faire op\'erer $G$ sur l'ensemble $E$ \`a droite, de fa\c con que les op\'erations de~$G$ commutent \`a celles de~$\pi$. On
constate alors aussit\^ot que $X$ est principal homog\`ene sous
$G$ si et seulement si l'ensemble $E$ est un espace principal
homog\`ene
\ifthenelse{\boolean{orig}}
{sous $G$}
{sous $G$,}
\ie si et seulement si $G$ y op\`ere de fa\c con simplement transitive. (D'ailleurs $X$ est formellement
principal homog\`ene \sss
$E$ principal homog\`ene \emph{ou} vide). Comparant
avec~\Ref{V.5.3}, on voit que si $X$ est principal homog\`ene sous
$G$ \emph{et connexe}, alors l'homomorphisme donn\'e de~$G$ dans le
groupe oppos\'e \`a $\Aut(X)$ est un \emph{isomorphisme}; et
d'ailleurs pour qu'un objet $X$
\ifthenelse{\boolean{orig}}{de~$C$}{de~$\cal{C}$}
soit connexe et principal homog\`ene sous le groupe oppos\'e \`a
$\Aut(X)$, il faut et il suffit, avec les notations du \No \Ref{V.4}, qu'il
soit isomorphe \`a un $P_i$ \emph{galoisien}. Dans le cas type
$\cal{C}=\cal{C}(\pi)$, cela signifie que $X$ est isomorphe \`a un
quotient de~$\pi$ par un sous-groupe invariant.

Supposons toujours donn\'e un foncteur fondamental~$F$. Alors la
donn\'ee d'un
\marginpar{134}
$X$ principal homog\`ene sous un groupe fini $G$ op\'erant \`a
droite, et d'un point $a\in F(X)$, est \'equivalente \`a la
donn\'ee d'un homomorphisme de~$\pi$ dans le groupe~$G$. En effet,
\`a un tel homomorphisme on fait correspondre l'ensemble $E=G$, en y
faisant op\'erer~$\pi$ \`a gauche gr\^ace \`a l'homomorphisme
donn\'e $\pi\to G$ et les translations \`a gauche de~$G$, et en y
faisant op\'erer $G$ \`a droite par translation \`a droite, le
point marqu\'e $a$ de~$E$ \'etant l'\'el\'ement unit\'e
de~$G$. Gr\^ace \`a ce qui pr\'ec\`ede, on obtient ainsi de
fa\c con essentiellement unique tout triple $(X, G, a)$ ayant les
propri\'et\'es envisag\'ees plus haut,
\ifthenelse{\boolean{orig}}
{puisque un}
{puisqu'un}
ensemble \`a
point marqu\'e principal homog\`ene sous un groupe $G$ s'identifie
\`a ce dernier. De cette fa\c con, on a une interpr\'etation
g\'eom\'etrique directe du foncteur $G\mto \Hom(\pi, G)$ de la
cat\'egorie des groupes finis dans~$\Ens$, foncteur qui est
pro-repr\'esentable \`a l'aide de~$\pi$, et dont la
consid\'eration fournirait donc une autre construction du groupe
$\pi$ associ\'e \`a~$F$.
\end{remarque}

\section{Foncteurs exacts d'une cat\'egorie galoisienne dans une autre}
\label{V.6}

\begin{proposition}
\label{V.6.1}
Soient $\cal{C}$, $\cal{C'}$ deux cat\'egories galoisiennes,
$H\colon \cal{C}\to \cal{C'}$ un foncteur covariant, $F'$ un foncteur
fondamental sur~$\cal{C'}$ et $F=F' \circ H$. Conditions
\'equivalentes:
\begin{enumerate}
\item[(i)] $H$ est \emph{exact}, \ie exact \`a gauche et exact
\`a droite.
\item[(ii)] $H$ est exact \`a gauche, transforme sommes finies en
sommes finies, et \'epimorphismes en \'epimorphismes (ou encore:
objets $\not\approx \emptyset_{\cal{C}}$ en objets $\not\approx
\emptyset_{\cal{C'}}$).
\item[(iii)] $F$ est un foncteur fondamental sur~$\cal{C}$.
\end{enumerate}
\end{proposition}

L'implication (i)$\To$(ii) est un fait g\'en\'eral aux
cat\'egories. D'ailleurs la premi\`ere forme donn\'ee \`a (ii)
implique la seconde, comme on voit en notant que si $X$ est un objet
de~$\cal{C}$, alors $X$ est $\not\approx \emptyset_{\cal{C}}$ \sss le
morphisme $X\to e_{\cal{C}}$ est un \'epimorphisme; on notera que
$F$ \'etant suppos\'e exact \`a gauche transforme $e_{\cal{C}}$
en~$e_{\cal{C'}}$. La deuxi\`eme forme donn\'ee \`a (ii)
implique~(iii), car $F$ \'etant exact \`a gauche donc
pro-repr\'esentable est justiciable de~\Ref{V.5.6}
crit\`ere~(iii). Enfin (iii) implique~(i), comme il r\'esulte du
fait que $F$ est exact, et \og conservatif\fg (\ie satisfait l'axiome
(G~6) de~\no \Ref{V.4}).

Soit alors $\Gamma$ le groupo\"ide fondamental de~$\cal{C}$,
$\Gamma'$ celui de~$\cal{C'}$. Donc si
\marginpar{135}
$H$ est exact, alors $F' \mto F' \circ H$ est un foncteur du
groupo\"ide~$\Gamma'$ dans le groupo\"ide~$\Gamma$, qu'on
d\'enotera par $\tH$:
$$
\tH (F')(X)=F' (H(X))
$$
qu'on peut aussi \'ecrire, avec la notation $F(X)=E_X (F)$
introduite dans~\no \Ref{V.6}:
$$
\ifthenelse{\boolean{orig}}
{E_{H(X)} (F')=E_X (\tH(F'))}
{E_{H(X)} (F')=E_X (\tH(F')).}
$$
Cette derni\`ere formule montre, compte tenu de~\Ref{V.5.8}
ou~\Ref{V.4.1}, que le foncteur exact $H$ est d\'etermin\'e (\`a
isomorphisme unique pr\`es) quand on conna\^it le foncteur
correspondant~$\tH$. Fixons nous un~$F'$, soit $F=\tH (F')$,
alors $\tH$ d\'efinit un homomorphisme de~$\pi_{F'} =\Aut(F')$
dans $\pi_F =\Aut(F)$:
$$
\tH\colon \pi_{F'} \to \pi_F \quad (F=\tH (F')=F' \circ H).
$$
D'ailleurs la formule plus haut montre (compte tenu de~\Ref{V.5.8})
que cet homomorphisme a la propri\'et\'e suivante: pour tout
ensemble fini $E$ o\`u $\pi_F$ op\`ere contin\^ument, le groupe
$\pi_{F'}$ op\`ere \'egalement \emph{contin\^ument} gr\^ace
\`a l'homomorphisme pr\'ec\'edent $\pi_{F'} \to
\pi_F$. \hbox{Appliquant} ceci aux quotients de~$\pi_F$ par ses sous-groupes
invariants ouverts, on voit que la condition pr\'ec\'edente
signifie aussi que l'homomorphisme consid\'er\'e est
continu. R\'eciproquement,
\ifthenelse{\boolean{orig}}
{donnons nous}
{donnons-nous}
un objet $F$ de~$\Gamma$, un
objet $F'$ de~$\Gamma'$ et un homomorphisme continu
$$
u\colon \pi_{F'} \to \pi_F ,
$$
il lui correspond donc un foncteur de~$\cal{C} (\pi)$
dans~$\cal{C}(\pi')$, manifestement exact, donc en vertu
de~\Ref{V.4.1} il lui correspond un foncteur $H$ de~$\cal{C}$ dans
$\cal{C'}$ qui est exact, et tel que $\tH\colon \pi_{F'} \to \pi_F$
soit pr\'ecis\'ement~$u$. On peut aussi, au lieu d'un
homomorphisme de groupes, partir d'un \emph{foncteur}
$$
U\colon \Gamma' \to \Gamma
$$
qui est tel que pour \emph{tout} $F' \in \Gamma'$ (ou \emph{un} $F'
\in \Gamma'$, cela revient au m\^eme) l'homomorphisme correspondant
$\pi_{F'} \to \pi_F$ soit continu: un tel foncteur est isomorphe \`a
un foncteur de la forme~$\tH$, o\`u $H\colon \cal{C} \to
\cal{C'}$ est un foncteur exact d\'etermin\'e \`a isomorphisme
unique pr\`es. Ainsi:

\begin{corollaire}
\label{V.6.2}
Pour
\marginpar{136}
qu'un foncteur $H\colon \cal{C} \to \cal{C'}$ de cat\'egories
galoisiennes soit exact, il faut et il suffit qu'il existe des
\'equivalences $\cal{C}(\pi)\to\cal{C}$ et
\ifthenelse{\boolean{orig}}
{$\cal{C'}\to\cal{C}(\pi)$}
{$\cal{C'}\to\cal{C}(\pi')$}
qui transforment le foncteur $H$ en le foncteur
$\cal{C}(\pi)\to\cal{C}(\pi')$ associ\'e \`a un homomorphisme de
groupes topologiques $\pi'\to\pi$.
\end{corollaire}

\begin{corollaire}
\label{V.6.3}
Soient $\cal{C}$, $\cal{C'}$ deux cat\'egories galoisiennes,
$\Gamma$, $\Gamma'$ leurs groupo\"ides fondamentaux. Alors la
cat\'egorie des foncteurs exacts de~$\cal{C}$ dans $\cal{C'}$ est
\'equivalente \`a la cat\'egorie des foncteurs $U\colon \Gamma'
\to \Gamma$ ayant la propri\'et\'e suivante: pour tout $F'$ dans~$\Gamma'$ (ou \emph{un} $F'$ dans~$\Gamma'$, cela revient au
m\^eme), posant $F=U(F')$, l'homomorphisme
$$
\pi_{F'} =\Aut (F')\to \pi_F =\Aut (F)
$$
d\'efini par $U$ est continu.
\end{corollaire}

Consid\'erons le pro-groupe fondamental $\Pi$ de~$\cal{C}$, alors un
foncteur exact $H$ le transforme en un pro-groupe $H(\Pi)$
de~$\cal{C'}$, nous allons d\'efinir un homomorphisme
$$
\Pi' \to H(\Pi)
$$
(o\`u $\Pi'$ est le pro-groupe fondamental de~$\cal{C'}$), par la
condition que pour tout objet $F'$ de~$\Gamma'$, l'homomorphisme
correspondant
$$
F' (\Pi')=\pi_{F'} \to F' (H(\Pi))=\pi_F \qquad
\ifthenelse{\boolean{orig}}
{\text{(o\`u $F=F'\circ H=\tH(F')$)};}
{\text{(o\`u $F=F'\circ H=\tH(F')$)}}
$$
soit l'homomorphisme naturel
$$
\Aut(F')\to \Aut(F' \circ H).
$$
Comme ce dernier est fonctoriel en~$F'$, il d\'efinit bien en vertu
de~\Ref{V.5.8} un homomorphisme de pro-objets, et en fait de
pro-groupes, de~$\cal{C'}$. Cet homomorphisme est dit
\emph{associ\'e} au foncteur~$H$.

Soit maintenant $H'$ un deuxi\`eme foncteur exact, de la
cat\'egorie galoisienne $\cal{C'}$ dans une cat\'egorie
galoisienne~$\cal{C''}$. Il est trivial qu'on a
$$
{}^t (H' H)=\tH\,\tH'
$$
(N.B. on a l\`a une identit\'e de foncteurs, et non seulement un
isomorphisme canonique). On a une propri\'et\'e de
transitivit\'e analogue pour les homomorphismes associ\'es des
pro-groupes fondamentaux.

Nous
\marginpar{137}
allons maintenant interpr\'eter les propri\'et\'es du foncteur
exact $H$ en termes de l'homomorphisme correspondant
$$
u\colon \pi_{F'} \to \pi_{F} \qquad \text{(o\`u $F' =F'\circ H$)}.
$$
Il est commode d'introduire la notion d'\emph{objet ponctu\'e}
\index{ponctu\'e (objet)|hyperpage}%
de la cat\'egorie galoisienne $\cal{C}$ (muni de son foncteur
fondamental $F$): c'est par d\'efinition un objet $X$ de~$\cal{C}$
muni d'un \'el\'ement $a$ de~$F(X)$. Il s'interpr\`ete donc
comme un ensemble fini o\`u $\pi_F$ op\`ere contin\^ument \`a
gauche, muni d'un point~$a$. Donc les objets ponctu\'es
\emph{connexes} de~$\cal{C}$ s'identifient en vertu de~\Ref{V.5.3} aux
sous-groupes ouverts de~$\pi_F$. Si $U$ et $V$ sont deux tels
sous-groupes, correspondants \`a des objets ponctu\'es connexes
$X$, $Y$ de~$\cal{C}$, alors il existe un homomorphisme ponctu\'e de
$X$ dans $Y$ si et seulement si $U\subset V$, et cet homomorphisme est
alors unique. Bien entendu, le foncteur $H$ transforme objets
ponctu\'es en objets ponctu\'es (puisque $F=F' \circ H$). Notons
d'autre part qu'un sous-groupe ferm\'e d'un groupe tel que $\pi_F$
est l'intersection des sous-groupes ouverts qui le contiennent; par
suite,
\ifthenelse{\boolean{orig}}{\ignorespaces}{si}
$M$, $N$ sont deux sous-groupes ferm\'es, alors $M\subset N$ si et
seulement si tout sous-groupe ouvert qui contient $N$ contient
\'egalement~$M$. Gr\^ace \`a ces remarques, on prouve facilement
les r\'esultats qui suivent:

\begin{proposition}
\label{V.6.4}
Soit $X$ un objet ponctu\'e connexe de~$\cal{C}$, associ\'e \`a
un sous-groupe ouvert $U$ de~$\pi_F$. Pour que
\ifthenelse{\boolean{orig}}{$u$}{$U$}
contienne $u(\pi_{F'})$ (\resp le sous-groupe invariant ferm\'e
engendr\'e par~$u(\pi_{F'})$) il faut et il suffit que $H(X)$
admette une section ponctu\'ee (\resp soit compl\`etement
d\'ecompos\'e).
\end{proposition}

On appellera \emph{section} --- sous-entendu: au-dessus de l'objet
final --- d'un objet $X$ d'une cat\'egorie galoisienne~$\cal{C}$, un
morphisme de l'objet final $e_{\cal{C}}$ dans~$X$, ce qui revient
\`a la donn\'ee d'un \'el\'ement $a$ de~$F(X)$ invariant par
$\pi_F$; si $X$ est ponctu\'e, on dit qu'on a une \emph{section
ponctu\'ee} si elle est compatible avec les structures ponctu\'ees
sur~$X$ et~$e_{\cal{C}}$, \ie si $a$ est pr\'ecis\'ement l'objet
marqu\'e de~$F(X)$. Une telle section est donc unique, et existe si
et seulement si l'objet marqu\'e de~$F(X)$ est invariant
par~$\pi_F$. Enfin, un objet d'une cat\'egorie galoisienne est dit
\emph{compl\`etement d\'ecompos\'e}
\index{compl\`etement d\'ecompos\'e (objet)|hyperpage}%
\index{decompose (objet completement)@d\'ecompos\'e (objet compl\`etement)|hyperpage}%
s'il est isomorphe \`a une somme d'objets finaux, \ie si $\pi_F$
op\`ere trivialement dans $F(X)$ --- condition \'evidemment plus
forte que l'existence d'une section ponctu\'ee, lorsque $X$ est
ponctu\'e. La proposition~\Ref{V.6.4} r\'esulte trivialement des
d\'efinitions et remarques pr\'ec\'edentes.

\begin{corollaire}
\label{V.6.5}
Pour
\marginpar{138}
que $u$ soit trivial, il faut et il suffit que pour tout objet
$X$ de~$\cal{C}$, $H(X)$ soit compl\`etement d\'ecompos\'e.
\end{corollaire}

\begin{proposition}
\label{V.6.6}
Soit $X'$ un objet ponctu\'e connexe de~${\cal{C}}'$, associ\'e
\`a un sous-groupe ouvert $U'$ de~$\pi_{F'}$. Pour que $U'$
contienne~$\Ker u$, il faut et il suffit qu'il existe un objet
ponctu\'e connexe $X$ de~$\cal{C}$ et un homomorphisme ponctu\'e
de la composante connexe ponctu\'ee $X'_0$ de~$H(X)$ dans~$X'$ (donc
que $X$ soit isomorphe comme objet ponctu\'e \`a un quotient de la
composante connexe neutre de l'image inverse d'un objet ponctu\'e
de~$\cal{C}$). Si $u$ est surjectif, la condition pr\'ec\'edente
\'equivaut aussi \`a la suivante: $X$ est isomorphe \`a
un~$H(X)$, o\`u $X$ est un objet ponctu\'e de~$\cal{C}$.
\end{proposition}

(On appelle \emph{composante connexe neutre}
\index{neutre (composante connexe d'un objet ponctu\'e)|hyperpage}%
d'un objet ponctu\'e $X$ d'une cat\'egorie galoisienne~$\cal{C}$,
l'unique sous-objet connexe ponctu\'e de~$X$; il correspond \`a la
trajectoire sous $\pi_F$ du point marqu\'e de~$F(X)$, en vertu
de~\Ref{V.5.3}). Comme le fait que $U'$ contienne $\Ker u$ ne
d\'epend pas de la ponctuation choisie de~$X'$ (car une autre
ponctuation revient \`a remplacer $U$ par un sous-groupe
conjugu\'e \`a~$U$), on voit:

\begin{corollaire}
\label{V.6.7}
Pour que $U'$ contienne~$\Ker u$, il faut et il suffit qu'il existe un
\hbox{objet~$X$} de~$\cal{C}$ (qu'on peut supposer connexe) et un morphisme
d'une composante connexe de~$H(X)$ dans~$X'$. Si $u$ est surjectif,
cela signifie aussi que $X'$ est isomorphe \`a un objet de la
forme~$H(X)$.
\end{corollaire}

\begin{corollaire}
\label{V.6.8}
Pour que $u$ soit injectif, il faut et il suffit que pour tout objet
$X'$ de~${\cal{C}}'$, il existe un objet $X$ de~$\cal{C}$ et un
homomorphisme d'une composante connexe de~$H(X)$ dans~$X'$.
\end{corollaire}

\begin{proposition}
\label{V.6.9}
Les conditions suivantes sont \'equivalentes:
\begin{enumerate}
\item[(i)] L'homomorphisme $u\colon \pi_{F'} \to \pi_F$ est surjectif.
\item[(ii)] Pour tout objet connexe $X$ de~$\cal{C}$, $H(X)$ est
connexe.
\item[(iii)] Le foncteur $H$ est pleinement fid\`ele.
\end{enumerate}
\end{proposition}

Ce dernier fait signifie que pour deux objets~$X$, $Y$ de~$\cal{C}$,
l'application naturelle
$$
\Hom(X,Y) \to \Hom(H(X),H(Y))
$$
est bijective.

\begin{corollaire}
\label{V.6.10}
Pour
\marginpar{139}
que $u$ soit un isomorphisme, il faut et suffit que $H$ soit une
\'equivalence de cat\'egories, ou encore que les deux conditions
suivantes soient v\'erifi\'ees:
\begin{enumerate}
\item[a)] pour tout objet connexe $X$ de~$\cal{C}$, $H(X)$ est connexe
\item[b)] tout objet de~$\cal{C}'$ est isomorphe \`a un objet de la
forme~$H(X)$.
\end{enumerate}
\end{corollaire}

\begin{proposition}
\label{V.6.11}
Soient $H\colon \cal{C}\to \cal{C}'$ et $H'\colon\cal{C}\to \cal{C}''$
des foncteurs exacts entre cat\'egories galoisiennes, $F''$ un
foncteur fondamental sur~$\cal{C}''$, posons $F'=F''H'$ et $F=F'H$, et
consid\'erons les homomorphismes associ\'es
$$
u'\colon\pi_{F''}\to \pi_{F'}\qquad u\colon\pi_{F'}\to \pi_F.
$$
Pour que $\Ker u\subset \Im u'$ \ie pour que $uu'$ soit
l'homomorphisme trivial, il faut et il suffit que pour tout objet $X$
de~$\cal{C}$, $H'(H(X))$ soit compl\`etement d\'ecompos\'e. Pour
que $\Ker u \supset \Im u' $, il faut et il suffit que pour tout objet
ponctu\'e connexe $X'$ de~$\cal{C}'$ tel que $H'(X')$ admette une
section ponctu\'ee, il existe un objet~$X$ de~$\cal{C}$ et un
homomorphisme d'une composante connexe de~$H(X)$ dans~$X'$.
\end{proposition}

La premi\`ere assertion r\'esulte de la derni\`ere affirmation
de~\Ref{V.6.4}. La deuxi\`eme r\'esulte de la conjonction
de~\Ref{V.6.4} et~\Ref{V.6.6}.

\begin{remarque}
\label{V.6.12}
Il n'est pas vrai en g\'en\'eral, sous les conditions
de~\Ref{V.6.8} que $X'$ soit isomorphe \`a un objet de la forme
$H(X)$. On peut montrer que pour que tout objet connexe (donc tout
objet) de~$\cal{C}'$ soit isomorphe \`a un objet de la forme $H(X)$,
il faut et il suffit que $u$ soit un isomorphisme de~$\pi_{F'}$ sur un
sous-groupe \emph{facteur direct} de~$\pi_F$. En pratique cependant,
on construit directement un homomorphisme $\pi_F\to\pi_{F'}$ inverse
\`a droite de~$u$ \`a l'aide d'un foncteur exact convenable
de~$\cal{C'}$ dans~$\cal{C}$.
\end{remarque}

\begin{proposition}
\label{V.6.13}
Soient $\cal{C}$ une cat\'egorie galoisienne munie d'un foncteur
fondamental $F$, $S$ un objet connexe de~$\cal{C}$, $\cal{C}'$ la
cat\'egorie des objets de~$\cal{C}$ au-dessus de~$S$. Alors
$\cal{C'}$ est une cat\'egorie galoisienne, et le foncteur $X\mto
H(X)=X\times S$ de~$\cal{C}$ dans~$\cal{C}'$ est exact. Soit $a\in
F(S)$, et soit $F'$ le foncteur de~$\cal{C}'$ dans la cat\'egorie
des ensembles finis d\'efini par
$$
F'(X')=\text{image inverse de~$a$ par }F(X')\to F(S).
$$
Alors on a un isomorphisme $F\cong F'\circ H$, et l'homomorphisme
correspondant
$$
u:\pi_{F'}\to \pi_F
$$
\marginpar{140}%
est un isomorphisme de~$\pi_{F'}$ sur le sous-groupe ouvert $U$ de
$\pi_{F}$ stabilisateur de l'\'el\'ement marqu\'e $a$ de~$F(X)$.
\end{proposition}

La d\'emonstration est laiss\'ee au lecteur.

\section{Cas des pr\'esch\'emas}
\label{V.7}
Soit $S$ un pr\'esch\'ema localement noeth\'erien et
\emph{connexe}, et soit
$$
a\colon \Spec(\Omega)\to S
$$
un point g\'eom\'etrique de~$S$, \`a valeurs dans un corps
alg\'ebriquement clos $\Omega$. On posera
$$
\cal{C}=\text{cat\'egorie des rev\^etements \'etales de~$S$}
$$
et pour un objet $X$ de~$\cal{C}$, \ie un rev\^etement \'etale
$X$ de~$S$, on pose
$$
F(X)=\text{ensemble des points g\'eom\'etriques de~$X$ au-dessus
de~$a$.}
$$
Ainsi, $F$ devient un foncteur sur~$\cal{C}$ \`a valeurs dans la
cat\'egorie des ensembles finis. Les propri\'et\'es (G~1) \`a
(G~6) sont satisfaites: pour (G~1), c'est contenu dans les sorites de
I~\Ref{I.4.6}, (G~2) r\'esulte de~\Ref{V.3.4}, (G~3) de~\Ref{V.3.5},
(G~4) est trivial par d\'efinition, (G~5) r\'esulte de~\Ref{V.3.5} et
du d\'ebut du \No \Ref{V.2}, enfin (G~6) est prouv\'e dans~\Ref{V.3.7}. On
peut donc appliquer les r\'esultats des \No \Ref{V.4}, \Ref{V.5}, \Ref{V.6}. Cela permet en
particulier de d\'efinir un pro-objet $P$ de~$\cal{C}$
repr\'esentant $F$, appel\'e \emph{rev\^etement universel de~$S$
au point~$a$},
\index{revetement universel d'un preschema en un point@rev\^etement universel d'un pr\'esch\'ema en un point|hyperpage}%
et un groupe topologique $\pi=\Aut(F)=\Aut^0 P$, appel\'e
\emph{groupe fondamental de~$S$ en~$a$},
\index{groupe fondamental d'un pr\'esch\'ema en un point|hyperpage}%
et not\'e~$\pi_1(S,a)$.
\label{indnot:eg}\oldindexnot{$\pi_1(S,a)$|hyperpage}%
Le foncteur $F$ d\'efinit alors une \'equivalence de la
cat\'egorie $\cal{C}$ avec la cat\'egorie des ensembles finis
o\`u $\pi=\pi_1(S,a)$ op\`ere contin\^ument. Cette
\'equivalence permet donc d'interpr\'eter les op\'erations
courantes de limites projectives et limites inductives finies sur des
rev\^etements (produits, produits fibr\'es, sommes, passage au
quotient, \ifthenelse{\boolean{orig}}{etc...}{etc.}) en termes des op\'erations analogues dans
$\cal{C}(\pi)$, \ie en termes des op\'erations \'evidentes sur
des ensembles finis o\`u $\pi$ op\`ere. Notons d'ailleurs, puisque
les composantes connexes topologiques d'un rev\^etement \'etale
sont \'egalement des rev\^etements \'etales, qu'\emph{un \hbox{objet $X$} de~$\cal{C}$ est connexe dans $\cal{C}$ si et seulement si il est
topologiquement connexe}; en vertu de~\Ref{V.5.3} cela signifie donc
que $\pi_1$ op\`ere transitivement dans~$F(X)$.
\marginpar{141}
Notons que pour qu'un objet $X$ de~$\cal{C}$ soit fid\`element plat
et quasi-compact sur~$S$ (comme il est d\'ej\`a plat et
quasi-compact sur~$S$), il faut et suffit que $X\to S$ soit surjectif
\ie soit un \'epimorphisme dans $\cal{C}$, ou encore que $X\neq
\emptyset$. On conclut alors du crit\`ere \Ref{V.2.6}~(iii) que~$X$
\emph{est un rev\^etement principal de
\ifthenelse{\boolean{orig}}{$X$}{$S$}
de groupe $G$ si et seulement si il est un espace principal
homog\`ene sous $G$ dans la cat\'egorie $\cal{C}$}, (tel qu'il a
\'et\'e d\'efini dans \No \Ref{V.5}).

Si $a'$ est un autre point g\'eom\'etrique de~$S$ (correspondant
\`a un corps alg\'ebriquement clos $\Omega'$, qui peut \^etre
diff\'erent de~$\Omega$ et qui peut m\^eme avoir une
caract\'eristique diff\'erente), il d\'efinit un foncteur fibre
$F'=F_{a'}$ de~$\cal{C}$ dans la cat\'egorie des ensembles finis,
qui est encore exact, donc isomorphe \`a $F=F_a$. Par suite, les
groupes fondamentaux $\pi_1(S;a)$
\label{indnot:eg1}\oldindexnot{$\pi_1(S,a)$|hyperpage}%
%
pour $a$ variable sont isomorphes entre eux. Si on d\'esigne par
$\pi_1(S;a,a')$
\label{indnot:eh}\oldindexnot{$\pi_1(S;a,a')$|hyperpage}%
l'ensemble des isomorphismes (ou ce qui revient au m\^eme,
l'ensembles des homomorphismes) $F_a\to F_{a'}$ des foncteurs fibres
associ\'es, on obtient ainsi un \emph{groupo\"ide} dont
l'ensemble des objets est l'ensemble des points g\'eom\'etriques
de~$S$, les groupes fondamentaux \'etant les groupes
d'automorphismes des objets dudit groupo\"ide. L'ensemble
$\pi_1(S;a',a)$ peut \^etre appel\'e l'\emph{ensemble des classes
de chemins de~$a$ \`a $a'$}. Ces classes se composent donc de fa\c con \'evidente. Enfin, on peut d\'efinir un pro-groupe $\Pi_1^S$
de~$\cal{C}$, qu'on pourra appeler \emph{pro-groupe fondamental
de}~$S$
\index{pro-groupe fondamental d'un pr\'esch\'ema|hyperpage}%
ou \emph{syst\`eme local des groupes fondamentaux sur}~$S$,
\index{systeme local des groupes fondamentaux@syst\`eme local des groupes fondamentaux|hyperpage}%
d\'efini \`a isomorphisme unique pr\`es par la condition qu'on
ait un isomorphisme, fonctoriel en le point g\'eom\'etrique $a$ de
$S$:
$$
F_a(\Pi_1^S)=\pi_1(S;a)
$$
(\cf remarque~\Ref{V.5.10}). En particulier, si $s$ est un point
ordinaire de~$S$, la fibre de~$G$ en~$s$ est un pro-groupe sur~$\kres(s)$,
limite projective de groupes finis \'etales sur~$\kres(s)$; on pourrait
appeler ce pro-groupe \emph{le groupe fondamental de~$S$ en le point
ordinaire $s$ de~$S$}, et le noter $\pi_1(S,s)$. Par d\'efinition, ces
points \`a valeurs dans une extension alg\'ebriquement close
$\Omega$ de~$\kres(s)$ sont les \'el\'ements de~$\pi_1(S;a)$, o\`u
$a$ est le point g\'eom\'etrique de~$S$ d\'efini par ladite
extension. En particulier, (prenant pour $S$ le spectre d'un corps)
\`a tout corps $k$ est associ\'e canoniquement et fonctoriellement
un pro-groupe sur~$k$, qu'on pourrait noter $\pi_1(k)$, limite
projective de groupes finis \'etales sur~$k$, et dont les points
dans une extension alg\'ebriquement close $\Omega$ de~$k$
s'identifient
\marginpar{142}
aux \'el\'ements du groupe de Galois topologique de~$\bar{k}/k$,
o\`u $\bar{k}$ est la cl\^oture galoisienne de~$k$ dans $\Omega$
(\cf \Ref{V.8.1}). Ce groupe $\pi_1(k)$ ne semble pas encore avoir
retenu l'attention des alg\'ebristes.

Soit maintenant
$$
f\colon S'\to S
$$
un morphisme d'un pr\'esch\'ema connexe localement noeth\'erien
dans un autre, soit $a'$ un point g\'eom\'etrique de~$S'$ et soit
$a=f(a')$ son image directe dans $S$. Alors le foncteur \og image
inverse\fg induit un foncteur de la cat\'egorie $\cal{C}(S)$
\label{indnot:ei}\oldindexnot{$\cal{C}(S)$|hyperpage}%
des rev\^etements \'etales de~$S$, dans la cat\'egorie
$\cal{C}(S')$ des rev\^etements \'etales de~$S'$:
$$
f^\bbullet:\cal{C}(S)\to \cal{C}(S')
$$
On a d'ailleurs un isomorphisme de foncteurs
$$
F_a\cong F_{a'}\circ f^\bbullet,
$$
de sorte que $f^\bbullet$ est un foncteur \emph{exact}, auquel
s'appliquent les r\'esultats du \No \Ref{V.6}. On a en particulier un
homomorphisme canonique
$$
u=\pi_1(f;a'):\pi_1(S',a')\to \pi_1(S;a) \qquad (a'=f(a))
$$
\label{indnot:ej}\oldindexnot{$\pi_1(f;a')$|hyperpage}%
qui permet de reconstituer le foncteur image inverse, comme une
op\'eration de restriction des groupes d'op\'erateurs. Les
propri\'et\'es du foncteur $f^\bbullet$ s'expriment donc de
fa\c con simple par les propri\'et\'es de l'homomorphisme de
groupes associ\'es, comme il a \'et\'e explicit\'e dans le
\No \Ref{V.6}. Si en particulier $S'$ est un rev\^etement \'etale de~$S$,
alors l'homomorphisme $u$ est un isomorphisme de~$\pi_1(S',a')$ sur le
sous-groupe ouvert de~$\pi_1(S,a)$ qui d\'efinit le rev\^etement
\'etale connexe ponctu\'e $S'$ de~$S$ (\ie le stabilisateur~$U$
de~$a'\in F_a(S')$ dans $\pi_1(S;A)$).

Si on d\'esire interpr\'eter les homomorphismes $\pi_1(f;a')$ pour
un point g\'eom\'etrique variable $a'$, on doit, conform\'ement
\`a ce qui a \'et\'e dit dans le \No \Ref{V.6}, consid\'erer un
homomorphisme
$$
\Pi_1(f):\Pi_1^{S'}\to f^\bbullet(\Pi_1^S)
$$
de pro-groupes sur~$S$, et prendre l'homomorphisme correspondant pour
les fibres g\'eom\'etriques.

\section{Cas d'un pr\'esch\'ema de base normale}
\label{V.8}
\marginpar{143}
\begin{proposition}
\label{V.8.1}
Soit $S$ le spectre d'un corps $k$, et soit $\Omega$ une extension
alg\'ebriquement close de~$k$, d\'efinissant un point
g\'eom\'etrique $a$ de~$S$ \`a valeurs dans $\Omega$. Soit
$\bar{k}$ la cl\^oture s\'eparable de~$k$ dans $\Omega$. Alors il
existe un isomorphisme canonique de~$\pi_1(S,a)$ sur le groupe de
Galois topologique de~$\bar{k}/k$.
\end{proposition}

Soit $k'$ la cl\^oture alg\'ebrique de~$k$ dans $\Omega$, il
correspond donc \`a un point g\'eom\'etrique~$b$ de~$S$, \`a
valeurs dans $k'$. L'homomorphisme naturel de foncteurs $F_b\to F_a$
est \'evidemment un isomorphisme, car un $k$-homomorphisme d'une
extension finie s\'eparable de~$k$ dans $\Omega$ prend
n\'ecessairement ses valeurs dans $\bar{k}$ et a fortiori dans
$k'$. D'autre part, le groupe $\pi'$ des $k$-automorphismes de~$k'/k$
op\`ere de fa\c con \'evidente sur~$F_b$, d'o\`u un
homomorphisme $\pi'\to \Aut(F_b)\isomto \Aut(F_a)=\pi_1(S;a)$.
D'autre part, il est bien connu que l'homomorphisme naturel de~$\pi'$
dans le groupe $\pi$ des automorphismes de~$\bar{k}/k$ est un
isomorphisme. On obtient ainsi un homomorphisme canonique $\pi\to
\pi_1(S;a)$, reste \`a montrer que c'est un isomorphisme. En effet,
cet homomorphisme est injectif, car un \'el\'ement du noyau est un
automorphisme de~$\bar{k}/k$ qui induit l'identit\'e sur toute
sous-extension s\'eparable finie, donc est trivial. Cet
homomorphisme est surjectif, car si $X$ est un rev\^etement
\'etale \emph{connexe} de~$S$, donc d\'efini par une
\emph{extension} finie s\'eparable~$L/k$, alors $\pi$ est transitif
sur l'ensemble des $k$-homomorphismes de~$L$ dans $k'$, comme bien
connu.

\begin{proposition}
\label{V.8.2}
Soient $S$ un pr\'esch\'ema connexe, localement noeth\'erien et
normal, $K=\kres(s)$ son corps des fonctions $=$ le corps r\'esiduel en
son point g\'en\'erique $s$, $\Omega$ une extension
alg\'ebriquement close de~$K$, d\'efinissant un point
g\'eom\'etrique $a'$ de~$S'=\Spec(K)$ et un point
g\'eom\'etrique $a$ de~$S$. Alors l'homomorphisme
$\pi_1(S';a')\to\pi_1(S;a)$ est surjectif. Lorsqu'on identifie le
premier groupe au groupe de Galois de la cl\^oture s\'eparable
$\bar{K}$ de~$K$ dans~$\Omega$ (\cf \Ref{V.8.1}) alors le noyau de
l'homomorphisme pr\'ec\'edent correspond par la th\'eorie de
Galois \`a la sous-extension de~$\bar{K}/K$ compos\'ee des
extensions finies de~$K$ dans $\Omega$ qui sont non ramifi\'ees sur~$S$.
\end{proposition}

La premi\`ere assertion signifie que l'image inverse sur~$S'$ d'un
rev\^etement \'etale connexe $X$ de~$S$ est connexe, \ie que $X$
est int\`egre, ce n'est autre que (I~\Ref{I.10.1}). Le noyau de
l'homomorphisme pr\'ec\'edent s'interpr\`ete alors comme
form\'e des automorphismes de~$\bar{K}/K$ qui induisent
l'identit\'e sur les ensembles~$F_a(X)$,
\marginpar{144}
o\`u on peut supposer le rev\^etement \'etale $X$ de~$S$
connexe. Mais cela signifie que cet automorphisme induit
l'identit\'e sur les sous-extensions finies de~$\bar{K}/K$ qui sont
non ramifi\'ees sur~$S$, ce qui prouve la derni\`ere assertion.

\begin{remarquestar}
Gr\^ace \`a cette interpr\'etation du groupe fondamental du
pr\'esch\'ema normal $S$ en termes de th\'eorie des Galois
habituelle, la d\'efinition \'etait connue dans ce cas depuis
longtemps.
\end{remarquestar}

\section{Cas des pr\'esch\'emas non connexes: cat\'egories
multigaloisiennes}
\label{V.9}

Soit $S$ un pr\'esch\'ema localement noeth\'erien, et soient
$(S_i)_{i\in I}$ ses composantes connexes. Alors la cat\'egorie
$\cal{C}(S)$ des rev\^etements \'etales de~$S$ est \'equivalente
\`a la cat\'egorie produit des $\cal{C}(S_i)$, qui
s'interpr\`etent en termes des groupes fondamentaux des $S_i$, une
fois choisie un point g\'eom\'etrique dans chaque $S_i$. Dans
l'application de la th\'eorie de la descente pour les morphismes
\'etales, il est parfois malcommode de faire choix pour tout $S_i$
d'un point g\'eom\'etrique de~$S_i$. Il est plus commode alors
de recourir \`a la g\'en\'eralisation naturelle de~\Ref{V.5.8}
pour interpr\'eter $\cal{C}(S)$ comme une cat\'egorie de foncteurs
sur le groupo\"ide des points g\'eom\'etriques de~$S$,
consid\'er\'e comme somme des groupo\"ides correspondants aux
composantes connexes de~$S$; les foncteurs en question sont les
foncteurs \`a valeurs dans la cat\'egorie des ensembles finis,
satisfaisant la propri\'et\'e de continuit\'e analogue \`a
celle invoqu\'ee dans~\Ref{V.5.8}. En pratique, on aura une famille
$(a_t)_{t\in E}$ de points g\'eom\'etriques de~$S$, telle que
toute composante connexe $S_i$ de~$S$ en contienne au moins un, et on
pourra alors, comme dans~\Ref{V.5.10}, remplacer le groupo\"ide de
tous les points g\'eom\'etriques de~$S$ par le groupo\"ide
analogue dont l'ensemble sous-jacent est~$E$. Bien entendu, ces
consid\'erations devraient s'ins\'erer dans des d\'efinitions
g\'en\'erales concernant les cat\'egories qui sont
\'equivalentes \`a des cat\'egories produits de cat\'egories
de la forme $\cal{C}(\pi)$, et qu'on pourra appeler
\emph{cat\'egories multi-galoisiennes}.
\index{multi-galoisienne (cat\'egorie)|hyperpage}%
Nous en laisserons le d\'etail au lecteur.

\chapter{Cat\'egories fibr\'ees et descente}
\label{VI}
\marginpar{145}

\setcounter{section}{-1}
\section{Introduction}

Contrairement \`a ce qui avait \'et\'e annonc\'e dans
l'introduction \`a l'expos\'e pr\'ec\'edent, il s'est
av\'er\'e impossible de faire de la descente dans la cat\'egorie
des pr\'esch\'emas, m\^eme dans des cas particuliers, sans avoir
d\'evelopp\'e au pr\'ealable avec assez de soin le langage de la
descente dans les cat\'egories g\'en\'erales.

La notion de \og descente\fg fournit le cadre g\'en\'eral pour tous
les proc\'ed\'es de \og recollement\fg d'objets, et par
cons\'equent de \og recollement\fg de cat\'egories. Le cas le plus
classique de recollement est relatif \`a la donn\'ee d'un espace
topologique $X$ et d'un recouvrement de~$X$ par des ouverts $X_i$; si
on se donne pour tout $i$ un espace fibr\'e (disons) $E_i$ au-dessus
de~$X_i$, et pour tout couple $(i,j)$ un isomorphisme $f_{ji}$ de
$E_i|X_{ij}$ sur~$E_j|X_{ij}$ (o\`u on pose $X_{ij}=X_i\cap X_j$),
satisfaisant une condition de transitivit\'e bien connue (qu'on
\'ecrit de fa\c con abr\'eg\'ee $f_{kj}f_{ji}=f_{ki}$), on
sait qu'il existe un espace fibr\'e $E$ sur~$X$, d\'efini \`a
isomorphisme pr\`es par la condition que l'on ait des isomorphismes
$f_i\colon E|X_i\isomto E_i$, satisfaisant les relations
$f_{ji}=f_jf_i^{-1}$ (avec l'abus d'\'ecriture habituel). Soit $X'$
l'espace somme des $X_i$, qui est donc un espace fibr\'e au-dessus
de~$X$ (\ie muni d'une application continue $X'\to X$). On peut
interpr\'eter de fa\c con plus concise la donn\'ee des $E_i$
comme un espace fibr\'e $E'$ sur~$X'$, et la donn\'ee des $f_{ji}$
comme un isomorphisme entre les deux images inverses (par les deux
projections canoniques) $E''_1$ et $E''_2$ de~$E'$ sur~$X''=X'\times_X
X'$, la condition de recollement pouvant alors s'\'ecrire comme une
identit\'e entre isomorphismes d'espaces fibr\'es $E'''_1$ et
$E'''_3$ sur le produit fibr\'e triple $X'''=X'\times_X X'\times_X
X'$ (o\`u $E'''_i$ d\'esigne l'image inverse de~$E'$ sur~$X'''$
par la projection canonique d'indice $i$). La construction de~$E$,
\`a partir de~$E'$ et de~$f$, est un cas typique
\marginpar{146}
de proc\'ed\'e de \og descente\fg. D'ailleurs, partant d'un espace
fibr\'e $E$ sur~$X$, on dit que $X$ est \og localement trivial\fg, de
fibre $F$, s'il existe un recouvrement ouvert $(X_i)$ de~$X$ tel que
les $E|X_i$ soient isomorphes \`a $F\times X_i$, ou ce qui revient
au m\^eme, tel que l'image inverse $E'$ de~$E$ sur~$X'=\coprod_i
X_i$ soit isomorphe \`a $X'\times F$.

Ainsi, la notion de \og recollement\fg d'objets comme celle de
\og localisation\fg d'une propri\'et\'e, sont li\'ees \`a
l'\'etude de certains types de \og changements de base\fg $X'\to X$.
En G\'eom\'etrie Alg\'ebrique, bien d'autres types de changement
de base, et notamment les morphismes $X'\to X$ fid\`element plats,
doivent \^etre consid\'er\'es comme correspondant \`a un
proc\'ed\'e de \og localisation\fg relativement aux
pr\'esch\'emas, ou autres objets, \og au-dessus\fg de~$X$. Ce type de
localisation est utilis\'e tout autant que la simple localisation
topologique (qui en est un cas particulier d'ailleurs). Il en est de
m\^eme (dans une moindre mesure) en G\'eom\'etrie Analytique.

La plupart des d\'emonstrations, se r\'eduisant \`a des
v\'erifications, sont omises ou simplement esquiss\'ees; le cas
\'ech\'eant nous pr\'ecisons les diagrammes moins \'evidents
qui s'introduisent dans une d\'emonstration.

\section{Univers, cat\'egories, \'equivalence de cat\'egories}
\label{VI.1}

Pour \'eviter certaines difficult\'es logiques, nous admettrons
ici la notion d'\emph{Univers}, qui est un ensemble \og assez gros\fg
pour qu'on n'en sorte pas par les op\'erations habituelles de la
th\'eorie des ensembles; un \og \emph{axiome des Univers}\fg garantit
que tout objet se trouve dans un Univers. Pour des d\'etails, voir
un livre en pr\'eparation par C\ptbl Chevalley et le conf\'erencier\footnote{%
\ifthenelse{\boolean{orig}}
{Les auteurs d\'efinitifs sont C\ptbl Chevalley et P\ptbl Gabriel. Le livre doit sortir en l'an 2000. \Cf aussi SGA~4 VI.7.1.}
{Les auteurs d\'efinitifs sont C\ptbl Chevalley et P\ptbl Gabriel. Le livre doit sortir en l'an 3000. En attendant, \cf aussi SGA~4 I.}}. Ainsi, le sigle~$\Ens$
\label{indnot:fb}\oldindexnot{$\Ens$|hyperpage}%
d\'esigne, non pas la cat\'egorie de tous les ensembles (notion
qui n'a pas de sens), mais la cat\'egorie des ensembles qui se
trouvent dans un Univers donn\'e (que nous ne pr\'eciserons pas
ici dans la notation). De m\^eme, $\Cat$
\label{indnot:fc}\oldindexnot{$\Cat$|hyperpage}%
d\'esignera la cat\'egorie des cat\'egories se trouvant dans
l'Univers en question, les \og morphismes\fg d'un objet $X$ de~$\Cat$
dans un autre $Y$, \'etant par d\'efinition les \emph{foncteurs}
de~$X$ dans $Y$.

Si
\marginpar{147}
$\cal{C}$ est une cat\'egorie, nous d\'esignons par $\Ob(\cal{C})$
\label{indnot:fd}\oldindexnot{$\Ob(\cal{C})$|hyperpage}%
\emph{l'ensemble des objets} de~$\cal{C}$, par~$\Fl(\cal{C})$
\label{indnot:fe}\oldindexnot{$\Fl(\cal{C})$|hyperpage}%
\emph{l'ensemble des fl\`eches} de~$\cal{C}$ (ou morphismes de
$\cal{C}$). Nous \'ecrirons donc $X\in\Ob(\cal{C})$ en \'evitant
l'abus de notation courant $X\in\cal{C}$. Si $\cal{C}$ et $\cal{C}'$
sont deux cat\'egories, un \emph{foncteur} de~$\cal{C}$ dans
$\cal{C}'$ sera toujours ce qu'on appelle commun\'ement un foncteur
\emph{covariant} de~$\cal{C}$ dans $\cal{C}'$; sa donn\'ee implique
celle de la cat\'egorie d'arriv\'ee et la cat\'egorie de
d\'epart, $\cal{C}$ et $\cal{C}'$. Les foncteurs de~$\cal{C}$ dans
$\cal{C}'$ forment un ensemble, not\'e $\Hom(\cal{C},\cal{C}')$, qui
est l'ensemble des objets d'une cat\'egorie not\'ee
$\SheafHom(\cal{C},\cal{C}')$.
\label{indnot:ff}\oldindexnot{$\SheafHom(\cal{C},\cal{C}')$|hyperpage}%
Par d\'efinition, un \emph{foncteur contravariant} de~$\cal{C}$ dans
$\cal{C}'$ est un foncteur de la \emph{cat\'egorie oppos\'ee}
$\cal{C}^\circ$
\label{indnot:fg}\oldindexnot{$\cal{C}^\circ$|hyperpage}%
de~$\cal{C}$ dans~$\cal{C}'$.

Nous admettrons la notion de \emph{limite projective} et de
\emph{limite inductive} d'un foncteur $F\colon \cal{I}\to \cal{C}$, et
en particulier les cas particuliers les plus courants de cette notion:
produits cart\'esiens et produits fibr\'es, notions duales de
sommes directes et de sommes amalgam\'ees, et les propri\'et\'es
formelles habituelles de ces op\'erations.

Par exemple, dans la cat\'egorie $\Cat$ introduite plus haut, les
limites projectives (relatives \`a des cat\'egories $\cal{I}$ se
trouvant dans l'Univers choisi) existent; l'ensemble d'objets
(\resp l'ensemble de fl\`eches) de la cat\'egorie limite
projective $\cal{C}$ des $\cal{C}_i$, s'obtient en prenant la limite
projective des ensembles d'objets (\resp des ensembles de fl\`eches)
des cat\'egories $\cal{C}_i$. Le cas le plus connu est celui du
produit d'une famille de cat\'egories. Nous utiliserons constamment
par la suite le produit fibr\'e de deux cat\'egories sur une
troisi\`eme.

Pour tout ce qui concerne les cat\'egories et foncteurs, en
attendant le livre en pr\'eparation d\'ej\`a signal\'e, voir
\cite{VI.1} (qui est n\'ecessairement fort incomplet, m\^eme en ce qui
concerne les g\'en\'eralit\'es esquiss\'ees dans le
pr\'esent num\'ero).

Prenons cette occasion pour expliciter la notion d'\'equivalence de
cat\'egories,
\index{equivalence@\'equivalence de cat\'egories|hyperpage}%
qui n'est pas expos\'ee de fa\c con satisfaisante
dans~\cite{VI.1}. Un foncteur $F\colon\cal{C}\to \cal{C}'$ est dit
\emph{fid\`ele}
\index{fidele (foncteur)@fid\`ele (foncteur)|hyperpage}%
(\resp \emph{pleinement fid\`ele})
\index{pleinement fid\`ele (foncteur)|hyperpage}%
si pour tout couple d'objets $S$, $T$ de~$\cal{C}$, l'application
$u\mto F(u)$ de~$\Hom(S,T)$ dans $\Hom(F(S),F(T))$ est injective
(\resp bijective). On dit que $F$ est une \emph{\'equivalence} de
cat\'egories si
\marginpar{148}
$F$ est pleinement fid\`ele, et si de plus tout objet $S'$ de
$\cal{C}'$ est isomorphe \`a un objet de la forme $F(S)$. On montre
qu'il revient au m\^eme de dire qu'il existe un foncteur $G$ de
$\cal{C}'$ dans $\cal{C}$ \emph{quasi-inverse de} $F$, \iev tel que
$GF$ soit isomorphe \`a $\id_{\cal{C}}$. Lorsqu'il en est ainsi, la
donn\'ee d'un foncteur $G\colon\cal{C}'\to \cal{C}$ et d'un
isomorphisme $\varphi\colon GF\to \id_{\cal{C}'}$ \'equivaut \`a
la donn\'ee, pour tout $S'\in\Ob(\cal{C}')$, d'un couple $(S,u)$
form\'e d'un objet $S$ de~$\cal{C}$ et d'un isomorphisme $u\colon
F(S)\to S'$, soit $(G(S),\varphi(S))$. (Avec ces notations, il existe
un foncteur unique $\cal{C}'\to\cal{C}$ ayant l'application donn\'ee
$S\mto G(S)$ comme application-objets, et tel que l'application
$S\mto \varphi(S)$ soit un homomorphisme de foncteurs $FG\to
\id_{\cal{C}'}$). Enfin, si $G$ est un foncteur quasi-inverse de~$F$,
et si on choisit des isomorphismes $\varphi\colon
FG\isomto\id_{\cal{C}'}$ et $\psi\colon GF\isomto\id_{\cal{C}}$, alors
les deux conditions de compatibilit\'es sur~$\varphi$, $\psi$
\'enonc\'ees dans \cite[I.1.2]{VI.1} sont en fait \'equivalentes l'une
\`a l'autre, et pour tout isomorphisme $\varphi$ choisi, il existe
un isomorphisme $\psi$ unique tel que lesdites conditions soient
satisfaites.

\section{Cat\'egories sur une autre}
\label{VI.2}

Soit $\cal{E}$ une cat\'egorie dans~$\Univ$, c'est donc un objet
de~$\Cat$, et on peut consid\'erer la cat\'egorie
$\Cat_{/\cal{E}}$
\label{indnot:fh}\oldindexnot{$\Cat_{/\cal{E}}$|hyperpage}%
des \og objets de~$\Cat$ au-dessus de~$\cal{E}$\fg. Un objet de cette
cat\'egorie est donc un foncteur
$$
p\colon \cal{F}\to\cal{E}
$$
On dit aussi que la cat\'egorie~$\cal{F}$, munie d'un tel foncteur,
est une \emph{cat\'egorie au-dessus de}~$\cal{E}$, ou une
$\cal{E}$\emph{-cat\'egorie}. On appellera donc $\cal{E}$-foncteur
d'une cat\'egorie $\cal{F}$ sur~$\cal{E}$ dans une cat\'egorie
$\cal{G}$ sur~$\cal{E}$, un foncteur
$$
f\colon\cal{F}\to\cal{G}
$$
tel que
$$
qf=p
$$
o\`u $p$ et $q$ sont les foncteurs-projection pour $\cal{F}$
\resp $\cal{G}$. L'ensemble des $\cal{E}$-foncteurs $f$ de~$\cal{F}$
dans $\cal{G}$ est donc en correspondance biunivoque avec l'ensemble
des fl\`eches d'origine $\cal{F}$ et d'extr\'emit\'e $\cal{G}$
dans~$\Cat_{/\cal{E}}$, sans
\marginpar{149}
pourtant qu'on ait l\`a une identit\'e (puisque la donn\'ee d'un
$f$ comme dessus ne d\'etermine pas $\cal{F}$ et $\cal{G}$ en tant
que cat\'egories sur~$\cal{E}$); mais bien entendu, comme dans toute
autre cat\'egorie~$\cal{C}_{/S}$, on fera couramment l'abus de
langage consistant \`a identifier les $\cal{E}$-foncteurs (au sens
explicit\'e plus haut) \`a des fl\`eches dans une
cat\'egorie~$\Cat_{/\cal{E}}$.

On d\'esignera par
$$
\Hom_{\cal{E}}(\cal{F},\cal{G})
$$
l'ensemble des $\cal{E}$-foncteurs de~$\cal{F}$ dans~$\cal{G}$. Bien
entendu, un compos\'e de~$\cal{E}$-foncteurs est un
$\cal{E}$-foncteur (la composition en question correspondant par
d\'efinition \`a la composition des fl\`eches
dans~$\Cat_{/\cal{E}}$).

Consid\'erons maintenant deux $\cal{E}$-foncteurs
$$
f,g\colon \cal{F}\to\cal{G}
$$
et un homomorphisme de foncteurs:
$$
u\colon f\to g
$$
On dit que $u$ est un $\cal{E}$-\emph{homomorphisme} ou un
\og \emph{homomorphisme de~$\cal{E}$-foncteurs}\fg, si pour tout
$\xi\in\Ob(\cal{F})$, on a
$$
q(u(\xi))=\id_{p(\xi)} \quoi\text{,}
$$
en paroles: posant
\ifthenelse{\boolean{orig}}
{$S=p(\xi)=qf(\xi)=qg(\xi)\in\Ob\cal{E}$,}
{$S=p(\xi)=qf(\xi)=qg(\xi)\in\Ob(\cal{E})$,}
le morphisme
$$
u(\xi)\colon f(\xi)\to g(\xi)
$$
dans $\cal{G}$ est un $\id_S$-morphisme. (De fa\c con
g\'en\'erale, pour tout morphisme $\alpha\colon T \to S$
dans~$\cal{E}$, et toute cat\'egorie $\cal{G}$ au-dessus
de~$\cal{E}$, un morphisme $v$ dans $\cal{G}$ est appel\'e un
$\alpha$-\emph{morphisme} si $q(v)=\alpha$, $q$ d\'esignant le
foncteur projection $\cal{G}\to \cal{E}$). Si on a un troisi\`eme
$\cal{E}$-foncteur $h\colon \cal{F}\to \cal{G}$ et un
$\cal{E}$-homomorphisme $v\colon g\to h$, alors $vu$ est \'egalement
un $\cal{E}$-homomorphisme. Ainsi,
\marginpar{150}
\emph{les $\cal{E}$-foncteurs de~$\cal{F}$ dans~$\cal{G}$, et
les $\cal{E}$-homomorphismes de tels, forment une sous-cat\'egorie
de la cat\'egorie $\SheafHom(\cal{F},\cal{G})$ de tous les foncteurs
de~$\cal{F}$ dans~$\cal{G}$, qu'on appellera la cat\'egorie des
$\cal{E}$-foncteurs de~$\cal{F}$ dans $\cal{G}$ et qu'on notera}
$$
\SheafHom_{\cal{E}/-}(\cal{F},\cal{G})
$$
\label{indnot:fi}\oldindexnot{$\SheafHom_{\cal{E}/-}(\cal{F},\cal{G})$|hyperpage}%
C'est aussi la sous-cat\'egorie noyau du couple de foncteurs
$$
\xymatrix@C=.5cm{R,S\colon\SheafHom(\cal{F},\cal{G})
\ar@<2pt>[r]\ar@<-2pt>[r]&\SheafHom(\cal{F},\cal{E})}\quoi\text{,}
$$
o\`u $R$ est le foncteur constant d\'efini par l'objet $p$ de
$\SheafHom(\cal{F},\cal{E})$, et o\`u $S$ est le foncteur $f\mto
q\circ f$ d\'efini par $q\colon\cal{G}\to\cal{E}$.

Pour finir ces g\'en\'eralit\'es, il reste \`a d\'efinir les
accouplements naturels entre les cat\'egories
$\SheafHom_{\cal{E}/-}(\cal{F},\cal{G})$ par composition de
$\cal{E}$-foncteurs. En d'autres termes, on veut d\'efinir un
\og foncteur composition\fg:
\begin{equation*}
\label{eq:VI.2.i}
\tag{i} {\SheafHom_{\cal{E}/-}(\cal{F},\cal{G})\times
\SheafHom_{\cal{E}/-}(\cal{G},\cal{H})\to
\SheafHom_{\cal{E}/-}(\cal{F},\cal{H})}
\end{equation*}
lorsque $\cal{F}$, $\cal{G}$, $\cal{H}$ sont trois cat\'egories
sur~$\cal{E}$, de telle fa\c con que ce foncteur induise, pour les
objets, l'application de composition $(f,g)\mto gf$ de
$\cal{E}$-foncteurs $f\colon\cal{F}\to\cal{G}$ et $g\colon \cal{G}\to
\cal{H}$. Pour ceci, rappelons qu'on d\'efinit un foncteur
canonique:
\begin{equation*}
\label{eq:VI.2.ii}
\tag{ii}
{\SheafHom(\cal{F},\cal{G})\times\SheafHom(\cal{G},\cal{H})\to
\SheafHom(\cal{F},\cal{H})}
\end{equation*}
qui, pour les objets, n'est autre que l'application de composition
$(f,g)\mto gf$ de foncteurs, et qui transforme une fl\`eche
$(u,v)$,~o\`u
$$
u\colon f\to f'\quoi,\quad v\colon g\to g'
$$
sont des fl\`eches dans $\SheafHom(\cal{F},\cal{G})$ \resp dans
$\SheafHom(\cal{G},\cal{H})$, en la fl\`eche
$$
v*u\colon gf\to g'f'
$$
\label{indnot:fj}\oldindexnot{$v*u$|hyperpage}%
d\'efinie par la relation:
\marginpar{151}
$$
v*u(\xi)=v(f'(\xi))\cdot g(u(\xi))=g'(u(\xi))\cdot v(f(\xi))
$$
Il est bien connu que l'on obtient bien ainsi un homomorphisme de~$gf$
dans~$g'f'$, et que (pour $f,g$ et $u,v$ variables) on obtient ainsi
un foncteur~\eqref{eq:VI.2.ii}, \ie qu'on a
\begin{equation*}
\label{eq:VI.2.I}
\tag{I} {\id_g*\id_f=\id_{gf}\quoi,}
\end{equation*}
\begin{equation*}
\label{eq:VI.2.II}
\tag{II} {(v'*u')\circ(v*u)=(v'\circ v)*(u'\circ u)}
\end{equation*}
Rappelons aussi qu'on a une formule d'associativit\'e pour les
accouplements canoniques~\eqref{eq:VI.2.ii}, qui s'exprime d'une part
par l'associativit\'e $(hg)f=h(gf)$ de la composition de foncteurs,
et d'autre part par la formule
\begin{equation*}
\label{eq:VI.2.III}
\tag{III} {(w*v)*u=w*(v*u)}
\end{equation*}
pour les produits de composition d'homomorphismes de foncteurs (o\`u
$u\colon f\to f'$ et $v\colon g\to g'$ sont comme dessus, et o\`u on
suppose donn\'e de plus un homomorphisme $w\colon h\to h'$ de
foncteurs $h,h'\colon \cal{H}\to\cal{K}$). Je dis maintenant que
\emph{lorsque~$\cal{F}$, $\cal{G}$ sont des $\cal{E}$-cat\'egories,
le foncteur composition canonique}~\eqref{eq:VI.2.ii} \emph{induit un
foncteur}~\eqref{eq:VI.2.i}. Comme on sait d\'ej\`a que le
compos\'e de deux $\cal{E}$-foncteurs est un $\cal{E}$-foncteur,
cela revient \`a dire que \emph{lorsque $u\colon f\to f'$ et
$v\colon g\to g'$ sont des homomorphismes de~$\cal{E}$-foncteurs,
alors $v*u\colon gf\to g'f'$ est \'egalement un homomorphisme de
$\cal{E}$-foncteurs.} Cela r\'esulte en effet trivialement des
d\'efinitions. Comme les accouplements~\eqref{eq:VI.2.i} sont
induits par les accouplements~\eqref{eq:VI.2.ii}, ils satisfont \`a
la m\^eme propri\'et\'e d'associativit\'e, exprim\'ee aussi
dans les formules $(hg)f=h(gf)$ et $(w*v)*u=w*(v*u)$ pour des
$\cal{E}$-foncteurs et des $\cal{E}$-homomorphismes de
$\cal{E}$-foncteurs.

Pour compl\'eter le formulaire~\eqref{eq:VI.2.I},
\eqref{eq:VI.2.II}, \eqref{eq:VI.2.III}, rappelons aussi les formules:
\begin{equation*}
\label{eq:VI.2.IV}
\tag{IV} {v*\id_{\cal{F}}=v\quad\text{et}\quad\id_{\cal{G}}*u=u}\quoi,
\end{equation*}
o\`u
\marginpar{152}
pour simplifier on \'ecrit $v*f$ ou $u*g$ au lieu de~$v*u$, lorsque
$u$ \resp $v$ est l'automorphisme identique de~$f$ \resp $g$.

De la d\'efinition des accouplements~\eqref{eq:VI.2.i} r\'esulte
que $\SheafHom_{\cal{E}/-}(\cal{F},\cal{G})$ \emph{est un foncteur en
\ifthenelse{\boolean{orig}}
{$\cal{E}$,~$\cal{G}$,}
{$\cal{F}$,~$\cal{G}$,}
de la cat\'egorie produit
${\Cat_{/\cal{E}}}^\circ\times\Cat_{/\cal{E}}$ dans la
cat\'egorie~$\Cat$}. Si en effet $g\colon \cal{G}\to\cal{G}_1$ est
un $\cal{E}$-foncteur, \ie un objet de
$\SheafHom_{\cal{E}/-}(\cal{G},\cal{G}_1)$, alors faisant
dans~\eqref{eq:VI.2.i} $\cal{H}=\cal{G}_1$, il lui correspond un
foncteur
$$
g_*\colon\SheafHom_{\cal{E}/-}(\cal{F},\cal{G})\to
\SheafHom_{\cal{E}/-}(\cal{F},\cal{G}_1)
$$
On d\'efinit de la fa\c con analogue, pour un $\cal{E}$-foncteur
$f\colon \cal{F}_1\to \cal{F}$, un foncteur
$$
f^*:\SheafHom_{\cal{E}/-}(\cal{F},\cal{G})\to
\SheafHom_{\cal{E}/-}(\cal{F}_1,\cal{G})
$$
Pour abr\'eger, on d\'esigne ces foncteurs aussi par les sigles
$f\mto g\circ f$ \resp $g\mto g\circ f$ (qui d\'esignent
seulement, en fait, les applications correspondantes sur les ensembles
d'objets). Il r\'esulte de la propri\'et\'e d'associativit\'e
indiqu\'ee plus haut que de cette fa\c con, on obtient bien comme
annonc\'e un foncteur
$\Cat_{/\cal{E}}^\circ\times\Cat_{/\cal{E}}\to\Cat$.

\section{Changement de base dans les cat\'egories sur~$\mathcal{E}$}
\label{VI.3}

Comme dans $\Cat$ les limites projectives (relativement \`a des
cat\'egories $\cal{I}$ \'el\'ements de~$\Univ$) existent, il en
est de m\^eme dans $\Cat_{/\cal{E}}$, en particulier les produits
cart\'esiens y existent, qui s'interpr\`etent comme des produits
fibr\'es dans~$\Cat$. Conform\'ement aux notations
g\'en\'erales, si $\cal{F}$ et $\cal{G}$ sont des cat\'egories
au-dessus de~$\cal{E}$, on d\'esigne par
$$
\cal{F}\times_{\cal{E}}\cal{G}
$$
\label{indnot:fk}\oldindexnot{$\cal{F}\times_{\cal{E}}\cal{G}$|hyperpage}%
leur produit dans~$\Cat_{/\cal{E}}$, \ie leur produit fibr\'e
au-dessus de~$\cal{E}$ dans~$\Cat$,
\marginpar{153}
en tant que cat\'egorie au-dessus de~$\cal{E}$. Ainsi,
$\cal{F}\times_{\cal{E}} \cal{G}$ est muni de deux $\cal{E}$-foncteurs
$\pr_1$ et $\pr_2$, qui d\'efinissent, pour toute cat\'egorie
$\cal{H}$ au-dessus de~$\cal{E}$, une bijection
$$
\Hom_{\cal{E}}(\cal{H},\cal{F}\times_{\cal{E}}\cal{G})\isomto
\Hom_{\cal{E}}(\cal{H},\cal{F})\times\Hom_{\cal{E}}(\cal{H},\cal{G}).
$$
Cette bijection provient d'ailleurs d'un isomorphisme de
cat\'egories
$$
\SheafHom_{\cal{E}/-}(\cal{H},\cal{F}\times_{\cal{E}}
\cal{G})\isomto\SheafHom_{\cal{E}/-}(\cal{H},\cal{F})\times
\SheafHom_{\cal{E}/-}(\cal{H},\cal{G})
$$
en prenant les ensembles d'objets des deux membres, o\`u le foncteur
\'ecrit est celui qui a pour composantes les foncteurs $h\mto
\pr_1\circ h$ et $h\mto \pr_2\circ h$ du premier membre dans les deux
facteurs du second. On laisse au lecteur le soin de v\'erifier qu'on
obtient bien ainsi un isomorphisme (le fait analogue \'etant vrai,
plus g\'en\'eralement, chaque fois qu'on a une limite projective
de cat\'egories --- et non seulement dans le cas d'un produit
fibr\'e).

Rappelons d'ailleurs qu'on a (comme il a \'et\'e dit dans
le~\No \Ref{VI.1}):
\begin{eqnarray*}
\Ob(\cal{F}\times_\cal{E}\cal{G})&=&
\Ob(\cal{F})\times_{\Ob(\cal{E})}\Ob(\cal{G})\\
\Fl(\cal{F}\times_\cal{E}\cal{G})&=&
\Fl(\cal{F})\times_{\Fl(\cal{E})}\Fl(\cal{G})\quoi,
\end{eqnarray*}
la composition des fl\`eches se faisant d'ailleurs composante par
composante.

Dans la suite, nous consid\'erons un foncteur
$$
\lambda\colon\cal{E}'\to\cal{E}
$$
et pour toute cat\'egorie $\cal{F}$ au-dessus de~$\cal{E}$, on
consid\`ere $\cal{F}\times_\cal{E}\cal{E}'$ comme une cat\'egorie
au-dessus de~$\cal{E}'$ gr\^ace \`a $\pr_2$; en d'autres termes,
nous interpr\'etons l'op\'eration \og produit fibr\'e\fg comme une
op\'eration \emph{\og changement de base\fg,} le foncteur
$\lambda\colon\cal{E}'\to\cal{E}$ prenant le nom de \emph{\og foncteur
de changement de base\fg.}
\index{changement de base (foncteur)|hyperpage}%
\index{foncteur changement de base|hyperpage}%
Conform\'ement aux faits g\'en\'eraux bien connus, on obtient
ainsi un foncteur, dit
\marginpar{154}
\emph{foncteur changement de base} pour $\lambda$:
$$
\lambda^*\colon\Cat_{/\cal{E}}\to\Cat_{/\cal{E}'}\quoi,
$$
\label{indnot:fl}\oldindexnot{$\lambda^*\colon\Cat_{/\cal{E}}\to\Cat_{/\cal{E}'}$|hyperpage}%
(adjoint du foncteur \og restriction de la base\fg qui, \`a toute
cat\'egorie $\cal{F}'$ au-dessus de~$\cal{E}'$, de foncteur
projection~$p'$, associe~$\cal{F}'$, consid\'er\'e comme
cat\'egorie au-dessus de~$\cal{E}$ par le foncteur~$p=\lambda
p'$). Comme il est bien connu dans le cas g\'en\'eral d'un
foncteur changement de base dans une cat\'egorie, le foncteur
changement de base \og commute aux limites projectives\fg, et en
particulier \og transforme\fg produits fibr\'es sur~$\cal{E}$ en
produits fibr\'es sur~$\cal{E}'$.

Soient $\cal{F}$ et $\cal{G}$ deux cat\'egories au-dessus
de~$\cal{E}$, on va d\'efinir un \emph{isomorphisme canonique}:
\begin{multline*}
\label{eq:VI.3.i}
\tag{i}
\SheafHom_{\cal{E}'/-}(\cal{F}',\cal{G}')\isomto
\SheafHom_{\cal{E}/-}(\cal{F}\times_{\cal{E}}\cal{E}',\cal{G})\\
\text{o\`u }\cal{F}'=\cal{F}\times_{\cal{E}}\cal{E}'\text{, }
\cal{G}'=\cal{G}\times_{\cal{E}}\cal{E}'.
\end{multline*}
Pour ceci, consid\'erons le foncteur
$$
\pr_1\colon\cal{G}'=\cal{G}\times_{\cal{E}}\cal{E}'\to \cal{G}\quoi,
$$
et d\'efinissons~\eqref{eq:VI.3.i} par
$$
F\mto \pr_1\circ F\quoi,
$$
qui a priori d\'esigne un foncteur
\begin{equation*}
\label{eq:VI.3.ii}
\tag{ii} {\SheafHom(\cal{F}',\cal{G}')\to\SheafHom(\cal{F}',\cal{G})}
\end{equation*}
Il faut donc v\'erifier seulement que ce dernier induit un foncteur
pour les sous-cat\'egories~\eqref{eq:VI.3.i}, et que ce dernier est
un isomorphisme. Que~\eqref{eq:VI.3.ii} induise une bijection
$$
\Hom_{\cal{E}'/-}(\cal{F}',\cal{G}')\isomto
\Hom_{\cal{E}/-}(\cal{F}\times_{\cal{E}}\cal{E}',\cal{G})
$$
est la propri\'et\'e caract\'eristique du foncteur changement de
base. Il reste donc \`a
\marginpar{155}
prouver que si $F$, $G$ sont des $\cal{E}'$-foncteurs
$\cal{F}'\to\cal{G}'$, alors \emph{l'application}
$$
u\mto \pr_1\circ u
$$
\emph{induit une bijection}
$$
\Hom_{\cal{E}'}(F,G)\isomto\Hom_{\cal{E}}(\pr_1\circ F,\pr_1\circ
G)\quoi.
$$
La v\'erification de ce fait est imm\'ediate, et laiss\'ee au
lecteur.

Il r\'esulte de cet isomorphisme~\eqref{eq:VI.3.i}, et de la fin du
num\'ero pr\'ec\'edent, que
$$
\SheafHom_{\cal{E}'/-}(\cal{F}\times_{\cal{E}}\cal{E}',\cal{G}
\times_{\cal{E}}\cal{E}')
$$
\emph{peut \^etre consid\'er\'e comme un foncteur en}
$\cal{E}',\cal{F},\cal{G}$, \emph{de la cat\'egorie}
$\Cat_{/\cal{E}}^\circ\times\Cat_{/\cal{E}}^\circ\allowbreak\times\Cat_{/\cal{E}}$
\emph{dans la cat\'egorie}~$\Cat$, isomorphe au foncteur d\'efini
par l'expression
$\SheafHom_{\cal{E}/-}(\cal{F}\times_{\cal{E}}\cal{E}',\cal{G})$. En
particulier, pour $\cal{F}$,~$\cal{G}$ fix\'es, on obtient un
foncteur en $\cal{E}'$, $\cal{E}'\to
\SheafHom_{\cal{E}''/-}(\cal{F}',\cal{G}')=
\SheafHom_{\cal{E}'/-}(\cal{F}\times_{\cal{E}}
\cal{E}',\cal{G}\times_{\cal{E}}\cal{E}')$, et en particulier le
$\cal{E}$-foncteur de projection $\lambda\colon \cal{E}'\to\cal{E}$
d\'efinit un morphisme \ie un foncteur
$$
\lambda^*_{\cal{F},\cal{G}}\colon
\SheafHom_{\cal{E}/-}(\cal{F},\cal{G})\to
\SheafHom_{\cal{E}'/-}(\cal{F}',\cal{G}')
$$
que nous allons expliciter. Pour les ensembles d'objets des deux
membres, c'est l'application
$$
f\mto f'=f\times_{\cal{E}}\cal{E}'
$$
qui exprime la d\'ependance fonctorielle de
$\cal{F}\times_{\cal{E}}\cal{E}'$ de l'objet $\cal{F}$
sur~$\cal{E}$. D'autre part, consid\'erons deux $\cal{E}$-foncteurs
$$
f,g\colon\cal{F}\to \cal{G}
$$
et un homomorphisme de~$\cal{E}$-foncteurs
$$
u\colon f\to g\quoi,
$$
on
\marginpar{156}
va expliciter l'homomorphisme de~$\cal{E}'$-foncteurs correspondant:
$$
u'\colon f'\to g'\quoi.
$$
Pour tout
$$
\xi'=(\xi,S')\in\Ob(\cal{F}')
$$
avec
$$
\xi\in\Ob(\cal{F}),\quad S'\in\Ob(\cal{E}'),\quad p(\xi)=\lambda(S')=S
$$
le morphisme
$$
u'(\xi')\colon f'(\xi')=(f(\xi),S')\to g'(\xi')=(g(\xi),S')\quad
\text{dans } \cal{G}'
$$
est d\'efini par la formule
$$
u'(\xi')=(u(\xi),\id_{S'})
$$
(ce qui est bien un $S'$-morphisme dans~$\cal{G}'$, car
$q(u(\xi))=\lambda(\id_{S'})=\id_S$).

Consid\'erons maintenant un $\cal{E}$-foncteur quelconque
$$
\lambda'\colon\cal{E}''\to\cal{E}'
$$
et le foncteur correspondant
$$
\SheafHom_{\cal{E}'/-}(\cal{F}\times_{\cal{E}}\cal{E}',
\cal{G}\times_{\cal{E}}\cal{E}')\to
\SheafHom_{\cal{E}''/-}(\cal{F}\times_{\cal{E}}\cal{E}'',
\cal{G}\times_{\cal{E}}\cal{E}'')\quoi,
$$
je dis que ce foncteur n'est autre que le foncteur qu'on obtient par
le proc\'ed\'e pr\'ec\'edent, en partant de~$\cal{F}'$ et
$\cal{G}'$ sur~$\cal{E}'$ et en consid\'erant $\cal{E}''$ comme une
cat\'egorie sur~$\cal{E}'$, compte tenu des isomorphismes de
\emph{\og transitivit\'e de changement de base\fg}:
$$
\ifthenelse{\boolean{orig}}
{\cal{F}'\times_{\cal{E}'}\cal{E}''\isomto \cal{F}''=
\cal{F}\times_{\cal{E}}\cal{E}''\quad
\text{et}\quad\cal{G}'\times_{\cal{E}}\cal{E}''\isomto\cal{G}''
=\cal{G}\times_{\cal{E}}\cal{E}''\quoi,}
{\cal{F}'\times_{\cal{E}'}\cal{E}''\isomto
\cal{F}''=\cal{F}\times_{\cal{E}}\cal{E}''\quad
\text{et}\quad\cal{G}'\times_{\cal{E}'}\cal{E}''\isomto\cal{G}''
=\cal G\times_{\cal E}\cal E''\quoi,}
$$
qui impliquent un isomorphisme canonique
$$
\SheafHom_{\cal{E}''/-}(\cal{F}'\times_{\cal{E}'}\cal{E}'',
\cal{G}'\times_{\cal{E}'}\cal{E}'')\isomto
\SheafHom_{\cal{E}''/-}(\cal{F}\times_{\cal{E}}\cal{E}'',
\cal{G}\times_{\cal{E}}\cal{E}'')
$$
La
\marginpar{157}
v\'erification de cette compatibilit\'e est imm\'ediate, et
laiss\'ee au lecteur.

Les foncteurs qu'on vient de d\'efinir sont compatibles avec les
accouplements d\'efinis au num\'ero pr\'ec\'edent, de fa\c con pr\'ecise, si $\cal{F}$, $\cal{G}$,~$\cal{H}$ sont des
cat\'egories au-dessus de~$\cal{E}$ et si on pose
$$
\cal{F}'=\cal{F}\times_{\cal{E}}\cal{E}'\quoi,\quad \cal{G}'=
\cal{G}\times_{\cal{E}}\cal{E}'\quoi,\quad \cal{H}'=
\cal{H}\times_{\cal{E}}\cal{E}'\quoi,
$$
on a commutativit\'e dans le diagramme de foncteurs suivant:
$$
\xymatrix{ \SheafHom_{\cal{E}/-}(\cal{F},\cal{G})\times
\SheafHom_{\cal{E}/-}(\cal{G},\cal{H}) \ar[r]
\ar[d]_{\lambda^*_{\cal{F},\cal{G}}\times \lambda^*_{\cal{G},\cal{H}}}
& \SheafHom_{\cal{E}/-}(\cal{F},\cal{H})
\ar[d]^{\lambda^*_{\cal{F},\cal{H}}} \\
\ifthenelse{\boolean{orig}}
{\SheafHom_{\cal{E}/-}(\cal{F}',\cal{G}')\times}
{\SheafHom_{\cal{E}'/-}(\cal{F}',\cal{G}')\times}
\SheafHom_{\cal{E}'/-}(\cal{G}',\cal{H}')\ar[r] &
\SheafHom_{\cal{E}'/-}(\cal{F}',\cal{H}')}
$$
o\`u les fl\`eches horizontales sont les foncteurs-composition
d\'efinis au num\'ero pr\'ec\'edent. Cette commutativit\'e
s'exprime par les formules
$$
(gf)'=g'f'
$$
pour $f\in\Hom_{\cal{E}}(\cal{F},\cal{G})$,
$g\in\Hom_{\cal{E}}(\cal{G},\cal{H})$, (formule qui exprime simplement
la fonctorialit\'e du changement de base), et la formule
$$
(v*u)'=v'*u'
$$
lorsque $u\colon f\to f_1$ est une fl\`eche de
$\SheafHom_{\cal{E}/-}(\cal{F},\cal{G})$ et $v\colon g\to g_1$ une
fl\`eche de~$\SheafHom_{\cal{E}/-}(\cal{G},\cal{H})$. La
v\'erification de cette formule r\'esulte facilement des
d\'efinitions.

Dans la suite, nous nous int\'eresserons surtout \`a
$\SheafHom_{\cal{E}}(\cal{F},\cal{G})$ (et certaines
sous-cat\'egories remarquables de celle-ci) lorsque
$\cal{F}=\cal{E}$, et introduisons pour cette raisons une notation
sp\'eciale:
$$
\bf{\Gamma}(\cal{G}/\cal{E})=
\SheafHom_{\cal{E}}(\cal{E},\cal{G})\quoi,\quad
\Gamma(\cal{G}/\cal{E})=\Ob(\bf{\Gamma}(\cal{G}/\cal{E}))=
\ifthenelse{\boolean{orig}}
{\Hom_{\cal{E}}(\cal{G},\cal{E})\quoi.}
{\Hom_{\cal{E}}(\cal{E},\cal{G})\quoi.}
$$
\label{indnot:fm}\oldindexnot{$\bf{\Gamma}(\cal{G}/\cal{E})$ et $\Gamma(\cal{G}/\cal{E})$|hyperpage}%

\begin{remarquesstar}
Lorsque
\marginpar{158}
$\cal{E}$ est une cat\'egorie ponctuelle,
\ie $\Ob(\cal{E})$ et $\Fl(\cal{E})$ r\'eduits \`a un seul
\'el\'ement, ce qui signifie aussi que $\cal{E}$ est un objet
final de la cat\'egorie~$\Cat$, alors la donn\'ee d'une
cat\'egorie sur~$\cal{E}$ est \'equivalente \`a la donn\'ee
d'une cat\'egorie tout court, (car il y aura un foncteur unique de
$\cal{F}$ dans~$\cal{E}$). De fa\c con plus pr\'ecise,
$\Cat_{/\cal{E}}$ est alors isomorphe \`a~$\Cat$. De plus, les
cat\'egories $\SheafHom_{\cal{E}/-}(\cal{F},\cal{G})$ ne sont alors
autres que les $\SheafHom(\cal{F},\cal{G})$. Rappelons alors que la
formule fondamentale
$$
\Hom(\cal{H},\SheafHom(\cal{F},\cal{G}))\isomto
\Hom(\cal{F}\times\cal{H},\cal{G})
$$
(isomorphisme fonctoriel en les trois arguments qui y figurent),
permet d'interpr\'eter $\SheafHom(\cal{F},\cal{G})$ axiomatiquement,
en termes internes \`a la cat\'egorie~$\Cat$, de sorte que le
formulaire connu des cat\'egories $\SheafHom$ appara\^it comme un
cas particulier d'un formulaire valable dans les cat\'egories telles
que~$\Cat$, o\`u des \og objets $\SheafHom$\fg (d\'efinis par la
formule pr\'ec\'edente) existent. Il y a une interpr\'etation
analogue de~$\SheafHom_{\cal{E}/-}(\cal{F},\cal{G})$ lorsqu'on suppose
\`a nouveau $\cal{E}$ quelconque, par la formule
$$
\Hom(\cal{H},\SheafHom_{\cal{E}/-}(\cal{F},\cal{G}))\isomto
\Hom_{\cal{E}}(\cal{F}\times\cal{H},\cal{G})
$$
(isomorphisme fonctoriel en les trois arguments). De cette fa\c con,
les propri\'et\'es formelles expos\'ees dans les \No \Ref{VI.2},~\Ref{VI.3} sont
des cas particuliers de r\'esultats plus g\'en\'eraux, valables
dans les cat\'egories o\`u les objets
$\SheafHom_{\cal{E}/-}(\cal{F},\cal{G})$ (lorsque $\cal{F}$,~$\cal{G}$
sont deux objets de la cat\'egorie au-dessus d'un troisi\`eme~$\cal{E}$) existent.
\end{remarquesstar}

\section[Cat\'egories-fibres]{Cat\'egories-fibres; \'equivalence de
$\mathcal{E}$-cat\'egories}
\label{VI.4}

Soit~$\cal{F}$ une cat\'egorie sur~$\cal{E}$, et soit $S \in
\Ob({\cal{E}})$. On appelle \emph{cat\'egorie-fibre de}~$\cal{F}$
en~$S$
\index{cat\'egorie-fibre|hyperpage}%
la sous-cat\'egorie~$\cal{F}_S$
\label{indnot:fn}\oldindexnot{$\cal{F}_S$|hyperpage}%
de~$\cal{F}$ image r\'eciproque de la
\marginpar{159}
sous-cat\'egorie ponctuelle de~$\cal{E}$ d\'efinie par $S$. Donc
les objets de~$\cal{F}_S$ sont les objets~$\xi$ de~$\cal{F}$ tels que
$p(\xi)=S$, ses morphismes sont les morphismes~$u$ de~$\cal{F}$ tels
que $p(u)=\id_S$, \ie les $S$-morphismes dans $\cal{F}$. Bien
entendu,~$\cal{F}_S$ est canoniquement isomorphe au produit fibr\'e
$\cal{F} \times_{\cal{E}} \{S\}$, o\`u~$\{S\}$ d\'esigne la
sous-cat\'egorie ponctuelle de~$\cal{E}$ d\'efinie par~$S$, munie
de son foncteur d'inclusion dans~$\cal{E}$. Il en r\'esulte (compte
tenu de la transitivit\'e du changement de base) que si on fait le
changement de base $\lambda\colon \cal{E}' \to \cal{E}$, alors pour
tout $S' \in \Ob(\cal{E}')$, \emph{la projection $\pr_1\colon
\cal{F}' = \cal{F} \times_{\cal{E}} \cal{E}' \to \cal{F}$ induit un
isomorphisme}
$$
\cal{F}'_{S'} \to \cal{F}_S \qquad (\text{o\`u $S=\lambda(S')$})%
\ifthenelse{\boolean{orig}}{}{.}
$$

\begin{proposition}
\label{VI.4.1} Soit $f\colon \cal{F} \to \cal{G}$ un
$\cal{E}$-foncteur. Si~$f$ est pleinement fid\`ele, alors pour tout
changement de base $\cal{E}' \to \cal{E}$, le foncteur correspondant
$f' \colon \cal{F}' = \cal{F} \times_{\cal{E}} \cal{E}' \to\allowbreak \cal{G}' =
\cal{G} \times_{\cal{E}} \cal{E}'$ est pleinement fid\`ele.
\end{proposition}

La v\'erification est imm\'ediate; plus g\'en\'eralement, on
peut montrer que toute limite projective de
\ifthenelse{\boolean{orig}}{foncteur}{foncteurs}
pleinement
\ifthenelse{\boolean{orig}}{fid\`ele}{fid\`eles}
(ici,~$f$ et les foncteurs identiques
dans~$\cal{E},\cal{E}'$) est un foncteur pleinement fid\`ele.

On notera que l'assertion analogue \`a 4.1,
\ifthenelse{\boolean{orig}}{ou}{o\`u}
\og pleinement fid\`ele\fg est remplac\'ee par \og \'equi\-va\-lence
de cat\'egories\fg, est fausse, d\'ej\`a pour~$\cal{G}=\cal{E}$.
Cependant:

\begin{proposition}
\label{VI.4.2} Soit $f \colon \cal{F} \to \cal{G}$ un
$\cal{E}$-foncteur. Les conditions suivantes sont \'equivalentes:
\begin{enumerate}
\item[(i)] Il existe un $\cal{E}$-foncteur $g \colon \cal{G} \to
\cal{F}$ et des $\cal{E}$-isomorphismes
$$
gf \isomto \id_{\cal{F}}, \quad fg \isomto \id_{\cal{G}}
\ifthenelse{\boolean{orig}}{}{.}
$$
\item[(ii)] Pout toute cat\'egorie $\cal{E}'$ sur~$\cal{E}$, le
foncteur
$$
f' = f \times_{\cal{E}} \cal{E}' \colon \cal{F}' = \cal{F}
\times_{\cal{E}} \cal{E}' \to \cal{G}'=\cal{G} \times_{\cal{E}}
\cal{E}'
$$
est une \'equivalence de cat\'egories.
\item[(iii)]
\marginpar{160}
$f$ est une \'equivalence de cat\'egories, et pout tout $S \in
\Ob(\cal{E})$, le foncteur $f_S \colon \cal{F}_S \to \cal{G}_S$ induit
par~$f$ est une \'equivalence de cat\'egories.
\item[(iii~bis)] $f$ est pleinement fid\`ele, et pour tout $S \in
\Ob(\cal{E})$ et tout
\ifthenelse{\boolean{orig}}
{$\eta\in\Ob \cal{G}_S$,}
{$\eta\in\Ob(\cal{G}_S)$,}
il existe un $\xi \in \Ob(\cal{F}_S)$ et un $S$-isomorphisme
\ifthenelse{\boolean{orig}}
{$u \colon f((\xi) \to \eta$.}
{$u \colon f(\xi) \to \eta$.}
\end{enumerate}
\end{proposition}

\subsubsection*{D\'emonstration} \'Evidemment (i) implique que~$f$ est une
\'equivalence de cat\'egories (notion qui se d\'efinit par la
m\^eme condition, mais o\`u les isomorphismes de foncteurs ne sont
pas astreints \`a \^etre des $\cal{E}$-morphismes). D'autre part,
il r\'esulte des fonctorialit\'es du num\'ero pr\'ec\'edent
que la condition (i) est conserv\'ee apr\`es changement de
base~$\cal{E}'\to\cal{E}$. Il s'ensuit que (i)$\To$(ii).
\'Evidemment (ii)$\To$(iii), car il suffit de faire $\cal{E}'
= \cal{E}$ et~$\cal{E}'=\{S\}$. Il est encore plus trivial que
(iii)$\To$(iii~bis), reste \`a prouver que (iii~bis)$\To$(i). Pour ceci, choisissons pour tout $\eta \in
\Ob(\cal{G})$ un $g(\eta) \in \Ob(\cal{F})$ et un isomorphisme
$u(\eta)\colon f(g(\eta))\to\eta$ qui soit tel que~$q(u(\eta)) =
\id_S$, o\`u~$S=q(\eta)$. C'est possible gr\^ace \`a la
deuxi\`eme condition (iii~bis). Comme il est connu et imm\'ediat,
le fait que~$f$ est pleinement fid\`ele implique que~$g$ peut de
fa\c con unique \^etre consid\'er\'e comme un foncteur
de~$\cal{G}$ dans~$\cal{F}$, de fa\c con que les~$u(\eta)$
d\'efinissent un homomorphisme (donc un isomorphisme) fonctoriel
$u\colon fg \isomto \id_{\cal{G}}$. De plus, par construction~$g$ est
un $\cal{E}$-foncteur et~$u$ un $\cal{E}$-homomorphisme. Aux
donn\'ees pr\'ec\'edentes correspond alors un isomorphisme
fonctoriel $v\colon gf \to \id_{\cal{F}}$, d\'efini par la condition
que $f*v=u*f$, et on constate tout de suite que c'est \'egalement un
$\cal{E}$-morphisme, cqfd.
\begin{definition}
\label{VI.4.3} Si les conditions pr\'ec\'edentes sont
v\'erifi\'ees, on dit que~$f$ est une \'equivalence de
cat\'egories sur~$\cal{E}$, ou une $\cal{E}$-\'equivalence.
\end{definition}

\begin{corollaire}
\label{VI.4.4} Supposons que le foncteur projection $p\colon
\cal{F}\to\cal{E}$ soit un foncteur transportable, \ie que pour tout
isomorphisme $\alpha\colon T \to S$ dans~$\cal{E}$ et tout objet~$\xi$
dans~$\cal{F}_T$, il existe un isomorphisme~$u$ dans~$\cal{F}$ de
source~$\xi$ tel que~$p(u)=\alpha$. Alors tout $\cal{E}$-foncteur
$f\colon \cal{F}\to\cal{G}$ qui est une \'equivalence de
cat\'egories, est une $\cal{E}$-\'equivalence.
\end{corollaire}

R\'esulte
\marginpar{161}
du crit\`ere (iii~bis).

\begin{corollaire}
\label{VI.4.5} Soit $f\colon \cal{F}\to\cal{G}$ une
$\cal{E}$-\'equivalence. Alors pour toute cat\'egorie~$\cal{H}$
sur~$\cal{E}$, les foncteurs correspondants:
\begin{align*}
\SheafHom_{\cal{E}/-}(\cal{G},\cal{H}) & \to
\SheafHom_{\cal{E}/-}(\cal{F},\cal{H}) \\
\SheafHom_{\cal{E}/-}(\cal{H},\cal{F})
& \to \SheafHom_{\cal{E}/-}(\cal{H},\cal{G})
\end{align*}
(\cf \No \Ref{VI.2}) sont des \'equivalences de cat\'egories.
\end{corollaire}

Cela r\'esulte du crit\`ere (i) par le raisonnement habituel.

\section{Morphismes cart\'esiens, images inverses, foncteurs
car\-t\'e\-siens}
\label{VI.5}

Soit~$\cal{F}$ une cat\'egorie sur~$\cal{E}$, de
foncteur-projection~$p$.
\ifthenelse{\boolean{orig}}
{}
{\enlargethispage{.5cm}}

\begin{definition}
\label{VI.5.1} Consid\'erons un morphisme
$$
\alpha\colon \eta\to\xi
$$
dans~$\cal{F}$, et soient
$$
S = p(\xi), \quad T = p(\eta), \quad f=p(\alpha)%
\ifthenelse{\boolean{orig}}{}{.}
$$
On dit que~$\alpha$ est un \emph{morphisme cart\'esien}
\index{cart\'esien (morphisme)|hyperpage}%
si pour tout $\eta' \in \Ob(\cal{F}_T)$ et tout $f$-morphisme $u\colon
\eta'\to\xi$, il existe un $T$-morphisme unique $\overline{u}\colon
\eta'\to\eta$ tel que~$u=\alpha\circ\overline{u}$.
\end{definition}

Cela signifie donc que pour tout $\eta'\in\Ob(\cal{F}_T)$,
l'application $v \mto \alpha\circ v$:
\begin{equation}
\tag{i} \Hom_T(\eta',\eta) \to \Hom_f(\eta',\xi)
\end{equation}
est bijective. Cela signifie aussi que le couple $(\eta,\alpha)$
\emph{repr\'esente le foncteur en~$\eta'$} $\cal{F}^{\circ}_T \to
\Ens$ du deuxi\`eme membre. Si pour un morphisme $f\colon T \to S$
dans~$\cal{E}$ donn\'e, et un $\xi \in \Ob(\cal{F}_S)$ donn\'e, il
existe un tel couple $(\eta,\alpha)$, \ie un morphisme
\marginpar{162}
cart\'esien~$\alpha$ dans~$\cal{F}$ de but~$\xi$, tel que
$p(\alpha)=f$, alors~$\eta$ est d\'etermin\'e dans $\cal{F}_T$
\`a isomorphisme unique pr\`es. On dit alors que \emph{l'image
inverse de~$\xi$ par~$f$ existe}, et un objet~$\eta$ de~$\cal{F}_T$
muni d'un $f$-morphisme cart\'esien $\alpha\colon \eta \to \xi$ est
appel\'e \emph{une image inverse de~$\xi$ par~$f$}.
\index{image inverse|hyperpage}%
Souvent on suppose
\ifthenelse{\boolean{orig}}{choisi}{choisie}
une telle image inverse chaque fois qu'elle existe ($\cal{F}$
\'etant fix\'e); on notera alors l'image inverse par des symboles
tels que~$f^*_{\cal{F}}(\xi)$, ou simplement~$f^*(\xi)$
\label{indnot:fo}\oldindexnot{$f^*_{\cal{F}}(\xi)$ ou $f^*(\xi)$|hyperpage}%
\ifthenelse{\boolean{orig}}{o\`u~$\xi \times_S T$}
{ou~$\xi \times_S T$}
lorsque ces notations n'entra\^inent pas
\ifthenelse{\boolean{orig}}{des confusions;}{de confusions;}
le morphisme canonique $\alpha\colon \eta \to \xi$ sera alors
not\'e, dans ce qui suit, par~$\alpha_f(\xi)$.
\label{indnot:fp}\oldindexnot{$\alpha_f(\xi)$|hyperpage}%
Si pour tout $\xi \in \Ob(\cal{F}_S)$, l'image inverse de~$\xi$
par~$f$ existe, on dira aussi que \emph{le foncteur image inverse
par~$f$ dans~$\cal{F}$ existe}, et~$f^*(\xi)$ devient alors un
\emph{foncteur covariant en~$\xi$}, de~$\cal{F}_S$ dans~$\cal{F}_T$.
Ceci provient du fait que le deuxi\`eme membre dans (i) d\'epend
de fa\c con covariante de~$\xi$, \ie de fa\c con pr\'ecise
d\'esigne un foncteur de~$\cal{F}^{\circ}_T \times \cal{F}_S$
dans~$\Ens$. Cette d\'ependance fonctorielle pour~$f^*(\xi)$
s'explicite ainsi: consid\'erons des $f$-morphismes cart\'esiens
$$
\alpha\colon \eta \to \xi, \quad \alpha'\colon \eta' \to \xi'
$$
et un $S$-morphisme $\lambda\colon \xi \to \xi'$, alors il existe un
$T$-morphisme et un seul $\mu\colon \eta\to\eta'$ tel que l'on ait
$$
\alpha' \mu = \lambda \alpha
$$
(comme il r\'esulte du fait que~$\alpha'$ est cart\'esien).

Notons aussi le fait imm\'ediat suivant: consid\'erons un
diagramme commutatif
$$
\xymatrix@C=1.2cm{
\xi\ar[d]_-{\lambda}& \ar[l]_-{\alpha}\eta \ar[d]^-{\mu}\\
\xi' &\ar[l]_-{\alpha'}\eta'
}
$$
dans~$\cal{F}$,
\marginpar{163}
o\`u~$\alpha$ et~$\alpha'$ sont des $f$-morphismes, et~$\lambda$ un
$S$-isomorphisme,~$\mu$ un $T$-isomorphisme. \emph{Pour que~$\alpha$
soit cart\'esien, il faut et il suffit que~$\alpha'$ le
soit}.

\begin{definition}
\label{VI.5.2} Un $\cal{E}$-foncteur $F\colon \cal{F}\to\cal{G}$ est
appel\'e un \emph{foncteur cart\'esien}
\index{cart\'esien (foncteur)|hyperpage}%
s'il transforme morphismes cart\'esiens en morphismes
cart\'esiens. On d\'esigne par $\SheafHom_\cart(\cal{F},\cal{G})$
\label{indnot:fq}\oldindexnot{$\SheafHom_\cart(\cal{F},\cal{G})$|hyperpage}%
la sous-cat\'egorie pleine de
$\SheafHom_{\cal{E}/-}(\cal{F},\cal{G})$ form\'ee des foncteurs
cart\'esiens.

Par exemple, consid\'erant~$\cal{E}$ comme une cat\'egorie
sur~$\cal{E}$ gr\^ace au foncteur identique, tout morphisme
de~$\cal{E}$ est cart\'esien, donc un foncteur cart\'esien
de~$\cal{E}$ dans~$\cal{F}$ est un foncteur section $F\colon
\cal{E}\to\cal{F}$ qui transforme tout morphisme de~$\cal{E}$ en un
morphisme cart\'esien; un tel foncteur s'appelle une \emph{section
cart\'esienne} de~$\cal{F}$ sur~$\cal{E}$.
\end{definition}

\begin{proposition}
\label{VI.5.3} \textup{(i)}~Un foncteur $F\colon \cal{F} \to \cal{G}$ qui est une
$\cal{E}$-\'equivalence, est un foncteur cart\'esien.
\textup{(ii)}~Soient~$F,G$ deux $\cal{E}$-foncteurs
\emph{isomorphes}~$\cal{F}\to\cal{G}$. Si l'un est cart\'esien,
l'autre l'est. \textup{(iii)}~Le compos\'e de deux foncteurs cart\'esiens
$\cal{F}\to\cal{G}$ et $\cal{G}\to\cal{H}$ est un foncteur
cart\'esien.
\end{proposition}

L'assertion (iii) est triviale sur la d\'efinition, (ii)
r\'esulte de la remarque pr\'ec\'edant~\Ref{VI.5.2}, (i) r\'esulte
facilement de la d\'efinition et du crit\`ere~\Ref{VI.4.2}~(iii); plus
pr\'ecis\'ement, un morphisme~$\alpha$ dans~$\cal{F}$ est
cart\'esien si et seulement si~$F(\alpha)$ l'est.

\begin{corollaire}
\label{VI.5.4}
Soit $F \colon \cal{F} \to \cal{G}$ une
$\cal{E}$-\'equivalence. Alors pour toute cat\'egorie~$\cal{H}$
sur~$\cal{E}$, les foncteurs correspondants $G \mto G \circ F$ et
$G \mto F \circ G$ induisent des \'equivalences de
cat\'egories:
\begin{eqnarray*}
\SheafHom_\cart(\cal{G},\cal{H})\lto{\approx}
\SheafHom_\cart(\cal{F},\cal{H}) \\
\SheafHom_\cart(\cal{H},\cal{F})\lto{\approx}
\SheafHom_\cart(\cal{H}, \cal{G})
\end{eqnarray*}
\end{corollaire}

Cela se d\'eduit de la fa\c con habituelle de~\Ref{VI.4.2}
crit\`ere~(i) et de~\Ref{VI.5.3} (i) (ii)~(iii).
\marginpar{164}
On peut pr\'eciser que \emph{le $\cal{E}$-foncteur
$G\colon \cal{G} \to \cal{H}$ est cart\'esien si et seulement si $G
\circ F$ l'est}, et de m\^eme \emph{un $\cal{E}$-foncteur $G\colon
\cal{H} \to \cal{F}$ est cart\'esien si et seulement si $F \circ G$
l'est}.

Il r\'esulte de~\Ref{VI.5.4}~(iii) que si on consid\`ere la
sous-cat\'egorie~$\Cat^\cart_{/\cal{E}}$
\label{indnot:fr}\oldindexnot{$\Cat^\cart_{/\cal{E}}$|hyperpage}%
de~$\Cat_{/\cal{E}}$ dont les objets sont les m\^emes que ceux de
$\Cat_{/\cal{E}}$, et dont les morphismes sont les foncteurs
\emph{cart\'esiens} alors on a comme au~\No \Ref{VI.2} des accouplements:
$$
\SheafHom_\cart(\cal{F},\cal{G})\times
\SheafHom_\cart(\cal{G},\cal{H}) \to \SheafHom_\cart(\cal{F},\cal{H})
$$
\ifthenelse{\boolean{orig}}{induit}{induits}
par ceux du~\No \Ref{VI.2}, permettant de consid\'erer
$\SheafHom_\cart(\cal{F},\cal{G})$ comme un foncteur en $\cal{F}$,
$\cal{G}$, de la cat\'egorie $\left(\Cat^\cart_{/\cal{E}}
\right)^\circ \times \Cat^\cart_{/\cal{E}}$
\ifthenelse{\boolean{orig}}{dans~$\Cat$).}{dans~$\Cat$.}
Nous aurons besoin de cette remarque surtout pour le cas o\`u
$\cal{F} = \cal{G}$:

\begin{definition}
\label{VI.5.5}
Soit $\cal{F}$ une cat\'egorie sur~$\cal{E}$. On d\'esigne par
$$
\varprojLim\cal{F}/\cal{E}
$$
\label{indnot:fs}\oldindexnot{$\varprojLim\cal{F}/\cal{E}$|hyperpage}%
la cat\'egorie des $\cal{E}$-foncteurs cart\'esiens
$\cal{E}\to\cal{F}$, \ie des sections cart\'esiennes de~$\cal{F}$
sur~$\cal{E}$.
\end{definition}

D'apr\`es ce qu'on vient de dire,
$\varprojLim\cal{F}/\cal{E}$ est un foncteur
en~$\cal{F}$, de la cat\'egorie $\Cat^\cart_{/\cal{E}}$ dans la
cat\'egorie $\Cat$.

Nous verrons plus bas les relations entre cette op\'eration
$\varprojLim$ et la notion de limite projective de
cat\'egories, ainsi que de nombreux exemples.

\section[Cat\'egories fibr\'ees et cat\'egories
pr\'efibr\'ees]{Cat\'egories fibr\'ees et cat\'egories
pr\'efibr\'ees. Produits et changement de base dans icelles}
\label{VI.6}

\begin{definition}
\label{VI.6.1}
Une cat\'egorie $\cal{F}$ sur~$\cal{E}$ est appel\'ee une
\emph{cat\'egorie fibr\'ee}
\index{fibree (categorie)@fibr\'ee (cat\'egorie)|hyperpage}%
(et on dit alors que le foncteur $\cal{F}\to\cal{E}$ est
\emph{fibrant})
\index{fibrant (foncteur)|hyperpage}%
si elle satisfait les deux axiomes suivants:
\begin{itemize}
\item{$\Fib_\rm{I}$}
\marginpar{165}
Pour tout morphisme $f\colon T\to S$ dans~$\cal{E}$, le foncteur image
inverse par $f$ dans $\cal{F}$ existe.
\item{$\Fib_\rm{II}$} Le compos\'e de deux morphismes cart\'esiens
est cart\'esien.
\end{itemize}
Une cat\'egorie $\cal{F}$ sur~$\cal{E}$ satisfaisant la condition
$\Fib_\rm{I}$ est appel\'ee une \emph{cat\'egorie
pr\'efibr\'ee sur~$\cal{E}$}.
\index{prefibree (categorie)@pr\'efibr\'ee (cat\'egorie)|hyperpage}%
\end{definition}

Si~$\cal{F}$ est une cat\'egorie fibr\'ee
(\resp pr\'efibr\'ee) sur~$\cal{E}$, une sous-cat\'egorie
$\cal{G}$ de~$\cal{F}$ est appel\'ee une \emph{sous-cat\'egorie
fibr\'ee}
\index{sous-cat\'egorie fibr\'ee (\resp pr\'efibr\'ee)|hyperpage}%
(\resp une \emph{sous-cat\'egorie pr\'efibr\'ee}) si c'est une
cat\'egorie fibr\'ee (\resp pr\'efibr\'ee) sur~$\cal{E}$, et
si de plus le foncteur d'inclusion est cart\'esien. Si par exemple
$\cal{G}$ est une sous-cat\'egorie \emph{pleine} de~$\cal{F}$, on
voit que cela signifie que pour tout morphisme $f\colon T\to S$ dans
$\cal{E}$ et pour tout $\xi \in \Ob(\cal{G}_S)$, $f^*_{\cal{F}}(\xi)$
est $T$-isomorphe \`a un objet de~$\cal{G}_T$. Un autre cas
int\'eressant est le suivant: $\cal{F}$ \'etant une cat\'egorie
fibr\'ee sur~$\cal{E}$, consid\'erons la sous-cat\'egorie
$\cal{G}$ de~$\cal{F}$ ayant m\^emes objets, et dont les morphismes
sont les morphismes \emph{cart\'esiens} de~$\cal{F}$; en particulier
les morphismes de~$\cal{G}_S$ sont les isomorphismes de
$\cal{F}_S$. On voit de suite que c'est bien une sous-cat\'egorie
fibr\'ee de~$\cal{F}$, car dans la bijection
$$
\Hom_T(\eta',\eta) \isomto \Hom_f(\eta',\xi)
$$
relative \`a un $f$-morphisme cart\'esien $\alpha$ dans $\cal{F}$,
aux $T$-isomorphismes du premier membre correspondent les morphismes
cart\'esiens du second. Par d\'efinition, les sections
cart\'esiennes $\cal{E}\to\cal{F}$
\ifthenelse{\boolean{orig}}{(correspondent}{correspondent}
alors biunivoquement aux $\cal{E}$-foncteurs quelconques
$\cal{E}\to\cal{G}$ (mais on notera que le foncteur naturel
$$
\SheafHom_{\cal{E}/-}(\cal{E},\cal{G})\to
\SheafHom_\cart(\cal{E},\cal{F})=
\varprojLim(\cal{F}/\cal{E})
$$
est fid\`ele, mais en g\'en\'eral n'est pas pleinement
fid\`ele, \ie n'est pas un isomorphisme).

\begin{remarquesstar}
Soit $\cal{F}$ une cat\'egorie sur~$\cal{E}$. Les conditions
suivantes sont \'equivalentes:~(i)~Tous les morphismes de~$\cal{F}$
sont cart\'esiens (ii)~$\cal{F}$~est une cat\'egorie fibr\'ee
sur~$\cal{E}$, et les $\cal{F}_S$ sont des groupo\"ides%
\ifthenelse{\boolean{orig}}{,}{}
(\ie tout morphisme dans $\cal{F}_S$ est un isomorphisme). On dit
alors que $\cal{F}$ est une cat\'egorie \emph{fibr\'ee en
groupo\"ides}
\index{fibree en groupoides@fibr\'ee en groupo\"ides|hyperpage}%
%
\marginpar{166}
sur~$\cal{E}$. Ce sont elles qu'on rencontre surtout en \og th\'eorie
des modules\fg. Si $\cal{E}$ est un groupo\"ide, on montre que les
conditions (i) et (ii) \'equivalent aussi \`a la suivante:
(iii)~$\cal{F}$ est un groupo\"ide, et le foncteur projection
$p\colon \cal{F}\to\cal{E}$ est transportable (\cf \Ref{VI.4.4}). Par
exemple, si $\cal{E}$ et $\cal{F}$ sont des groupo\"ides tels que
\ifthenelse{\boolean{orig}}
{$\Ob \cal{E}$ et $\Ob \cal{F}$}
{$\Ob(\cal{E})$ et $\Ob(\cal{F})$}
soient r\'eduits \`a un point, de sorte que $\cal{E}$ et $\cal{F}$
sont d\'efinis, \`a isomorphisme pr\`es, par des groupes $E$ et
$F$, et le foncteur $p\colon \cal{F}\to\cal{E}$ est d\'efini par un
homomorphisme de groupes $p\colon F\to E$, alors $\cal{F}$ est
fibr\'e sur~$\cal{E}$ si et seulement si $p$ est surjectif, \ie si
$p$ d\'efinit une \emph{extension} du groupe $E$ par le groupe $G =
\Ker p$.
\end{remarquesstar}

\begin{proposition}
\label{VI.6.2}
Soit $F \colon \cal{F} \to \cal{G}$ une
$\cal{E}$-\'equivalence. Pour que $\cal{F}$ soit une cat\'egorie
fibr\'ee (\resp pr\'efibr\'ee) sur~$\cal{E}$, il faut et il
suffit que $\cal{G}$ le soit.
\end{proposition}

R\'esulte facilement des d\'efinitions et de la remarque
signal\'ee plus haut qu'un morphisme
\ifthenelse{\boolean{orig}}{$\alpha x$}{$\alpha$}
dans $\cal{F}$ est cart\'esien si et seulement si $F(\alpha)$ l'est.

\begin{proposition}
\label{VI.6.3}
Soient $\cal{F}_1$, $\cal{F}_2$ deux cat\'egories sur~$\cal{E}$, et
soit $\alpha = (\alpha_1,\alpha_2)$ un morphisme dans
$\cal{F}=\cal{F}_1 \times_{\cal{E}} \cal{F}_2$. Pour que $\alpha$ soit
cart\'esien, il faut et il suffit que ses composantes le soient.
\end{proposition}

Soit en effet $\xi_i$ le but et $\eta_i$ la source de~$\alpha_i$, et
soit $f\colon T\to S$ le morphisme de~$\cal{E}$ tel que $\alpha_1$ et
$\alpha_2$ soient des $f$-morphismes. Pour tout
$\eta'=(\eta'_1,\eta'_2)$ dans $\cal{F}_T$, on a un diagramme
commutatif
$$
\xymatrix{ \Hom_T(\eta',\eta) \ar[r] \ar[d] & \Hom_f(\eta',\xi) \ar[d]
\\ \Hom_T(\eta'_1,\eta_1)\times\Hom_T(\eta'_2,\eta_2) \ar[r] &
\Hom_f(\eta'_1,\xi_1)\times \Hom_f(\eta'_2,\xi_2) }
$$
o\`u les fl\`eches verticales sont des bijections. Donc si l'une
des fl\`eches horizontales est une bijection, il en est de m\^eme
de l'autre. Cela montre d\'ej\`a que si $\alpha_1$, $\alpha_2$
sont cart\'esiens (donc la deuxi\`eme fl\`eche horizontale
bijective) alors $\alpha$ l'est. La r\'eciproque se voit en faisant
dans le diagramme ci-dessus $\eta'_i=\eta_i$ d'o\`u
$\Hom_T(\eta'_i,\eta_i) \neq \emptyset$, d'abord pour $i=2$ ce qui
prouve que $\alpha_1$ est cart\'esien, puis pour $i=1$ ce qui prouve
que $\alpha_2$ est cart\'esien.

\begin{corollaire}
\label{VI.6.4}
Soit
\marginpar{167}
$\cal{F} = \cal{F}_1\times_{\cal{E}}\cal{F}_2$, et soit
$F=(F_1,F_2)$ un $\cal{E}$-foncteur $\cal{G}\to\cal{F}$. Pour que
\ifthenelse{\boolean{orig}}{$\cal{F}$}{$F$}
soit cart\'esien, il faut et il suffit que $F_1$ et $F_2$ le
soient. On obtient ainsi un isomorphisme de cat\'egories
$$
\SheafHom_\cart(\cal{G}, \cal{F}_1 \times_{\cal{E}} \cal{F}_2) \isomto
\SheafHom_\cart(\cal{G},\cal{F}_1)
\times\SheafHom_\cart(\cal{G},\cal{F}_2)
$$
et en particulier (faisant $\cal{G}=\cal{E}$) un isomorphisme de
cat\'egories
$$
\varprojLim(\cal{F}_1 \times_{\cal{E}}
\cal{F}_2/\cal{E}) \isomto
\varprojLim(\cal{F}_1/\cal{E}) \times
\varprojLim(\cal{F}_2/\cal{E})
$$
\end{corollaire}

\begin{corollaire}
\label{VI.6.5}
Soient $\cal{F}_1$ et $\cal{F}_2$ deux cat\'egories fibr\'ees
(\resp pr\'efibr\'ees) au-dessus de~$\cal{E}$, alors leur produit
fibr\'e $\cal{F}=\cal{F}_1 \times_{\cal{E}}\cal{F}_2$ est une
cat\'egorie fibr\'ee (\resp pr\'efibr\'ee) sur~$\cal{E}$.
\end{corollaire}

Ces r\'esultats s'\'etendent d'ailleurs au cas du produit
fibr\'e d'une famille quelconque de cat\'egories sur~$\cal{E}$.

\begin{proposition}
\label{VI.6.6}
Soient $\cal{F}$ une cat\'egorie sur~$\cal{E}$, de
foncteur-projection $p$, et soit $\lambda \colon \cal{E}' \to \cal{E}$
un foncteur, consid\'erons
$\cal{F}'=\cal{F}\times_{\cal{E}}\cal{E}'$ comme une cat\'egorie sur~$\cal{E}'$ par le foncteur-projection $p'=p \times_{\cal{E}}
\id_{\cal{E}'}$. Soit $\alpha'$ un morphisme de~$\cal{F}'$, pour que
$\alpha'$ soit un morphisme cart\'esien, il faut et il suffit que
son image $\alpha$ dans $\cal{F}$ le soit.
\end{proposition}

La d\'emonstration est imm\'ediate et laiss\'ee au lecteur.

\begin{corollaire}
\label{VI.6.7}
Pour tout foncteur cart\'esien $F\colon \cal{F} \to \cal{G}$ de
cat\'egories sur~$\cal{E}$, le foncteur
$F'=F\times_{\cal{E}}\cal{E}'$ de
$\cal{F}'=\cal{F}\times_{\cal{E}}\cal{E}'$ dans
$\cal{G}'=\cal{G}\times_{\cal{E}}\cal{E}'$ est cart\'esien.
\end{corollaire}

Par suite, le foncteur $\SheafHom_{\cal{E}}(\cal{F},\cal{G}) \to
\SheafHom_{\cal{E}'}(\cal{F}',\cal{G}')$ consid\'er\'e dans
\ifthenelse{\boolean{orig}}{\ignorespaces}{le}
\No \Ref{VI.3} induit un foncteur
$$
\SheafHom_\cart(\cal{F},\cal{G}) \to
\SheafHom_\cart(\cal{F}',\cal{G}') \rm{;}
$$
en d'autres termes, pour $\cal{F}$, $\cal{G}$ fix\'es, \emph{on peut
consid\'erer}
$$
\SheafHom_\cart(\cal{F} \times_{\cal{E}} \cal{E}',
\cal{G}\times_{\cal{E}} \cal{E}')
$$
\emph{comme
\marginpar{168}
un foncteur en~$\cal{E}'$, de la cat\'egorie
${\Cat_{/\cal{E}}}^\circ$ dans $\Cat$}. Si on laisse varier
\'egalement $\cal{F}$, $\cal{G}$, on trouve%
\ifthenelse{\boolean{orig}}{,}{}
un foncteur de la cat\'egorie ${\Cat_{/\cal{E}}}^\circ \times \big(
\Cat^\cart_{/\cal{E}} \big)^\circ \times \Cat^\cart_{/\cal{E}}$
dans~$\Cat$. Lorsqu'on tient compte de l'isomorphisme
$$
\Hom_{\cal{E}'}(\cal{F}',\cal{G}') \isomto
\Hom_\cal{E}(\cal{F}\times_\cal{E}\cal{E}',\cal{G})
$$
envisag\'e au \No \Ref{VI.3}, alors les $\cal{E}'$-foncteurs cart\'esiens
de~$\cal{F}'$ dans~$\cal{G}'$ correspondent aux $\cal{E}$-foncteurs
$\cal{F} \times_\cal{E} \cal{E}' \to \cal{G}$ qui transforment tout
morphisme dont la premi\`ere projection est un morphisme
cart\'esien de~$\cal{F}$, en un morphisme cart\'esien de
$\cal{G}$. Faisant $\cal{F}=\cal{E}$, on trouve (apr\`es changement
de notation):

\begin{corollaire}
\label{VI.6.8}
$\varprojLim(\cal{F}'/\cal{E}')$ est isomorphe \`a
la sous-cat\'egorie pleine de
$\SheafHom_{\cal{E}/-}(\cal{E}',\cal{F})$ form\'ee des
$\cal{E}$-foncteurs $\cal{E}'\to\cal{F}$ qui transforment morphismes
quelconques en morphismes cart\'esiens. En particulier, si $\cal{F}$
est une cat\'egorie fibr\'ee et si $\tilde{\cal{F}}$ est la
sous-cat\'egorie de~$\cal{F}$ dont les morphismes sont les
morphismes cart\'esiens de~$\cal{F}$, alors on a une bijection
$$
\Ob \varprojLim(\cal{F}'/\cal{E}') \isomto
\Hom_{\cal{E}/-}(\cal{E}',\tilde{\cal{F}}) \rm{.}
$$
\end{corollaire}

Cela pr\'ecise la fa\c con dont l'expression
$\varprojLim(\cal{F} \times_{\cal{E}} \cal{E}' /
\cal{E}')$ doit \^etre consid\'er\'ee comme un foncteur en
$\cal{E}'$ et en $\cal{F}$, de la cat\'egorie
${\Cat_{/\cal{E}}}^\circ \times \Cat^\cart_{/\cal{E}}$ dans la
cat\'egorie $\Cat$. On verra ult\'erieurement une d\'ependance
fonctorielle plus compl\`ete par rapport \`a $\cal{E}'$, lorsque
$\cal{F}$ est astreint \`a \^etre une cat\'egorie fibr\'ee.

\begin{corollaire}
\label{VI.6.9}
Soit $\cal{F}$ une cat\'egorie fibr\'ee (\resp pr\'efibr\'ee)
sur~$\cal{E}$, alors $\cal{F}'=\cal{F}\times_{\cal{E}}\cal{E}'$ est
une cat\'egorie fibr\'ee (\resp pr\'efibr\'ee) sur~$\cal{E}'$.
\end{corollaire}

\begin{proposition}
\label{VI.6.10} Soient $\cal{F}$, $\cal{G}$ des cat\'egories
pr\'efibr\'ees sur~$\cal{E}$, $F$ un $\cal{E}$-foncteur
cart\'esien de~$\cal{F}$ dans $\cal{G}$. Pour que $F$ soit
fid\`ele,
\marginpar{169}
\ifthenelse{\boolean{orig}}
{\resp pleinement fid\`ele, (\resp une $\cal{E}$-\'equivalence)}
{(\resp pleinement fid\`ele, \resp une $\cal{E}$-\'equivalence)}
il faut et il suffit que pour tout
\ifthenelse{\boolean{orig}}
{$S\in \Ob \cal{E}$,}
{$S\in \Ob(\cal{E})$,}
le foncteur induit $F_S\colon \cal{F}_S\to \cal{G}_S$ soit
fid\`ele (\resp pleinement fid\`ele, \resp une \'equivalence).
\end{proposition}

D\'emonstration imm\'ediate \`a partir des d\'efinitions.

Pour finir ce num\'ero, nous donnons quelques propri\'et\'es des
cat\'egories fibr\'ees, utilisant l'axiome $\Fib_{\mathrm{II}}$.
\begin{proposition}
\label{VI.6.11} Soit $\cal{F}$ une cat\'egorie pr\'efibr\'ee
sur~$\cal{E}$. Pour que $\cal{F}$ soit fibr\'ee, il faut et il
suffit qu'elle satisfasse la condition suivante:

$\Fib_{\mathrm{II}}'$: Soit $\alpha\colon\eta\to\xi$ un morphisme
cart\'esien dans $\cal{F}$ au-dessus du morphisme $f\colon T\to S$
de~$\cal{E}$. Pour tout morphisme $g\colon U\to T$ dans $\cal{E}$, et
tout
\ifthenelse{\boolean{orig}}
{$\zeta\in\Ob \cal{F}_U$,}
{$\zeta\in\Ob(\cal{F}_U)$,}
l'application $u\mto\alpha\circ u$:
$$
\Hom_g(\zeta,\eta)\to\Hom_{fg}(\zeta,\xi)
$$
est bijective.
\end{proposition}

En d'autres termes, dans une cat\'egorie \emph{fibr\'ee}
sur~$\cal{E}$, les diagrammes cart\'esiens sont caract\'eris\'es
par une propri\'et\'e, plus forte a priori que celle de la
\ifthenelse{\boolean{orig}}
{d\'efinition, (qu'on obtient en faisant $g=\id_T$ dans}
{d\'efinition (qu'on obtient en faisant $g=\id_T$ dans}
l'\'enonc\'e qui pr\'ec\`ede).
\begin{corollaire}
\label{VI.6.12} Soient $\cal{F}$ une cat\'egorie sur~$\cal{E}$,
$\alpha$ un morphisme dans $\cal{F}$. Pour que $\alpha$ soit un
isomorphisme, il faut que $p(\alpha)=f$ soit un isomorphisme et que
$\alpha$ soit cart\'esien; la r\'eciproque est vraie si $\cal{F}$
est fibr\'ee sur~$\cal{E}$.
\end{corollaire}

En effet, si $\alpha$ est un isomorphisme il en est \'evidemment de
m\^eme de~$f=p(\alpha)$; pour tout
\ifthenelse{\boolean{orig}}
{$\eta'\in\Ob\cal{F}_T$,}
{$\eta'\in\Ob(\cal{F}_T)$,}
l'application $u\mto \alpha\circ u$
$$
\Hom(\eta',\eta)\to \Hom(\eta',\xi)
$$
est bijective; comme $f$ est un isomorphisme, on voit de suite
\ifthenelse{\boolean{orig}}
{que un}
{qu'un}
\'el\'ement du premier membre est un $T$-morphisme si et seulement
si son image dans le second est un $f$-morphisme, donc on obtient
ainsi une bijection
$$
\Hom_T(\eta',\eta)\to\Hom_f(\eta',\xi)
$$
ce
\marginpar{170}
qui prouve la premi\`ere assertion. R\'eciproquement, supposons
que $f$ soit un isomorphisme et que $\alpha$ satisfasse la condition
\'enonc\'ee dans $\Fib_{\mathrm{II}'}$ (ce qui signifie donc,
lorsque $\cal{F}$ est
\ifthenelse{\boolean{orig}}
{fibr\'e}
{fibr\'ee}
sur~$\cal{E}$, que $\alpha$ est cart\'esien), alors on voit tout de
suite que pour tout $\zeta\in\Ob\cal{F}$, l'application $u\mto
\alpha\circ u$ de~$\Hom(\zeta,\eta)$ dans $\Hom(\zeta,\xi)$ est
bijective, donc $\alpha$ est un isomorphisme.
\begin{corollaire}
\label{VI.6.13} Soient $\alpha\colon\eta\to\xi$ et
$\beta\colon\zeta\to \eta$ deux morphismes composables dans la
cat\'egorie $\cal{F}$ fibr\'ee sur~$\cal{E}$. Si $\alpha$ est
cart\'esien alors $\beta$ l'est si et seulement si $\alpha\beta$
l'est.
\end{corollaire}

On utilise la d\'efinition des morphismes cart\'esiens sous la
forme renforc\'ee de~\Ref{VI.6.11}.

\section{Cat\'egories cliv\'ees sur~$\mathcal{E}$}
\label{VI.7}

\begin{definition}
\label{VI.7.1} Soit $\cal{F}$ une cat\'egorie sur~$\cal{E}$. On
appelle \emph{clivage}
\index{clivage|hyperpage}%
de~$\cal{F}$ sur~$\cal{E}$ une fonction qui attache \`a tout
$f\in\Fl(\cal{E})$ un foncteur image inverse pour $f$ dans $\cal{F}$,
soit~$f^*$. Le clivage est dit \emph{normalis\'e}
\index{clivage normalis\'e|hyperpage}%
\index{normalis\'e (clivage)|hyperpage}%
si $f=\id_S$ implique $f^*=\id_{\cal{F}_S}$. On appelle
\emph{cat\'egorie cliv\'ee} (\resp \emph{cat\'egorie cliv\'ee
normalis\'ee})
\index{cat\'egorie cliv\'ee (\resp cliv\'ee normalis\'ee)|hyperpage}%
une cat\'egorie $\cal{F}$ sur~$\cal{E}$ munie d'un clivage
(\resp d'un clivage normalis\'e).
\end{definition}

Il est \'evident que $\cal{F}$ admet un clivage si et seulement si
$\cal{F}$ est pr\'efibr\'ee sur~$\cal{E}$, et alors $\cal{F}$
admet un clivage normalis\'e. L'ensemble des clivages sur~$\cal{F}$
est en correspondance biunivoque avec l'ensemble des parties $K$ de
$\Fl(\cal{F})$ satisfaisant les conditions suivantes:
\begin{enumerate}
\item[a)] Les $\alpha\in K$ sont des morphismes cart\'esiens.
\item[b)] Pour tout morphisme $f\colon T\to S$ dans $\cal{E}$ et tout
$\xi\in\Ob(\cal{F}_S)$, il existe un $f$-morphisme unique dans $K$, de
but $\xi$.
\end{enumerate}

Pour que le clivage d\'efini par $K$ soit normalis\'e, il faut et
il suffit que $K$ satisfasse de plus la condition
\begin{enumerate}
\item[c)] Les morphismes identiques dans $\cal{F}$ appartiennent
\`a~$K$.
\end{enumerate}

Les
\marginpar{171}
morphismes \'el\'ements de~$K$ pourront \^etre appel\'es les
\emph{\og morphismes de transport\fg}
\index{transport (morphisme de)|hyperpage}%
pour le clivage envisag\'e.

La notion d'isomorphisme de cat\'egories cliv\'ees sur~$\cal{E}$
est claire. Plus g\'en\'eralement, on peut d\'efinir les
morphismes de~$\cal{E}$-cat\'egories cliv\'ees comme les foncteurs
de~$\cal{E}$-cat\'egories $\cal{F}\to \cal{G}$ qui appliquent
morphismes de transport en morphismes de transport. (Ce sont en
particulier des foncteurs cart\'esiens). De cette fa\c con les
cat\'egories cliv\'ees sur~$\cal{E}$ sont les objets d'une
cat\'egorie, la \emph{cat\'egorie des cat\'egories cliv\'ees
sur} $\cal{E}$. Le lecteur explicitera l'existence de produits,
li\'ee au fait que si une cat\'egorie sur~$\cal{E}$ est produit de
cat\'egories $\cal{F}_i$ sur~$\cal{E}$ munies chacune d'un clivage,
alors $\cal{F}$ est muni d'un clivage naturel correspondant. On laisse
\'egalement au lecteur d'expliciter la notion de changement de base
dans les cat\'egories cliv\'ees.

Nous d\'esignerons par $\alpha_f(\xi)$ le morphisme canonique
$$
\alpha_f(\xi)\colon f^*(\xi)\to\xi.
$$
Il est, on l'a dit, fonctoriel en $\xi$, \ie on a un homomorphisme
fonctoriel
$$
\alpha_f\colon i_Tf^*\to i_S,
$$
o\`u pour tout $S\in\Ob(\cal{E})$, $i_S$ d\'esigne le foncteur
d'inclusion
$$
i_S\colon\cal{F}_S\to \cal{F}
$$

Consid\'erons maintenant des morphismes
$$
f\colon T\to S\quad\text{ et }\quad g\colon U\to T
$$
dans $\cal{E}$, et soit $\xi\in\Ob(\cal{F}_S)$, il existe alors un
unique $U$-morphisme
$$
c_{f,g}(\xi)\colon g^*f^*(\xi)\to (fg)^*(\xi)
$$
rendant
\marginpar{172}
commutatif le diagramme
$$
\xymatrix{
f^*(\xi) \ar[d]_{\alpha_f(\xi)} & & \ar[ll]_{\alpha_g(f^*(\xi))}
g^*(f^*(\xi)) \ar[d]^{c_{f,g}(\xi)} \\
\xi & & \ar[ll]^{\alpha_{fg}(\xi)} (fg)^*(\xi) }
$$
(en vertu de la d\'efinition de~$(fg)^*(\xi)$). Pour $\xi$ variable,
cet homomorphisme est fonctoriel, \ie : on a un homomorphisme
$$
c_{f,g}\colon g^*f^*\to (fg)^*
$$
de foncteurs $\cal{F}_S\to\cal{F}_U$. Notons tout de suite:
\begin{proposition}
\label{VI.7.2} Pour que la cat\'egorie cliv\'ee $\cal{F}$
sur~$\cal{E}$ soit fibr\'ee, il faut et il suffit que les $c_{f,g}$
soient des isomorphismes.
\end{proposition}

On en conclut, prenant pour $f$ un isomorphisme, pour $g$ son inverse,
et en consid\'erant les isomorphismes $c_{f,g}$ et $c_{g,f}$:
\begin{corollaire}
\label{VI.7.3} Si $\cal{F}$ est une cat\'egorie \emph{fibr\'ee}
cliv\'ee sur~$\cal{E}$, alors pour tout isomorphisme $f\colon T\to
S$ dans $\cal{E}$, $f^*$ est une \'equivalence de cat\'egories
$\cal{F}_S\to \cal{F}_T$.
\end{corollaire}

\begin{proposition}
\label{VI.7.4} Soit $\cal{F}$ une cat\'egorie cliv\'ee sur~$\cal{E}$.
On a

$\qquad\qquad A)\qquad
\left\{\begin{array}{ccc}c_{f,\id_T}(\xi)&=&\alpha_{\id_T}(f^*(\xi))\\
c_{\id_S,f}(\xi)&=&f^*(\alpha_{\id_S}(\xi))
\end{array}\right.$

\smallskip
$\qquad\qquad B)\qquad c_{f,gh}(\xi)\cdot
c_{g,h}(f^*(\xi))=c_{fg,h}(\xi)\cdot h^*(c_{f,g}(\xi))$

\noindent(Dans
\marginpar{173}
ces formules, $f,g,h$ d\'esignent des morphismes
$$
V\to U\to T\to S
$$
et $\xi$ un objet de~$\cal{F}_S$).
\end{proposition}

La premi\`ere et seconde relation, dans le cas d'un clivage
normalis\'e, prennent la forme plus simple
\begin{equation*}
\label{eq:VI.7.A'}\quad A')\qquad {c_{f,\id_T}=\id_{f^*},\quad
c_{\id_S,f}=\id_{f^*}.\qquad\qquad\qquad\qquad\qquad\qquad}
\end{equation*}

Quant \`a la troisi\`eme, elle se visualise par la
commutativit\'e du diagramme
\begin{equation*}
\label{eq:VI.7.D} \tag{$D$}
\begin{split}
\xymatrix@C=2.2cm{
h^*g^*f^*(\xi) \ar[r]^-{c_{g,h}(f^*(\xi))} \ar[d]_{h^*(c_{f,g}(\xi))}
& (gh)^*(f^*(\xi)) \ar[d]^{c_{f,gh}(\xi)} \\
h^*(fg)^*(\xi) \ar[r]_-{c_{fg,h}(\xi)} & (fgh)^*(\xi) }
\end{split}
\end{equation*}
Dans le cas des cat\'egories
\ifthenelse{\boolean{orig}}
{fibr\'ees,}
{fibr\'ees}
(o\`u les $c_{f,g}$ sont des isomorphismes), cette commutativit\'e
peut s'exprimer intuitivement par le fait que \emph{l'utilisation
successive des isomorphismes de la forme $c_{f,g}$ ne conduit pas
\`a des \og identifications contradictoires\fg.} On peut \'ecrire
\'egalement cette formule sans argument $\xi$, par l'utilisation du
produit de convolution d'homomorphismes de foncteurs:
$$
c_{fg,h}\circ (h^**c_{f,g})=c_{f,gh}\circ(c_{g,h}*f^*).
$$

\ifthenelse{\boolean{orig}}
{\ignorespaces}
{\enlargethispage{.5cm}}
La d\'emonstration des deux premi\`eres formules~\Ref{VI.7.4} est
triviale, esquissons celle de la troisi\`eme. Pour ceci,
consid\'erons, en plus du carr\'e $(D)$, le carr\'e
d'homomorphismes:
\begin{equation*}
\label{eq:VI.7.D'} \tag{$D'$}
\begin{split}
\xymatrix@C=2cm{
g^*f^*(\xi) \ar[r]^-{\alpha_g(f^*(\xi))}
\ar[d]_{c_{f,g}(\xi)} & f^*(\xi) \ar[d]^{\alpha_f(\xi)} \\
(fg)^*(\xi) \ar[r]_-{\alpha_{fg}(\xi)} & \xi }
\end{split}
\end{equation*}
qui
\marginpar{174}
est commutatif par d\'efinition de
$c_{f,g}(\xi)$. Consid\'erons le diagramme obtenu en joignant les
sommets de~$(D)$ aux sommets correspondants de~$(D')$ par les
homomorphismes de la forme $\alpha$:
$$
\begin{array}{ll}
\alpha_h(g^*f^*(\xi)),&\alpha_{gh}(f^*(\xi))\\
\alpha_h((fg)^*(\xi)),&\alpha_{fgh}(\xi).
\end{array}
$$
Les quatre faces lat\'erales du cube ainsi obtenu sont \'egalement
commutatives: pour celle de gauche, cela provient du fait que la
colonne gauche de~$(D)$ se d\'eduit de la colonne gauche de~$(D')$
par application de~$h$, et que $\alpha_h$ est un homomorphisme
fonctoriel; pour les trois autres, ce n'est autre que la
d\'efinition des op\'erations $c$ des trois c\^ot\'es restants
de~$(D)$. Ainsi les cinq faces du cube autres que la face
sup\'erieure sont commutatives. Il en r\'esulte que les deux
$(fgh)$-morphismes $h^*g^*f^*(\xi)\to(fgh)^*(\xi)$ d\'efinis par
$(D)$ ont un m\^eme compos\'e avec $\alpha_{fgh}(\xi)\colon
(fgh)^*(\xi)\to\xi$, donc ils sont \'egaux par d\'efinition de
$(fgh)^*$.

Bornons-nous pour la suite aux cat\'egories cliv\'ees
\emph{normalis\'ees}. Une telle cat\'egorie donne naissance aux
objets suivants:
\begin{enumerate}
\item[a)] Une application $S\mto \cal{F}_S$ de~$\Ob(\cal{E})$ dans
$\Cat$.
\item[b)] Une application $f\mto f^*$, associant \`a toute
$f\in\Fl(\cal{E})$, de source $T$ et de but $S$,
\ifthenelse{\boolean{orig}}
{une}
{un}
foncteur $f^*\colon \cal{F}_S\to\cal{F}_T$.
\item[c)] Une application $(f,g)\mto c_{f,g}$, associant \`a tout
couple de fl\`eches $(f,g)$ de~$\cal{E}$, un homomorphisme
fonctoriel $c_{f,g}\colon g^*f^*\to (fg)^*$.
\end{enumerate}

D'ailleurs ces donn\'ees satisfont aux conditions exprim\'ees dans
les formules $A')$ et~$B)$ donn\'ees plus haut. (N.B. Si on ne
s'\'etait pas born\'e au cas d'un clivage normalis\'e, il aurait
fallu introduire un objet suppl\'ementaire, savoir une fonction
$S\mto \alpha_S$ qui associe \`a tout objet $S$ de~$\cal{E}$ un
homomorphisme fonctoriel $\alpha_S\colon (\id_S)^*\to
\id_{\cal{F}_S}$; la condition $A')$ se remplacerait alors par la
condition $A)$).

Nous
\marginpar{175}
allons montrer maintenant comment on peut reconstituer (\`a
isomorphisme unique pr\`es) la cat\'egorie cliv\'ee
normalis\'ee $\cal{F}$ sur~$\cal{E}$ \`a l'aide des objets
pr\'ec\'edents.

\section[Cat\'egorie cliv\'ee d\'efinie par un pseudo-foncteur]{Cat\'egorie cliv\'ee d\'efinie par un pseudo-foncteur
$\mathcal{E}^\circ\to\Cat$} \label{VI.8}

Appelons, pour abr\'eger, \emph{pseudo-foncteur} de~$\cal{E}^\circ$
dans $\Cat$ (il faudrait dire, pseudo-foncteur \emph{normalis\'e}),
un ensemble de donn\'ees a), b), c) comme ci-dessus, satisfaisant les
conditions $A')$ et $B)$. Au num\'ero pr\'ec\'edent, nous avons
associ\'e, \`a une cat\'egorie cliv\'ee normalis\'ee
sur~$\cal{E}$, un pseudo-foncteur $\cal{E}^\circ\to\Cat$, ici nous
allons indiquer la construction inverse. Nous laisserons au lecteur la
v\'erification de la plupart des d\'etails, ainsi que du fait que
ces constructions sont bien \og inverses\fg l'une de l'autre. De
fa\c con pr\'ecise, il y aurait lieu de consid\'erer les
pseudo-foncteurs $\cal{E}^\circ\to\Cat$ comme les objets d'une
nouvelle cat\'egorie, et de montrer que nos constructions
fournissent des \'equivalences, quasi-inverses l'une de l'autre,
entre cette derni\`ere et la cat\'egorie des cat\'egories
cliv\'ees au-dessus de~$\cal{E}$, d\'efinie au num\'ero
pr\'ec\'edent.

On pose
$$
\cal{F}_\circ=\underset{S\in\Ob(\cal{E})}\coprod
\ifthenelse{\boolean{orig}}
{\Ob\cal{F}(S),}
{\Ob(\cal{F}(S)),}
$$
ensemble somme des ensembles
\ifthenelse{\boolean{orig}}
{$\Ob\cal{F}(S)$}
{$\Ob(\cal{F}(S))$}
(N.B. nous noterons ici $\cal{F}(S)$ et non $\cal{F}_S$ la valeur en
l'objet $S$ de~$\cal{E}$ du pseudo-foncteur donn\'e, pour \'eviter
des confusions de notation par la suite). On a donc une application
\'evidente:
$$
p_\circ\colon\cal{F}_\circ\to\Ob\cal{E}.
$$
Soient
$$
\ifthenelse{\boolean{orig}}
{\xi}
{\overline{\xi}}
=(S,\xi),\quad \overline{\eta}=(T,\eta)\qquad\qquad
\text{(avec $\xi\in\Ob\cal{F}(S)$, $\eta\in\Ob\cal{F}(T)$)}
$$
deux \'el\'ements de~$\cal{F}_\circ$, et soit $f\in\Hom(T,S)$, on
posera
$$
h_f(\overline{\eta},\overline{\xi})=\Hom_{\cal{F}(T)}(\eta,f^*(\xi)).
$$
\marginpar{176}%
Si on a de plus un morphisme $g\colon U\to T$ dans $\cal{E}$, et un
$\zeta\in\Ob\cal{F}(U)$, on d\'efinit une application, not\'ee
$(u,v)\mto u\circ v$:
$$
h_f(\overline{\eta},\overline{\xi})\times
h_g(\overline{\zeta},\overline{\eta})\to
h_{fg}(\overline{\zeta},\overline{\xi}),
$$
\ie une application
$$
\Hom_{\cal{F}(T)}(\eta,f^*(\xi))\times\Hom_{\cal{F}(U)}(\zeta,g^*(\eta))
\to \Hom_{\cal{F}(U)}(\zeta,(fg)^*(\xi)),
$$
par la formule
$$
u\circ v=c_{f,g}(\xi)\cdot g^*(u)\cdot v,
$$
\ie $u\circ v$ est le compos\'e de la s\'equence
$$
\zeta\lto{u}
g^*(\eta)\lto{g^*(u)}
g^*f^*(\xi)\lto{c_{f,g}(\xi)} (fg)^*(\xi).
$$
On posera d'autre part
$$
h(\overline{\eta},\overline{\xi})=\underset{f\in\Hom(T,S)}\coprod
h_f(\overline{\eta},\overline{\xi}),
$$
et les accouplements pr\'ec\'edents d\'efinissent des
accouplements
$$
h(\overline{\eta},\overline{\xi})\times
h(\overline{\zeta},\overline{\eta})\to
h(\overline{\zeta},\overline{\xi}),
$$
tandis que la d\'efinition des $h(\overline{\eta},\overline{\xi})$
implique une application \'evidente:
$$
p_{\overline{\eta},\overline{\xi}}\colon
h(\overline{\eta},\overline{\xi})\to \Hom(T,S).
$$
Ceci dit, on v\'erifie les points suivants:
\begin{enumerate}
\item[1)] La composition entre \'el\'ements des
\ifthenelse{\boolean{orig}}
{$h(\eta,\xi)$}
{$h(\overline{\eta},\overline{\xi})$}
est \emph{associative}.
\item[2)]
\marginpar{177}
Pour tout $\overline{\xi}=(\xi,S)$ dans $\cal{F}_\circ$,
consid\'erons l'\'el\'ement
\ifthenelse{\boolean{orig}}
{}
{identit\'e}
de
$$
\ifthenelse{\boolean{orig}}
{h_{\id_S}(\overline{\xi},\overline{\xi})=
\Hom_{\cal{F}_S}(\id_S^*(\xi),\xi)=\Hom_{\cal{F}_S}(\xi,\xi),}
{h_{\id_S}(\overline{\xi},\overline{\xi})=
\Hom_{\cal{F}(S)}(\id_S^*(\xi),\xi)=\Hom_{\cal{F}(S)}(\xi,\xi),}
$$
et son image dans $h(\overline{\xi},\overline{\xi})$. Cet objet est
une \emph{unit\'e} \`a gauche et \`a droite pour la composition
entre \'el\'ements des $h(\overline{\eta},\overline{\xi})$.

Cela montre d\'ej\`a que \emph{l'on obtient une cat\'egorie}
$\cal{F}$, en posant
$$
\Ob\cal{F}=\cal{F}_\circ,\quad \Fl{\cal{F}}=
\ifthenelse{\boolean{orig}}
{\underset{\xi,\overline{\eta}\in\cal{F}_\circ} \coprod
h(\overline{\eta},\overline{\xi}).}
{\underset{\overline{\xi},\overline{\eta}\in\cal{F}_\circ}
\coprod h(\overline{\eta},\overline{\xi}).}
$$
(N.B. on ne peut prendre simplement pour $\Fl\cal{F}$ la
\emph{r\'eunion} des ensembles $h(\overline{\eta},\overline{\xi})$,
car ces derniers ne sont pas n\'ecessairement disjoints). De plus:
\item[3)] Les applications $p_\circ\colon \Ob\cal{F}\to\Ob\cal{E}$ et
$p_1=(p_{\overline{\eta},\overline{\xi}})\colon\Fl\cal{F}\to\Fl{\cal{E}}$
d\'efinissent un \emph{foncteur} $p\colon\cal{F}\to\cal{E}$. De
cette fa\c con, $\cal{F}$ devient une cat\'egorie sur~$\cal{E}$, de
plus l'application \'evidente
$h_f(\overline{\eta},\overline{\xi})\to\Hom(\overline{\eta},\overline{\xi})$
induit une \emph{bijection}
$$
h_f(\overline{\eta},\overline{\xi})\isomto
\Hom_f(\overline{\eta},\overline{\xi}).
$$
\item[4)] Les applications \'evidentes
$$ \Ob\cal{F}(S)\to\cal{F}_\circ=\Ob\cal{F},\qquad \Fl\cal{F}(S)\to
\Fl\cal{F},
$$
o\`u la deuxi\`eme est d\'efinie par les applications
\'evidentes
$$
\Hom_{\cal{F}(S)}(\xi,\xi')=
h_{\id_S}(\overline{\xi},\overline{\xi}')\to
\Hom(\overline{\xi},\overline{\xi}')
$$
d\'efinissent un \emph{isomorphisme}
$$
i_S \colon \cal{F}(S)\isomto\cal{F}_S.
$$
\item[5)] Pour tout objet $\overline{\xi}=(S,\xi)$ de~$\cal{F}$, et
tout morphisme $f\colon T\to S$ de~$\cal{E}$, consid\'erons
\marginpar{178}
l'\'el\'ement $\overline{\eta}=(T,\eta)$ de~$\cal{F}_T$, avec
$\eta=f^*(\xi)$, et l'\'el\'ement $\alpha_f(\xi)$ de
$\Hom(\overline{\eta},\overline{\xi})$, image de~$\id_{f^*(\xi)}$ par
le morphisme
$\Hom_{\cal{F}(T)}(f^*(\xi),f^*(\xi))=h_f(\overline{\eta},\overline{\xi})
\to \Hom_f(\overline{\eta},\overline{\xi})$. \emph{Cet \'el\'ement
est cart\'esien, et c'est l'identit\'e dans $\overline{\xi}$ si}
$f=\id_S$, en d'autres termes, l'ensemble des $\alpha_f(\xi)$
d\'efinit un \emph{ clivage normalis\'e de~$\cal{F}$ sur}
$\cal{E}$. De plus, par construction, on a commutativit\'e dans le
diagramme de foncteurs
$$
\xymatrix{
\cal{F}(S) \ar[r]^{f^*} \ar[d]_{i_S} & \cal{F}(T) \ar[d]^{i_T} \\
\cal{F}_S \ar[r]^{f^*_{\cal{F}}} & \cal{F}_T }
$$
o\`u $f^*_{\cal{F}}$ est le foncteur image inverse par $f$, relatif
au clivage consid\'er\'e sur~$\cal{F}$. Enfin:
\item[6)] les homomorphismes $c_{f,g}$ donn\'es avec le
pseudo-foncteur sont transform\'es, par les isomorphismes $i_S$, en
les homomorphismes fonctoriels $c_{f,g}$ associ\'es au clivage
de~$\cal{F}$.
\end{enumerate}

Nous nous bornons \`a donner la v\'erification de 1) (qui est, si
possible, moins triviale que les autres). Il suffit de prouver
l'associativit\'e de la composition entre les objets d'ensembles de
la forme $h_f(\overline{\eta},\overline{\xi})$. Consid\'erons donc
dans $\cal{E}$ des morphismes
$$
S\lfrom{f} T\lfrom{g} U\lfrom{h} V
$$
et des objets
$$
\xi,\ \eta,\ \zeta,\ \tau
$$
dans $\cal{F}(S),\cal{F}(T),\cal{F}(U),\cal{F}(V)$, enfin des
\'el\'ements
$$
\begin{array}{ccc}
u\in
h_f(\overline{\eta},\overline{\xi})&=&\Hom_{\cal{F}(T)}(\eta,f^*(\xi))\\
v\in
h_g(\overline{\zeta},\overline{\eta})&=&\Hom_{\cal{F}(U)}(\zeta,g^*(\eta))\\
w\in
h_h(\overline{\tau},\overline{\zeta})&=&\Hom_{\cal{F}(V)}(\tau,h^*(\zeta)).
\end{array}
$$
On
\marginpar{179}
veut prouver la formule
$$
(u\circ v)\circ w = u\circ (v\circ w),
$$
qui est une \'egalit\'e dans
$\Hom_{\cal{F}(V)}(\tau,(fgh)^*(\xi))$. En vertu des d\'efinitions
les deux membres de cette \'egalit\'e s'obtiennent par composition
suivant le contour sup\'erieur et inf\'erieur du diagramme
ci-dessous:
$$
\xymatrix{
\tau \ar[r]^w \ar[drr]_{v\circ w} & h^*(\zeta)
\ar[r]^{h^*(v)} \ar@<-2pt>
`u[rrrrr] `[rrrrr]^{h^*(u\circ v)} [rrrrr] & h^*g^*(\eta)
\ar[rr]^{h^*g^*(u)}
\ar[d]_{c_{g,h}(\eta)} & & h^*g^*f^*(\xi) \ar[rr]^{h^*(c_{f,g}(\xi))}
\ar[d]^{c_{g,h}(f^*(\xi))} & & h^*(fg)^*(\xi) \ar[d]^{c_{fg,h}(\xi)} \\
& & (gh)^*(\eta) \ar[rr]^{(gh)^*(u)} & & (gh)^*f^*(\xi)
\ar[rr]^{c_{f,gh}(\xi)} & & (fgh)^*(\xi) }
$$
Or le carr\'e m\'edian est commutatif parce que $c_{g,h}$ est un
homomorphisme fonctoriel, et le carr\'e de droite est commutatif en
vertu de la condition~$B)$ pour un pseudo-foncteur. D'o\`u le
r\'esultat annonc\'e.

Bien entendu, il reste \`a pr\'eciser, lorsque le pseudo-foncteur
envisag\'e provient d\'ej\`a d'une cat\'egorie cliv\'ee
normalis\'ee $\cal{F}'$ sur~$\cal{E}$, comment on obtient un
isomorphisme naturel entre $\cal{F}'$ et $\cal{F}$. Nous en laissons
le d\'etail au lecteur.

Nous laissons \'egalement au lecteur d'interpr\'eter, en termes de
pseudo-foncteurs, la notion d'image inverse d'une cat\'egorie
cliv\'ee $\cal{F}$ sur~$\cal{E}$ par un foncteur changement de base
$\cal{E}' \to\cal{E}$.

\section[Exemple]{Exemple: cat\'egorie cliv\'ee d\'efinie par un foncteur
$\cal{E}^\circ\to\Cat$; cat\'egories scind\'ees sur~$\cal{E}$}
\label{VI.9}
Supposons qu'on ait un foncteur
$$
\ifthenelse{\boolean{orig}}
{\phi\colon\cal{E}^\circ\to\Cat,}
{\varphi\colon\cal{E}^\circ\to\Cat,}
$$
il d\'efinit alors un pseudo-foncteur en posant
$$
\ifthenelse{\boolean{orig}}
{\cal{F}(S)=\phi(S),}
{\cal{F}(S)=\varphi(S),}
\quad f^*=\varphi(f),\quad c_{f,g}=\id_{(fg)^*}
$$
\marginpar{180}%
Donc la construction du num\'ero pr\'ec\'edent nous donne une
cat\'egorie $\cal{F}$ cliv\'ee sur~$\cal{E}$, dite associ\'ee au
foncteur $\varphi$. Pour qu'une cat\'egorie cliv\'ee sur~$\cal{E}$
soit isomorphe \`a une cat\'egorie cliv\'ee d\'efinie par un
foncteur $\varphi:\cal{E}^\circ\to\Cat$, il faut et il suffit
manifestement qu'elle satisfasse les conditions:
$$
(fg)^*=g^*f^*,\quad c_{f,g}=\id_{(fg)^*}\quoi.
$$
En termes de l'ensemble $K$ des morphismes de transport, cela signifie
aussi simplement que \emph{le compos\'e de deux morphismes de
transport est un morphisme de transport}. Un clivage d'une
cat\'egorie $\cal{F}$ sur~$\cal{E}$ satisfaisant la condition
pr\'ec\'edente est appel\'e un \emph{scindage}
\index{scindage|hyperpage}%
de~$\cal{F}$ sur~$\cal{E}$, et une cat\'egorie $\cal{F}$
sur~$\cal{E}$ munie d'un scindage est appel\'ee une
\emph{cat\'egorie scind\'ee sur~$\cal{E}$}.
\index{scindee (categorie)@scind\'ee (cat\'egorie)|hyperpage}%
C'est donc un cas particulier de la notion de cat\'egorie
cliv\'ee. La cat\'egorie des cat\'egories scind\'ees
sur~$\cal{E}$ est donc \'equivalente \`a
$\SheafHom(\cal{E}^\circ,\Cat)$. Noter qu'une cat\'egorie
scind\'ee sur~$\cal{E}$ est a fortiori une cat\'egorie cliv\'ee
sur~$\cal{E}$.

Si $\cal{F}$ est une cat\'egorie fibr\'ee sur~$\cal{E}$, il
n'existe pas toujours de scindage sur~$\cal{F}$. Supposons par exemple
que $\Ob\cal{E}$ et $\Ob\cal{F}$ soient r\'eduits \`a un
\'el\'ement, et que l'ensemble des endomorphismes dudit est un
groupe $E$ \resp $F$, de sorte que le foncteur projection $p$ est
donn\'e par un homomorphisme de groupes $p\colon F\to E$, surjectif
puisque $p$ est fibrant. On v\'erifie alors aussit\^ot que
l'ensemble des clivages de~$\cal{F}$ sur~$\cal{E}$ est en
correspondance biunivoque avec l'ensemble des applications $s\colon
E\to F$ telles que $ps=\id_E$ (\ie l'ensemble des \og syst\`emes de
repr\'esentants\fg pour les classes mod le sous-groupe $G$ noyau de
l'homomorphisme surjectif $p\colon F\to E$). Un clivage est un
scindage si et seulement si $s$ est un homomorphisme de groupes. Dire
qu'il existe un scindage signifie donc que l'extension de groupes $F$
de~$E$ par $G$ est triviale, ce qui s'exprime, lorsque $G$ est
commutatif, par la nullit\'e d'une certaine classe de cohomologie
dans $\H^2(E,G)$ (o\`u $G$ est consid\'er\'e comme un groupe
o\`u $E$ op\`ere).

Supposons cependant que $\cal{F}$ soit une cat\'egorie fibr\'ee
sur~$\cal{E}$ telle que les
\marginpar{181}
$\cal{F}_S$ soient des cat\'egories \emph{rigides}, \ie le groupe
des automorphismes de tout objet de~$\cal{F}_S$ est r\'eduit \`a
l'identit\'e. Il est facile alors de prouver que $\cal{F}$ admet un
scindage sur~$\cal{E}$. En effet, on constate d'abord que la question
d'existence d'un scindage n'est pas modifi\'ee si on remplace
$\cal{F}$ par une cat\'egorie $\cal{E}$-\'equivalente, ce qui nous
ram\`ene en l'occurrence au cas o\`u les $\cal{F}_S$ sont des
cat\'egories rigides \emph{et r\'eduites} (\ie deux objets
isomorphes dans $\cal{F}_S$ sont identiques). Mais si $G$ est une
cat\'egorie rigide et r\'eduite, tout isomorphisme de deux
foncteurs $H\to G$ (o\`u $H$ est une cat\'egorie quelconque) est
une identit\'e. Il s'ensuit que si $\cal{F}$ est une cat\'egorie
fibr\'ee sur~$\cal{E}$, telle que les cat\'egories-fibres soient
rigides et r\'eduites, alors il existe un clivage \emph{unique} de
$\cal{F}$ sur~$\cal{E}$, qui est n\'ecessairement un scindage. Donc
$\cal{F}$ est isomorphe \`a la cat\'egorie d\'efinie par un
foncteur $\varphi\colon \cal{E}^\circ\to \Cat$, tel que les
$\varphi(S)$ soient des cat\'egories rigides et discr\`etes, et le
foncteur $\varphi$ est d\'efini \`a isomorphisme pr\`es.

\section{Cat\'egories co-fibr\'ees, cat\'egories bi-fibr\'ees}
\label{VI.10}
Consid\'erons une cat\'egorie $\cal{F}$ au-dessus de~$\cal{E}$,
avec le foncteur projection
$$
p\colon \cal{F}\to \cal{E},
$$
elle d\'efinit une cat\'egorie $\cal{F}^\circ$ au-dessus de
$\cal{E}^\circ$, par le foncteur projection
$$
p^\circ\colon \cal{F}^\circ\to \cal{E}^\circ.
$$
Un morphisme $\alpha\colon \eta\to\xi$ dans $\cal{F}$ est dit
\emph{co-cart\'esien} si c'est un morphisme cart\'esien pour
$\cal{F}^\circ$ sur~$\cal{E}^\circ$. Explicitant, on voit que cela
signifie que pour tout objet $\xi'$ de~$\cal{F}_S$, l'application
$u\mto u\circ \alpha$
$$
\Hom_S(\xi,\xi')\to\Hom_f(\eta,\xi')
$$
est bijective. On dit alors aussi que $(\xi,\alpha)$ est une
\emph{image directe}
\index{image directe|hyperpage}%
de~$\eta$ par $f$, dans la cat\'egorie $\cal{F}$ sur~$\cal{E}$. Si
elle existe pour tout $\eta$ dans $\cal{F}_T$, on dit que le foncteur
image directe par $f$ existe, et on note ce foncteur
\marginpar{182}
$f_*^{\cal{F}}$ ou~$f_*$,
\label{indnot:ft}\oldindexnot{$f_*^{\cal{F}}$ ou~$f_*$|hyperpage}%
une fois choisi. Il est donc d\'efini par
un isomorphisme de bifoncteurs sur~$\cal{F}_T^\circ\times\cal{F}_S$:
$$
\Hom_S(f_*(\eta),\xi)\isomto \Hom_f(\eta,\xi).
$$
Si donc $f_*$ existe, pour que $f^*$ existe, il faut et il suffit que
\ifthenelse{\boolean{orig}}
{$f$}
{$f_*$}
admette un foncteur adjoint, \ie qu'il existe un foncteur
$f^*\colon \cal{F}_S\to\cal{F}_T$ et un isomorphisme de bifoncteurs
$$
\ifthenelse{\boolean{orig}}
{\Hom_E(f_*(\eta),\xi)\isomto \Hom_T(\eta,f^*(\xi)).}
{\Hom_S(f_*(\eta),\xi)\isomto \Hom_T(\eta,f^*(\xi)).}
$$
Soit $g\colon U\to T$ un autre morphisme dans $\cal{E}$, et supposons
que les images inverses et directes par $f,g$ et $fg$
existent. Consid\'erons alors les homomorphismes fonctoriels
$$
\begin{array}{cccc}
c^{f,g}\colon&f_*g_*&\dpl\from&(fg)_*\\
c_{f,g}\colon&g^*f^*&\dpl\to&(fg)^*.
\end{array}
$$
On constate que si on consid\`ere $f_*g_*$ et $g^*f^*$ comme un
couple de foncteurs adjoints, ainsi que $(fg)_*$ et $(fg)^*$, les deux
homomorphismes pr\'ec\'edents sont adjoints l'un de l'autre. Donc
l'un est un isomorphisme si et seulement si l'autre l'est. En
particulier:
\begin{proposition}
\label{VI.10.1} Supposons que la cat\'egorie $\cal{F}$ sur~$\cal{E}$
soit pr\'efibr\'ee et co-pr\'efibr\'ee. Pour qu'elle soit
fibr\'ee, il faut et il suffit qu'elle soit
\ifthenelse{\boolean{orig}}
{cofibr\'ee.}
{co-fibr\'ee.}
\end{proposition}

Bien entendu, on dit que $\cal{F}$ est co-pr\'efibr\'ee \resp co-fibr\'ee
\index{co-fibr\'ee (cat\'egorie)|hyperpage}%
sur~$\cal{E}$, si $\cal{F}^\circ$ est pr\'efibr\'ee \resp fibr\'ee sur~$\cal{E}$. Nous dirons que $\cal{F}$ est bi-fibr\'ee
sur~$\cal{E}$,
\index{bi-fibr\'ee (cat\'egorie)|hyperpage}%
si elle est \`a la fois fibr\'ee et co-fibr\'ee sur~$\cal{E}$.

\section{Exemples divers}
\label{VI.11}

\begin{enumerate}

\item[a)] \textbf{Cat\'egories des fl\`eches de}~$\cal{E}$. Soit
$\cal{E}$ une cat\'egorie. D\'esignons par $\mathbf{\Delta^1}$ la
cat\'egorie associ\'ee \`a l'ensemble totalement ordonn\'e
\`a deux \'el\'ements $[0,1]$;
\marginpar{183}
elle a donc deux objets 0 et 1, et en plus des deux morphismes
identiques une fl\`eche $(0,1)$ de source 0 et but 1. Soit
$$
\CatFl(\cal{E})=\SheafHom(\mathbf{\Delta^1},\cal{E})
$$
on l'appelle la \emph{cat\'egorie des fl\`eches de} $\cal{E}$.
L'objet 1 de~$\mathbf{\Delta^1}$ d\'efinit un foncteur canonique,
appel\'e \emph{foncteur-but}
$$
\CatFl(\cal{E})\to \cal{E}
$$
(le foncteur d\'efini par l'objet 0 de~$\mathbf{\Delta^1}$ est
appel\'e \emph{foncteur-source}). Pour tout objet $S$ de~$\cal{E}$,
la cat\'egorie-fibre $\CatFl(\cal{E})_S$ est canoniquement isomorphe
\`a la cat\'egorie $\cal{E}_{/S}$ des objets de~$\cal{E}$
au-dessus de~$S$.

Consid\'erons un morphisme $f\colon T\to S$ dans $\cal{E}$, alors il
lui correspond un foncteur canonique
$$
f_*\colon \cal{E}_{/T}=\cal{F}_T\to \cal{E}_{/S}=\cal{F}_S
$$
et un isomorphisme fonctoriel
$$
\Hom_S(f_*(\eta),\xi)\isomto\Hom_f(\eta,\xi)
$$
qui fait donc de~$f_*$ un foncteur image directe pour $f$ dans
$\cal{F}$. On a d'ailleurs ici
$$
(\id_S)_*=\id_{\cal{F}_S},\quad (fg)_*=f_*g_*,\quad
c^{f,g}=\id_{(fg)},
$$
\ie $\cal{F}$ est muni d'un co-scindage sur~$\cal{E}$. A fortiori,
$\cal{F}$ est co-fibr\'ee sur~$\cal{E}$. Notons maintenant que
l'ensemble des morphismes dans $\cal{F}$ est en correspondance
biunivoque avec l'ensemble des diagrammes carr\'es commutatifs dans
$\cal{E}$.
$$
\xymatrix{
X\ar[d]_-u&\ar[l]_-{f'}Y\ar[d]^-v\\
S&\ar[l]_-{f}T
}
$$
Par
\marginpar{184}
d\'efinition, le morphisme en question est cart\'esien si le
carr\'e est cart\'esien dans $\cal{E}$, \ie s'il fait de~$Y$ un
produit fibr\'e de~$X$ et $T$ sur~$S$. Le foncteur image inverse
\ifthenelse{\boolean{orig}}
{$f$}
{$f^*$}
existe donc si et seulement si pour tout objet $X$ sur~$S$, le produit
fibr\'e $X\times_ST$ existe. Il r\'esulte de~\Ref{VI.10.1} que si
le produit de deux objets sur un troisi\`eme existe toujours dans~$\cal{E}$, \ie si $\cal{F}$ est pr\'efibr\'ee sur~$\cal{E}$,
alors $\cal{F}$ est m\^eme fibr\'ee sur~$\cal{E}$.

\medskip
\item[b)] \textbf{Cat\'egorie des pr\'efaisceaux ou faisceaux sur
des espaces variables}

Soit $\cal{E}=\mathbf{Top}$ la cat\'egorie des espaces topologiques.
Si $T$ est un espace topologique, nous noterons $\cal{U}(T)$ la
cat\'egorie des ouverts de~$T$, o\`u les morphismes sont les
applications d'inclusion. Si $\cal{C}$ est une cat\'egorie, un
foncteur $\cal{U}(T)^\circ\to \cal{C}$ s'appelle un
\emph{pr\'efaisceau} sur~$T$ \`a valeurs dans $\cal{C}$, et un
\emph{faisceau} s'il satisfait une condition d'exactitude \`a gauche
que nous ne r\'ep\'etons pas ici. La \emph{cat\'egorie
$\cal{P}(T)$ des pr\'efaisceaux sur~$T$ \`a valeurs dans
$\cal{C}$}, est par d\'efinition la cat\'egorie
$\SheafHom(\cal{U}(T)^\circ,\cal{C})$, et la cat\'egorie
$\cal{F}(T)$ des faisceaux sur~$T$ \`a valeurs dans $\cal{C}$ est la
sous-cat\'egorie pleine dont les objets sont les objets de
$\SheafHom(\cal{U}(T)^\circ,\cal{C})$ qui sont des faisceaux. Si
$f\colon T\to S$ est un morphisme dans $\cal{E}$, \ie une
application continue d'espaces topologiques, il lui correspond par
l'application croissante $U\mto f^{-1}(U)$ un foncteur
$\cal{U}(S)\to \cal{U}(T)$, d'o\`u un foncteur
$$
f_*\colon\SheafHom(\cal{U}(T)^\circ,\cal{C})\to
\SheafHom(\cal{U}(S)^\circ,\cal{C})
$$
appel\'e \emph{foncteur image directe de pr\'efaisceaux par}
$f$. On voit aussit\^ot que l'image directe d'un faisceau est un
faisceau, donc le foncteur
\ifthenelse{\boolean{orig}}
{$f_*\cal{P}(T)\to\cal{P}(S)$}
{$f_*\colon \cal{P}(T)\to\cal{P}(S)$}
induit un foncteur, \'egalement not\'e
\ifthenelse{\boolean{orig}}
{$f_*\colon\cal{P}(T)\to\cal{P}(S)$.}
{$f_*\colon\cal{F}(T)\to\cal{F}(S)$.}
On v\'erifie de plus trivialement (par l'associativit\'e de la
composition des foncteurs) qu'on a, pour une deuxi\`eme application
continue $g\colon U\to T$, l'identit\'e
$$
(gf)_* =g_*f_*,\qquad
\text{de m\^eme }\ (\id_S)_*=\id_{\cal{P}(S)}.
$$
De cette fa\c con, on a obtenu un foncteur
$$
S\mto \cal{P}(S)
$$
\resp $$
S\mto \cal{F}(S)
$$
de
\marginpar{185}
$\cal{E}$ dans~$\Cat$. En fait, nous nous int\'eressons au
foncteur correspondant
$$
S\mto \cal{P}(S)^\circ,\qquad\text{\resp~}S\mto
\cal{F}(S)^\circ.
$$
Il d\'efinit une cat\'egorie co-fibr\'ee, et m\^eme
co-scind\'ee, sur la cat\'egorie des espaces topologiques qu'on
appelle la \emph{cat\'egorie co-fibr\'ee des pr\'efaisceaux}
(\resp \emph{faisceaux}) \emph{\`a valeurs dans} $\cal{C}$
(sous-entendu: sur des espaces variables). Explicitant la construction
du \No \Ref{VI.8}, on voit qu'un morphisme d'un pr\'efaisceau $B$ sur~$T$
dans un pr\'efaisceau $A$ sur~$S$ est un couple $(f,u)$ form\'e
d'une application continue de~$T$ dans $S$, et d'un morphisme $u\colon
A\to f_*(B)$ dans la cat\'egorie $\cal{P}(S)$. Cette description
vaut \'egalement pour les morphismes de faisceaux, $\cal{F}$
\'etant une sous-cat\'egorie pleine de~$\cal{P}$.

Dans les cas les plus importants, la cat\'egorie $\cal{P}$ et la
cat\'egorie $\cal{F}$ au-dessus de~$\cal{E}$ sont aussi des
cat\'egories fibr\'ees, \ie pour toute application continue, les
foncteurs image directe $\cal{P}(T)\to \cal{P}(S)$ et
$\cal{F}(T)\to\cal{F}(S)$ ont un foncteur adjoint, qui est alors
not\'e $f^*$ et appel\'e foncteur image inverse de
pr\'efaisceaux \resp foncteur image inverse de faisceaux, par
l'application continue $f$. Ce foncteur existe par exemple si
$\cal{C}=\Ens$. On peut montrer que le foncteur $f^*\colon
\cal{P}(S)\to \cal{P}(T)$ existe chaque fois que dans $\cal{C}$ les
limites inductives (relatives \`a des diagrammes dans l'Univers
consid\'er\'e) existent. La question est moins facile pour
$\cal{F}$; on notera en effet que (m\^eme dans le cas
$\cal{C}=\Ens$) l'image inverse d'un pr\'efaisceau qui est un
faisceau n'est en g\'en\'eral pas un faisceau, en d'autres termes le
foncteur image inverse de faisceau n'est pas isomorphe au foncteur
induit par le foncteur image inverse de pr\'efaisceaux (malgr\'e
la notation commune $f^*$). Ainsi, $\cal{F}$ est une
sous-cat\'egorie co-fibr\'ee de~$\cal{P}$, mais pas une
sous-cat\'egorie fibr\'ee, \ie \emph{le foncteur d'inclusion
$\cal{F}\to\cal{P}$ n'est pas fibrant}.

La cat\'egorie co-fibr\'ee $\cal{P}$ peut se d\'eduire d'une
cat\'egorie co-fibr\'ee (ou plut\^ot fibr\'ee) plus
g\'en\'erale, obtenue ainsi. Pour toute cat\'egorie $\cal{U}$ (dans
l'Univers fix\'e), on pose
$$
\cal{P}(\cal{U})=\SheafHom(\cal{U},\cal{C})
$$
et
\marginpar{186}
on note que $\cal{U}\mto\cal{P}(\cal{U})$ est de fa\c con
naturelle un foncteur contravariant en~$\cal{U}$, de la cat\'egorie
$\Cat$ dans $\Cat$. Il d\'efinit donc une cat\'egorie scind\'ee
au-dessus de~$\cal{E}=\Cat$, que nous noterons $\Cat_{/\!/\cal{C}}$. Les
objets de cette cat\'egorie sont les couples $(\cal{U},p)$ d'une
cat\'egorie $\cal{U}$ et d'un foncteur $p\colon\cal{U}\to\cal{C}$,
et un morphisme de~$(\cal{U},p)$ dans $(\cal{V},q)$ est
essentiellement un couple $(f,u)$, o\`u $f$ est un foncteur
$\cal{U}\to\cal{V}$ et $u$ un homomorphisme de foncteurs $u\colon p\to
qf$. Nous laissons au lecteur le soin d'expliciter la composition des
morphismes dans $\Cat_{/\!/\cal{C}}$. Le foncteur-projection
$$
\cal{F}=\Cat_{/\!/\cal{C}}\to\cal{E}=\Cat
$$
associe au couple $(\cal{U},p)$ l'objet $\cal{U}$; la
cat\'egorie-fibre en $\cal{U}$ est la cat\'egorie
$\SheafHom(\cal{U},\cal{C})$ (\`a isomorphisme pr\`es). Lorsque
dans $\cal{C}$ les limites inductives existent, on montre facilement
que la cat\'egorie fibr\'ee $\Cat_{/\!/\cal{C}}$ sur~$\Cat$ est
\'egalement co-fibr\'ee sur~$\Cat$, \ie on peut d\'efinir la
notion d'\emph{image directe d'un foncteur} $p\colon\cal{U}\to\cal{C}$
par un foncteur $f\colon\cal{U}\to\cal{V}$. La cat\'egorie des
pr\'efaisceaux se d\'eduit de la cat\'egorie fibr\'ee
pr\'ec\'edente par le changement de base
$$
\mathbf{Top}^\circ\to\Cat
$$
(foncteur $S\mto\cal{U}(S)$ d\'efini plus haut), ce qui donne une
cat\'egorie fibr\'ee sur~$\mathbf{Top}^\circ$, et en passant \`a
la cat\'egorie oppos\'ee, on obtient la cat\'egorie
co-fibr\'ee $\cal{P}$ des pr\'efaisceaux au-dessus de
$\mathbf{Top}$. La notion d'image inverse d'un foncteur correspond
\`a celle d'image directe de pr\'efaisceau, la notion d'image
directe d'un foncteur \`a celle d'image inverse d'un
pr\'efaisceau.

\medskip
\item[c)] \textbf{Objets \`a op\'erateurs au-dessus d'un objet
\`a op\'erateurs}

Soit $\cal{F}$ une cat\'egorie sur~$\cal{E}$, et soit $S$ un objet de
$\cal{E}$ o\`u un groupe $G$ op\`ere, \`a gauche pour fixer les
id\'ees. Cet objet \`a op\'erateurs peut s'interpr\'eter comme
correspondant \`a un foncteur $\lambda\colon\cal{E}'\to \cal{E}$ de
la cat\'egorie (\`a un seul objet, ayant $G$ comme groupe
d'endomorphismes) $\cal{E}'$ d\'efinie par $G$, dans la
cat\'egorie $\cal{E}$, et d\'efinie donc par changement de base
une cat\'egorie $\cal{F}'$ au-dessus de~$\cal{E}'$, qui est
fibr\'ee \resp co-fibr\'ee lorsque $\cal{F}$ l'est
sur~$\cal{E}$. Une
\marginpar{187}
section de~$\cal{E}'$ sur~$\cal{F}'$ (n\'ecessairement
cart\'esienne, car $\cal{E}'$ est un groupo\"ide, et tout
isomorphisme dans $\cal{F}'$ est cart\'esien en vertu
de~\Ref{VI.6.12}), peut aussi s'interpr\'eter comme un
$\cal{E}$-foncteur $\cal{E}'\to\cal{F}$ au-dessus de~$\lambda$, ou
aussi comme un objet \`a op\'erateurs $\xi$ dans $\cal{F}$
\og au-dessus\fg de l'objet \`a op\'erateurs $S$.

\medskip
\item[d)] \textbf{Couples de foncteurs adjoints quasi-inverses;
autodualit\'es}

Lorsque la cat\'egorie-base $\cal{E}$ est r\'eduite \`a deux
objets $a,b$ et, en plus des fl\`eches identiques, \`a deux
isomorphismes $f\colon a\to b$ et $g\colon b\to a$ inverse l'un de
l'autre (\ie $\cal{E}$ est un groupo\"ide connexe rigide avec deux
objets), une cat\'egorie cliv\'ee normalis\'ee sur~$\cal{E}$ est
essentiellement la m\^eme chose que le syst\`eme form\'e par
deux cat\'egories $\cal{F}_a$ et $\cal{F}_b$ et un \emph{couple de
foncteurs adjoints} $G\colon \cal{F}_a\to\cal{F}_b$ et $F\colon
\cal{F}_b\to\cal{F}_a$, qui soient des \'equivalences de
cat\'egories (donc quasi-inverses l'un de l'autre). On prendra pour
$\cal{F}_a$ et $\cal{F}_b$ les cat\'egories fibres de~$\cal{F}$,
pour $F$ et $G$ les foncteurs $f^*$ et $g^*$, et les deux
isomorphismes
$$
u\colon FG\isomto\id_{\cal{F}_a}\qquad v\colon
GF\isomto\id_{\cal{F}_b}
$$
sont $c_{g,f}$ et $c_{f,g}$. Les deux conditions usuelles de
compatibilit\'e entre $u$ et $v$ ne sont autres que la condition
\Ref{VI.7.4}~$B)$ pour les compos\'es $fgf$ et $gfg$. Il est facile
de montrer que ces conditions suffisent \`a impliquer qu'on a bien
un pseudo-foncteur $\cal{E}^\circ\to\Cat$.

Un cas int\'eressant est celui o\`u l'on a
$$
\cal{F}_b=\cal{F}_a^\circ,\quad G=F^\circ,\quad v=u^\circ.
$$
On appelle \emph{autodualit\'e}
\index{autodualit\'e|hyperpage}%
dans une cat\'egorie $\cal{C}$, la
donn\'ee d'un foncteur $D\colon \cal{C}\to\cal{C}^\circ$, et d'un
isomorphisme $u\colon DD^\circ\isomto\id_{\cal{C}}$, tels que $u$ et
l'isomorphisme $u^\circ\colon D^\circ D\isomto\id_{\cal{C}^\circ}$
fassent de~$(D,D^\circ)$ un couple de foncteurs adjoints,
(n\'ecessairement quasi-inverses l'un de l'autre). Cette condition
s'\'ecrit:
$$
D(u(x))=u(D(x))\qquad \text{pour tout $x\in\Ob(\cal{C})$.}
$$
\item[e)]
\marginpar{188}
\textbf{Cat\'egories au-dessus d'une cat\'egorie
discr\`ete $\cal{E}$.} On dit que $\cal{E}$ est une
\emph{cat\'egorie discr\`ete} si toute fl\`eche y est une
fl\`eche identique, de sorte que $\cal{E}$ est d\'efini \`a
isomorphisme unique pr\`es par la connaissance de l'ensemble
$I=\Ob(\cal{E})$. La donn\'ee d'une cat\'egorie $\cal{F}$
au-dessus de~$\cal{E}$ \'equivaut donc (\`a isomorphisme unique
pr\`es) \`a la donn\'ee d'une famille de cat\'egories
$\cal{F}_i$ ($i\in I$), les cat\'egories fibres. Toute cat\'egorie
$\cal{F}$ sur~$\cal{E}$ est fibr\'ee, tout $\cal{E}$-foncteur
$\cal{F}\to\cal{G}$ est cart\'esien, on a un isomorphisme canonique
$$ \SheafHom_{\cal{E}/-}(\cal{F},\cal{G}) \isomto
\underset{i}\prod\SheafHom(\cal{F}_i,\cal{G}_i).
$$
En particulier, on obtient
$$
\bf{\Gamma}(\cal{F}/\cal{E})=\varprojLim\cal{F}/\cal{E}\isomto\underset{i}\prod \cal{E}_i.
$$
\item[f)] \textbf{Supposons que $\cal{E}$ ait exactement deux objets
$S$ et $T$,} et en plus des \textbf{morphismes identiques, un
morphisme }$f\colon T\to S$. Alors une cat\'egorie $\cal{F}$
au-dessus de~$\cal{E}$ est d\'efinie, \`a $\cal{E}$-isomorphisme
unique pr\`es, par la donn\'ee de deux cat\'egories $\cal{F}_S$
et $\cal{F}_T$ et d'un bifoncteur $H(\eta,\xi)$ sur~$\cal{F}_T^\circ\times\cal{F}_S$, \`a valeurs dans $\Ens$. En effet,
si~$\cal{F}$ est une cat\'egorie au-dessus de~$\cal{E}$, on lui
associe les deux cat\'egories-fibres $\cal{F}_S$ et~$\cal{F}_T$ et
le bifoncteur $H(\eta,\xi)=\Hom_f(\eta,\xi)$. On laisse au lecteur le
soin d'expliciter la construction en sens inverse. Pour que la
cat\'egorie envisag\'ee soit fibr\'ee (ou pr\'efibr\'ee,
cela revient au m\^eme) il faut et il suffit que le foncteur $H$ soit
repr\'esentable par rapport \`a l'argument $\xi$. Pour qu'elle
soit co-fibr\'ee, il faut et il suffit que $H$ soit
repr\'esentable par rapport \`a l'argument~$\eta$.
\item[g)] Soit $\cal{F}=\cal{C}\times\cal{E}$, consid\'er\'ee
comme cat\'egorie au-dessus de~$\cal{E}$ gr\^ace \`a
$\pr_2$. Alors $\cal{F}$ est fibr\'ee et co-fibr\'ee sur~$\cal{E}$,
et est m\^eme munie d'un scindage et d'un co-scindage canonique,
correspondant au foncteur constant sur~$\cal{E}$, \resp sur~$\cal{E}^\circ$, \`a valeurs dans $\Cat$, de valeur $\cal{C}$. On a
$$
\bf{\Gamma}(\cal{F}/\cal{E})\simeq\SheafHom(\cal{E},\cal{C})
$$
et
\marginpar{189}
$\varprojLim\cal{F}/\cal{E}$ correspond \`a la
sous-cat\'egorie pleine form\'ee des foncteurs
$F\colon\cal{E\to\cal{C}}$ transformant morphismes quelconques en
isomorphismes.
\end{enumerate}

\section{Foncteurs sur une cat\'egorie cliv\'ee}
\label{VI.12}

Soit $\cal{F}$ une cat\'egorie cliv\'ee normalis\'ee
sur~$\cal{E}$. Pour tout objet $S$ de~$\cal{E}$ on d\'esigne par
$$
i_S\colon \cal{F}_S\to\cal{F}
$$
le foncteur d'inclusion. On a donc un homomorphisme fonctoriel, pour
tout morphisme $f\colon T\to S$ dans~$\cal{E}$:
$$
\alpha_f\colon i_T f^*\to i_S,
$$
o\`u $f^*$ est le foncteur changement de base
$\cal{F}_S\to\cal{F}_T$ pour $f$ d\'efini par le clivage. Soit
maintenant
$$
F\colon\cal{F}\to\cal{C}
$$
un foncteur de~$\cal{F}$ dans une cat\'egorie $\cal{C}$, posons,
pour tout~$S\in\Ob(\cal{E})$,
$$
F_S=F\circ i_S\colon\cal{F}_S\to\cal{C}
$$
et pour tout~$f\colon T\to S$ dans~$\cal{E}$,
$$
\varphi_f=F*\alpha_f\colon F_T f^*\to F_S
$$
On a ainsi, \`a tout foncteur $F\colon \cal{F}\to\cal{C}$,
associ\'e une famille $(F_S)$ de foncteurs $\cal{F}_S\to\cal{C}$, et
une famille $(\varphi_f)$ d'homomorphismes de foncteurs $F_T f^*\to
F_S$. Ces familles satisfont aux conditions suivantes:

a) $\varphi_{\id_S}=\id_{F_S}$.

b) Pour deux morphismes $f\colon T\to S$ et $g\colon U\to T$
dans~$\cal{E}$, on a commutativit\'e dans le carr\'e
d'homomorphismes fonctoriels:
$$
\xymatrix@C=2cm{ F_U g^* f^* \ar[r]^-{F_U*c_{f,g}} \ar[d]_{\varphi_g * f^*} &
F_U(fg)^* \ar[d]^{\varphi_{fg}} \\ F_T f^*\ar[r]^-{\varphi_f} & F_S.}
$$
\marginpar{190}%
La premi\`ere relation est triviale, et la deuxi\`eme relation
s'obtient en appliquant le foncteur~$F$ au diagramme commutatif
$$
\xymatrix@C=2cm{ g^* f^* (\xi) \ar[r]^-{c_{f,g}(\xi)}
\ar[d]_{\alpha_g(f^*(\xi))} & (fg)^*(\xi) \ar[d]^{\alpha_{fg}(\xi)}\\
f^*(\xi) \ar[r]^-{\alpha_f(\xi)}& \xi }
$$
pour un objet variable $\xi$ dans~$\cal{F}_S$.

Si $G$ est un deuxi\`eme foncteur $\cal{F}\to\cal{C}$, donnant
naissance \`a des foncteurs $G_S\colon\cal{F}_S\to\cal{C}$ et des
homomorphismes fonctoriels $\psi_f\colon G_Tf\to G_S$, et si $u\colon
F\to G$ est un homomorphisme fonctoriel, alors il lui correspond des
homomorphismes fonctoriels $u*i_S$:
$$
u_S\colon F_S\to G_S
$$
et on constate aussit\^ot que pour tout morphisme $f\colon T\to S$
dans $\cal{E}$, on a commutativit\'e dans les carr\'es
$$
\begin{array}{c}
\xymatrix@C=1.5cm{ F_Tf^* \ar[r]^-{\phi_f} \ar[d]_{u_T * f^*} & F_S
\ar[d]^{u_S} \\ G_T f^* \ar[r]^-{\psi_f}& G_S}
\end{array}
\leqno\quad \textup{c)}
$$

\begin{proposition}
\label{VI.12.1}
Soit $\cal{H}(\cal{F},\cal{C})$ la cat\'egorie dont les objets sont
les couples de familles $(F_S)$ ($S\in\Ob(\cal{F})$) de foncteurs
$\cal{F}_S\to\cal{C}$, et de familles $(\varphi_f)$
($f\in\Fl(\cal{F})$) d'homomorphismes fonctoriels $F_Tf^*\to F_S$,
satisfaisant les conditions \emph{a) et~b)}, et o\`u les morphismes
sont les familles $(u_S)$ ($S\in\Ob(\cal{F})$) d'homomorphismes
\marginpar{191}
$F_S\to G_S$, v\'erifiant la condition de commutativit\'e
\emph{c)} \'ecrite plus haut, (la composition des morphismes se
faisant par la composition des homomorphismes de foncteurs
$\cal{F}_S\to\cal{C}$). Alors les deux lois explicit\'ees plus haut
d\'efinissent un \emph{isomorphisme} $K$ de la cat\'egorie
$\SheafHom(\cal{F},\cal{C})$ avec la cat\'egorie
$\cal{H}(\cal{F},\cal{C})$.
\end{proposition}

Il est trivial qu'on a bien l\`a un \emph{foncteur} de la
premi\`ere cat\'egorie dans la seconde. Ce foncteur est pleinement
fid\`ele, car pour $F,G$ donn\'es, $\Hom(F,G)\to\allowbreak\Hom(K(F),K(G))$
est trivialement injectif; pour montrer que c'est surjectif, il suffit
de noter que la condition de commutativit\'e c) exprime la
fonctorialit\'e des applications $u(\xi)=u_S(\xi)\colon
F(S)=F_S(\xi)\to G(\xi)=G_S(\xi)$ pour les homomorphismes de la forme
$\alpha_f(\xi)$ dans $\cal{F}$, d'autre part on a la fonctorialit\'e
sur chaque cat\'egorie fibre \ie pour les morphismes dans $\cal{F}$
qui sont des $T$-morphismes ($T\in\Ob(\cal{E})$), d'o\`u la
fonctorialit\'e pour tout morphisme dans $\cal{F}$,
\ifthenelse{\boolean{orig}}
{puisque un}
{puisqu'un}
$f$-morphisme (o\`u $f\colon T\to S$ est un morphisme dans
\ifthenelse{\boolean{orig}}
{($\cal{E}$)}
{$\cal{E}$}
est de fa\c con unique) un compos\'e d'un morphisme
$\alpha_f(\xi)$ et d'un $T$-morphisme. Il reste donc \`a prouver que
le foncteur~$K$ est bijectif pour les objets. L'argument
pr\'ec\'edent montre d\'ej\`a que $K$ est injectif pour les
objets, reste \`a prouver qu'il est surjectif, \ie que si on part
d'un syst\`eme $(F_S),(\varphi_f)$, satisfaisant a) et b); et si on
d\'efinit une application $\Ob\cal{F}\to\Ob\cal{C}$ par
$$
F(\xi)=F_S(\xi) \quad\text{pour }
\xi\in\Ob\cal{F}_S\subset\Ob\cal{F}
$$
et une application $\Fl(\cal{F})\to\Fl(\cal{C})$ par
$$
F(\alpha_f(\xi)u')=\varphi_f(\xi)\,F_T(u')
$$
pour tout morphisme $f\colon T\to S$ dans $\cal{E}$, tout objet $\xi$
de~$\cal{F}_S$ et tout $T$-morphisme $u'$ de but $f^*(\xi)$, alors on
obtient un \emph{foncteur} $F$ de~$\cal{F}$ dans $\cal{C}$. En effet,
la relation $F(\id_\xi)=\id_{F(\xi)}$ est triviale, il reste \`a
prouver la multiplicativit\'e $F(uv)=F(u)F(v)$ lorsqu'on a un
$f$-morphisme $u\colon\eta\to\xi$ et un $g$-morphisme
$v\colon\zeta\to\nu$, avec $f\colon T\to S$ et $g\colon U\to T$ des
morphismes de~$\cal{E}$. Posant $w=uv$, on aura
$$
u=\alpha_f(\xi)u'\quoi,\quad v=\alpha_g(\eta)v'\quoi,\quad w=\alpha_{fg}(\xi)w'
$$
\marginpar{192}%
avec
$$
w'=c_{f,g}(\xi)g^*(u')v'\quad \hbox{(\cf \No \Ref{VI.8})}.
$$

Avec ces notations, il faut prouver la commutativit\'e du contour
ext\'erieur du diagramme ci-dessous:
$$
\xymatrix{ F_U(\zeta) \ar@{-<} `u[r] `[rrrrrr]^{F(w')} [rrrrrr]
\ar[rr]^{F_u(v')} \ar[drr]_{F(v)} & & F_Ug^*(\eta)
\ar[rr]^{F_ug^*(u')} \ar[d]^{\varphi_g(\eta)} & & F_Ug^*f^*(\xi)
\ar[rr]^{F_U(c_{f,g}(\xi))}\ar[d]^{\varphi_g(f^*(\xi))}& &
F_U(fg)^*(\xi) \ar[d]^{\varphi_{fg}(\xi)} \\ & & F_T(\eta) \ar@{-<}
`d[r] `[rrrr]_{F(u)} [rrrr] \ar[rr]_{F_T(u')} & & F_Tf^*(\xi)
\ar[rr]_{\varphi_f(\xi)} & & F_S(\xi) }
$$
Or le triangle gauche est commutatif par d\'efinition de~$F(v)$, le
carr\'e m\'edian est commutatif car d\'eduit de l'homomorphisme
$u'\colon \xi\to f(\eta)$ par l'homomorphisme fonctoriel $\alpha_g$,
enfin le carr\'e de droite est commutatif en vertu de la condition
b). La conclusion voulue en r\'esulte.

Supposons maintenant que $\cal{C}$ soit \'egalement une
cat\'egorie cliv\'ee normalis\'ee sur~$\cal{E}$, que nous
appellerons dor\'enavant $\cal{G}$, et que nous nous int\'eressons
aux $\cal{E}$-foncteurs de~$\cal{F}$ dans $\cal{G}$. Si $F$ est un tel
foncteur, il induit des foncteurs
$$
F_S\colon \cal{F}_S\to\cal{G}_S
$$
pour les cat\'egories fibres. D'autre part, pour tout morphisme
$f\colon T\to S$ dans~$\cal{E}$, et tout objet $\xi$ dans~$\cal{F}_S$,
le $f$-morphisme $F(\alpha_f(\xi))$ se factorise de fa\c con unique
par un $T$-morphisme
$$
\varphi_f(\xi)\colon F_T(f^*_{\cal{F}}(\xi))\to
f^*_{\cal{G}}(F_S(\xi))
$$
(o\`u le $\cal{F}$ ou le $\cal{G}$ en indice indique la
cat\'egorie cliv\'ee pour laquelle on prend le foncteur image
inverse), d'o\`u un homomorphisme fonctoriel de foncteurs de
$\cal{F}_S$ dans $\cal{G}_T$:
$$
\varphi_f\colon F_T f^*_{\cal{F}}\to f^*_{\cal{G}}F_S.
$$
\marginpar{193}%
Les deux syst\`emes $(F_S)$ et $(\varphi_f)$ satisfont les
conditions suivantes:

a') $\varphi_{\id_S}=\id_{F_S}$.

b') Pour deux morphismes $f\colon T\to S$ et $g\colon U\to T$ dans
$\cal{E}$, on a commutativit\'e dans le diagramme d'homomorphismes
fonctoriels suivant:
$$
\xymatrix{ F_Ug_F^*f_F^* \ar[rr]^{F_U *c_{f,g}^\cal{F}}
\ar[d]_{\varphi_{g^*}f^*_{\cal{F}}} & & F_U(fg)^*_F
\ar[dd]^{\varphi_{fg}} \\ g^*_{\cal{G}}F_Tf^*_F
\ar[d]_{g^*_{\cal{G}}*\varphi_f} & & \\ g^*_{\cal{G}}f^*_{\cal{G}}F_S
\ar[rr]^{c_{f,g}^\cal{G}*F_S} & & (fg)^*_{\cal{G}}F_S. }
$$
Nous en laissons la v\'erification au lecteur, ainsi que
l'\'enonc\'e et la d\'emonstration de l'analogue de la
proposition~\Ref{VI.12.1}, impliquant que l'on obtient ainsi une
correspondance biunivoque entre l'ensemble des $\cal{E}$-foncteurs de
$\cal{F}$ dans $\cal{G}$, et l'ensemble des syst\`emes
$(F_S)$,$(\varphi_f)$ satisfaisant les conditions a') et b')
ci-dessus. Bien entendu, dans cette correspondance, les foncteurs
cart\'esiens sont caract\'eris\'es par la propri\'et\'e que
les homomorphismes $\varphi_f$ sont des isomorphismes.

\begin{remarquestar}
Bien entendu, il y a int\'er\^et le plus souvent \`a raisonner
directement sur des cat\'egories fibr\'ees sans utiliser des
clivages explicites, ce qui dispense en particulier de faire appel,
pour la notion simple de~$\cal{E}$-foncteur ou de~$\cal{E}$-foncteur
cart\'esien, \`a une interpr\'etation pesante comme
ci-dessus. C'est pour \'eviter des lourdeurs insupportables, et pour
obtenir des \'enonc\'es plus intrins\`eques,
\marginpar{194}
que nous avons d\^u renoncer \`a partir comme dans \cite{VI.2} de la
notion de cat\'egorie cliv\'ee (appel\'ee \og cat\'egorie
fibr\'ee\fg dans \loccit), qui passe au second rang au profit de
celle de cat\'egorie fibr\'ee. Il est d'ailleurs probable que,
contrairement \`a l'usage encore pr\'epond\'erant maintenant,
li\'e \`a d'anciennes habitudes de pens\'ee, il finira par
s'av\'erer plus commode dans les probl\`emes universels, de ne pas
mettre l'accent sur \emph{une} solution suppos\'ee choisie une fois
pour toutes, mais de mettre toutes les solutions sur un pied
d'\'egalit\'e.
\end{remarquestar}



\refstepcounter{chapter}
\addtocontents{toc}{\protect\medskip\protect\par\protect\noindent\protect\textbf{VII: n'existe pas}}

\chapterspace{-2}
\chapter{Descente fid\`element plate}
\label{VIII}
\marginpar{195}

\section{Descente des Modules quasi-coh\'erents}
\label{VIII.1}

Soit $\Sch$
\label{indnot:hb}\oldindexnot{$\Sch$|hyperpage}%
la cat\'egorie des pr\'esch\'emas. Proc\'edant comme
dans~VI~\Ref{VI.11}.b, on trouve que la cat\'egorie des
couples~$(X,F)$ d'un pr\'esch\'ema~$X$ et d'un Module~$F$ sur~$X$,
(o\`u les morphismes sont d\'efinis comme dans \loccit \`a
l'aide de la notion d'image directe de Module par un morphisme
d'espaces annel\'es) peut \^etre consid\'er\'ee comme une
cat\'egorie fibr\'ee au-dessus de~$\Sch$, le foncteur de
changement de base relativement \`a un morphisme~$f\colon X\to Y$
dans $\Sch$ \'etant le foncteur image inverse de Modules
par~$f$. (On notera que la cat\'egorie fibre en~$X\in \Ob(\Sch)$ de
la cat\'egorie fibr\'ee pr\'ec\'edente est la cat\'egorie
\emph{oppos\'ee} \`a la cat\'egorie des Modules sur~$X$). Comme
l'image inverse d'un Module quasi-coh\'erent est quasi-coh\'erent,
on voit que la sous-cat\'egorie pleine de la cat\'egorie des
couples~$(X,F)$, form\'ee des couples pour lesquels $F$ est
quasi-coh\'erent, est une sous-cat\'egorie fibr\'ee de la
cat\'egorie fibr\'ee pr\'ec\'edente. (Par contre, si on ne
fait pas d'hypoth\`eses sur~$f$, l'image directe d'un Module
quasi-coh\'erent n'est pas en g\'en\'eral un Module
quasi-coh\'erent). On appellera simplement cette cat\'egorie
fibr\'ee la \emph{cat\'egorie fibr\'ee des Modules
quasi-coh\'erents sur les pr\'esch\'emas}.

Rappelons d'autre part qu'un morphisme~$f\colon X\to Y$ d'espaces
annel\'es est dit \emph{fid\`element plat}
\index{fidelement plat (morphisme)@fid\`element plat (morphisme)|hyperpage}%
\index{plat (morphisme fid\`element)|hyperpage}%
s'il est \emph{plat} (\ie pour tout~$x\in X$, $\mathcal{O}_{X,x}$ est
un module plat sur~$\mathcal{O}_{Y,f(x)}$, \eqref{IV}), et
\emph{surjectif}. On dit que $f$ est un morphisme \emph{quasi-compact}
si l'image inverse par $f$ de toute partie quasi-compacte est
quasi-compacte; lorsque $f$ est un morphisme de pr\'esch\'emas,
cela signifie aussi que l'image inverse par~$f$ d'un ouvert affine
de~$Y$ est r\'eunion \emph{finie} d'ouverts affines de~$X$.

\begin{theoreme}
\label{VIII.1.1}
Soit
\marginpar{196}
$\cal{F}$ la cat\'egorie fibr\'ee des Modules
quasi-coh\'erents sur les pr\'esch\'emas. Soit
\ifthenelse{\boolean{orig}}
{$g\colon S''\to S$}
{$g\colon S'\to S$}
un morphisme de pr\'esch\'emas, fid\`element plat et
quasi-compact. Alors $g$ est un morphisme de~$\cal{F}$-descente
effective.
\index{descente effective (morphisme de)|hyperpage}%
\end{theoreme}

Rappelons\footnote{Nous admettrons ici la th\'eorie
g\'en\'erale de la descente expos\'ee expos\'ee en d\'etail
dans l'article de J\ptbl \textsc{Giraud} cit\'e dans la
note~\eqref{footnotegiraud} en bas de page de l'Avertissement, travail
que nous citerons \cite{VIII.D} par la suite. \Cf aussi \cite{VIII.2} pour un expos\'e
succinct.} que cela signifie deux choses:

\begin{corollaire}[Descente d'homomorphismes de Modules]
\label{VIII.1.2}
\index{descente d'homomorphismes de modules|hyperpage}%
Soient $g\colon S'\to S$ un morphisme de pr\'esch\'emas,
fid\`element plat et quasi-compact, $F$ et $G$ deux Modules
quasi-coh\'erents sur~$S$, $F'$ et $G'$ leurs images inverses
sur~$S'$, enfin $F''$ et $G''$ leurs images inverses
sur~$S''=S'\times_{S}S'$. Consid\'erons le diagramme d'applications
d'ensembles d\'efini par les foncteurs changement de base par~$g$,
$p_{1}$, $p_{2}$ (o\`u $\xymatrix@C=.5cm{p_{1},p_{2}\colon
S'\times_{S}S'\ar@<2pt>[r]\ar@<-2pt>[r] & S'}$ sont les deux
projections):
$$
\xymatrix@C=.5cm{\Hom_{S}(F,G)\ar[r] &
\Hom_{S'}(F',G')\ar@<2pt>[r]\ar@<-2pt>[r] & \Hom_{S''}(F'',G'').}
$$
Ce diagramme est exact, \ie d\'efinit une bijection du premier
ensemble sur l'ensemble des co\"incidences des deux applications
\'ecrites du deuxi\`eme dans le troisi\`eme.
\end{corollaire}

En d'autres termes, le foncteur changement de base
\ifthenelse{\boolean{orig}}
{par~$f$,}
{par~$g$,}
$F\mto F'$, d\'efinit un foncteur \emph{pleinement
fid\`ele} de la cat\'egorie des Modules quasi-coh\'erents sur~$S$ dans la cat\'egorie des Modules quasi-coh\'erents sur~$S'$
munis d'une donn\'ee de descente relativement
\ifthenelse{\boolean{orig}}
{\`a~$f$.}
{\`a~$g$.}
De plus:

\begin{corollaire}[Descente de Modules]
\label{VIII.1.3}
\index{descente de modules|hyperpage}%
Pour tout Module quasi-coh\'erent~$F'$ sur~$S'$, toute donn\'ee de
descente sur~$F'$ relativement \`a $g$ est \emph{effective},
\index{donn\'ee de descente effective|hyperpage}%
\ie $F'$ est isomorphe avec sa donn\'ee de descente \`a l'image
inverse par~$g$ d'un Module quasi-coh\'erent sur~$S$
(d\'etermin\'e \`a isomorphisme unique pr\`es en vertu
de~\Ref{VIII.1.2}).
\end{corollaire}

En d'autres termes, le foncteur pleinement fid\`ele
pr\'ec\'edent est m\^eme une
\emph{\'equivalence}. Pratiquement, cela signifie qu'il revient au
m\^eme de se donner un Module quasi-coh\'erent sur~$S$, ou un
Module quasi-coh\'erent sur~$S'$ muni d'une donn\'ee de descente
relativement \`a~$g$.

\subsubsection*{D\'emonstration de~\Ref{VIII.1.1}}
Soit
\marginpar{197}
d'abord $T$ un~$S$-pr\'esch\'ema qui est $S$-isomorphe \`a
la somme d'une famille
\ifthenelse{\boolean{orig}}
{de ouverts}
{d'ouverts}
induits~$S_{i}$ de~$S$ qui recouvrent~$S$. Alors il est \'evident
que le morphisme structural~$T\to S$ est un morphisme de
$\cal{F}$-descente effective (cela signifie pr\'ecis\'ement que la
donn\'ee d'un Module quasi-coh\'erent~$F$ sur~$S$ \'equivaut
\`a la donn\'ee de Modules quasi-coh\'erents~$F_{i}$ sur
les~$S_{i}$, et d'isomorphismes de recollement~$\varphi_{ji}\colon
F_{i}|S_{i}\cap S_{j}\to F_{j}|S_{i}\cap S_{j}$ satisfaisant la
condition de cocha\^ines bien connue). En vertu de~VII,~8
il s'ensuit que pour v\'erifier que $g\colon S'\to S$ est un
morphisme de~$\cal{F}$-descente effective, il suffit de le
v\'erifier pour le morphisme~$g_{T}\colon T'=T\times_{S}S'\to T$
d\'eduit de~$g$ par le changement de base~$T\to S$. (Remarquer que
l'hypoth\`ese sur~$T\to S$ reste stable par changement de base
quelconque, donc que $T\to S$ est en fait un morphisme de
$\cal{F}$-descente effective \emph{universel}). Prenant pour~$S_{i}$
des ouverts affines qui recouvrent~$S$, on est donc ramen\'e au cas
o\`u~$S$ est affine.

Alors $S'$ est r\'eunion finie d'ouverts affines, et prenant le
$S$-sch\'ema somme de ces derniers, on trouve un $S$-sch\'ema
affine~$S_{1}$ et un $S$-morphisme~$S_{1}\to S'$ plat et
surjectif. Donc $S_{1}$ est aussi fid\`element plat sur~$S$. Si donc
on prouve qu'un morphisme fid\`element plat et affine est un
morphisme de~$\cal{F}$-descente effective, donc un morphisme de
$\cal{F}$-descente strict universel, (l'hypoth\`ese \'etant en
effet stable par changement de base), on en conclut en particulier que
le morphisme structural~$S_{1}\to S$ est un morphisme de
$\cal{F}$-descente strict universel, et comme il existe un
$S$-morphisme~$S_{1}\to S'$, il en r\'esultera bien, par \cite{VIII.D}, que
$g\colon S'\to S$ est un morphisme de~$\cal{F}$-descente strict.

\ifthenelse{\boolean{orig}}{}
{\enlargethispage{.5cm}}%
Cela nous ram\`ene donc au cas o\`u $g$ est un morphisme affine,
et comme on a vu on peut alors de plus supposer $S$ affine, donc
\emph{on peut supposer $S$ et $S'$ affines}. Dans ce cas,
\Ref{VIII.1.2} \'equivaut au

\begin{lemme}
\label{VIII.1.4}
Soient $A$ un anneau, $A'$ une~$A$-alg\`ebre fid\`element plate,
$M$ et $N$ deux~$A$-modules, $M'$ et $N'$ les~$A'$-modules d\'eduits
par changement d'anneau~$A\to A'$, et $M''$,$N''$
les~$A''=A'\otimes_{A}A'$-modules d\'eduits par changement
d'anneau~$A\to A''$. Alors la suite d'applications ensemblistes
$$
\xymatrix@C=.5cm{\Hom_{A}(M,N)\ar[r] &
\Hom_{A'}(M',N')\ar@<2pt>[r]\ar@<-2pt>[r] & \Hom_{A''}(M'',N'')}
$$
est exacte.
\end{lemme}

Comme
\marginpar{198}
l'homomorphisme~$N\to N'$ est injectif ($A'$ \'etant fid\`element
plat sur~$A$) on voit que la premi\`ere fl\`eche est injective. Il
reste \`a prouver que si un
\ifthenelse{\boolean{orig}}
{$A'$ homomorphisme}
{$A'$-homomorphisme}
$u'\colon M'\to N'$ est compatible avec les donn\'ees de descente,
alors il provient d'un $A$-homomorphisme~$u\colon M\to N$. Or cela
signifie aussi simplement que $u'$ applique le sous-ensemble~$M$
de~$M'$ dans le sous-ensemble~$N$ de~$N'$ (l'application~$u\colon M\to
N$ induite sera alors automatiquement $A$-lin\'eaire puisque $u'$
est $A'$-lin\'eaire, et on voit de m\^eme que $u'$ est
n\'ecessairement \'egal \`a~$u\otimes_{A}A'$). Or si $x\in M$,
alors $u'(x)$ est un \'el\'ement dans le noyau du couple
d'applications~$\xymatrix@C=.5cm{N'\ar@<2pt>[r]\ar@<-2pt>[r] & N''}$. On est
donc ramen\'e pour prouver~\Ref{VIII.1.4}
au cas particulier suivant (correspondant au cas o\`u~$M=A$):

\begin{corollaire}
\label{VIII.1.5}
Soit $N$ un~$A$-module, alors la suite d'applications ensemblistes
$$
\xymatrix@C=.5cm{N\ar[r] & N'\ar@<2pt>[r]\ar@<-2pt>[r] & N''}
$$
est exacte.
\end{corollaire}

Soit en effet $A_{1}$ une~$A$-alg\`ebre fid\`element plate. Pour
montrer que la suite envisag\'ee est exacte, il suffit de prouver
que la suite qui s'en d\'eduit par le changement d'anneau~$A\to
A_{1}$ l'est. Or cette derni\`ere, comme on voit de suite, est celle
relative au
\ifthenelse{\boolean{orig}}
{$A$-module}
{$A_1$-module}
$N_{1}=N\otimes_{A}A_{1}$ et \`a la $A_{1}$-alg\`ebre
$A_{1}'=A_{1}\otimes_{A}A'$. Il suffit donc de trouver un $A_{1}$
fid\`element plat sur~$A$, tel que $\Spec(A_{1}')\to \Spec(A_{1})$
soit un morphisme de~$\cal{F}$-descente strict. Or il suffit en effet
de prendre $A_{1}=A'$, car alors le morphisme pr\'ec\'edent admet
un morphisme inverse \`a droite, donc en vertu de~\cite{VIII.D} c'est un
morphisme de descente effective pour n'importe quelle cat\'egorie
fibr\'ee sur~$\Sch$.

Il reste enfin \`a montrer que si $N'$ est un~$A'$-module muni d'une
donn\'ee de descente pour~$A\to A'$, \ie muni d'un isomorphisme
$$
\varphi\colon N_{1}'\isomto N_{2}'
$$
entre les deux modules d\'eduits de~$N$ par les changements
d'anneaux $\xymatrix@C=.5cm{A'\ar@<2pt>[r]\ar@<-2pt>[r] & A'\otimes_{A}A'}$,
alors
\marginpar{199}
$N'$ est isomorphe avec sa donn\'ee de descente \`a un
module~$N\otimes_{A}A'$. On voit facilement, compte tenu
de~\Ref{VIII.1.5}, que cet \'enonc\'e \'equivaut au suivant:

\begin{lemme}
\label{VIII.1.6}
Soit $N'$ un~$A'$-module muni d'une donn\'ee de descente
relativement \`a~$A\to A'$ (o\`u $A'$ est une
$A$-alg\`ebre). Soit $N$ le sous-$A$-module de~$N'$ form\'e
des~$x$ tels que $\varphi(x\otimes_{A'}1_{A'})=1_{A'}\otimes_{A'}x$,
et consid\'erons l'homomorphisme canonique
$$
N\otimes_{A}A'\to N'\quoi,
$$
(qui est alors compatible avec les donn\'ees de descente). Si $A'$
est fid\`element plat sur~$A$, cet homomorphisme est un
isomorphisme.
\end{lemme}

D\'emontrons ce lemme. Soit encore $A_{1}$ une $A$-alg\`ebre
fid\`element plate, pour montrer que le morphisme envisag\'e est
un isomorphisme, il suffit de prouver qu'il le devient apr\`es le
changement d'anneau
\ifthenelse{\boolean{orig}}
{$A_{1}\to A$.}
{$A\to A_1$.}
Or, utilisant la platitude de~$A_{1}$ sur~$A$, on voit que
l'homomorphisme ainsi obtenu n'est autre que celui qu'on obtiendrait
directement en termes du module~$N'\otimes_{A}A_{1}$ sur~$A_{1}'=A'\otimes_{A}A_{1}$, muni de la donn\'ee de descente
relativement \`a~$A_{1}\to A_{1}'$ qui se d\'eduit canoniquement
par changement d'anneau de celle qui \'etait donn\'ee
sur~$N'$. Ainsi il suffit de trouver un~$A_{1}$ fid\`element plat
sur~$A$ tel que $\Spec(A_{1}')\to \Spec(A_{1})$ soit un morphisme de
$\cal{F}$-descente effective. On prend alors comme ci-dessus
$A_{1}=A'$. Cela ach\`eve la d\'emonstration de~\Ref{VIII.1.6}, et
par l\`a la d\'emonstration de~\Ref{VIII.1.1}.

\begin{corollaire}[Descente de sections de Modules]
\label{VIII.1.7}
\index{descente de sections de modules|hyperpage}%
Soit $g\colon S'\to S$ un morphisme de pr\'esch\'emas,
fid\`element plat et quasi-compact. Pour tout Module
quasi-coh\'erent $G$ sur~$S$, soient $G'$ et $G''$ ses images
inverses sur~$S'$ et $S''=S'\times_{S}S'$, et consid\'erons le
diagramme d'homomorphismes de Modules sur~$S$:
$$
\xymatrix@C=.5cm{G\ar[r] & g_*(G')\ar@<2pt>[r]\ar@<-2pt>[r] & h_*(G'')}
$$
(o\`u $h\colon S''\to S$ est le morphisme structural). Ce diagramme
est \emph{exact}.
\end{corollaire}

En
\marginpar{200}
effet, cela signifie que pour tout ouvert~$U$ dans~$S$, le diagramme
correspondant form\'e par les sections sur~$U$ est exact. On peut
\'evidemment supposer alors $U=S$, et l'exactitude en question est
alors un cas particulier de~\Ref{VIII.1.2}, obtenu en faisant
$F=\mathcal{O}_{S}$.

Comme le foncteur image inverse de Modules est exact \`a droite, on
conclut formellement \`a partir de~\Ref{VIII.1.1}:

\begin{corollaire}[Descente de Modules quotients]
\label{VIII.1.8}
\index{descente de modules quotients|hyperpage}%
Avec les notations de~\Ref{VIII.1.7}, soit de plus, pour tout Module
quasi-coh\'erent~$F$ sur un pr\'esch\'ema, $\Quot(F)$ l'ensemble
des Modules quasi-coh\'erents quotients de~$F$. Avec cette
convention, le diagramme d'applications d'ensembles:
$$
\xymatrix@C=.5cm{\Quot(G)\ar[r] & \Quot(G')\ar@<2pt>[r]\ar@<-2pt>[r] &
\Quot(G'')}
$$
est exact.
\end{corollaire}

(On aurait \'evidemment le m\^eme \'enonc\'e avec les
sous-Modules au lieu de Modules quotients, puisque les deux se
correspondent biunivoquement). Faisant en particulier
$G=\mathcal{O}_{S}$, on trouve:

\begin{corollaire}[Descente des sous-pr\'esch\'emas ferm\'es]
\label{VIII.1.9}
\index{descente des sous-pr\'esch\'emas ferm\'es|hyperpage}%
Pour tout pr\'esch\'ema~$X$, soit $H(X)$ l'ensemble des
sous-pr\'esch\'emas ferm\'es de~$X$. Avec cette notation, et
sous les conditions de~\Ref{VIII.1.7}, le diagramme d'applications
d'ensembles suivant
$$
\xymatrix@C=.5cm{H(S)\ar[r] & H(S')\ar@<2pt>[r]\ar@<-2pt>[r] & H(S'')}
$$
\emph{est exact}.
\end{corollaire}

Il y a lieu de compl\'eter le th\'eor\`eme~\Ref{VIII.1.1} par le
r\'esultat suivant:

\begin{proposition}[Descente de propri\'et\'es de Modules]
\label{VIII.1.10}
\index{descente de propri\'et\'es de modules|hyperpage}%
Soient $g\colon S'\to S$ un morphisme fid\`element plat et
quasi-compact, $F$ un Module quasi-coh\'erent sur~$S$. Pour que~$F$
soit de type fini, \resp de pr\'esentation finie, \resp localement
libre et de type fini, il faut et il suffit que son image inverse $F'$
sur~$S'$ le soit.
\end{proposition}

Il n'y a qu'\`a prouver le \og il suffit\fg. On peut \'evidemment
supposer~$S$ affine,
\marginpar{201}
et rempla\c cant alors~$S'$ par
\ifthenelse{\boolean{orig}}
{la}
{une}
somme d'ouverts affines recouvrant~$S'$, on est ramen\'e au cas
o\`u~$S'$ est \'egalement affine. Alors notre \'enonc\'e
\'equivaut au suivant:

\begin{corollaire}
\label{VIII.1.11}
Soient $A$ un anneau, $A'$ une~$A$-alg\`ebre fid\`element plate,
$M$ un~$A$-module, $M'$ le~$A'$-module $M\otimes_{A}A'$. Pour que~$M$
soit de type fini (\resp de pr\'esentation finie, \resp localement
libre de type fini) il faut et il suffit que~$M'$ le soit.
\end{corollaire}

En effet, on a $M=\varinjlim_i M_{i}$, o\`u les~$M_{i}$
sont les sous-modules de type fini de~$M$. Par suite
$M'=\varinjlim_i M_{i}'$, et si $M'$ est de type fini, $M'$
est \'egal \`a l'un des~$M_{i}'$, donc par fid\`ele platitude
$M$ est \'egal \`a~$M_{i}$, donc $M$ est de type fini. Par suite
il existe une suite exacte
$$
0\to R\to L\to M\to 0\quoi,
$$
avec $L$ libre de type fini, d'o\`u une suite exacte
$$
0\to R'\to L'\to M'\to 0\quoi,
$$
avec $L'$ libre de type fini. Si donc $M'$ est de pr\'esentation
finie, $R'$ est de type fini, donc d'apr\`es ce qui pr\'ec\`ede
$R$ est de type fini, donc $M$ est de pr\'esentation finie. Enfin,
dire que~$M$ est localement libre et de type fini, signifie qu'il est
de pr\'esentation finie et plat (\cf \Ref{IV} dans le cas
noeth\'erien; le cas g\'en\'eral est laiss\'e au
lecteur). Comme chacune de ces propri\'et\'es se descend bien, il
en est de m\^eme de leur conjonction, ce qui ach\`eve la
d\'emonstration.

\begin{remarque}
\label{VIII.1.12}
La conjonction de~\Ref{VIII.1.1} et de~\Ref{VIII.1.10} montre que
l'\'enonc\'e~\Ref{VIII.1.1} reste encore valable, quand on
remplace la cat\'egorie fibr\'ee~$\cal{F}$ par la
sous-cat\'egorie fibr\'ee form\'ee des Modules
quasi-coh\'erents de type fini, \resp de pr\'esentation finie,
\resp localement libres de type fini, \resp localement libres de rang
donn\'e $n$.
\end{remarque}

\section{Descente des pr\'esch\'emas affines sur un autre}
\label{VIII.2}
\marginpar{202}

Comme le foncteur image inverse de Modules est compatible avec le
produit tensoriel et d'autres op\'erations tensorielles, le
th\'eor\`eme~\Ref{VIII.1.1} implique diverses variantes, obtenues
en envisageant, au lieu d'un seul Module quasi-coh\'erent, un Module
quasi-coh\'erent ou un syst\`eme de Modules quasi-coh\'erents
muni de structures suppl\'ementaires diverses s'exprimant \`a
l'aide d'op\'erations tensorielles. Par exemple, la donn\'ee de
trois Modules quasi-coh\'erents~$F$, $G$, $H$ sur~$S$ et d'un
accouplement
$$
F\otimes G\to H\quoi,
$$
\'equivaut \`a la donn\'ee de trois Modules
quasi-coh\'erents~$F'$, $G'$, $H'$ sur~$S'$, munis de donn\'ees de
descente relativement \`a~$g\colon S'\to S$, et munis d'un
accouplement
$$
F'\otimes G'\to H'\quoi,
$$
\og compatible\fg avec ces donn\'ees de descente, au sens \'evident
du terme. Par exemple, si $F=G=H$, on voit que la donn\'ee d'un
Module quasi-coh\'erent~$F$ sur~$S$ muni d'une loi d'alg\`ebre
(dont pour l'instant nous ne supposons pas qu'elle satisfasse \`a
aucun axiome d'associativit\'e, commutativit\'e ou d'existence
d'une section unit\'e), \'equivaut \`a la m\^eme donn\'ee
sur~$S'$, munie en plus d'une donn\'ee de descente. Et utilisant les
r\'esultats du num\'ero pr\'ec\'edent, on constate
aussit\^ot que pour que~$F$ satisfasse \`a l'un des axiomes
habituels auxquels on vient de faire allusion, il faut et il suffit
qu'il en soit ainsi pour~$F'$. Par exemple, la donn\'ee d'une
Alg\`ebre quasi-coh\'erente~$\mathcal{A}$ sur~$S$ (par quoi nous
sous-entendons d\'esormais: associative, commutative, \`a section
unit\'e) \'equivaut \`a la donn\'ee d'une Alg\`ebre
quasi-coh\'erente~$\mathcal{A}'$ sur~$S'$, munie d'une donn\'ee de
descente relativement \`a
\ifthenelse{\boolean{orig}}
{$G\colon S'\to S$.}
{$g\colon S'\to S$.}
Si on se rappelle l'\'equivalence entre la cat\'egorie duale des
Alg\`ebres quasi-coh\'erentes sur~$S$, et de la cat\'egorie
des~$S$-pr\'esch\'emas affines sur~$S$, (EGA~II par.1), on trouve
aussit\^ot:

\begin{theoreme}
\label{VIII.2.1}
Soit $\cal{F}'$ la cat\'egorie fibr\'ee des morphismes affines de
pr\'esch\'emas $f\colon X\to S$, consid\'er\'ee comme
sous-cat\'egorie fibr\'ee de la cat\'egorie fibr\'ee
\marginpar{203}
des fl\`eches dans la cat\'egorie des pr\'esch\'emas $\Sch$
\textup{(VI~\Ref{VI.11}.a)}. Soit $g\colon S'\to S$ un morphisme de
pr\'esch\'emas fid\`element plat et quasi-compact. Alors $g$ est
un morphisme de~$\cal{F}'$-descente effective.
\end{theoreme}

\section[Descente de propri\'et\'es ensemblistes]{Descente de propri\'et\'es ensemblistes et de propri\'et\'es de finitude de morphismes\protect\footnotemark}
\label{VIII.3}


\footnotetext{Pour d'autres r\'esultats comme ceux trait\'es dans les
num\'eros 3 et~4, \Cf EGA~IV~2.3, 2.6,~2.7.}

\begin{proposition}
\label{VIII.3.1}
Soient $f\colon X\to Y$ un~$S$-morphisme, $g\colon S'\to S$ un
morphisme surjectif, $f'\colon X'=X\times_{S}S'\to Y'=Y\times_{S}S'$
le morphisme d\'eduit de~$f$ par le changement de base \`a l'aide
de~$g\colon S'\to S$. Pour que $f$ soit surjectif (\resp radiciel), il
faut et il suffit que $f'$ le soit.
\end{proposition}

On note que $f'$ peut aussi s'obtenir par le changement de base~$Y'\to
Y$, qui est \'egalement surjectif puisque d\'eduit de~$g\colon
S'\to S$ qui l'est. D'autre part, pour tout~$y\in Y$ et tout~$y'\in Y'$
au-dessus de~$y$, on a un isomorphisme
$$
X'_{y'}\isomfrom X_{y}\otimes_{\kres(y)}\kres(y')\quoi,
$$
o\`u $X_{y}$ d\'esigne la fibre de~$X$~en~$y$, et $X'_{y'}$ celle
de~$X'$~en $y'$. Il en r\'esulte que~$X_{y}$ est non vide (\resp a
au plus un point, et ce dernier correspond \`a une extension
r\'esiduelle radicielle) si et seulement si
\ifthenelse{\boolean{orig}}
{$X_{y'}$}
{$X'_{y'}$}
a la m\^eme propri\'et\'e. Cela prouve~\Ref{VIII.3.1}.

\begin{corollaire}
\label{VIII.3.2}
Sous les conditions de~\Ref{VIII.3.1}, si $f'$ est injectif
(\resp bijectif), alors $f$ l'est \'egalement.
\end{corollaire}

Cela provient du fait que si $X'_{y'}$ a au plus un point
(\resp exactement un point) il en est de m\^eme de~$X_{y}$; il en
est bien ainsi, puisque le morphisme $X'_{y'}\to X_{y}$ est surjectif
(car d\'eduit de~$\Spec(\kres(y'))\to \Spec(\kres(y))$ qui l'est).

\begin{proposition}
\label{VIII.3.3}
Avec les notations de~\Ref{VIII.3.1}. Supposons que~$g\colon S'\to S$
soit surjectif et quasi-compact (\resp fid\`element plat et
quasi-compact). Pour que~$f$ soit quasi-compact (\resp de type fini),
il faut et il suffit que~$f'$ le soit.
\end{proposition}

Il
\marginpar{204}
y a \`a prouver seulement le~\og il suffit\fg. On peut \'evidemment
supposer~$S=Y$, puisque l'hypoth\`ese faite sur~$g\colon S'\to S$
est conserv\'ee pour~$Y'\to Y$. De plus on peut supposer~$Y$
affine. Alors $Y'$ est quasi-compact, donc $X'$ est quasi-compact
(puisque~$f'$ l'est par hypoth\`ese). Soit $(X_{i})_{i\in I}$ une
famille d'ouverts affines de~$X$ recouvrant~$X$, alors les~$X_{i}'$
sont des ouverts de~$X'$ recouvrant~$X'$, donc il y a une sous-famille
finie qui recouvre~$X'$. Comme $X'\to X$ est surjectif, il s'ensuit
que les~$X_{i}$ correspondants recouvrent d\'ej\`a~$X$, donc~$X$ est
quasi-compact, \ie $f$ est quasi-compact. Supposons maintenant $f'$
de type fini, prouvons que~$f$ l'est, en supposant~$g$ fid\`element
plat. Rempla\c cant $Y'$ par la somme d'une famille d'ouverts
affines le recouvrant, on peut supposer $Y'$ affine. Enfin, $X$
\'etant recouvert par un nombre fini d'ouverts affines $X_{i}$ par
ce qui pr\'ec\`ede, il faut montrer qu'ils sont de type fini sur~$Y$ sachant que $X_{i}'$ est de type fini sur~$Y'$. Cela nous
ram\`ene alors au

\begin{corollaire}
\label{VIII.3.4}
Soient $B$ une $A$-alg\`ebre, $A'$ une~$A$-alg\`ebre
fid\`element plate, $B'=B\otimes_{A}A'$ la~$A$-alg\`ebre
d\'eduite de~$B$ par changement d'anneau. Pour que~$B$ soit de type
fini, il faut et il suffit que~$B'$ le soit.
\end{corollaire}

Il y a \`a prouver seulement le~\og il suffit\fg. On a
$B=\varinjlim_i B_{i}$, o\`u les $B_{i}$ sont les
sous-alg\`ebres de type fini de~$B$. On a donc
$B'=\varinjlim_i B_{i}'$, et si $B'$ est de type fini
sur~$A'$, alors $B'$ est \'egal \`a un des~$B_{i}'$, donc $B$ est
\'egal \`a~$B_{i}$, donc est de type fini.

\begin{corollaire}
\label{VIII.3.5}
Supposons encore le morphisme de changement de base~$g\colon S'\to S$
fid\`element plat et quasi-compact. Pour que~$f$ soit quasi-fini, il
faut et il suffit que~$f'$ le soit.
\end{corollaire}

En effet, la propri\'et\'e \og quasi-fini\fg est par d\'efinition
la conjonction de \og type fini\fg et \og \`a fibres finies\fg, dont
chacune se descend bien par~$g$, la premi\`ere en vertu
de~\Ref{VIII.3.3}, la seconde par le raisonnement de~\Ref{VIII.3.1}
(n'utilisant que le fait que $g$ soit surjectif).

\begin{remarques}
\label{VIII.3.6}
Soient $A$ un anneau, $X$ un~$A$-pr\'esch\'ema. On voit facilement
que
\marginpar{205}
les conditions suivantes sont \'equivalentes:

\begin{enumerate}
\item[(i)] Il existe un anneau noeth\'erien~$A_{0}$ (qu'on peut si
on veut supposer un sous-anneau de type fini de~$A$), un
$A_{0}$-pr\'esch\'ema de type fini~$X_{0}$, un homomorphisme
\ifthenelse{\boolean{orig}}
{$A\to A_{0}$}
{$A_0\to A$}
et un $A$-isomorphisme~$X\isomto X_{0}\times_{A_{0}}A$.
\item[(ii)] Le morphisme diagonal~$X\to X\times_{\Spec(A)}X$ est
quasi-compact (condition vide si $X$ est s\'epar\'e sur~$A$), $X$
est r\'eunion finie d'ouverts affines~$X_{i}$ dont les
anneaux~$B_{i}$ sont des alg\`ebres de pr\'esentation finie
sur~$A$, \ie quotients d'alg\`ebres de polyn\^omes \`a un
nombre fini d'ind\'etermin\'ees, par des id\'eaux de type fini.
\end{enumerate}

Si $X$ est lui-m\^eme affine, d'anneau~$B$, ces conditions
signifient simplement que~$B$ est une alg\`ebre de pr\'esentation
finie sur~$A$.

Un morphisme~$f\colon X\to Y$ est dit \emph{morphisme de
pr\'esentation finie},
\index{presentation finie@pr\'esentation finie (morphisme de)|hyperpage}%
et on dit encore que~$X$ est de pr\'esentation finie sur~$Y$, si $Y$
est r\'eunion d'ouverts affines~$Y_{i}$, tel que~$X|Y_{i}$ en tant
que $Y_{i}$-pr\'esch\'ema satisfasse aux conditions
\'equivalentes pr\'ec\'edentes. Il en est alors de m\^eme pour
$X|Y'$ pour \emph{tout} ouvert affine $Y'$ dans $Y$. C'est l\`a une
propri\'et\'e stable par changement de base, et d'ailleurs le
compos\'e de deux morphismes de pr\'esentation finie est de
pr\'esentation finie.

Ces notions pos\'ees, on voit sur (ii), proc\'edant comme
dans~\Ref{VIII.1.10}, que cet \'enonc\'e reste valable en y
rempla\c cant les mots \og de type fini\fg par \og de pr\'esentation
finie\fg.
\end{remarques}

\section{Descente de propri\'et\'es topologiques}
\label{VIII.4}

\begin{theoreme}
\label{VIII.4.1}
Soient $g\colon Y'\to Y$ un morphisme, et $Z$ une partie de~$Y$. On
suppose que~$g$ est plat, et qu'il existe un morphisme
quasi-compact~$f\colon X\to Y$ tel que~$Z=f(X)$ (N.B. si $Y$ est
noeth\'erien, cette derni\`ere condition est impliqu\'ee par
\ifthenelse{\boolean{orig}}
{$Z$ est constructible}
{\og $Z$ est constructible\fg}). Alors on a
$$
g^{-1}(\overline{Z})=\overline{g^{-1}(Z)}\quoi.
$$
\end{theoreme}

\ifthenelse{\boolean{orig}}
{On peut supposer~$Y$ affine, puis $Y'$ affine. Comme $Y$ est affine,
$X$ est \emph{r\'eunion
\marginpar{206}
finie d'ouverts affines~$X_{i}$, et rempla\c cant~$X$ par la somme
des~$X_{i}$, on peut supposer \'egalement~$X$
affine. Soient~$A,A',B$ les anneaux de~$Y,Y',X$, $B'=B\otimes_{A}A'$
celui de~$X'=X\times_{Y}Y'$, $I$ le noyau de~$A\to B$, $I'$ le noyau
de~$A'\to B'$, donc les parties ferm\'ees de~$Y$ et~$Y'$
d\'efinies par ces id\'eaux sont respectivement l'adh\'erence
de~$Z=f(X)$ et l'adh\'erence de~$Z'=f'(X')=g^{-1}(Z)$. On veut
\'etablir que cette derni\`ere est \'egale \`a
$g^{-1}(\overline{Z})$, ce qui r\'esultera de~$I'=IA'$, lui-m\^eme
cons\'equence de la platitude de~$A'$ sur~$A$}.}
{On peut supposer~$Y$ affine, puis $Y'$ affine. Comme $Y$ est affine,
$X$ est r\'eunion
\marginpar{206}
finie d'ouverts affines~$X_{i}$, et rempla\c cant~$X$ par la somme
des~$X_{i}$, on peut supposer \'egalement~$X$
affine. Soient~$A,A',B$ les anneaux de~$Y,Y',X$, $B'=B\otimes_{A}A'$
celui de~$X'=X\times_{Y}Y'$, $I$ le noyau de~$A\to B$, $I'$ le noyau
de~$A'\to B'$, donc les parties ferm\'ees de~$Y$ et~$Y'$
d\'efinies par ces id\'eaux sont respectivement l'adh\'erence
de~$Z=f(X)$ et l'adh\'erence de~$Z'=f'(X')=g^{-1}(Z)$. On veut
\'etablir que cette derni\`ere est \'egale \`a
$g^{-1}(\overline{Z})$, ce qui r\'esultera de~$I'=IA'$, lui-m\^eme
cons\'equence de la platitude de~$A'$ sur~$A$.}

\begin{corollaire}
\label{VIII.4.2}
Soient $g\colon Y'\to Y$ un morphisme plat et quasi-compact, et $Z'$
une partie ferm\'ee de~$Y'$ satur\'ee pour la relation
d'\'equivalence ensembliste d\'efinie par $g$. Alors on a
$Z'=g^{-1}(\overline{g(Z')})$.
\end{corollaire}

\ifthenelse{\boolean{orig}}{}
{\enlargethispage{.5cm}}%
On a en effet $Z'=g^{-1}(Z)$, avec~$Z=g(Z')$. On peut alors
appliquer~\Ref{VIII.4.1}, en notant que la condition mise sur~$Z$
dans~\Ref{VIII.4.1} est bien v\'erifi\'ee en prenant pour~$X$ le
pr\'esch\'ema~$Z'$ muni de la structure r\'eduite induite
par~$Y'$. (Le fait que~$g$ soit quasi-compact assure alors que~le
morphisme induit~$f\colon Z'\to Y$ est quasi-compact).

L'\'enonc\'e~\Ref{VIII.4.2} signifie aussi que \emph{la topologie
de~$g(Y')$ induite par~$Y$ est quotient de celle de~$Y'$}. En
particulier:

\begin{corollaire}
\label{VIII.4.3}
Soit $g\colon Y'\to Y$ un morphisme fid\`element plat et
quasi-compact. Alors $g$ fait de~$Y$ un espace topologique quotient
de~$Y'$, \ie pour une partie~$Z$ de~$Y$, $Z$ est ferm\'ee
(\resp ouverte) si et seulement si $Z'=g^{-1}(Z)$ l'est.
\end{corollaire}

Rappelons maintenant que deux \'el\'ements $a,b$, de~$Y'$ ont
m\^eme image dans $Y$ si et seulement si ils sont de la forme
$p_{1}(c),p_{2}(c)$ pour un \'el\'ement convenable $c$
dans~$Y''=Y'\times_{Y}Y'$. Il en r\'esulte que si $g$ est surjectif,
on a un diagramme \emph{exact} d'ensembles
$$
\xymatrix@C=.5cm{\cal{P}(Y)\ar[r] & \cal{P}(Y')\ar@<2pt>[r]\ar@<-2pt>[r] &
\cal{P}(Y'')}\quoi,
$$
o\`u pour tout ensemble~$E$, on d\'esigne par~$\cal{P}(E)$
l'ensemble de ses parties. Ceci pos\'e, \Ref{VIII.4.3} peut aussi
s'interpr\'eter ainsi:

\begin{corollaire}[Descente des parties ouvertes \resp ferm\'ees]
\label{VIII.4.4}
Soit
\marginpar{207}%
\index{descente des parties ouvertes (\resp ferm\'ees)|hyperpage}%
$g\colon Y'{\to} Y$ comme dans~\Ref{VIII.4.3}. Pour tout
pr\'esch\'ema~$X$, soit $\Ouv(X)$ \resp $\Fer(X)$ l'ensemble de
ses parties ouvertes \resp l'ensemble de ses parties
ferm\'ees. Alors on a des diagrammes \emph{exacts} d'applications
ensemblistes (d\'eduits de~$g$ et des deux projections de
$Y''=Y'\times_{Y}Y'$):
\vspace*{-5pt}
\begin{gather*}
\xymatrix@C=.5cm{\Ouv(Y)\ar[r] & \Ouv(Y')\ar@<2pt>[r]\ar@<-2pt>[r] &
\Ouv(Y'')}\\[-3pt]
\xymatrix@C=.5cm{\Fer(Y)\ar[r] & \Fer(Y')\ar@<2pt>[r]\ar@<-2pt>[r] &
\Fer(Y'')}
\end{gather*}
\end{corollaire}

On a le compl\'ement suivant \`a~\Ref{VIII.4.3}:

\begin{corollaire}
\label{VIII.4.5}
Soit $g\colon Y'\to Y$ comme dans~\Ref{VIII.4.3}, et soit~$Z$ une
partie de~$Y$ telle qu'il existe un morphisme quasi-compact~$f\colon
X\to Y$ d'image~$Z$ (par exemple $Z$ constructible, $Y$
noeth\'erien). Pour que $Z$ soit une partie localement ferm\'ee
de~$Y$, il faut et il suffit que $Z'=g^{-1}(Z)$ soit une partie
localement ferm\'ee de~$Y'$.
\end{corollaire}

Il suffit de prouver le \og il suffit\fg. Soit $Y_{1}$ le
sous-pr\'esch\'ema ferm\'e de~$Y$, adh\'erence de~$Z$ muni de
la structure r\'eduite induite, et soit~$Y_{1}'=Y_{1}\times_{Y}Y'$
le sous-pr\'esch\'ema ferm\'e de~$Y'$ image r\'eciproque
de~$Y_{1}$. Son ensemble sous-jacent est l'image
inverse~$g^{-1}(Y_{1})=g^{-1}(\overline{Z})$, donc est \'egal en
vertu de~\Ref{VIII.4.1} \`a~$\overline{Z'}$. Comme $Z'$ est
localement ferm\'e dans~$Y'$, il est ouvert dans~$\overline{Z'}$
donc ouvert dans~$Y_{1}'$. Or ce dernier est fid\`element plat et
quasi-compact sur~$Y_{1}$, donc en vertu de~\Ref{VIII.4.3} on en
conclut que~$Z$ est ouvert dans~$Y_{1}$, \ie dans~$\overline{Z}$, ce
qui signifie que $Z$ est localement ferm\'e.

\begin{corollaire}
\label{VIII.4.6}
Soit $g\colon S'\to S$ un morphisme fid\`element plat et
quasi-compact, $f\colon X\to Y$ un $S$-morphisme, $f'\colon X'\to Y'$
le $S'$-morphisme qui s'en d\'eduit par changement de
base. Supposons que~$f'$ soit une application ouverte (\resp une
application ferm\'ee, \resp \ifthenelse{\boolean{orig}}
{quasi-compact}
{quasi-compacte}
et un hom\'eomorphisme dans, \resp un hom\'eomorphisme sur); alors
$f$ a la m\^eme propri\'et\'e.
\end{corollaire}

Comme $Y'$ est fid\`element plat et quasi-compact sur~$Y$, on peut
supposer~$Y=S$. Soit~$Q$
\marginpar{208}
une partie de~$X$, on a alors (en d\'esignant par $h$ le morphisme
de projection $X'\to X$):
$$
g^{-1}(f(Q))=f'(h^{-1}(Q))\quoi.
$$

Si $Q$ est ouvert (\resp ferm\'e), il en est de m\^eme
de~$h^{-1}(Q)$, donc aussi de~$f'(h^{-1}(Q))$ si on suppose $f'$ une
application ouverte (\resp ferm\'ee), donc il en est de m\^eme
de~$f(Q)$ en vertu de la formule pr\'ec\'edente et
de~\Ref{VIII.4.3}. Cela prouve les deux premi\`eres assertions
dans~\Ref{VIII.4.6}, il reste \`a examiner le cas o\`u~$f'$ est un
hom\'eomorphisme dans, et prouver alors que~$f$ est un
hom\'eomorphisme dans. (Le cas d'un hom\'eomorphisme sur
r\'esultera alors de~\Ref{VIII.3.1}). En vertu
\ifthenelse{\boolean{orig}} {de~\Ref{VIII.3.1}} {de~\Ref{VIII.3.2}}
$f$ est injectif, il reste \`a prouver que l'application $X\to f(X)$
est ouverte. On sait d\'ej\`a que~$f$ est quasi-compact en vertu
de~\Ref{VIII.3.3}. Il suffit de prouver que pour toute partie
ferm\'ee~$Z$ de~$X$, on a $Z=f^{-1}(\overline{f(Z)})$, ce qui
\'equivaut \`a la formule analogue pour les images inverses par
l'application surjective~$h\colon X'\to X$, \ie \`a
$$
Z'=f'^{-1}(g^{-1}(\overline{f(Z)})\quoi,
$$
o\`u on pose $Z'=h^{-1}(Z)$. Or en vertu de~\Ref{VIII.4.1}
appliqu\'e \`a la partie~$f(Z)$ de~$Y$, on a
$g^{-1}(\overline{f(Z)})=\overline{g^{-1}(f(Z))}$, et la formule \`a
prouver \'equivaut \`a
$$
Z'=f'^{-1}(\overline{f'(Z')})\quoi,
$$
qui r\'esulte de l'hypoth\`ese que $f'$ est un hom\'eomorphisme
dans.

N.B. Dans ce dernier raisonnement, supposant d\'ej\`a que~$f$ est
quasi-compact, on n'a pas utilis\'e que~$g$ est quasi-compact, mais
seulement que
\ifthenelse{\boolean{orig}}
{$g'$}
{$g$}
est fid\`element plat. Donc c'est sous cette hypoth\`ese qu'on
peut descendre la propri\'et\'e \og hom\'eomorphisme dans\fg, ou
\og hom\'eomorphisme sur\fg, ou encore gr\^ace au raisonnement
pr\'ec\'edent, la propri\'et\'e \og $f'$ est quasi-compact et
fait de~$f'(X')$ un espace topologique quotient de~$X'$\fg.

Nous dirons qu'un morphisme~$f\colon X\to Y$ de pr\'esch\'emas est
\emph{universellement ouvert} (\resp \emph{universellement ferm\'e},
\resp \emph{universellement bicontinu}, \ifthenelse{\boolean{orig}}{etc...}{etc.})
\index{universellement ouvert (\resp universellement ferm\'e, \resp
universellement bicontinu, etc.)|hyperpage}%
si pour tout changement de base~$Y'\to Y$, $f'\colon X'\to Y'$ est un
morphisme ouvert (\resp ferm\'e, \resp un
\ifthenelse{\boolean{orig}}
{hom\'eom.}
{hom\'eomorphisme}
sur l'espace image). On tire alors de~\Ref{VIII.4.6}:

\begin{corollaire}
\label{VIII.4.7}
Sous
\marginpar{209}
les conditions de~\Ref{VIII.4.6}, pour que~$f$ soit
universellement ouvert, (\resp universellement ferm\'e, \resp un
hom\'eomorphisme dans universel, \resp un hom\'eomorphisme
universel), il faut et il suffit que $f'$ le soit.
\end{corollaire}

\begin{corollaire}
\label{VIII.4.8}
Sous les conditions de~\Ref{VIII.4.6}, pour que~$f$ soit
s\'epar\'e (\resp propre) il faut et il suffit que $f'$ le soit.
\end{corollaire}

Dire que~$f$ est s\'epar\'e signifie que le morphisme
diagonal~$X\to X\times_{Y}X$ est ferm\'e ou aussi universellement
ferm\'e, et la premi\`ere assertion
\ifthenelse{\boolean{orig}}
{}
{de}
\Ref{VIII.4.8} r\'esulte donc de~\Ref{VIII.4.7}. Dire que~$f$ est
propre signifie que $f$ satisfait les conditions a) $f$ est de type
fini b) $f$ est s\'epar\'e c) $f$ est universellement
ferm\'e. La condition a) se descend bien en vertu de~\Ref{VIII.3.3},
b) aussi d'apr\`es ce qu'on vient de voir, enfin c) \'egalement
par~\Ref{VIII.4.7}.

\begin{remarques}
\label{VIII.4.9}
Rappelons que lorsque $g\colon Y'\to Y$ est un morphisme plat de type
fini, avec~$Y$ localement noeth\'erien, alors $g$ est un morphisme
ouvert (VI~\Ref{IV.6.6}) ce qui est un r\'esultat plus
pr\'ecis que~\Ref{VIII.4.3}. On notera cependant que si~$f$ est un
morphisme fid\`element plat et quasi-compact de pr\'esch\'emas
noeth\'eriens, alors $f$ n'est pas en g\'en\'eral un morphisme
ouvert. Soit par exemple $Y$ un sch\'ema irr\'eductible dont le
point g\'en\'erique~$y$ n'est pas ouvert (par exemple une courbe
alg\'ebrique), et prenons pour~$Y'$ le sch\'ema somme~$Y\amalg
\Spec(\kres(y))$, alors l'image par le morphisme structural~$Y'\to Y$ de
la partie ouverte~$\Spec(\kres(y))$ n'est pas une partie ouverte
de~$Y$. Le lecteur remarquera \'egalement que divers \'enonc\'es
du pr\'esent expos\'e deviennent faux si on y abandonne
l'hypoth\`ese que le morphisme fid\`element plat envisag\'e est
aussi quasi-compact, le cas type mettant les \'enonc\'es en
d\'efaut \'etant celui o\`u on prend pour~$Y'$ le sch\'ema
somme des spectres des anneaux locaux des points de~$Y$. Par exemple,
prenant encore pour~$Y$ une courbe alg\'ebrique irr\'eductible, et
pour~$Z$ la partie de~$Y$ r\'eduite au point g\'en\'erique, son
image inverse dans~$Y'$ est ouverte, sans que~$Z$ soit ouverte.
\end{remarques}

\subsection{}
\label{VIII.4.10}
Divers
\marginpar{210}
\'enonc\'es du pr\'esent expos\'e restent valables en y
rempla\c cant l'hypoth\`ese que~$Y'$ soit plat sur~$Y$ par la
suivante: il existe un Module de type fini~$F$ sur~$Y'$, de
support~$Y'$, plat par rapport \`a~$Y$; l'hypoth\`ese de
fid\`ele platitude sera alors remplac\'ee par la
pr\'ec\'edente, plus l'hypoth\`ese que $Y'\to Y$ est
surjectif. Ceci s'applique aux deux premi\`eres assertions
dans~\Ref{VIII.1.10}, \`a~\Ref{VIII.3.3}, \Ref{VIII.3.5},
\Ref{VIII.4.1} et par suite \`a tous les r\'esultats du
pr\'esent num\'ero.

\section{Descente de morphismes de pr\'esch\'emas}
\label{VIII.5}

\begin{proposition}
\label{VIII.5.1}
Soit $g\colon S'\to S$ un morphisme de pr\'esch\'emas.
\begin{enumerate}
\item[a)] Supposons que~$g$ soit surjectif, et que l'homomorphisme
$$
g^*\colon \mathcal{O}_{S}\to g_*(\mathcal{O}_{S'})
$$
soit injectif, alors~$g$ est un \'epimorphisme dans la cat\'egorie
des pr\'esch\'emas, et m\^eme dans la cat\'egorie des espaces
annel\'es.
\item[b)] Supposons que~$g$ soit surjectif et fasse de~$S$ un espace
topologique quotient de~$S'$. Soit~$S''=S'\times_{S}S'$ et
soit~$h\colon S''\to S$ le morphisme structural, consid\'erons le
diagramme d'homomorphismes canonique:
$$
\xymatrix@C=.5cm{\mathcal{O}_{S} \ar[r] &
g_*(\mathcal{O}_{S'})\ar@<2pt>[r]\ar@<-2pt>[r] &
h_*(\mathcal{O}_{S''})}\quoi,
$$
et supposons ce diagramme \emph{exact}. Alors $g$ est un
\'epimorphisme effectif
\index{epimorphisme effectif@\'epimorphisme effectif|hyperpage}%
dans la cat\'egorie des pr\'esch\'emas (et aussi dans la
cat\'egorie des espaces annel\'es), \ie le diagramme
$$
\xymatrix@C=.5cm{S & S'\ar[l] & S''\ar@<2pt>[l]\ar@<-2pt>[l]}
$$
est exact.
\end{enumerate}
\end{proposition}

\subsubsection*{D\'emonstration} a) Il faut montrer qu'un morphisme d'espaces
annel\'es~$f\colon S\to Z$ est connu quand on
conna\^it~$fg$. Or comme~$g$ est surjectif, on conna\^it
l'application ensembliste~$f_{0}$ sous-jacente \`a~$f$, reste \`a
d\'eterminer l'homomorphisme de faisceaux
d'anneaux~$\mathcal{O}_{Z}\to \mathcal{O}_{S}$, ou ce qui revient au
m\^eme l'homomorphisme
\marginpar{211}
de faisceaux d'anneaux
$$
u\colon f_{0}^{-1}(\mathcal{O}_{Z})\to \mathcal{O}_{S}
$$
d\'efini par~$f$. On conna\^it d\'ej\`a l'homomorphisme
$$
(fg)_{0}^{-1}(\mathcal{O}_{Z})=
g_{0}^{-1}(f_{0}^{-1}(\mathcal{O}_{Z}))\to \mathcal{O}_{S'}
$$
d\'efini par~$fg$, ou ce qui revient au m\^eme, on dispose d'un
homomorphisme
$$
f_{0}^{-1}(\mathcal{O}_{Z})\to
g_{0*}(\mathcal{O}_{S'})=g_*(\mathcal{O}_{S'})\quoi.
$$
On constate aussit\^ot que ce dernier n'est autre que le compos\'e
de~$g^*\colon \mathcal{O}_{S}\to g_*(\mathcal{O}_{S'})$ et
de~$u$, et comme $g^*$ est injectif, $u$ est connu quand on
conna\^it~$g^*u$. [N.B. on n'a pas utilis\'e \'evidemment
que~$g\colon S'\to S$ est un morphisme de pr\'esch\'emas,
l'\'enonc\'e vaudrait pour un morphisme quelconque d'espaces
annel\'es; la m\^eme remarque vaut pour b), tant dans la
cat\'egorie des espaces annel\'es, que dans la cat\'egorie des
espaces annel\'es en anneaux locaux. Noter aussi que si~$g$ est un
morphisme de pr\'esch\'emas pas n\'ecessairement surjectif, mais
tel que~$g^*\colon \mathcal{O}_{S}\to g_*(\mathcal{O}_{S'})$
soit injectif, alors pour deux morphismes~$f_{1},f_{2}$ de~$S$ dans un
\emph{sch\'ema}~$Z$ tel que~$f_{1}g=f_{2}g$, on a $f_{1}=f_{2}$; en
effet, si~$I$ est l'Id\'eal sur~$S$ qui d\'efinit le
sous-pr\'esch\'ema de~$S$ des co\"incidences de~$f_{1},f_{2}$
(image inverse du sous-pr\'esch\'ema diagonal de~$Z\times Z$
par~$(f_{1},f_{2})$), on voit que~$I$ est contenu
dans~$\Ker(g^*)$].

b) On doit montrer que pour tout espace annel\'e~$Z$, le diagramme
suivant d'applications
$$
\xymatrix@C=.5cm{\Hom(S,Z)\ar[r] & \Hom(S',Z)\ar@<2pt>[r]\ar@<-2pt>[r] &
\Hom(S'',Z)}
$$
est exact, et qu'il en est de m\^eme lorsque~$Z$ est un espace
annel\'e en anneaux locaux et qu'on se borne aux homomorphismes
d'espaces annel\'es en anneaux locaux. Comme on sait d\'ej\`a par a)
que la premi\`ere application est injective, il reste \`a voir que
si~$f'\colon S'\to Z$ est un homomorphisme d'espaces annel\'es tel
que~$f'p_{1}=f'p_{2}$,
\marginpar{212}
alors~$f'$ est de la forme~$fg$, o\`u~$f\colon S\to Z$ est un
homomorphisme d'espaces annel\'es. Comme~$g$ est surjectif, il est
alors \'evident que si~$f'$ est un morphisme d'espaces annel\'es
en anneaux locaux, il en sera de m\^eme pour~$f$.

De l'hypoth\`ese sur~$f'$ r\'esulte que l'application ensembliste
sous-jacente~$f_{0}'$ est constante sur les fibres de
l'application~$g_{0}$, donc comme cette derni\`ere est surjective,
$f_{0}'$ se factorise de fa\c con unique en~$f_{0}'=f_{0}g_{0}$,
o\`u~$f_{0}\colon S\to Z$ est une application, n\'ecessairement
continue puisque~$g_{0}$ identifie~$S$ \`a un espace topologique
quotient de~$S'$. Consid\'erons maintenant l'homomorphisme
$$
f_{0}^{-1}(\mathcal{O}_{Z})\to g_*(\mathcal{O}_{S'})
$$
d\'eduit de l'homomorphisme~$(f_{0}g_{0})^{-1}(\mathcal{O}_{Z})\to
\mathcal{O}_{S'}$ correspondant
\`a~$f'$. L'hypoth\`ese~$f'p_{1}=f'p_{2}$ s'interpr\`ete alors
en disant que les compos\'es de l'homomorphisme pr\'ec\'edent
avec les deux homomorphismes
$\xymatrix@C=.5cm{g_*(\mathcal{O}_{S'})\ar@<2pt>[r]\ar@<-2pt>[r] &
h_*(\mathcal{O}_{S''})}$ sont les m\^emes donc d'apr\`es
l'hypoth\`ese b) il se factorise par un morphisme
$$
f_{0}^{-1}(\mathcal{O}_{Z})\to \mathcal{O}_{S}\quoi.
$$
Ce dernier d\'efinit un morphisme d'espaces annel\'es~$f\colon
S\to Z$, qui est le morphisme cherch\'e.

\begin{theoreme}
\label{VIII.5.2}
Soit $\cal{F}$ la cat\'egorie fibr\'ee des fl\`eches dans la
cat\'egorie~$\Sch$ des pr\'esch\'emas
\textup{(VI~\Ref{VI.11}.a)}. Alors tout morphisme~$g\colon S'\to S$
fid\`element plat et quasi-compact est un morphisme
de~$\cal{F}$-descente (ou encore, comme on dit, un morphisme de
descente dans~$\Sch$).
\index{descente (morphisme de)|hyperpage}%
\end{theoreme}

Cela signifie donc ceci: soit~$S''=S'\times_{S}S'$, et pour deux
pr\'esch\'emas~$X,Y$ sur~$S$, consid\'erons leurs images
inverses~$X',Y'$ sur~$S'$ et leurs images inverses~$X'',Y''$
sur~$S''$, d'o\`u un diagramme d'applications
$$
\xymatrix@C=.5cm{\Hom_{S}(X,Y)\ar[r] &
\Hom_{S'}(X',Y')\ar@<2pt>[r]\ar@<-2pt>[r] & \Hom_{S''}(X'',Y'')}\quoi;
$$
ces notations pos\'ees, \Ref{VIII.5.2} signifie que ce diagramme est
exact. On notera qu'il
\marginpar{213}
n'est pas vrai en g\'en\'eral que~$g$ soit un morphisme de
descente effective, \ie que pour tout pr\'esch\'ema~$X'$
sur~$S'$, toute donn\'ee de descente sur~$X'$ relativement
\`a~$g\colon S'\to S$ soit effective. La question de
l'effectivit\'e, souvent d\'elicate, sera examin\'ee au
\No \Ref{VIII.7}.

On a vu \cite{VIII.D}, (compte tenu que dans~$\Sch$ les produits fibr\'es
existent) que l'\'enonc\'e~\Ref{VIII.5.2} \'equivaut au suivant:

\begin{corollaire}
\label{VIII.5.3}
Un morphisme fid\`element plat et quasi-compact de
pr\'esch\'emas est un \'epi\-mor\-phisme effectif universel.
\end{corollaire}

Comme un morphisme fid\`element plat et quasi-compact reste tel par
toute extension de la base, on est ramen\'e \`a prouver que c'est
un \'epimorphisme effectif. On applique alors le
crit\`ere~\Ref{VIII.5.1} b), qui donne le r\'esultat voulu, compte
tenu de~\Ref{VIII.4.3} et~\Ref{VIII.1.7}.

\begin{corollaire}
\label{VIII.5.4}
Soit $g\colon S'\to S$ un morphisme fid\`element plat et
quasi-compact, $f\colon X\to Y$ un~$S$-morphisme, $f'\colon X'\to Y'$
le~$S'$-morphisme qui s'en d\'eduit par le changement de base~$S'\to
S$. Pour que~$f$ soit un isomorphisme, il faut et il suffit que~$f'$
le soit.
\end{corollaire}

En effet, si $f'$ est un isomorphisme, c'est aussi un isomorphisme
pour les structures de descente naturelles sur~$X',Y'$, et comme le
foncteur~$X\mto X'$ de~$\Sch_{/S}$ dans la cat\'egorie des
pr\'esch\'emas sur~$S'$ munis d'une donn\'ee de descente
relativement \`a~$g$ est pleinement fid\`ele par~\Ref{VIII.5.2},
la conclusion voulue appara\^it.

\begin{corollaire}
\label{VIII.5.5}
Sous les conditions de~\Ref{VIII.5.4}, pour que~$f$ soit une immersion
ferm\'ee (\resp une immersion ouverte, \resp une immersion
quasi-compacte) il faut et suffit que~$f'$ le soit.
\end{corollaire}

On peut supposer comme d'habitude~$Y=S$, et il y a \`a prouver
seulement le \og il suffit\fg. Notons que le fait que~$X'/Y'$ soit muni
d'une donn\'ee de descente relativement \`a~$g\colon Y'\to Y$, et
que le morphisme structural~$f'\colon X'\to Y'$ soit une immersion
donc un monomorphisme, implique que les deux sous-objets de~$Y''$
images
\marginpar{214}
inverses de~$X'/Y'$ par l'une et l'autre projection de~$S''$ dans~$S'$
sont les m\^emes. Si~$f'$ est une immersion ferm\'ee, il en
r\'esulte en vertu de~\Ref{VIII.1.9} qu'il existe un
sous-pr\'esch\'ema ferm\'e $X_{1}$ de~$Y$ dont l'image inverse
par~$g\colon Y'\to Y$ est~$X'$. Donc par l'unicit\'e de la solution
d'un probl\`eme de descente relativement \`a un morphisme de
$\cal{F}$-descente,
\ifthenelse{\boolean{orig}}
{}
{il}
r\'esulte que~$X_{1}$ est $Y$-isomorphe \`a~$X$, donc $f\colon
X\to Y$ est une immersion ferm\'ee. On proc\`ede de m\^eme pour
une immersion ouverte, en utilisant~\Ref{VIII.4.4}. Si enfin $f'$ est
une immersion quasi-compacte, $f$ est quasi-compact en vertu
de~\Ref{VIII.3.3}, donc on peut appliquer \`a la partie~$f(X)$
de~$Y$ le crit\`ere~\Ref{VIII.4.5}, qui prouve que~$f(X)$ est
localement ferm\'e puisque son image inverse~$f'(X')$ dans~$Y'$
l'est. Rempla\c cant alors~$Y$ par une partie ouverte dans
laquelle~$f(X)$ soit ferm\'ee, on est ramen\'e au cas o\`u~$f'$
est une immersion ferm\'ee, donc $f$ l'est en vertu de ce qui
pr\'ec\`ede.

\begin{corollaire}
\label{VIII.5.6}
Sous les conditions de~\Ref{VIII.5.4}, pour que~$f$ soit affine, il
faut et il suffit que~$f'$ le soit.
\end{corollaire}

On proc\`ede comme dans~\Ref{VIII.5.5}, en utilisant~\Ref{VIII.2.1}
(On peut aussi utiliser le crit\`ere cohomologique de Serre [EGA II
5.2], qui d\'emontre~\Ref{VIII.5.6} sans utiliser de technique de
descente).

\begin{corollaire}
\label{VIII.5.7}
Sous les conditions de~\Ref{VIII.5.4}, pour que~$f$ soit entier
(\resp fini, \resp fini et localement libre) il faut et il suffit
que~$f'$ le soit.
\end{corollaire}

Il y a \`a prouver seulement le \og il suffit\fg, et comme d'habitude
on peut supposer $Y=S$, $Y$ affine, et $Y'$ affine. Comme
l'hypoth\`ese implique que~$f'$ est affine, il en est de m\^eme
de~$f$ d'apr\`es~\Ref{VIII.5.6}, donc~$X$ et par suite~$X'$ sont
affines. Soient~$A,A',B,B'=B\otimes_{A}A'$ les anneaux
de~$Y,Y',X,X'$. On a $B=\varinjlim_i B_{i}$, o\`u~$B_{i}$
parcourt les sous-$A$-alg\`ebres de~$B$ qui sont de type fini
sur~$A$, d'o\`u $B'=\varinjlim_i B_{i}'$, o\`u
les~$B_{i}'$ sont des sous-alg\`ebres de type fini de la
$A'$-alg\`ebre~$B'$. Si $B'$ est entier
\ifthenelse{\boolean{orig}}
{sur~$A$,}
{sur~$A'$,}
les~$B_{i}'$ sont des modules de type fini sur~$A'$, donc $A'$
\'etant fid\`element plat sur~$A$, les~$B_{i}$ sont des modules de
type fini sur~$A$, \ie $B$ est entier sur~$A$. On voit de m\^eme
que si~$B'$ est fini sur~$A'$, $B$ l'est sur~$A$. M\^eme conclusion
pour \og localement libre de type fini\fg, \cf \Ref{VIII.1.11}.

\begin{corollaire}
\label{VIII.5.8}
Sous
\marginpar{215}
les conditions de~\Ref{VIII.5.4}, supposons~$f$ quasi-compact et
soient~$\mathcal{L}$ un Module inversible sur~$X$, et~$\mathcal{L}'$
son image inverse sur~$X'$. Pour que~$\mathcal{L}$ soit ample
(\resp tr\`es ample) relativement \`a~$f$, il faut et il suffit
\ifthenelse{\boolean{orig}}
{que~$f'$}
{que~$\cal{L}'$}
soit ample (\resp tr\`es ample) relativement \`a~$f'$.
\end{corollaire}

Il y a \`a prouver seulement le \og il suffit\fg. L'hypoth\`ese
sur~$\mathcal{L}$ implique en tous cas que~$f'$ est s\'epar\'e,
donc $f$ est s\'epar\'e par~\Ref{VIII.4.8}, et comme~$f$ est
quasi-compact, et~$g\colon Y'\to Y$ est plat, le calcul des images
directes par recouvrements affines montre qu'on a des isomorphismes
$$
g^*(f_*(\mathcal{L}^{\otimes n}))\isomto
f_*'(\mathcal{L}'^{\otimes n})
$$
pour tout entier~$n$, donc on a un isomorphisme
$$
g^*(\mathcal{S})\isomto \mathcal{S}'\quoi,
$$
o\`u~$\mathcal{S}$ (\resp $\mathcal{S}'$) d\'esigne l'Alg\`ebre
gradu\'ee quasi-coh\'erente sur~$Y$ (\resp sur~$Y'$) somme directe
des~$f_*(\mathcal{L}^{\otimes n})$
(\resp des~$f_*'(\mathcal{L}'^{\otimes n})$) pour~$n\geq
0$. Notons que pour tout~$n\geq 0$, le conoyau de l'homomorphisme
canonique~$f'^*(\mathcal{S}_{n}')\to \mathcal{L}'^{\otimes n}$
est l'image inverse par~$X'\to X$ du conoyau
de~$f^*(\mathcal{S}_{n})\to \mathcal{L}^{\otimes n}$, donc son
support~$Z_{n}'$ est l'image inverse du support~$Z_{n}$. Si
\ifthenelse{\boolean{orig}}
{$f$}
{$\cal{L}'$}
est ample l'intersection des~$Z_{n}'$ est vide, donc comme~$X'\to X$
est surjectif, l'intersection des~$Z_{n}$ est vide, \ie on a un
morphisme canonique
$$
j\colon X\to \Proj(\mathcal{S})
$$
(EGA~II 3). D'ailleurs, le morphisme analogue
$$
j'\colon X'\to \Proj(\mathcal{S}')
$$
n'est autre que celui qui est d\'eduit du pr\'ec\'edent par le
changement de base~$Y'\to Y$ (\loccit). Ceci pos\'e, dire que
$\mathcal{L}'$ est ample relativement \`a~$f'$ signifie que~$j'$ est
une immersion, d'ailleurs n\'ecessairement quasi-compacte
puisque~$f'$ est quasi-compact. Donc en vertu de~\Ref{VIII.5.5} $j$
est une immersion, \ie $\mathcal{L}$ est ample relativement
\`a~$f$. --- On proc\`ede de fa\c con toute analogue dans le cas
de \og tr\`es
\marginpar{216}
ample\fg, en se bornant ci-dessus \`a~$n=1$, et en rempla\c cant la
consid\'eration de~$\Proj(\mathcal{S})$ par celle du fibr\'e
projectif~$\mathcal{P}(\mathcal{S}_{1})$ associ\'e
\`a~$\mathcal{S}_{1}$.

Rappelons (EGA~II 5.1.1) qu'un morphisme quasi-compact~$f$ est dit
\emph{quasi-affine} si pour tout ouvert affine~$U$ dans~$Y$,
$f^{-1}(U)$ est un pr\'esch\'ema isomorphe \`a un
sous-sch\'ema ouvert d'un sch\'ema affine. On montre (\loccit)
qu'il revient au m\^eme de dire que~$\mathcal{O}_{X}$ est ample (ou
aussi: tr\`es ample) relativement \`a~$f$. Donc~\Ref{VIII.5.8}
implique:

\begin{corollaire}
\label{VIII.5.9}
Sous les conditions de~\Ref{VIII.5.4}, et supposant~$f$ quasi-compact,
pour que~$f$ soit quasi-affine, il faut et il suffit que~$f'$ le soit.
\end{corollaire}

\begin{remarques}
\label{VIII.5.10}
L'exemple de vari\'et\'e non projective de Hironaka montre qu'on
peut avoir un morphisme propre~$f\colon X\to Y$ de vari\'et\'es
alg\'ebriques non singuli\`eres (avec~$Y$ projective), tel que~$Y$
soit r\'eunion de deux ouverts~$Y_{i}$ tels
que~$X_{i}=X\times_{Y}Y_{i}$ soit projectif sur~$Y_{i}$, mais~$f$
n'\'etant pas projectif. Donc posant~$Y'=Y_{1}\amalg Y_{2}$, $Y'$
est fid\`element plat et quasi-compact (et m\^eme quasi-fini)
sur~$Y$, $f'\colon X'\to Y'$ est projectif, mais~$f$ n'est pas
projectif. Il faut donc faire attention que pour
appliquer~\Ref{VIII.5.8}, et d\'eduire du fait que~$f'$ est
projectif la m\^eme conclusion sur~$f$, il faut disposer d\'ej\`a
sur~$X'$ d'un Module inversible~$\mathcal{L}'$ ample pour~$f'$,
\emph{muni d'une donn\'ee de descente relativement \`a}~$X'\to X$,
(ce qui permet de consid\'erer~$\mathcal{L}'$ comme l'image inverse
d'un Module inversible~$\mathcal{L}$ sur~$X$, qui sera alors ample
pour~$f$ gr\^ace \`a~\Ref{VIII.5.8}). Lorsque~$g\colon S'\to S$
est fini et localement libre, voir cependant~\Ref{VIII.7.7}.
\end{remarques}

\section[Application aux morphismes finis et quasi-finis]{Application aux morphismes finis et quasi-finis\protect\footnotemark}
\label{VIII.6}
\footnotetext{\Cf EGA IV 18.12 pour des g\'en\'eralisations \`a des pr\'esch\'emas non n\'ecessairement localement noeth\'eriens}

Nous allons d\'emontrer les deux th\'eor\`emes suivants:

\begin{theoreme}
\label{VIII.6.1}
Soit~$f\colon X\to Y$ un morphisme \emph{propre \`a fibres finies},
avec~$Y$ localement noeth\'erien. Alors~$f$ est fini.
\end{theoreme}

\begin{theoreme}
\label{VIII.6.2}
Soit~$f\colon X\to Y$ un morphisme \emph{quasi-fini} et
\emph{s\'epar\'e}, avec~$Y$ localement noeth\'erien. Alors~$f$
est quasi-affine, et a fortiori quasi-projectif.
\end{theoreme}

\begin{remarques}
\label{VIII.6.3}
Le
\marginpar{217}
th\'eor\`eme~\Ref{VIII.6.1} est bien connu, et d\^u \`a
Chevalley dans le cas de vari\'et\'es alg\'ebriques; on en
trouvera aussi une d\'emonstration simple dans [EGA III 4],
utilisant le \og th\'eor\`eme des fonctions holomorphes\fg. La
d\'emonstration donn\'ee ici n'utilise pas ce dernier
th\'eor\`eme, mais par contre la th\'eorie de descente; nous
la donnons comme prime au lecteur, car on l'a \og pour rien\fg en
m\^eme temps que celle de~\Ref{VIII.6.2}. Rappelons aussi ([EGA III
4] ou \cite{VIII.1}) que la forme globale du \og Main Theorem\fg de Zariski,
d\'eduit du \og th\'eor\`eme des fonctions holomorphes\fg, affirme
que si~$f\colon X\to Y$ est quasi-fini et \emph{quasi-projectif}, $Y$
\'etant noeth\'erien, alors~$X$ est $Y$-isomorphe \`a un
sous-pr\'esch\'ema ouvert d'un~$Y$-pr\'esch\'ema
\emph{fini}~$Z$. La conjonction du \og Main Theorem\fg et
de~\Ref{VIII.6.2} s'\'enonce donc ainsi:
\end{remarques}

\begin{corollaire}
\label{VIII.6.4}
Soit~$f\colon X\to Y$ un morphisme quasi-fini et s\'epar\'e,
avec~$Y$ noeth\'erien. Alors~$X$ est $Y$-isomorphe \`a un
sous-pr\'esch\'ema ouvert d'un~$Y$-pr\'esch\'ema fini~$Z$.
\end{corollaire}

Une autre cons\'equence int\'eressante de~\Ref{VIII.6.2} pour la
th\'eorie de la descente sera donn\'ee avec~\Ref{VIII.7.9}.

\subsubsection*{D\'emonstration de~\Ref{VIII.6.1} et~\Ref{VIII.6.2}}
Nous admettrons le fait suivant, dont la d\'emonstration est
facile\footnote{\Cf EGA IV 8}:

\begin{lemme}
\label{VIII.6.5}
Soit~$X$ un pr\'esch\'ema de type fini sur~$Y$ localement
noeth\'erien, et soit~$y\in Y$. Pour qu'il existe un voisinage
ouvert~$U$ de~$y$ tel que~$X|U$ soit fini (\resp quasi-affine,
\resp\ ...) sur~$U$, il faut et il suffit
que~$X\times_{Y}\Spec(\mathcal{O}_{y})$ soit fini (\resp quasi-affine,
\resp\ ...) sur~$\Spec(\mathcal{O}_{y})$.
\end{lemme}

Comme d'autre part la propri\'et\'e pour~$f\colon X\to Y$
d'\^etre fini, \resp quasi-affine, est locale sur~$Y$, on est
ramen\'e pour prouver~\Ref{VIII.6.1} et~\Ref{VIII.6.2} au cas
o\`u~$Y$ est le spectre d'un anneau local, et est a fortiori de
dimension finie. (N.B. on appelle \emph{dimension d'un
pr\'esch\'ema}~$Y$ le sup des dimensions de Krull de ses anneaux
locaux). Nous proc\'edons par r\'ecurrence sur~$$
n=\dim(Y)\quoi,
$$
l'assertion
\marginpar{218}
\'etant triviale pour~$n<0$. Nous pouvons donc supposer~$n\geq 0$,
et l'assertion d\'emontr\'ee pour les dimensions~$n'<n$. On peut
\`a nouveau supposer que~$Y$ est le spectre d'un anneau local
noeth\'erien~$A$, de dimension~$n$. Notons que les hypoth\`eses
faites dans~\Ref{VIII.6.1} et~\Ref{VIII.6.2} sont stables par
changement de base (on s'en est d\'ej\`a servi dans la r\'eduction
du d\'ebut), elles resteront vraies apr\`es le changement de
base~$\Spec(\widehat{A})\to \Spec(A)$. Comme ce dernier est
fid\`element plat et quasi-compact, les
\'enonc\'es~\Ref{VIII.5.7} et~\Ref{VIII.5.9} nous ram\`enent au
cas o\`u~$A$ est de plus complet. Utilisant alors le fait que tout
anneau local noeth\'erien~$B$ sur~$A$ quasi-fini sur~$A$ est fini
sur~$A$, et le fait que~$X$ est s\'epar\'e sur~$Y$ et la fibre
de~$y$ est form\'ee de points isol\'es, on trouve une
d\'ecomposition
$$
X=X'\amalg X''
$$
o\`u~$X'$ est \emph{fini} sur~$Y$, et o\`u la fibre de~$X''$
en~$y$ est vide. Si~$X$ est propre sur~$Y$, il en est de m\^eme
de~$X''$, donc son image dans~$Y$ est ferm\'ee, et comme elle ne
contient pas~$y$, elle est vide, donc~$X''=\emptyset$ donc~$X=X'$,
ce qui montre que~$X$ est fini sur~$Y$ et d\'emontre~\Ref{VIII.6.1}
(N.B. l'hypoth\`ese de r\'ecurrence est ici inutile). Si~$X$ est
quasi-fini sur~$Y$, $X''$ l'est aussi, or~$X''$ se trouve en fait sur
\ifthenelse{\boolean{orig}}
{l'ouvert~$Y-(y)$}
{l'ouvert~$Y-\{ y\}$}
de~$Y$, \emph{qui est de dimension} $<n$. En vertu de l'hypoth\`ese
de r\'ecurrence, $X''$ est quasi-affine
\ifthenelse{\boolean{orig}}
{sur~$Y-(y)$,}
{sur~$Y-\{ y\}$,}
donc aussi sur~$Y$, il en est \'evidemment de m\^eme de~$X'$, donc
aussi de leur somme~$X$, ce qui prouve~\Ref{VIII.6.2}.

\begin{remarque}
\label{VIII.6.6}
Les th\'eor\`emes~\Ref{VIII.6.1} et~\Ref{VIII.6.2} restent
valables si on ne suppose plus~$Y$ localement noeth\'erien, \`a
condition de sp\'ecifier que l'on suppose~$f$ de pr\'esentation
finie
\ifthenelse{\boolean{orig}}
{(\cf \Ref{VIII.3.6}).}
{(\cf \Ref{VIII.3.6}).}
En effet, on peut encore supposer~$Y$ affine, et alors on v\'erifie
sans difficult\'e que la situation~$f\colon X\to Y$ est
\ifthenelse{\boolean{orig}}
{d\'eduit,}
{d\'eduite,}
par un changement de base~$Y\to Y_{0}$, d'une situation~$f_{0}\colon
X_{0}\to Y_{0}$ satisfaisant les m\^emes hypoth\`eses que~$f$,
avec~$Y_{0}$ \emph{noeth\'erien}. Donc d'apr\`es le
r\'esultat~\Ref{VIII.6.1} \resp \Ref{VIII.6.2}, $f_{0}$ est fini
\resp quasi-affine, donc il en est de m\^eme de~$f$. Ce genre de
raisonnement est souvent utile pour se d\'ebarrasser
d'hypoth\`eses noeth\'eriennes,
\marginpar{219}
(qui finissent toujours par \^etre g\^enantes dans les
applications).
\end{remarque}

\section{Crit\`eres d'effectivit\'e pour une donn\'ee de descente}
\label{VIII.7}

Consid\'erons comme d'habitude un morphisme de pr\'esch\'emas
$$
g\colon S'\to S
$$
et un $S'$-pr\'esch\'ema~$X'$. Conform\'ement aux faits
\ifthenelse{\boolean{orig}}
{g\'en\'eraux~VII,~9,}
{g\'en\'eraux~VII,~9\footnote{n'existe pas; voir expos\'e 212 du
S\'eminaire Bourbaki (1962)},}
la donn\'ee d'une donn\'ee de descente sur~$X'$, relativement
\`a~$g$, est \'equivalente \`a la donn\'ee d'un couple
d'\'equivalence \cite{VIII.3}:
$$
\xymatrix@C=.5cm{ q_{1},q_{2}\colon X''\ar@<2pt>[r]\ar@<-2pt>[r] & X'}
$$
tel que le morphisme structural~$X'\to S'$ soit compatible avec ce
couple et le couple d'\'equivalence
$$
\xymatrix@C=.5cm{p_{1},p_{2}\colon
S''=S'\times_{S}S'\ar@<2pt>[r]\ar@<-2pt>[r] & S'}
$$
d\'efini par~$g$, et tel que les deux carr\'es (ou l'un des deux,
cela revient au m\^eme par raison de sym\'etrie) extraits du
diagramme correspondant
$$
\xymatrix{X'\ar[d] & X''\ar@<2pt>[l]\ar@<-2pt>[l]\ar[d]\\ S' &
S''\ar@<2pt>[l]\ar@<-2pt>[l]}
$$
(en utilisant soit~$p_{1},q_{1}$, soit~$p_{2},q_{2}$), soit
\emph{cart\'esien}. Une solution du probl\`eme de descente
pos\'e par cette donn\'ee de descente, \ie un objet~$X$ sur~$S$
muni d'un isomorphisme~$X\times_{S}S'\isomfrom X'$ compatible avec les
donn\'ees de descente, \'equivaut \`a la donn\'ee d'un
carr\'e \emph{cart\'esien}
$$
\xymatrix{X\ar[d] & X'\ar[l]_{h}\ar[d]\\ S & S'\ar[l]_{g}}
$$
satisfaisant~$hq_{1}=hq_{2}$.

Comme
\marginpar{220}
l'ensemble des morphismes fid\`element plats et quasi-compacts est
stable par changement de base, et qu'un morphisme fid\`element plat
et quasi-compact est un \'epimorphisme effectif en vertu
de~\Ref{VIII.5.3}, il r\'esulte de la th\'eorie g\'en\'erale
\cite{VIII.D}:

\begin{proposition}
\label{VIII.7.1}
Supposons~$g\colon S'\to S$ fid\`element plat et quasi-compact. Pour
qu'une don\-n\'ee de descente sur~$X'$ relativement \`a~$g$ soit
effective, il faut et il suffit que la relation
d'\'equivalence~$R=(q_{1},q_{2})$ qu'elle d\'efinit soit effective
(\ie le quotient~$X'/R$ existe et $X''$ devient le carr\'e
fibr\'e de~$X'$ sur~$X'/R$), et que le morphisme canonique~$X'\to
X'/R$ soit fid\`element plat et quasi-compact.
\end{proposition}

Ainsi la question de l'effectivit\'e d'une donn\'ee de descente
est un cas particulier de la question d'effectivit\'e d'un graphe
d'\'equivalence, et divers crit\`eres d'effectivit\'e donn\'es
dans ce num\'ero peuvent s'obtenir de cette fa\c con. N\'eanmoins on dispose dans le contexte de la descente du
th\'eor\`eme~\Ref{VIII.2.1}, qui implique que \emph{si $X'$ est
affine sur~$S'$, toute donn\'ee de descente sur~$X'$ relativement
\`a~$g$ est effective}, \'enonc\'e qui n'a pas d'analogue pour
le passage au quotient par un graphe d'\'equivalence plat
g\'en\'eral. Tous les crit\`eres d'effectivit\'e que nous
donnons ici peuvent aussi \^etre consid\'er\'es comme
d\'eduits du pr\'ec\'edent.

Soit~$U'$ un sous-pr\'esch\'ema de~$X'$ (ou plus
g\'en\'eralement un sous-objet de~$X'$ dans la
cat\'egorie~$\Sch$); on dit que~$U'$ est \emph{stable par la
donn\'ee de descente} sur~$X'$, si on peut trouver sur~$U'$ une
donn\'ee de descente relativement \`a~$g$, telle que
l'immersion~$U'\to X'$ soit compatible avec les donn\'ees de
descente. Cela signifie aussi que les images inverses de~$U'$
dans~$X''$ par~$q_{1}$ et~$q_{2}$ sont les m\^emes (ou aussi, comme
on dit, que~$U'$ est \emph{stable par la relation
d'\'equivalence}~$R$), et bien entendu la donn\'ee de descente en
question sur~$U'$ est alors unique, et dite \emph{donn\'ee de
descente induite}
\index{donn\'ee de descente induite|hyperpage}%
par celle de~$X'$. Ceci pos\'e:

\begin{proposition}
\label{VIII.7.2}
Soit~$(X_{i}')$ un recouvrement de~$X'$ par des ouverts~$X_{i}'$
stables
\marginpar{221}
%
par la donn\'ee de descente. Pour que la donn\'ee de descente
sur~$X'$ soit effective, il faut et il suffit qu'il en soit de
m\^eme des donn\'ees de descente induites dans les~$X_{i}'$.
\end{proposition}

C'est l\`a une cons\'equence facile de~\Ref{VIII.7.1} par exemple,
et le d\'etail de la d\'emonstration est laiss\'e au lecteur.

\begin{corollaire}
\label{VIII.7.3}
Soit~$(S_{i})$ un recouvrement ouvert de~$S$, et pour tout~$i$
soient~$S_{i}'$ et~$X_{i}'$ d\'eduits de~$S'$
\ifthenelse{\boolean{orig}}
{de}
{et}
$X'$ par le changement de base~$S_{i}\to S$. Pour que la donn\'ee de
descente sur~$X'$ soit effective, il faut et il suffit que, pour
tout~$i$, la donn\'ee de descente sur~$X_{i}'$ relativement
\`a~$g_{i}\colon S_{i}'\to S_{i}$ soit effective.
\end{corollaire}

Ce crit\`ere nous ram\`ene toujours pratiquement au cas o\`u~$S$
est affine. Dans le cas o\`u~$S'$ est \'egalement affine, ce qui
est le cas le plus fr\'equent dans les applications, on~a:

\begin{corollaire}
\label{VIII.7.4}
Supposons~$S$ et~$S'$ affines. Pour que la donn\'ee de descente
sur~$X'$ soit effective, il faut et il suffit que~$X'$ soit
r\'eunion d'ouverts~$X_{i}'$ affines et stables par la donn\'ee de
descente.
\end{corollaire}

La suffisance provient de~\Ref{VIII.7.2} et du fait que si~$X_{i}'$
est affine, il est affine sur~$S'$ et on peut
appliquer~\Ref{VIII.2.1}. Pour la n\'ecessit\'e, on note que
si~$X'$ provient de~$X$, et si~$X$ est recouvert par des ouverts
affines~$X_{i}$, alors les~$X_{i}'=X_{i}\times_{S}S'$ sont des ouverts
affines stables par les donn\'ees de descente et recouvrant~$X'$.

\begin{corollaire}
\label{VIII.7.5}
Soit~$g\colon S'\to S$ un morphisme fid\`element plat, quasi-compact
et \emph{radiciel}. Alors~$g$ est un morphisme de descente
\emph{effective}, \ie pour tout~$X'$ sur~$S'$, toute donn\'ee de
descente sur~$X'$ relativement \`a~$g\colon S'\to S$ est effective.
\end{corollaire}

En effet, en vertu de~\Ref{VIII.7.3} on peut supposer~$S$ affine, donc
comme~$S'$ est radiciel sur~$S$ donc s\'epar\'e, $S'$ est
s\'epar\'e. D'ailleurs pour tout~$x'\in X'$, la
fibre~$R(x')=q_{2}(q_{1}^{-1}(x'))$ de la relation d'\'equivalence
ensembliste d\'efinie
\marginpar{222}
par la relation d'\'equivalence~$R$ est r\'eduite \`a un point,
car~$g$ \'etant radiciel, il en est de m\^eme de~$p_{1},p_{2}$ qui
s'en d\'eduisent par changement de base~$S'\to S$, donc aussi
de~$q_{1},q_{2}$ qui se d\'eduisent des pr\'ec\'edents par
changement de base~$X'\to S''$. Donc \emph{tout ouvert de~$X'$ est
stable} par la donn\'ee de descente. Recouvrons alors~$X'$ par des
ouverts affines~$X_{i}'$, ils sont affines sur~$S$ puisque~$S'$ est
s\'epar\'e, donc la donn\'ee de descente induite est effective
par~\Ref{VIII.2.1}. On conclut alors par~\Ref{VIII.7.2}.

On notera que~\Ref{VIII.7.5} donne le seul cas connu d'un morphisme de
descente effective dans la cat\'egorie des pr\'esch\'emas, et
c'est probablement le seul cas en effet, m\^eme en se limitant aux
sch\'emas noeth\'eriens, ou aux sch\'emas de type fini sur un
corps.

Lorsqu'on suppose~$S$ localement noeth\'erien et~$S'$ de type fini
sur~$S$, l'\'enonc\'e~\Ref{VIII.7.5} est aussi un cas particulier
du suivant (qui g\'en\'eralise la descente galoisienne de Weil et
la descente ins\'eparable de Cartier):

\begin{corollaire}
\label{VIII.7.6}
Soit~$g\colon S'\to S$ un morphisme fini localement libre
(\ie d\'efini par une Alg\`ebre sur~$S$ qui est un module
localement libre de type fini) et surjectif (donc~$g$ est
fid\`element plat et quasi-compact, donc un morphisme de
descente). Soit~$X'$ un~$S'$-pr\'esch\'ema muni d'une donn\'ee
de descente. Pour que cette donn\'ee soit effective, il faut et il
suffit que pour tout~$x'\in X'$, la
fibre~$R(x')=q_{2}(q_{1}^{-1}(x'))$ soit contenue dans un ouvert
affine (condition automatiquement v\'erifi\'ee si~$X'$ est
quasi-projectif sur~$S'$).
\end{corollaire}

La remarque entre parenth\`eses provient du fait que si~$s$ est le
point de~$S$ au-dessous de~$x'$, alors $R(x')$ est fini et contenu
dans la
\ifthenelse{\boolean{orig}}
{fibre~$X_{s}$,}
{fibre~$X'_{s}$,}
d'autre part comme~$X'$ est quasi-projectif sur~$S'$ et~$S'$ est fini
sur~$S$, $X'$ est quasi-projectif sur~$S$, ce qui implique qu'une
fibre
\ifthenelse{\boolean{orig}}
{de~$X/S$}
{de~$X'$ sur~$S$}
est contenue dans un ouvert affine.

Comme toute partie finie d'un sch\'ema affine admet un syst\`eme
fondamental de voisinages affines, on voit qu'on ne perd pas
l'hypoth\`ese en se restreignant au-dessus d'un ouvert affine
de~$S$, ce qui en vertu de~\Ref{VIII.7.3} nous ram\`ene au cas
o\`u~$S$ est affine. En vertu de~\Ref{VIII.7.4}, on est ramen\'e
\`a montrer
\ifthenelse{\boolean{orig}}
{que~$X'$}
{que~$x'$}
est contenu
\marginpar{223}
dans un ouvert affine \emph{stable} par la donn\'ee de
descente. Soit en effet~$U$ un ouvert affine contenant~$R(x')$, alors
le satur\'e
$$
\ifthenelse{\boolean{orig}}
{R(X-U)=q_{2}(q_{1}^{-1}(X-U))}
{R(X'-U)=q_{2}(q_{1}^{-1}(X'-U))}
$$
ne rencontre pas~$R(x')$, d'autre part comme~$q_{2}$ est fini (car~$g$
donc~$p_{2}$ l'est) donc ferm\'e, le deuxi\`eme membre est une
partie ferm\'ee de~$X'$. Soit~$U'$ son compl\'ementaire dans~$X'$,
c'est donc un ouvert \emph{satur\'e} et on a
$$
R(x')\subset U'\subset U
$$
avec~$U$ affine, mais~$U'$ pas affine a priori. Comme une partie
finie~$R(x')$ dans un sch\'ema affine~$U$ a un syst\`eme
fondamental de voisinages affines de la forme~$U_{f}$, on voit,
rempla\c cant~$f$ par sa restriction \`a~$U'$, qu'il existe une
section~$f$ de~$\mathcal{O}_{U}$ telle que:
$$
R(x')\subset U_{f}',\quad U_{f}'\textrm{ est affine}.
$$
Soit alors~$U''=q_{1}^{-1}(U')=q_{2}^{-1}(U')$, d\'esignons encore
par~$q_{1},q_{2}$ les morphismes induits~$U''\to U'$, et
consid\'erons
$$
f'=\Norm_{q_{2}}(q_{1}^*(f))\quoi,
$$
o\`u~$\Norm_{q_{2}}$ d\'esigne la \emph{norme} relativement au
morphisme fini localement libre~$q_{2}\colon U''\to U'$. La
compatibilit\'e de la formation de la norme avec le changement de
base implique facilement que~$f'$ est une section \emph{invariante}:
$$
q_{1}^*(f')=q_{2}^*(f')
$$
ce qui implique que~$U_{f'}'$ est un ouvert satur\'e de~$U'$. De
fa\c con plus pr\'ecise d'ailleurs, d\'esignant par~$Z(f')$
l'ensemble des z\'eros d'une section~$f'$, on trouve en vertu des
propri\'et\'es des normes:
$$
Z(f')=q_{2}(Z(q_{1}^*(f)))=q_{2}(q_{1}^*(Z(f)))=
R(U'-U_{f}')\quoi,
$$
ce
\marginpar{224}
qui implique que~$U_{f'}'=U'-Z(f')$ est satur\'e, contient~$R(x')$,
et est contenu dans~$U_{f}'$. Comme ce dernier est affine, il s'ensuit
que~$U_{f'}'$ l'est aussi (car \'egal \`a~$(U_{f}')_{f''}$,
avec~$f''=f'|U_{f'}'$). C'est donc un ouvert affine satur\'e
contenant~$R(x')$ donc~$x'$, ce qui ach\`eve la d\'emonstration.

On notera que ce raisonnement s'applique chaque fois qu'on a une
relation d'\'equivalence (ou m\^eme seulement de
pr\'e\'equivalence, \cf \cite{VIII.3}) dans un pr\'esch\'ema~$X'$, fini
et localement libre et d'ailleurs~\Ref{VIII.7.6} est aussi un cas
particulier du r\'esultat analogue pour les pr\'e\'equivalences
finies et localement libres, \cf \loccit M\^eme remarque
pour~\Ref{VIII.7.7} ci-dessous.

On peut aussi, une fois obtenue l'existence d'un ouvert quasi-affine
satur\'e~$U'$ contenant~$x'$, faire appel \`a~\Ref{VIII.7.9}
et~\Ref{VIII.7.2} ce qui \'evite le recours aux normes.

Notons d'ailleurs que sous les conditions de~\Ref{VIII.7.6}, si la
donn\'ee de descente sur~$X'$ est effective, $X'$ provenant de~$X$
sur~$S$, alors le morphisme~$X'\to X$ est fini, localement libre et
surjectif, car d\'eduit de~$g$ par le changement de base~$X\to
S$. Il s'ensuit (EGA~II 6.6.4) que si~$X'$ est quasi-projectif
sur~$S'$ donc sur~$S$, alors~$X$ est quasi-projectif sur~$S$, (un
faisceau inversible relativement ample sur~$X$ \'etant obtenu en
prenant la \emph{norme} d'un faisceau inversible sur~$X'$ relativement
ample sur~$S$, ou sur~$S'$, cela revient au m\^eme). On obtient
ainsi:

\begin{corollaire}
\label{VIII.7.7}
Un morphisme~$g\colon S'\to S$ fini localement libre et surjectif est
un morphisme de descente effective pour la cat\'egorie fibr\'ee
des pr\'esch\'emas quasi-projectifs sur d'autres, \ie pour
tout~$X'$ quasi-projectif sur~$S'$, toute donn\'ee de descente
sur~$X'$ relativement \`a~$g$ est effective, et le
$S$-pr\'esch\'ema descendu~$X$ est quasi-projectif sur~$S$.
\end{corollaire}

\begin{proposition}
\label{VIII.7.8}
Soit~$g\colon S'\to S$ un morphisme fid\`element plat et
quasi-compact. Alors~$g$ est un morphisme de descente effective pour
la cat\'egorie fibr\'ee des pr\'esch\'emas~$Z$ quasi-compacts
sur un pr\'esch\'ema~$T$, munis d'un faisceau inversible ample
\marginpar{225}
relativement \`a~$T$. En particulier, pour tout
pr\'esch\'ema~$X'$ sur~$S'$, muni d'une donn\'ee de descente
relativement \`a~$g\colon S'\to S$, et tout faisceau
inversible~$\mathcal{L}'$ sur~$X'$ ample relativement \`a~$S'$, muni
\'egalement d'une donn\'ee de descente relativement \`a celle
donn\'ee sur~$X'$, (\ie muni d'un isomorphisme
de~$q_{1}^*(\mathcal{L}')$ avec~$q_{2}^*(\mathcal{L}')$,
satisfaisant la condition de transitivit\'e habituelle), la
donn\'ee de descente sur~$X'$ est effective, et le faisceau
inversible~$\mathcal{L}$ sur le pr\'esch\'ema descendu~$X$,
d\'eduit de~$\mathcal{L}'$ par descente, est ample relativement
\`a~$S$.
\end{proposition}

La d\'emonstration est toute analogue \`a celle de~\Ref{VIII.5.8},
en notant que sur l'Alg\`ebre gradu\'ee
quasi-coh\'erente~$\mathcal{S}'$ sur~$S'$ d\'efinie
par~$\mathcal{L}'$, il y a une donn\'ee de descente, permettant de
construire une Alg\`ebre gradu\'ee
quasi-coh\'erente~$\mathcal{S}$ sur~$S$ gr\^ace
\`a~\Ref{VIII.1.1} d'o\`u un~$P=\Proj(\mathcal{S})$ sur~$S$ tel
que~$P'=\Proj(\mathcal{S}')$ s'identifie avec sa donn\'ee de
descente \`a~$P\times_{S}S'$. Comme par hypoth\`ese $X'$
s'identifie \`a un ouvert de~$P'$, n\'ecessairement stable par la
donn\'ee de descente sur~$P'$, la donn\'ee de descente sur~$X'$
est \'egalement effective, et on obtient le pr\'esch\'ema
descendu comme un ouvert dans~$P$. Le d\'etail est laiss\'e au
lecteur. --- En particulier, faisant~$\mathcal{L}'=\mathcal{O}_{X'}$, on
trouve:

\begin{corollaire}
\label{VIII.7.9}
Soit~$g\colon S'\to S$ un morphisme fid\`element plat et
quasi-compact, et soit~$X'$ un pr\'esch\'ema \emph{quasi-affine}
au-dessus de~$S'$, alors toute donn\'ee de descente sur~$X'$
relativement \`a~$g$ est effective, et le pr\'esch\'ema
descendu~$X$ est quasi-affine sur~$S$.
\end{corollaire}

En vertu de~\Ref{VIII.6.2}, ce r\'esultat s'applique en particulier
si~$S'$ est localement noeth\'erien et~$X'$ est quasi-fini et
s\'epar\'e sur~$S'$, plus g\'en\'eralement si~$S'$ est
quelconque et~$X'$ est de pr\'esentation finie, quasi-fini et
s\'epar\'e sur~$S'$ (\cf \Ref{VIII.6.6}).

\begin{remarques}
\label{VIII.7.10}
Les r\'esultats donn\'es dans ce num\'ero \'epuisent les
crit\`eres actuellement connus d'effectivit\'e, et probablement
m\^eme les crit\`eres utiles existants\footnote{Cette opinion
s'est trouv\'ee partiellement erron\'ee, voir p\ptbl ex. J.P\ptbl \textsc{Murre},
S\'em. Bourbaki 294 (Appendix), Mai 1965, et des r\'esultats
sp\'eciaux (notamment de \textsc{N\'eron} et \textsc{Raynaud}). Pour la descente des
sch\'emas en groupes; \cf M\ptbl Raynaud, th\`ese (1968).}. On notera
les contre-exemples suivants \`a l'appui de cette assertion:
\begin{enumerate}
\item[(i)]
\marginpar{226}
Si~$S$ est le spectre d'un corps, et~$S'$ est le spectre d'une
extension quadratique galoisienne, on peut trouver un~$X'$ sur~$S'$
propre et lisse sur~$S'$, de dimension 3, muni d'une donn\'ee de
descente qui n'est pas effective (Serre).
\item[(ii)] On peut trouver un~$S$ spectre d'un anneau local
r\'egulier de dimension 3 (si on veut, l'anneau local d'un
sch\'ema alg\'ebrique sur un corps de caract\'eristique
donn\'ee), un~$T$ rev\^etement principal de~$S$ de
groupe~$\ZZ/2\ZZ$, tel que, si~$t$ d\'esigne l'un des points de~$T$
au-dessus du point ferm\'e~$s$ de~$S$, et $S'=T-s$, on puisse
trouver un~$X'$ \emph{projectif} sur~$S'$, r\'egulier, muni d'une
donn\'ee de descente relativement \`a~$g\colon S'\to S$, cette
donn\'ee de descente n'\'etant pas effective.
\end{enumerate}

On utilise pour ces constructions l'exemple de Hironaka de
vari\'et\'es non projectives. Pour~(i), il suffit d'utiliser le
fait qu'on peut trouver au-dessus de~$k$ un sch\'ema propre et
lisse~$X_{0}$ de dimension 3, sur lequel~$G=\ZZ/2\ZZ$ op\`ere sans
inertie, et dans lequel il existe deux points~$a,b$ rationnels
\ifthenelse{\boolean{orig}}
{sur~$P$,}
{sur~$S$,}
congrus sous~$G$, qui ne sont pas contenus dans un ouvert affine. On
pose alors~$X'=X_{0}\times_{k}k'$, on fait op\'erer~$G$ sur~$X'$
gr\^ace aux op\'erations de~$G$ sur les deux facteurs, ce qui
donne une donn\'ee de descente sur~$X'$ relativement \`a~$g\colon
\Spec(k')\to \Spec(k)$. Au-dessus de~$a$ \resp $b$, il y a exactement
un point~$a'$ \resp $b'$, (\`a extension r\'esiduelle
quadratique), et~$a'$ et~$b'$ sont congrus sous~$G$, puisque~$X'\to
X_{0}$ est compatible avec les op\'erations de~$G$. Alors~$a'$
et~$b'$ ne peuvent \^etre contenus dans un ouvert affine, soit~$U'$,
car alors~$U=X_{0}-\Im(X'-U')$ serait un ouvert de~$X_{0}$
contenant~$(a,b)$ et dont l'image inverse dans~$X'$ serait contenue
dans~$U'$, donc quasi-affine, donc~$U$ serait quasi-affine, et par
suite~$(a,b)$ aurait un voisinage affine dans~$U$.

Pour~(ii), on utilise le fait que dans l'exemple de Hironaka, $X_{0}$
est obtenu comme pr\'esch\'ema propre au-dessus
d'un~$k$-sch\'ema projectif~$Y$, lisse sur~$k$ (le
morphisme~$f\colon X_{0}\to Y$ \'etant d'ailleurs birationnel, mais
peu importe), le groupe~$G$ op\'erant \'egalement sur~$Y$ de
fa\c con compatible avec ses op\'erations sur~$X_{0}$, enfin
posant~$S'=Y-f(b)$, $X'=X_{0}|S'$, $X'$ est projectif
sur~$S'$. Alors~$X_{0}$ est muni d'une donn\'ee de descente
naturelle relative au morphisme
\marginpar{227}
canonique~$Y\to S=Y/G$, gr\^ace aux op\'erations de~$G$
sur~$X_{0}$ compatibles avec ses op\'erations sur~$Y$. Cette
donn\'ee de descente n'est pas effective, puisque~$(a,b)$ n'est pas
contenu dans un ouvert affine. La donn\'e de descente induite
sur~$X'$ relativement \`a~$g\colon S'\to S$ n'est alors pas
effective, comme on v\'erifie facilement.
\end{remarques}


\chapter[Descente des morphismes \'etales]{Descente des morphismes \'etales.\\ Application au groupe fondamental}
\label{IX}
\marginpar{228}

\section{Rappels sur les morphismes \'etales}
\label{IX.1}

Nous allons ici passer en revue les propri\'et\'es des morphismes
\'etales (d\'evelopp\'es dans~\Ref{I}) qui vont nous servir,
en profitant de cette occasion pour \'eliminer de la th\'eorie les
hypoth\`eses noeth\'eriennes superflues. Le lecteur notera que
m\^eme si on ne s'int\'eresse qu'aux sch\'emas noeth\'eriens, la
technique de descente conduit \`a introduire des sch\'emas non
noeth\'eriens (tels que~$\Spec(\widehat{A}\otimes_{A}\widehat{A})$,
o\`u~$A$ est un anneau local noeth\'erien), et pour pouvoir
appliquer le langage des cat\'egories fibr\'ees, il importe de
d\'efinir les notions telles que morphisme \'etale \ifthenelse{\boolean{orig}}{etc...}{etc.}, sans y
introduire de restriction noeth\'erienne. Le lecteur qui
r\'epugnerait \`a v\'erifier ou \`a admettre que les
\'enonc\'es ci-dessous sont vrais sans hypoth\`eses
noeth\'eriennes, pourra se contenter de les admettre sous les
hypoth\`eses noeth\'eriennes de~\Ref{I}, \`a condition
d'introduire ces m\^emes hypoth\`eses noeth\'eriennes dans les
\'enonc\'es des num\'eros suivants, et d'utiliser la
d\'efinition~\Ref{IX.1.1} ci-dessous pour les sch\'emas non
noeth\'eriens qui s'introduisent dans les raisonnements.

\begin{definition}
\label{IX.1.1}
Soient~$f\colon X\to S$ un morphisme de pr\'esch\'emas, et~$x$ un
point de~$X$. On dit que~$f$ est \emph{\'etale en}~$x$,
\index{etale (morphisme, algebre)@\'etale (morphisme, alg\`ebre)|hyperpage}%
ou que~$X$ est \'etale sur~$S$ en~$x$, s'il existe un voisinage
ouvert affine~$U$ de~$s=f(x)$, un voisinage ouvert affine~$V$ de~$x$
au-dessus de~$U$, un sch\'ema \emph{noeth\'erien} affine~$U_{0}$,
un~$U_{0}$-sch\'ema \emph{\'etale}~\eqref{I} et affine~$V_{0}$, un
morphisme~$U\to U_{0}$, et un~$U$-isomorphisme
$$
V\isomto V_{0}\times_{U_{0}}U\quoi.
$$
\end{definition}

On
\marginpar{229}
notera que lorsque~$S$ est localement noeth\'erien, cette
terminologie co\"incide avec celle de \loccit On dira de m\^eme
que~$f$ \emph{est \'etale}, ou que~$X$ \emph{est \'etale sur}~$S$,
si~$f$ est \'etale en tout point~$x$ de~$X$. Avec ces
d\'efinitions, les propositions ci-dessous se ram\`enent sans
difficult\'e au cas noeth\'erien, o\`u elles sont
d\'emontr\'ees dans~I~\No \Ref{I.4},
\Ref{I.5},~\Ref{I.7}. Pour des d\'etails, le lecteur pourra
consulter EGA~IV\footnote{De fa\c con pr\'ecise, EGA~IV 17,~18.}.

\begin{remarques}
\label{IX.1.2}
Si $f$ est \'etale en $x$, alors $f$ est \og \emph{de pr\'esentation
finie en}~$x$\fg (VIII~\Ref{VIII.3.5}), l'anneau local de~$x$ dans la
fibre $f^{-1}(s)$ est une \emph{extension finie s\'eparable} de
$\kres(s)$, enfin $f$ est \emph{plat en}~$x$. On peut montrer que la
r\'eciproque est vraie, donc que la d\'efinition~\Ref{IX.1.1} est
la m\^eme que dans le cas o\`u $S$ est localement noeth\'erien,
sauf qu'il faut remplacer la condition \og de type fini en~$x$\fg par
\og de pr\'esentation finie en~$x$\fg. Comme ce r\'esultat est de
d\'emonstration d\'elicate, nous n'avons pas voulu ici donner
cette d\'efinition de la notion de morphisme \'etale, qui ne se
pr\^ete pas directement \`a la d\'emonstration des
propri\'et\'es qui vont suivre.
\end{remarques}

Notons d'abord qu'on a trivialement:
\begin{proposition}
\label{IX.1.3}
Si $f\colon X\to S$ est \'etale, alors tout morphisme $f'\colon
X'\to S'$ qui s'en d\'eduit par changement de base $S'\to S$ est
\'egalement \'etale.
\end{proposition}

On peut donc dire que les morphismes \'etales forment une
\emph{sous-cat\'egorie fibr\'ee} de la cat\'egorie des
fl\`eches dans~$\Sch$
\ifthenelse{\boolean{orig}}
{(\cf VI~\Ref{VI.11}~a)}
{(\cf VI~\Ref{VI.11}~a))}. L'objet du pr\'esent
expos\'e est l'\'etude des propri\'et\'es d'exactitude de
cette cat\'egorie fibr\'ee sur~$\Sch$.
\begin{proposition}
\label{IX.1.4}
Soit $f\colon X\to S$ un morphisme de pr\'esch\'emas. Pour que ce
soit une immersion ouverte, il faut et il suffit qu'il soit \'etale
et radiciel.
\end{proposition}

\Cf I~\Ref{I.5.1}. On en conclut que si $X$ est \'etale sur~$S$, toute section de~$X$ sur~$S$ est une immersion ouverte, donc,
utilisant encore~\Ref{IX.1.4}, on trouve:
\begin{corollaire}
\label{IX.1.5}
Soit $X$ un $S$-pr\'esch\'ema \'etale. Alors il y a une
correspondance biunivoque entre l'ensemble des sections de~$X$
sur~$S$, et l'ensemble des parties
\marginpar{230}
ouvertes $\Gamma$ de~$X$ telles que le morphisme $\Gamma \to S$ induit
par le morphisme structural soit \emph{radiciel} et \emph{surjectif}.
\end{corollaire}

Si d'ailleurs $X$ est s\'epar\'e sur~$S$, $\Gamma$ sera une partie
de~$X$ \`a la fois ouverte et ferm\'ee, mais peu importe. --- Faisant un changement de base \'evident, on peut mettre~\Ref{IX.1.5}
sous la forme en apparence plus g\'en\'erale:
\begin{corollaire}
\label{IX.1.6}
Soient $X$ et $Y$ deux $S$-pr\'esch\'emas, $Y$ \'etant \'etale
sur~$S$. Alors l'application $f \mto \Gamma_f$ qui associe \`a
tout $S$-morphisme $f$ de~$X$ dans $Y$ la partie de~$X\times_S Y$
sous-jacente au graphe de~$f$, est une bijection de~$\Hom_S(X,Y)$ sur
l'ensemble des parties ouvertes de~$\Gamma$ de~$X\times_S Y$ telles
que le morphisme $\Gamma\to X$ induit par $\pr_1$ soit
\emph{radiciel} et \emph{surjectif}.
\end{corollaire}

\begin{proposition}
\label{IX.1.7}
Soit $S_0$ le sous-pr\'esch\'ema de~$S$ d\'efini par un
Nil-id\'eal quasi-coh\'erent, \ie tel que $S_0$ ait m\^eme
ensemble sous-jacent que~$S$. Alors le foncteur $X \mto X\times_S
S_0$ de la cat\'egorie des pr\'esch\'emas \'etales sur~$S$
dans la cat\'egorie des pr\'esch\'emas \'etales sur~$S_0$, est
une \'equivalence de cat\'egories.
\end{proposition}

Le fait que ce foncteur soit pleinement fid\`ele est une
cons\'equence imm\'ediate de~\Ref{IX.1.6}. Le fait qu'il soit
essentiellement surjectif est contenu dans~I~\Ref{I.8.3}. On notera
que dans l'\'equivalence pr\'ec\'edente,
$X$ est de type fini \ie quasi-fini sur~$S$, (\resp fini \ie un rev\^etement \'etale
de~$S$), si et seulement si $X_0$ satisfait \`a la condition
analogue sur~$S_0$; m\^eme remarque pour la condition de
s\'eparation. Ces faits sont imm\'ediats, et aussi contenus
dans~\Ref{IX.2.4} plus bas.

\begin{corollaire}
\label{IX.1.8}
Soit $A$ un anneau local noeth\'erien complet de corps r\'esiduel
$k$. Alors le foncteur $B \mto B\otimes_A k$ est une
\'equivalence de la cat\'egorie des alg\`ebres finies et
\'etales sur~$A$, avec la cat\'egorie des alg\`ebres finies et
\'etales sur~$k$, (\ie compos\'ees d'un nombre fini d'extensions
finies s\'eparables de~$k$).
\end{corollaire}

\begin{proposition}
\label{IX.1.9}
Pour que $X$ soit un \emph{rev\^etement \'etale} de~$S$ \ie fini
et \'etale sur~$S$, il faut et il suffit que $X$ soit $S$-isomorphe
au spectre d'une Alg\`ebre $\cal{A}$ sur~$S$, qui soit un Module
localement libre de type fini, et telle que pour tout
\marginpar{231}
$s \in S$, $\cal{A}_s\otimes_{\cal{O}_s}\kres(s)$ soit une alg\`ebre
s\'eparable sur~$\kres(s)$, donc en l'occurrence compos\'ee directe
d'extensions finies s\'eparables de~$\kres(s)$.
\end{proposition}

Enfin le r\'esultat suivant est de nature moins \'el\'ementaire,
\'etant la conjonction de~I~\Ref{I.8.4} et du
\emph{th\'eor\`eme d'existence de faisceaux en g\'eom\'etrie
alg\'ebrique}, (EGA~III~5; \cf aussi \cite{IX.1} th\ptbl 3).

\begin{theoreme}
\label{IX.1.10}
Soient $S$ le spectre d'un anneau local noeth\'erien complet, $X$ un
$S$-sch\'ema propre, $X_0$ la fibre de~$X$ au point ferm\'e de
$S$, (de sorte que $X_0$ est un sous-sch\'ema ferm\'e de
$X$). Alors le foncteur restriction $X' \mto X'\times_X X_0$ est
une \'equivalence de la cat\'egorie des rev\^etements \'etales
de~$X$ avec la cat\'egorie des rev\^etements \'etales de~$X_0$.
\end{theoreme}

\section{Morphismes submersifs et universellement submersifs}
\label{IX.2}

\begin{definition}
\label{IX.2.1}
Un morphisme $g\colon S' \to S$ de pr\'esch\'emas est dit
\emph{submersif}
\index{submersif (morphisme)|hyperpage}%
s'il est surjectif, et fait de~$S$ un espace topologique quotient
de~$S'$ (\ie une partie $U$ de~$S$ telle que $f^{-1}(U)$ soit
ouverte, est ouverte). On dit que $f$ est \emph{universellement
submersif}
\index{universellement submersif (morphisme)|hyperpage}%
si pour tout morphisme $T \to S$, le morphisme $f'\colon T'=S'\times_S
T \to T$ d\'eduit de~$f$ par changement de base est submersif.
\end{definition}

Il est imm\'ediat que le compos\'e de deux morphismes submersifs
(\resp universellement submersifs) est submersif
(\resp universellement submersif), et qu'un changement de base dans un
morphisme universellement submersif donne un morphisme universellement
submersif (vu qu'on fait ce qu'il faut pour cela). Si $fg$ est
submersif (\resp universellement submersif), $f$ l'est.
\begin{exemples}
\label{IX.2.2}
a) Un morphisme surjectif qui est ouvert, ou ferm\'e, est submersif,
donc un morphisme surjectif universellement ferm\'e ou
universellement ouvert est universellement submersif. Par exemple
\emph{un morphisme propre surjectif est universellement
submersif}. D'autre part \emph{un morphisme fid\`element plat et
quasi-compact est universellement submersif} (VIII~\Ref{VIII.4.3}). Ce
seront les deux cas les plus importants pour nous.
\end{exemples}

On
\marginpar{232}
peut appliquer \`a un morphisme submersif ou universellement
submersif $g\colon S' \to S$ les raisonnements de~VIII~\Ref{VIII.4.3},
on trouve en particulier:
\begin{proposition}
\label{IX.2.3}
Supposons $g\colon S' \to S$ submersif. Alors le diagramme suivant
d'applications est exact:
$$
\xymatrix@C=.5cm{\Ouv (S)\ar[r] &
\Ouv (S')\ar@<2pt>[r]\ar@<-2pt>[r] & \Ouv(S''),}
$$
o\`u $S'' = S'\times_S S'$, et o\`u $\Ouv(X)$ d\'esigne
l'ensemble des parties ouvertes du pr\'esch\'ema~$X$.
\end{proposition}

\begin{proposition}
\label{IX.2.4}
Soient $g\colon S' \to S$ un morphisme universellement submersif,
$f\colon X \to Y$ un $S$-morphisme, et $f'\colon X' \to Y'$ le
$S'$-morphisme qui s'en d\'eduit par changement de base. Pour que
$f$ soit ouverte (\resp ferm\'ee), il suffit que $f'$ le soit. Pour
que $f$ soit universellement ouverte, \resp universellement
ferm\'ee, \resp s\'epar\'ee, il faut et il suffit que $f'$ le
soit. Si de plus $g$ est quasi-compact, et $f$ localement de type
fini, pour que $f$ soit propre, il faut et il suffit que $f'$ le soit.
\end{proposition}

Pour ce dernier point, on note que si $f'$ est propre donc
quasi-compact alors $f$ est quasi-compact (VIII~\Ref{VIII.3.3}) donc
de type fini puisqu'il est localement de type fini. D'autre part il
est s\'epar\'e et universellement ferm\'e d'apr\`es ce qui
pr\'ec\`ede, donc il est propre.
\begin{proposition}
\label{IX.2.5}
Soit $S'$ un pr\'esch\'ema de type fini sur le spectre $S$ d'un
anneau local noeth\'erien complet, supposons que la fibre du point
ferm\'e $s$ de~$S$ soit finie, donc les anneaux locaux dans $S'$ des
points $s'$ de cette fibre sont finis sur~$A=\cal{O}_s$. Soit~$S''$ le
sch\'ema somme des spectres des $\cal{O}_{S',s'}$ en question,
consid\'er\'e comme $S$-sch\'ema fini. Pour que $g\colon S' \to
S$ soit universellement submersif, il faut et il suffit que le
morphisme structural $S''\to S'$ soit surjectif.
\end{proposition}

Comme il y a un $S$-morphisme naturel $S'' \to S'$, et qu'un morphisme
fini surjectif est universellement submersif d'apr\`es~\Ref{IX.2.2},
la condition \'enonc\'ee est suffisante. Montrons donc que si
$S''\to S'$ n'est pas surjectif, alors $g$ n'est pas universellement
submersif. En effet, soit $t$ un point de~$S$ qui n'est pas dans
\marginpar{233}
l'image de~$S''$; il existe alors un $S$-sch\'ema $T$, spectre d'un
anneau de valuation discr\`ete, dont l'image dans $S$ est
$\{s,t\}$. Notons que l'image de~$S''$ dans $S'$ est ouverte, car le
morphisme $S'' \to S'$ est un isomorphisme local, et d'autre part
cette image contient $S'_s$, et ne rencontre pas $S'_t$. Il s'ensuit
que l'image inverse de cette derni\`ere dans $T'=S'\times_S T$ est
\emph{ouverte}, et identique \`a l'image inverse du point ferm\'e
de~$T$. Cela montre que $T'\to T$ n'est pas submersif, donc $S'\to S$
n'est pas universellement submersif.
\begin{remarque}
\label{IX.2.6}
Utilisant le crit\`ere IV~\Ref{IV.6.3} pour qu'une partie
constructible d'un espace noeth\'erien soit ouverte, on trouve
facilement le crit\`ere valuatif suivant pour qu'un morphisme
$g\colon S'\to S$ \emph{de type fini}, avec $S$ localement
noeth\'erien, soit universellement submersif: il faut et il suffit
que pour tout $S$-sch\'ema $T$, spectre d'un anneau de valuation
discr\`ete, posant $T'=S'\times_S T$, l'image inverse dans $T'$ du
point ferm\'e de~$T$ soit non ouverte.
\end{remarque}

\section{Descente de morphismes de pr\'esch\'emas \'etales}
\label{IX.3}

\begin{proposition}
\label{IX.3.1}
Soient $g\colon S'\to S$ un morphisme \emph{surjectif} de
pr\'esch\'emas, $X$ et $Y$ deux pr\'esch\'emas sur~$S$, $X'$,
$Y'$ leurs images inverses sur~$S'$. Si $Y$ est non ramifi\'e sur~$S$, alors l'application canonique
$$
\Hom_S(X,Y) \to \Hom_{S'}(X',Y')
$$
est injective.
\end{proposition}

En effet, en vertu de~\Ref{IX.1.6}, un $S$-morphisme $f\colon X\to Y$
est connu quand on conna\^it l'ensemble sous-jacent \`a son
graphe $\Gamma$, qui est une partie de~$Z=X\times_SY$. Comme
$$
Z'=Z\times_SS'=X'\times_{S'}Y' \to Z
$$
est surjectif, (puisque $S'\to S$ l'est), cette partie $\Gamma$ est
connue quand on conna\^it son image inverse dans $X'\times_{S'}Y'$,
qui n'est autre que l'ensemble sous-jacent au graphe de~$f'$. D'o\`u
la conclusion.

Une partie~$\Gamma$ de~$Z$ est le graphe d'un $S$-morphisme $f\colon
X\to Y$ si et seulement si elle est ouverte, et si le morphisme induit
par~$\pr_1$ de~$\Gamma$ dans~$Y$ est radiciel
\marginpar{234}
et surjectif \cf \Ref{IX.1}. Lorsque la premi\`ere
propri\'et\'e est v\'erifi\'ee, la deuxi\`eme l'est si et
seulement si l'image inverse~$\Gamma'$ de~$\Gamma$ dans~$Z'$ satisfait
la m\^eme condition VIII~\Ref{VIII.3.1}. Si on sait enfin que
$Z'\to Z$ est submersif, ce qui sera le cas en particulier si $S'\to
S$ est universellement submersif, alors~$\Gamma$ est ouvert si et
seulement si~$\Gamma'$ l'est. Ainsi, l'ensemble $\Hom_S(X,Y)$ est
alors en correspondance biunivoque avec l'ensemble des parties
ouvertes~$\Gamma'$ de~$Z'$ telles que le morphisme projection
$\pr_1\colon Z'\to X'$ soit radiciel et surjectif, (\ie correspondant \`a un $S'$-morphisme $f'\colon X'\to Y'$), et qui
sont satur\'ees pour la relation d'\'equivalence d\'efinie par
$Z'\to Z$, \ie dont les deux images inverses dans $Z''=Z'\times_Z
Z'=Z\times_S S''$ (o\`u $S''=S'\times_S S'$), par l'une et l'autre
projection, sont \'egales. Or ces derni\`eres sont les graphes
des deux $S''$-morphismes $X''\to Y''$ d\'eduits de~$f'$ par
changement de base, par l'une et l'autre projection $S''\to S'$. On a
ainsi obtenu:

\begin{proposition}
\label{IX.3.2}
Soient $g\colon S'\to S$ un morphisme \emph{universellement submersif}
de pr\'esch\'emas, $S''=S'\times_S S'$, $X$ et~$Y$ deux
$S$-pr\'esch\'emas, $X'$ et~$Y'$ leurs images inverses sur~$S'$,
et $X''$, $Y''$ leurs images inverses sur~$S''$. Si~$Y$ est \'etale
sur~$S$ le diagramme canonique suivant d'applications est exact:
$$
\xymatrix@C=.5cm{ {\Hom_S(X,Y)} \ar[r] & {\Hom_{S'}(X',Y')}
\ar@<2pt>[r]\ar@<-2pt>[r] & {\Hom_{S''}(X'',Y'')}. }
$$
\end{proposition}

Prenant~$X$ et~$Y$ \'etales sur~$S$, on trouve l'\'enonc\'e
suivant, qui d'ailleurs redonne~\Ref{IX.3.2} (m\^eme en se
restreignant \`a $X=S$, auquel cas on peut en effet toujours se
ramener dans~\Ref{IX.3.2}, par le changement de base $X\to S$);

\begin{corollaire}
\label{IX.3.3}
Un morphisme universellement submersif de pr\'esch\'emas est un
morphisme de descente pour la cat\'egorie fibr\'ee des
pr\'esch\'emas \'etales sur d'autres.
\end{corollaire}

J'ignore d'ailleurs si c'est n\'ecessairement un morphisme de
descente \emph{effective} pour la cat\'egorie fibr\'ee en
question, m\^eme en faisant de plus l'hypoth\`ese que~$S$ est
noeth\'erien et~$g$ de type fini, et en se bornant aux
rev\^etements \'etales. Nous donnerons n\'eanmoins au
num\'ero suivant des crit\`eres utiles d'effectivit\'e.

\begin{corollaire}
\label{IX.3.4}
Soit
\marginpar{235}
$g\colon S'\to S$ un morphisme universellement submersif, dont
les fibres $g^{-1}(s)$ sont \og g\'eom\'etriquement connexes\fg,
\index{geometriquement connexe@g\'eom\'etriquement connexe|hyperpage}%
\ie pour toute extension $K/\kres(s)$, $g^{-1}(s)\otimes_{\kres(s)} K$ est
connexe. Alors~$S'$ est connexe si~$S$ l'est. Le foncteur de la
cat\'egorie des pr\'esch\'emas \'etales sur~$S$ dans la
\ifthenelse{\boolean{orig}}
{cat\'egories}
{cat\'egorie}
des pr\'esch\'emas \'etales sur~$S'$ d\'efini par~$g$ est
pleinement fid\`ele.
\end{corollaire}

Une partie de~$S'$ qui est \`a la fois ouverte et ferm\'ee est
satur\'ee pour la relation d'\'equivalence ensembliste d\'efinie
par~$g$, puisque les fibres sont connexes, donc est l'image inverse
d'une partie de~$S$, qui est n\'ecessairement ouverte et ferm\'ee
puisque~$g$ est submersif. Si donc~$S$ est connexe, $S'$ l'est. Cela
implique aussi le r\'esultat suivant: le compos\'e~$fg$ de deux
morphismes \`a fibres universellement connexes,
\index{universellement connexe: synonyme de g\'eom\'etriquement connexe|hyperpage}%
$f$ \'etant universellement submersif, est \`a fibres
universellement connexes; si~$S_1'$ et~$S_2'$ sur~$S$ ont des fibres
universellement connexes, il en est de m\^eme de~$S_1'\times_S
S_2'$. En particulier, sous les conditions de~\Ref{IX.3.4}, $S''$ a
des fibres universellement connexes sur~$S$. Soient alors~$X$ et~$Y$
\'etales sur~$S$, et soit~$u'$ un $S'$-morphisme de~$X'$ dans~$Y'$,
prouvons qu'il est compatible avec les donn\'ees de descente (ce qui
entra\^ine la conclusion voulue gr\^ace \`a~\Ref{IX.3.3}). Or
soient~$u_1''$ et~$u_2''$ les deux $S''$-morphismes $X''\to Y''$
d\'eduits de~$u'$. Le sous-pr\'esch\'ema de~$S''$ des
co\"incidences de~$u_1''$ et~$u_2''$ est un sous-pr\'esch\'ema
ouvert induit, ferm\'e fibre par fibre, comme image inverse du
pr\'esch\'ema diagonal de~$Y''$ sur~$S''$\kern1pt\footnote{noter que les
fibres de~$S'$ sur~$S$ sont s\'epar\'ees!}. C'est donc l'image
inverse d'une partie de~$S$. Comme elle contient la diagonale
dans~$S''$, elle est identique \`a~$S''$, d'o\`u $u_1''=u_2''$
cqfd.

\section{Descente de pr\'esch\'emas \'etales: crit\`eres
d'effectivit\'e}
\label{IX.4}

\ifthenelse{\boolean{orig}}{}
{\enlargethispage{.7cm}}%
\begin{proposition}
\label{IX.4.1}
Soit $g\colon S'\to S$ un morphisme fid\`element plat et
quasi-compact. Alors~$g$ est un morphisme de descente effective pour
la cat\'egorie fibr\'ee des pr\'esch\'emas \'etales,
s\'epar\'es et de type fini sur d'autres.
\end{proposition}

C'est en effet un morphisme de descente pour la cat\'egorie
fibr\'ee en question, en vertu de~\Ref{IX.3.3} ou de
VIII~\Ref{VIII.5.2} au choix. Reste \`a montrer que si~$X'$ est
\'etale, s\'epar\'e et de type fini sur~$S'$, et muni d'une
donn\'ee de descente relativement \`a $g\colon S'\to S$, cette
derni\`ere est effective dans la cat\'egorie fibr\'ee en
question. Or on voit facilement que si~$X$ est un pr\'esch\'ema
sur~$S$, alors
\marginpar{236}
il est \'etale sur~$S$ si et seulement si il est \'etale sur~$S'$
(en vertu de la d\'efinition~\Ref{IX.1.1} et de \loccit \Ref{VIII.3.6}). Donc il est \'etale, s\'epar\'e et de
type fini sur~$S$ si et seulement si~$X'$ l'est sur~$S'$, \cf par
exemple~\Ref{IX.2.4}. Donc il suffit de s'assurer de
l'effectivit\'e de la donn\'ee de descente sur~$X$ pour la
cat\'egorie fibr\'ee des fl\`eches de~$\Sch$. Or~$X'$ est
quasi-affine sur~$S'$ en vertu de VIII~\Ref{VIII.6.2}
et~\Ref{VIII.6.6}. On peut alors conclure en utilisant
VIII~\Ref{VIII.7.9}. Le lecteur notera d'ailleurs que la
d\'emonstration demande moins si on se borne aux pr\'esch\'emas
\'etales et \emph{finis} sur d'autres, car on peut alors invoquer
directement VIII~\Ref{VIII.2.1}.

\begin{corollaire}
\label{IX.4.2}
Soient $g\colon S'\to S$ un morphisme universellement submersif, $X'$
un
\ifthenelse{\boolean{orig}}
{$S'$-pr\'esch\'emas}
{$S'$-pr\'esch\'ema}
\'etale s\'epar\'e et de type fini, muni d'une donn\'ee de
descente relativement \`a~$g$, $S_1\to S$ un morphisme
fid\`element plat et quasi-compact, $S_1'$ et~$X_1'$ d\'eduits
de~$S'$ et~$X'$ par le changement de base, de sorte que $S_1'\to S_1$
est universellement submersif, $X_1'$ est \'etale s\'epar\'e et
de type fini sur~$S_1'$, et muni d'une donn\'ee de descente
relativement \`a $g_1\colon S_1'\to S_1$. Pour que la donn\'ee de
descente sur~$X'$ soit effective, il faut et il suffit que la
donn\'ee de descente sur~$X_1'$ le soit.
\end{corollaire}

Cela r\'esulte de la th\'eorie de la descente dans les
cat\'egories \cite{IX.D}, compte tenu de~\Ref{IX.4.1} et~\Ref{IX.3.3}.

On prouve de fa\c con analogue:

\begin{corollaire}
\label{IX.4.3}
Soient $g\colon S'\to S$ un morphisme universellement submersif, $X'$
un~$S'$-pr\'esch\'ema \'etale muni d'une donn\'ee de descente
relativement \`a~$g$, $(S_i)$ un recouvrement de~$S$ par des
ouverts. Pour que la donn\'ee de descente soit effective, il faut
et il suffit que pour tout~$i$, la donn\'ee de descente
correspondante sur~$X_i'=X\times_S S_i$, relativement au morphisme
$g_i\colon S_i'=S'\times_S S_i\to S_i$, le soit.
\end{corollaire}

Ce dernier r\'esultat conduit \`a d\'egager un crit\`ere
d'effectivit\'e local:

\enlargethispage{.3cm}
\begin{proposition}
\label{IX.4.4}
Soit $g\colon S'\to S$ un morphisme de pr\'esentation finie
\textup{VIII}~\Ref{VIII.3.6}
\marginpar{237}
et universellement submersif, $X'$ un pr\'esch\'ema \'etale et
de pr\'esentation finie sur~$S'$, muni d'une donn\'ee de descente
relativement \`a~$g$, enfin~$a$ un point de~$S$. Pour qu'il existe
un voisinage ouvert~$U$ de~$a$, tel que la donn\'ee de descente
correspondante sur~$X_U'=X'\times_S U$ relativement au morphisme
$$g_U\colon S_U'=S'\times_S U\to S_U=U$$ soit effective, il faut et il
suffit que la donn\'ee de descente correspondante sur $X_a'=X'\times_S\Spec(\cal{O}_a)$, relativement au morphisme
$$g_a\colon S_a'=S'\times_S\Spec(\cal{O}_a)\to S_a=\Spec(\cal{O}_a),$$
soit effective.
\end{proposition}

La n\'ecessit\'e \'etant triviale, montrons la suffisance. On
dispose donc d'un pr\'esch\'ema \'etale de type fini~$X_a$
sur~$S_a$, et d'un isomorphisme
\begin{equation*}
\label{eq:IX.4.*}
\tag{$*$} {X_a'\isomto X_a\times_{S_a}S_a'}
\end{equation*}
compatible avec les donn\'ees de descente. Conform\'ement \`a
un sorite g\'en\'eral facile sur les pr\'esch\'emas
d\'efinis sur une limite inductive d'anneaux (ici les anneaux~$A_f$,
o\`u~$A$ est l'anneau d'un voisinage ouvert affine de~$a$, et
o\`u~$f$ parcourt les \'el\'ements de~$A$ qui ne sont pas dans
l'id\'eal premier correspondant \`a~$a$), on peut trouver un
voisinage ouvert~$U$ de~$a$, un pr\'esch\'ema \'etale de type
fini~$X_U$ sur~$U=S_U$, et un $S_a$-isomorphisme $X_a\isomto
X_U\times_{S_U}S_a$. De plus, prenant~$U$ assez petit, on peut alors
supposer que l'isomorphisme~\eqref{eq:IX.4.*} provient d'un
isomorphisme:
$$
X_U'\isomto X_U\times_{S_U} S_U';
$$
ce dernier pourrait ne pas \^etre compatible avec les donn\'ees de
descente, cependant \`a condition de r\'etr\'ecir~$U$, il sera
compatible avec les donn\'ees de descente. Cela ach\`eve la
d\'emonstration.

\begin{corollaire}
\label{IX.4.5}
Sous les conditions de~\Ref{IX.4.4}, pour que la donn\'ee de
descente sur~$X'$ soit effective, il faut et il suffit que pour tout
$a\in S$, la donn\'ee de descente correspondante sur~$X_a'$,
relativement au morphisme
$S_a'=S'\times_S\Spec(\cal{O}_a)\to\Spec(\cal{O}_a)$, le soit.
Lorsque~$S$ est localement noeth\'erien, et~$X'$ s\'epar\'e
sur~$S'$, on peut dans le crit\`ere pr\'ec\'edent remplacer
aussi~$\cal{O}_a$ par son compl\'et\'e.
\end{corollaire}

La
\marginpar{238}
premi\`ere assertion r\'esulte de~\Ref{IX.4.4} et~\Ref{IX.4.3}, la
deuxi\`eme est alors cons\'equence de~\Ref{IX.4.2}. Utilisant
encore~\Ref{IX.4.2} et le fait que pour tout anneau local
noeth\'erien~$A$, on peut trouver un anneau local noeth\'erien
complet~$B$, et un homomorphisme local $A\to B$, tel que~$B$ soit plat
sur~$A$ et que $B/\goth{m}B$ soit une extension donn\'ee du corps
r\'esiduel $k=A/\goth{m}$ de~$A$, on trouve:

\begin{corollaire}
\label{IX.4.6}
Sous les conditions de~\Ref{IX.4.4}, supposons de plus~$X'$
s\'epar\'e sur~$S'$, et~$S$ localement noeth\'erien. Pour que la
donn\'ee de descente sur~$X'$ soit effective, il faut et il suffit
que pour tout pr\'esch\'ema~$S_1$ sur~$S$, spectre d'un anneau
local complet \`a corps r\'esiduel alg\'ebriquement clos, la
donn\'ee de descente correspondante sur~$X_1'=X'\times_S S_1$,
relativement au morphisme $g_1\colon S_1'\to S_1$, soit effective.
\end{corollaire}

\begin{theoreme}
\label{IX.4.7}
Soit $g\colon S'\to S$ un morphisme fini et surjectif, et de
pr\'esentation finie (cette derni\`ere hypoth\`ese \'etant
cons\'equence des autres si~$S$ est localement
noeth\'erien)\footnote{On peut montrer qu'il suffit en fait que~$g$
soit un morphisme \emph{entier}, en se ramenant au cas du texte par un
proc\'ed\'e de passage \`a la limite dans le style EGA IV~8.}.
Alors~$g$ est un morphisme de descente effective pour la cat\'egorie
fibr\'ee des pr\'esch\'emas \'etales, s\'epar\'es, de type
fini sur d'autres.
\end{theoreme}

Il faut montrer que si~$X'$ est \'etale, s\'epar\'e, de type
fini sur~$S'$, et muni d'une donn\'ee de descente relativement
\`a~$g$, alors cette donn\'ee est effective.
Utilisant~\Ref{IX.4.3}, on se ram\`ene facilement au cas o\`u~$S$
est noeth\'erien. Gr\^ace \`a~\Ref{IX.4.5}, on peut donc
supposer que~$S$ est le spectre d'un anneau local noeth\'erien, a
fortiori que
$$
\dim S=n<+\infty.
$$
On raisonne alors par r\'ecurrence sur~$\dim S=n$, l'assertion
\'etant triviale pour $n<0$. Supposons donc $n\geq 0$ et le
th\'eor\`eme d\'emontr\'e pour les dimensions $n'<n$. En
vertu de~\Ref{IX.4.6} on est ramen\'e au cas o\`u~$S$ est le
spectre d'un anneau local complet, donc~$S'$ est une r\'eunion finie
de spectres d'anneaux locaux complets. On a donc
$$
X'=X_1'\sqcup X_2'
$$
o\`u~$X_1'$ est \emph{fini} sur~$S'$, et o\`u~$X_2'$ n'a aucun
point au-dessus d'un des points ferm\'es
\marginpar{239}
de~$S'$. Consid\'erons les morphismes
$$
\xymatrix@C=.5cm{ q_1,q_2\colon X'' \ar@<2pt>[r] \ar@<-2pt>[r] & X' }
$$
correspondants \`a la donn\'ee de descente, compatibles avec
$\xymatrix@C=.5cm{p_1,p_2\colon S'' \ar@<2pt>[r] \ar@<-2pt>[r] & S'}$. On
voit aussit\^ot que
$$
X''=q_i^{-1}(X_1')\sqcup q_i^{-1}(X_2')\quad i=1,2
$$
est la d\'ecomposition canonique analogue de~$X''$ sur~$S''$, ce qui
implique $q_1^{-1}(X_1')=\allowbreak q_2^{-1}(X_1')$ et par suite~$X_1'$ et~$X_2'$
sont munis de donn\'ees de descente induites. Or soit~$T$ l'ouvert
de~$S$ compl\'ementaire de son point ferm\'e, donc $T'=S'\times_S
T$ est la partie de~$S'$ compl\'ementaire de l'ensemble des points
ferm\'es, et~$X_2'$, qui se trouve tout entier au-dessus de~$T'$,
est muni d'une donn\'ee de descente relativement au morphisme $T'\to
T$ induit par~$g$. Comme ce dernier est fini surjectif, et que $\dim
T<\dim S=n$, cette donn\'ee de descente est effective par
l'hypoth\`ese de r\'ecurrence. On voit donc qu'il suffit de
prouver que la donn\'ee de descente sur~$X_1'$ est effective, donc
on peut maintenant supposer~$X'$ \'etale et \emph{fini} sur~$S'$.
(N.B. le raisonnement par r\'ecurrence est inutile si on se borne
aux rev\^etements \'etales dans l'\'enonc\'e~\Ref{IX.4.7}).
Soit alors~$S_0$ le spectre du corps r\'esiduel de~$A$, soit
$S_0'=S'\times_S S_0$ et d\'efinissons de m\^eme~$S_0''$, $S_0'''$
\`a partir des carr\'es et cubes fibr\'es~$S''$ et~$S'''$
de~$S'$ sur~$S$. En vertu de~\Ref{IX.1.8}, les morphismes $S_0\to S$,
$S_0'\to S'$, etc. induisent des \'equivalences pour les
cat\'egories des rev\^etements \'etales de~$S$ et~$S_0$ d'une
part, $S'$ et~$S_0'$ d'autre part, \ifthenelse{\boolean{orig}}{etc...}{etc.} D'apr\`es les sorites de
la th\'eorie de la descente dans les cat\'egories \cite{IX.D}, il s'ensuit
que pour que $g\colon S'\to S$ soit un morphisme de descente effective
pour la cat\'egorie fibr\'ee des rev\^etements \'etales, il
faut et il suffit qu'il en soit ainsi de~$g_0\colon S_0'\to S_0$.
Mais c'est bien le cas, comme cas particulier de~\Ref{IX.4.1}, par
exemple. Cela ach\`eve la d\'emonstration.

\begin{corollaire}
\label{IX.4.8}
La conclusion de~\Ref{IX.4.7} subsiste si on suppose seulement que
$S'\to S$ est universellement submersif, de type fini et quasi-fini,
pourvu qu'on suppose~$S$ localement noeth\'erien.
\end{corollaire}

En
\marginpar{240}
vertu de~\Ref{IX.4.6}, on peut en effet supposer que~$S$ est le
spectre d'un anneau local noeth\'erien complet. Alors en vertu
de~\Ref{IX.2.5}, il existe un morphisme fini et surjectif $S_1\to S$,
et un $S$-morphisme $S_1\to S'$. Comme $S_1\to S$ est un morphisme de
descente strict universel pour la cat\'egorie fibr\'ee
envisag\'ee, en vertu de~\Ref{IX.4.7}, et que $S'\to S$ est un
morphisme de descente universel pour ladite, \Ref{IX.4.8} r\'esulte
des sorites g\'en\'eraux~\cite{IX.D}.

\begin{corollaire}
\label{IX.4.9}
Soit $g\colon S'\to S$ un morphisme de type fini, surjectif et
universellement ouvert, avec~$S$ localement noeth\'erien. Alors~$g$
est un morphisme de descente effective pour la cat\'egorie
fibr\'ee des pr\'esch\'emas \'etales, s\'epar\'es et de
type fini sur d'autres.
\end{corollaire}

Proc\'edant comme dans~\Ref{IX.4.7}, on est ramen\'e au cas
o\`u~$S$ est le spectre d'un anneau local noeth\'erien et
complet~$A$. Soit~$A_1$ une alg\`ebre finie sur~$A$, de
spectre~$S_1$, telle que $S_1\to S$ soit fini et \emph{surjectif},
donc un morphisme de descente effective universel pour la
cat\'egorie fibr\'ee envisag\'ee, gr\^ace \`a~\Ref{IX.4.7}.
Il r\'esulte alors des th\'eor\`emes g\'en\'eraux \cite{IX.D}
que~$g$ est un morphisme de descente effective pour la cat\'egorie
fibr\'ee envisag\'ee, si et seulement si le morphisme
correspondant $g_1\colon S_1'=S'\times_S S_1\to S_1$ l'est. Comme ce
dernier satisfait aux m\^emes hypoth\`eses que~$g$, on est
ramen\'e \`a prouver~\Ref{IX.4.9} pour~$S_1$ au lieu de~$S$.
Prenant d'abord pour~$A_1$ le compos\'e direct des $A/\goth{p}_i$,
pour les id\'eaux premiers minimaux~$\goth{p}_i$ de~$A$, on est
ramen\'e au cas o\`u~$A$ est \emph{int\`egre}. On montre
alors\footnote{\Cf EGA IV~14.3.13 et~14.5.4.} qu'il existe un
sous-sch\'ema int\`egre~$S_1$ de~$S'$, quasi-fini sur~$S$ et
dominant~$S$, passant par un point de la fibre de~$S'$ en le point
ferm\'e~$y$ de~$S$ (gr\^ace au fait que~$S'$ est universellement
ouvert de type fini sur~$S$ local noeth\'erien int\`egre, et
$S_y'\neq\emptyset$). Comme~$A$ est complet, $S_1$ est fini sur~$S$,
et comme il domine~$S$, le morphisme $S_1\to S$ est surjectif.
Rempla\c cant encore une fois~$S$ par~$S_1$, on est ramen\'e au
cas o\`u~$S'$ a une section sur~$S$, o\`u l'\'enonc\'e est
trivial.

\begin{theoreme}
\label{IX.4.10}
Soit $g\colon S'\to S$ un morphisme fini radiciel surjectif, de
pr\'esentation finie (cette derni\`ere condition \'etant
superflue si~$S$ est localement noeth\'erien\footnote{Il suffit
m\^eme que~$g$ soit entier, radiciel surjectif, comme on voit par
une r\'eduction facile au cas du texte, style EGA~IV~8, \cf SGA~4
VIII~1.1.}).
\marginpar{241}
Alors le foncteur image inverse induit une \'equivalence de la
cat\'egorie des pr\'esch\'emas \'etales sur~$S$ avec la
cat\'egorie des pr\'esch\'emas \'etales sur~$S'$.
\end{theoreme}

Comme les morphismes diagonaux de~$S'$ dans $S'\times_S S'$ et
$S'\times_S S'\times_S S'$ sont des immersions surjectives, donc
induisent en vertu de~\Ref{IX.1.9} des \'equivalences des
cat\'egories des pr\'esch\'emas \'etales sur~$S'\times_S S'$
\resp $S'\times_S S'\times_S S'$ avec la cat\'egorie des
pr\'esch\'emas \'etales sur~$S'$, il r\'esulte des sorites de
la descente \cite{IX.D} que tout~$X'$ \'etale sur~$S'$ est muni d'une
donn\'ee de descente et d'une seule relativement \`a $g\colon
S'\to S$. Donc~\Ref{IX.3.3} implique que le foncteur image inverse
par~$g$, de la cat\'egorie des pr\'esch\'emas \'etales sur~$S$
dans la cat\'egorie des pr\'esch\'emas \'etales sur~$S'$, est
\emph{pleinement fid\`ele}. Reste \`a montrer qu'il est
essentiellement surjectif, \ie que tout~$X'$ \'etale sur~$S'$ est
isomorphe \`a l'image inverse d'un~$X$ \'etale sur~$S$. La
question \'etant \'evidemment locale sur~$S$ \emph{et sur}~$X'$,
on peut supposer $S$,~$S'$,~$X'$ affines. Mais alors~$X'$ est
s\'epar\'e de type fini sur~$S'$, et on peut appliquer le
crit\`ere d'effectivit\'e~\Ref{IX.4.7}.

\begin{corollaire}
\label{IX.4.11}
La conclusion~\Ref{IX.4.9} subsiste en rempla\c cant l'hypoth\`ese
sur~$g$ par: $g$~est fid\`element plat, quasi-compact et radiciel.
\end{corollaire}

M\^eme d\'emonstration, en invoquant~\Ref{IX.4.1} au lieu
de~\Ref{IX.4.7}.

On notera que la d\'emonstration de~\Ref{IX.4.7} est
\og \'el\'ementaire\fg en ce qu'elle n'utilise pas les
th\'eor\`emes
\ifthenelse{\boolean{orig}}
{}
{de}
finitude et de comparaison pour les morphismes propres (EGA~III~3,~4,~5). Il n'en est plus de m\^eme du r\'esultat suivant:

\begin{theoreme}
\label{IX.4.12}
Soit $g\colon S'\to S$ un morphisme propre, surjectif, de
pr\'esentation finie (cette derni\`ere hypoth\`ese \'etant
cons\'equence de la premi\`ere si~$S$ est localement
noeth\'erien). Alors~$g$ est un morphisme de descente effective
pour la cat\'egorie fibr\'ee des rev\^etements \'etales de
pr\'esch\'emas.
\end{theoreme}

En vertu de~\Ref{IX.3.3} et de~\Ref{IX.2.2}, on est ramen\'e \`a
prouver que pour tout rev\^etement \'etale~$X'$ sur~$S'$, muni
d'une donn\'ee de descente relativement \`a $g\colon S'\to S$,
cette donn\'ee de descente est effective. Utilisant~\Ref{IX.4.3},
on est ramen\'e facilement
\marginpar{242}
au cas o\`u~$S$ est noeth\'erien, et utilisant~\Ref{IX.4.6}, on
peut donc supposer que~$S$ est le spectre d'un anneau local
noeth\'erien \emph{complet}~$A$. Introduisons~$S''$ et~$S'''$ comme
d'habitude, soit~$S_0$ le spectre du corps r\'esiduel de~$A$, et
soient~$S_0'$,~$S_0''$,~$S_0'''$ d\'eduits de~$S'$,~$S''$,~$S'''$
par le changement de base $S_0\to S$, \ie les fibres
de~$S'$,~$S''$,~$S'''$ au point ferm\'e de~$S$.
D'apr\`es~\Ref{IX.1.10}, les morphismes $S_0\to S$, $S_0'\to S'$,
etc. induisent des \'equivalences de la cat\'egorie des
rev\^etements \'etales sur le sch\'ema-but avec la cat\'egorie
des rev\^etements \'etales sur le sch\'ema-source. Par suite,
$g\colon S'\to S$ est un morphisme de descente strict pour la
cat\'egorie fibr\'ee des rev\^etements \'etales de
pr\'esch\'emas, si et seulement si $g_0\colon S_0'\to S_0$ l'est,
ce qui est bien le cas en vertu de~\Ref{IX.4.1}. Cela ach\`eve la
d\'emonstration de~\Ref{IX.4.12}. (Dans ce raisonnement, on n'avait
besoin de~\Ref{IX.1.10} que le fait que le foncteur envisag\'e
dans~\Ref{IX.1.10} est \emph{pleinement fid\`ele}, ce qui n'utilise
\emph{pas} le th\'eor\`eme d'existence de faisceaux coh\'erents
en g\'eom\'etrie alg\'ebrique.)

\section{Traduction en termes du groupe fondamental}
\label{IX.5}

Soit
$$
g\colon S'\to S
$$
un \emph{morphisme de descente effective} pour la cat\'egorie
fibr\'ee des \emph{rev\^etements \'etales} de
pr\'e\-sch\'e\-mas, par exemple un morphisme propre, surjectif, de
pr\'esentation finie \eqref{IX.4.12}, ou un morphisme fid\`element
plat et quasi-compact. Introduisant comme d'habitude~$S''$,~$S'''$,
et d\'esignant par~$\cal{C}$,~$\cal{C}'$,~$\cal{C}''$,~$\cal{C}'''$
la cat\'egorie des rev\^etements \'etales
de~$S$,~$S'$,~$S''$,~$S'''$ respectivement, on a donc un diagramme
$2$-exact de cat\'egories
\begin{equation*}
\label{eq:IX.5.*}
\tag{$*$} 
\xymatrix@C=.5cm{
{\cal{C}} \ar[r]^-{~p^*} & {\cal{C}'}
\ar@<2pt>[rr]^{p_1^*,p_2^*} \ar@<-2pt>[rr] && {\cal{C}''}
\ar@<4pt>[rrr]^{p_{21}^*,p_{32}^*,p_{31}^*} \ar[rrr] \ar@<-4pt>[rrr] &&&
{\cal{C}'''}
}
\end{equation*}
correspondant au diagramme
$$
\xymatrix@C=.5cm{ S & \ar[l]_-{~p} S' && \ar@<2pt>[ll] \ar@<-2pt>[ll]_-{p_1,p_2}
S'' &&& \ar@<4pt>[lll] \ar[lll] \ar@<-4pt>[lll]_{p_{21},p_{32},p_{31}} S'''.
}
$$
Supposons
\marginpar{243}
les pr\'esch\'emas~$S$,~$S'$,~$S''$,~$S'''$ sommes disjointes de
pr\'esch\'emas connexes, ce qui sera le cas en particulier si ce
sont des pr\'esch\'emas localement connexes, a fortiori s'ils sont
localement noeth\'eriens (par exemple si~$S'$ est de type fini
sur~$S$ localement noeth\'erien). Alors les
cat\'egories~$\cal{C}$,~$\cal{C}'$ ... dans~\eqref{eq:IX.5.*} sont
des cat\'egories multigaloisiennes V~\Ref{V.9}, d\'ecrites
chacune par une collection de groupes topologiques compacts totalement
disconnexes, savoir les groupes fondamentaux des composantes connexes
des pr\'esch\'emas~$S$,~$S'$,~$S''$,~$S'''$. Nous supposons pour
simplifier~$S$ connexe, et allons donner alors un proc\'ed\'e de
calcul pour son groupe fondamental, en termes de la cat\'egorie
fibr\'ee form\'ee avec~$\cal{C}'$,~$\cal{C}''$,~$\cal{C}'''$,
convenablement explicit\'ee \`a l'aide des groupes fondamentaux
exprimant ces cat\'egories. Le lecteur notera que le
proc\'ed\'e esquiss\'e est valable en fait dans le cadre
g\'en\'eral des cat\'egories multigaloisiennes (qui n'ont pas
\`a provenir de pr\'esch\'emas
donn\'es~$S$,~$S'$,~$S''$,~$S'''$). C'est d'ailleurs l'analogue du
proc\'ed\'e bien connu pour calculer le groupe fondamental d'un
espace topologique~$S$, r\'eunion localement finie de sous-espaces
ferm\'es~$S_i$ (ou r\'eunion quelconque de sous-espaces
ouverts~$S_i$), \`a l'aide des groupes fondamentaux des composantes
\ifthenelse{\boolean{orig}}
{des~$X_i$}
{des~$S_i$}
et des composantes des $S_i\cap S_j$. Bien entendu, la
situation analogue dans le cadre des pr\'esch\'emas tombe bien
dans le cadre g\'en\'eral de la descente, en introduisant le
pr\'esch\'ema~$S'$ somme des~$S_i$ et le morphisme canonique
$g\colon S'\to S$.

Posons
$$
E'=\pi_0(S'), \quad E''=\pi_0(S''), \quad E'''=\pi_0(S'''),
$$
o\`u~$\pi_0$ d\'esigne le foncteur \og ensemble des composantes
connexes\fg. Comme les produits fibr\'es de~$S'$ sur~$S$ forment un
objet simplicial de~$\Sch$, il est transform\'e par le
foncteur~$\pi_0$ en un ensemble simplicial dont~$E'$,~$E''$,~$E'''$
sont les composantes de dimension~$0$,~$1$,~$2$. Nous aurons \`a
utiliser les applications simpliciales
$$
q_i=\pi_0(p_i),\quad(i=1,2)\quad\text{et}
\quad q_{ij}=\pi_0(p_{ij}),\quad (i,j)=(2,1),(3,2),(3,1),
$$
mises en \'evidence dans le diagramme
\begin{equation*}
\label{eq:IX.5.1}
\tag{1} {\xymatrix@C=.5cm{ E' && \ar@<2pt>[ll] \ar@<-2pt>[ll]_{q_1,q_2} E'' &&&
\ar@<4pt>[lll] \ar[lll] \ar@<-4pt>[lll]_{q_{21},q_{32},q_{31}} E'''. }}
\end{equation*}
Les objets
\marginpar{244}
de~$E'$ seront not\'es avec un accent, comme~$s'$, ceux
de~$E''$ \resp $E'''$ seront not\'es avec un~$''$ \resp un~$'''$.
Le fait que~$S$ soit connexe se traduit par $\pi_0(K)=0$, o\`u~$K$
est l'ensemble simplicial d\'efini par $g\colon S'\to S$, ou encore
par le fait que la relation d'\'equivalence dans~$E'$ engendr\'ee
par le couple d'applications $(q_1,q_2)$ est transitive.

Nous choisirons une fois pour toutes un \'el\'ement~$s_0'$
dans~$E'$, et pour tout~$s'$ dans~$E'$, un
\'el\'ement~$\overline{s'}\in E''$ tel que\footnote{On fera
attention que l'\'el\'ement~$\overline{s'}$ dont l'existence est
admise implicitement, n'existe pas dans tous les cas. Par suite, le
th\'eor\`eme~\Ref{IX.5.1} tel qu'il est \'enonc\'e ne
s'applique pas dans tous les cas. Il n'est pas difficile cependant,
en s'inspirant du texte \'ecrit, de modifier l'\'enonc\'e de ce
th\'eor\`eme de telle fa\c con qu'il donne une m\'ethode de
calcul qui s'applique dans tous les cas. En particulier, les
corollaires dudit th\'eor\`eme sont valables tels quels.}
$$
q_1(\overline{s'})=s_0', \quad q_2(\overline{s'})=s',
$$
mettant ainsi en \'evidence la connexit\'e de~$S$. Pour
tout~$s'\in E'$, choisissons un point
g\'eom\'etrique~$\underline{s}'$ dans la composante connexe~$s'$
de~$S'$; ce point interviendra en fait par le foncteur-fibre~$F_{s'}'$
correspondant sur la cat\'egorie multigaloisienne~$\cal{C}'$. Le
groupe des automorphismes de ce foncteur, \ie le groupe fondamental
de~$S'$ en~$\underline{s}'$, sera not\'e~$\pi_{s'}$. On choisit de
m\^eme des~$\underline{s}''$ et des~$\underline{s}'''$, donc des
foncteurs~$F_{s''}''$ et~$F_{s'''}'''$, d'o\`u des groupes
fondamentaux~$\pi_{s''}$ et~$\pi_{s'''}$. Ainsi
$$
\pi_{s'}=\pi_1(S',\underline{s}'), \quad
\pi_{s''}=\pi_1(S'',\underline{s}''), \quad
\pi_{s'''}=\pi_1(S''',\underline{s}''').
$$
Pour tout $s''\in E''$, $p_1(\underline{s}'')$ se trouve dans la
m\^eme composante connexe que $\underline{q_1(s'')}$, donc il existe
un isomorphisme de foncteurs
\ifthenelse{\boolean{orig}}
{$F''_{s''}\circ p_1^* F'_{s'}$}
{$F''_{s''}\circ p_1^*\isomto F'_{s'}$}
(\ie une \og classe de chemins\fg de~$p_1(\underline{s}'')$ \`a
$\underline{q_1(s'')}$). Cette remarque se r\'ep\`ete pour~$q_2$,
et les~$q_{ij}$. Choisissons toutes ces classes de chemins:
$$
F''_{s''}\circ p_i^*\isomto F'_{q_i(s'')},\quad F_{s'''}'''\circ
p_{ij}^*\isomto F'''_{q_{ij}(s''')},
$$
(pour
\marginpar{245}
$i=1,2$ et $(i,j)=(2,1),(3,2),(3,1)$). Il en r\'esulte en
particulier des homomorphismes de groupes:
\begin{equation*}
\label{eq:IX.5.2}
\tag{2} {q_i^{s''}\colon \pi_{s''}\to\pi_{q_i(s'')},\quad
q_{ij}^{s'''}\colon \pi_{s'''}\to\pi_{q_{ij}(s''')},}
\end{equation*}
(m\^emes valeurs de~$i$ et $(i,j)$). Enfin, rappelons-nous que dans
la structure de cat\'egorie cliv\'ee de
fibres~$\cal{C}'$,~$\cal{C''}$,~$\cal{C'''}$ figurent aussi des
isomorphismes de foncteurs:
$$
p_{21}^*p_1^*\isomto p_{31}^*p_1^*,\quad p_{21}^*p_2^*\isomto
p_{32}^*p_1^*,\quad p_{31}^*p_2^*\isomto p_{32}^*p_2^*,
$$
d\'eduits d'isomorphismes des deux membres respectivement avec les
$u_i^*$~($i=1,2,3$), o\`u les~$u_i$ sont les trois projections
de~$S'''$ dans~$S'$. Quand on explicite ces donn\'ees, on trouve
pour tout~$s'''$ un \'el\'ement bien d\'etermin\'e
\begin{equation*}
\label{eq:IX.5.3}
\tag{3} {a_i^{s'''}\in\pi_{v_i(s''')},}
\end{equation*}
(o\`u les~$v_i$,~$i=1,2,3$ sont les trois applications $E'''\to E'$
d\'efinies par $v_i=\pi_0(u_i)$), d'ailleurs soumis aux conditions:
$$
q_1^{s_1''}q_{21}^{s'''}=\mathrm{int}(a_1^{s'''})
q_1^{s_2''}q_{31}^{s'''}\qquad (s_1''=q_{21}(s'''),\
s_2''=q_{31}(s''')),
$$
et les deux conditions analogues, faisant intervenir les~$a_2$
et~$a_3$. Le lecteur notera d'ailleurs que les
donn\'ees~\eqref{eq:IX.5.1},~\eqref{eq:IX.5.2},~\eqref{eq:IX.5.3}
permettent de reconstituer, \`a une \'equivalence de
cat\'egories fibr\'ees pr\`es, la cat\'egorie fibr\'ee
envisag\'ee de fibres~$\cal{C}'$,~$\cal{C}''$,~$\cal{C}'''$. Elles
doivent donc permettre en principe de reconstituer~$\cal{C}$ \`a
\'equivalence pr\`es, donc son groupe fondamental \`a
isomorphisme pr\`es. Nous d\'eterminerons en fait le groupe
fondamental en le point g\'eom\'etrique $p(\underline{s}_0')$
de~$S$, \ie le groupe des automorphismes de~$F'_{s_0'}\circ p^*$.

On note que la donn\'ee d'un objet~$X'$ de~$\cal{C}'$ est
\'equivalente essentiellement \`a la donn\'ee d'ensembles
finis~$X'_{s'}$ ($s'\in E'$) o\`u les~$\pi_{s'}$ op\`erent
contin\^ument. Une
\marginpar{246}
donn\'ee de recollement sur un tel objet revient alors \`a la
donn\'ee, pour tout $s''\in E''$, d'une bijection:
$$
\varphi_{s''}\colon X'_{q_1(s'')}\isomto X'_{q_2(s'')}
$$
compatible avec les op\'erations de~$\pi_{s''}$, op\'erant sur
l'un et l'autre membre gr\^ace aux homomorphismes $q_i^{s''}\colon
\pi_{s''}\to\pi_{q_i(s'')}$. Prenant d'abord les~$s''$ de la
forme~$\overline{s'}$, on voit qu'une telle donn\'ee d\'efinit des
bijections
$$
\psi_{s'}\colon X'_{s_0'}=F_0'(X')\to X_{s'}'
$$
ce qui permet d'identifier les~$X'_{s'}$ au m\^eme ensemble
$F_0'(X')=X'_{s_0'}$, sur lesquels tous les groupes~$\pi_{s'}$ vont
d\`es lors op\'erer. Cela pos\'e, les
bijections~$\varphi_{s''}$ vont correspondre \`a des bijections
$$
g_{s''}\colon F_0'(X')\isomto F_0'(X'),
$$
soumis d'une part aux relations de commutation avec~$\pi_{s''}$:
\begin{equation*}
\label{eq:IX.5.a}
\tag*{a)} {g_{s''}q_1^{s''}(g'')=q_2^{s''}(g'')g_{s''}\qquad (s''\in E'',\ g''\in\pi_{s''}),}
\end{equation*}
d'autre part aux relations
\begin{equation*}
\label{eq:IX.5.b}
\tag*{b)} {g_{\overline{s'}}=g_{\overline{s_0'}}\qquad (s'\in E'),}
\end{equation*}
exprimant la fa\c con dont nous avions identifi\'e entre eux
les~$X_{s'}$. Quand on explicite la condition pour qu'une telle
donn\'ee de recollement soit en fait une donn\'ee de descente, on
trouve les relations:
\begin{equation*}
\label{eq:IX.5.c}
\tag*{c)} {a_3^{s'''}g_{q_{31}(s''')}a_1^{s'''}=
g_{q_{32}(s''')}a_2^{s'''}g_{q_{21}(s''')}\qquad (s'''\in E''').}
\end{equation*}

Cela nous donne une \'equivalence entre la cat\'egorie des objets
de~$\cal{C}'$ munis d'une
\marginpar{247}
donn\'ee de descente, et la cat\'egorie des ensembles finis o\`u
les groupes~$\pi_{s'}$ op\`erent contin\^ument, munis de plus de
bijections~$g_{s''}$, satisfaisant les
relations~\Ref{eq:IX.5.a},~\Ref{eq:IX.5.b},~\Ref{eq:IX.5.c}. Soit
alors~$G$ le groupe engendr\'e par les groupes~$\pi_{s'}$ et les
nouveaux g\'en\'erateurs~$g_{s''}$, soumis aux
relations~\Ref{eq:IX.5.a},~\Ref{eq:IX.5.b},~\Ref{eq:IX.5.c}, et
soit~$\pi$ le groupe limite projective des quotients de~$G$ par les
sous-groupes d'indice fini dont les images inverses dans les
groupes~$\pi_{s'}$ soient des sous-groupes ouverts. On dit aussi
que~$\pi$ est le \emph{groupe de type galoisien engendr\'e par
les~$\pi_{s'}$ et les~$g_{s''}$, soumis aux
relations~\Ref{eq:IX.5.a},~\Ref{eq:IX.5.b},~\Ref{eq:IX.5.c}}. On
constate aussit\^ot que la cat\'egorie envisag\'ee est aussi
\'equivalente \`a la cat\'egorie des ensembles finis o\`u le
groupe topologique~$\pi$ op\`ere contin\^ument. Cela \'etablit
l'\'enonc\'e suivant:
\begin{theoreme}
\label{IX.5.1}
Soit $g\colon S'\to S$ un morphisme de pr\'esch\'emas qui soit un
morphisme de descente effective pour la cat\'egorie fibr\'ee de
rev\^etements \'etales de pr\'esch\'emas (\cf \Ref{IX.4.9}
et~\Ref{IX.4.12}). Supposons~$S$ connexe, et~$S'$, son carr\'e
fibr\'e~$S''$ et son cube fibr\'e~$S'''$, sommes de
pr\'esch\'emas connexes (ce qui est le cas par exemple si~$S'$ est
de type fini sur~$S$ localement noeth\'erien et connexe).
Choisissons comme dessus: un point g\'eom\'etrique dans toute
composante connexe de~$S'$,~$S''$,~$S'''$, certaines classes de
chemins, un $s_0'\in E'$, et pour tout $s'\in E'$ un
\ifthenelse{\boolean{orig}}
{$s'\in E''$}
{$s''\in E''$}
dont les deux images dans~$E'$ soient~$s_0'$ et~$s'$.
($E'$,~$E''$,~$E'''$ d\'esignent respectivement l'ensemble des
composantes connexes de~$S'$,~$S''$,~$S'''$). Alors le groupe
fondamental de~$S$ en le point g\'eom\'etrique image de~$s_0'$ est
canoniquement isomorphe au groupe de type galoisien engendr\'e par
les $\pi_{s'}=\pi_1(S',\underline{s}')$ ($s'\in E'$) et des
g\'en\'erateurs~$g_{s''}$ ($s''\in E''$), soumis aux
relations~\Ref{eq:IX.5.a},~\Ref{eq:IX.5.b},~\Ref{eq:IX.5.c} ci-dessus
faisant intervenir les \'el\'ements des groupes
$\pi_{s''}=\pi_1(S'',\underline{s}'')$, et les
\'el\'ements~$a_i^{s'''}$ ($i=1,2,3$, $s'''\in E'''$) introduits
plus haut.
\end{theoreme}

\begin{corollaire}
\label{IX.5.2}
Supposons que~$S'$ et~$S''$ n'aient qu'un nombre fini de composantes
connexes, et que les groupes fondamentaux des composantes connexes
de~$S'$ soient topologiquement de g\'en\'eration finie. Alors le
groupe fondamental de~$S$ est topologiquement de g\'en\'eration
finie.
\end{corollaire}

Ainsi,
\marginpar{248}
nous prouverons plus tard que le groupe fondamental d'un sch\'ema
projectif normal sur un corps alg\'ebriquement clos est
topologiquement de g\'en\'eration finie. Utilisant le lemme de
Chow et la normalisation des sch\'emas alg\'ebriques, il
s'ensuivra que le m\^eme r\'esultat est vrai pour tout sch\'ema
propre sur un corps alg\'ebriquement clos.

\begin{corollaire}
\label{IX.5.3}
Supposons que~$S'$,~$S''$,~$S'''$ n'aient qu'un nombre fini de
composantes connexes, que les groupes fondamentaux des composantes
connexes de~$S'$ soient topologiquement de \emph{pr\'esentation
finie}, et les groupes fondamentaux des composantes connexes de~$S''$
topologiquement de \emph{g\'en\'eration finie}. Alors le groupe
fondamental de~$S$ est topologiquement de \emph{pr\'esentation
finie}.
\end{corollaire}

On notera qu'on peut exprimer~\Ref{IX.4.9} (restreint aux
\emph{rev\^etements} \'etales) en disant qu'\emph{un morphisme
fini radiciel surjectif de pr\'esch\'emas noeth\'eriens induit
un isomorphisme des groupes fondamentaux}; de fa\c con imag\'ee,
on peut donc dire que le groupe fondamental est un \emph{invariant
topologique} pour les pr\'esch\'emas. On peut expliciter plus
g\'en\'eralement, \`a l'aide de~\Ref{IX.5.1}, l'effet sur le
groupe fondamental d'op\'erations sur les pr\'esch\'emas, telles
que le \og pincement\fg du pr\'esch\'ema suivant un ensemble fini de
points, ayant une signification topologique simple. On trouve par
exemple:

\begin{corollaire}
\label{IX.5.4}
Soient $g\colon S'\to S$ un morphisme fini de pr\'esentation finie,
$T$ une partie discr\`ete de~$S$. Pour tout $s\in S$, soit $n(s)$ le
\og nombre g\'eom\'etrique de points\fg dans la fibre $g^{-1}(s)$
(qui s'explicite aussi comme le degr\'e s\'eparable de~$g^{-1}(s)$
sur~$\kres(s)$, somme des degr\'es s\'eparables de ses extensions
r\'esiduelles). On suppose que pour $s\in S-T$, on a $n(s)=1$. Pour
tout $s\in T$, soit $K_s$ une extension alg\'ebriquement close de
$\kres(s)$, $I_s$ l'ensemble des points g\'eom\'etriques de~$S'$ \`a
valeurs dans~$K_s$ (c'est un ensemble \`a $n(s)$ \'el\'ements),
$I_s'$ le compl\'ementaire d'un point choisi de~$I_s$, et enfin $I'$
l'ensemble r\'eunion des $I_s'$. On suppose $S'$ connexe. Alors le
groupe fondamental de~$S$ est isomorphe au groupe de type galoisien
engendr\'e par le groupe fondamental de~$S'$, et des
g\'en\'erateurs $g_i$ ($i\in I'$), soumis \`a aucune condition
suppl\'ementaire.
\end{corollaire}

Le
\marginpar{249}
d\'etail de la d\'emonstration est laiss\'e au lecteur;
l'\'enonc\'e obtenu n'est que la traduction, en langage de la
th\'eorie des groupes, du fait qu'on a une \'equivalence de la
cat\'egorie $C$ des rev\^etements \'etales de~$S$, et de la
cat\'egorie des rev\^etements \'etales $X'$ de~$S'$,
\emph{munis} pour tout $s\in T$ d'un syst\`eme transitif de
bijections entre les $n(s)$ fibres de~$X'$ aux points de~$g^{-1}(s)$
\`a valeurs dans~$K_s$. (Sous cette forme intrins\`eque bien
entendu, il n'est plus n\'ecessaire de supposer $S'$ connexe).
\begin{exemple}
\label{IX.5.5}
On prouve facilement que la courbe rationnelle $\PP^1_k$ sur un corps
alg\'ebriquement clos $k$ est simplement
connexe\footnote{\Cf Exp\ptbl XI~\Ref{XI.1.1}.}. Donc le groupe
fondamental d'une courbe rationnelle compl\`ete ayant exactement un
point double, \`a $n$ branches analytiques, est le groupe de type
galoisien libre engendr\'e par $n-1$ g\'en\'erateurs. Par
exemple, dans le cas d'un point double ordinaire, on trouve le groupe
fondamental $\hat\ZZ$, comme annonc\'e dans (I~\Ref{I.11}~a)). Par
contre, l'existence d'un point de rebroussement (qui est un point
\og g\'eom\'etriquement unibranche\fg) n'a pas d'influence sur le
groupe fondamental.
\end{exemple}

\begin{corollaire}
\label{IX.5.6}
Soit $g\colon S'\to S$ un morphisme de pr\'esch\'emas
universellement submersif, \`a fibres g\'eom\'etriquement
connexes, $S$ \'etant connexe. Alors $S'$ est connexe, et
choisissant un point g\'eom\'etrique $s'$ dans~$S'$ et
d\'esignant par $s$ son image dans~$S$, l'homomorphisme
$$
\pi_1(S',s')\to \pi_1(S,s)
$$
est \emph{surjectif}. Si $g$ est un morphisme de descente effective
pour la cat\'egorie fibr\'ee des rev\^etements \'etales de
pr\'esch\'emas (\cf \emph{\Ref{IX.4.12}}), introduisant le point
g\'eom\'etrique $s^{\prime\prime}=\diag(s')$
de~$S^{\prime\prime}=S'\times_S S'$, et les deux homomorphismes
$$
p_{1*},p_{2*}\colon \pi_1(S^{\prime\prime},s^{\prime\prime}) \to
\pi_1(S',s')
$$
induits par les deux projections, $\pi_1(S, s)$ est isomorphe au
conoyau de ce couple de morphismes dans la cat\'egorie des groupes
de type galoisien, \ie au quotient de~$\pi_1(S',s')$ par le
sous-groupe invariant ferm\'e engendr\'e par les
\marginpar{250}
\'el\'ements de la forme
$p_{1*}(g^{\prime\prime})p_{2*}(g^{\prime\prime})^{-1}$, avec
$g^{\prime\prime}\in \pi_1(S^{\prime\prime},s^{\prime\prime})$.
\end{corollaire}

On sait en effet par~\Ref{IX.3.4} que le foncteur $X\mto
X\times_SS'$ des rev\^etements \'etales sur~$S$ dans les
rev\^etements \'etales sur~$S'$ est pleinement fid\`ele, ce qui
\'equivaut au fait que l'homomorphisme sur les groupes fondamentaux
est un \'epimorphisme (V~\Ref{V.6.9}). La derni\`ere assertion est
une cons\'equence imm\'ediate de la description~\Ref{IX.5.1}.

\begin{remarque}
\label{IX.5.7}
Il n'est pas connu \`a l'heure actuelle si le groupe fondamental
d'un sch\'ema propre sur un corps alg\'ebriquement clos $k$ est
topologiquement de pr\'esentation finie\footnote{Cela semble
tr\`es improbable dans le cas des courbes lisses de genre $g\geq 2$,
en caract\'eristique $p>0$. Quand on remplace $\pi_1$ par son plus
grand quotient premier \`a $p$, par contre, il semble que les
techniques bien connues permettent de donner une r\'eponse
affirmative, m\^eme sans hypoth\`ese de propret\'e. \Cf un
travail en pr\'eparation de J.P\ptbl Murre.}. Utilisant~\Ref{IX.5.3}, une
technique bien connue de sections hyperplanes, et la
d\'esingularisation des surfaces normales, on est ramen\'e au cas
d'une \emph{surface lisse sur} $k$. Cela permet du moins de montrer,
par voie transcendante, que la r\'eponse est affirmative en
caract\'eristique $0$ (et ceci sans \^etre oblig\'e d'admettre
la triangulabilit\'e de vari\'et\'es alg\'ebriques
singuli\`eres). En caract\'eristique $p>0$, la difficult\'e
principale semble dans le cas des courbes, dont on sait seulement que
le groupe fondamental est un quotient de celui qui se pr\'esente
dans le cas classique (\cf expos\'e suivant), le noyau par lequel on
divise \'etant cependant fort mal connu.
\end{remarque}

\begin{remarque}
\label{IX.5.8}
On pourrait expliciter d'autres cas particuliers que~\Ref{IX.5.4}
et~\Ref{IX.5.6} o\`u~\Ref{IX.5.1} prend une forme
particuli\`erement simple. Un cas int\'eressant est celui o\`u
$S$ est le quotient de~$S'$ par un groupe fini d'automorphismes
$\Gamma$. \emph{Alors la cat\'egorie des rev\^etements \'etales
de~$S'$ est \'equivalente \`a la cat\'egorie des rev\^etements
\'etales $X'$ de~$S'$, o\`u le groupe $\Gamma$ op\`ere de fa\c con compatible avec ses op\'erations sur~$S'$, de telle fa\c con
que pour tout $s'\in S'$ et tout $g\in \Gamma_{s'}$} (o\`u
$\Gamma_{s'}$ d\'esigne \emph{le groupe d'inertie} de~$s'$
dans~$\Gamma$), \emph{$g$~op\`ere trivialement dans la fibre}
$X_{s'}^\prime$. Si $S'$ est connexe cet \'enonc\'e
s'interpr\`ete de la fa\c con suivante. Soit $C'_0$ la
cat\'egorie des rev\^etements \'etales de~$S'$ o\`u $\Gamma$
op\`ere de fa\c con compatible avec ses op\'erations sur~$S'$
(mais sans satisfaire n\'ecessairement la condition ci-dessus sur
les groupes d'inertie des points de~$S'$). On voit facilement que
c'est une cat\'egorie galoisienne (V~\Ref{V.5}), et que pour
\marginpar{251}
tout point g\'eom\'etrique $a'$ de~$S'$, le foncteur fibre
$X'\mto X_{a'}'$ sur~$C'_0$ est un foncteur fondamental. Soit
$\pi_1(S',\Gamma;a')=G$ le groupe des automorphismes de ce foncteur,
muni de sa topologie habituelle. On a alors une suite exacte canonique
$$
e\to \pi_1(S',a')\to G\to \Gamma\to e
$$
(cas particulier de (V~\Ref{V.6.13}), o\`u on prend pour $S$ le
rev\^etement trivial $S'\times \Gamma$ de~$S'$ d\'efini par
$\Gamma$, o\`u on fait op\'erer $\Gamma$ de fa\c con
\'evidente). D'ailleurs pour tout point g\'eom\'etrique $b'$
de~$S'$, on a un isomorphisme $\pi_1(S', \Gamma ; b' )\to G =
\pi_1(S', \Gamma; a')$ d\'efini \`a automorphisme int\'erieur
pr\`es provenant de~$\pi_1(S',a')$, et comme $\Gamma_{b'}$
s'applique de fa\c con \'evidente dans le premier membre, on
obtient un homomorphisme
$$
u_{b'}\colon \Gamma_{b'}\to G \quoi,
$$
d\'efini \`a automorphisme int\'erieur pr\`es (provenant
\ifthenelse{\boolean{orig}}
{de~$\pi_1(S',a')$,}
{de~$\pi_1(S',a')$),}
dont le compos\'e avec l'homomorphisme canonique
$G\to \Gamma$ est d'ailleurs l'immersion canonique $\Gamma_{b'} \to
\Gamma$. Ceci pos\'e, \emph{le groupe fondamental $\pi_1(S,a)$
\emph{est} canoniquement isomorphe au groupe quotient de~$G =
\pi_1(S', \Gamma; a')$ par le sous-groupe invariant ferm\'e
engendr\'e par les images des homomorphismes $\Gamma_{b'} \to
G$}. En particulier, l'image de~$\pi_1(S',a')$ dans $\pi_1(S,a)$ est
un sous-groupe invariant, et le quotient correspondant est isomorphe
\`a un quotient de~$\Gamma$. On peut d'ailleurs r\'eduire le
nombre des \og relations\fg introduites en introduisant, pour tout $g\in
\Gamma$, $ g\neq e$, le sous-pr\'esch\'ema $S'_g$ des
co\"incidences des automorphismes $\id_S$ et $g$ de~$S$, en
choisissant un point g\'eom\'etrique $b'_{g,i}$ dans chaque
composante connexe de~$S'_g$, puis un des homomorphismes
correspondants $\pi_1(S',\Gamma; b'_{g,i})\to G$, d'o\`u des
rel\`evements $\bar{g}_i$ de~$g$ dans~$G$. Il suffit alors de
prendre le quotient de~$G$ par le sous-groupe invariant ferm\'e
de~$G$ engendr\'e par les $\bar{g}_i$.

Lorsque $a'$ est invariant par $\Gamma$, on voit ais\'ement que
$\Gamma$ op\`ere de fa\c con naturelle sur~$\pi_1(S',a')$, et $G$
s'identifie au produit semi-direct correspondant. Identifiant alors~$\Gamma$
\ifthenelse{\boolean{orig}}
{a}
{\`a}
un sous-groupe de~$G$, on voit que dans les relations
\marginpar{252}
introduites plus haut, faisant $b'=a'$, on trouve \og $g = e$\fg pour
$g\in I$. Donc \emph{si $S'$ a un point g\'eom\'etrique $a'$ fixe
par}~$I$ (\ie un point $s'$ dont le groupe d'inertie est $\Gamma$),
\emph{alors $\pi_1(S,a)$ est un groupe quotient du groupe quotient de
type galoisien de~$\pi_1(S',a')$ obtenu en \emph{\og rendant
triviales\fg} les op\'erations de~$\Gamma$ sur~$\pi_1(S',a')$; et il
est m\^eme isomorphe \`a ce dernier groupe si on suppose que pour
tout $g\in G$, l'ensemble d'inertie $S_g'$ est connexe, donc passe par
la localit\'e de~$a'$}. Cette derni\`ere assertion est en effet
contenue dans la deuxi\`eme description donn\'ee plus haut pour
les relations \`a introduire dans $G$.

Ce dernier r\'esultat s'applique en particulier si l'on prend pour
$S'$ la puissance cart\'esienne $X^n$ d'un pr\'esch\'ema connexe
sur un corps alg\'ebriquement clos, pour $\Gamma$ le groupe
sym\'etrique $\Gamma={\goth {S}}_n$, op\'erant de la fa\c con
habituelle, d'o\`u pour $S$ la puissance sym\'etrique \nieme
de~$S$. Prenant alors pour $a'$ un point g\'eom\'etrique
localis\'e en la diagonale, on est sous les conditions
pr\'ec\'edentes, les ensembles d'inertie $S_g'$ contenant en effet
tous la diagonale. Utilisant le fait, prouv\'e dans l'expos\'e
suivant, que si $X$ est propre connexe sur~$k$, le groupe fondamental
de~$X^n$ s'identifie \`a $\pi_1(X)^n$, on trouve le r\'esultat
amusant suivant: \emph{Si $X$ est propre connexe sur~$k$
alg\'ebriquement clos, le groupe fondamental de sa puissance
sym\'etrique \nieme, $n \geq 2$, est isomorphe au groupe
fondamental de~$X$ rendu ab\'elien}. (J'ignore si le fait analogue
en Topologie alg\'ebrique est connu; il devrait pouvoir
s'\'etablir par la m\^eme m\'ethode de descente). Prenons par
exemple pour $X$ une courbe rationnelle $X=\PP^1_k$, on trouve une
\Nieme d\'emonstration du fait que $\PP^r_k$ est simplement
connexe, utilisant le fait que $\PP^1_k$ l'est. Prenons maintenant
pour $X$ une courbe simple sur~$k$, et $n\geq 2g-1$, de sorte que
$\Symm^n(X)$ est fibr\'e sur la jacobienne $J$, de fibres des
espaces projectifs, donc (comme on verra \`a l'aide des
r\'esultats des deux expos\'es suivants) a m\^eme groupe
fondamental que $J$. On retrouve alors sans d\'evissage le fait bien
connu que \emph{le groupe fondamental de la jacobienne de~$X$ est
isomorphe au groupe fondamental de~$X$ rendu ab\'elien}.
\end{remarque}

\section[Une suite exacte fondamentale]{Une suite exacte fondamentale. Descente par morphismes \`a
fibres relativement connexes}
\label{IX.6}
\marginpar{253}

\begin{theoreme}
\label{IX.6.1}
Soient $S$ le spectre d'un anneau artinien $A$ de corps r\'esiduel
$k$, $\overline{k}$ une cl\^oture alg\'ebrique de~$k$, $X$ un $S$
pr\'esch\'ema, $X_0=X\otimes_A k$, $\overline{X}_0=X\otimes_A
\overline{k}$, $\overline{a}$ un point g\'eom\'etrique
de~$\overline{X}$, $a$ son image dans~$X$, $b$ son image dans~$S$. On
suppose que $X_0$ est quasi-compact et \emph{g\'eom\'etriquement
connexe} sur~$k$ (N.B. si $X$ est propre sur~$S$, cela signifie que
$\H^0(X_0,\cal{O}_{X_0})$ est un anneau artinien \emph{local} de corps
r\'esiduel \emph{radiciel} sur~$k$). Alors la suite d'homomorphismes
canoniques
$$
e\mto \pi_1(\overline{X}_0,\overline{a})\to
\pi_1(X,a)\to\pi_1(S,b)\to e
$$
est exacte, et on a
$$
\pi_1(S,b)\isomfrom \pi_1(k,\overline{k})=\text{groupe de Galois
de~$\overline{k}$ sur~$k$.}
$$
\end{theoreme}

Comme les groupes fondamentaux ne changent pas en ch\^atrant par les
\'el\'ements nilpotents, on peut supposer $A=k$, ce qui rend
d\'ej\`a \'evident le dernier isomorphisme. Soit $k'$ la
cl\^oture s\'eparable de~$k$ dans~$\overline{k}$, et
consid\'erons $X'=X\otimes_k k'$, et l'image $a'$ de~$\overline{a}$
dans~$X'$. On a une suite d'homomorphismes canoniques
$$
e\mto \pi_1(\overline{X}_0,\overline{a})\to
\pi_1(X',a')\to\pi_1(S',b')\to e
$$
(o\`u $S'=\Spec(k')$). Enfin, on a un homomorphisme canonique de
cette suite dans celle relative \`a~$X/k$, gr\^ace au~diagramme
$$
\xymatrix{ S & X \ar[l] & \overline{X}_0 \ar[l] \\ S' \ar[u] &
X'\ar[u] \ar[l] & \overline{X}_0 \ar[l] \ar@{=}[u]. }
$$
On voit d'autre part que cet homomorphisme de suites de groupes est un
isomorphisme, comme il r\'esulte de~\Ref{IX.4.11}. On est donc
ramen\'e \`a prouver que la deuxi\`eme suite est exacte, \ie on
peut supposer que $k$ est \emph{parfait}. Soient alors $k_i$ les
sous-extensions galoisiennes finies de~$k$ dans~$\overline{k}$, posons
$X_i=X\otimes_k k_i$, et soit $a_i$ l'image de~$\overline{a}$
dans~$X_i$. On laisse au lecteur
\ifthenelse{\boolean{orig}}
{\ignorespaces}
{le soin}
de v\'erifier que l'homomorphisme
naturel
$$
\pi_1(\overline{X}_0,a)\isomto \varprojlim\pi_1(X_i,a_i)
$$
\marginpar{254}%
est un isomorphisme, ce qui signifie simplement qu'un rev\^etement
\'etale de~$\overline{X}$ provient d'un rev\^etement \'etale
d'un~$X_i$, et que ce dernier est essentiellement unique, modulo
passage \`a un~$X_j$, $j\geq i$. D'autre part, soit $\pi_i$ le
groupe de Galois de~$k_i$ sur~$k$, \ie le groupe oppos\'e au groupe
des $S$-automorphismes de~$S_i=\Spec(k_i)$. Comme le foncteur
$S'\mto X\times_S S'$ des rev\^etements \'etales de~$S$ dans
les rev\^etements \'etales de~$X$ est pleinement fid\`ele
\eqref{IX.3.4}, il s'ensuit que $\pi_i$ est aussi isomorphe au groupe
oppos\'e aux groupes des $X$-automorphismes du rev\^etement
principal connexe $X_i$ de~$X$. Il r\'esulte donc de
V~\Ref{V.6.13} que l'on a une suite exacte
$$
e\mto \pi_1(X_i,a_i)\to \pi_1(X,a)\to\pi_i\to e.
$$
Passant \`a la limite projective sur~$i$ dans ces suites exactes, on
trouve une suite exacte (puisqu'on est dans la cat\'egorie des
groupes de type galoisien), qui n'est autre que la suite envisag\'ee
dans~\Ref{IX.6.1}. Cela ach\`eve la d\'emonstration.

La traduction de l'exactitude \`a droite dans~\Ref{IX.6.1} en
langage g\'eom\'etrique est la suivante:

\begin{corollaire}
\label{IX.6.2}
Avec les notations pr\'ec\'edentes, soit $X'$ un rev\^etement
\'etale de~$X$, et soit $\overline{X'}_0$ le rev\^etement
\'etale correspondant de~$\overline{X}_0$. Les conditions suivantes
sont \'equivalentes:
\begin{enumerate}
\item[(i)] Il existe un $S'$ \'etale sur~$S$, et un $X$-isomorphisme
$X'\isomto X\times_S S'$, ($S'$ est alors d\'etermin\'e \`a
isomorphisme unique pr\`es en vertu de~\Ref{IX.3.4}).
\item[(ii)] $\overline{X}_0'$ est compl\`etement d\'ecompos\'e
sur~$\overline{X}_0$.

Si $X'$ est connexe, ces conditions \'equivalent aussi \`a:
\item[(ii~bis)] $\overline{X}_0'$ a une section sur~$\overline{X}_0$.
\end{enumerate}
\end{corollaire}

(N.B. Ce dernier compl\'ement est essentiel; l'\'equivalence
de~(i) et (ii) signifie seulement que $\pi_1(S,b)$ est le groupe
quotient de~$\pi_1(X,a)$ par le sous-groupe
\marginpar{255}
invariant ferm\'e engendr\'e par l'image
de~$\pi_1(\overline{X}_0,\overline{a})$, et non par cette image
elle-m\^eme). Sous les conditions pr\'ec\'edentes, nous dirons
que $X'$ est un rev\^etement \emph{g\'eom\'etriquement trivial}
\index{geometriquement trivial (revetement)@g\'eom\'etriquement trivial (rev\^etement)|hyperpage}%
de~$X$.

\begin{remarque}
\label{IX.6.3}
On ne peut dans l'\'enonc\'e~\Ref{IX.6.1} remplacer $\overline{k}$
par une extension alg\'ebriquement close quelconque de~$k$, m\^eme
si $k$ est d\'ej\`a suppos\'e alg\'ebriquement clos. En
d'autres termes, il n'est pas vrai en g\'en\'eral que si $X$ est
un sch\'ema alg\'ebrique connexe sur un corps alg\'ebriquement
clos $k$, son groupe fondamental ne change pas en rempla\c cant $k$
par une extension alg\'ebriquement close; c'est d\'ej\`a faux
par exemple en caract\'eristique $p>0$ pour la droite affine
sur~$k$, \`a cause des ph\'enom\`enes de \og ramification
sup\'erieure\fg au point \`a l'infini, impliquant une structure
\og continue\fg pour le groupe fondamental. Nous verrons cependant dans
l'expos\'e suivant que de tels ph\'enom\`enes ne peuvent se
produire si $X$ est \emph{propre} sur~$k$. Nous montrerons aussi par
voie transcendante qu'il en est de m\^eme si $k$ est de
caract\'eristique nulle.
\end{remarque}

\begin{corollaire}
\label{IX.6.4}
Supposons que $a$ soit localis\'e en un $x\in X$ qui est rationnel
sur~$k$ (ou plus g\'en\'eralement, ayant un corps r\'esiduel
radiciel sur~$k$). Alors la suite exacte~\Ref{IX.6.1} est scind\'ee.
\end{corollaire}

On peut supposer $S=\Spec(k)$. Si $x$ est rationnel sur~$k$, il
correspond \`a une section $S\to X$ de~$X$ sur~$S$, transformant $b$
en~$a$, et d\'efinissant un homomorphisme $\pi_1(S,b)\to \pi_1(X,a)$
qui est le scindage cherch\'e. Si $\kres(x)$ est radiciel sur~$k$, on se
ram\`ene au cas pr\'ec\'edent en faisant l'extension de la base
$\Spec(\kres(x))\to\Spec(k)$.

\begin{theoreme}
\label{IX.6.5}
Soient $f\colon X\to S$ un morphisme propre et surjectif de
pr\'esentation finie, \`a fibres g\'eom\'etriquement connexes,
$X'$ un pr\'esch\'ema de pr\'esentation finie et propre sur~$X$,
$s$ un point de~$S$, $F=X_s$ la fibre de~$X$ en~$s$, et $F_1'$ une
composante connexe de la fibre $F'=X'_s$ de~$X'$ en~$s$. Pour qu'il
existe un voisinage ouvert $X_1'$ de~$F_1'$ dans~$X'$, un
$S$-sch\'ema \'etale $S_1'$ et un $X$-isomorphisme $X_1'\isomto
S_1'\times_S X$, il faut et il suffit que $X'$ soit \'etale sur~$X$
aux points de~$F_1'$, et que $F_1'$ soit un rev\^etement
g\'eom\'etriquement trivial de~$F$.
\end{theoreme}

La
\marginpar{256}
n\'ecessit\'e de la condition \'etant triviale, il reste
\`a prouver la suffisance. On se ram\`ene facilement au cas o\`u
$S$ est noeth\'erien. Consid\'erons la factorisation de Stein
$X\to T\to S$ de~$f$, o\`u $T$ est le spectre de l'Alg\`ebre
$f_*(\cal{O}_X)$ sur~$S$. Comme les fibres de~$X$ sur~$S$ sont
g\'eom\'etriquement connexes, et $f$ est surjectif, le morphisme
$T\to S$ est fini surjectif et radiciel, donc \eqref{IX.4.10} tout
$T'$ \'etale sur~$T$ provient par image inverse d'un $S'$ \'etale
sur~$S$. Cela nous ram\`ene \`a
\ifthenelse{\boolean{orig}}
{prouver~\Ref{IX.6.5}.}
{prouver~\Ref{IX.6.5}}
en y rempla\c cant $S$ par~$T$, \ie dans le cas o\`u on suppose
$f_*(\cal{O}_X)=\cal{O}_S$. Consid\'erons alors la factorisation de
Stein $X'\to S'\to S$ du morphisme propre $h\colon X'\to S$, o\`u
$S'$ est le spectre de l'Alg\`ebre $h_*(\cal{O}_{X'})$. Les
morphismes $X'\to X$ et $X'\to S'$ d\'efinissent un morphisme
canonique
$$
X'\to X\times_S S',
$$
et notre assertion est contenue dans la suivante:

\begin{corollaire}
\label{IX.6.6}
Soit $f\colon X\to S$ un morphisme propre de pr\'esch\'emas
localement noeth\'eriens, tel que $f_*(\cal{O}_X)=\cal{O}_S$, et
soit $X'$ un pr\'esch\'ema propre sur~$X$. Consid\'erons la
factorisation de \emph{Stein} $X'\to S'\to S$ pour $X'\to S$, et le
morphisme canonique $X'\to X\times_S S'$. Soient $s$ un point de~$S$,
$s'$ un point de~$S'$ au-dessus de~$s$, correspondant \`a une
composante connexe $F_1'$ de la fibre $X'_s$ de~$X'$ en~$s$. Pour que
le morphisme $X'\to X\times_S S'$ soit un isomorphisme au-dessus d'un
voisinage ouvert $U'$ de~$s'$ \'etale sur~$S$, il faut et il suffit
que $X'$ soit \'etale sur~$X$ en les points de~$F_1'$, et que $F_1'$
soit un rev\^etement g\'eom\'etriquement trivial de la fibre
$F=X_s$.
\end{corollaire}

La n\'ecessit\'e \'etant encore triviale, il reste \`a prouver
la suffisance. La conclusion signifie aussi que a) le morphisme
d\'eduit de~$X'\to X\times_S S'$ par le changement de base
$\Spec(\widehat{\cal{O}}_{s'})\to S'$ est un isomorphisme, b) $S'$ est
\'etale sur~$S$ en~$s'$, \ie $\widehat{\cal{O}}_{s'}$ est \'etale
sur~$\widehat{\cal{O}}_s$. Sous cette forme, on voit que la conclusion
est invariante par changement de base $\Spec(\widehat{\cal{O}}_s)\to
S$. Les hypoth\`eses \'etant \'egalement stables par ce
changement de base, on peut donc supposer que $S$ est le spectre d'un
anneau local noeth\'erien complet. On peut de plus \'evidemment
supposer $X'$ connexe, ce qui implique ici que
$S'=\Spec(\cal{O}_{s'})$, $F'=F_1'$.

Comme
\marginpar{257}
l'ensemble des points de~$X'$ o\`u $X'$ est \'etale sur~$X$ est
ouvert et contient la fibre $X'_{s'}=F'$, $X'$ \'etant propre
sur~$S$, il s'ensuit que $X'$ est \'etale sur~$X$. Comme il induit
sur~$F=X_s$ un rev\^etement \'etale isomorphe \`a un
$F\otimes_{\kres(s)}L$, o\`u $L$ est \'etale sur~$\kres(s)$, il
r\'esulte de~\Ref{IX.1.10} qu'il est isomorphe \`a un
rev\^etement de la forme $X\times_S T$, avec $T$ \'etale sur~$S$.
(N.B. ici encore, il suffit d'utiliser que le foncteur
de~\Ref{IX.1.10} est pleinement fid\`ele, qui r\'esulte du fait
qu'un isomorphisme formel de faisceaux coh\'erents sur~$X$ provient
d'un isomorphisme de ces faisceaux). Donc, si $T$ est d\'efini par
l'alg\`ebre~$B$ finie sur~$A$, $X'$ s'identifie au spectre de
l'Alg\`ebre $\cal{O}_X\otimes_A B$ sur~$X$, d'o\`u r\'esulte
aussit\^ot, puisque $f_*(\cal{O}_X)=\cal{O}_S$, que
$h_*(\cal{O}_{X'})$ est d\'efini par~$B$, donc l'homomorphisme
canonique $X'\to X\times_S S'$ n'est autre que l'isomorphisme
envisag\'e $X'\isomto X\times_S T$. Cela ach\`eve la
d\'emonstration.

\begin{corollaire}
\label{IX.6.7}
Sous les conditions de~\Ref{IX.6.5}, pour qu'il existe un
pr\'esch\'ema $S'$ \'etale sur~$S$ et un $X$-isomorphisme
$X'\isomto X\times_S S'$, il faut et il suffit que $X'$ soit \'etale
sur~$X$ et que pour tout $s\in S$, $X'_{s'}$ soit un rev\^etement
g\'eom\'etriquement trivial de~$X_s$.
\end{corollaire}

En effet, s'il en est ainsi, $X'$ est r\'eunion d'ouverts $X_i'$ qui
sont isomorphes \`a des images inverses de~$S_i'$ \'etales
sur~$S$. On voit alors facilement que ces $S_i'$ se recollent en un
$S'$ \'etale sur~$S$, et qu'on obtient un isomorphisme $X'\isomto
X\times_S S'$. On peut dire par exemple que les $X_i'$ sont munis de
donn\'ees de descente relativement \`a~$X\to S$, qui se recollent
n\'ecessairement en une donn\'ee de descente sur~$X'$ tout entier
relativement \`a~$X\to S$. Et comme cette derni\`ere est effective
sur les~$X_i'$, il s'ensuit facilement (gr\^ace \`a un sorite
oubli\'e au num\'ero \Ref{IX.4}) qu'elle est effective. On peut aussi
\'enoncer \Ref{IX.6.7} sous la forme suivante:

\begin{corollaire}
\label{IX.6.8}
Soit $f\colon X\to S$ un morphisme propre surjectif de
pr\'esentation finie, \`a fibres g\'eom\'etriquement
connexes. Alors $f$ est un morphisme de descente effective pour la
cat\'egorie fibr\'ee des pr\'esch\'emas \'etales finis sur
d'autres. Le foncteur $S'\mto X\times_S S'$ induit une
\'equivalence de la cat\'egorie des pr\'esch\'emas \'etales
et finis sur~$S$ avec la cat\'egorie des pr\'esch\'emas
\marginpar{258}
\'etales et finis sur~$X$ qui induisent sur chaque fibre $X_s$ un
rev\^etement g\'eom\'etriquement trivial.
\end{corollaire}

\begin{remarque}
\label{IX.6.9}
Soit $f\colon X\to S$ un morphisme propre et surjectif, avec $Y$
localement noeth\'erien, alors $f$ se factorise en un morphisme
$X\to S'$ satisfaisant l'hypoth\`ese de~\Ref{IX.6.8}, et un
morphisme fini surjectif $S'\to S$ justiciable de~\Ref{IX.4.7}, donc
en produit de deux morphismes qui sont des \emph{morphismes de
descente effective universels} pour la cat\'egorie fibr\'ee des
pr\'esch\'emas \'etales et finis sur d'autres. On peut en
conclure que $f$ lui-m\^eme est un morphisme de descente effective
universel pour la cat\'egorie fibr\'ee envisag\'ee. On retrouve
ainsi \Ref{IX.4.12} par une m\'ethode diff\'erente.
\end{remarque}

\begin{remarque}
\label{IX.6.10}
La conclusion de~\Ref{IX.6.7} ne reste pas valable si on remplace
l'hypoth\`ese que $f$ est propre par: $X$ est de type fini sur~$S$
et admet une section sur~$S$ (donc $f$ est universellement submersif
et un morphisme de descente pour la cat\'egorie fibr\'ee des
pr\'esch\'emas \'etales sur d'autres), m\^eme lorsque $S$ est
le spectre d'un anneau de valuation discr\`ete et lorsque $X'$ est
un rev\^etement \'etale de~$X$. Pour le voir, on part d'un $Z$
propre sur~$S$, dont la fibre g\'en\'erale
est une courbe rationnelle non singuli\`ere, et la fibre
sp\'eciale $Z_0$ consiste en deux droites concourantes. Par exemple,
si $t$ est une uniformisante de l'anneau de valuation~$A$, on prend le
sous-sch\'ema ferm\'e $Z$ de~$\PP_A^2$ d\'efini par
l'\'equation homog\`ene $x^2+y^2+tz^2=0$ (coordonn\'ees
homog\`enes $x$, $y$, $z$). On prend pour $X$ le compl\'ementaire
du point singulier $a$ de~$Z_0$ dans la r\'eunion $Z\cup
\PP_k^2$. Les fibres de~$X$ sont $\PP_k^1$ et $\PP_k^2 - a$, donc sont
g\'eom\'etriquement simplement connexes (\ie tout rev\^etement
\'etale d'une telle fibre est g\'eom\'etriquement trivial).
Cependant
\marginpar{259}
on construit facilement, en proc\'edant comme dans le \No \Ref{IX.4}, des
rev\^etements \'etales de~$X$ qui ne proviennent pas de
rev\^etements \'etales de~$S$, par recollement de rev\^etements
triviaux de~$Z-a$ et de~$\PP_k^2-a$. Il est possible par contre que la
conclusion de~\Ref{IX.6.7} subsiste si on y remplace l'hypoth\`ese
de propret\'e par celle que $X$ soit \emph{universellement ouvert}
de pr\'esentation finie sur~$S$\kern1pt\footnote{C'est maintenant
prouv\'e, $g$ \'etant seulement universellement ouvert et
surjectif; \cf SGA~4 XV~1.15.}. C'est vrai du moins si on suppose que
les fibres de~$X$ sur~$S$ sont g\'eom\'etriquement
irr\'eductibles, et non seulement g\'eom\'etriquement connexes.
Signalons seulement qu'on peut dans cette question se ramener au cas
o\`u $S$ est le spectre d'un anneau de valuation discr\`ete
complet \`a corps r\'esiduel alg\'ebriquement clos.
\end{remarque}

L'interpr\'etation de~\Ref{IX.6.7} en termes du groupe fondamental
est la suivante:

\begin{corollaire}
\label{IX.6.11}
Soit $f\colon X\to S$ un morphisme propre surjectif de
pr\'esentation finie, \`a fibres g\'eom\'etriquement
connexes. On suppose $X$ donc $S$ connexe. Soient $a$ un point
g\'eom\'etrique de~$X$, $b$ son image dans~$S$, et pour tout $s\in
S$, choisissons une cl\^oture alg\'ebrique $\overline{\kres(s)}$ de
$\kres(s)$, un point g\'eom\'etrique $a_s$ de~$S_s$ \`a valeurs
dans cette extension, et une classe de chemins de~$a_s$ dans~$a$,
d'o\`u un homomorphisme $\pi_1(\overline X_s,a_s)\to \pi_1(X,a)$,
o\`u $\overline X_s=X_s\otimes_{\kres(s)}\overline{\kres(s)}$. Alors
l'homomorphisme $\pi_1(X,a)\to \pi_1(S,b)$ est surjectif, et son noyau
est le sous-groupe invariant ferm\'e de~$\pi_1(X,a)$ engendr\'e
par les images des $\pi_1(\overline X_s,a_s)$.
\end{corollaire}

\begin{remarque}
\label{IX.6.12}
Sous les conditions de~\Ref{IX.6.7}, supposant $S$ noeth\'erien,
on voit facilement que l'ensemble des points $s\in S$ tels que
$S'_{s'}$ soit g\'eom\'etriquement trivial sur~$X_s$ est un
ensemble constructible; si $X'$ est propre sur~$X$, il est m\^eme
ouvert, comme on voit sur~\Ref{IX.6.6}. Cela permet donc, si $S$ et
un pr\'esch\'ema de Jacobson (par exemple de type fini sur un
corps), ou
\marginpar{260}
$X'$ est propre sur~$X$, de se borner pour v\'erifier les conditions
de~\Ref{IX.6.7} aux points $s$ de~$S$ qui sont \emph{ferm\'es}. De
m\^eme, dans~\Ref{IX.6.11} il suffit alors de prendre les
$\pi_1(\overline X_s,a_s)$ pour les points de~$S$ qui sont ferm\'es.
\end{remarque}


\chapter{Th\'eorie de la sp\'ecialisation du~groupe~fondamental}
\label{X}
\marginpar{261}
Dans le pr\'esent expos\'e, nous nous bornons \`a l'\'etude du
groupe fondamental des fibres g\'eom\'etriques dans un morphisme
\emph{propre}, \ie du groupe fondamental d'un sch\'ema
alg\'ebrique \emph{propre} variable. Dans un expos\'e
ult\'erieur, nous g\'en\'eraliserons la technique employ\'ee
aux rev\^etements \'etales \emph{mod\'er\'ement ramifi\'es}
\og \`a l'infini\fg. Cela nous donnera par exemple une solution du
\og Probl\`eme des trois points\fg dans le cas des rev\^etements
galoisiens d'ordre premier \`a la caract\'eristique (\ie une
d\'etermination des rev\^etements galoisiens de la droite
$\PP_k^1$, ramifi\'es au plus en trois points donn\'es et
mod\'er\'ement ramifi\'es en ces points), et de ses variantes
\'evidentes.

\section[La suite exacte d'homotopie pour un morphisme propre]
{La suite exacte d'homotopie pour un morphisme propre et s\'eparable}
\label{X.1}

\begin{definition}
\label{X.1.1}
Un pr\'esch\'ema $X$ sur un corps $k$ est dit
\emph{s\'eparable}, ou \emph{s\'eparable sur~$k$}, si pour toute
extension $K$ de~$k$, $X\otimes_k K$ est r\'eduit. Si $f\colon X\to
Y$ est un morphisme de pr\'esch\'emas, on dit que $f$ est
\emph{s\'eparable}, ou que $X$ \emph{est s\'eparable sur~$Y$}, si
$X$ est plat sur~$Y$ et si pour tout $y\in Y$, la fibre $X\otimes_Y
\kres(y)$ est s\'eparable sur~$\kres(y)$.
\end{definition}

Si $X$ est un pr\'esch\'ema sur un corps~$k$, dire qu'il est
s\'eparable signifie aussi qu'il est \emph{r\'eduit}, et que les
corps~$\kres(x)$, pour $x$ point g\'en\'erique d'une composante
irr\'eductible de~$X$, sont des extensions s\'eparables de~$k$. Si
$k$ est parfait, il revient donc au m\^eme de dire que $X$ est
r\'eduit. Notons que si $X$ est s\'eparable sur~$Y$, alors pour
tout changement de base $Y'\to Y$, $X'=X\times_Y Y'$ est s\'eparable
sur~$Y'$. On peut prouver aussi, moyennant des hypoth\`eses de
finitude convenables, que le compos\'e de deux morphismes
s\'eparables est
\marginpar{262}
un morphisme s\'eparable. Nous en aurons besoin seulement sous la
forme suivante: \emph{si $X$ est s\'eparable sur~$Y$, et $X'$
\'etale sur~$X$, $X'$ est s\'eparable sur~$Y$.} C'est en effet
une cons\'equence imm\'ediate des d\'efinitions
et~I~\Ref{I.9.2}. Par ailleurs, l'hypoth\`ese \og morphisme
s\'eparable\fg nous servira par l'interm\'ediaire de la proposition
suivante:
\begin{proposition}
\label{X.1.2}
Soit $f\colon X\to Y$ un morphisme propre et s\'eparable, avec $Y$
localement noeth\'erien, et consid\'erons sa factorisation de
\emph{Stein} $X\lto{f'}Y'\dpl\to Y$, (o\`u
$f_*'(\cal{O}_X)=\cal{O}_{Y'}$, $Y'$ \'etant fini sur~$Y$ et
isomorphe au spectre de l'alg\`ebre $f_*(\cal{O}_X)$). Alors $Y'$
est un \emph{rev\^etement \'etale} de~$Y$.
\end{proposition}

Cette proposition figurera dans (EGA~III 7)\footnote{\Cf EGA~III
7.8.10 (i)}. Indiquons le principe de la d\'emonstration. On se
ram\`ene facilement au cas o\`u $Y$ est le spectre d'un anneau
local complet~$A$, et faisant encore une extension finie plate
convenable de ce dernier (correspondant \`a une extension
r\'esiduelle convenable), on peut supposer que les composantes
connexes de la fibre du point ferm\'e $y$ sont
g\'eom\'etriquement connexes, ce qui signifie aussi que
$\H^0(X_y,\cal{O}_{X_y})$ se d\'ecompose en un produit de corps
identiques \`a $k=\kres(y)$. Supposant alors $X$ connexe, ce qui est
loisible, on aura $\H^0(X_y,\cal{O}_{X_y})=k$, donc l'homomorphisme
$A\to \H^0(X_y,\cal{O}_{X_y})$ est \emph{surjectif}. On en conclut par
une proposition g\'en\'erale (du type K\"unneth) que
$f_*(\cal{O}_X)$ est d\'efini par un module $B$ sur~$A$ qui est
libre sur~$A$, et que $B/\goth{m} B\to \H^0(X_y,\cal{O}_{X_y})=k$ est
bijectif. Donc en l'occurrence $B$ est une alg\`ebre \'etale
sur~$A$, ce qui ach\`eve la d\'emonstration.

\begin{theoreme}
\label{X.1.3}
Soit $f\colon X\to Y$ un morphisme propre et s\'eparable, avec $Y$
localement noeth\'erien et connexe, et supposons
$f_*(\cal{O}_X)=\cal{O}_Y$ (ce qui implique que les fibres de~$X$ sur~$Y$ sont g\'eom\'etriquement connexes, et r\'eciproquement
gr\^ace \`a~\Ref{X.1.2}). Soient~$y$ un point de~$Y$,
$\overline{\kres(y)}$ une cl\^oture alg\'ebrique de~$\kres(y)$,
$\overline X_y =X_y\otimes_{\kres(y)} \overline{\kres(y)}$. Soient \hbox{enfin~$X'$} un rev\^etement \'etale \emph{connexe} de~$X$, et $\overline
X_y' =X_y'\otimes_{\kres(y)} \overline{\kres(y)}$. Pour qu'il existe un
rev\^etement \'etale $Y'$ de~$Y$ et un $X$-isomorphisme
\marginpar{263}
$X'\isomto X\times_Y Y'$, il faut et il suffit que $\overline{X'}_y$
admette une section sur~$\overline X_y$.
\end{theoreme}

Posant $Y'=\Spec(h_*(\cal{O}_{X'}))$ (o\`u $h\colon X'\to Y$ est le
morphisme compos\'e $X'\to X\to Y$), il suffit de prouver que le
$Y$-morphisme canonique
$$
X'\to X\times_Y Y'
$$
est un \emph{isomorphisme}, et que $Y'$ est \'etale sur~$Y$. Or nous
savons d\'ej\`a par~\Ref{X.1.2} que $Y'$ est \'etale sur~$Y$,
donc $X\times_Y Y'$ est \'etale sur~$X$, donc le morphisme $X'\to
X\times_Y Y'$ est \'egalement \'etale (I~\Ref{I.4.8}). D'ailleurs,
$Y'$ est connexe comme image de~$X'$ qui l'est, donc $X\times_Y Y'$
est connexe puisque $X$ est \`a fibres connexes sur~$Y$
(IX~\Ref{IX.3.4} et V~\Ref{V.6.9}~(iii)). Donc pour prouver que
$X'\to X\times_Y Y'$ est un isomorphisme, il suffit de voir que son
degr\'e de projection en \emph{un} point de~$X\times_Y Y'$ est
\'egal \`a~1. Or ceci r\'esulte facilement de l'hypoth\`ese
que $\overline{X'}_y$ admet une section sur~$\overline X_y$, soit
par utilisation de IX~\Ref{IX.6.6}, soit plus simplement en notant qu'il
suffit de prouver l'existence d'un tel point dans $X\times_Y Y'$
apr\`es changement de base
\ifthenelse{\boolean{orig}}
{$\Spec(\overline k)\to Y$}
{$\Spec(\overline{\kres(y)})\to Y$}, o\`u cela
est \'evident. Cela ach\`eve la d\'emonstration de~\Ref{X.1.3}.

Tenant compte de IX~\Ref{IX.3.4} et du dictionnaire V~\Ref{V.6.9} et
V~\Ref{V.6.11}, on peut mettre~\Ref{X.1.3} sous la forme
\'equivalente suivante

\begin{corollaire}
\label{X.1.4}
Avec les notations pr\'ec\'edentes pour $f\colon X\to Y$ et
$\overline X_y$, soit $\bar{a}$ un point g\'eom\'etrique de
$\overline X_y$, $a$ son image dans $X$ et $b$ son image
dans~$Y$. Alors la suite suivante d'homomorphismes de groupes est
\emph{exacte}:
$$
\pi_1(\overline X_y,\bar{a})\to \pi_1(X,a)\to \pi_1(Y,b)\to e.
$$
\end{corollaire}

\begin{remarques}
\label{X.1.5}
On notera que la d\'emonstration de~\Ref{X.1.3} fait intervenir de
fa\c con essentielle~\Ref{X.1.2} et par l\`a le \og premier
th\'eor\`eme de comparaison\fg en g\'eom\'etrie
alg\'ebrico-formelle. Par contre, la th\'eorie de la descente de
l'expos\'e~\Ref{IX}
\marginpar{264}
n'est intervenue que par l'interm\'ediaire de IX~\Ref{IX.3.4}, dont
une d\'emonstration directe dans le cas d'un morphisme \emph{propre}
$f\colon X\to Y$ tel que $f_*(\cal{O}_X)=\cal{O}_Y$ est facile. Soit
en effet $Y'$ \'etale sur~$Y$ et supposons que $X' = X\times_Y Y'$
soit somme disjointe de deux ouverts non vides, prouvons qu'il en est
de m\^eme de~$Y'$. En effet, on aura $Y'=\Spec(\cal{A})$, donc
$X'=\Spec(\cal{B})$ avec $\cal{B}=\cal{A}\otimes_{\cal{O}_Y}
\cal{O}_X$, et la d\'ecomposition de~$X'$ en somme directe
correspond \`a une d\'ecomposition de~$\cal{B}$ en produit de deux
Alg\`ebres non nulles $\cal{B}_1$ et $\cal{B}_2$. Comme
$f_*(\cal{O}_X)=\cal{O}_Y$ on conclut facilement
$f_*(\cal{B})=\cal{A}$, donc $\cal{A}$ sera somme de deux Alg\`ebres
(non nulles \'egalement, car leurs sections unit\'e sont non
nulles) $f_*(\cal{B}_1)$ et $f_*(\cal{B}_2)$, cqfd.
\end{remarques}

\subsection{}
\label{X.I.6}
Supposons encore que $f$ soit propre et s\'eparable, mais ne faisons
plus d'hypoth\`ese sur~$f_*(\cal{O}_X)$, qui correspondra \`a un
rev\^etement \'etale $Y'$ bien d\'etermin\'e de~$Y$,
d'ailleurs ponctu\'e au-dessus de~$b$ par l'image $b'$
de~$a$. Appliquant alors \Ref{X.1.4} au morphisme canonique $X\to Y'$,
et supposant $f$ surjectif, la suite exacte~\Ref{X.1.4} est
remplac\'ee par la suivante, analogue de la suite exacte d'homotopie
des espaces fibr\'es en topologies alg\'ebriques:
$$
\pi_1(\overline X_y,\overline a)\to \pi_1(X,a)\to \pi_1(Y,b)\to
\pi_0(\overline X_y,\overline a)\to \pi_0(X,a)\to \pi_0(Y,b)\to e.
$$
Bien entendu, dans~\Ref{X.1.4} on ne peut pas en g\'en\'eral
affirmer que l'homomorphisme
$\pi_1(\overline X_y,\overline a)\to \pi_1(X,a)$
soit injectif; en topologie alg\'ebrique, son noyau est l'image de
$\pi_2(Y,b)$, et il y aurait lieu en g\'eom\'etrie alg\'ebrique
\'egalement d'introduire des groupes d'homotopie en toutes
dimensions, et la suite exacte d'homotopie compl\`ete pour un
morphisme propre satisfaisant des hypoth\`eses convenables (par
exemple d'\^etre un morphisme
\ifthenelse{\boolean{orig}}
{lisse.}
{lisse).}
On ne dispose \`a l'heure
actuelle d'aucun r\'esultat dans ce sens, \`a l'exception d'une
d\'efinition raisonnable (sinon d\'efinitive) des groupes
d'homotopie sup\'erieure.

\begin{corollaire}
\label{X.1.7}
Soient $k$ un corps alg\'ebriquement clos, $X$ et $Y$ deux
pr\'esch\'emas connexes sur~$k$, on suppose $X$ propre sur~$k$ et
$Y$ localement
\marginpar{265}
noeth\'erien. Soient $a$ un point g\'eom\'etrique de~$X$, $b$
un point g\'eom\'etrique de~$Y$ \`a valeurs dans la m\^eme
extension alg\'ebriquement close $K$ de~$k$, consid\'erons le
point g\'eom\'etrique $c=(a,b)$ de~$X\times_k Y$, et
l'homomorphisme $\pi_1(X\times_k Y,c)\to \pi_1(X,a)\times\pi_1(Y,b)$
d\'eduit des homomorphismes sur les groupes fondamentaux
associ\'es aux deux projections $X\times_k Y\to X$ et $X\times_k
Y\to Y$. L'homomorphisme pr\'ec\'edent est un \emph{isomorphisme}.
\end{corollaire}

Supposons d'abord $K=k$. Posons $Z=X\times_k Y$, consid\'erons la
projection $f\colon Z\to Y$ et la localit\'e $y$ du point
g\'eom\'etrique $b$ de~$Y$, appliquons \`a la situation le
r\'esultat~\Ref{X.1.4}. On notera pour ceci que quitte \`a passer
\`a $X_\red$ (ce qui ne change pas les groupes fondamentaux
envisag\'es), on peut supposer d\'ej\`a $X$ r\'eduit donc
s\'eparable sur~$k$, donc~$Z$ est s\'eparable sur~$k$, et
\'evidemment \`a fibres g\'eom\'etriquement connexes (puisque~$X$ est connexe). La fibre g\'eom\'etrique de~$Z$ en $b$ est
canoniquement isomorphe \`a $X\otimes_k K=X$. D'autre part, comme le
compos\'e des morphismes $X\to Z\to X$ est l'identit\'e, on trouve
que $\pi_1(X,a)\to \pi_1(Z,c)$ est injectif et 1.4 nous donne une
suite exacte:
$$
e\to \pi_1(X,a)\to \pi_1(Z,c)\to \pi_1(Y,b)\to e.
$$
D'autre part, on a la suite exacte canonique:
$$
e\to \pi_1(X,a)\to \pi_1(X,a)\times \pi_1(Y,b)\to \pi_1(Y,b)\to e,
$$
o\`u les deux homomorphismes \'ecrits sont l'injection canonique
et la projection canonique. On a enfin un homomorphisme de la
premi\`ere suite exacte dans la deuxi\`eme, \`a l'aide des
morphismes identiques sur les termes extr\^emes, et l'homomorphisme
canonique $\pi_1(Z,c)\to \pi_1(X,a)\times\pi_1(Y,b)$ pour les termes
m\'edians. La commutativit\'e du diagramme ainsi obtenu se
v\'erifie trivialement. Comme les homomorphismes sur les termes
extr\^emes sont des isomorphismes, il en est de m\^eme pour les
termes m\'edians, ce qui prouve \Ref{X.1.7} dans ce cas.

Lorsqu'on ne suppose plus $K=k$, on trouve seulement un isomorphisme
$$
\pi_1(Z,c)\to \pi_1(X\otimes_k K,a)\times \pi_1(Y,b),
$$
\marginpar{266}%
et~\Ref{X.1.7} \'equivaut alors au cas particulier suivant:

\begin{corollaire}
\label{X.1.8}
Soient $X$ un sch\'ema propre et connexe sur un corps
alg\'ebriquement clos~$k$, $k'$ une extension alg\'ebriquement
close de~$k$, $a'$ un point g\'eom\'etrique de~$X\otimes_k k'$ et
$a$ son image dans~$X$. Alors l'homomorphisme canonique
$\pi_1(X\otimes_k k',a')\to \pi_1(X,a)$ est un \emph{isomorphisme.}
\end{corollaire}

Le fait que cet homomorphisme soit surjectif \'equivaut \`a dire
que si $X'$ est un rev\^etement \'etale connexe de~$X$, alors
$X'\otimes_k k'$ est \'egalement connexe, et r\'esulte
aussit\^ot du fait que $k$ est alg\'ebriquement clos; c'est aussi
un cas particulier
\ifthenelse{\boolean{orig}}
{IX.\Ref{IX.3.4}.}
{de~IX~\Ref{IX.3.4}.}
L'hypoth\`ese de propret\'e sur~$X$ n'a pas encore servi. Ceci
dit, dire que l'homomorphisme envisag\'e est injectif signifie aussi
ceci: \emph{tout rev\^etement \'etale de~$X\otimes_k k'$ est
isomorphe \`a l'image inverse d'un rev\^etement \'etale de~$X$.}
Il est essentiellement sorital qu'on peut trouver une
sous-$k$-alg\`ebre $A$ de
\ifthenelse{\boolean{orig}}
{$K$}
{$k'$},
de type fini sur~$k$, et un
rev\^etement \'etale de~$X\otimes_k A$ dont l'image inverse sur~$X\otimes_k k'$ soit isomorphe au rev\^etement donn\'e. Soit donc
$Y=\Spec(A)$, qui est un $k$-sch\'ema int\`egre de type fini, donc
a des points rationnels sur~$k$. Appliquons alors \Ref{X.1.7} au
groupe fondamental de~$X\times Y$ en un point $(a,b)$ rationnel
sur~$k$: on trouve que tout rev\^etement \'etale connexe de
$X\times Y$ est isomorphe \`a un quotient d'un rev\^etement
$X'\times Y'$, o\`u $X'$, $Y'$ sont des rev\^etements galoisiens
\'etales de~$X$ et $Y$ de groupes~$G$, $G'$ par un sous-groupe $H$
de~$G\times G'$. Cela implique que l'image inverse de ce
rev\^etement de~$X\times Y$ sur~$X\times Y'$ est isomorphe \`a un
rev\^etement de la forme $X_1'\times Y'$, o\`u $X_1'$ est un
rev\^etement \'etale de~$X$. Si donc $L$ est le corps des
fonctions de~$Y$, \'egal au corps des fractions de~$A$ dans~$k'$, le
rev\^etement \'etale de~$X\otimes_k L$ induit par le
rev\^etement donn\'e de~$X\times_k Y$ est tel qu'il existe une
extension finie s\'eparable $L'$ de~$L$, telle que l'image inverse
dudit rev\^etement sur~$X\otimes_k L'$ est isomorphe \`a
$X_1'\otimes_k L'$. Or $k'$ \'etant alg\'ebriquement clos, on peut
supposer que l'extension $L'$ de~$L$ est contenue dans~$k'$. Cela
prouve que le rev\^etement \'etale donn\'e de~$X\otimes_k k'$
est isomorphe \`a $X_1'\otimes_k k'$, cqfd.

La forme explicite signal\'ee en passant pour les rev\^etements
\'etales d'un produit
\marginpar{267}
$X\times_k Y$ implique aussit\^ot le r\'esultat suivant:

\begin{corollaire}
\label{X.1.9}
Soient $k$ un corps alg\'ebriquement clos, $X$ et $Y$ deux
pr\'esch\'emas localement noeth\'eriens sur~$k$, $Z=X\times_k
Y$ leur produit, $Z'$ un rev\^etement \'etale de~$Z$. Pour tout
point $y\in Y$ rationnel sur~$k$, soit $i_y \colon \Spec(k)\to Y$ le
morphisme canoniquement associ\'e, $j_y = \id_X\times_k i_y$ le
morphisme $X\to Z$ correspondant. Soit enfin $X_y'$ le rev\^etement
\'etale de~$X$ image inverse de~$Z'$ par~$j_y$. On suppose $Y$
connexe, et $X$ ou~$Y$ propre sur~$k$. Alors les rev\^etements
$X_y'$ de~$X$ sont tous isomorphes.
\end{corollaire}

De fa\c con imag\'ee, on peut dire qu'\emph{une famille de
rev\^etements \'etales de~$X$, param\'etr\'ee par un
pr\'esch\'ema connexe~$Y$, est constante si $X$ ou le
pr\'esch\'ema de param\`etres~$Y$ est propre sur~$k$.}

\begin{remarques}
\label{X.1.10}
Les corollaires~\Ref{X.1.7} \`a~\Ref{X.1.9} sont d\^us \`a
Lang-Serre~\cite{X.2} dans le cas des sch\'emas alg\'ebriques normaux.
(Leur travail a \'et\'e la motivation initiale pour la th\'eorie
du groupe fondamental d\'evelopp\'ee dans ce S\'eminaire). Comme
l'ont remarqu\'e ces auteurs, ces r\'esultats deviennent inexacts
lorsqu'on y abandonne l'hypoth\`ese de propret\'e, du moins en
caract\'eristique $p>0$. Prenant par exemple pour $X$ la droite
affine $X=\Spec(k[t])$, il n'est pas difficile de voir que les
rev\^etements de~$X$, param\'etr\'es par la droite affine
$Y=\Spec(k[s])$, d\'efinis par les \'equations
$$
x^p - x = st,
$$
sont \'etales et deux \`a deux non isomorphes. Cela met en
d\'efaut~\Ref{X.1.9} et a fortiori~\Ref{X.1.7}, et on voit de
m\^eme que si $s$ est consid\'er\'e comme un \'el\'ement
transcendant sur~$k$ dans une extension alg\'ebriquement close $K$
de~$k$, on trouve un rev\^etement \'etale $X'$ de~$X\otimes_kK$
qui ne provient pas d'un rev\^etement \'etale de~$X$.
\end{remarques}

\section[Application du th\'eor\`eme d'existence de faisceaux]{Application du th\'eor\`eme d'existence de faisceaux:
th\'eo\-r\`eme de semi-con\-ti\-nu\-i\-t\'e pour les groupes
fondamentaux des fibres d'un morphisme propre et s\'eparable}
\label{X.2}
\marginpar{268}

\begin{theoreme}
\label{X.2.1}
Soient $Y$ le spectre d'un anneau local noeth\'erien
\emph{complet}, de corps r\'esiduel~$\kres$, $X$ un $Y$-sch\'ema
propre, $X_0=X\otimes_A \kres$, $a_0$ un point g\'eom\'etrique de
$X_0$ et $a$ le point g\'eom\'etrique correspondant de~$X$. Alors
l'homomorphisme canonique $\pi_1(X_0,a_0)\to \pi_1(X,a)$ est un
\emph{isomorphisme.}
\end{theoreme}

Ce n'est qu'une traduction, dans le langage du groupe fondamental, du
r\'esultat rappel\'e dans IX~\Ref{IX.1.10}. C'est ici que le
th\'eor\`eme d'existence
\ifthenelse{\boolean{orig}}
{des}
{de}
faisceaux en g\'eom\'etrie
alg\'ebrico-formelle s'introduit de fa\c con essentielle dans la
th\'eorie du groupe fondamental.

Introduisons maintenant une cl\^oture alg\'ebrique $\bar{k}$ du
corps r\'esiduel $k$, et la fibre g\'eom\'etrique $\overline X_0
= X_0\otimes_k \bar{k}$. On a donc la suite exacte (IX~\Ref{IX.6.1})
$$
\ifthenelse{\boolean{orig}}
{e\to \pi_1(\overline X_0,\overline a_0)\to \pi_1(X_0,a_0)\to
\pi_1(k,\bar{k})\to e,}
{e\to \pi_1(\overline X_0,\overline a)\to \pi_1(X_0,a_0)\to
\pi_1(k,\bar{k})\to e,}
$$
d'autre part on a l'isomorphisme~\Ref{X.2.1} et l'isomorphisme
analogue, plus \'el\'ementaire,
$$
\pi_1(k,\bar{k})\to \pi_1(Y,b),
$$
o\`u $b$ est l'image de~$a$ dans~$Y$. On trouve ainsi:

\begin{corollaire}
\label{X.2.2}
\ifthenelse{\boolean{orig}}
{Avec les notations pr\'ec\'edentes, supposons $\overline X_0$
connexe, et soient $\overline a_0$ un point g\'eom\'etrique de
$\overline X_0 = X_0\otimes_k \bar{k}$, $a_0$ son image dans~$X$,}
{Avec les notations pr\'ec\'edentes, supposons
$\overline X_0=X_0\otimes_k \bar{k}$
connexe, et soient $\overline a_0$ un point g\'eom\'etrique de
$\overline X_0$, $a_0$ son image dans~$X$,}
$b_0$ son image dans~$Y$. Alors la suite d'homomorphismes canoniques
suivante
$$
\ifthenelse{\boolean{orig}}
{e\to \pi_1(\overline X_0,\overline a_0 )\to \pi_1(X,a_0)\to
\pi_1(Y,b_0)\to e}
{e\to \pi_1(\overline X_0,\overline a)\to \pi_1(X,a_0)\to
\pi_1(Y,b_0)\to e}
$$
est exacte.
\end{corollaire}

On comparera cette suite \`a la suite exacte~\Ref{X.1.4}, mais on
notera que~a)~on n'a pas eu \`a faire d'hypoth\`ese de platitude,
ou de s\'eparabilit\'e sur les fibres, pour $X\to Y$; b)~on a le
compl\'ement important que \emph{le morphisme $\pi_1(\overline
X_0,\overline a_0)\to \pi_1(X,a_0)$ est injectif.}

Ce dernier fait nous permettra de comparer le groupe fondamental des
autres
\marginpar{269}
fibres g\'eom\'etriques de~$X$ sur~$Y$ \`a celui de~$\overline
X_0$. Soit en effet $y_1$ un point quelconque de~$Y$, $X_1$ sa fibre
et $\overline X_1$ sa fibre g\'eom\'etrique, relativement \`a
une extension alg\'ebriquement close de~$\kres(y_1)$, $\overline a_1$
un point g\'eom\'etrique de~$\overline X_1$, $a_1$ son image dans
$X$ et $b_1$ son image dans~$Y$. Choisissons une \og classe de chemins\fg
de~$a_1$ \`a $a_0$, d'o\`u une classe de chemins de~$b_1$ \`a
$b_0$, d'o\`u un diagramme commutatif d'homomorphismes:
$$
\xymatrix@=.5cm{
&\pi_1(\overline X_1,\overline a_1) \ar[r] & \pi_1(X,a_1) \ar[r]\ar[d] & \pi_1(Y,b_1) \ar[r]\ar[d] & e \\
e \ar[r]& \pi_1(\overline X_0,\overline a_0) \ar[r] & \pi_1(X,a_0) \ar[r] &\pi_1(Y,b_0) \ar[r] & e,
}
$$
o\`u les deux fl\`eches verticales \'ecrites sont des
isomorphismes. Comme la deuxi\`eme ligne est exacte, on trouve donc
un homomorphisme canonique, que nous appellerons \emph{l'homomorphisme
de sp\'ecialisation pour le groupe fondamental} (ne d\'ependant
que de la classe de chemins choisis de~$a_1$ \`a~$a_0$, donc
\emph{d\'efini modulo automorphisme int\'erieur de}
$\pi_1(X,a_0)$):
$$
\pi_1(\overline X_1,\overline a_1)\to \pi_1(\overline X_0,\overline
a_0).
$$
Lorsque la premi\`ere ligne ci-dessus est \'egalement exacte, il
s'ensuit aussit\^ot que l'homo\-mor\-phisme de sp\'ecialisation
est surjectif. On trouve donc, compte tenu de~\Ref{X.1.4}:

\begin{corollaire}
\label{X.2.3}
Sous les conditions de~\Ref{X.2.1}, supposons de plus que le morphisme
$f\colon X\to Y$ soit \emph{s\'eparable} \eqref{X.1.1} et $\overline
X_0$ connexe (donc en vertu de~\Ref{X.1.2} on a
$f_*(\cal{O}_X)=\cal{O}_Y$). Alors pour toute fibre
g\'eom\'etrique $\overline X_1$ de~$X$ sur~$Y$, munie d'un point
g\'eom\'etrique $\overline a_1$, l'homomorphisme de
sp\'ecialisation
\index{homomorphisme de sp\'ecialisation pour le groupe fondamental|hyperpage}%
\index{sp\'ecialisation pour le groupe fondamental (homomorphisme de)|hyperpage}%
d\'efini ci-dessus est un homomorphisme \emph{surjectif.}
\end{corollaire}
C'est l\`a un r\'esultat de \emph{semi-continuit\'e} pour le
groupe fondamental, qui ne semble pas encore avoir d'analogue en
topologie alg\'ebrique. On peut d'ailleurs l'\'enoncer sous des
conditions plus g\'en\'erales:

\begin{corollaire}
\label{X.2.4}
Soient
\marginpar{270}
$f\colon X\to Y$ un morphisme propre \`a fibres
g\'eom\'etriquement connexes, avec $Y$ localement noeth\'erien,
$y_0$ et $y_1$ deux points de~$Y$ tels que $y_0\in\overline{\{y_1\}}$,
$\overline X_0$ et $\overline X_1$ les fibres g\'eom\'etrique de
$X$ correspondant \`a des extensions alg\'ebriquement closes
donn\'ees de~$\kres(y_0)$ et $\kres(y_1)$, $\overline a_0$
\resp $\overline a_1$ un point g\'eom\'etrique de~$\overline X_0$
\resp $\overline X_1$. Alors on peut d\'efinir de fa\c con
naturelle un homomorphisme de sp\'ecialisation:
$$
\pi_1(\overline X_1,\overline a_1)\to \pi_1(\overline X_0,\overline
a_0),
$$
d\'efini \`a automorphisme int\'erieur pr\`es, et c'est l\`a
un homomorphisme \emph{surjectif} si $f$ est un morphisme
s\'eparable \eqref{X.1.1}.
\end{corollaire}

En effet, il r\'esulte d'abord de~\Ref{X.1.8} que \Ref{X.2.4} est
essentiellement ind\'ependant des extensions alg\'ebriquement
closes choisies pour les corps r\'esiduels $\kres(y_0)$
et~$\kres(y_1)$. Cela nous permet de remplacer $Y$ par un
pr\'esch\'ema $Y'$ sur~$Y$ ayant un point $y_0'$ (\resp~$y_1'$)
au-dessus de~$y_0$ (\resp $y_1$). On prendra alors pour $Y'$ le
spectre du compl\'et\'e de l'anneau local de~$y_0$ dans~$Y$, et on
applique~\Ref{X.2.3}.

\begin{remarques}
\label{X.2.5}
La conclusion finale de~\Ref{X.2.4} sur la surjectivit\'e de
l'homomorphisme de sp\'e\-cia\-li\-sa\-tion, et a fortiori les
r\'esultats~\Ref{X.1.3} et~\Ref{X.1.4} dont elle est une
cons\'equence, devient inexacte si on ne suppose plus que $f\colon
X\to Y$ est un morphisme s\'eparable, m\^eme pour des sch\'emas
projectifs sur un corps alg\'ebriquement clos de
caract\'eristique~$0$. Nous en verrons plus loin des exemples, tant
dans le cas o\`u $f$ est plat mais o\`u $f$ admet une fibre non
s\'eparable ($X$ et $Y$ \'etant cependant lisses sur~$k$), que
dans le cas o\`u les fibres de~$f$ sont bien s\'eparables mais
o\`u $f$ n'est pas plat (par exemple $f\colon X\to Y$ \'etant un
morphisme birationnel de sch\'emas int\`egres normaux),
\cf XI~\Ref{XI.3}. Dans ces exemples, il peut arriver que le groupe
fondamental de la fibre g\'eom\'etrique g\'en\'erique soit
nul, mais non celui d'une fibre g\'eom\'etrique sp\'eciale
convenable. D'autre part, m\^eme si $f\colon X\to Y$ est un
morphisme propre s\'eparable comme dans~\Ref{X.2.4}, il arrive
couramment que le morphisme de sp\'ecialisation ne soit
\marginpar{271}
pas un isomorphisme. Ainsi, il est facile de donner des exemples
o\`u $\overline X_1$ est une courbe elliptique non singuli\`ere,
(donc son groupe fondamental est commutatif, et sa composante
$\ell$-primaire pour un nombre premier $\ell$ premier \`a la
caract\'eristique est isomorphe \`a $\ZZ_\ell^2$, \cf \Ref{XI}),
tandis que $\overline X_0$ est form\'e, soit de deux courbes
rationnelles non singuli\`eres se coupant en deux points, soit de
deux courbes rationnelles tangentes en un point, soit enfin d'une
courbe rationnelle ayant un point singulier qui est un point de
rebroussement. (Pour la classification compl\`ete des courbes
elliptiques d\'eg\'en\'er\'ees, voir les travaux r\'ecents
de Kodaira~\cite{X.1} et N\'eron). On voit alors que dans ces cas, le
groupe fondamental de~$\overline X_0$ est respectivement $\hat\ZZ$,
$e$, $e$, donc \og strictement plus petit\fg que celui de~$\overline
X_1$. Nous verrons cependant plus loin, lorsque $f$ est un morphisme
\emph{lisse}, une majoration du noyau de l'homomorphisme de
sp\'ecialisation, qui implique en particulier que si $\kres(y_0)$ est
de \emph{caract\'eristique} 0, l'homomorphisme de sp\'ecialisation
est un isomorphisme. Mais m\^eme pour un morphisme lisse, si la
caract\'eristique de~$\kres(y_0)$ est $>0$, il peut arriver que
l'homomorphisme de sp\'ecialisation ne soit pas un isomorphisme,
comme on voit par exemple dans le cas o\`u $X$ est un sch\'ema
ab\'elien sur~$Y$ (de dimension relative~$1$, si on veut),
\ifthenelse{\boolean{orig}}
{\cf XI~\Ref{XI.2}).}
{\cf XI~\Ref{XI.2}.}
Une th\'eorie satisfaisante de la
sp\'ecialisation du groupe fondamental doit tenir compte de la
\og composante continue\fg du \og vrai\fg groupe fondamental, correspondant
\`a la classification des rev\^etements principaux de groupe
structural des groupes infinit\'esimaux; moyennant quoi on serait en
droit \`a s'attendre que les \og vrais\fg groupes fondamentaux des
fibres g\'eom\'etriques d'un morphisme lisse et propre $f\colon
X\to Y$ forment un joli syst\`eme local sur~$X$, limite projective
de sch\'emas en groupes finis et plats sur~$X$\kern1pt\footnote{Cette
conjecture extr\^emement s\'eduisante est malheureusement mise en
d\'efaut par un exemple in\'edit de M\ptbl Artin, d\'ej\`a lorsque
les fibres de~$f$ sont des courbes alg\'ebriques de genre donn\'e
$g\ge 2$.}. Nous reviendrons ult\'erieurement sur ce point de vue,
notre objet pr\'esent \'etant au contraire de pousser aussi loin
que possible les ph\'enom\`enes communs \`a la th\'eorie
topologique et la th\'eorie sch\'ematique du groupe fondamental.
\end{remarques}

Soit maintenant $X_0$ une courbe propre, lisse et connexe de genre $g$
sur un corps al\-g\'e\-bri\-quement clos~$k$. Si $k$ est de
caract\'eristique z\'ero, son groupe fondamental peut se
d\'eterminer par voie transcendante de la fa\c con suivante. On
sait que $X_0$ provient par extension de la base d'une courbe
d\'efinie
\marginpar{272}
sur une extension alg\'ebriquement close de degr\'e de
transcendance fini du corps premier $\QQ$, et compte tenu
de~\Ref{X.1.8}, on peut supposer que $k$ est lui-m\^eme de degr\'e
de transcendance fini sur~$\QQ$. On peut donc supposer que $k$ est un
sous-corps du corps $\CC$ des nombres complexes, et une nouvelle
application de~\Ref{X.1.8} nous permet de supposer que $k=\CC$. Il
n'est pas difficile alors de v\'erifier que le groupe fondamental de
$X$ est isomorphe au compactifi\'e du groupe fondamental de l'espace
topologique associ\'e $\tilde X$ (surface compacte orient\'ee de
genre~$g$), pour la topologie d\'efinie par les sous-groupes
d'indice fini\footnote{Cette d\'eduction \'etait explicit\'ee
dans un des expos\'es oraux qui n'ont pas \'et\'e
r\'edig\'es.}. Il est d'autre part classique que le groupe
fondamental topologique est engendr\'e par $2g$ g\'en\'erateurs
$s_i,t_i$ ($1\le i\le g$), soumis \`a une seule relation:
$$
(s_1t_1s_1^{-1}t_1^{-1})\cdots(s_gt_gs_g^{-1}t_g^{-1})=1.
$$
Donc le groupe fondamental de~$X$ admet $2g$ g\'en\'erateurs
\emph{topologiques} $s_i,t_i$ ($1\le i\le g$), li\'es par la seule
relation pr\'ec\'edente. Si maintenant la caract\'eristique de
$k$ est $p>0$, d\'esignons par $A$ l'anneau des vecteurs de Witt
construit avec~$k$, par $K$ une extension alg\'ebriquement close de
son corps des fractions. On a vu dans III~\Ref{III.7.4} qu'il
existe un sch\'ema $X$ sur~$Y=\Spec(A)$, propre et lisse sur~$Y$, se
r\'eduisant suivant $X_0$. Appliquons-lui~\Ref{X.2.3}, on trouve un
morphisme \emph{surjectif}
$$
\pi_1(X_1)\to \pi_1(X_0),
$$
o\`u $X_1=X\otimes_A K$. Il est imm\'ediat\footnote{\cf EGA
IV~12.2.} que $X_1$ est lisse sur~$K$, connexe \eqref{X.1.2}, de
dimension~$1$, et son genre est \'egal \`a~$g$ (d'apr\`es
l'invariance de la caract\'eristique d'Euler-Poincar\'e, \cf EGA
III~7). Comme $K$ est de caract\'eristique~$0$, on peut lui
appliquer le r\'esultat pr\'ec\'edent. On a ainsi prouv\'e par
\emph{voie transcendante}:

\begin{theoreme}
\label{X.2.6}
Soit $X_0$ une courbe alg\'ebrique lisse propre et connexe sur un
\marginpar{273}
corps al\-g\'e\-bri\-quement clos $k$, et soit $g$ son genre. Alors
$\pi_1(X_0)$ admet un syst\`eme de~$2g$ g\'en\'erateurs
topologiques, li\'es par la relation \'ecrite plus haut. Lorsque la
caract\'eristique de~$k$ est~$0$, $\pi_1(X_0)$ est m\^eme le
groupe de type galoisien libre pour les g\'en\'erateurs et la
relation qui pr\'ec\`edent.
\end{theoreme}

\begin{remarques}
\label{X.2.7}
Il n'existe pas \`a l'heure actuelle, \`a la connaissance du
r\'edacteur, de d\'e\-mon\-stra\-tion par voie purement
alg\'ebrique du r\'esultat pr\'ec\'edent, (sauf pour les genres
$0,1$). Pour commencer, on ne voit gu\`ere comment distinguer dans
$\pi_1(X_0)$ $2g$~\'el\'ements, dont on pourrait attendre ensuite
qu'ils forment un syst\`eme de g\'en\'erateurs
topologiques. \`A cet \'egard, la situation de la droite
rationnelle priv\'ee de~$n$ points, et l'\'etude des
rev\^etements d'icelle mod\'er\'ement ramifi\'es en ces
points, est plus sympathique, puisque la consid\'eration des groupes
de ramification en ces $n$ points fournit $n$ \'el\'ements du
groupe fondamental \`a \'etudier, dont on montre en effet qu'ils
engendrent topologiquement ce groupe fondamental, comme nous verrons
ult\'erieurement\footnote{\Cf Exp\ptbl \Ref{XII}. Encore ces
\'el\'ements ne sont-ils d\'etermin\'es vraiment que modulo
conjugaison, et il convient de faire un choix \emph{simultan\'e}
judicieux de ces \'el\'ements dans leurs classes.}. Mais m\^eme
dans ce cas particuli\`erement concret, il ne semble pas exister de
d\'emonstration purement alg\'ebrique. Une telle
d\'emonstration serait \'evidemment extr\^emement
int\'eressante. Le seul fait concernant le groupe fondamental d'une
courbe qu'on sache d\'emontrer par voie purement alg\'ebrique
(exception faite du th\'eor\`eme de finitude faible~\Ref{X.2.12}
ci-dessous, prouv\'e par voie alg\'ebrique par Lang-Serre~\cite{X.2}),
semble la d\'etermination du groupe fondamental rendu ab\'elien
via la jacobienne (signal\'ee dans IX~\Ref{IX.5.8} derni\`ere
ligne).
\end{remarques}

\subsection{}
\label{X.2.8}
La derni\`ere assertion~\Ref{X.2.6} n'est plus valable en
caract\'eristique $p>0$, comme on voit d\'ej\`a dans le cas des
courbes elliptiques. Comme nous l'avons d\'ej\`a signal\'e, nous
ne savons pas si le groupe fondamental de~$X_0$ est topologiquement de
pr\'esentation finie; cela semble tout \`a fait improbable.

\begin{theoreme}
\label{X.2.9}
Soient $k$ un corps alg\'ebriquement clos, et $X$ un sch\'ema
propre et connexe sur~$k$. Alors le groupe fondamental de~$X$ est
topologiquement de g\'en\'eration finie.
\end{theoreme}

Nous
\marginpar{274}
proc\'edons par r\'ecurrence sur~$n=\dim X$, l'assertion
\'etant triviale pour~$n\le 0$. Supposons donc~$n>0$, et le
th\'eor\`eme d\'emontr\'e pour les
dimensions~$n'<n$. D'apr\`es le lemme de Chow (EGA~II~5.6.2) il
existe un sch\'ema projectif $X'$ sur~$k$ et un morphisme surjectif
$X'\to X$. On peut \'evidemment supposer $X'$ r\'eduit, et en
passant au normalis\'e, normal. Gr\^ace \`a la th\'eorie de
la descente, il suffit de prouver que les groupes fondamentaux des
composantes connexes de~$X'$ sont topologiquement de
g\'en\'eration finie (IX~\Ref{IX.5.2}). Cela nous ram\`ene
donc au cas o\`u $X$ est \emph{projectif} et \emph{normal}. Si
alors $n=1$, il suffit d'appliquer~\Ref{X.2.6}. Si $n\ge 2$, on
consid\`ere une immersion projective $X\to \PP_k^r$, et une section
hyperplane $Y=X\cdot H$
\ifthenelse{\boolean{orig}}
{(muni}
{(munie}
de la structure r\'eduite induite), telle que
\ifthenelse{\boolean{orig}}
{$Y\ne X$}
{$Y\ne X$,}
\ie $H\not\supset X$. On aura alors $\dim Y<n$, et tenant compte de
l'hypoth\`ese de r\'ecurrence, il suffit de prouver que
$\pi_1(Y)\to \pi_1(X)$ est \emph{surjectif.} Or plus
g\'en\'eralement:

\begin{lemme}
\label{X.2.10}
Soient $X$ un pr\'esch\'ema propre sur un corps alg\'ebriquement
clos~$k$, $g\colon X\to \PP_k^r$ un morphisme. On suppose $X$
irr\'eductible et normal et $\dim g(X)\ge 2$. Soient $H$ un
hyperplan de~$\PP_k^r$ et $Y=X\times_{\PP_k^r} H$. Alors $Y$ est
connexe, et l'homomorphisme $\pi_1(Y)\to\pi_1(X)$ est surjectif.
\end{lemme}

Ces assertions r\'esultent en effet de la suivante:

\begin{corollaire}
\label{X.2.11}
Sous les conditions pr\'ec\'edentes, soient $X'$ un rev\^etement
\'etale connexe de~$X$, et $Y'=X'\times_X Y=X'\times_{\PP_k^r} H$ le
rev\^etement induit sur~$Y$. Alors $Y'$ est connexe.
\end{corollaire}

Comme $X$ est normal, $X'$ est normal, donc \'etant connexe, $X'$
est irr\'eductible; de plus son image dans $\PP_k^r$ est de
dimension~$\ge 2$. Un lemme bien connu d\^u \`a Zariski (et
appel\'e \og \emph{th\'eor\`eme de Bertini}\fg) implique donc que
si $H_1'$ est l'hyperplan g\'en\'erique dans~$\PP_k^r$,
d\'efini sur une extension $K$ de~$k$, alors $X'\times_{\PP^r} H_1$
est universellement irr\'eductible donc universellement connexe
sur~$K$. Le th\'eor\`eme de connexion de Zariski (EGA~III~4)
implique alors que pour \emph{tout} hyperplan $H$ (d\'efini sur une
extension quelconque de~$k$) $X'\times_{\PP^r} H$ est
g\'eom\'etriquement connexe. Cela ach\`eve la d\'emonstration
de~\Ref{X.2.11}, donc de~\Ref{X.2.9}.

\begin{corollaire}[Lang-Serre]
\label{X.2.12}
Sous
\marginpar{275}
les conditions de~\Ref{X.2.9}, pour tout
groupe fini~$G$, l'ensemble des classes, \`a isomorphisme pr\`es,
de rev\^etements principaux de~$X$ de groupe~$G$, est fini.
\end{corollaire}

\begin{remarque}
\label{X.2.13}
Sous les conditions de~\Ref{X.2.10} nous prouverons lorsque $\dim
g(X)\ge 3$ (du moins lorsque $g$ est une immersion et $X$
r\'egulier), un r\'esultat plus pr\'ecis, connu en
g\'eom\'etrie alg\'ebrique sous le nom de
\og \emph{th\'eor\`eme de Lefschetz}\fg: $\pi_1(Y)\to \pi_1(X)$
\emph{est un isomorphisme.}\footnote{\Cf le s\'eminaire SGA~2 (1962)
faisant suite \`a celui-ci%
\ifthenelse{\boolean{orig}}
{.}
{, th\'eor\`eme X~3.10.}}
Il y a dans les cas classiques des
\'enonc\'es analogues pour les groupes d'homologie et les groupes
d'homotopie sup\'erieure, qui t\^ot ou tard devront \^etre
englob\'es dans la g\'eom\'etrie alg\'ebrique
abstraite. M\^eme pour la cohomologie de Hodge $\H^p(X,\Omega^q)$,
il ne semble pas que la question ait encore \'et\'e
\'etudi\'ee; il n'est d'ailleurs gu\`ere probable que pour cette
derni\`ere, les th\'eor\`emes de Lefschetz subsistent tels quels
en caract\'eristique~$p>0$.
\end{remarque}

\ifthenelse{\boolean{orig}}
{}
{\begin{remarqueMR}
\label{X.2.14}
Soit $R$ un anneau de valuation discr\`ete complet, de
corps r\'esiduel $\kres$ alg\'ebriquement clos de caract\'eristique $p>0$, de
corps des fractions~$K$ et soit $Y$ une courbe propre et lisse,
connexe, de genre $g$ sur~$R$. On a donc un morphisme de
sp\'ecialisation surjectif $\sp : \pi_1(Y_{\overline{K}}) 
\to\pi_1(Y_{\kres})$.
On a d\'ej\`a not\'e que si $K$ est de caract\'eristique 0, $\sp$ n'est pas un
isomorphisme d\`es que $g\geq1$. Supposons~$K$ de caract\'eristique $p$,
de sorte que $R$ est isomorphe \`a l'anneau de s\'eries 
formelles~$\kres[[T]]$.
En \'egale caract\'eristique $p >0$, le groupe fondamental n'est pas
d\'etermin\'e par le genre~$g$ comme on le voit d\'ej\`a sur les courbes
elliptiques qui peuvent \^etre ordinaires ou supersinguli\`eres. Citons
le r\'esultat r\'ecent de A\ptbl Tamagawa (non encore publi\'e). Si~$G$ est un
groupe profini, on note $G^{\res}$, le groupe profini quotient de~$G$,
limite projective des quotients topologiques de~$G$ qui sont finis
r\'esolubles.

\begin{theoreme*}[A\ptbl Tamagawa]
Supposons $g\geq2$, que la fibre sp\'eciale $Y_{\kres}$
soit d\'efinissable sur un corps fini et que le morphisme
$\sp^{\res}:\pi_1(Y_{\overline{K}})^{\res}\to\pi_1(Y_{\kres})^{\res}$ d\'eduit de
$\sp$ par passage au quotient, soit bijectif. Alors la courbe $Y$ est
constante sur~$R$.
\end{theoreme*}

Notons que le groupe de Galois de~$\overline{K}/K$ est r\'esoluble. L'\'enonc\'e
pr\'ec\'edent peut encore se traduire de la fa\c{c}on suivante : supposons que
tout rev\^etement fini \'etale $X_K\to Y_K$ de la fibre g\'en\'erique $Y_K$,
galoisien de groupe de Galois r\'esoluble, ait potentiellement bonne
r\'eduction (c'est-\`a-dire s'\'etende en un rev\^etement fini \'etale de~$Y$
apr\`es remplacement \'eventuel de~$R$ par son normalis\'e dans une extension
finie de~$K$). Alors $Y$ est constante sur~$R$.
\end{remarqueMR}}

\section[Application du th\'eor\`eme de puret\'e]{Application du th\'eor\`eme de puret\'e:
th\'eor\`eme de continuit\'e pour les groupes fondamentaux des
fibres d'un morphisme propre et lisse}
\label{X.3}

Rappelons sans d\'emonstration\footnote{Pour une d\'emonstration, \cf SGA~2 X~\Ref{X.3.4}.} le

\begin{theoremedepurete}[Zariski-Nagata]
\label{X.3.1}
\index{theoreme de purete@th\'eor\`eme de puret\'e|hyperpage}%
Soit $f\colon X\to Y$ un morphisme
quasi-fini et dominant de pr\'esch\'emas int\`egres, avec $X$
normal, $Y$ r\'egulier localement noeth\'erien, et soit $Z$
l'ensemble des points de~$X$ o\`u $f$ n'est pas \'etale,
\ie o\`u~$f$ est ramifi\'e (cela revient au m\^eme,
\textup{I~\Ref{I.9.5}~(ii)}). Si $Z\ne X$, $Z$ est de codimension~$1$ dans
$X$ en tous ses points, \ie pour toute composante irr\'eductible
$Z'$ de~$Z$ de point g\'en\'erique~$z$, la dimension de Krull de
$\cal{O}_{X,z}$ est \'egale \`a~$1$.
\end{theoremedepurete}

Rappelons qu'un pr\'esch\'ema est dit \emph{normal}
\resp \emph{r\'egulier} si ses anneaux locaux sont normaux
\resp r\'eguliers, et que la relation $Z\ne X$ signifie aussi que
l'extension finie $R(Z)/R(X)$ (o\`u $R$ d\'esigne le corps des
fonctions rationnelles) est \emph{s\'eparable}. Se pla\c cant en
le point g\'en\'erique $z$ d'une composante $Z'$ de~$Z$, et
localisant en le point $y$ de~$Y$ en-dessous de~$z$, on trouve
\marginpar{276}
l'\'enonc\'e \'equivalent:

\begin{corollaire}
\label{X.3.2}
Soient $A$ un anneau local noeth\'erien r\'egulier, $A\to B$ un
homomorphisme local injectif tel que $B$ soit normal, localis\'e
d'une alg\`ebre de type fini sur~$A$, et \emph{quasi-fini} sur~$A$;
on suppose de plus que $\dim A (=\dim B)\ge 2$, et que pour tout
id\'eal premier $\goth{p}$ de~$B$ distinct de l'id\'eal maximal,
$B$ est \'etale sur~$A$ en~$\goth{p}$, \ie $B_\goth{p}$ est
\'etale sur~$A_\goth{q}$ (o\`u $\goth{q}=A\cap\goth{p}$). Alors
$B$ est \'etale sur~$A$.
\end{corollaire}

D'ailleurs, il n'est pas difficile de r\'eduire ce dernier
\'enonc\'e au cas o\`u $A$ est un anneau local \emph{complet},
donc o\`u $B$ est \emph{fini} sur~$A$. Zariski~\cite{X.5} donne une
d\'emonstration simple de ce r\'esultat, valable dans le cas
d'\'egales caract\'eristiques; le cas g\'en\'eral est d\^u
\`a Nagata~\cite{X.3}, qui s'appuie sur un r\'esultat d\'elicat de
Chow; ce dernier n'a \'et\'e v\'erifi\'e par aucun des
participants du S\'eminaire, et devrait faire l'objet d'un
expos\'e ult\'erieur. Signalons seulement ici la d\'emonstration
tr\`es simple dans le cas particulier o\`u $\dim A=2$, qui est
suffisant pour l'application la plus importante que nous en ferons
dans le pr\'esent num\'ero. Comme $B$ est normal, il est un
$B$-module de profondeur (ancienne terminologie: codimension
cohomologique) $\ge 2$, donc c'est un $A$-module de profondeur~$\ge
2$, et comme $A$ est r\'egulier de dimension $2$, il en r\'esulte
que $B$ est un \emph{module libre} sur~$A$\kern1pt\footnote{\Cf EGA $0_{\textup{IV}}$
17.3.4}. Il r\'esulte alors de I~\Ref{I.4.10} que l'ensemble
des id\'eaux premiers~$\goth{q}$ de~$A$ en lesquels $B$ est
ramifi\'e sur~$A$ est la partie de~$\Spec(A)$ d\'efinie par un
id\'eal principal (engendr\'e par le discriminant d'une base de
$B$ sur~$A$), donc est vide si elle est contenue dans le point
ferm\'e de~$\Spec(A)$, ce qui prouve~\Ref{X.3.2} lorsque $\dim A=2$.

Nous utiliserons surtout~\Ref{X.3.1} sous la forme \'equivalente:

\begin{corollaire}
\label{X.3.3}
Soient $X$ un pr\'esch\'ema localement noeth\'erien
r\'egulier, $U$ une partie ouverte de~$X$ com\-pl\'e\-men\-taire
d'une partie ferm\'ee $Z$ de~$X$ de codimension~$\ge 2$. Alors le
foncteur $X'\mto X'\times_X U$ de la cat\'egorie des
rev\^etements \'etales de~$X$ dans la cat\'egorie des
rev\^etements \'etales de~$U$ est une \'equivalence
\marginpar{277}
de cat\'egories; en particulier, si $a$ est un point
g\'eom\'etrique de~$U$, l'homomorphisme canonique $\pi_1(U,a)\to
\pi_1(X,a)$ est un isomorphisme.
\end{corollaire}

La derni\`ere assertion est \'evidemment cons\'equence de la
premi\`ere, et pour celle-ci on peut \'evidemment supposer que $X$
est connexe donc irr\'eductible. De la normalit\'e de~$X$
r\'esulte d\'ej\`a que le foncteur $X'\mto X'\times_X U$ de
la cat\'egorie des rev\^etements localement libres (pas
n\'ecessairement \'etales) de~$X'$ dans la cat\'egorie des
rev\^etements de~$U$ est pleinement fid\`ele, car le foncteur
$\cal{E}\mto\cal{E}\vert U$ de la cat\'egorie des Modules
localement libres sur~$X$ dans la cat\'egorie des Modules localement
libres sur~$U$ l'est. Il reste donc \`a prouver que pour tout
rev\^etement \emph{\'etale} $U'$ de~$U$, il existe un
rev\^etement \'etale $X'$ de~$X$ (n\'ecessairement unique
d'apr\`es ce qui pr\'ec\`ede), tel que $U'$ soit isomorphe \`a
$X'\times_X U$. On peut \'evidemment supposer $U'$ connexe, donc
irr\'eductible puisque ($U$~\'etant normal) $U'$ est normal. Soit
$K$ le corps des fonctions rationnelles sur~$X$, ou sur~$U$ (c'est
pareil), $K'$ celui de~$U'$, alors $U'$ s'identifie au normalis\'e
de~$U$ dans~$K'$ I~\Ref{I.10.3}. Soit $X'$ le normalis\'e de
$X$ dans~$K'$ (EGA~II~6.3), alors $X'=X\times_X U \cong U'$, d'autre
part $X'$ est normal, int\`egre, le morphisme structural $f\colon
X'\to X$ est \emph{fini} et dominant (car $X$ est normal et $K'/K$ est
une extension finie s\'eparable). Il est \'etale dans
$U'=f^{-1}(U)=X'=f^{-1}(Z)$, et comme $Z$ est de codimension $\ge 2$
dans $X$, $f^{-1}(Z)$ est de codimension $\ge 2$ dans~$X'$. On conclut
alors de~\Ref{X.3.1} que $X'$ est \'etale sur~$X$, ce qui ach\`eve
la d\'emonstration.

\medskip
Soit maintenant $f\colon X\to Y$ une application rationnelle d'un
pr\'esch\'ema localement noeth\'erien et r\'egulier $X$
dans un pr\'esch\'ema~$Y$, et supposons que $f$ soit d\'efini
dans un ouvert $U$ compl\'ementaire d'une partie ferm\'ee de
codimension~$\ge 2$. Alors on d\'eduit de~\Ref{X.3.3} un foncteur,
d\'efini \`a isomorphisme pr\`es, de la cat\'egorie des
rev\^etements \'etales de~$Y$ dans la cat\'egorie des
rev\^etements \'etales de~$X$, d'o\`u pour tout point
g\'eom\'etrique $a$ de~$U$, d'image $b$ dans~$Y$, un homomorphisme
canonique
$$
\pi_1(X,a)\to \pi_1(Y,b)
$$
(d\'eduit
\marginpar{278}
de l'homomorphisme canonique $\pi_1(U,a)\to \pi_1(Y,b)$ gr\^ace
\`a l'isomorphisme $\pi_1(U,a)\isomto \pi_1(X,a))$. Lorsque $f$ est
un morphisme dominant, $X$ et $Y$ \'etant int\`egres de corps $K$
et~$L$, de sorte que $K$ est une extension de~$L$, et que $Y$ est
normal, ces correspondances se pr\'ecisent en termes d'extensions de
corps en notant que pour toute extension finie $L'$ de~$L$, non
ramifi\'ee sur~$Y$, l'alg\`ebre $K'=L'\otimes_L K$ sur~$K$ est non
ramifi\'ee sur~$X$.

En particulier, ces r\'eflexions montrent que le groupe fondamental
des pr\'esch\'emas localement noeth\'eriens connexes
r\'eguliers, ponctu\'es par des points g\'eom\'etriques
localis\'es en codimension~$\le 1$, est un \emph{foncteur} lorsque
l'on prend comme morphismes dans cette cat\'egorie les applications
rationnelles dominantes d\'efinies dans des compl\'ementaires de
parties ferm\'ees de codimension~$\ge 2$. Se rappelant par exemple
qu'une application rationnelle d'un sch\'ema normal sur un corps $k$
dans un sch\'ema propre sur~$k$ est d\'efinie dans le
compl\'ementaire d'un ensemble de codimension~$\ge 2$, on trouve:

\begin{corollaire}[Invariance birationnelle du groupe fondamental]
\label{X.3.4}
\index{invariance birationnelle du groupe fondamental|hyperpage}%
\index{theoreme d'invariance birationnelle du groupe fondamental@th\'eor\`eme d'invariance birationnelle du groupe fondamental|hyperpage}%
Soient $k$ un corps, $X$ et $Y$ deux sch\'emas propres sur~$k$ et
r\'eguliers, $f\colon X\to Y$ une application birationnelle de~$X$
dans~$Y$, $\Omega$ une extension alg\'ebriquement close du corps des
fonctions $K$ de~$X$, permettant de d\'efinir le groupe fondamental
de~$X$ et le groupe fondamental de~$Y$. Ces derniers sont alors
canoniquement isomorphes.
\end{corollaire}

Cela signifie aussi que pour une extension finie $K'$ de~$K$, si elle
est non ramifi\'ee sur un \og mod\`ele\fg propre non singulier $X$
de~$K$, elle l'est sur tout autre mod\`ele propre non singulier.

\begin{remarque}
\label{X.3.5}
Pour d'autre applications du th\'eor\`eme de puret\'e, voir les
travaux
\ifthenelse{\boolean{orig}}
{de \textsc{Abhyankar}}
{d'\textsc{Abhyankar}}
expos\'es dans~\cite{X.4}, inspir\'es par les
r\'esultats de \textsc{Zariski} \cite[chap\ptbl VIII]{X.6}, d\'emontr\'es par voie
topologique. Ces derniers sont loin d'avoir \'et\'e assimil\'es
par la g\'eom\'etrie alg\'ebrique \og abstraite\fg et m\'eritent
de nouveaux efforts.
\end{remarque}

Nous
\marginpar{279}
aurons besoin de quelques faits \'el\'ementaires de la th\'eorie
de la ramification. Soient $V$ un anneau de valuation discr\`ete de
corps des fractions~$K$, corps r\'esiduel~$k$, $L$ une extension
galoisienne de~$K$ de groupe~$G$, $V'$ le normalis\'e de~$V$
dans~$L$, qui est un module libre de rang $n=[L{\colon }K]$ sur~$V$,
$\goth{m'}$ un id\'eal maximal de~$V'$, $G_d$ le sous-groupe de~$G$
form\'e des \'el\'ements laissant $\goth{m'}$ invariant, de
sorte que $G_d$ op\`ere dans l'extension r\'esiduelle
$k'=V'/\goth{m'}$ de~$k$, et $G_i$ le sous-groupe des \'el\'ements
de~$G_d$ op\'erant trivialement (rappelons que $G_d$ et $G_i$ sont
appel\'es respectivement sous-groupes de \emph{d\'ecomposition} et
d'\emph{inertie} de~$G$). On dit que $L$ est \emph{mod\'er\'ement
ramifi\'e} sur~$V$ si $n_i=[G_i{\colon }e]$ est d'ordre $p$ premier
\`a la caract\'eristique de~$k$ (condition toujours
v\'erifi\'ee si $k$ est de caract\'eristique~0). Il est bien
connu que $G_i$ se plonge alors canoniquement dans le groupe ${k'}^*$,
donc est isomorphe au groupe des racines
\ifthenelse{\boolean{orig}}
{$n_i$.\`emes}
{$n_i$-i\`emes} de l'unit\'e
dans~${k'}^*$, ce qui implique en particulier que $G_i$ \emph{est
cyclique}. Le cas type de cette situation est celui o\`u on pose
$L=K[t]/(t^n - u)$, $u$ \'etant une uniformisante de~$V$ et $n$ un
entier premier \`a~$p$: si $K$ contient les racines \niemes de
l'unit\'e, $L$ est une extension galoisienne totalement ramifi\'ee
de~$K$, de groupe de Galois $G=G_i$ isomorphe \`a~$\ZZ/n\ZZ$.

\begin{lemme}[%
\ifthenelse{\boolean{orig}}
{Lemme de \textsc{Abhyankar}}
{Lemme d'\textsc{Abhyankar}}]
\label{X.3.6}
\index{Abhyankar (lemme d')|hyperpage}%
\index{lemme d'Abhyankar|hyperpage}%
Soit $V$ un anneau de valuation discr\`ete de
corps des fractions~$K$, $L$ et $K'$ deux extensions galoisiennes
de~$K$ \emph{mod\'er\'ement ramifi\'ees} sur~$V$, $n$ et $m$ les
ordres des groupes d'inertie
\ifthenelse{\boolean{orig}}
{correspondants.}
{correspondants,}
$L'$ une extension
compos\'ee de~$L$ et $K'$ sur~$K$. Si $m$ est un multiple de~$n$,
alors $L'$ est non ramifi\'ee sur les localis\'es de la
cl\^oture normale $V'$ de~$V$ dans~$K'$.
\end{lemme}

Soient en effet $W'$ le normalis\'e de~$V'$ dans~$L'$, $\goth{m'}$
un id\'eal maximal de~$V'$, $\goth{n'}$ un id\'eal maximal de~$W'$
au-dessus de
\ifthenelse{\boolean{orig}}
{$\goth{n'}$,}
{$\goth{m'}$,}
$\goth{n}$ l'id\'eal maximal qu'il induit sur le normalis\'e $W$
de~$V$ dans~$L$, $G,\ H,\ M$ les groupes de Galois de~$L,\ K',\ L'$
sur~$K$, et $G_i,\ H_i,\ L_i'$ les groupes d'inertie correspondant aux
id\'eaux maximaux choisis. Alors $M$ se plonge dans le produit
$G\times H$ et $M_i$ dans le produit $G_i\times H_i$, de telle fa\c con que les projections $M\to G$ et $M\to H$, $M_i\to G_i$
\marginpar{280}
et $M_i\to H_i$ soient surjectives (sorite du corps
interm\'ediaire). Il en r\'esulte d\'ej\`a, puisque $G_i$ et
$H_i$ sont par hypoth\`ese cycliques d'ordres $m$ et $n$ premiers
\`a~$p$, que $M_i$ est d'ordre premier \`a~$p$, donc cyclique, et
comme $m$ est multiple de~$n$ donc les \'el\'ements de~$G_i\times
H_i$ sont de puissance
\ifthenelse{\boolean{orig}}
{$m$.\`eme}
{$m$-i\`eme} nulle, $M_i$ est d'ordre
divisant~$m$, donc d'ordre \'egal \`a $m$ puisque $M_i\to H_i$ est
surjectif. Ce dernier homomorphisme est donc \'egalement
injectif. Or son noyau est le groupe d'inertie de~$\goth{n'}$
au-dessus de~$\goth{m'}$, ce qui prouve que $L'$ est non ramifi\'e
sur~$K'$ en~$\goth{n'}$. D'o\`u le lemme.

Pla\c cons-nous maintenant sous les conditions de~\Ref{X.2.4},
o\`u on a un homomorphisme de sp\'ecialisation
$$
\pi_1(\overline X_1,\overline a_1)\to\pi_1(\overline X_0,\overline
a_0)
$$
qui est \emph{surjectif}, relativement \`a un morphisme propre et
s\'eparable $f\colon X\to Y$. Nous voulons pr\'eciser le noyau de
cet homomorphisme. Proc\'edant comme dans la d\'emonstration
de~\Ref{X.2.4}, on voit que dans cette question, on peut toujours
supposer que $Y$ est le spectre d'un \emph{anneau de valuation
discr\`ete~$V$, complet et \`a corps r\'esiduel
alg\'ebriquement clos} (car on peut toujours trouver un tel anneau
et un morphisme de son spectre $Y'$ dans $Y$ dont l'image soit
$\{y_0,y_1\}$). Alors on a $X_0=\overline X_0$, $\kres(y_0)=k=$ corps
r\'esiduel de~$V$, $\kres(y_1)=K=$ corps des fractions de~$V$. Soit
$K_s$ la cl\^oture s\'eparable de~$K$, $\overline K$ sa
cl\^oture alg\'ebrique, et pour tout sous-anneau $W$ de~$\overline
K$ contenant~$V$, posons $X_W=X\otimes_V W$. En particulier on a
$$
X_V=\ifthenelse{\boolean{orig}}{V}{X}, \qquad X_K=X_1, \qquad X_{\overline K}=\overline X_1.
$$
D'ailleurs le morphisme canonique $\overline X_1=X_{\overline K}\to
X_{K_s}$ induit un isomorphisme sur les groupes fondamentaux
(IX~\Ref{IX.4.11}) de sorte que, compte tenu de
l'isomorphisme~\Ref{X.2.1} $\pi_1(X_0)\to\pi_1(X)$, on est ramen\'e
\`a \'etudier l'homomorphisme surjectif
$$
\pi_1(X_{K_s})\to\pi_1(X)
$$
associ\'e
\marginpar{281}
au morphisme canonique $X_{K_s}\to X$. La d\'etermination du noyau
de ce dernier revient \`a la solution du probl\`eme suivant:
\emph{on a un rev\^etement principal connexe $Z_{K_s}$ de
$X_{K_s}$, de groupe~$G$} (donc associ\'e \`a un homomorphisme de
$\pi_1(X_{K_s})$ dans~$G$), \emph{d\'eterminer sous quelles
conditions il est isomorphe \`a l'image r\'eciproque d'un
rev\^etement principal $Z$ de~$X$ de groupe~$G$.}

Notons d'abord que $K_s$ est r\'eunion filtrante croissante de ses
sous-extensions finies $K'$ sur~$K$, et que par suite, $Z_{K_s}$
est isomorphe \`a l'image inverse d'un rev\^etement principal
$Z_{K'}$ de~$X_{K'}$ pour un $K'$ convenable (on fera attention
cependant que pour $K'$ fix\'e, $Z_{K'}$ n'est pas d\'etermin\'e
de fa\c con unique). Dire que $Z_{K_s}$ est isomorphe \`a l'image
inverse d'un rev\^etement principal $Z$ de~$X$ signifie qu'il existe
une sous-extension finie $K''\supset K'$ de~$K_s$ telle que
$Z_{K''}=Z_{K'}\otimes_{K'} K''$ soit isomorphe \`a $Z\otimes_V
K''$. D\'esignons maintenant pour une sous-extension finie $K'$
de~$K_s$, par $V'$ le normalis\'e de~$V$ dans~$K'$, qui est un
anneau de valuation discr\`ete, complet, de corps
r\'esiduel~$k$. Donc le morphisme canonique $X_{V'}\to X_V$ induit
un isomorphisme pour les fibres au-dessus des points ferm\'es de
$Y=\Spec(V)$ et $Y'=\Spec(V')$ et il r\'esulte alors de~\Ref{X.2.1}
appliqu\'e \`a $X_V$ et $X_{V'}$ que l'homomorphisme induit pour
les groupes fondamentaux $\pi_1(X_{V'})\to \pi_1(X_V)$ est un
isomorphisme, ou encore que tout rev\^etement principal de~$X_{V'}$
est l'image inverse d'un rev\^etement principal de~$X_V$
d\'etermin\'e \`a isomorphisme pr\`es. Cela implique donc le

\begin{lemme}
\label{X.3.7}
Soit $Z_{K'}$ un rev\^etement principal connexe de~$X_{K'}$ de
groupe~$G$, $Z_{K}$ son image inverse sur~$X_{K_s}$\quoi. Pour que ce
dernier soit isomorphe \`a l'image inverse d'un rev\^etement
principal $Z$ de~$X$, il faut et il suffit qu'il existe une extension
finie $K''\supset K'$ de~$K$ dans $K_s$ telle que le rev\^etement
principal $Z_{K''}$ de~$X_{K''}$ soit induit par un rev\^etement
principal de~$X_{V''}$.
\end{lemme}

Supposons en particulier que les $X_{V''}$ soient normaux (il suffit
par exemple pour cela que $X_0$ soit normal, et a fortiori que $X_0$
soit simple, \cf I~\Ref{I.9.1}).
\marginpar{282}
Comme ils sont connexes, ils sont donc irr\'eductibles. Soient $L$ le
corps des fonctions rationnelles pour $X$ et $X_{K}$, $L'$ celui pour
$X_{V'}$ et $X_{K'}$, $L''$ celui pour $X_{V''}$ et~$X_{K''}$\quoi. Alors
sous les conditions de~\Ref{X.3.7}, $Z_{K'}$ d\'efinit une extension
finie s\'eparable $R'$ de~$L'$, et $Z_{K''}$ d\'efinit l'extension
$R''=R'\otimes_{L'} L'' = R'\otimes_{K'} K''$. La condition
envisag\'ee dans~\Ref{X.3.7} signifie donc aussi qu'il existe une
extension finie s\'eparable $K''$ de~$K'$ telle que
$R''=R'\otimes_{K'} K''$ soit \emph{non ramifi\'ee} au-dessus du
sch\'ema normal $X_{V''}$ de corps $L''=L'\otimes_{K'} K''$, et non
seulement au-dessus de la partie ouverte $X_{K''}$ de~$X_{V''}$.

Nous supposons dor\'enavant que $f\colon X\to Y$ est un morphisme
\emph{lisse,} donc les morphismes $X_{V'}\to \Spec(V')$ sont lisses,
donc les sch\'emas $X_{V'}$ sont \emph{r\'eguliers.} Noter que la
fibre du point ferm\'e de~$\Spec(V')$ dans $X_{V'}$ est
irr\'eductible et de codimension~$1$. Soit $\goth{o'}$ son anneau
local, qui est donc un anneau de valuation discr\`ete de corps
$L'$, de corps r\'esiduel isomorphe au corps des fonctions
rationnelles de~$X_0$, donc ayant m\^eme caract\'eristique que
$k$. D\'efinissons de m\^eme $\goth{o''}$ dans $L''$, qui est
\'evidemment le normalis\'e de
\ifthenelse{\boolean{orig}}
{$\goth{o''}$}
{$\goth{o'}$}
dans~$L''$. Il
r\'esulte alors du th\'eor\`eme de puret\'e~\Ref{X.3.1}
ou~\Ref{X.3.3} que pour que $R''$ soit non ramifi\'e sur~$X_{V''}$,
il faut et il suffit que $R''$ soit non ramifi\'e sur~$\goth{o''}$,
normalis\'e de~$\goth{o'}$ dans~$L''$.

Notons maintenant que si $u'$ est une uniformisante de~$V'$, c'est
aussi une uniformisante de~$\goth{o'}$. Si alors $n$ est un entier
premier \`a la caract\'eristique $p$ de~$k$, et si on prend
$K''=K'[t]/(t^n-u)$, alors
\ifthenelse{\boolean{orig}}
{$K'$}
{$K''$}
est une extension galoisienne finie
de~$K'$ et $L''$ est isomorphe \`a $L'[t]/(t^n-u')$, donc est
mod\'er\'ement ramifi\'e sur~$\goth{o'}$ et de groupe d'inertie
d'ordre~$n$. Supposons alors $G$ \emph{d'ordre premier \`a~$p$,} ce
qui implique que $R'$ est mod\'er\'ement ramifi\'e sur~$\goth{o'}$, et prenons pour $n$ un multiple premier \`a $p$ de
l'ordre du groupe d'inertie de~$R'$ sur~$\goth{o'}$ (par exemple
$n=[G\colon e]$). Appliquant le
\ifthenelse{\boolean{orig}}
{lemme de \textsc{Abhyankar}}
{lemme d'\textsc{Abhyankar}}~\Ref{X.3.6}, on
voit que la condition envisag\'ee dans~\Ref{X.3.7} est
v\'erifi\'ee.

Cela prouve le th\'eor\`eme suivant:

\begin{theoreme}
\label{X.3.8}
Soient $f\colon X\to Y$ un morphisme propre et lisse, \`a fibres
\marginpar{283}%
g\'eo\-m\'e\-tri\-que\-ment connexes, avec $Y$ localement noeth\'erien,
$y_0$ et $y_1$ deux points de~$Y$ tels que
\ifthenelse{\boolean{orig}}
{$y_0\in\overline{y_1}$}
{$y_0\in\overline{\{y_1\}}$},
$\overline{X}_0$ et $\overline{X}_1$ les fibres g\'eom\'etriques
correspondantes, consid\'erons l'homomorphisme de sp\'ecialisation
\eqref{X.2.4} $\pi_1(\overline X_1)\to\pi_1(\overline X_0)$. Cet
homomorphisme est surjectif, et tout homomorphisme continu de
$\pi_1(\overline X_1)$ dans un groupe fini $G$ d'ordre premier \`a
la caract\'eristique $p$ de~$\kres(y_0)$ provient d'un homomorphisme de
$\pi_1(\overline{X}_0)$ dans~$G$.
\end{theoreme}

En d'autres termes:

\begin{corollaire}
\label{X.3.9}
Si $\kres(y_0)$ est de caract\'eristique nulle, alors l'homomorphisme
de sp\'ecialisation est un isomorphisme. Si
\ifthenelse{\boolean{orig}}
{$p>0$}
{$\kres(y_0)$ est de caract\'eristique $p>0$},
alors le noyau de
l'homomorphisme de sp\'ecialisation est contenu dans l'intersection
des noyaux des homomorphismes continus de~$\pi_1(\overline X_1)$ dans
des groupes finis d'ordre premier \`a~$p$ (ou encore, le sous-groupe
invariant ferm\'e engendr\'e par un $p$-sous-groupe de Sylow du
groupe de type galoisien $\pi_1(\overline X_1)$); si donc
\ifthenelse{\boolean{orig}}
{$\pi_1(X_1)^{(p)}$}
{$\pi_1(\overline X_1)^{(p)}$}
d\'esigne le groupe quotient de
$\pi_1(\overline X_1)$ par le sous-groupe ferm\'e pr\'ec\'edent,
et si on d\'efinit de m\^eme $\pi_1(\overline X_0)^{(p)}$, alors
l'homomorphisme de sp\'ecialisation induit un \emph{isomorphisme}
$$
\pi_1(\overline X_1)^{(p)}\isomto \pi_1(\overline X_0)^{(p)}
$$
\end{corollaire}

On notera que la d\'emonstration de \Ref{X.3.8} est purement alg\'ebrique.
Proc\'edant comme dans \Ref{X.2.6}, on en conclut \emph{par voie
transcendante:}

\begin{corollaire}
\label{X.3.10}
Soit $X_0$ une courbe propre, lisse et connexe de genre $g$ sur un
corps al\-g\'e\-bri\-quement clos de caract\'eristique~$p$. Avec
la notation introduite dans~\Ref{X.3.9}, le groupe $\pi_1(X_0)^{(p)}$
est isomorphe \`a $\Gamma^{(p)}$, o\`u $\Gamma$ est le groupe de
type galoisien engendr\'e par des g\'en\'erateurs $s_i,~t_i$
($1\le i\le g$) li\'es par la relation
$$
(s_1t_1s_1^{-1}t_1^{-1})\cdots (s_gt_gs_g^{-1}t_g^{-1})=1
$$
\end{corollaire}

\begin{remarques}
\label{X.3.11}
\marginpar{284}
Dans le cas o\`u $\kres(y_0)$ est de caract\'eristique nulle, le
r\'esultat~\Ref{X.3.9} est bien connu par voie transcendante. On
notera que la d\'emonstration de~\Ref{X.3.10} fait appel au
th\'eor\`eme de puret\'e dans le cas d'in\'egales
caract\'eristiques, mais dans le cas d'anneaux de dimension 2
seulement, o\`u la d\'emonstration dudit th\'eor\`eme est
facile et a \'et\'e rappel\'ee dans le texte.
\end{remarques}


\chapterspace{-2}
\chapter{Exemples et compl\'ements}
\label{XI}
\marginpar{285}

\section{Espaces projectifs, vari\'et\'es unirationnelles}
\label{XI.1}

\begin{proposition}
\label{XI.1.1}
Soient $k$ un corps alg\'ebriquement clos, $X=\PP_k^r$ l'espace
projectif de dimension $r$ sur~$k$. Alors $X$ est \emph{simplement
connexe,} \ie $\pi_1(X)=0$.
\end{proposition}

Pour $r=0$, c'est trivial. Si $r=1$, il faut montrer que si $X'$ est
un rev\^etement \'etale connexe non vide de~$X=\PP_k^1$, alors
$X'\isomto X$. La formule du genre nous donne ici, si~$g$ et $g'$ sont
les genres de~$X$ et~$X'$:
$$
1-g' = d(1-g)\quoi,
$$
o\`u $d$ est le degr\'e de~$X'$ sur~$X$. Comme $g=0$, on aura donc
$1-g'=d$, ce qui exige $d=1$ puisque $g'\ge 0$, ce qui prouve
$X'\isomto X$. Lorsque $r\ge 2$, on proc\`ede par r\'ecurrence
sur~$r$, en supposant que $\PP^{r'}$ est simplement connexe pour
$r'<r$. Appliquant ceci \`a un hyperplan de~$\PP^r$ et utilisant
(X~\Ref{X.2.10}), il en r\'esulte bien que $\PP^r$ est simplement
connexe. Autre d\'emonstration: on aura
$\pi_1(\PP^1\times\cdots\times \PP^1) = \pi_1(\PP^1)\times\cdots
\times \pi_1(\PP^1)$ en vertu de (X~\Ref{X.1.7}), donc $(\PP^1)^r$
est simplement connexe puisque $\PP^1$ l'est, donc $\PP^r$ est
simplement connexe en vertu de l'invariance birationnelle du groupe
fondamental (X~\Ref{X.3.4}). Cette d\'emonstration montre plus
g\'en\'eralement:

\begin{corollaire}
\label{XI.1.2}
Soit $X$ un sch\'ema propre et normal sur un corps
alg\'ebriquement clos~$k$; si $X$ est une vari\'et\'e
rationnelle, \ie int\`egre et son corps des fonctions est une
extension transcendante pure de~$k$, alors $X$ est simplement connexe.
\end{corollaire}

Ce r\'esultat s'applique en particulier aux vari\'et\'es
grassmanniennes et plus
\marginpar{286}%
g\'en\'e\-ra\-le\-ment aux vari\'et\'es~$G/H$, o\`u $G$ est un
groupe lin\'eaire connexe sur~$k$ et $H$ un sous-groupe
alg\'ebrique contenant un sous-groupe de Borel
\ifthenelse{\boolean{orig}}
{de~$B$.}
{de~$G$.}

Rappelons qu'on appelle vari\'et\'e \emph{unirationnelle sur~$k$}
\index{unirationnelle (vari\'et\'e)|hyperpage}%
un sch\'ema propre et int\`egre sur~$k$ dont le corps des
fonctions $K$ est contenu dans une extension transcendante pure~$K'$
de~$k$, finie sur~$K$ (\ie ayant m\^eme degr\'e de transcendance
sur~$k$ que~$K$), \ie s'il existe une application rationnelle
dominante $f\colon\PP_k^r\to X$, avec $r=\dim X$. Si $X$ est normale,
on voit donc par les r\'eflexions pr\'ec\'edant X~\Ref{X.3.4}
que pour tout rev\^etement \'etale connexe $X'$ de~$X$, de
corps~$L/K$, l'alg\`ebre $L\otimes_K K'$ sur~$K'$ est non
ramifi\'ee sur le mod\`ele~$\PP^r$, donc compl\`etement
d\'ecompos\'ee en vertu de~\Ref{X.1.1} ce qui montre que $L$ est
$K$-isomorphe \`a une sous-extension de~$K'/K$. Cela prouve donc,
compte tenu de V~\Ref{V.8.2}:

\begin{corollaire}
\label{XI.1.3}
Le groupe fondamental d'une vari\'et\'e unirationnelle normale sur
un corps al\-g\'e\-bri\-quement clos, est fini.
\end{corollaire}

(N.B. On notera que dans la d\'efinition de \og unirationnelle\fg, on
n'avait pas besoin que $K'/K$ soit finie).

\ifthenelse{\boolean{orig}}
{\begin{remarques}
\label{XI.1.4}
Bien entendu, les r\'esultats de ce num\'ero sont bien
connus. D'autre part,
\ifthenelse{\boolean{orig}}{J.P\ptbl}{J-P\ptbl}%
Serre a montr\'e~\cite{XI.10} que lorsque $X$ est
une vari\'et\'e projective unirationnelle lisse sur un corps
alg\'ebriquement clos de \emph{caract\'eristique nulle,} $X$ est
simplement connexe. Sa d\'emonstration est transcendante en ce
qu'elle utilise le th\'eor\`eme de sym\'etrie de Hodge, et qu'on
ignore si ce r\'esultat s'\'etend \`a la
caract\'eristique~$p>0$. Il semble d'ailleurs qu'on ne connaisse pas
d'exemple de vari\'et\'e unirationnelle lisse sur~$k$ qui ne soit
d\'ej\`a rationnelle.
\end{remarques}}
{\begin{remarques}
\label{XI.1.4}
Bien entendu, les r\'esultats de ce num\'ero sont bien
connus. D'autre part,
\ifthenelse{\boolean{orig}}{J.P\ptbl}{J-P\ptbl}%
Serre a montr\'e~\cite{XI.10} que lorsque $X$ est
une vari\'et\'e projective unirationnelle lisse sur un corps
alg\'ebriquement clos de \emph{caract\'eristique nulle,} $X$ est
simplement connexe. Sa d\'emonstration est transcendante en ce
qu'elle utilise le th\'eor\`eme de sym\'etrie de Hodge.

\lcrochetbf Ajout\'e en 2003 (MR) Bornons-nous au cas d'un corps alg\'ebriquement clos $k$ de
caract\'eristique $p\geq 0$. En caract\'eristique $p>0$, la d\'efinition de
$k$-vari\'et\'e unirationnelle donn\'ee plus haut est celle de~$k$-vari\'et\'e
faiblement unirationnelle, par opposition \`a celle de~$k$-vari\'et\'e
fortement unirationnelle, o\`u l'on suppose de plus que $K'$ est une
extension s\'eparable de~$K$.

Dans le cas de la dimension $2$, il existe des surfaces projectives et
lisses, faiblement unirationnelles qui ont un groupe fondamental (fini
d'apr\`es 1.3) non trivial et donc qui ne sont pas rationnelles (T\ptbl Shioda, On unirationality of supersingular surfaces, Math.~Ann. \textbf{225}
(1977), p\ptbl155--159.). Par contre une surface fortement
unirationnelle est rationnelle: cela r\'esulte du crit\`ere de rationalit\'e de
Castelnuovo, \'etendu \`a la caract\'eristique $p>0$ par O\ptbl Zariski (\cf J-P\ptbl Serre, S\'em. Bourbaki \no146, 1957 et {\OE}uvres--Collected
papers, vol\ptbl 1, p\ptbl 467).

Sur le corps $\CC$ des nombres complexes, on conna\^it des exemples de vari\'et\'es projectives et
lisses de \hbox{dimension $\geq 3$} qui sont unirationnelles et non
rationnelles (\cf P\ptbl Deligne, Vari\'et\'es unirationnelles non
rationnelles [d'apr\`es M\ptbl Artin et D\ptbl Mumford],  S\'em.\ Bourbaki
\no402, 1971-72, Lecture Notes vol\ptbl 317). C'est le cas en particulier d'une
hypersurface cubique lisse dans l'espace projectif $\PP^4$.
Certains de ces exemples s'\'etendent en caract\'eristique $>0$ pour donner
des vari\'et\'es fortement unirationnelles non rationnelles (\cf J.P\ptbl Murre, Reduction of the proof of the non-rationality of a non
singular cubic threefold to a result of Mumford, Compositio Math. \textbf{27}
(1973), p\ptbl63--82).

Soit $V$ une $k$-vari\'et\'e projective normale, int\`egre. On dit que $V$ est
rationnellement connexe, s'il existe un $k$-sch\'ema de type fini int\`egre
$T$ et un $k$-morphisme $F : T\times\PP^1\to V$ tel que le morphisme :
\begin{equation}\tag*{(1)}
\begin{split}
T \times\PP^1\times\PP^1&\to V\times V\\
(t, u, u')&\mto(F(t,u), F(t,u'))
\end{split}
\end{equation}
soit dominant. En particulier, par deux points rationnels de~$V$,
assez g\'en\'eraux, passe une courbe rationnelle. La terminologie est
justifi\'ee par le fait que si $V$ est rationnellement connexe, on peut
joindre deux points rationnels par une cha\^ine finie de courbes
rationnelles. Si $k$ est de caract\'eristique $p>0$, on renforce la
d\'efinition pr\'ec\'edente en demandant que l'application (1) soit
g\'en\'eriquement lisse. On obtient alors la notion de vari\'et\'e
s\'eparablement rationnellement connexe. Ainsi les vari\'et\'es
unirationnelles sont rationnellement connexes et en caract\'eristique
$p>0$, les vari\'et\'es fortement unirationnelles sont s\'eparablement
rationnellement connexes. J\ptbl Koll\'ar a montr\'e que les vari\'et\'es
s\'eparablement rationnellement connexes ont un groupe fondamental
alg\'ebrique trivial (\cf O\ptbl Debarre, Vari\'et\'es rationnellement
connexes [d'apr\`es T\ptbl Graber, J\ptbl Harris, J\ptbl Starr et A.J\ptbl de~Jong],
S\'em. Bourbaki \no906, 2001-2002).\rcrochetbf
\end{remarques}}

\section{Vari\'et\'es ab\'eliennes}
\label{XI.2}
Soient $k$ un corps alg\'ebriquement clos, $A$ une vari\'et\'e
ab\'elienne sur~$k$, \ie un sch\'ema en groupes sur~$k$, propre
sur~$k$, lisse sur~$k$, et connexe, enfin $G$ un sch\'ema en groupes
commutatifs de type fini sur~$k$. D\'esignons par $\Ext(A,G)$ le
groupe des classes d'extensions commutatives de~$A$ par~$G$, par
$\H^1(A,G)$ le groupe des classes de fibr\'es principaux sur~$A$ de
groupe~$G$ (comparer \No \Ref{XI.4} plus bas), et consid\'erons
l'homomorphisme canonique
$$
\Ext(A,G)\to \H^1(A,G)\quoi.
$$
Un
\marginpar{287}
raisonnement de Serre \cite[chap\ptbl VII, th\ptbl 5]{XI.5} montre que c'est un
homomorphisme injectif, qui a pour image l'ensemble des
\og \'el\'ements primitifs\fg de~$\H^1(A,G)$, \ie des
\'el\'ements $\xi$ pour lesquels on~a:
$$
\pi^*(\xi) = \pr_1^*(\xi) + \pr_2^*(\xi)\quoi,
$$
o\`u $\pr_i$ sont les deux projections de~$A\times A$
sur~$A$, et $\pi\colon A\times A\to A$ la loi de composition de~$A$
(N.B. Serre n'\'enonce son th\'eor\`eme que pour $G$
lin\'eaire et connexe, et bien entendu lisse sur~$k$, mais en
simplifiant la premi\`ere partie de son raisonnement, on voit que
ces restrictions sont inutiles: il suffit de noter que tout morphisme
de~$A$ dans un sch\'ema en groupes $E$ de type fini sur~$k$, qui
transforme unit\'e en unit\'e, est un homomorphisme de groupes, et
d'appliquer ceci aux sections au-dessus de~$A$ d'une extension $E$ de
$A$ par~$G$).

Nous allons appliquer ce r\'esultat au cas o\`u $G$ est un groupe
fini s\'eparable sur~$k$, \ie un groupe fini ordinaire, suppos\'e
commutatif. Utilisant alors $\pi_1(A\times A)\isomto\allowbreak \pi_1(A)\times
\pi_1(A)$ (X~\Ref{X.1.7}) et interpr\'etant $\H^1(X,G)$ comme
$\Hom(\pi_1(X),G)$ pour tout sch\'ema alg\'ebrique~$X$, en
particulier pour $X=A$ ou $X=A\times A$, on voit que toute classe de
$\H^1(A,G)$ est primitive, donc on a un isomorphisme
$$
\Ext(A,G)\isomto \H^1(A,G)\quoi,
$$
en d'autres termes \emph{tout rev\^etement principal de~$A$ de
groupe structural commutatif~$G$, ponctu\'e au-dessus de l'origine
de~$A$, est muni de fa\c con unique d'une structure de groupe
alg\'ebrique admettant le point marqu\'e comme origine, et tel que
$A'\to A$ soit un homomorphisme de groupes alg\'ebriques.} En
particulier, si $A'$ est connexe, c'est \'egalement une
vari\'et\'e ab\'elienne, isog\`ene \`a~$A$.

D'autre part, comme le foncteur $X\mto \pi_1(X)$ des sch\'emas
alg\'ebriques ponctu\'es $X$ dans les groupes commute au produit
(IX~\Ref{IX.1.7}), il transforme un groupe dans la premi\`ere
cat\'egorie en un groupe dans la cat\'egorie des groupes, \ie en
un groupe \emph{commutatif}. Donc \emph{si $A$ est une vari\'et\'e
ab\'elienne, $\pi_1(A)$
\marginpar{288}
est un groupe commutatif.} Donc pour conna\^itre $\pi_1(A)$, il
suffit de conna\^itre le foncteur $G\mto
\H^1(A,G)=\Hom(\pi_1(A),G)$ pour $G$ variant dans les groupes finis
\emph{commutatifs}. Enfin, rappelons que pour tout entier $n>0$,
l'homomorphisme de multiplication par $n$ dans~$A$:
$$
A\lto{n} A
$$
est surjectif, donc \`a noyau fini, \ie c'est une isog\'enie, et
qu'il en r\'esulte que toute isog\'enie $A'\to A$ est quotient
d'une isog\'enie du type pr\'ec\'edent. De ceci, et de
raisonnements standard (\cf par exemple~\cite{XI.6}) on tire:

\begin{theoreme}[Serre-Lang]
\label{XI.2.1}
Soit $A$ une vari\'et\'e ab\'elienne sur un corps
alg\'ebriquement clos~$k$, et pour tout entier $n>0$ consid\'erons
le groupe fini ordinaire $K_n$ sous-jacent au noyau ${}_nA$ de la
multiplication par $n$ dans~$A$, enfin posons pour tout nombre
premier~$\ell$:
$$
T_\ell(A)=\varprojlim_r K_{\ell^r}
$$
et
$$
T_{.}(A)=\underset{\ell}{\prod}T_\ell(A)=\varprojlim_n K_n
$$
(o\`u pour $m$ multiple de~$n$, $m=ns$, on envoie $K_m$ dans $K_n$
par la multiplication par~$s$). Alors le groupe $\pi_1(A)$ est
canoniquement isomorphe \`a $T_{.}(A)$, donc pour tout nombre
premier $\ell$, la composante $\ell$-primaire de~$\pi_1(A)$ est
canoniquement isomorphe \`a $T_\ell(A)$.
\end{theoreme}

On notera que ces isomorphismes sont fonctoriels pour $A$ variable. Le
module $T_\ell(A)$ est appel\'e le \emph{module $\ell$-adique de
Tate} de la vari\'et\'e ab\'elienne~$A$. C'est un foncteur
additif en~$A$, en particulier il donne lieu \`a une
repr\'esentation de l'anneau $\Hom(A,A)$ des endomorphismes de~$A$
dans $T_\ell(A)$, appel\'ee \emph{repr\'esentation $\ell$-adique
de Weil}, et qui joue un r\^ole important dans la th\'eorie des
vari\'et\'es ab\'eliennes (\cf par exemple \hbox{\cite[chap~VII]{XI.4}}). Le
th\'eor\`eme~\Ref{XI.2.1} en donne une interpr\'etation en
termes de la repr\'esentation naturelle dans le \emph{groupe
d'homologie $\ell$-adique} de~$A$, $\H_1(A,\ZZ_\ell)=\pi_1(A)_\ell$,
ce
\marginpar{289}
qui est \'evidemment plus satisfaisant a priori, du point de vue
notamment de la formule de Lefschetz, \cite[Chap~V]{XI.4}. Rappelons ici les
r\'esultats de Weil sur la structure de~$T_\ell(A)$:
\begin{enumerate}
\item[a)] Si $n$ est premier \`a $\car(k)$, alors $K_n$ est un module
libre de rang \'egal \`a $2\dim A$ sur~$\ZZ/n\ZZ$, donc si $\ell$
est un nombre premier $\neq \car(k)$, $T_\ell(A)$ est un module libre
de rang \'egal \`a $2\dim A$ sur l'anneau $\ZZ_\ell$ des entiers
$\ell$-adiques;
\item[b)] Si $n$ est une puissance de~$\car(k)=p$, alors $K_n$ est un
module libre de rang $\nu\leq\dim A$ sur~$\ZZ/n\ZZ$, $\nu$
ind\'ependant de~$n$, donc $T_p(A)$ est un module libre de rang
$\nu\leq \dim A$ sur l'anneau $\ZZ_p$ des entiers $p$-adiques.
\end{enumerate}

Cela montre que dans la th\'eorie du groupe fondamental
d\'evelopp\'ee ici, le groupe fondamental d'une vari\'et\'e
ab\'elienne variable ne varie pas de fa\c con r\'eguli\`ere
avec le param\`etre, sa composante $p$-primaire pouvant diminuer
brusquement pour des valeurs du param\`etre $t$ correspondant \`a
une caract\'eristique r\'esiduelle~$p$; le cas le mieux connu de
ce ph\'enom\`ene est celui des courbes elliptiques. On notera
cependant que, quel que soit $n$ (premier ou non \`a la
caract\'eristique) le vrai noyau ${}_nA$ dans $A$ pour la
multiplication par $n$ est un sch\'ema en groupe fini sur~$k$ de
degr\'e $n^{2g}$, o\`u $g=\dim A$, qui sera non s\'eparable sur~$k$ si $n$ est un multiple de~$p=\car(k)$. D'ailleurs, lorsque $A$
varie dans une famille de vari\'et\'es ab\'eliennes, \ie si on
a un sch\'ema ab\'elien $A$ sur un sch\'ema de base $S$, on
montre plus g\'en\'eralement que ${}_n A$ est un sch\'ema en
groupes fini et plat sur~$S$, de degr\'e $n^{2g}$ sur~$S$,
c'est-\`a-dire \`a condition de tenir compte des parties
infinit\'esimales des noyaux~${}_n A$, ils se comportent de fa\c con r\'eguli\`ere quel que soit~$n$. Cela sugg\`ere que le
\og vrai\fg groupe fondamental d'une vari\'et\'e ab\'elienne $A$
est le pro-groupe alg\'ebrique (limite projective formelle de
groupes finis sur~$k$) $\varprojlim_n{{}_nA}$, o\`u
par \og vrai groupe fondamental\fg d'un sch\'ema alg\'ebrique~$X$, il
faut entendre: le pro-groupe qui classifie les rev\^etements
principaux de~$X$ de groupe structural un groupe fini quelconque $G$
sur~$k$ (pas n\'ecessairement s\'eparable sur~$k$). De cette
fa\c con par exemple, on r\'ecup\`ere par les
repr\'esentations de~$\Hom(A,A)$ dans la composante $p$-primaire du
vrai groupe fondamental de~$A$, le polyn\^ome caract\'eristique de
Weil d\'efini par ce dernier \`a l'aide des $\ell\neq p$, de fa\c con plus naturelle que la construction de Serre~\cite{XI.8}.

\section{C\^ones projetants, exemple de Zariski}
\label{XI.3}
\marginpar{290}

Soit toujours $k$ un corps alg\'ebriquement clos pour simplifier, et
soit $V$ un $k$-sch\'ema projectif connexe, sous-sch\'ema
ferm\'e de~$\PP^r_k$, qu'on pourra si on veut supposer non
singulier. Soient $Y=\hat C$ le c\^one projetant projectif de~$V$,
$y_0$ son sommet, $X=\hat C_V$ la fermeture projective habituelle du
fibr\'e vectoriel $C_V=\VV(\cal{O}_V(1))$ associ\'e \`a
$\cal{O}_V(1)$, enfin
$$
f\colon X\to Y
$$
le morphisme canonique, contractant la section nulle $X_0$ de~$C_V$
sur~$X$ en un point (EGA~II~8.6.4). Comme $X$ est un fibr\'e
localement trivial sur~$V$, de fibres $\PP^1$ donc de fibres
simplement connexes, le morphisme $p\colon X\to V$ induit en vertu
\ifthenelse{\boolean{orig}}
{de~(XI~\Ref{XI.4}.)}
{de~(X~\Ref{X.1.4})}
un isomorphisme:
$$
\pi_1(X)\isomto\pi_1(V)\quoi.
$$
Comme $p$ induit un isomorphisme $X_0\to V$, on en conclut qu'\emph{un
rev\^etement \'etale de~$X$ est compl\`etement d\'ecompos\'e
si et seulement si sa restriction \`a $X_0$ l'est}. Or pour tout
rev\^etement \'etale $Y'$ de~$Y$, l'image inverse $X'=X\times_Y
Y'$ est un rev\^etement \'etale de~$X$ compl\`etement
d\'ecompos\'e sur la fibre~$X_0$, donc trivial. Comme
l'homomorphisme $\pi_1(X)\to\pi_1(Y)$ est surjectif (IX~\Ref{IX.3.4}),
on en conclut que
$$
\pi_1(Y)=(e)
$$
en d'autres termes \emph{tout c\^one projetant projectif est
simplement connexe.} (N.B. en ca\-rac\-t\'e\-ris\-tique~$0$, le
m\^eme r\'esultat sera valable en prenant pour $Y$ le c\^one
projetant affine).

Supposons maintenant $V$ r\'eguli\`ere \ie lisse sur~$k$; alors
$X$ est r\'eguli\`ere, et pour une immersion projective convenable
de~$V$, on trouve alors un c\^one projetant $Y$ \emph{normal}. Si
$V$ n'est pas simplement connexe, donc $X$ non simplement connexe,
soit $X'$ un rev\^etement \'etale connexe non trivial
de~$X$. Comme les fibres de~$X$ en les points $y\in Y$ distincts de
$y_0$ sont r\'eduites \`a un point, on voit que la restriction de
$X'$ \`a ses fibres (en particulier \`a la
fibre
\marginpar{291}
g\'en\'erique) est triviale; cependant $X$ ne provient pas par
image inverse d'un rev\^etement \'etale de~$Y$, puisque $Y$ est
simplement connexe et que $X'$ serait compl\`etement
d\'ecompos\'e. Cela montre que (X~\Ref{X.1.3} et~\Ref{X.1.4})
deviennent faux si on remplace l'hypoth\`ese que $f$ est
s\'eparable par celle plus faible que ses fibres sont des
sch\'emas alg\'ebriques s\'eparables (ou m\^eme lisses) sur
les~$\kres(s)$. On notera de m\^eme que les groupes fondamentaux des
fibres g\'eom\'etriques $\overline X_y$ des $y\neq y_0$ sont
\'evidemment r\'eduits \`a $(e)$ puisque ces fibres sont
r\'eduites \`a un point, alors que $\pi_1(X_0)\neq e$, donc le
th\'eor\`eme de semi-continuit\'e (X~\Ref{X.2.4}) est
\'egalement en d\'efaut pour~$f$.

Indiquons enfin l'exemple, signal\'e par Zariski, mettant en
d\'efaut ces m\^emes th\'eor\`emes, lorsqu'on y remplace
l'hypoth\`ese que $f$ est s\'eparable par celle que $f$ est
plate. Soit $f\colon X\to Y$ un morphisme d'une surface non
singuli\`ere projective dans la droite rationnelle $Y=\PP^1$, tel
que
\ifthenelse{\boolean{orig}}
{$K=k(x)$}
{$K=k(X)$}
soit une extension \og r\'eguli\`ere\fg \ie primaire et
s\'eparable de~$k(f)$, (\ie la fibre g\'en\'erique
g\'eom\'etrique est connexe et s\'eparable), et telle que le
diviseur $(f)=X_0-X_\infty$ soit un multiple
\nieme
d'un diviseur (o\`u $n$ est un entier premier \`a la
caract\'eristique). Il est possible de construire de tels exemples
en toute caract\'eristique. Soit $X'$ le normalis\'e de~$X$ dans
$K(f^{1/n})$, o\`u $K=k(X)$ est le corps des fonctions de~$X$. Il
r\'esulte de l'hypoth\`ese sur~$(f)$ que $X'$ est \'etale
sur~$X$. Soit $Y'$ le normalis\'e de~$Y$ dans $k(t)(t^{1/n})$, il
est ramifi\'e sur~$Y$ en les points $t=0$ et $t=\infty$ exactement,
et la restriction $X'|f^{-1}(U)$ est isomorphe \`a l'image inverse
de~$Y'|U$. En particulier, la restriction de~$X'$ \`a la fibre
g\'en\'erique \emph{g\'eom\'etrique} de~$X$ sur~$Y$ se
d\'ecompose compl\`etement. Cependant, $X'$ n'est pas isomorphe
\`a l'image inverse d'un rev\^etement \'etale de~$Y$, car on
voit tout de suite que ce dernier serait n\'ecessairement~$Y'$, ce
qui est absurde puisque $Y'$ est ramifi\'e sur~$Y$\kern2pt\footnote{On peut
remarquer, du point de vue de la \og topologie \'etale\fg (SGA~4~VII),
que dans cet exemple $\R^1(f_{\et})_*(\ZZ/n\ZZ)$ est \og non
s\'epar\'e\fg sur~$S$.}.

Voici (d'apr\`es Serre) une fa\c con simple de r\'ealiser les
conditions de cet exemple, en s'inspirant de \cite[\No 20]{XI.5}: on prend
pour $n$ un nombre premier $\geq 5$, distinct de la
caract\'eristique, et on fait op\'erer $G=\ZZ/n\ZZ$ dans
$\underline{k}^4$%
\ifthenelse{\boolean{orig}}
{}
{{\renewcommand{\thefootnote}{*}\addtocounter{footnote}{-1}\kern1pt\footnote{\lcrochetbf Ajout\'e en 2003: $\underline{k}^4$ d\'esigne l'espace affine de dimension $4$ sur $k$.\rcrochetbf}}}
en multipliant les coordonn\'ees par quatre
caract\`eres distincts de~$G$ (ce qui est possible puisque $n\geq
5$). Alors $G$ op\`ere sur l'espace projectif $\PP^3_k$, et les
seuls points fixes sous $G$ sont les quatre points correspondants aux
axes de coordonn\'ees. La surface $X'$ d'\'equation
\ifthenelse{\boolean{orig}}
{$x^4+y^4+z^4+t^4=0$}
{$x^n+y^n+z^n+t^n=0$}
est lisse sur
\marginpar{292}
$k$ (crit\`ere jacobien), et ne contient aucun des points fixes,
donc $G$ \'etant d'ordre premier, op\`ere sur~$X$ \og sans points
fixes\fg \ie $X$ est un rev\^etement principal de~$X=X'/G$ de groupe
$G$. Soit $g=x/y$ dans $k(X')=K'$, c'est un g\'en\'erateur
kumm\'erien de~$K'$ sur~$K=k(X)$ si les caract\`eres choisis
\'etaient~$\chi^i$, $i=0,1,2,3$, avec $\chi$ un caract\`ere
primitif, soit $f$ sa puissance \nieme, qui est un
\'el\'ement de~$K$. On voit tout de suite que $K'$ est une
extension r\'eguli\`ere de~$k(g)$, ce qui r\'esulte du fait que
la courbe plane d'\'equation homog\`ene en $U,T,Z$:
$T^n+Z^n+(1+g^n)U^n=0$ est lisse sur~$k(g)$ (crit\`ere jacobien), et
qu'on sait que toute courbe plane est connexe. D'autre part on a
$k(f)=K\cap k(g)$, puisque le deuxi\`eme membre est une extension de
$k(f)$ contenue dans l'extension $k(g)$ de degr\'e premier, et
distincts de~$k(g)$ (puisque $g\not\in K$). Cela implique que $K$ est
une extension r\'eguli\`ere de~$k(f)$. Enfin le diviseur de~$f$
sur~$X$ est un multiple \nieme d'un diviseur, car son image
inverse sur~$X'$ est le diviseur de~$g^n$, donc un multiple
\nieme, et on peut redescendre parce que $X'$ est \'etale sur~$X$. On aurait fini si l'application rationnelle $f\colon X\to \PP^1$
\'etait un morphisme, c'est-\`a-dire si les diviseurs des
z\'eros et des p\^oles de~$f$ ne se rencontraient point. En fait,
on v\'erifie ais\'ement (en regardant encore sur~$X'$) que les
deux diviseurs en question sont les produits par $n$ de deux courbes
lisses sur~$k$, se coupant transversalement en un point~$a$. Rempla\c cant maintenant $X$ par le sch\'ema obtenu en faisant \'eclater~$a$, les conditions pr\'ec\'edentes ($\divisor(f)$ divisible par
$n$, et $k(X_1)=k(X)$ extension r\'eguli\`ere de~$k(X)$) restent
v\'erifi\'ees, mais de plus $f$ est un \emph{morphisme} $X_1\to
\PP^1$, donc on est sous les conditions voulues.

\section{La suite exacte de cohomologie}
\label{XI.4}

Soit $S$ un pr\'esch\'ema, de sorte que la cat\'egorie
$\Sch_{/S}$ des pr\'esch\'emas sur~$X$ est d\'etermin\'ee,
donc aussi la notion de groupe dans icelle, qu'on appellera aussi
\emph{pr\'esch\'ema en groupes sur~$S$}, ou simplement
\emph{$S$-groupe}. Pour simplifier l'exposition et fixer les
id\'ees, nous nous bornerons le plus souvent par la suite \`a des
groupes qui sont \emph{affines} et \emph{plats} sur~$S$\kern1pt\footnote{En
fait, pour ce qui va suivre, l'hypoth\`ese quasi-affine au lieu
d'affine suffirait, \cf note de bas de page~\eqref{note-296}, p\ptbl\pageref{note-296}.}, ce qui
suffira pour les applications que nous avons en vue. (Bien entendu, on
rencontre de nombreux cas o\`u ni l'une ni l'autre hypoth\`ese
n'est v\'erifi\'ee). Soit
\marginpar{293}
$G$ un tel $S$-groupe, et soit $P$ un pr\'esch\'ema sur~$S$ sur
lequel $G$ op\`ere \`a droite, ce qui implique en particulier un
morphisme
$$
\pi:P\times_SG\to P
$$
dans la cat\'egorie $\Sch_{/S}$, satisfaisant les axiomes bien
connus. On dit que $P$ est \emph{formellement principal homog\`ene
sous} $G$ si le morphisme
$$
P\times_S G\to P\times_S P
$$
de composantes $\pr_1$ et $\pi$ est un isomorphisme; il revient au
m\^eme de dire que pour tout objet $S'$ de~$\Sch_{/S'}$,
$P(S')=\Hom_S(S',P)$ consid\'er\'e comme ensemble \`a groupes
d'op\'erateurs $G(S')=\Hom_S(S',G)$, est vide ou principal
homog\`ene (\ie vide ou isomorphe \`a $G(S')$ sur lequel le
groupe $G(S')$ op\`ere par translations \`a droite). On dit que
$P$ est \emph{trivial} si $P$ est isomorphe \`a $G$, sur lequel $G$
op\`ere par translations \`a droite, ou ce qui revient au
m\^eme, si chacun des ensembles \`a op\'erateurs $P(S')$ sous
$G(S')$ est trivial. On v\'erifie, par exemple par le
proc\'ed\'e brevet\'e de passage au cas ensembliste, que $P$
\emph{est trivial si et seulement
\ifthenelse{\boolean{orig}}
{si il}
{s'il}
est formellement principal homog\`ene, et admet une section sur~$S$}
(ce dernier fait s'\'enon\c cant en termes cat\'egoriques en
disant que $P$ a une section sur l'objet final $e=S$ de~$\Sch_{/S}$,
\ie qu'il existe un morphisme de~$e$ dans $P$). Pour d\'efinir la
notion de fibr\'e principal homog\`ene $P$ sous $G$, plus forte
que celle de fibr\'e formellement principal homog\`ene, il faut
pr\'eciser d'abord dans $\Sch_{/S}$ un ensemble de morphismes qui
seront utilis\'es pour la \og descente\fg, et joueront le r\^ole de
\og morphismes de localisation\fg pour \og trivialiser\fg des fibr\'es. Le
choix le plus ad\'equat varie suivant le contexte, aucun ne
contenant tous les autres\footnote{\label{note-293}Voir \`a ce sujet
SGA~3~IV, notamment~\S 4.}. Ici, il sera commode d'adopter la
d\'efinition suivante:

\begin{definition}
\label{XI.4.1}
Soit $G$ un $S$-groupe. On appelle \emph{fibr\'e principal
homog\`ene} (\`a droite) sous~$G$, un $S$-pr\'esch\'ema $P$
\`a~$S$-groupe \`a droite~$G$, tel qu'il existe un recouvrement de~$S$ par des ouverts~$U_i$, et pour tout $i$ un morphisme de changement
de base $S'_i\to U_i$ fid\`element plat et quasi-compact, tel que
\ifthenelse{\boolean{orig}}
{$P'_i=P\times_S S'_i$ soit un pr\'esch\'ema \`a op\'erateurs
trivial sous $G'=G\times_S S'$.}
{$P'=P\times_S S'$ soit un pr\'esch\'ema \`a op\'erateurs
trivial sous $G'=G\times_S S'$, o\`u $S'$ est le
$S$-pr\'esch\'ema somme disjointe des~$S'_i$.}
\end{definition}

(On
\marginpar{294}
notera que le foncteur changement de base $X\mto X'=X\times_S
S'$
\'etant exact \`a gauche, transforme groupes en groupes, objets
\`a groupe d'op\'erateurs en objets \`a groupes
d'op\'erateurs). Notons que \Ref{XI.4.1} est \emph{stable par
changement de base}. Notons aussi:

\begin{proposition}
\label{XI.4.2}
Soient $G$ un $S$-groupe, plat et quasi-compact sur~$S$, $P$ un
$S$-pr\'esch\'ema o\`u $G$ op\`ere \`a droite. Conditions
\'equivalentes:
\begin{enumerate}
\item[(i)] $P$ est un fibr\'e principal homog\`ene
\ifthenelse{\boolean{orig}}
{sous~$G$.}
{sous~$G$;}
\item[(ii)] $P$ est formellement principal homog\`ene sous~$G$,
\index{formellement principal homog\`ene (pr\'esch\'ema)|hyperpage}%
et le morphisme structural $P\to S$ est fid\`element plat et
quasi-compact.
\end{enumerate}
\end{proposition}

Si $P$ est principal homog\`ene sous~$G$, alors avec les notations
de~\Ref{XI.4.1} $P'$ est fid\`element plat et quasi-compact sur~$S'$
(puisque $G'$ l'est, et $P'$ lui est $S'$-isomorphe), donc $P$ a les
m\^emes propri\'et\'es au-dessus de~$S$, (pour \og surjectif\fg et
\og quasi-compact\fg, \cf VIII~\Ref{VIII.3.1}, pour \og plat\fg c'est un
oubli dans les sorites de l'expos\'e~\Ref{VIII}). Inversement, si
(ii) est v\'erifi\'e, prenons le changement de base $S'=P$, qui
est bien fid\`element plat et quasi-compact sur~$S$; alors $P'$ sera
formellement principal homog\`ene sur~$S'$ puisque $P$ l'est sur~$S$
et que le foncteur changement de base est exact \`a gauche, d'autre
part $P'$ a une section sur~$S'$,
\ifthenelse{\boolean{orig}}
{savoir}
{\`a savoir}
la section diagonale, donc
c'est un fibr\'e principal trivial, ce qui ach\`eve la
d\'emonstration.

\begin{corollaire}
\label{XI.4.3}
Si $G$ est affine et plat sur~$S$, tout fibr\'e principal
homog\`ene $P$ sous $G$ est affine et plat sur~$S$.
\end{corollaire}

En effet, il le devient par extension fid\`element plate et
quasi-compacte de la base, et on applique (VIII~\Ref{VIII.5.6}).

L'utilit\'e de la d\'efinition~\Ref{XI.4.1} pour des $S$-groupes
\emph{plats} et \emph{affines} sur~$S$ tient \`a
(VIII~\Ref{VIII.2.1}), \ie au fait que les morphismes $S'\to S$
envisag\'es dans~\Ref{XI.4.1} sont des morphismes de descente
effective pour la cat\'egorie des pr\'esch\'emas affines sur
d'autres. Gr\^ace \`a ce fait, la v\'erification des faits
esquiss\'es ci-dessous se fait de fa\c con essentiellement
\og cat\'egorique\fg\footnote{\Cf \loccit dans note de bas de
page~\eqref{note-293}, p\ptbl\pageref{note-293}.}. Soit $E$ un $S$-pr\'esch\'ema sur lequel le
$S$-groupe $G$ op\`ere \`a gauche, et soit $P$ un fibr\'e
principal homog\`ene (\`a droite) sous~$G$, nous voulons
d\'efinir un fibr\'e associ\'e $E^{(P)}$, \og localement\fg
isomorphe \`a~$E$. Pour ceci, faisons op\'erer
\marginpar{295}
\`a droite $G$ dans $P\times_S E$ suivant la loi $(x,y)\mto
(xg,g^{-1}y)$, qui d\'ecrit de telles op\'erations dans le
contexte ensembliste, et s'\'etend aux cat\'egories par le
proc\'ed\'e brevet\'e. On posera, sous-r\'eserve d'existence:
$$
E^{(P)}=(P\times_S E)/G
$$
moyennant quoi on constate que $P\times_S E$ sera un pr\'esch\'ema
au-dessus de~$T=E^{(P)}$, \`a groupes d'op\'erateurs \`a droite
$G_T=G\times_S T$; pour \^etre \`a l'aise, on aimerait que de plus
$P\times_S E$ soit un fibr\'e principal homog\`ene sur~$T$ de
groupe $G_T$. Pour v\'erifier l'existence de~$E^{(P)}$ et la
propri\'et\'e pr\'ec\'edente, reprenons le $S'$ de la
d\'efinition~\Ref{XI.4.1} et regardons la situation image inverse
sur~$S'$ de la situation initiale. Du fait que $P'$ est
\ifthenelse{\boolean{orig}}
{trivial}
{trivial,}
\ie isomorphe \`a $G'_d$, on voit tout de suite que $E^{\prime(P')}$
existe, et a la propri\'et\'e voulue d'exactitude. En fait,
$E'\times_{S'} P'$ est $G'$-isomorphe au produit $E'\times_{S'} G'$,
donc $E^{\prime(P')}$ est isomorphe \`a~$E'$. De plus, la formation du
\og fibr\'e associ\'e\fg dans le cas d'un espace \`a op\'erateurs
trivial commute \`a toute extension de la base, et prenant en
occurrence les diverses extensions de la base
$\xymatrix@C=.5cm{S''\ar@<2pt>[r]\ar@<-2pt>[r]&S'}$
et
$\xymatrix@C=.5cm{S'''\ar@<4pt>[r]\ar@<0pt>[r]\ar@<-4pt>[r]&S'}$, o\`u
$S''$ et $S'''$ sont les produits fibr\'es double et triple de~$S'$
sur~$S$, on constate que $E^{\prime(P')}$ \emph{est muni d'une donn\'ee
de descente relativement au morphisme $S'\to S$, et que $E^{(P)}$
existe avec les propri\'et\'es requises si et seulement si cette
donn\'ee de descente est effective}; bien entendu $E^{(P)}$ n'est
autre alors que l'objet descendu. (Utiliser le fait que $S'\to S$ est
un morphisme de descente dans la cat\'egorie des
$S$-pr\'esch\'emas, \cf VIII~\Ref{VIII.5.2}). Il s'ensuit~que
\emph{le fibr\'e associ\'e existe si $E$ est affine sur~$S$}. Nous
appliquerons cette construction au cas o\`u on~a un homomorphisme de~$S$-groupes $G\to H$, et qu'on prend pour $E$ le $S$-pr\'esch\'ema~$H$ muni des op\'erations de~$G$ sur~$H$ \`a gauche r\'esultant
du morphisme donn\'e; comme $H$~op\`ere \`a droite sur
lui-m\^eme de fa\c con \`a commuter aux op\'erations de~$G$
sur~$H$, et que (sous r\'eserve d'existence au-dessus de~$S$) la
formation du fibr\'e associ\'e commute \`a l'extension de la
base, on constate ais\'ement que $H$ va op\'erer \`a droite sur~$P^{(H)}$, qui est d\`es lors un fibr\'e principal homog\`ene
sous $H$ au sens de~\Ref{XI.4.1}, et
\marginpar{296}
de fa\c con pr\'ecise est trivialis\'ee par le m\^eme
morphisme $S'\to S$ que~$P$. En particulier, \emph{\`a tout
fibr\'e principal homog\`ene $P$ sous $G$ et tout homomorphisme de~$S$-groupes $G\to H$, avec $H$ affine sur~$S$, est associ\'e un
fibr\'e principal homog\`ene de groupe~$H$}, de fa\c con
fonctorielle en $(G\to H)$, et compatible avec les changements de base
quelconques $T\to S$.

\begin{definition}
\label{XI.4.4}
Soit $G$ un $S$-pr\'esch\'ema. On note $\H^0(S,G)$ l'ensemble des
sections de~$G$ sur~$S$, qu'on consid\'erera comme un groupe lorsque
$G$ est un $S$-groupe. Dans ce cas, on note $\H^1(S,G)$ l'ensemble des
classes, \`a isomorphisme pr\`es, de fibr\'es principaux
homog\`enes sur~$S$ de groupe $S$, en consid\'erant $\H^1(S,G)$
comme muni du \og point marqu\'e\fg qui correspond aux fibr\'es
triviaux\footnote{\label{note-296}Cette notation n'est coh\'erente
\ifthenelse{\boolean{orig}}
{vis \`a vis}
{vis-\`a-vis}
des notations cohomologiques g\'en\'erales (SGA~4~V)
que lorsqu'on dispose de crit\`eres d'effectivit\'e de descente,
qui ne sont gu\`ere assur\'es que si $G$ est affine (ou seulement
quasi-affine, \cf (VIII~\Ref{VIII.7.9})).}.
\end{definition}

Ainsi, $\H^0(S,G)$ est un foncteur en le $S$-pr\'esch\'ema~$G$,
\`a valeurs dans la cat\'egorie des ensembles. Ce foncteur est
exact \`a gauche, a fortiori commute aux produits finis, ce qui
implique en effet qu'il transforme groupes en groupes, groupes
commutatifs en groupes commutatifs. De fa\c con analogue,
$\H^1(S,G)$ est un foncteur en le $S$-groupe \emph{affine}~$G$, \`a
valeurs dans la cat\'egorie des ensembles gr\^ace \`a la
formation des fibr\'es associ\'es; on constate facilement que ce
foncteur commute aux produits finis. En particulier il transforme les
groupes dans la cat\'egorie des $S$-groupes affines, \ie les
$S$-groupes \emph{affines commutatifs}, en des groupes de la seconde,
et m\^eme en des groupes commutatifs (puisque les groupes de la
premi\`ere cat\'egorie sont commutatifs). Ainsi, \emph{si $G$ est
un $S$-groupe affine commutatif, $\H^1(S,G)$ est un groupe
commutatif}, et un homomorphisme $G\to H$ de~$S$-groupes affines
commutatifs donne naissance \`a un homomorphisme de groupes
$\H^1(S,G)\to\H^1(S,H)$.

Pour simplifier, nous nous bornons pour la suite \`a la
consid\'eration de~$S$-groupes \emph{affines et commutatifs}. Soit
$$
0\to G'\lto{u} G\lto{v}G''\to 0
$$
une suite de morphismes de tels groupes, \emph{nous dirons que cette
suite est exacte si $vu=0$} (ce qui permet de consid\'erer $G$
comme un pr\'esch\'ema sur~$G''$, \`a groupes d'op\'erateurs
\`a droite~$G'_{G''}$) \emph{et si $G$ est un fibr\'e principal
homog\`ene
\marginpar{297}%
sur~$G''$ de groupe $G'_{G''}=G'\times_S G''$}. Cela implique en
particulier que $u\colon G'\to G$ est un noyau de~$v$, et a fortiori
cela implique l'exactitude de la suite
$$0\to\H^0(X,G')\to\H^0(X,G)\to\H^0(X,G'').$$ Cela implique de plus la
possibilit\'e de d\'efinir une application
$$
\partial\colon\H^0(X,G'')\to\H^1(X,G')\quoi,
$$
en associant \`a toute section de~$G''$ sur~$S$, \ie \`a tout
$S$-morphisme $f\colon S\to G''$, le fibr\'e principal homog\`ene
$P_f$ de groupe $G'\simeq f^*(G'_{G''})$ sur~$S$, image inverse du
fibr\'e principal homog\`ene $G$ sur~$G''$. Du point de vue
$S$-pr\'esch\'emas, ce n'est donc autre que l'image inverse par
$v\colon G\to G''$ du sous-pr\'esch\'ema image de~$S$ par
l'immersion~$f$, et les op\'erations de~$G'$ sur~$P_f$ sont induites
par les op\'erations \`a droite de~$G'$
\ifthenelse{\boolean{orig}}
{sur~$G$).}
{sur~$G$.}
Nous laissons \'egalement au lecteur la v\'erification de la
proposition suivante, qui ne pr\'esente pas de difficult\'es
autres que de r\'edaction:

\begin{proposition}
\label{XI.4.5}
L'application $\partial\colon \H^0(X,G'')\to\H^1(X,G')$ est un
homomorphisme de groupes. La suite d'homomorphismes suivante est
exacte:
\begin{multline*}
0\to\H^0(X,G')\lto{u^0}\H^0(X,G)\lto{v^0}\H^0(X,G')\lto{\partial}\H^1(X,G'')\\
\lto{u^1}\H^1(X,G) \lto{v^1}\H^1(X,G'')
\end{multline*}
(o\`u les homomorphismes autres que $\partial$ proviennent de la loi
fonctorielle de~$\H^0$ \resp~$\H^1$).
\end{proposition}

\begin{remarques}
\label{XI.4.6}
Le point de vue expos\'e ici pour l'\'etude des fibr\'es
principaux homog\`enes est visiblement inspir\'e de Serre~\cite{XI.7}, que
le lecteur aura tout int\'er\^et \`a consulter. Lorsqu'on veut
un formalisme qui s'applique \'egalement \`a des $S$-groupes
structuraux qui sont quasi-projectifs sur~$S$ (de fa\c con \`a
englober les sch\'emas ab\'eliens projectifs en particulier), on a
int\'er\^et \`a modifier \Ref{XI.4.1} en y demandant que $S'$
soit somme de pr\'esch\'emas $S'_i$ qui sont finis et localement
libres sur des ouverts $S_i$ de~$S$ recouvrant $S$. Les
d\'eveloppements pr\'ec\'edents sont alors valables, y inclus
notamment \Ref{XI.4.5}, en rempla\c cant partout l'hypoth\`ese
affine par l'hypoth\`ese quasi-projective, et en interpr\'etant de
fa\c con correspondante la d\'efinition donn\'ee plus faut d'une
suite exacte de~$S$-groupes. Il suffit en effet de remplacer la
r\'ef\'erence \`a (VIII~\Ref{VIII.2.1})
\marginpar{298}
par (VIII~\Ref{VIII.7.7}): les morphisme utilis\'es $S'\to S$ sont
des morphismes de descente effective pour la cat\'egorie fibr\'ee
des pr\'esch\'emas quasi-projectifs sur d'autres. On fera
attention cependant que cette deuxi\`eme notion de fibr\'e
principal homog\`ene est plus restrictive que la
premi\`ere~\Ref{XI.4.1}.
\end{remarques}

\subsection{}
\label{XI.4.7}
On obtient une notion encore plus restrictive de fibr\'e principal
homog\`ene en demandant que $S$ soit recouvert par des ouverts $S_i$
tels que pour tout~$i$, $P|S_i$ soit un fibr\'e \`a op\'erateurs
trivial sous~$G|S_i$: on dira alors que $P$ est un fibr\'e principal
homog\`ene \emph{localement trivial}.
\index{principal homog\`ene localement trivial (fibr\'e)|hyperpage}%
Les classes de ces fibr\'es, pour $G$ donn\'e, forment une partie
de~$\H^1(X,G)$, qui est en correspondance biunivoque avec
$\H^1(X,\cal{O}_X(G))$, o\`u $\cal{O}_X(G)$ est le faisceau (au sens
ordinaire) des sections de~$G$ sur~$S$, \cf \cite{XI.2}. Pour que ces $\H^1$
donnent encore lieu \`a une suite exacte de
cohomologie~\Ref{XI.4.5}, il faut \'evidemment supposer que la suite
$0\to G'\to G\to G''\to 0$ soit exacte au sens raisonnable pour ce
nouveau contexte, \ie que $G$ soit fibr\'e localement trivial
sur~$G''$, de groupe $G'_{G''}$; cela signifie aussi que $u\colon
G'\to G$ est un noyau de~$v\colon G\to G''$, et que $G$ admet
localement une section sur~$G''$.

\subsection{}
\label{XI.4.8}
Il est \'evidemment tr\`es d\'esirable de continuer la suite
exacte~\Ref{XI.4.5} en introduisant les groupes de cohomologie
sup\'erieurs $\H^i(X,G)$. Cela est possible en se pla\c cant dans
le cadre de la \og Cohomologie de Weil\fg: on consid\`ere la
cat\'egorie $\cal{B}$ des pr\'esch\'emas quasi-compacts sur~$S$,
muni de l'ensemble $\cal{M}$ des morphismes fid\`element plats et
quasi-compacts, qu'on appellera morphismes localisants. Un \og faisceau
de Weil\fg ab\'elien sur~$S$ (ou mieux, sur~$(\cal{B},\cal{M})$) est
alors un foncteur contravariant $\cal{F}$ de~$\cal{B}$ dans la
cat\'egorie des groupes ab\'eliens, transformant sommes en
produits, et une suite
${\def\labelstyle{\scriptstyle}\xymatrix@C=.5cm{T''=T'\times_TT'\ar@<+2pt>[rr]^-{\pr_1,\pr_2}\ar@<-2pt>[rr] &&T'\ar[r]^f & T}}$, avec $f\in\cal{M}$, en un diagramme \emph{exact}
d'ensembles $\xymatrix@C=.5cm{\cal{F}(T)\ar[r] & \cal{F}(T')
\ar@<2pt>[r]\ar@<-2pt>[r] & \cal{F}(T'')}$. Les faisceaux de Weil
forment une cat\'egorie ab\'elienne \`a limites inductives
exactes admettant un g\'en\'erateur, donc admettant suffisamment
d'objets injectifs~\cite{XI.1}. Les foncteurs d\'eriv\'es droits du
foncteur $\Gamma(\cal{F})=\cal{F}(S)$ sont alors not\'es
$\H^i(S,\cal{F})$. D'autre part, tout $S$-groupe commutatif
d\'efinit \'evidemment un faisceau de Weil (VIII~\Ref{VIII.5.2}),
dont le $\H^0$ et $\H^1$ ne sont autres que $\H^0(S,G)$ et
$\H^1(S,G)$, ce qui permet de d\'efinir les autres $\H^i(S,G)$ de
fa\c con raisonnable. On
\marginpar{299}
montre d'ailleurs qu'une suite exacte de~$S$-groupes d\'efinit une
suite exacte de faisceaux de Weil, ce qui permet de retrouver et de
prolonger la suite exacte~\Ref{XI.4.5}\footnote{Pour une \'etude
syst\'ematique de ce point de vue, \cf SGA~4~I~\`a~IX.}.

\subsection{}
\label{XI.4.9}
Il serait indiqu\'e de d\'evelopper les variantes non commutatives
de~\Ref{XI.4.5} comme dans~\cite{XI.2}. Pour un d\'eveloppement
syst\'ematique, dans le cadre qui convient, des diverses notions
cohomologiques esquiss\'ees dans le pr\'esent num\'ero, nous
renvoyons \`a un travail en pr\'eparation de
J\ptbl \textsc{Giraud}\footnote{\Cf J\ptbl \textsc{Giraud}, 
\ifthenelse{\boolean{orig}}
{\emph{Alg\`ebre homologique non ab\'elienne}, \`a para\^itre dans}
{\emph{Cohomologie non ab\'elienne},}
Springer-Verlag 1971.}.

\section{Cas particuliers de fibr\'es principaux}
\label{XI.5}

Supposons maintenant que $S$ soit connexe, et muni d'un point
g\'eom\'etrique~$a$, d'o\`u un groupe fondamental $\pi_1(S,a)$
permettant de classifier les rev\^etements \'etales de~$S$: la
cat\'egorie des rev\^etements \'etales de~$S$ est
\'equivalente \`a la cat\'egorie des ensembles finis o\`u
$\pi_1$ op\`ere contin\^ument. Il s'ensuit qu'un sch\'ema en
groupes fini et \'etale $G$ sur~$S$ est d\'etermin\'e
essentiellement par un groupe fini ordinaire $\cal{G}$, sur lequel
$\pi_1$ op\`ere contin\^ument par automorphismes de groupe. Un
rev\^etement \'etale $P$ de~$S$ o\`u $G$ op\`ere \`a droite
est d\'etermin\'e essentiellement par un ensemble fini $\cal{P}$
o\`u $\pi_1$ op\`ere contin\^ument (\`a gauche), et sur lequel
$\cal{G}$ op\`ere \`a droite de fa\c con compatible avec les
op\'erations de~$\pi_1$:
$$
s(p\cdot g)=(sp)\cdot(sg)\qquad \text{pour }s\in\pi_1,\
p\in\cal{P},\ g\in\cal{G}.
$$
On v\'erifie que $P$ est un fibr\'e principal homog\`ene au sens
de~\Ref{XI.4.1} si et seulement si $\cal{P}$ est un ensemble principal
homog\`ene sous~$\cal{G}$ (utiliser par exemple le
crit\`ere~\Ref{XI.4.2}). En d'autres termes, \emph{la cat\'egorie
des fibr\'es principaux homog\`enes sur~$S$ de groupe $G$ est
\'equivalente \`a la cat\'egorie des fibr\'es principaux
homog\`enes de groupe $G$ dans la cat\'egorie des ensembles finis
o\`u $\pi_1$ op\`ere contin\^ument.} On en d\'eduit en
particulier une bijection canonique, fonctorielle en $G$:
\begin{equation*}
\label{eq:XI.5.etoile}
\tag{$*$} \H^1(S,G)\isomto \H^1(\pi_1,\cal{G})\quoi,
\end{equation*}
o\`u
\marginpar{300}
le deuxi\`eme membre d\'esigne l'ensemble des classes, \`a
isomorphisme pr\`es, des fibr\'es principaux homog\`enes sous
$\cal{G}$ dans la cat\'egorie des ensembles finis o\`u $\pi_1$
op\`ere (inutile d'ailleurs de pr\'eciser: contin\^ument),
ensemble qui s'explicite de fa\c con bien connue comme ensemble
quotient de l'ensemble $Z^1(\pi_1,G)$ des $1$-cocycles
$\varphi:\pi_1\to \cal{G}$ (satisfaisant $\varphi(1)=1$,
$\varphi(st)=\varphi(s)(s.\varphi(t))$) par le groupe $\cal{G}$ qui y
op\`ere de fa\c con naturelle).

\label{page-300}
Un cas important est celui o\`u $\pi_1$ op\`ere trivialement dans~$\cal{G}$, \ie lorsque $G$ est un rev\^etement compl\`etement
d\'ecompos\'e de~$S$, isomorphe \`a la somme de~$\cal{G}$
exemplaires de~$S$; on \'ecrit alors aussi $\H^1(S,\cal{G})$ au lieu
de~$\H^1(S,G)$, et cet ensemble est en correspondance biunivoque
\eqref{eq:XI.5.etoile} avec $\H^1(\pi_1,\cal G)= \Hom(\pi_1,\cal{G})/$automorphismes int\'erieurs de~$\cal{G}$. On notera d'ailleurs que dans ce cas, un fibr\'e
principal homog\`ene sur~$S$ de groupe $G$ n'est autre chose qu'un
\emph{rev\^etement principal} de~$S$ de groupe $\cal{G}$
(V~\Ref{V.2.7}), et la correspondance biunivoque pr\'ec\'edente
est celle qui se d\'eduit de la correspondance entre rev\^etements
principaux de~$S$ de groupe~$\cal{G}$, \emph{ponctu\'es} au-dessus
de~$a$, et les homomorphismes continus de~$\pi_1(S,a)$ dans $\cal{G}$
(\Ref{V}~fin~du~\No \Ref{V.5}).

L'int\'er\^et de relier la th\'eorie des rev\^etements
\'etales avec celle des fibr\'es principaux (d\'ej\`a
implicite dans \textsc{A\ptbl Weil}, G\'en\'eralisation des Fonctions
Ab\'eliennes, et explicit\'ee par \textsc{S\ptbl Lang} dans sa th\'eorie
g\'eom\'etrique du corps de classes, \cf Serre~\cite{XI.5}), vient du fait
que tout $S$-groupe qui est fini et \'etale sur~$S$ peut se plonger
dans un $S$-groupe~$H$, affine et lisse sur~$S$, \`a fibres
connexes, commutatif lorsque $G$ l'est de sorte que par la suite
exacte~\Ref{XI.4.5} (et \'eventuellement ses variantes non
commutatives), la classification \og discr\`ete\fg des rev\^etements
principaux de groupe $G$ peut s'\'etudier \`a l'aide de la
classification \og continue\fg des fibr\'es principaux de groupe~$H$,
et r\'eciproquement d'ailleurs. Pour l'id\'ee de la construction
g\'en\'erale de l'immersion de~$G$ dans $H$ (assez peu
utilis\'ee en pratique semble-t-il), se reporter \`a
\cite[VI~2.8]{XI.5}. Nous nous contentons de d\'evelopper au \No suivant
deux cas particuliers importants, d'ailleurs classiques. Nous y aurons
besoin d'un r\'esultat auxiliaire:

\begin{proposition}
\label{XI.5.1}
Soit
\marginpar{301}
$S$ un pr\'esch\'ema, $G$ un $S$-groupe isomorphe \`a
$\Gl(n)_S$ (par exemple $\GG_{m\,S}$) ou $\GG_{a\,S}$, alors tout
fibr\'e principal homog\`ene sous $G$ est localement trivial.
\end{proposition}

Pr\'ecisons que $\Gl(n)_S$ ($n$ entier $\geq 0$) d\'esigne le
$S$-groupe qui repr\'esente le foncteur contravariant $T\mto
\Gl(n,\Gamma(T,\cal{O}_T))$ en le $S$-pr\'esch\'ema~$T$, en
particulier $\GG_{m\,S}$ (\og groupe multiplicatif sur~$S$\fg)
repr\'esente le foncteur contravariant $T\mto
\Gamma(T,\cal{O}_T^*)$, donc comme pr\'esch\'ema sur~$S$ est
isomorphe \`a $\Spec {\cal{O}_S[t,t^{-1}]}$, o\`u $t$ est une
ind\'etermin\'ee. De m\^eme $\GG_{a\,S}$ repr\'esente le
foncteur contravariant $T\mto \Gamma(T,\cal{O}_T)$, il est donc
isomorphe comme $S$-pr\'esch\'ema \`a $\Spec(\cal{O}_S[t])$,
o\`u $t$ est une ind\'etermin\'ee. Notons que par d\'evissage,
\Ref{XI.5.1} redonne le r\'esultat de locale trivialit\'e de
Rosenlicht, relatif au cas o\`u $G$ admet une \og suite de
composition\fg dont les facteurs cons\'ecutifs sont des groupes du
type envisag\'e ici. (Pour une \'etude plus fine des questions de
locale trivialit\'e des fibr\'es principaux homog\`enes,
\cf \cite{XI.7}~et~\cite{XI.3}).

La premi\`ere assertion se d\'emontre en remarquant que
$G(T)=\Aut(\cal{O}_T^n)$, et que les morphismes $S'\to S$ intervenant
dans~\Ref{XI.4.1} (\ie qui sont fid\`element plats et
quasi-compacts) sont des morphismes de descente effective pour la
cat\'egorie fibr\'ee des Modules localement isomorphes \`a
$\cal{O}_T^n$, \ie localement libres de rang~$n$
(VIII~\Ref{VIII.1.12}). La deuxi\`eme se d\'emontre de fa\c con
analogue, en notant que dans ce cas on a $G(T)=\Aut(\cal{E}_T)$,
o\`u $\cal{E}_T$ est l'\emph{extension} triviale de~$\cal{O}_T$ par
$\cal{O}_T$ (et o\`u les automorphismes bien entendu doivent
respecter la structure d'extension), et que les morphismes $S'\to S$
intervenant dans~\Ref{XI.4.1} sont des morphismes de descente
effective pour la cat\'egorie fibr\'ee des extensions de~$\cal{O}_T$ par $\cal{O}_T$ (comme il r\'esulte facilement de
VIII~\Ref{VIII.1.1}), et que de telles extensions sont automatiquement
localement triviales.

\begin{remarque}
\label{XI.5.2}
On notera que le m\^eme type de d\'emonstration s'applique au
groupe symplectique $\Symp(2n)_S$, compte tenu qu'une forme
altern\'ee sur un module localement isomorphe \`a
$\cal{O}_S^{2n}$, qui est \og non d\'eg\'en\'er\'ee\fg
\ie d\'efinit un isomorphisme de ce Module sur son dual, est
localement isomorphe \`a la forme standard. Le r\'esultat pour le
groupe orthogonal est par contre faux, d\'ej\`a si $S$ est le
spectre d'un corps, car il peut y avoir des formes
\marginpar{302}
quadratiques sur un corps qui ne sont pas isomorphes \`a la forme
standard. D'ailleurs on montre essentiellement dans \cite{XI.3} que les
groupes $\Gl$, $\Symp$, $\GG_a$ et ceux qui se d\'evissent en tels
groupes, sont \`a peu de choses pr\`es les seuls pour lesquels on
ait un r\'esultat de trivialit\'e locale du type consid\'er\'e
ici.
\end{remarque}

\begin{corollaire}
\label{XI.5.3}
On a des bijections canoniques
$$
\H^1(S,\GL(n)_S)\isomfrom \H^1(S,\GL(n,\cal{O}_S))\quoi,
$$
en particulier
$$
\H^1(S,\GG_{m\,S})\isomfrom \H^1(S,\cal{O}_S^*)\quoi,
$$
et
$$
\H^1(S,\GG_{a\,S})\isomfrom \H^1(S,\cal{O}_S)\quoi,
$$
o\`u les deuxi\`emes membres d\'esignent des cohomologies de
l'espace topologique $S$ \`a coefficients dans des faisceaux
ordinaires.
\end{corollaire}

En particulier, $\H^1(S,\GL(n)_S)$ s'identifie \`a l'ensemble des
classes, \`a isomorphisme pr\`es, de Modules localement libres de
rang $n$ sur~$S$, et $\H^1(S,\GG_{a\,S})$ s'identifie \`a l'ensemble
des classes d'extensions du Module $\cal{O}_S$ par lui-m\^eme.

\section[Application aux rev\^etements principaux]{Application aux rev\^etements principaux: th\'eories de
Kummer et d'Artin-Schreier}
\label{XI.6}


\begin{proposition}
\label{XI.6.1}
Soient $S$ un pr\'esch\'ema, $n$ un entier $>0$, soit $u_n\colon
\GG_{m\,S}\to\GG_{m\,S}$ l'homomorphisme de puissance \nieme, et
$\bbmu_{n\,S}$
\label{indnot:kb}\oldindexnot{$\bbmu_{n\,S}$|hyperpage}%
son noyau. Alors $\bbmu_{n\,S}$ est fini et localement libre de
rang $n$ sur~$S$, et il est \'etale sur~$S$ si et seulement si pour
tout $s\in S$, la caract\'eristique de~$s$ est premi\`ere \`a
$n$. La suite d'homomorphismes
$$
0\to \bbmu_{n\,S}\to\GG_{m\,S}\lto{u_n}\GG_{m\,S}\to 0
$$
est exacte au sens du \No \Ref{XI.4}. (On l'appellera la \emph{suite exacte de
Kummer}
\index{Kummer (suite exacte de)|hyperpage}%
sur~$S$, relativement \`a l'entier~$n$).
\end{proposition}

On
\marginpar{303}
a
$$
\GG_m=\Spec\cal{O}_S[t,t^{-1}],
$$
et $u_n$ correspond \`a l'homomorphisme $u_n$ sur les
$\cal{O}_S$-alg\`ebres affines, donn\'e par
$$
u_n(t)=t^n,
$$
d'autre part la section unit\'e de~$\GG_{m\,S}$ correspond \`a
l'homomorphisme d'augmentation de~$\cal{O}_S$-alg\`ebres, donn\'e
par
$$
\varepsilon(t)=1,
$$
dont le noyau est donc l'Id\'eal principal $(t-1)$. L'image de ce
dernier par $u_n$ est donc l'Id\'eal principal $(1-t^n)$, et on
trouve:
$$
\bbmu_{n\,S}=\Spec\cal{O}_S[t]/(1-t^n),
$$
ce qui montre en particulier que $\bbmu_{n\,S}$ est fini sur~$S$,
et d\'efini par une Alg\`ebre sur~$S$ qui est libre de rang $n$,
ayant la base form\'ee des $t^i$ ($0\leq i\leq n-1$). Pour qu'il soit
\'etale en $s\in S$, il faut et il suffit que l'Alg\`ebre
r\'eduite $k[t]/(1-t^n)$, o\`u $k=\kres(s)$, obtenue par adjonction
formelle des racines \niemes de l'unit\'e \`a $k$, soit
s\'eparable sur~$k$, \ie les racines de~$1-t^n$ dans une
cl\^oture alg\'ebrique de~$k$ sont toutes distinctes, ce qui
\'equivaut au fait que $n$ soit premier \`a la
caract\'eristique. Enfin, pour montrer que la suite
d'homomorphismes dans~\Ref{XI.6.1} est exacte, on est ramen\'e en
vertu du crit\`ere~\Ref{XI.4.2} \`a prouver que
\ifthenelse{\boolean{orig}}
{$v$}
{$u_n$}
est fid\`element plat. On peut \'evidemment supposer $S$ affine
d'anneau~$A$, donc $\GG_{m\,S}$ affine d'anneau $B=A[t,t^{-1}]$, et il
suffit de v\'erifier que $u_n$ fait de~$B$ un module libre de rang
$n$ sur~$B$, ou ce qui revient au m\^eme, que $u_n$ est injectif, et
que $A[t,t^{-1}]$ est un module libre de rang $n$ sur~$A[t^n,t^{-n}]$. En effet, on v\'erifie facilement que les
$t^i$ ($0\leq i\leq n-1$) forment une base de l'un sur l'autre, ce qui
ach\`eve la d\'emonstration.

\begin{definition}
\label{XI.6.2}
On
\marginpar{304}
appelle $\bbmu_{n\,S}$ le \emph{groupe de Kummer de rang $n$
sur~$S$}, et on appelle \emph{rev\^etement principal kumm\'erien
de rang
\ifthenelse{\boolean{orig}}
{$n\,S$}}
{$n$ sur~$S$}
tout fibr\'e principal homog\`ene sur~$S$ de groupe le groupe de
Kummer de rang~$n$.
\end{definition}

L'ensemble de ces rev\^etements est un groupe, not\'e
$\H^1(S,\bbmu_{n\,S})$ ou simplement $\H^1(S,\bbmu_n)$. On
notera que la formation du groupe de Kummer de rang $n$ sur~$S$ est
compatible avec l'extension de la base, donc que $\bbmu_{n\,S}$
provient par extension de la base du \emph{groupe de Kummer absolu
$\bbmu_n$} sur~$\Spec(\ZZ)$.

D\'esignons par $(\ZZ/n\ZZ)_S$ le $S$-groupe d\'efini par le
groupe fini ordinaire $\ZZ/n\ZZ$. Si~$G$ est un $S$-groupe quelconque,
les homomorphismes de~$S$-groupes $u$ de~$(\ZZ/n\ZZ)_S$ dans $G$
correspondent biunivoquement, et de fa\c con compatible avec le
changement de base, aux sections de~$G$ sur~$S$ dont la puissance
\nieme est la section unit\'e, en faisant correspondre \`a
$u$ l'image par $u$ de la section de~$(\ZZ/n\ZZ)_S$ sur~$S$ d\'efini
par le g\'en\'erateur $1\bmod n\ZZ$ de~$\ZZ/n\ZZ$. Ceci pos\'e:

\begin{corollaire}
\label{XI.6.3} Si $\bbmu_{n\,S}$ est \'etale sur~$S$, on
obtient ainsi une correspondance biunivoque entre les isomorphismes de~$S$-groupes $(\ZZ/n\ZZ)_S\isomto\bbmu_{n\,S}$, et les sections de~$\cal{O}_S$ qui sont d'ordre $n$ exactement sur chaque composante
connexe de~$S$ (une telle section s'appellera \og racine primitive
\nieme de l'unit\'e sur~$S$\fg). Donc pour que
$\bbmu_{n\,S}$ soit isomorphe en tant que $S$-groupe \`a
$(\ZZ/n\ZZ)_S$, il faut et il suffit qu'il soit \'etale
\ifthenelse{\boolean{orig}}
{sur~$S$}
{sur~$S$,}
\ie que les caract\'eristiques r\'esiduelles de~$S$ soient
premi\`eres \`a $n$, et qu'il existe une racine primitive
\nieme de l'unit\'e sur~$S$.
\end{corollaire}

Cela explique le r\^ole jou\'e dans la th\'eorie kumm\'erienne
classique par l'hypoth\`ese que le corps de base (jouant le r\^ole
de~$S$) soit de caract\'eristique premi\`ere \`a $n$ et
contienne les racines \niemes de l'unit\'e, et par le choix
d'une racine primitive \nieme de l'unit\'e. Une fois qu'on
dispose du langage des sch\'emas, il n'y a plus lieu de s'embarrasser
de ces hypoth\`eses, et il convient de raisonner directement sur~$\bbmu_n$ au lieu de~$\ZZ/n\ZZ$. Ainsi, la conjonction
de~\Ref{XI.6.1}, \Ref{XI.4.5} et~\Ref{XI.5.3} nous donne la relation
g\'en\'erale suivante entre la th\'eorie des rev\^etements
principaux kumm\'eriens et celle des groupes de Picard:

\begin{proposition}
\label{XI.6.4}
Soient
\marginpar{305}
$S$ un pr\'esch\'ema, on a une suite exacte canonique
$$
\textstyle
0\to\H^0(S,\bbmu_n)\to\H^0(S,\cal{O}^*_S)\to
\H^0(S,\cal{O}^*_S)\to\H^1(S,\bbmu_n)\to
\H^1(S,\cal{O}^*_S)\to\H^1(S,\cal{O}^*_S),
$$
d'o\`u, en posant
\ifthenelse{\boolean{orig}}
{$\H^1(S,\cal{O}_S)=\Pic(S)$,}
{$\H^1(S,\cal{O}_S^*)=\Pic(S)$,}
et en d\'esignant pour tout groupe ab\'elien $A$, par ${_n}A$ et
$A_n$ les noyau et conoyau de la multiplication par $n$ dans $A$, la
suite exacte:
$$
0\to\H^0(S,\cal{O}_S)^*_n\to\H^1(S,\bbmu_n)\to
{\sideset{_n}{}\Pic}(S)\to 0.
$$
\end{proposition}

Nous allons expliciter deux cas importants, o\`u l'un ou l'autre
terme extr\^eme de cette suite exacte sont nuls:

\begin{corollaire}
\label{XI.6.5} Supposons $\sideset{_n}{}\Pic(S)=0$, (par exemple
que $S$ soit le spectre d'un anneau local, ou d'un anneau factoriel),
et soit $A$ l'anneau $\H^0(S,\cal{O}_S)$. Alors on a un isomorphisme
canonique
$$
\H^1(S,\bbmu_n)\isomto A^*/{A^*}^n.
$$
\end{corollaire}

C'est essentiellement l'\'enonc\'e classique de la th\'eorie de
Kummer, lorsque $S$ est le spectre d'un corps.

\begin{corollaire}
\label{XI.6.6} Supposons que tout \'el\'ement de
$\H^0(S,\cal{O}_S)$ soit une puissance \nieme, par exemple que
$\H^0(S,\cal{O}_S)$ soit un compos\'e de corps alg\'ebriquement
clos ou que $S$ soit r\'eduit et propre sur un corps
alg\'ebriquement clos $k$. Alors on a un isomorphisme canonique
$$
\H^1(S,\bbmu_n)\isomto\sideset{_n}{}\Pic(S).
$$
\end{corollaire}

En particulier, lorsque $S$ est propre et connexe sur un corps
alg\'ebriquement clos $k$, cela met en relation le groupe
fondamental de~$S$ avec les points d'ordre fini du sch\'ema de
Picard $P$ de~$S$ sur~$k$; ainsi on aura un isomorphisme
$$
\Hom(\pi_1(S),\ZZ/n\ZZ)\simeq {_n P}(k)
$$
pour $n$ premier \`a la caract\'eristique, qui est souvent
utilis\'e en g\'eom\'etrie
\marginpar{306}
alg\'ebrique. Comme application, lorsque la composante connexe $P^0$
de~$P$ est un sch\'ema en groupes complet, de dimension $g$, on voit
en utilisant les r\'esultats rappel\'es dans le \No \Ref{XI.2}, et la
finitude du groupe de torsion de
\ifthenelse{\boolean{orig}}{N\'eron-S\'ev\'eri}{N\'eron-Severi},
que pour
tout nombre premier $\ell$ premier \`a la caract\'eristique, la
composante $\ell$-primaire du groupe fondamental $\pi_1(S)$ rendu
ab\'elien est un module de type fini et de rang $2g$ sur l'anneau
$\ZZ_\ell$ des entiers $\ell$-adiques (et d'ailleurs libre sauf pour
un nombre fini au plus de valeurs de~$\ell$). Comme l'a remarqu\'e
Serre, cela permet de prouver sous certaines conditions que lorsque
$X$ est un sch\'ema plat et projectif sur~$S$ connexe, alors les
sch\'emas de Picard des fibres de~$X$ ont toutes la m\^eme
dimension, en appliquant le th\'eor\`eme de semicontinuit\'e
(X~\Ref{X.2.3}); l'argument de Serre s'applique d\`es que le
sch\'ema de Picard de~$X$ sur~$S$ existe, et que les Picards
connexes des fibres de~$X$ sur~$S$ sont des sch\'emas en groupes
propres, par exemple lorsque les fibres g\'eom\'etriques de~$X$
sur~$S$ sont normales ($X$ \'etant toujours plat et projectif
sur~$S$), en particulier si $X$ est lisse et projectif sur~$S$.

Soit maintenant $p$ un nombre premier, et supposons que $S$ soit un
pr\'esch\'ema de caract\'eristique $p$, \ie tel que
$p\cdot\cal{O}_S=0$. Alors l'homomorphisme de puissance \pieme
dans $\cal{O}_S$ est additif, et le morphisme correspondant, obtenu en
rempla\c cant $S$ par un~$T$ variable sur~$S$:
$$
F\colon\GG_{a\,S}\to\GG_{a\,S}
$$
est donc un homomorphisme de~$S$-groupes, qu'on appelle
l'\emph{homomorphisme de Frobenius} (N.B. Un tel morphisme est
d\'efini pour tout $S$-pr\'esch\'ema $G$ qui provient par
extension de la base d'un pr\'esch\'ema $G_0$ sur le corps premier
$\ZZ/p\ZZ$, et ce morphisme est un homomorphisme de groupes si $G_0$
est un pr\'esch\'ema en groupes). Nous poserons:
$$
\wp=\id-F\colon\GG_{a\,S}\to\GG_{a\,S}.
$$
Consid\'erons d'autre part le $S$-groupe $(\ZZ/p\ZZ)_S$ d\'efini
par le groupe fini ordinaire~$\ZZ/p\ZZ$, nous avons dit que pour tout
$S$-groupe~$G$, les homomorphismes de~$S$-groupe de~$(\ZZ/p\ZZ)_S$
dans $G$ correspondent biunivoquement aux sections de~$G$ sur~$S$ dont
la puissance \pieme est la section unit\'e. Lorsque
$G=\GG_{a\,S}$, elles correspondent
\marginpar{307}
donc aux sections quelconques de~$G$ sur~$S$. Prenant en particulier
la section de~$\GG_{a\,S}$ sur~$S$ correspondant \`a la section
unit\'e du faisceau d'anneaux $\cal{O}_S$, on trouve un
homomorphisme de~$S$-groupes
$$
\ifthenelse{\boolean{orig}}
{i\colon(\ZZ/p\ZZ)\to\GG_{a\,S}}
{i\colon(\ZZ/p\ZZ)_S\to\GG_{a\,S}}
$$

\begin{proposition}
\label{XI.6.7} La suite d'homomorphismes de~$S$-groupes
$$
\ifthenelse{\boolean{orig}}
{0\to(\ZZ/p\ZZ)_S\to\GG_{a\,S}\to G_{a\,S}\to 0}
{0\to(\ZZ/p\ZZ)_S\to\GG_{a\,S}\to \GG_{a\,S}\to 0}
$$
est exacte (au sens du \No \Ref{XI.4}). (On l'appelle la suite exacte
\emph{d'Artin-Schreier} sur~$S$).
\index{Artin-Schreier (suite exacte d')|hyperpage}%
\end{proposition}

Il suffit de le prouver sur le corps premier $k=\ZZ/p\ZZ$. Il suffit
de remarquer que l'homomorphisme $\wp^*\colon k[t]\to k[t]$ d\'efini
par $\wp^*(t)=t-t^p$ fait de~$k[t]$ un module libre de rang $p$ sur~$k[t]$, de fa\c con pr\'ecise que $k[t]$ est un module libre sur~$k[s]$, o\`u $s=t-t^p$, ayant la base form\'ee des $t^i$ ($0\leq i\leq p-1$).

On en conclut, utilisant \Ref{XI.4.5} et~\Ref{XI.5.3}:

\begin{proposition}
\label{XI.6.8} On a une suite exacte canonique:
\begin{multline*}
0\to\H^0(S,\ZZ/p\ZZ)\to\H^0(S,\cal{O}_S)\to\H^0(S,\cal{O}_S)\to
\H^1(S,\ZZ/p\ZZ)\\
\to\H^1(S,\cal{O}_S)\to\H^1(S,\cal{O}_S),
\end{multline*}
d'o\`u une suite exacte:
$$
0\to\H^0(S,\cal{O}_S)/\wp\H^0(S,\cal{O}_S)\to
\H^1(S,\ZZ/p\ZZ)\to\H^1(S,\cal{O}_S)^F\to 0,
$$
(o\`u l'exposant $F$ dans le dernier terme signifie le sous-groupe
des invariants par l'endomorphisme $F$, \'egal au noyau de
$\wp=\id-F$).
\end{proposition}

Explicitons encore deux cas extr\^emes:

\begin{corollaire}
\label{XI.6.9} Supposons que $\H^1(S,\cal{O}_S)^F=0$, par exemple
que $S$ soit un sch\'ema affine. Alors, posant
$A=\H^0(S,\cal{O}_S)$, on a un isomorphisme canonique
$$
\H^1(S,\ZZ/p\ZZ)\isomto A/\wp A.
$$
\end{corollaire}

C'est la \emph{th\'eorie d'Artin-Schreier} dans la forme classique,
du moins lorsque $A$ est le spectre d'un corps.

\begin{corollaire}
\label{XI.6.10}
Supposons
\marginpar{308}%
que $\wp\H^0(S,\cal{O}_S)\mkern1mu{=}\mkern1mu\H^0(S,\cal{O}_S)$, par exemple
que $\H^0(S,\cal{O}_S)$ soit un compos\'e de corps
alg\'ebriquement clos, ou que $S$ soit propre sur un corps
alg\'ebriquement clos. Alors on a un isomorphisme canonique:
$$
\H^1(S,\ZZ/p\ZZ)\isomto\H^1(S,\cal{O}_S)^F
$$
\end{corollaire}

\begin{remarques}
\label{XI.6.11}
Le dernier \'enonc\'e est d\^u \`a
\ifthenelse{\boolean{orig}}{J.P\ptbl}{J-P\ptbl}%
Serre
\cite{XI.9}. Il est
possible \'egalement de d\'evelopper une th\'eorie analogue pour
le groupe structural $\ZZ/p^n\ZZ$ pour $n$ quelconque, en utilisant au
lieu de~$\GG_a$ le sch\'ema en groupes de Witt $\WW_n$, \cf \loccit On notera qu'en caract\'eristique $p>0$, la th\'eorie de
Kummer ne donne plus de renseignement sur les rev\^etements
principaux d'ordre $p$, puisque $\bbmu_p$ est alors un groupe
\og infinit\'esimal\fg \ie radiciel sur la base, donc sans rapport
direct avec $\ZZ/p\ZZ$; aussi \`a premi\`ere vue, la th\'eorie
de ces rev\^etements n'est plus justiciable (lorsque $S$ est un
sch\'ema propre sur un corps alg\'ebriquement clos pour fixer les
id\'ees), de la th\'eorie du sch\'ema de Picard comme
dans~\Ref{XI.6.6}. N\'eanmoins, si on se rappelle que l'espace
tangent de Zariski \`a l'origine dans
$\mathbf{Pic}_{S/K}$\kern1pt\footnote{Pour la d\'efinition de
$\mathbf{Pic}_{S/K}$, \cf A. Grothendieck, S\'em. Bourbaki \No 232
(F\'evrier 1962).} s'identifie \`a $\H^1(S,\cal{O}_S)$, on
constate que \emph{la connaissance du sch\'ema en groupes
${_p}\mathbf{Pic}_{S/k}$, noyau de la multiplication par $p$ dans
$\mathbf{Pic}_{S/k}$, implique celle de~$\H^1(S,\ZZ/p\ZZ)$ aussi bien
que celle de~$\H^1(S,\bbmu_p)$; on notera qu'elle implique aussi
celle de~$\H^1(S,\bbalpha_p)$}, o\`u $\bbalpha_p$
d\'esigne le sch\'ema en groupes infinit\'esimal sur le corps
premier, noyau de~$F\colon\GG_a\to\GG_a$ (qui peut se d\'ecrire
aussi comme le spectre de l'alg\`ebre enveloppante restreinte de la
$p$-alg\`ebre de Lie triviale de dimension~$1$): en effet la suite
exacte~\Ref{XI.4.5} donne ici:
$$
\H^1(S,\bbalpha_p)\simeq\Ker
(F\colon\H^1(S,\cal{O}_S)\to\H^1(S,\cal{O}_S)),
$$
et plus g\'en\'eralement, d\'esignant par $\bbalpha_{p^n}$
le noyau dans $\GG_a$ du \nieme it\'er\'e de~$F$, on aura
$$
\H^1(S,\bbalpha_{p^n})\simeq\Ker(F^n\colon\H^1(S,\cal{O}_S)\to
\H^1(S,\cal{O}_S)).
$$

En fait, la connaissance de~${_p}\mathbf{Pic}_{S/k}$ \'equivaut \`a
celle de~$\H^1(S,G)$ pour tout groupe alg\'ebrique
\marginpar{309}
commutatif fini annul\'e par $p$, plus g\'en\'eralement, la
connaissance de~${_{p^n}}\mathbf{Pic}_{S/k}$ \'equivaut \`a celle
de~$\H^1(S,G)$ pour tout groupe alg\'ebrique commutatif fini $G$
annul\'e par $p^n$, en vertu du th\'eor\`eme suivant qui englobe
dans le cas envisag\'e \`a la fois la th\'eorie de Kummer et
celle de Artin-Schreier:

Soit $G$ un groupe alg\'ebrique fini sur~$k$,
$D(G)=\SheafHom_{k\textup{-groupes}}(G,\GG_m)$ son \emph{dual de
Cartier} (dont l'alg\`ebre affine est port\'ee par l'espace
vectoriel dual de l'alg\`ebre affine de~$G$, \ie par
l'hyperalg\`ebre de~$G$ au sens de Dieudonn\'e-Cartier), alors on
a un isomorphisme canonique:
\begin{equation}\label{eqXI.6.11}\tag{$*$}
\H^1(S,G)\simeq\Hom_{k\textup{-groupes}}(D(G),\mathbf{Pic}_{S/k}).
\end{equation}
(N.B. $S$ est un sch\'ema propre sur~$k$ alg\'ebriquement clos, tel
que $\H^0(S,\cal{O}_S)=k$). Cette formule peut encore s'exprimer en
disant que le \og vrai groupe fondamental\fg de~$S$ auquel il \'etait
fait allusion au \No \Ref{XI.2}, rendu ab\'elien, est isomorphe \`a la
limite projective des $D(P_i)$, o\`u $P_i$ parcourt les sous-groupes
alg\'ebriques \emph{finis} de~$\mathbf{Pic}_{S/k}$, qu'on notera
$T^\bbullet(\mathbf{Pic}_{X/k})$. Lorsque $S$ est une vari\'et\'e
ab\'elienne, on a vu dans~\Ref{XI.2.1} que ce groupe est
\'egalement isomorphe au \og vrai\fg module de Tate
$T_{\bbullet}(S)=\varprojlim\, {{_n}S}$, et l'isomorphisme
\eqref{eqXI.6.11} s'\'ecrit alors de fa\c con plus frappante
$$
\Ext^1(A,G)\simeq\Hom(D(G),B),
$$
$A$ \'etant une vari\'et\'e ab\'elienne, $B$ sa duale, $G$ un
groupe alg\'ebrique fini sur~$k$. Les r\'esultats qu'on vient
d'indiquer peuvent se g\'en\'eraliser d'ailleurs au cas o\`u $k$
est remplac\'e par un pr\'esch\'ema de base quelconque, et \`a
d'autres groupes de coefficients $G$ que des groupes finis.
\end{remarques}


\chapterspace{-4}
\chapter{G\'eom\'etrie alg\'ebrique et~g\'eom\'etrie~analytique}
\label{XII}
\marginpar{311}

\begin{center}
{Mme M.\ \textsc{Raynaud}\footnote{D'apr\`es des notes in\'edites de
A\ptbl Grothendieck.}}
\end{center}
\vspace*{1cm}

Proc\'edant comme dans \cite{XII.10}, on associe \`a tout sch\'ema $X$
localement de type fini sur le corps des nombres complexes $\CC$ un
espace analytique~$X^\an$
\label{indnot:lb}\oldindexnot{$X^\an$|hyperpage}%
dont l'ensemble sous-jacent est $X(\CC)$.

Dans les \No \Ref{XII.2} et~\Ref{XII.3} de cet expos\'e, nous donnons un
\og dictionnaire\fg entre les propri\'et\'es usuelles de~$X$ et de
$X^\an$ et entre les propri\'et\'es d'un morphisme $f\colon X\to
Y$ et du morphisme associ\'e $f^\an\colon X^\an \to Y^\an$.
\label{indnot:lc}\oldindexnot{$f^\an$|hyperpage}%

Nous montrons ensuite que les th\'eor\`emes de comparaison entre
faisceaux coh\'erents sur~$X$ et $X^\an$, \'etablis dans \cite[\No
12]{XII.10} pour une vari\'et\'e projective, sont encore valables lorsque
$X$ est un sch\'ema propre.

Enfin nous prouvons au \No \Ref{XII.5} l'\'equivalence de la cat\'egorie des
rev\^etements \'etales finis de~$X$ et de la cat\'egorie des
rev\^etements \'etales finis de~$X^\an$. En prime au lecteur, nous
donnons une nouvelle d\'emonstration du th\'eor\`eme de
Grauert-Remmert \cite{XII.6}, utilisant la r\'esolution des
singularit\'es~\cite{XII.8}.

\section{Espace analytique associ\'e \`a un sch\'ema}
\label{XII.1}
\marginpar{312}

Soit $X$ un sch\'ema localement de type fini sur~$\CC$. Soit $\Phi$
le foncteur de la cat\'egorie des espaces analytiques \cite[\No 9]{XII.4} dans
la cat\'egorie des ensembles, qui \`a un espace analytique~$\othercal{X}$ associe l'ensemble des morphismes d'espaces annel\'es en
$\CC$-alg\`ebres $\Hom_\CC (\othercal{X},X)$. On a le th\'eor\`eme
suivant:

\begin{theoremedefinition}
\label{XII.1.1} Le foncteur $\Phi$ est repr\'esentable par un
espace analytique $X^\an$ et un morphisme $\varphi\colon X^\an\to
X$. On dit que $X^\an$ est l'espace analytique associ\'e \`a~$X$.
\index{espace analytique associ\'e|hyperpage}%

Si $|X^\an|$ est l'ensemble sous-jacent \`a $X^\an$, $\varphi$
induit une bijection de~$|X^\an|$ sur l'ensemble $X(\CC)$ des points
de~$X$ \`a valeurs dans $\CC$. De plus, pour chaque point $x$ de~$X^\an$, le morphisme
$$
\varphi_x\colon\cal{O}_{X,\varphi(x)}\to\cal{O}_{X^\an,x},
$$
qui est n\'ecessairement local, donne par passage aux
compl\'et\'es un isomorphisme
$$
\widehat{\varphi}_x\colon\widehat{\cal{O}}_{X,\varphi(x)}\isomto
\widehat{\cal{O}}_{X^\an,x},
$$
En particulier le morphisme $\varphi$ est plat.
\end{theoremedefinition}

Notons que le fait que $\varphi$ induise une bijection de~$X^\an$ sur~$X(\CC)$ r\'esulte de la propri\'et\'e universelle de~$X^\an$.
D'autre part on a les assertions suivantes:

\begin{enumerate}
\item [a)] Si le th\'eor\`eme est vrai pour un sch\'ema $Y$, il
en est de m\^eme pour tout sous-sch\'ema $X$ de~$Y$. Supposons
d'abord que $X$ soit un sous-sch\'ema ouvert de~$Y$; si $\psi\colon
Y^\an\to Y$ est le morphisme canonique, $\psi^{-1}(X)$ est un ouvert
\marginpar{313}
de~$Y^\an$ que l'on muni de la structure d'espace analytique induite
par celle de~$Y^\an$. Comme tout morphisme d'un espace analytique
$\othercal{X}$ dans $X$ se factorise \`a travers $Y^\an$ d'apr\`es la
propri\'et\'e universelle de ce dernier, donc \`a travers
$X^\an$ qui est le produit fibr\'e $Y^\an\times_Y X$, $X^\an$ est
l'espace analytique associ\'e \`a $X$. Enfin l'assertion
concernant les $\varphi_x$ est \'evidente.

Il suffit maintenant de consid\'erer le cas o\`u $X$ est un
sous-sch\'ema ferm\'e de~$Y$. Soit~$I$ le $\cal{O}_Y$-Id\'eal
coh\'erent d\'efinissant $X$; alors $I\cdot\cal{O}_{Y^\an}$ est un
faisceau coh\'erent d'id\'eaux sur~$\cal{O}_{Y^\an}$ qui
d\'efinit un sous-espace analytique ferm\'e $X^\an$ de~$Y^\an$; on
voit comme dans le cas d'un sous-sch\'ema ouvert que $X^\an$ est
l'espace analytique associ\'e \`a $X$. Soit $\varphi\colon
X^\an\to X$ le morphisme canonique. Pour tout point $x$ de~$X^\an$, le
morphisme $\varphi_x$ n'est autre que le morphisme
$$
\cal{O}_{Y,\psi(x)}/I_{\psi(x)}\to
\cal{O}_{Y^\an,x}/I_{\psi(x)}\cdot\cal{O}_{Y^\an,x}
$$
induit par $\psi_x$; son compl\'et\'e
$$
\widehat{\varphi}_x\colon
\widehat{\cal{O}}_{Y,\psi(x)}/I_{\psi(x)}\cdot\widehat{\cal{O}}_{Y,\psi(x)}
\to\widehat{\cal{O}}_{Y^\an,x}/I_{\psi(x)}\cdot\widehat{\cal{O}}_{Y^\an,x}
$$
est un isomorphisme puisque $\widehat\psi_x$ en est un, ce qui
d\'emontre a).

\item [b)] Si l'on a deux $\CC$-sch\'emas $X_1$, $X_2$, tels que
$X_1^\an$ et $X_2^\an$ existent, alors il en est de m\^eme de
$(X_1\times X_2)^\an$. Soient en effet $\varphi_1\colon X_1^\an\to
X_1$, $\varphi_2\colon X_2^\an\to X_2$ les morphismes canoniques,
$p_1$, $p_2$ les deux projections de~$X_1^\an\times X_2^\an$. On
d\'eduit formellement de EGA~I~1.8.1 que $X_1\times X_2$ est le
produit de~$X_1$ et $X_2$ dans la cat\'egorie des espaces
annel\'es en anneaux locaux; il en r\'esulte que les morphismes
$\varphi_1\cdot p_1$ et $\varphi_2\cdot p_2$ d\'efinissent un
\marginpar{314}
morphisme $\varphi\colon X_1^\an\times X_2^\an\to X_1\times X_2$ et
que le couple $(X_1^\an\times X_2^\an,\varphi)$ repr\'esente le
foncteur $\othercal{X}\mto\Hom_\CC(\othercal{X},X_1\times X_2)$.

\item [c)] Si $\othercal{E}^1$ d\'esigne l'espace affine de dimension 1,
\ie l'espace topologique $\CC$ muni du faisceau des fonctions
holomorphes, le foncteur%
\ifthenelse{\boolean{orig}}
{}
{{\renewcommand{\thefootnote}{*}\addtocounter{footnote}{-1}\footnote{\lcrochetbf Ajout\'e en 2003: $E^1_\CC$ d\'esigne la droite affine (alg\'ebrique) sur $\CC$.\rcrochetbf}}}
$\othercal{X}\mto\Hom_\CC (\othercal{X},E^1_\CC)$
est repr\'esentable par $\othercal{E}^1$, le morphisme canonique
$\varphi\colon\othercal{E}^1\to E^1_\CC$ \'etant le morphisme
\'evident. En effet se donner un morphisme d'un espace analytique
$\othercal{X}$ dans $E^1_\CC$ \'equivaut \`a se donner un
\'el\'ement de~$\Gamma(\othercal{X},\cal{O}_\othercal{X})$, ce qui revient
aussi \`a se donner un morphisme de~$\othercal{X}$ dans~$\othercal{E}^1$. On
a \'evidemment une bijection $|\othercal{E}^1|\simeq E^1(\CC)$, et, pour
chaque point $x\in\othercal{E}^1$, le morphisme $\widehat\varphi_x$ n'est
autre que le morphisme identique d'un anneau de s\'eries formelles
\`a une variable sur~$\CC$.
\end{enumerate}

\ifthenelse{\boolean{orig}}{}
{\enlargethispage{.5cm}}%
On d\'eduit de b) et c) que le th\'eor\`eme est vrai pour
l'espace affine $E^n_\CC$, $n\geq 0$. Utilisant a), on voit qu'il en
est de m\^eme pour tout sch\'ema affine $X$, localement de type
fini sur~$\CC$. Si l'on ne suppose plus $X$ affine et si $(X_i)$ est
un recouvrement de~$X$ par des ouverts affines, il r\'esulte de la
propri\'et\'e universelle et de a) que les $X_i^\an$ se recollent
et d\'efinissent ainsi l'espace analytique $X^\an$ associ\'e
\`a~$X$.

\subsection{}
\label{XII.1.2}
Soit $f\colon X\to Y$ un morphisme de~$\CC$-sch\'emas localement de
type fini. Si $\varphi\colon X^{\an}\to X$ et $\psi\colon Y^{\an}\to
Y$ sont les morphismes canoniques, il r\'esulte de la
propri\'et\'e universelle de~$Y^{\an}$ qu'il existe un unique
morphisme $f^{\an}\colon X^{\an}\to Y^{\an}$ tel que le diagramme
$$
\xymatrix{
X^{\an} \ar[r]\ar[d]_-{f^{\an}}& X\ar[d]^-f\\ Y^{\an} \ar[r]& Y
}
$$
soit
\marginpar{315}
commutatif. On a donc d\'efini un foncteur $\Phi$ de la
cat\'egorie des $\CC$-sch\'emas localement de type fini dans la
cat\'egorie des espaces analytiques.

Le foncteur $\Phi$ commute aux limites projectives finies. Il suffit
en effet de voir que $\Phi$ commute aux produits fibr\'es. Or, si
$X$, $Y$, $Z$, sont des sch\'emas localement de type fini sur~$\CC$,
il r\'esulte du fait que $X\times_Z Y$ est le produit fibr\'e de
$X$ et $Y$ au-dessus de~$Z$ dans la cat\'egorie des espaces
annel\'es en anneaux locaux que $X^{\an}\times_{Z^{\an}}Y^{\an}$
satisfait \`a la propri\'et\'e universelle qui caract\'erise
$(X\times_Z Y)^{\an}$.

\subsection{}
\label{XII.1.3}
Soient $X$ un $\CC$-sch\'ema localement de type fini, $X^{\an}$
l'espace analytique associ\'e, $\varphi\colon X^{\an}\to X$ le
morphisme canonique. Si $F$ est un ${\cal{O}}_X$-Module, l'image
inverse $\varphi^*F=F^{\an}$ est un faisceau de modules sur~${\cal{O}}_{X^{\an}}$. On d\'efinit ainsi un foncteur de la
cat\'egorie des $\cal{O}_X$-Modules dans la cat\'egorie des
Modules sur~$X^{\an}$. Ce foncteur commute aux limites inductives
(EGA~0~4.3.2). Le faisceau $\cal{O}_{X^{\an}}$ \'etant coh\'erent
\cite[\no18~\S 2~th\ptbl 2]{XII.4}, il transforme faisceaux coh\'erents en
faisceaux coh\'erents (EGA~0~5.3.11). On a de plus:

\begin{subproposition}
\label{XII.1.3.1}
Le foncteur qui \`a un $\cal{O}_{X^{\an}}$-Module $F$ associe son
image inverse $F^{\an}$ sur~$X^{\an}$ est exact, fid\`ele,
conservatif.
\end{subproposition}

L'exactitude r\'esulte du fait que le morphisme $\varphi\colon
X^{\an}\to X$ est plat \eqref{XIII.1.1}. Prouvons que le foncteur $F\mto
F^{\an}$ est fid\`ele. Compte tenu de l'exactitude, il suffit de
montrer que, si $F^{\an}$ est nul, il en est de m\^eme de~$F$. Or,
pour tout point $x$ de~$X^{\an}$, on a alors $F_{\varphi
(x)}\otimes_{\cal{O}_{X,\varphi (x)}}\cal{O}_{X^{\an},x}=0$. Le
morphisme $\cal{O}_{X,\varphi (x)}\to \cal{O}_{X^{\an},x}$ \'etant
fid\`element plat, on a $F_{\varphi (x)}=0$ pour tout point
ferm\'e $\varphi (x)$ de~$X$, et, comme $X$ est de Jacobson
(EGA~IV~10.4.8), ceci implique que $F$ est nul.

Le
\marginpar{316}
fait que le foncteur $F\mto F^{\an}$ soit conservatif est formel
\`a partir de l'exactitude et de la fid\'elit\'e.

\section[Comparaison des propri\'et\'es d'un sch\'ema]{Comparaison des propri\'et\'es d'un sch\'ema et de l'espace
analytique \hbox{associ\'e}}
\label{XII.2}

\begin{proposition}
\label{XII.2.1}
Soient $X$ un $\CC$-sch\'ema localement de type fini, $X^{\an}$
l'espace analytique associ\'e, $n$ un entier. Consid\'erons la
propri\'et\'e $P$ d'\^etre
\begin{itemize}
\item[(i)] non vide
\item[(i')] discret
\item[(ii)] de Cohen-Macaulay
\item[(iii)] $(\mathrm{S}_n)$
\item[(iv)] r\'egulier
\item[(v)] $(\mathrm{R}_n)$
\item[(vi)] normal
\item[(vii)] r\'eduit
\item[(viii)] de dimension $n$.
\end{itemize}
Alors, pour que $X$ poss\`ede la propri\'et\'e $P$, il faut et
il suffit qu'il en soit ainsi de~$X^{\an}$.
\end{proposition}

Soit $\varphi\colon X^{\an}\to X$ le morphisme
canonique. (i)~r\'esulte du fait que l'on a $|X^{\an}|=X(\CC )$
\eqref{XIII.1.1} et du fait que $X$ est de Jacobson (EGA~IV~10.4.8). Dire que~$X$
(\resp $X^{\an}$) est discret \'equivaut \`a dire que l'on a $\dim
X=0$ (\resp $\dim X^{\an}=0$ d'apr\`es \cite[\no19~\S 4~cor\ptbl 6]{XII.4}); (i')
r\'esulte donc de~(viii).

Soit $P$ l'une des propri\'et\'es (ii) \`a (vii). Pour que $X$
poss\`ede la propri\'et\'e $P$, il faut et il suffit que $P$
soit v\'erifi\'ee en chaque point ferm\'e de~$X$; en effet, $X$
\'etant excellent (EGA~IV~7.8.6~(iii)), l'ensemble des points
\marginpar{317}
o\`u $X$ v\'erifie $P$ est un ouvert (\loccit) et, si cet ouvert
contient tous les points ferm\'es, il est \'egal \`a $X$ tout
entier. Dire que $X$ (\resp $X^{\an}$) a la propri\'et\'e $P$
\'equivaut donc \`a dire que, pour tout point $x$ de~$X^{\an}$,
l'anneau local $\cal{O}_{X,\varphi (x)}$ (\resp $\cal{O}_{X^{\an},x}$)
a la propri\'et\'e $P$. Comme le fait qu'un anneau local excellent
ait la propri\'et\'e $P$ se voit apr\`es passage au
compl\'et\'e, la proposition r\'esulte des isomorphismes
$\widehat{\cal{O}}_{X,\varphi(x)}\isomto\widehat{\cal{O}}_{X^{\an},x}$
dans les cas (ii) \`a (vii). Il en est de m\^eme dans le cas
(viii), compte tenu des relations
$$
\dim X=\sup_x~\dim \cal{O}_{X,\varphi (x)}\qquad \dim
X^{\an}=\sup_x~\dim \cal{O}_{X^{\an},x}
$$
o\`u $x\in X^{\an}$. Ceci ach\`eve la d\'emonstration.

\begin{proposition}
\label{XII.2.2}
Soient $X$ un $\CC$-sch\'ema localement de type fini,
$\varphi\colon X^{\an}\to X$ le morphisme canonique, $T$ une partie
localement constructible de~$X$. Alors on a la relation
$$
\varphi^{-1}(\overline{T})=\overline{\varphi^{-1}(T)}.
$$
\end{proposition}

On peut supposer que $T$ est un ouvert dense de~$X$. Soit $H$ le
sous-sch\'ema ferm\'e r\'eduit de~$X$ d'espace sous-jacent
$X-T$; l'espace associ\'e $H^{\an}$ est un sous-espace analytique
ferm\'e de~$X^{\an}$ d'espace sous-jacent
$X^{\an}-\varphi^{-1}(T)$. On doit montrer que tout point $x$ de
$H^{\an}$ appartient \`a $\overline{\varphi^{-1}(T)}$. Or, en un tel
point $x$, le germe d'espace analytique $(X^{\an},x)$ contient le
sous-germe $(H^{\an},x)$, et celui-ci est d\'efini par un Id\'eal
non nilpotent de~$\cal{O}_{X^{\an},x}$. Il r\'esulte alors du
Nullstellensatz \cite[\no19~\S 4~cor\ptbl 3]{XII.4} que tout voisinage ouvert de~$x$
contient des points de~$X^{\an}$ qui n'appartiennent pas \`a
$H^{\an}$, ce qui prouve bien que l'on a
$x\in\overline{\varphi^{-1}(T)}$.

\begin{corollaire}
\label{XII.2.3}
Soient
\marginpar{318}
$X$ un $\CC$-sch\'ema localement de type fini, $\varphi
\colon X^\an \to X$ le morphisme canonique, $T$ une partie localement
constructible de~$X$. Pour que $T$ soit une partie ouverte (\resp une
partie ferm\'ee, \resp une partie dense), il faut et il suffit qu'il
en soit ainsi de~$\varphi^{-1} (T)$.
\end{corollaire}

Le corollaire r\'esulte de~\Ref{XII.2.2} et du fait que, $X$
\'etant un sch\'ema de Jacobson (EGA~IV~10.4.8), deux parties
localement constructibles de~$X$ qui ont m\^eme trace sur l'ensemble
tr\`es dense $X (\CC )$ sont \'egales.

\begin{proposition}
\label{XII.2.4} Soit $X$ un $\CC$-sch\'ema localement de type
fini. Pour que $X$ soit connexe (\resp irr\'eductible), il faut et
il suffit qu'il en soit ainsi de~$X^\an$.
\end{proposition}

Supposons $X^\an $ connexe (\resp irr\'eductible). L'image $X (\CC)
$ de~$X^\an$ dans $X$ est alors connexe (\resp irr\'eductible). Il
en r\'esulte que $X$ est connexe (\resp irr\'eductible) car les
parties ferm\'ees de~$X$ et $X (\CC)$ se correspondent bijectivement
(EGA~IV~10.1.2).

Inversement supposons $X$ connexe (\resp irr\'eductible), et
montrons qu'il en est de m\^eme de~$X^\an$. On peut se borner au cas
o\`u $X$ est irr\'eductible. Supposons en effet $X$ connexe.
\'Etant donn\'e un point $x$ de~$X$, l'ensemble des points $y \in X$
tels qu'il existe une suite finie de sous-sch\'emas ferm\'es
irr\'eductibles $X_1, \dots, X_n$ de~$X$, avec $x \in X_1$, $y \in
X_n$, $X_i \cap X_{i+1} \neq \emptyset $ pour $1 \leq i \leq n-1$,
est un ensemble \`a la fois ouvert et ferm\'e, donc \'egal \`a
$X$ tout entier. Pour une suite $X_1, \dots, X_n$ telle que
pr\'ec\'edemment, on a aussi $X_{i}^\an \cap X_{i+1}^\an \neq
\emptyset $ pour $1 \leq i \leq n-1 $; si l'on suppose
d\'emontr\'e que les $X_{i}^\an$ sont connexes, il en est alors de
m\^eme de~$X^\an$.

On
\marginpar{319}
suppose d\'esormais $X$ irr\'eductible. On peut supposer de plus
$X$ affine. En effet, si $(U_i)_{i \in I}$ est un recouvrement de~$X$
par des ouverts affines, deux de ces ouverts ont une intersection non
vide, et la m\^eme propri\'et\'e est donc vraie pour le
recouvrement $(U_{i}^\an)_{i \in I} $ de~$X^\an$; si l'on suppose
d\'emontr\'e que les $U_{i}^\an $ sont irr\'eductibles, il en
est alors de m\^eme de~$X^\an$.

On peut supposer de plus que $X$ est normal. Soit en effet $\tilde X$
le normalis\'e de~$X$; comme le morphisme $\tilde X \to X$ est
surjectif, il est de m\^eme de~$\tilde X^\an \to X^\an$, ce qui
prouve que, si $\tilde X^\an $ est irr\'eductible, il est de
m\^eme de~$X^\an$.

On suppose d\'esormais $X$ affine normal. Comme les anneaux locaux
de~$X^\an$ sont int\`egres, il revient au m\^eme de dire que
$X^\an$ est irr\'eductible ou qu'il est connexe. En effet, si
$\cal{F}$ est une partie analytique ferm\'ee de~$X^\an$, l'ensemble
des points $x$ de~$X^\an$ o\`u l'on a $\codim_x (\cal{F}, X^\an ) =
0$ est un sous-ensemble analytique ferm\'e de~$X^\an$ \cite[\no 20 A
Cor\ptbl 1]{XII.4} qui est aussi ouvert; si $X^\an$ est connexe, ceci prouve que,
si l'on a $\cal{F} \neq X^\an$, $\cal{F}$ est rare, donc que $X^\an$
est irr\'eductible. On est ainsi ramen\'e \`a montrer que
$X^\an$ est connexe.

Soit
$$
i \colon X \to P
$$
une compactification de~$X$, o\`u $P$ est un $\CC$-sch\'ema
projectif normal et $i$ une immersion ouverte dominante. Il
r\'esulte alors de \cite[\no 12 th\ptbl 1]{XII.10} que $P^\an$ est connexe. Comme
$X^\an$ est obtenu en enlevant \`a $P^\an$ une partie analytique
ferm\'ee rare, il r\'esulte de~\Ref{XII.2.5} ci-dessous que $X^\an$
est aussi connexe.

\begin{lemme}
\label{XII.2.5}
Soient
\marginpar{320}
$\cal{P} $ un espace analytique normal connexe, $\cal{Y}$ une
partie analytique ferm\'ee rare, alors $\cal{X} = \cal{P} - \cal{Y}$
est connexe.
\end{lemme}

Lorsque $\cal{Y}$ est de codimension $\geq 2$, la proposition
r\'esulte de \cite[\no 3 prop\ptbl 4]{XII.11}. Dans le cas g\'en\'eral on peut
supposer, quitte \`a enlever \`a $\cal{P}$ une partie analytique
ferm\'ee de codimension $\geq 2$, que $\cal{P}$ et $\cal{Y}$
(consid\'er\'e comme sous-espace analytique r\'eduit de~$\cal{P}$) sont r\'eguliers. D'apr\`es le th\'eor\`eme des
fonctions implicites, tout point $y$ de~$\cal{Y}$ poss\`ede un
voisinage $\cal{U}$ isomorphe \`a une boule d'un espace affine
$\cal{E}^{n}$, de sorte que $\cal{U} \cap \cal{Y}$ soit d\'efini par
l'annulation d'un certain nombre de fonctions coordonn\'ees. Ceci
prouve que $\cal{U} - \cal{U} \cap \cal{Y}$ est connexe, et il en est
donc de m\^eme de~$\cal{X}$.

\begin{corollaire}
\label{XII.2.6} Soit $X$ un $\CC$-sch\'ema localement de type
fini; le morphisme
$$
\pi_0 (X^{\an}) \to \pi_0 (X)
$$
induit par le morphisme canonique $X^\an \to X$ est bijectif.
\end{corollaire}

\section{Comparaison des propri\'et\'es des morphismes}
\label{XII.3}

\begin{proposition}
\label{XII.3.1} Soient $f \colon X \to Y$ un morphisme de
$\CC$-sch\'emas localement de type fini, $f^\an \colon X^\an \to
Y^\an$ le morphisme d\'eduit de~$f$ sur les espaces analytiques
associ\'es. Soit $P$ la propri\'et\'e d'\^etre
\begin{enumerate}
\item[(i)] plat
\item[(ii)] net (\ie non ramifi\'e)
\item[(iii)] \'etale
\item[(iv)] lisse
\item[(v)] normal
\item[(vi)]
\marginpar{321}
r\'eduit
\item[(vii)] injectif
\item[(viii)] s\'epar\'e
\item[(ix)] un isomorphisme
\item[(x)] un monomorphisme
\item[(xi)] une immersion ouverte.
\end{enumerate}
Alors, pour que $f$ poss\`ede la propri\'et\'e $P$, il faut et
il suffit qu'il en soit ainsi de~$f^\an$.
\end{proposition}

Notons $\varphi \colon X^\an \to X $ et $\psi \colon Y^\an \to Y$ les
morphismes canoniques. Soient $x$ un point de~$X^\an$, $y = f^\an
(x)$. Les morphismes $\cal{O}_{Y^\an, y} \to {\mathcal O}_{X^\an, x}$
et $\cal{O}_{Y, \psi(y)} \to \cal{O}_{X, \varphi(x)} $ d\'eduits de
$f$ et $f^{\an}$ donnent le m\^eme morphisme par passage aux
compl\'et\'es \eqref{XII.1.1}. D'apr\`es \cite[ch\ptbl 3 \S 5 prop\ptbl 4]{XII.2}
(\resp EGA~IV 17.4.4) il revient donc au m\^eme de dire que $f^\an $
v\'erifie la propri\'et\'e (i) (\resp (ii)) ou de dire que $f$
v\'erifie (i) (\resp (ii)) en chaque point ferm\'e de~$X$. Comme
l'ensemble des points de~$X$ o\`u (i) (\resp (ii)) est
v\'erifi\'e est un ouvert (EGA~IV 11.1.1 et I 3.3), ceci
d\'emontre (i) et (ii), donc aussi (iii).

\ifthenelse{\boolean{orig}}{}
{\enlargethispage{.5cm}}%
Soit $P$ la propri\'et\'e (iv) (\resp (v), \resp (vi)). Compte tenu
de~\Ref{XII.2.1} ((v), (vi), (vii)), il revient au m\^eme de dire que
les fibres g\'eom\'etriques de~$f^\an$ aux diff\'erents points
$y$ de~$Y^\an$ sont r\'eguli\`eres (\resp normales, \resp r\'eduites) ou qu'il en est ainsi des fibres g\'eom\'etriques de
$f$ aux diff\'erents points ferm\'es $\psi (y)$ de~$Y$. Les cas
(iv) (\resp (v), \resp (vi)) r\'esultent alors de (i) et du fait que
l'ensemble des points de~$Y$ o\`u les fibres g\'eom\'etriques de
$f$ sont r\'eguli\`eres est un ouvert (EGA~IV 12.1.7).

(vii). Si $f$ est injectif, il en est de m\^eme de~$f^\an$.
Inversement supposons $f^\an$ injectif et montrons qu'il en est de
m\^eme de~$f$. On peut supposer
\marginpar{322}
$f$ de type fini. Le morphisme $f^\an$ \'etant injectif, les fibres
de~$f$ aux points ferm\'es de~$Y$ sont radicielles; comme l'ensemble
des points de~$Y$ dont la fibre est radicielle est localement
constructible (EGA~IV 9.6.1) et comme $Y$ est un sch\'ema de
Jacobson, $f$ a toutes ses fibres radicielles donc est injectif.

(viii). Soient $\Delta \colon X \to X \times_Y X $ et $\Delta^\an
\colon X^\an \to X^\an \times_{Y^\an} X^\an $ les immersions
diagonales, $\Theta \colon X^\an \times_{Y^\an} X^\an \to X \times_Y X
$ le morphisme canonique. En vertu de~\Ref{XII.2.3} il revient au
m\^eme de dire que $\Delta (X)$ est ferm\'e dans $X \times_Y X $
ou que $\Delta^\an (X^\an)$ est ferm\'e dans $X^\an \times_{Y^\an}
X^\an$.

Comme une immersion ouverte n'est autre qu'un morphisme \'etale
injectif (EGA~IV 17.9.1 et \cite[\no 13 \S 1]{XII.4}), (xi) r\'esulte de (iii)
et de (vii). Un isomorphisme \'etant la m\^eme chose qu'une
immersion ouverte surjective, (ix) r\'esulte de (xi) et de
(\Ref{XII.3.2}~(i)) ci-dessous. Dire que $f$ est un monomorphisme
\'equivaut \`a dire que la morphisme diagonal $\Delta \colon X \to
X \times_Y X$ est un isomorphisme, donc (x) r\'esulte de (ix).

\begin{proposition}
\label{XII.3.2}
Soient $X$ et $Y$ deux $\CC$-sch\'emas localement de type fini,
$f\colon X \to Y$ un morphisme de type fini, $f^\an \colon X^\an \to
Y^\an $ le morphisme d\'eduit de~$f$ sur les espaces analytiques
associ\'es. Soit $P$ la propri\'et\'e d'\^etre
\begin{enumerate}
\item[(i)] surjectif
\item[(ii)] dominant
\item[(iii)] une immersion ferm\'ee
\item[(iv)] une immersion
\item[(v)] propre\footnote{Nous dirons qu'un morphisme d'espaces
analytiques est propre s'il l'est au sens de \cite[ch\ptbl 1 \S 10 \no 1]{XII.1} et
s'il est s\'epar\'e.}
\item[(vi)]
\marginpar{323}
fini.
\end{enumerate}
Alors, pour que $f$ poss\`ede la propri\'et\'e $P$, il faut et
il suffit qu'il en soit ainsi de~$f^\an$.
\end{proposition}

Soient $\varphi \colon X^\an \to X$ et $\psi \colon Y^\an \to Y$ les
morphismes canoniques.

(i). Si $f$ est surjectif, pour tout point $y$ de~$Y^\an$, $f^{-1}
(\psi (y)) $ est une partie ferm\'ee non vide de~$X$; elle contient
donc au moins un point ferm\'e, ce qui prouve que $f^\an$ est
surjectif. Inversement, si $f^\an$ est surjectif, $f(X)$ est une
partie localement constructible de~$Y$ (EGA~IV 1.8.4) qui contient
tous les points ferm\'es de~$Y$; on a donc $f(X)=Y$.

(ii) r\'esulte de~\Ref{XII.2.2}.

(iii). Si $f$ est une immersion ferm\'ee, il en est de m\^eme de
$f^\an$ d'apr\`es \Ref{XII.1.1}~a). Inversement, si $f^\an$ est une
immersion ferm\'ee, il en est de m\^eme de~$f$ d'apr\`es
\Ref{XII.3.1}~(x) et \Ref{XII.3.2}~(v), car cela revient \`a dire
que $f$ est un monomorphisme propre (EGA~IV 8.11.5).

(iv). Il est clair que, si $f$ est une immersion, il en est de
m\^eme de~$f^\an$. Inversement supposons que $f^\an$ soit une
immersion, et soient $T$ l'image de~$X$ dans $Y$, $\overline{T}$
l'adh\'erence sch\'ematique de~$f$. On a une factorisation de~$f$,
$$
X \lto{i} \overline{T} \lto{j} J \quoi,
$$
o\`u $j$ est une immersion ferm\'ee, $i$ le morphisme canonique,
et on en d\'eduit la factorisation suivante de~$f^\an$
$$
X^\an \lto{i^\an} \overline{T}^{\mkern1mu \an} \lto{j^\an} Y^\an \
.
$$
Comme
\marginpar{324}
$T = f(X)$ est une partie localement constructible de~$Y$ (EGA~IV
1.8.4), on~a, d'apr\`es~\Ref{XII.2.2}, $\overline{T}^{\mkern1mu \an} =
\overline{f^\an (X^\an )}$. Il en r\'esulte que $i^\an (X^\an )$ est
un ouvert de~$\overline{T}^{\mkern1mu \an}$, donc que $i(X) $ est un ouvert de
$\overline{T}$. On consid\`ere la factorisation canonique de~$i$
$$
X \lto{i_1} i (X) \lto{i_2} \overline{T}.
$$
Le morphisme $i_1^\an $ est un monomorphisme propre, donc il en est de
m\^eme de~$i_1$ d'apr\`es \Ref{XII.3.2}~(v) et \Ref{XII.3.1}~(x);
ceci prouve que $i_1$ donc aussi $f$ est une immersion.

(v). Supposons que $f$ soit propre et montrons qu'il en est de
m\^eme de~$f^\an$. Le fait que $f^\an $ soit propre \'etant local
sur~$Y^\an$, on peut supposer $Y$ affine. D'apr\`es le lemme de Chow
(EGA~II 5.6.1), on peut trouver un $Y$-sch\'ema projectif $X'$ et un
morphisme projectif surjectif
$$
g \colon X' \to X.
$$
Le morphisme $(fg)^\an = f^\an g^\an $ est projectif donc propre,
$g^\an $ est surjectif, et il r\'esulte de \cite[ch\ptbl 1 \S 10]{XII.1} que
$f^\an $ est propre.

Inversement supposons $f^\an $ propre et montrons qu'il en est de
m\^eme de~$f$. D'apr\`es \Ref{XII.3.1}~(viii) $f$ est
s\'epar\'e. Il reste \`a prouver que $f$ est universellement
ferm\'e, et il suffit m\^eme de montrer que $f$ est ferm\'e; en
effet, pour tout $Y$-sch\'ema $Y'$ localement de type fini, le
morphisme
$$
f_{(Y')} = h \colon X \times_Y Y' \to Y'
$$
sera aussi ferm\'e puisque $h^\an $ est propre. Soit $T$ une partie
ferm\'ee de~$X$; $f(T)$ est un ensemble localement constructible, et
l'on a
$$
f^\an (\varphi^{-1} (T)) = \psi^{-1} (f(T)).
$$
Comme
\marginpar{325}
$f^\an $ est propre, $\psi^{-1} (f (T))$ est une partie ferm\'ee de
$Y^\an $, et il r\'esulte donc de~\Ref{XII.2.2} que l'on a
$$
\psi^{-1} (\overline{f(T)}) = \psi^{-1} (f(T))\quoi .
$$
Cela entra\^ine que l'on a $f(T) = \overline{f(T)}$, \ie que $f$
est ferm\'ee donc que $f$ est propre.

(vi). Il revient au m\^eme de dire qu'un morphisme est fini ou qu'il
est propre \`a fibres finies (EGA~III 4.4.2 et \cite[\no 19 \S 5]{XII.4}).
Comme l'ensemble des points o\`u les fibres de~$f$ sont finies est
localement constructible (EGA~IV 9.7.9), les fibres de~$f$ sont finies
si et seulement si il en est ainsi des fibres de~$f^\an $; (vi)
r\'esulte donc de (v).

\begin{remarque}
\label{XII.3.3}
\begin{enumerate}
\item[a)] soit $f \colon X \to Y$ un morphisme de~$\CC$-sch\'emas
localement de type fini. Le fait que $f^\an $ soit un isomorphisme
local n'entra\^ine pas qu'il en soit de m\^eme de~$f$. En effet,
si $f$ est \'etale, $f^\an $ est \'etale donc est un isomorphisme
local \hbox{\cite[\no 13 \S 1]{XII.4}}, mais il n'en est pas n\'ecessairement ainsi
de~$f$.
\item[b)]
\ifthenelse{\boolean{orig}}
{l'\'enonc\'e~\Ref{XII.3.2}}
{L'\'enonc\'e~\Ref{XII.3.2}}
n'est pas vrai si l'on ne suppose
pas $f$ de type fini. Montrons par exemple que $f^\an $ peut \^etre
une immersion ferm\'ee sans qu'il en soit de m\^eme de~$f$. Il
suffit en effet de prendre pour $X$ la somme de~$\ZZ $ copies de
$\Spec \CC$, et pour $Y$ la droite affine, et pour $f$ le morphisme
obtenu en envoyant les points de~$X$ sur des points distincts de~$Y$
formant une partie discr\`ete.
\end{enumerate}
\end{remarque}

\section[Th\'eor\`emes de comparaison cohomologique]{Th\'eor\`emes de comparaison cohomologique et th\'eor\`emes d'existence}
\label{XII.4}
\marginpar{326}

L'objet de ce num\'ero est de red\'emontrer les r\'esultats de
\cite[\no 2 th\ptbl 5 et th\ptbl 6]{XII.3}; ces derniers g\'en\'eralisent au cas
d'un sch\'ema propre les th\'eor\`emes \'etablis dans \cite[\no 12]{XII.10} lorsque $X$ est projectif, et les \'etendent au cas
relatif. Des r\'esultats plus g\'en\'eraux, concernant les
sch\'emas relatifs propres sur un espace analytique, sont
prouv\'es dans \cite[ch\ptbl VIII \no 3]{XII.7}.

Rappelons que la cohomologie de {\v{C}}ech utilis\'ee dans \cite[\no 12]{XII.10} co\"incide avec la cohomologie usuelle dans le cas
alg\'ebrique comme dans le cas analytique (EGA~III 1.4.1. et \cite[II 5.10]{XII.5}).
\subsection{}
\label{XII.4.1}
Soient $f \colon X \to Y $ un morphisme de~$\CC$-sch\'emas
localement de type fini et consid\'erons le diagramme commutatif
$$
\xymatrix{ X^\an \ar[r]^\varphi \ar[d]_{f^\an} & X \ar[d]^f
\\ Y^\an \ar[r]^\psi & Y \quoi.}
$$
Si $F$ est un $\cal{O}_X$-Module, on a, pour tout entier $p \geq
0$, des morphismes
$$
\R^p f_* F \lto{i} \R^p f_* (\varphi_* F^\an )
\lto{j} \R^p (f\cdot\varphi)_* F^\an \lto{k} \psi_*
(\R^p f_*^\an F^\an ) \quoi,
$$
o\`u $i$ se d\'eduit du morphisme canonique $F \to \varphi_*
F^\an $, et $j, k$ sont des \og edge-homomorphismes\fg de suites
spectrales de Leray. Au compos\'e $k.j.i$ est associ\'e un
morphisme canonique
\begin{equation*}
\label{eq:XII.4.1.1}
\tag{4.1.1} {\theta_p \colon (\R^p f_* F)^\an \to \R^p f_*^\an
(F^\an )}
\end{equation*}

\begin{theoreme}
\label{XII.4.2}
Soient
\marginpar{327}
$f \colon X \to Y $ un morphisme propre de~$\CC$-sch\'emas
localement de type fini, $F$ un $\cal{O}_X$-Module
coh\'erent. Alors, pour tout entier $p \geq 0$ le morphisme
\eqref{eq:XII.4.1.1}
$$
\theta_p \colon (\R^p f_* F)^\an \to \R^p f_*^\an (F^\an)
$$
est un isomorphisme.
\end{theoreme}

1) \emph{Cas o\`u $f$ est projectif}. La d\'emonstration est
analogue \`a celle de \cite[\no 13]{XII.10}. Rappelons-la bri\`evement. On
se ram\`ene au cas o\`u $X$ est un espace projectif type $\PP^r_Y
$ au-dessus de~$Y$. Soit $\cal{Y} = Y^\an$, $\cal{P} = \PP^r_\cal{Y}$;
on prouve d'abord que l'on a
$$
f_*^\an \cal{O}_\cal{P} = \cal{O}_\cal{Y} \quoi, \qquad \R^p f_*^\an
(\cal{O}_\cal{P}) = 0 \qquad \text{pour } p>0
$$
Pour v\'erifier les relations pr\'ec\'edentes, on peut en effet
se ramener au cas o\`u $\cal{Y}$ est une boule $\cal{B}$ d'un espace
affine $\cal{E}^n$. On consid\`ere le \og recouvrement standard\fg $\{
\cal{U}_i \}$ de~$\cal{P}$ par $r+1$ ouverts isomorphes \`a $\cal{B}
\times \cal{E}^r$. Comme ces ouverts sont de Stein, on~a, pour tout
entier $p \geq 0$, des isomorphismes
$$
\H^p (\{ \cal{U}_i \}, \cal{O}_\cal{P}) \isomto \H^p (\cal{P},
\cal{O}_\cal{P})
$$
On peut alors exprimer les sections du faisceau structural
$\cal{O}_\cal{P}$ sur les ouverts $\cal{U}_i $ et sur leurs
intersections en termes de s\'eries de Laurent; un calcul facile
prouve que l'on~a
$$
\H^0 (\cal{P}, \cal{O}_\cal{P} ) \isomto \H^0 (\cal{Y},
\cal{O}_\cal{Y} ) \quoi, \qquad \H^p (\cal{P}, \cal{O}_\cal{P} ) = 0
\qquad \text{pour } p>0
$$
La d\'emonstration s'ach\`eve alors en recopiant \cite[\no 12
lemme~5]{XII.10}, les groupes de cohomologie \'etant remplac\'es par les
faisceaux de cohomologie.

2)
\marginpar{328}
\emph{Cas o\`u $f$ est propre}. On utilise EGA~III 3.1.2 pour se
ramener au cas projectif. Soit $\cal{K}$ la cat\'egorie des
$\cal{O}_X$-Modules coh\'erents tels que $\theta_p$ soit un
isomorphisme pour tout $p \geq 0$. Il suffit de prouver que, pour
toute suite exacte $0 \to F' \to F \to F'' \to 0$
dont deux termes sont dans $\cal{K}$, il en est de m\^eme du
troisi\`eme, qu'un facteur direct d'un objet de~$\cal{K}$ est dans
$\cal{K}$, et que, pour tout point $x$ de~$X$, on peut trouver un
objet $F$ de~$\cal{K}$ tel que l'on ait $F_x \neq 0$.

Le premi\`ere condition r\'esulte par application du lemme des
cinq du diagramme commutatif suivant dont les lignes sont exactes:
$$
\xymatrix@=.7cm{ \ar[r] & (\R^p f_* F')^\an \ar[d] \ar[r] &(\R^p f_*
F)^\an \ar[d] \ar[r] & (\R^p f_* F'')^\an \ar[d] \ar[r] &
(\R^{p+1} f_* F' )^\an \ar[r] \ar[d] &\, \\ \ar[r] & \R^p
f_*^\an F'^\an \ar[r] &\R^p f_*^\an F^\an \ar[r] & \R^p
f_*^\an F''^\an \ar[r] & \R^{p+1} f_*^\an F'^\an \ar[r] &
\quoi,}
$$
et on v\'erifie de fa\c{c}on analogue la deuxi\`eme condition.

Pour v\'erifier la troisi\`eme condition, on peut se borner au cas
o\`u $X$ est un sch\'ema irr\'eductible de point
g\'en\'erique $x$. On pouvait supposer $Y$ noeth\'erien d\`es
le d\'ebut. D'apr\`es le lemme de Chow (EGA~II 5.6.1), on peut
trouver un $Y$-sch\'ema projectif~$X'$ et un morphisme projectif
surjectif $g \colon X' \to X$. D'autre part il existe un entier~$n$
tel que l'on ait $\R^p g_* (\cal{O}_{X'} (n)) = 0 $ pour tout $p>0$ et
que le morphisme canonique $g^* g_* (\cal{O}_{X'} (n)) \to
\cal{O}_{X'} (n)$ soit surjectif (EGA~III 2.2.1). Si l'on pose
$F = g_* (\cal{O}_{X'} (n)) $, le faisceau $F$ r\'epond
\`a la question. En effet on a $F_x \neq 0$; de plus la suite
spectrale de Leray
$$
\R^p f_* (\R^q g_* (\cal{O}_{X'} (n))) \To \R^{p+q}
(\ifthenelse{\boolean{orig}}{f.g}{f\cdot g})_* (\cal{O}_{X'} (n))
$$
\'etant d\'eg\'en\'er\'ee, on a un isomorphisme
$$
\R^p f_* F \isomto \R^p (\ifthenelse{\boolean{orig}}{f.g}{f\cdot g})_* (\cal{O}_{X'} (n))
$$
\marginpar{329}%
Comme dans le cas alg\'ebrique on a un isomorphisme canonique
$$\R^p f_*^\an F^\an \isomto
\R^p (\ifthenelse{\boolean{orig}}{f.g}{f\cdot g})_*^\an (\cal{O}_{X'}(n)^\an ),$$
et le diagramme
$$
\xymatrix@C=1cm{
(\R^p f_* F)^\an \ar[d]_{\theta_p} \ar[r]^-{\sim} & (\R^p (\ifthenelse{\boolean{orig}}{f.g}{f\cdot g})_* (\cal{O}_{X'} (n)))^\an \ar[d]_{\psi_p} \\
\R^p f_*^\an F^\an \ar[r]^-{\sim} & \R^p (\ifthenelse{\boolean{orig}}{f.g}{f\cdot g})_*^\an
(\cal{O}_{X'}(n)^\an)}
$$
est commutatif. D'apr\`es~1) $\psi_p $ est un isomorphisme; il en
est donc de m\^eme de~$\theta_p$, ce qui ach\`eve la
d\'emonstration.

\begin{corollaire}
\label{XII.4.3}
Soient $X$ un $\CC$-sch\'ema propre, $F$ un $\cal{O}_X$-Module
coh\'erent. Alors, pour tout entier $p\geq 0$, le morphisme
canonique
$$
\H^p(X,F)\to \H^p(X^{\an},F^{\an})
$$
est un isomorphisme.
\end{corollaire}

\begin{theoreme}
\label{XII.4.4}
Soit $X$ un $\CC$-sch\'ema propre. Le foncteur qui, \`a tout
$\cal{O}_X$-Module coh\'erent $F$, associe son image inverse
$F^{\an}$ sur~$X^{\an}$ est une \'equivalence de cat\'egories.
\end{theoreme}

1) \emph{Le foncteur est pleinement fid\`ele}. Soient en effet $F$
et $G$ deux $\cal{O}_X$-Modules coh\'erents. Le morphisme canonique
$$
\Hom_{\cal{O}_X}(F,G)\to \Hom_{\cal{O}_{X^{\an}}}(F^{\an},G^{\an})
$$
s'identifie au morphisme canonique
$$
\H^0(X,\SheafHom_{\cal{O}_X}(F,G))\to
\H^0(X^{\an},\SheafHom_{\cal{O}_X}(F,G))
$$
\marginpar{330}%
(EGA~$0_{\textup{I}}$~6.7.6). Comme $\SheafHom_{\cal{O}_X}(F,G)$ est
coh\'erent, il r\'esulte de~\Ref{XII.4.3} que ce morphisme est
bijectif.

2) \emph{Le foncteur est essentiellement surjectif}. Lorsque $X$ est
projectif l'assertion r\'esulte de \cite[\no 12~th\ptbl 3]{XII.10}. Le cas
g\'en\'eral se ram\`ene au pr\'ec\'edent en utilisant le
lemme de Chow (EGA~II~5.6.1). Soient en effet $X'$ un $\CC$-sch\'ema
projectif, $f\colon X'\to X$ un morphisme projectif surjectif, $U$ un
ouvert dense de~$X$ tel que $f$ induise un isomorphisme
$f^{-1}(U)\simeq U$. On raisonne par r\'ecurrence noeth\'erienne
sur~$X$; on peut donc supposer que, pour tout faisceau coh\'erent
$\cal{G}$ sur~$X^{\an}$ tel que l'on puisse trouver une partie
ferm\'ee $Y$ de~$X$ distincte de~$X$, satisfaisant \`a la relation
$Y^{\an}\supset \Supp \cal{G}$, il existe un faisceau coh\'erent $G$
sur~$X$ tel que l'on ait un isomorphisme $G^{\an}\simeq \cal{G}$.

Soit $\cal{F}$ un faisceau de modules coh\'erent sur~$\cal{O}_{X^{\an}}$, $\cal{K}$ et $\cal{L}$ les faisceaux
coh\'erents d\'efinis par la condition que la suite
$$
0\to\cal{K}\to\cal{F}\to f_*^{\an}f^{\an *}\cal{F}\to\cal{L}\to 0
$$
soit exacte. Comme $X'$ est projectif, il existe un
$\cal{O}_{X'}$-Module coh\'erent $F'$ tel que l'on ait
$F'^{\an}\simeq {f^{\an}}^*\cal{F}$; on d\'eduit alors
de~\Ref{XII.4.2} que l'on a un isomorphisme $(f_*F')^{\an}\simeq
f_*^{\an}{f^{\an}}^*\cal{F}$. Comme $\cal{K}|U^{\an}$ et
$\cal{L}|U^{\an}$ sont nuls, il existe des $\cal{O}_X$-Modules
coh\'erents $K$ et $L$ tel que l'on ait des isomorphismes
$K^{\an}\simeq \cal{K}$, $L^{\an}\simeq \cal{L}$. D'apr\`es~1) le
morphisme $f_*^{\an}{f^{\an}}^*\cal{F}\to\cal{L}$ provient d'un unique
morphisme $f_*F'\to L$; soit $I=\Ker(f_*F'\to L)$. Le faisceau
$\cal{F}$ est alors extension de~$I^{\an}$ par~$K^{\an}$, et il suffit
de voir que cette extension provient
\marginpar{331}
par image inverse d'une extension de~$I$ par~$K$. Il suffit donc de
prouver que le morphisme canonique
\begin{equation*}
\label{eq:XII.4.*}
\tag{$*$} \Ext_{\cal{O}_X}^q(I,K)^{\an}\isomto
\Ext_{\cal{O}_{X^{\an}}}^q(I^{\an},K^{\an})\quad q\neq 1
\end{equation*}
est bijectif. Or on a des isomorphismes
$\mathbf{Ext}_{\cal{O}_X}^q(I,K)^{\an}\isomto
\mathbf{Ext}_{\cal{O}_{X^{\an}}}^q(I^{\an},K^{\an})$ pour tout entier
$q\geq 0$ (EGA~$0_{\textup{III}}$~12.3.5), et un morphisme de suites
spectrales
$$
\xymatrix@R=.7cm@C=.5cm{
\H^p(X,\mathbf{Ext}_{\cal{O}_X}^q(I,K))\ar@{=>}[r]\ar[d] &
\Ext_{\cal{O}_X}^{p+q}(I,K)\ar[d]\\
\H^p(X^{\an},\mathbf{Ext}_{\cal{O}_{X^{\an}}}^q(I^{\an},K^{\an}))\ar@{=>}[r] & \Ext_{\cal{O}_{X^{\an}}}^{p+q}(I^{\an},K^{\an}).
}
$$
Ce morphisme est un isomorphisme car, d'apr\`es~\Ref{XII.4.3}, il en
est ainsi sur les termes $E^{pq}_2$, et ceci d\'emontre la
bijectivit\'e de~$(*)$.

\begin{corollaire}
\label{XII.4.5}
Le foncteur qui \`a tout $\CC$-sch\'ema propre $X$ associe
$X^{\an}$ est pleinement fid\`ele.
\end{corollaire}

On doit montrer que, si $X$ et $Y$ sont deux $\CC$-sch\'emas
propres, l'application canonique
$$
\Hom_{\CC}(X,Y)\to\Hom(X^{\an},Y^{\an})
$$
est bijective. Or se donner un morphisme de~$X$ dans~$Y$ (\resp de
$X^{\an}$ dans $Y^{\an}$) \'equivaut \`a se donner son graphe,
\ie un sous-sch\'ema ferm\'e $Z$ de~$X\times Y$ (\resp un
sous-espace analytique ferm\'e $\othercal{Z}$ de~$X^{\an}\times
Y^{\an}$), tel que la restriction de la premi\`ere projection
$X\times Y\to X$ \`a~$Z$ (\resp de~$X^{\an}\times Y^{\an}\to
X^{\an}$ \`a~$\othercal{Z}$) soit
\marginpar{332}
un isomorphisme. Comme la donn\'ee d'un sous-sch\'ema ferm\'e
de~$X\times Y$ (\resp d'un sous-espace analytique ferm\'e
de~$X^{\an}\times Y^{\an}$) \'equivaut \`a celle d'un faisceau
coh\'erent d'id\'eaux sur~$\cal{O}_{X\times Y}$ (\resp sur~$\cal{O}_{X^{\an}\times Y^{\an}}$), le corollaire r\'esulte
de~\Ref{XII.4.4}.

\begin{corollaire}
\label{XII.4.6}
Soit $X$ un $\CC$-sch\'ema propre. Le foncteur qui, \`a tout
sch\'ema fini (\resp \'etale fini) $X'$ au-dessus de~$X$, associe
$X'^{\an}$ est une \'equivalence de la cat\'egorie des sch\'emas
finis (\resp \'etales finis) au-dessus de~$X$ dans la cat\'egorie
des espaces analytiques finis (\resp \'etales finis) au-dessus
de~$X^{\an}$.
\end{corollaire}

En effet se donner un morphisme fini $X'\to X$ (\resp $X'^{\an}\to
X^{\an}$) \'equivaut \`a se donner un faisceau coh\'erent
d'alg\`ebres sur~$\cal{O}_X$ (\resp sur~$\cal{O}_{X^{\an}}$)
\cite[\no 19~\S 5~th\ptbl 2]{XII.4}. Le corollaire r\'esulte donc de~\Ref{XII.4.4}
dans le cas non resp\'e, et le cas resp\'e s'en d\'eduit compte
tenu de~\Ref{XII.3.1}~(iii).

\section{Th\'eor\`emes de comparaison des rev\^etements \'etales}
\label{XII.5}
\setcounter{subsection}{-1}
\subsection{}
\label{XII.5.0}
Pr\'ecisons la notion de rev\^etement fini d'un espace analytique.
\index{revetement fini d'un espace analytique@rev\^etement fini d'un espace analytique|hyperpage}%
Si $\othercal{X}$ est un espace analytique, on dit qu'un espace analytique
$\othercal{X}'$ fini au-dessus de~$\othercal{X}$ est un rev\^etement fini de
$\othercal{X}$ si toute composante irr\'eductible de~$\othercal{X}'$ domine
une composante irr\'eductible de~$\othercal{X}$.


\begin{theoreme}[\og Th\'eor\`eme d'existence de Riemann\fg]
\label{XII.5.1}
\index{Riemann (theoreme d'existence de)@Riemann (th\'eor\`eme d'existence de)|hyperpage}%
\index{theoreme d'existence de Riemann@ th\'eor\`eme d'existence de Riemann|hyperpage}%
Soient $X$ un $\CC$-sch\'ema localement de type fini, $X^{\an}$
l'espace analytique associ\'e \`a $X$. Le foncteur $\Psi$ qui,
\`a tout rev\^etement \'etale fini $X'$ de~$X$, associe
$X^{\prime\an}$ est une \'equivalence
\marginpar{333}
de la cat\'egorie des rev\^etements \'etales finis de~$X$ dans
la cat\'egorie des rev\^etements \'etales finis de~$X^{\an}$.
\end{theoreme}

1) \emph{Le foncteur $\Psi$ est pleinement fid\`ele.} Soient $X'$
et $X''$ deux rev\^etements \'etales finis de~$X$, et prouvons que
l'application canonique
\begin{equation*}
\label{eq:XII.5.*}
\tag{$*$} { \Hom_{X}(X',X'') \to
\Hom_{X^{\an}}(X^{\prime\an},X^{\prime\prime\an}) }
\end{equation*}
est bijective. On peut supposer $X'$ connexe. Se donner un
$X$-morphisme de~$X'$ dans $X''$ \'equivaut \`a se donner une
composante connexe $X_{i}$ de~$X'\times_{X}X''$ telle que le morphisme
$X_{i}\to X'$ induit par la premi\`ere projection soit un
isomorphisme. Comme les composantes connexes de~$X'\times_{X}X''$
correspondent bijectivement aux composantes connexes de
$X^{\prime\an}\times_{X^{\an}}X^{\prime\prime\an}$~\eqref{XII.2.6} et qu'un
morphisme $X_{i}\to X'$ est un isomorphisme si et seulement si il en
est ainsi de~$X_{i}^{\an}\to X^{\prime\an}$, ceci d\'emontre la
bijectivit\'e de~\eqref{eq:XII.5.*}.

2) \emph{Le foncteur $\Psi$ est essentiellement surjectif.} Soit
$\othercal{X}'$ un rev\^etement \'etale fini de~$X^{\an}$ et prouvons
qu'il existe un rev\^etement \'etale $X'$ de~$X$ tel que l'on ait
un isomorphisme $X^{\prime\an}\isomto\othercal{X}'$. Compte tenu de~1) la
question est locale sur~$X$, et on peut donc supposer $X$ affine.

a) R\'eduction au cas o\`u $X$ est normal. On peut supposer $X$
r\'eduit. Supposons en effet le th\'eor\`eme d\'emontr\'e
pour $X_{\red}$. Le foncteur qui, \`a un rev\^etement \'etale
fini $X'$ de~$X$ fait correspondre le rev\^etement \'etale fini
$X^{\prime\an}_{\red}$ de~$X_{\red}^{\an}$ est alors une \'equivalence.
Comme il s'obtient en composant $\Psi$ avec le foncteur $\Theta$ qui,
\`a un rev\^etement \'etale fini de~$X^{\an}$ associe son image
inverse sur~$X_{\red}^{\an}$, et que $\Theta$ est pleinement
fid\`ele, ceci montre que $\Psi$ est une \'equivalence de
cat\'egories.

On
\marginpar{334}
peut supposer $X$ normal. Soit en effet $\widetilde{X}$ le
normalis\'e de~$X$, $p\colon\widetilde{X}\to X$ le morphisme
canonique. Comme $p$ est fini, $p$ est un morphisme de descente
effective pour la cat\'egorie des rev\^etements
\'etales~(IX~\Ref{IX.4.7}). Le th\'eor\`eme \'etant
suppos\'e d\'emontr\'e pour $\widetilde{X}$, si l'on pose
$\widetilde{\othercal{X}'}=\othercal{X}'\times_{X^{\an}}\widetilde{X}^{\an}$,
il existe un rev\^etement \'etale $\widetilde{X}'$ de
$\widetilde{X}$ et un isomorphisme $\widetilde{X}'{}^{\an}\simeq
\widetilde{\othercal{X}'}$. Il r\'esulte alors de~1) que la donn\'ee
de descente naturelle que l'on a sur~$\widetilde{\othercal{X}'}$ se
rel\`eve en une donn\'ee de descente sur~$\widetilde{X}'$
relativement \`a $\widetilde{X}\to X$; ceci prouve l'existence d'un
rev\^etement \'etale $X'$ de~$X$ tel que l'on ait un isomorphisme
$i\colon X^{\prime\an}\times_{X^{\an}}\widetilde{X}^{\an} \simeq
\widetilde{\othercal{X}'}$, dont les images inverses par les deux
projections de
$\widetilde{X}^{\an}\times_{X^{\an}}\widetilde{X}^{\an}$ soient les
m\^emes. D'apr\`es~IX~\Ref{IX.3.2}, dont la d\'emonstration est
valable dans le cas analytique, le morphisme $\widetilde{X}^{\an}\to
X^{\an}$ est un morphisme de descente pour la cat\'egorie des
rev\^etements \'etales, et par suite $i$ provient d'un
isomorphisme $X^{\prime\an}\simeq
\othercal{X}'$.

b) R\'eduction au cas o\`u $X$ est r\'egulier. Soient $U$
l'ouvert des points r\'eguliers de~$X$, $i\colon U\to X$,
$i^{\an}\colon U^{\an}\to X^{\an}$ les morphismes canoniques; comme
$X$ est normal, on a $\codim(X-U,X)\geq 2$. Supposons qu'il existe un
rev\^etement \'etale $U'$ de~$U$ tel que l'on ait
$U^{\prime\an}\simeq \othercal{X}'|U^{\an}$ et montrons qu'alors $U'$ se
prolonge en un rev\^etement \'etale $X'$ de~$X$ tel que l'on ait
$X^{\prime\an}\simeq \othercal{X}'$. Il suffit de voir que $U'$ se prolonge
en un rev\^etement \'etale $X'$ de~$X$; en effet on aura alors un
isomorphisme $X^{\prime\an}|U^{\an}\simeq \othercal{X}'|U^{\an}$; mais, si
$\cal{F}$ et $\cal{G}$ sont les faisceaux coh\'erents d'alg\`ebres
sur~$\cal{O}_{X^{\an}}$ d\'efinissant respectivement $\othercal{X}'$ et
$X^{\prime\an}$, le fait que $X$ soit normal et que l'on ait
$\codim(X-U,X)\geq 2$ entra\^ine que les morphismes canoniques
$$
\cal{F}\to i_*^{\an}(\cal{F}|U^{\an}) \qquad
\cal{G}\to i_*^{\an}(\cal{G}|U^{\an})
$$
sont
\marginpar{335}
des isomorphismes \cite[\no 3~prop\ptbl 4]{XII.11}. Il en r\'esulte que $\cal{F}$
et $\cal{G}$ donc aussi $X^{\prime\an}$ et $\othercal{X}'$ sont isomorphes.

Soit $\varphi\colon X^{\an}\to X$ le morphisme canonique. Comme le
probl\`eme de prolonger $U'$ \`a $X$ est local sur~$X$, il suffit
de prouver que, pour tout point $y$ de~$X^{\an}- U^{\an}$, le rev\^etement \'etale
$U'_{\varphi(y)}=U'\times_{X}\Spec\cal{O}_{X,\varphi(y)}$ de
$U_{\varphi(y)}=U\times_{X}\Spec\cal{O}_{X,\varphi(y)}$ se prolonge
\`a $\Spec\cal{O}_{X,\varphi(y)}$. Soit $H$ la
$\cal{O}_{U}$-Alg\`ebre coh\'erente d\'efinissant $U'$. Le
morphisme canonique
$$
\alpha\colon (i_*H)^{\an} \to i_*^{\an}(H^{\an})
=\cal{F}
$$
d\'efinit un morphisme de faisceaux de modules sur~$\Spec\cal{O}_{X^{\an},y}$:
$$
\alpha_{y}\colon (i_*H)^{\an}_{y} \to \cal{F}_{y}\quoi,
$$
dont la restriction \`a $U_{y}=U_{\varphi(y)}
\times_{\Spec\cal{O}_{X,\varphi(y)}}\Spec\cal{O}_{X^{\an},y}$ est un
isomorphisme. Mais ceci prouve que $H|U_{y}$ est trivial, donc que
$U'_{\varphi(y)}$ se prolonge \`a $\Spec\cal{O}_{X,\varphi(y)}$.

c) Cas o\`u $X$ est affine r\'egulier. Soit
$$
j\colon X\to P
$$
une compactification de~$X$, o\`u $P$ est un $\CC$-sch\'ema
projectif et $j$ une immersion ouverte dominante. Gr\^ace au
th\'eor\`eme de r\'esolution des singularit\'es~\cite{XII.8}, on peut
trouver un sch\'ema r\'egulier $R$, un morphisme projectif
$r\colon R\to P$, tel que $r$ induise un isomorphisme $r^{-1}(X)\simeq
X$ et que $r^{-1}(X)$ soit le compl\'ementaire dans $R$ d'un
diviseur \`a croisements normaux. Soit
\marginpar{336}
$$
k\colon X\to R
$$
l'immersion canonique. On va montrer qu'il existe un rev\^etement
fini normal~\eqref{XII.5.0} $\othercal{R}'$ de~$R^{\an}$ qui prolonge le
rev\^etement \'etale $X^{\prime\an}$. D'apr\`es la
proposition~\Ref{XII.5.3} ci-dessous, un tel rev\^etement est
unique; le probl\`eme de prolonger $X^{\prime\an}$ est donc local sur~$R^{\an}$ au voisinage de~$R^{\an}- X^{\an}$. Or chaque
point de~$R^{\an}- X^{\an}$ a un voisinage ouvert
$\othercal{V}$ isomorphe \`a une boule d'un espace affine
$\othercal{E}^{n}$, tel que $\othercal{V}- \othercal{V}\cap X^{\an}$
soit d\'efini par l'annulation des $p$ premi\`eres fonctions
coordonn\'ees $z_{1},\dots,z_{p}$, avec $0\leq p\leq n$. Le groupe
fondamental de~$\othercal{U}=\othercal{V}\cap X^{\an}$ est isomorphe \`a
$\ZZ^{p}$, et tout rev\^etement \'etale de~$\othercal{U}$ est quotient
d'un rev\^etement de la forme
$$
\othercal{U}''= \othercal{U}[T_{1},\dots,T_{p}]/
(T_{1}^{n_{1}}-z_{1},\dots,T_{p}^{n_{p}}-z_{p})\quoi,
$$
o\`u les $n_{i}$ sont des entiers $>0$, par un sous-groupe $H$ du
groupe de Galois $\ZZ/n_{1}\ZZ\times\cdots\times\ZZ/n_{p}\ZZ$ de
$\othercal{U}''$. Or $\othercal{U}''$ se prolonge en le rev\^etement
r\'egulier
$$
\othercal{V}''= \othercal{V}[T_{1},\dots,T_{p}]/
(T_{1}^{n_{1}}-z_{1},\dots,T_{p}^{n_{p}}-z_{p})\quoi,
$$
de~$\othercal{V}$ sur lequel $H$ op\`ere, et le quotient de~$\othercal{V}''$
par $H$ est le prolongement cherch\'e.

La d\'emonstration s'ach\`eve alors gr\^ace \`a~\Ref{XII.4.6}.
Le rev\^etement $\othercal{R}'$ provient d'un rev\^etement fini $R'$
de~$R$; la restriction de~$R'$ \`a $X$ est un rev\^etement $X'$ de
$X$ tel que l'on ait $X^{\prime\an}\simeq \othercal{X}'$, et
d'apr\`es~\Ref{XII.3.1}(iii) $X'$ est un rev\^etement \'etale
de~$X$.

\begin{corollaire}
\label{XII.5.2}
Soient
\marginpar{337}%
$X$ un $\CC$-sch\'ema localement de type fini connexe,
$\varphi\colon X^{\an}{\to} X$ le morphisme canonique, $x$ un point de
$X^{\an}$. Soit $\pi_{1}(X^{\an},x)$ le groupe fondamental de
l'espace topologique $X^{\an}$ au point $x$, $\pi_{1}(X,\varphi(x))$
le groupe fondamental du sch\'ema $X$ au point $\varphi(x)$
\textup{(V~\Ref{V.7})}. Alors $\pi_{1}(X,\varphi(x))$ est canoniquement
isomorphe au compl\'et\'e de~$\pi_{1}(X^{\an},x)$, pour la
topologie des sous-groupes d'indice fini.
\end{corollaire}

Soit en effet $\cal{C}$ la cat\'egorie des rev\^etements
\'etales finis de~$X^{\an}$, $F$ le foncteur de~$\cal{C}$ dans
$\Ens$ qui, \`a tout rev\^etement \'etale fini $\othercal{X}'$ de
$X^{\an}$ associe l'ensemble des points de~$\othercal{X}'$ au-dessus de
$x$, et soit $\widehat{\pi}_{1}(X^{\an},x)$ le groupe profini
associ\'e \`a $\cal{C}$ et $F$ comme il est dit dans~V~\Ref{V.4}.
Comme tout rev\^etement \'etale fini de~$X^{\an}$ est quotient du
rev\^etement universel par un sous-groupe d'indice fini,
$\widehat{\pi}_{1}(X^{\an},x)$ n'est autre que le compl\'et\'e de
$\pi_{1}(X^{\an},x)$ pour la topologie des sous-groupes d'indice fini.
Le corollaire r\'esulte donc de~\Ref{XII.5.1} et~V~\Ref{V.6.10}.

\begin{proposition}
\label{XII.5.3}
Soient $\othercal{X}$ un espace analytique normal, $\othercal{Y}$ un
sous-ensemble analytique ferm\'e tel que
$\othercal{U}=\othercal{X}-\othercal{Y}$ soit dense dans $\othercal{X}$.
Alors le foncteur qui, \`a tout rev\^etement normal
fini~\eqref{XII.5.0} $\othercal{X}'$ de~$\othercal{X}$ associe sa restriction
\`a $\othercal{U}$ est pleinement fid\`ele.
\end{proposition}

Soient $\othercal{X}'$ and $\othercal{X}''$ deux rev\^etements finis normaux
de~$\othercal{X}$. On doit montrer que l'application canonique
$$
\Hom_{\othercal{X}}(\othercal{X}',\othercal{X}'') \to
\Hom_{\othercal{U}}(\othercal{X}'|\othercal{U},\othercal{X}''|\othercal{U})
$$
est bijective. Soient $u$, $v$ deux $\othercal{X}$-morphismes de
$\othercal{X}'$ dans $\othercal{X}''$ dont les restrictions \`a $\othercal{U}$
sont les m\^emes et prouvons que $u=v$. Les morphismes $u$ et
\marginpar{338}
$v$ co\"incident sur l'ouvert dense
$\othercal{U}\times_{\othercal{X}}\othercal{X}'$, donc sur les espaces
topologiques sous-jacents. D'apr\`es \hbox{\cite[\no 19~\S 4~cor\ptbl 5]{XII.4}} ceci
prouve que l'on a $u=v$.

Soit maintenant $u$ un $\othercal{U}$-morphisme de~$\othercal{X}'|\othercal{U}$
dans $\othercal{X}''|\othercal{U}$ et montrons qu'il se prolonge \`a
$\othercal{X}'$ tout entier. On peut supposer $\othercal{X}'$ r\'egulier.
En effet, $\othercal{X}'$ \'etant normal, on peut trouver un ouvert
$\othercal{V}$ de~$\othercal{X}$ dont le compl\'ementaire soit une partie
analytique de codimension $\geq 2$, tel que
$\othercal{X}'\times_{\othercal{X}}\othercal{V}=\othercal{V}'$ soit r\'egulier.
Soit $\othercal{V}''=\othercal{X}''\times_{\othercal{X}}\othercal{V}$ et supposons la
proposition d\'emontr\'ee pour $\othercal{V}$. On consid\`ere le
diagramme commutatif
$$
\xymatrix{ \othercal{V}' \ar[dd]_{g'}\ar[rd]\ar[rrr]^{i'} && &\othercal{X}'
\ar[dd]^{\hbox{\raise10mm\hbox{$f'$}}} \cr &\othercal{V}'' \ar[ld]_{g''\mkern-10mu}\ar[rrr]^-{i''~~} & &
&\othercal{X}'' \ar[ld]_{f''\mkern-12mu} \cr \othercal{V} \ar[rrr]^{i} && &\othercal{X} \cr
}
$$
\`A $u$ est associ\'e un morphisme de
$\cal{O}_{\othercal{V}}$-Alg\`ebres $g''_*\cal{O}_{\othercal{V}''}\to
g'_*\cal{O}_{\othercal{V}'}$, d'o\`u l'on d\'eduit un morphisme
$$
i_*g''_*\cal{O}_{\othercal{V}''} \to
i_*g'_*\cal{O}_{\othercal{V}'}\quoi.
$$
Compte tenu des isomorphismes
$i'_*\cal{O}_{\othercal{V}'}\simeq\cal{O}_{\othercal{X}'}$,
$i''_*\cal{O}_{\othercal{V}''}\simeq\cal{O}_{\othercal{X}''}$
\cite[\no 3~prop\ptbl 4]{XII.11} on en d\'eduit un morphisme de
$\cal{O}_{\othercal{X}}$-Alg\`ebres
$$
f''_*\cal{O}_{\othercal{X}''} \to
f'_*\cal{O}_{\othercal{X}'}\quoi,
$$
d'o\`u le morphisme $\othercal{X}'\to\othercal{X}''$ cherch\'e.

On
\marginpar{339}
suppose d\'esormais $\othercal{X}'$ r\'egulier. Soient
$\othercal{U}'=\othercal{U}\times_{\othercal{X}}\othercal{X}'$,
$\othercal{Y}'=\othercal{X}'-\othercal{U}'$. On consid\`ere $\othercal{Y}'$ comme
sous-espace analytique r\'eduit de~$\othercal{X}'$; si $\othercal{Y}'_{1}$
est le ferm\'e singulier de~$\othercal{Y}'$, on a $\dim \othercal{Y}'_{1} <
\dim \othercal{Y}'$~\cite[\no 20~D\ptbl Th\ptbl 3]{XII.4}. On voit donc par r\'ecurrence
sur la dimension de~$\othercal{Y}'$ que l'on peut supposer $\othercal{Y}'$
lisse. Comme il suffit de prolonger $u$ \`a un voisinage ouvert de
chaque point de~$\othercal{Y}'$, on peut supposer par le th\'eor\`eme
des fonctions implicites que $\othercal{X}'$ est une boule d'un espace
affine $\othercal{E}^{n}$ et $\othercal{Y}'$ le ferm\'e d\'efini par
l'annulation des $p$ premi\`eres fonctions coordonn\'ees
$z_{1},\dots,z_{p}$, avec $0\leq p\leq n$.

On associe \`a $u$ une section $s$ de
$p\colon\othercal{X}'\times_{\othercal{X}}\othercal{X}''\to\othercal{X}'$ au-dessus de
$\othercal{U}'$; quitte \`a restreindre $\othercal{X}'$, on peut supposer
$p_*(\cal{O}_{\othercal{X}'\times_{\othercal{X}}\othercal{X}''})$
engendr\'e par des \'el\'ements $x_{1},\dots,x_{q}$ de
$\Gamma(\othercal{X}',
p_*\cal{O}_{\othercal{X}'\times_{\othercal{X}}\othercal{X}''})$; soient
$u_{1},\dots,u_{q}\in\Gamma(\othercal{U}',\cal{O}_{\othercal{X}'})$ les images
par $s$ de~${x_{1}}|\othercal{U}',\dots,{x_{q}}|\othercal{U}'$. Dire
que $s$ se prolonge \`a $\othercal{X}'$ revient \`a dire que les
$u_{1},\dots,u_{q}$ se prolongent en des sections de
$\Gamma(\othercal{X}',\cal{O}_{\othercal{X}'})$. Mais, puisque $f$ est fini,
chaque $u_{i}$ est une s\'erie de Laurent en $z_{1},\dots,z_{p}$,
\`a coefficients des s\'eries enti\`eres en
$z_{p{+}1},\dots,z_{n}$, qui satisfont \`a des relations de
d\'ependance int\'egrale. Il en r\'esulte que $u_{i}$ est
born\'e donc est une s\'erie enti\`ere en $z_{1},\dots,z_{n}$,
et par suite se prolonge \`a $\othercal{X}'$.

On peut se demander si le foncteur introduit dans~\Ref{XII.5.3} est
une \'equivalence de cat\'egories. On a une r\'eponse \`a
cette question gr\^ace au th\'eor\`eme de \textsc{Grauert}-\textsc{Remmert} \cite{XII.6}
dont nous donnons une d\'emonstration ci-dessous utilisant la
r\'esolution des singularit\'es. On aurait aussi pu utiliser le
th\'eor\`eme de \textsc{Grauert}-\textsc{Remmert} pour d\'emontrer~\Ref{XII.5.1};
c'est ce que l'on faisait avant de disposer de~\cite{XII.8}.

\ifthenelse{\boolean{orig}}
{}
{\enlargethispage{.5cm}}
\begin{theoreme}[Th\'eor\`eme de \textsc{Grauert}-\textsc{Remmert}]
\label{XII.5.4}
Soient
\marginpar{340}
\index{Grauert-Remmert (th\'eor\`eme de)|hyperpage}%
\index{theoreme de Grauert-Remmert@th\'eor\`eme de Grauert-Remmert|hyperpage}%
$\othercal{X}$ un espace analytique normal, $\othercal{Y}$ un
sous-ensemble analytique ferm\'e tel que
$\othercal{U}=\othercal{X}-\othercal{Y}$ soit dense dans~$\othercal{X}$. Soit
$\othercal{U}'$ un rev\^etement normal fini de~$\othercal{U}$; on suppose
qu'il existe une partie analytique ferm\'ee rare $\othercal{S}$
de~$\othercal{X}$ telle que la restriction de~$\othercal{U}'$ \`a
$\othercal{U}-\othercal{U}\cap\othercal{S}$ soit \'etale. Alors il existe un
rev\^etement fini normal $\othercal{X}'$ de~$\othercal{X}$ qui prolonge
$\othercal{U}'$, et $\othercal{X}'$ est unique \`a isomorphisme pr\`es.
\end{theoreme}

L'unicit\'e r\'esulte de~\Ref{XII.5.3}. Le probl\`eme de
prolonger $\othercal{U}'$ est donc local sur~$\othercal{X}$. On peut supposer
$\othercal{U}$ r\'egulier et $\othercal{U}'$ \'etale sur~$\othercal{U}$. En
effet l'ensemble des points r\'eguliers de~$\othercal{U}$ est un ouvert
$\othercal{V}$ dense dans~$\othercal{X}$ dont le compl\'ementaire est une
partie analytique \hbox{\cite[\no20~D th\ptbl 2]{XII.4}} et il suffit de remplacer
$\othercal{U}$ par l'ouvert $\othercal{V}-\othercal{V}\cap\othercal{S}$.

Soit $y$ un point de~$\othercal{X}-\othercal{U}$ et montrons que l'on peut
prolonger $\othercal{U}'$ \`a un voisinage de~$y$. Quitte \`a
restreindre $\othercal{X}$ \`a un voisinage ouvert de~$y$, il
r\'esulte du th\'eor\`eme de r\'esolution des singularit\'es
\cite{XII.8} que l'on peut trouver un espace analytique r\'egulier
$\othercal{X}_1$, un morphisme projectif $f\colon\othercal{X}_1\to\othercal{X}$
induisant par restriction \`a~$\othercal{U}$ un isomorphisme
$\othercal{U}_1=f^{-1}(\othercal{U})\simeq\othercal{U}$, tel que $\othercal{U}_1$ soit
le compl\'ementaire dans~$\othercal{X}_1$ d'un diviseur \`a
croisements normaux. Montrons que $\othercal{U}'$ se prolonge en un
rev\^etement fini normal de~$\othercal{X}_1$. Comme la question est
locale sur~$\othercal{X}_1$, on peut supposer que $\othercal{X}_1$ est une
boule d'un espace affine $\othercal{E}^n$ et que $\othercal{X}_1-\othercal{U}_1$
est d\'efini par l'annulation des $p$ premi\`eres fonctions
coordonn\'ees
\ifthenelse{\boolean{orig}}
{$z_1,\dots,z_q$,}
{$z_1,\dots,z_p$,}
avec $0\leq p\leq n$. Le rev\^etement \'etale $\othercal{U}'$
de~$\othercal{U}_1$ est quotient d'un rev\^etement de la forme
$$
\othercal{U}_2=\othercal{U}_1[T_1,\dots,T_p]/\left(T_1^{n_1}-z_1,\dots,
T_p^{n_p}-z_p\right)
$$
par un sous-groupe $H$ du groupe de Galois de~$\othercal{U}_2$. Le
rev\^etement $\othercal{U}_2$ se prolonge
\marginpar{341}
en le rev\^etement
$$
\othercal{X}_2=\othercal{X}_1[T_1,\dots,T_p]/\left(T_1^{n_1}-z_1,\dots,
T_p^{n_p}-z_p\right)
$$
de~$\othercal{X}_1$ sur lequel $H$ op\`ere, et $\othercal{X}_2/H$ prolonge
$\othercal{U}'$ \`a~$\othercal{X}_1$.

Notons $\othercal{X}_1'$ le rev\^etement fini normal de~$\othercal{X}_1$ qui
prolonge $\othercal{U}'$, $\othercal{F}_1$ la
$\cal{O}_{\othercal{X}_1}$-Alg\`ebre coh\'erente d\'efinie
par~$\othercal{F}_1$. D'apr\`es le th\'eor\`eme de finitude de
\textsc{Grauert}-\textsc{Remmert} \hbox{\cite[\no15 th\ptbl 1.1]{XII.4}}, $f_*\othercal{F}_1$ est une
$\cal{O}_{\othercal{X}}$-Alg\`ebre coh\'erente. Il lui correspond donc
un rev\^etement fini $\othercal{X}'$ de~$\othercal{X}$ qui est d'ailleurs
normal puisque $\othercal{X}_1'$ l'est, et $\othercal{X}'$ est le prolongement
de~$\othercal{U}'$ cherch\'e.

\begin{remarque}
\label{XII.5.5}
Dans l'\'enonc\'e~\Ref{XII.5.4}, on ne peut supprimer
l'hypoth\`ese sur le lieu des points o\`u le morphisme
$\othercal{U}'\to\othercal{U}$ n'est pas \'etale. Soit par exemple
$\othercal{X}$ le disque unit\'e du plan complexe, $\othercal{U}$ le
compl\'ementaire de l'origine dans~$\othercal{X}$,
$\othercal{U}'=\othercal{U}[T]/(T^2-\sin 1/z)$, o\`u~$z$ est la fonction
coordonn\'ee sur~$\othercal{X}$. Alors $\othercal{U}'$ est un rev\^etement
fini normal de~$\othercal{U}$ qui ne se prolonge pas
\`a~$\othercal{X}$. Supposons en effet que $\othercal{U}'$ se prolonge en un
rev\^etement fini $\othercal{X}'$ de~$\othercal{X}$; le lieu des points
de~$\othercal{X}$ o\`u le morphisme $\othercal{X}'\to\othercal{X}$ n'est pas
\'etale est alors un ferm\'e analytique qui contient tous les
points $z$ tels que l'on ait $\sin 1/z=0$, ce qui est absurde.

On peut cependant supprimer l'hypoth\`ese sur le lieu singulier du
morphisme $\othercal{U}'\to\othercal{U}$ lorsque l'on a
$\codim(\othercal{X}-\othercal{U},\othercal{X})\geq~2$. On peut en effet supposer
$\othercal{U}$ r\'egulier. Le lieu des points de~$\othercal{U}$ o\`u
$\othercal{U}'\to\othercal{U}$ n'est pas \'etale est un diviseur
de~$\othercal{U}$, et il r\'esulte du th\'eor\`eme de \textsc{Remmert}-\textsc{Stein}
\cite[th\ptbl 3]{XII.9} qu'il est
\marginpar{342}
la trace sur~$\othercal{U}$ d'un diviseur de~$\othercal{X}$. Or, dans ce cas,
si $\othercal{A}$ est une $\cal{O}_{\othercal{U}}$-Alg\`ebre coh\'erente
telle que $\othercal{U}'=\Spec \an(\othercal{A})$, si $i:\othercal{U}\to\othercal{X}$
est le morphisme canonique, il suffit de prendre
$\othercal{X}'=\Spec\an(i_*\othercal{A})$; on sait en effet que $i_*\othercal{A}$
est
\ifthenelse{\boolean{orig}}
{coh\'erents}
{coh\'erente}
\cite[\no 1 th\ptbl 1]{XII.11}.
\end{remarque}

\ifthenelse{\boolean{orig}}
{}
{\begin{remarqueMR}
\label{XII.5.6}
Il existe des groupes de pr\'esentation finie $G$, non
triviaux, qui ne poss\`edent pas de sous-groupes d'indice fini,
distincts de~$G$, par exemple le groupe de G\ptbl Higman (\cf J-P\ptbl Serre, Arbres et amalgames, Ast\'erisque N$^\circ$46, prop\ptbl 6, chap\ptbl I, \S 1).
D\`es lors, si un tel groupe se r\'ealise comme groupe fondamental
topologique d'un sch\'ema $V$ sur~$\CC$, disons lisse et
projectif, l'espace topologique $V_{\mathrm{top}}$ sous-jacent \`a $V$ n'est pas
simplement connexe, mais $V$ est alg\'ebriquement simplement connexe.
\`A l'heure actuelle, on ne conna\^it pas de tels $V$. Signalons toutefois
que D\ptbl Toledo a construit des sch\'emas lisses et projectifs $V$ sur~$\CC$ dont le groupe fondamental topologique n'est pas s\'epar\'e
pour la topologie des sous-groupes d'indice fini (le morphisme naturel
$\pi_1(V_{\mathrm{top}})\to\pi_1(V_{\alg})$ n'est pas injectif). [D\ptbl Toledo, Projective varieties with non-residually finite fundamental group, Publ.\ Math.\ IHES \textbf{77} (1993), p\ptbl103--119].
\end{remarqueMR}}


\chapterspace{-4}
\chapter[Propret\'e cohomologique des faisceaux d'ensembles...]{Propret\'e cohomologique\\ des faisceaux d'ensembles et\\ des faisceaux de groupes non commutatifs}
\label{XIII}
\marginpar{344}

\begin{center}
{par Mme M.\ \textsc{Raynaud}\footnote{D'apr\`es des notes in\'edites de
A\ptbl Grothendieck.}}
\end{center}
\vspace*{1cm}

Cet expos\'e se propose d'utiliser la cohomologie \'etale pour
g\'en\'eraliser certains r\'esultats de \Ref{IX} et \Ref{X}. Il montre aussi
comment on peut \'etendre aux faisceaux en groupes non
n\'ecessairement commutatifs les r\'esultats de~SGA~5~II qui ont
encore un sens pour de tels faisceaux. On suppose connues les notions
de cohomologie \'etale expos\'ees dans~SGA~4.

Le r\'esultat principal \eqref{XIII.2.4} donne un exemple important
de morphisme non propre $f\colon~U\to~S$, qui soit \og cohomologiquement
propre en dimension~$\leq1$\fg, c'est-\`a-dire tel que, pour
certains faisceaux en groupes $F$ sur~$U$ (au sens de la topologie
\'etale), la formation de~$f_*F$ et $\R^1f_*F$ commute \`a tout
changement de base~$S'\to~S$. Cette propri\'et\'e est en effet
satisfaite par l'ouvert $U$ d'un sch\'ema $X$ propre sur~$S$,
compl\'ementaire d'un diviseur $D$ \`a croisements normaux
relativement \`a~$S$, du moins si l'on impose \`a~$F$ d'\^etre
constant fini, d'ordre premier aux caract\'eristiques
r\'esiduelles de~$S$. Si l'on ne suppose plus $F$ d'ordre premier
aux caract\'eristiques r\'esiduelles de~$S$, on a un r\'esultat
analogue en rempla\c cant $\R^1f_*F$ par le sous-faisceau
$\R^1_{\tame}f_*F$ obtenu en se bornant \`a consid\'erer les torseurs
sous~$F$ \og mod\'er\'ement ramifi\'es sur~$X$ relativement
\`a~$S$\fg. En particulier cela permet de montrer que le groupe
fondamental mod\'er\'ement ramifi\'e d'une courbe alg\'ebrique
propre et lisse sur un corps s\'eparablement clos, priv\'ee d'un
nombre fini de points ferm\'es, est topologiquement de type fini
\eqref{XIII.2.12}.

Le
\marginpar{345}
\no \Ref{XIII.4} est consacr\'e \`a la suite exacte d'homotopie et \`a
la formule de K\"unneth.

Enfin un appendice donne des variantes utiles du lemme d'Abhyankar
d\'emontr\'e dans X.\Ref{X.3.6}.

\setcounter{section}{-1}

\section{Rappels sur la th\'eorie des champs}
\label{XIII.0}

Nous utiliserons dans ce qui suit la th\'eorie des champs
expos\'ee dans \cite{XIII.1} et \cite{XIII.2}. Nous nous bornons au cas du site
\'etale d'un sch\'ema. \'Etant donn\'e un sch\'ema $X$, notons
$X_\et$ le site \'etale de~$X$. Rappelons qu'un
\ifthenelse{\boolean{orig}}
{champ}
{champ $F$}
\index{champ|hyperpage}%
sur~$X$ est une cat\'egorie fibr\'ee au-dessus de~$X_\et$ telle
que, pour tout sch\'ema $X'$ \'etale sur~$X$ et pour tout couple
\ifthenelse{\boolean{orig}}
{d'objet}
{d'objets}
$x$, $y$ de la fibre $F_{X'}$, le pr\'efaisceau
$\SheafHom_{X'}(x,y)$ soit un faisceau, et telle que, pour tout
morphisme \'etale surjectif $X'' \to X'$, tout objet de~$F_{X''}$
muni d'une donn\'ee de descente relativement \`a $X'' \to X'$ soit
image inverse d'un objet de~$F_{X'}$.

On note $F(X')$ la cat\'egorie des sections cart\'esiennes de
$F/X'$. Plus g\'en\'eralement, si $\Sch_X$ est la cat\'egorie
des sch\'emas au-dessus de~$X$ munie de la topologie \'etale, le
champ $F$ peut s'\'etendre en un champ $\cal{F}$ sur~$\Sch_X$ et,
pour tout morphisme $f\colon X' \to X$, on note encore $F(X')$ la
cat\'egorie des sections cart\'esiennes de ce champ $\cal{F}$
au-dessus de~$X'$.

Une gerbe
\index{gerbe|hyperpage}%
est un champ tel que, pour tout sch\'ema $X'$ \'etale
sur~$X$ et pour tout couple d'objets $x$, $y$ de~$F_{X'}$, tout
morphisme de~$x$ dans $y$ soit un isomorphisme, que $x$ et $y$ soient
localement isomorphes, et tel que l'ensemble des objets $X'$ de
$X_\et$ tels que $F_{X'}$ soit non vide est un raffinement de
$X_\et$. Par exemple le champ des torseurs sous un faisceau en groupes
est une gerbe qui, de plus, a une section cart\'esienne.
R\'eciproquement une gerbe qui a une section, \ie telle qu'il
existe un objet $x$ de~$F_X$, est \'equivalente au champ des
torseurs sous le faisceau
\marginpar{346}
en groupes $\SheafAut_X(x)$.

On a une notion \'evidente de sous-gerbe et de sous-gerbe maximale
d'un champ~$F$. \'Etant donn\'e une section cart\'esienne $x$ de
$F(X)$, il existe une unique sous-gerbe maximale $G_x$ de~$F$ telle
que $x$ se factorise \`a travers~$G_x$. On appelle $G_x$ la
sous-gerbe engendr\'ee par~$x$; c'est par d\'efinition une gerbe
triviale.
\index{gerbe triviale|hyperpage}%
Le pr\'efaisceau $SF$
\label{indnot:mb}\oldindexnot{$SF$ ou $S(F)$|hyperpage}%
d\'efini par
$$
SF(X') = \{ \text{sous-gerbes maximales de }F|X' \}
$$
est un faisceau appel\'e, faisceau des sous-gerbes maximales
de~$F$.
\index{faisceau des sous-gerbes maximales|hyperpage}%
Soit $O$ le pr\'efaisceau d\'efini par
$$
O(X') = \{ \text{classes d'objets de }F_{X'}\text{ mod. isomorphisme}
\} \text{.}
$$
En associant \`a tout objet $x$ de~$F_{X'}$ la sous-gerbe maximale
de~$F|X'$, engendr\'ee par~$x$, on obtient un morphisme
$$
O \to SF \text{;}
$$
d'apr\`es \cite[III 2.1.4]{XIII.2}, ce morphisme fait de~$SF$ un faisceau
associ\'e \`a~$O$.

Un champ $F$ est dit \emph{constructible}
\index{constructible (champ)|hyperpage}%
(\emph{\resp ind-$\LL$-fini},
\index{champ ind-fini (\resp ind-$l$-fini)|hyperpage}%
$\LL$ \'etant un ensemble de nombres premiers) si, pour tout
sch\'ema $X'$ \'etale sur~$X$ et pour tout objet $x$ de~$F_{X'}$,
il en est ainsi du faisceau $\SheafAut_{X'}(x)$ \cite[VII 2.2.1]{XIII.2}. On dit
qu'un champ est $1$-constructible
\index{champ $1$-constructible (\resp constructible)|hyperpage}%
s'il est constructible et si le faisceau des sous-gerbes maximales est
constructible.

\section{Propret\'e cohomologique}
\label{XIII.1}
\setcounter{subsection}{-1}
\subsection{}
\label{XIII.1.0}
Soient $S$ un sch\'ema, $f \colon X \to Y$ un morphisme de
$S$-sch\'emas. Si $S'$ est un $S$-sch\'ema, on consid\`ere le
diagramme suivant, dont tous les carr\'es sont cart\'esiens:
\marginpar{347}
\begin{equation*}
\label{eq:XIII.1.0.1}
\tag{\thesubsection.1}
\begin{array}{c}
\xymatrix{ X \ar[d]_f & \ar[l]_h X'
\ar[d]^{f'} \\ Y \ar[d] & \ar[l]_g Y' \ar[d] \\ S & \ar[l] \,S'. }
\end{array}
\end{equation*}

Si $Y_1$ est un sch\'ema \'etale au-dessus de~$Y$, on pose $X_1 =
X \times_Y Y_1$, $Y_1' = Y' \times_Y Y_1$, et on consid\`ere le
carr\'e cart\'esien
\begin{equation*}
\label{eq:XIII.1.0.2}
\tag{\thesubsection.2}
\begin{array}{c}
\xymatrix{ X_1 \ar[d]_{f_1} & \ar[l]_{h_1}
X_1' \ar[d]^{f_1'} \\ Y_1 & \ar[l]^{g_1} \,Y_1'. }
\end{array}
\end{equation*}

\begin{definition}
\label{XIII.1.1}
Soit $F$ un champ sur~$X$. On dit que $(F,f)$ est
\emph{cohomologiquement propre relativement \`a~$S$ en dimension
$\leqslant -1$}
\index{cohomologiquement propre|hyperpage}%
\index{propret\'e cohomologique|hyperpage}%
(\resp en dimension $\leqslant 0$, \resp en dimension $\leqslant 1$)
si, pour tout $S$-sch\'ema $S'$, le foncteur canonique (d\'efini
de fa\c con \'evidente par la propri\'et\'e universelle de
l'image inverse de champs):
$$
g^* f_* F \to f'_* h^* F \qquad\text{ (\cf 1.0.1)}
$$
est fid\`ele (\resp pleinement fid\`ele, \resp une \'equivalence
de cat\'egories).
\end{definition}

S'il n'y a pas de confusion possible sur~$S$, en particulier si $S=Y$,
on dit cohomologiquement propre au lieu de cohomologiquement propre
relativement \`a~$S$.

\subsection{}
\label{XIII.1.2}
Soit $F$ un faisceau d'ensembles sur~$X$; soit $\Phi$ le champ en
cat\'egories discr\`etes associ\'e \`a $F$, \ie le champ dont
la fibre au-dessus de tout sch\'ema $X_1$ \'etale sur~$X$ est la
cat\'egorie discr\`ete ayant pour ensemble d'objets $F(X_1)$. On
dit que $(F,f)$ est cohomologiquement propre relativement \`a~$S$ en
dimension $\leqslant -1$ (\resp en dimension $\leqslant 0$) si $(\Phi,
f)$ est cohomologiquement propre relativement \`a~$S$ en dimension
$\leqslant 0$ (\resp en dimension $\leqslant 1$).

Le
\marginpar{348}
morphisme canonique
\begin{equation*}
\label{eq:XIII.1.2.1}
\tag{\thesubsection.1} { g^* f_* F \to f'_* h^* F }
\end{equation*}
donne par passage aux champs en cat\'egories discr\`etes
associ\'ees le morphisme canonique
$$
g^* f_* \Phi \to f'_* h^* \Phi \text{.}
$$
Par suite dire que $(F,f)$ est cohomologiquement propre relativement
\`a~$S$ en dimension $\leqslant -1$ (\resp en dimension $\leqslant
0$) \'equivaut \`a dire que, pour tout $S$-sch\'ema $S'$, le
morphisme \eqref{eq:XIII.1.2.1} est injectif (\resp bijectif).

\subsection{}
\label{XIII.1.3}
Soit $F$ un faisceau en groupes sur~$X$ et $\Phi$ le champ des
torseurs sur~$X$ de groupe $F$ \cite[II 2.3.2]{XIII.1}. On dit que $(F,f)$ est
cohomologiquement propre relativement \`a~$S$ en dimension
$\leqslant -1$ (\resp $\leqslant 0$, \resp $\leqslant 1$) si
$(\Phi,f)$ est cohomologiquement propre relativement \`a~$S$ en
dimension $\leqslant -1$ (\resp $\leqslant 0$, \resp $\leqslant
1$). La condition de propret\'e cohomologique peut s'expliciter
comme suit.

\begin{subproposition}
\label{XIII.1.3.1}
Les notations sont celles de \eqref{eq:XIII.1.0.1} et
\eqref{eq:XIII.1.0.2}. Soit $F$ un faisceau en groupes sur~$X$. On
d\'esigne par $F'$ (\resp $F_1$, \resp $F_1'$, etc.), l'image
inverse de~$F$ sur~$X'$ (\resp sur~$X_1$, \resp sur~$X_1'$,
etc.). Alors les conditions suivantes sont \'equivalentes:
\begin{enumerate}
\item[(i)] $(F,f)$ est cohomologiquement propre relativement \`a~$S$
en dimension $\leqslant -1$ (\resp $\leqslant 0$, \resp $\leqslant
1$).
\item[(ii)] Pour tout morphisme $S' \to S$, pour tout sch\'ema $Y_1$
\'etale au-dessus de~$Y$, et pour tout torseur $P$ sur~$X_1$ de
groupe $F_1$, si $\leftexp{P\mkern-4mu}{F}_1$ d\'esigne le groupe tordu de
$F_1$ par $P$ \textup{\cite[II 4.1.2.3]{XIII.1}}, le morphisme canonique
$$
a_0 \colon g_1^*(f_{1 *}(\leftexp{P\mkern-4mu}{F}_1)) \to f'_{1
*}(\leftexp{P'\mkern-6mu}{F}'_1)
$$
est injectif (\resp $a_0$ est bijectif et le morphisme canonique
$$
a_1 \colon g^*(\R^1 f_* F) \to \R^1 f'_* F'
$$
est
\marginpar{349}
injectif, \resp $a_0$ et $a_1$ sont bijectifs).
\item[(ii~bis)] Pour tout morphisme $S' \to S$, pour tout sch\'ema
$Y_1$ \'etale au-dessus de~$Y$, pour tout torseur $P$ sur~$X_1$ de
groupe $F_1$, et pour tout torseur $R$ sous $\leftexp{P\mkern-4mu}{F}_1$, le
morphisme canonique
$$
\alpha_0 \colon g_1^*(f_{1 *}R) \to f'_{1 *} R'
$$
est injectif (\resp $\alpha_0$ est bijectif, \resp les morphismes
$\alpha_0$ et
$$
\alpha_1 \colon g_1^*(\R^1 f_{1 *}(\leftexp{P\mkern-4mu}{F}_1)) \to \R^1
f'_{1 *}(\leftexp{P'\mkern-6mu}{F}'_1)
$$
sont bijectifs).
\end{enumerate}
\end{subproposition}

\subsubsection*{D\'emonstration}
(i)$\To$(ii~bis). D'apr\`es \cite[II 4.2.5]{XIII.1} tout torseur $R$
de groupe $\leftexp{P\mkern-4mu}{F}_1$ est de la forme $R = Q
\overset{F_1}{\wedge} P^\circ$, o\`u $Q$ est un torseur de groupe
$F_1$ et $P^\circ$ l'oppos\'e de~$P$. On a alors $R' \simeq Q'
\overset{F'_1}{\wedge} {P'}^\circ$. Soit $\Phi$ le champ des torseurs
sous $F$ et soient $x$, $y$ (\resp~$x'$,~$y'$) les objets de la
cat\'egorie fibre $(g^* f_* \Phi)_{Y_1'}$ (\resp $(f'_*
\Phi')_{Y_1'}$) associ\'es \`a $P$, $Q$ (\resp $P'$,~$Q'$). On a
la relation
$$
Q\overset{F_1}{\wedge} P^\circ \simeq \SheafHom_{F_1}(P,Q) \text{,}
$$
et il en r\'esulte que l'on a des isomorphismes canoniques
$$
\SheafHom_{Y'_1}(x,y) \simeq g_1^* f_{1 *}(Q
\overset{F_1}{\wedge} P^\circ) \text{,}\quad \SheafHom_{Y_1'}(x',y')
\simeq f'_{1 *}(Q' \overset{F'_1}{\wedge}{P'}^\circ)\text{.}
$$
Par suite le morphisme $\alpha_0$ s'identifie au morphisme
$$
\SheafHom_{Y'_1}(x,y) \to \SheafHom_{Y'_1}(x',y') \text{,}
$$
d'o\`u il r\'esulte que, si $(F,f)$ est cohomologiquement propre
relativement \`a~$S$ en dimension $\leqslant -1$, $\alpha_0$ est
injectif et que, si $(F,f)$ est cohomologiquement propre en dimension
$\leqslant 0$, $\alpha_0$ est bijectif.

Supposons maintenant que $(F,f)$ soit cohomologiquement propre
relativement \`a~$S$ en dimension $\leqslant 1$, \ie que le
morphisme canonique
\marginpar{350}
$$
\varphi \colon g^* f_* \Phi \to f'_* \Phi'
$$
soit une \'equivalence. Soit $G$ le faisceau des sous-gerbes
maximales du champ $f_* \Phi$ \hbox{\cite[III~2.1.8]{XIII.1}}; on a alors un
isomorphisme $G \simeq \R^1 f_* F$. Comme $g^* G$ est le
faisceau des sous-gerbes maximales de~$g^* f_* \Phi$ \cite[III
2.1.5.5]{XIII.2}, le morphisme $\alpha_1$ est obtenu \`a partir de
$\varphi|Y'_1$ en prenant les faisceaux des sous-gerbes maximales,
donc est un isomorphisme.

(ii~bis)$\To$(ii). Il suffit de montrer que, si les
morphismes $\alpha_0$ sont bijectifs, alors les morphismes $a_1$ sont
injectifs. Soient $Y'_1$ un sch\'ema \'etale au-dessus de~$Y'$,
$s$ et $t$ deux \'el\'ements de~$g^*(\R^1 f_* F)(Y'_1)$
ayant m\^eme image dans $\R^1 f'_* F'(Y'_1)$ et montrons que
l'on a $s=t$. L'assertion est locale pour la topologie \'etale de
$Y'_1$ et, compte tenu de la d\'efinition de l'image inverse
$g^*(\R^1 f_* F)$, on peut supposer que $Y'_1$ est l'image
inverse d'un sch\'ema $Y_1$ \'etale au-dessus de~$Y$ et que $s$ et
$t$ proviennent de torseurs $P$ et $Q$ sur~$X_1$. L'hypoth\`ese
faite sur~$s$ et $t$ signifie alors que les images inverses $P'$ et
$Q'$ de~$P$ et $Q$ sur~$X'_1$ sont isomorphes localement pour la
topologie \'etale de~$Y'_1$. Si l'on pose $R =
\SheafHom_{F_1}(P,Q)$, le fait que le morphisme
$$
g_1^* f_{1 *} R \to f'_{1 *} R'
$$
soit bijectif prouve que $P$ et $Q$ sont isomorphes localement pour la
topologie \'etale de~$Y_1$, donc que l'on a $s=t$.

(ii)$\To$(i). Pour prouver que $\varphi$ est fid\`ele
(\resp pleinement fid\`ele), il suffit de montrer que, si $Y_1$ est
un sch\'ema \'etale sur~$Y$, si $P$, $Q$ sont deux torseurs sur~$X_1$ de groupe~$F_1$, si $x$, $y$ (\resp $x'$, $y'$) sont les objets
de~$(g^* f_* \Phi)_{Y'_1}$ (\resp $(f'_* \Phi')_{Y'_1}$)
associ\'es \`a $P$, $Q$ (\resp $P'$, $Q'$), le morphisme
$$
a \colon \Hom(x,y) \to \Hom(x',y')
$$
est injectif (\resp bijectif). Or $a$ s'identifie au morphisme
canonique
$$
\H^0(Y'_1,g^*_1 f_{1 *} (Q \overset{F_1}{\wedge} P^\circ)) \to
\H^0(Y'_1, f'_{1 *}(Q' \overset{F'_1}{\wedge} {P'}^\circ))\text{.}
$$
Si
\marginpar{351}
l'on a $\Hom(x,y) \neq \emptyset$, alors $Q \overset{F_1}{\wedge}
P^\circ$ est un torseur sous $\leftexp{P\mkern-4mu}{F}_1$ localement trivial sur~$Y_1$; il en r\'esulte que $f_{1 *}(Q \overset{F_1}{\wedge}
P^\circ)$ est un torseur sous $f_{1*}(\leftexp{P\mkern-4mu}{F}_1)$, et que
$g^*_1 f_{1 *} (Q \overset{F_1}{\wedge} P^\circ)$ est un torseur
trivial. Le morphisme $a$ s'identifie alors au morphisme canonique
$$
\H^0(Y_1, g^*_1 f_{1 *}(\leftexp{P\mkern-4mu}{F}_1)) \to \H^0 (Y'_1,
f'_{1 *}(\leftexp{P'\mkern-6mu}{F}'_1)) \text{.}
$$
Il en est de m\^eme si l'on a $\Hom(x',y') \neq \emptyset$ et si
$a_1$ est injectif car alors $Q' \overset{F'_1}{\wedge} {P'}^\circ$
est trivial, et il r\'esulte de l'injectivit\'e de~$a_1$ que $P$
et $Q$ sont localement isomorphes sur~$Y_1$. On en conclut que, si
$a_0$ est injectif (\resp si $a_0$ est bijectif et $a_1$ injectif),
$\varphi$ est fid\`ele (\resp pleinement fid\`ele).

Il reste \`a montrer que, si $a_0$ et $a_1$ sont bijectifs, le
foncteur $\varphi$ est essentiellement surjectif. Soient $Y''$ un
sch\'ema \'etale au-dessus de~$Y'$, $X'' = X' \times_{Y'} Y''$ et
soit $P''$ un torseur sur~$X''$ de groupe $F'' = F'|X''$. On va
montrer qu'il existe un \'el\'ement~$x$ de~$(g^* f_*
\Phi)_{Y''}$ dont l'image dans $(f'_* \Phi')_{Y''}$ est isomorphe
\`a $P''$. Soit $p''$ la classe de~$P''$. Du fait que $a_1$ est
surjectif r\'esulte que l'on peut trouver un morphisme \'etale
surjectif $Y''_1 \to Y''$, un morphisme \'etale $Y_1 \to Y$ tel que
l'on ait un morphisme $Y''_1 \to Y'_1$ et un torseur $P_1$ sur~$X_1$
de groupe $F_1$ dont l'image inverse $P''_1$ sur~$X''_1$ soit
isomorphe \`a l'image inverse de~$P''$. Utilisant le fait que
$\varphi$ est pleinement fid\`ele, on voit que l'objet $x_1$ de
$(g^* f_* \Phi)_{Y''_1}$ qui correspond \`a $P''_1$ est muni
d'une donn\'ee de descente relativement \`a $Y''_1 \to Y''$, donc
provient d'un \'el\'ement $x$ de~$(g^* f_*
\Phi)_{Y''}$. Comme l'image de~$x$ dans $(f'_* \Phi')_{Y''}$ est
$P''$, ceci prouve que $\varphi$ est essentiellement surjectif et
ach\`eve la d\'emonstration.

\begin{exemple}
\label{XIII.1.4}
Soit $f \colon X \to Y$ un morphisme propre. Il r\'esulte de \cite[VII
2.2.2]{XIII.2} que, pour tout champ ind-fini $F$ sur~$X$, le couple $(F,f)$
est cohomologiquement propre (relativement \`a $Y$) en dimension
$\leqslant 1$. En particulier, pour tout faisceau d'ensembles
(\resp tout faisceau de groupes, \resp tout faisceau en groupes
ind-fini) $F$ sur~$X$, $(F,f)$ est cohomologiquement
\marginpar{352}
propre en dimension $\leqslant 0$ (\resp en dimension $\leqslant 0$,
\resp en dimension $\leqslant 1$).
\end{exemple}

\begin{remarques}
\label{XIII.1.5}
a) Soit $F$ un faisceau en groupes sur~$X$ tel que $(F,f)$ soit
cohomologiquement propre relativement \`a~$S$ en dimension
$\leqslant -1$ (\resp $\leqslant 0$). Si l'on consid\`ere $F$ comme
faisceau d'ensembles, $(F,f)$ est cohomologiquement propre
relativement \`a~$S$ en dimension $\leqslant -1$ (\resp $\leqslant
0$), mais la r\'eciproque est fausse.

Soit par exemple $Y$ le spectre d'un anneau de valuation discr\`ete
strictement local de point ferm\'e $t$, de point g\'en\'erique
$s$, $f \colon X \to Y$ un sch\'ema non vide sur~$Y$ dont la fibre
ferm\'ee est vide, $F$ un faisceau en groupes constant non trivial
sur~$X$ et $P$ un torseur sous $F$ tel que l'on ait $\H^0(X_s,
\leftexp{P\mkern-4mu}{F}|X_s)=1$. Alors $(\leftexp{P\mkern-4mu}{F},f)$ est
cohomologiquement propre relativement \`a $Y$ en dimension
$\leqslant -1$ lorsque l'on consid\`ere $\leftexp{P\mkern-4mu}{F}$ comme
faisceau d'ensembles. Si l'on consid\`ere $\leftexp{P\mkern-4mu}{F}$ comme faisceau
en groupes, on a un isomorphisme ${\vphantom{F}}^{P^\circ\mkern-7mu}(\leftexp{P\mkern-4mu}{F})
\simeq F$; comme le morphisme canonique
$$
\H^0(X,F)\to\H^0(X_t,F|X_t)=1
$$
n'est pas injectif, ceci prouve que $(\leftexp{P\mkern-4mu}{F},f)$ n'est pas
cohomologiquement propre relativement \`a $Y$ en dimension
$\leqslant -1$.

b) Supposons $f$ coh\'erent (\ie quasi-compact et
quasi-s\'epar\'e). Soit $F$ un champ sur~$X$. Pour tout point
g\'eom\'etrique $\overline{y}$ de~$Y'$, on note $\overline{Y}$
(\resp $\overline{Y'}$) le localis\'e strict de~$Y$ (\resp $Y'$) en
$\overline{y}$, et on pose $\overline{X} = X \times_Y \overline{Y}$,
$\overline{X'} = X' \times_{Y'} \overline{Y'}$, etc. Pour que $(F,f)$
soit cohomologiquement propre relativement \`a~$S$ en dimension
$\leqslant -1$ (\resp $\leqslant 0$, \resp $\leqslant 1$), il faut et
il suffit que, pour tout $S$-sch\'ema $S'$ et pour tout point
g\'eom\'etrique $\overline{y}$ de~$Y'$, le foncteur canonique
$$
\overline{F}(\overline{X})\to\overline{F'}(\overline{X'})
$$
soit fid\`ele (\resp pleinement fid\`ele, \resp une
\'equivalence).

En
\marginpar{353}
effet, si $S'$ est un $S$-sch\'ema, pour que le foncteur
$$
g^* f_* F \to f'_* F'
$$
soit fid\`ele (\resp pleinement fid\`ele, \resp une
\'equivalence), il faut et il suffit qu'il en soit ainsi du foncteur
induit sur les fibres aux diff\'erents points g\'eom\'etriques
$\overline{y}'$ de~$\overline{Y'}$ \cite[III 2.1.5.9]{XIII.2}. L'assertion
r\'esulte donc du calcul des fibres g\'eom\'etriques de l'image
directe d'un champ par un morphisme coh\'erent \cite[VII 2.1.5]{XIII.2}.

c) Soit $F$ un champ sur~$X$. Le fait que $(F,f)$ soit
cohomologiquement propre relativement \`a~$S$ en dimension
$\leqslant -1$ (\resp $\leqslant 0$, \resp $\leqslant 1$) est local
sur~$Y$ pour la topologie \'etale.

Soit $S'$ un $S$-sch\'ema, $F'$ l'image inverse de~$F$ sur~$X'$
(\cf \eqref{eq:XIII.1.0.1}). Si $(F,f)$ est cohomologiquement propre
relativement \`a~$S$ en dimension $\leqslant 1$, il en est de
m\^eme de~$(F', f')$. Mais, si $(F,f)$ est cohomologiquement propre
relativement \`a~$S$ en dimension $\leqslant -1$ (\resp $\leqslant
0$), il n'en est pas n\'ecessairement de m\^eme de~$(F', f')$.

Soit par exemple $S'$ un anneau de valuation discr\`ete, $f' \colon
E_{S'} \to S'$ l'espace affine au-dessus de~$S'$, $x$ un point
ferm\'e de~$E_{S'}$ au-dessus du point g\'en\'erique de~$S'$ et
$F'$ le faisceau d'ensembles sur~$E_{S'}$ dont la restriction \`a
$E_{S'}-\{x\}$ est le faisceau constant \`a un \'el\'ement et
dont la fibre en un point g\'eom\'etrique au-dessus de~$x$ a deux
\'el\'ements. Alors $(F',f')$ n'est pas cohomologiquement propre
relativement \`a~$S'$ en dimension $\leqslant -1$. Soient $S =
S'[Z]$, $f \colon E_S \to S$ l'espace affine sur~$S$ et $T$ une partie
ferm\'ee de~$X = E_S$ qui ne rencontre par le ferm\'e $Z=0$ et
telle que $f(T)$ contienne le point g\'en\'erique de~$S$. Soient
$G$ l'image inverse de~$F'$ sur~$X$ et $F$ le faisceau sur~$X$ obtenu
en prolongeant $G|X-T$ par le vide. Alors $(F,f)$ est
cohomologiquement propre relativement \`a~$S$ en dimension
$\leqslant -1$, mais il n'en est plus de m\^eme apr\`es le
changement de base $S' \to S$ d\'efini par $Z=0$.

d)
\marginpar{354}
Soit $F$ un champ sur~$X$ tel que $(F,f)$ soit cohomologiquement
propre relativement \`a $Y$ en dimension $\leqslant -1$ (\resp $\leqslant 0$, \resp $\leqslant 1$). Alors, pour tout point
g\'eom\'etrique $\overline{y}$ de~$Y$, le foncteur canonique
$$
(f_* F)_{\overline{y}} \to F(X_{\overline{y}})
$$
est fid\`ele (\resp pleinement fid\`ele, \resp une \'equivalence
de cat\'egories).
\end{remarques}

\begin{proposition}
\label{XIII.1.6}
Soient $f \colon X \to Y$ et $g \colon Y \to Z$ deux $S$-morphismes,
$\Phi$ un champ sur~$X$.
\begin{enumerate}
\item[1)] Supposons que $(\Phi,f)$ et $(f_* \Phi, g)$ soient
cohomologiquement propres relativement \`a~$S$ en dimension
$\leqslant -1$ (\resp $\leqslant 0$, \resp $\leqslant 1$). Alors il en
est de m\^eme de~$(\Phi, gf)$.
\item[2)] Supposons que $(\Phi, gf)$ soit cohomologiquement propre
relativement \`a~$S$ en dimension $\leqslant -1$ (\resp que $(\Phi,
gf)$ soit cohomologiquement propre relativement \`a~$S$ en dimension
$\leqslant 0$ et $(\Phi,f)$ cohomologiquement propre relativement
\`a~$S$ en dimension $\leqslant -1$, \resp que $(\Phi, gf)$ soit
cohomologiquement propre relativement \`a~$S$ en dimension
$\leqslant 1$ et $(\Phi, f)$ cohomologiquement propre relativement
\`a~$S$ en dimension $\leqslant 0$). Alors $(f_* \Phi, g)$ est
cohomologiquement propre relativement \`a~$S$ en dimension
$\leqslant -1$ (\resp en dimension $\leqslant 0$, \resp en dimension
$\leqslant 1$).
\end{enumerate}
\end{proposition}

Pour tout $S$-sch\'ema $S'$, on consid\`ere le diagramme suivant,
dont tous les carr\'es sont cart\'esiens:
\begin{equation*}
\label{eq:XIII.1.6.1}
\tag{\thesubsection.1}
\begin{array}{c}
\xymatrix{ X' \ar[r]^{f'} \ar[d]_h & Y'
\ar[r]^{g'} \ar[d]_k & Z' \ar[r] \ar[d]_{m} & S' \ar[d] \\ X \ar[r]^f
& Y \ar[r]^g & Z \ar[r] & \,S. }
\end{array}
\end{equation*}
\ifthenelse{\boolean{orig}}
{le}
{Le}
morphisme canonique
$$
m^* (g_* f_* \Phi) \to g'_* f'_* (h^* \Phi)
$$
s'identifie
\marginpar{355}
au compos\'e des morphismes canoniques
$$
m^*(g_* f_* \Phi) \lto{i} g'_*(k^*f_* \Phi) \lto{j} g'_* f'_*(h^* \Phi).
$$
\begin{enumerate}
\item[1)] L'hypoth\`ese entra\^ine que $i$ et $j$ sont
fid\`eles (\resp pleinement fid\`eles, \resp des
\'equi\-va\-lences); il en est donc de m\^eme de~$ji$.
\item[2)] L'hypoth\`ese entra\^ine que $ji$ est fid\`ele
(\resp que $ji$ est
\ifthenelse{\boolean{orig}}
{pelinement}
{pleinement}
fid\`ele et $j$ fid\`ele, \resp que $ji$ est une \'equivalence
et $j$ pleinement fid\`ele); il en r\'esulte que $i$ est
fid\`ele (\resp pleinement fid\`ele, \resp une \'equivalence).
\end{enumerate}

\begin{corollaire}
\label{XIII.1.7}
Soient $f \colon X \to Y$ et $g \colon Y \to Z$ deux $S$-morphismes,
et soit $F$ un faisceau en groupes sur~$X$. Supposons que $(F,gf)$
soit cohomologiquement propre relativement \`a~$S$ en dimension
$\leqslant -1$ (\resp que $(F,gf)$ soit cohomologiquement propre
relativement \`a~$S$ en dimension $\leqslant 0$ et que $(F,f)$ soit
cohomologiquement propre relativement \`a~$S$ en dimension
$\leqslant -1$). Alors $(f_* F, g)$ est cohomologiquement propre
relativement \`a~$S$ en dimension $\leqslant -1$ (\resp en dimension
$\leqslant 0$).
\end{corollaire}

Reprenons les notations de \eqref{eq:XIII.1.6.1} et, pour tout
sch\'ema $Y_1$ \'etale au-dessus de~$Y$, notons $f_1$, $F_1$ les
images inverses respectives de~$f$, $F$ par le morphisme $Y_1 \to
Y$. Soient $\Phi$ le champ des torseurs sous $F$ et $\Psi$ le champ
des torseurs sous $f_* F$. On a un foncteur canonique
$$
\varphi \colon \Psi \to f_* \Phi \text{,}
$$
obtenu en associant \`a tout sch\'ema $Y_1$ \'etale sur~$Y$ et
\`a tout torseur $P$ sur~$Y_1$ de groupe $f_{1 *}F_1$ le torseur
$\tilde{P}$ sur~$X_1$ d\'eduit de~$f_1^* P$ par l'extension du
groupe structural $f_1^* f_{1 *}F_1 \to F_1$. Le foncteur
$\varphi$ est pleinement fid\`ele. En effet, si $P$ et $Q$ sont deux
torseurs sur~$Y_1$ de groupe $f_{1 *} F_1$, on a un morphisme
canonique
$$
\SheafIsom_{f_{1 *}F_1}(P,Q) \to
f_{1*}(\SheafIsom_{F_1}(\tilde{P},\tilde{Q}))
$$
qui
\marginpar{356}
est un isomorphisme car il en est ainsi localement. On en d\'eduit
que le morphisme canonique
$$
\Isom_{f_{1 *}F_1}(P,Q) \to \Isom_{F_1}(\tilde{P},\tilde{Q})
$$
est un isomorphisme, donc que $\varphi$ est pleinement fid\`ele.

On a un diagramme commutatif
$$
\xymatrix{ g'_* k^* \Psi \ar[r]^i \ar[d]_{g'_*
k^*(\varphi)} & m^*(g_* \Psi) \ar[d]^{m^*
g_*(\varphi)} \\ g'_* k^*(f_* \Phi) \ar[r]^j &
m^*(g_* f_* \Phi)\text{,} }
$$
o\`u $i$ et $j$ sont les morphismes de changement de base. Il
r\'esulte de~\Ref{XIII.1.6} 2) que $j$ est fid\`ele
(\resp pleinement fid\`ele). Comme $g'_* k^* (\varphi)$ et
$m^* g_* (\varphi)$ sont pleinement fid\`eles, on d\'eduit
du diagramme ci-dessus que $i$ est fid\`ele (\resp pleinement
fid\`ele).

\begin{corollaire}
\label{XIII.1.8}
Soient $f \colon X \to Y$ un $S$-morphisme coh\'erent, $g \colon Y
\to Z$ un $S$-morphisme propre, $\Phi$ un champ ind-fini sur~$X$ \textup{\cite[VII 2.2.1]{XIII.2}}. Supposons que $(\Phi,f)$ soit cohomologiquement propre
relativement \`a~$S$ en dimension $\leqslant -1$ (\resp $\leqslant
0$, \resp $\leqslant 1$), alors il en est de m\^eme de~$(\Phi, gf)$.
\end{corollaire}

Comme $f$ est coh\'erent, $f_* \Phi$ est un champ ind-fini (SGA~4 IX
1.6~(ii)). Le corollaire r\'esulte donc de \Ref{XIII.1.6}~1)
et~\Ref{XIII.1.4}.

\begin{corollaire}
\label{XIII.1.9}
Soient $f \colon X \to Y$ un $S$-morphisme entier, $g \colon Y \to Z$
un $S$-morphisme. Si $F$ est un faisceau d'ensembles sur~$X$, pour que
$(f_* F, g)$ soit cohomologiquement propre relativement \`a~$S$
en dimension $\leqslant -1$ (\resp $\leqslant 0$), il faut et il
suffit qu'il en soit ainsi de~$(F, gf)$. Si $F$ est un faisceau en
groupes sur~$X$, pour que $(f_* F, g)$ soit cohomologiquement
propre relativement \`a~$S$ en dimension $\leqslant -1$
(\resp $\leqslant 0$, \resp $\leqslant 1$), il faut et il suffit qu'il
en soit ainsi de~$(F, gf)$.
\end{corollaire}

L'assertion
\marginpar{357}
relative au cas d'un faisceau d'ensembles r\'esulte
de~\Ref{XIII.1.6} et du fait que $(F,f)$ est cohomologiquement propre
relativement \`a~$S$ en dimension $\leqslant 0$. Soient $F$ un
faisceau en groupes sur~$X$ et $\Phi$ le champ des torseurs sous $F$.
D'apr\`es SGA~4 VIII~5.8, tout torseur sous $F$ est localement
trivial sur~$Y$. Il en r\'esulte que le champ $f_* \Phi$ est
\'equivalent au champ des torseurs sous $f_* F$, l'\'equivalence
\'etant obtenue en associant \`a tout sch\'ema $Y_1$ \'etale
sur~$Y$ et \`a tout torseur $P$ sur~$X_1 = X \times_Y Y_1$ de groupe
$F|X_1$ le torseur $f_* P$ de groupe $f_* F | Y_1$. Comme $(F,f)$ est
cohomologiquement propre relativement \`a~$S$ en dimension
$\leqslant 1$, le corollaire r\'esulte donc de~\Ref{XIII.1.6}.

\enlargethispage{.5cm}
\skpt
\begin{definitions}
\label{XIII.1.10}

\subsubsection{}
\label{XIII.1.10.1}
Soit $E$ une cat\'egorie et consid\'erons un diagramme
$$
\xymatrix@C=.6cm{ \Phi \ar[r]^p & \Phi_1 \ar@<2pt>[r]^{p_1}
\ar@<-2pt>[r]_{p_2} & \Phi_2 } \text{,}
$$
o\`u $\Phi$, $\Phi_1$, $\Phi_2$ sont des cat\'egories fibr\'ees
au-dessus de~$E$ et les fl\`eches des morphismes de cat\'egories
fibr\'ees, et soit
$$
a \colon p_1 p \isomto p_2 p
$$
un isomorphisme de foncteurs.

On dit que le diagramme ci-dessus est exact si la condition suivante
est satisfaite

a) Pour tout couple d'objets $x$, $y$ de~$\Phi$ et tout morphisme $f
\colon p(x) \to p(y)$ tel que l'on ait $p_1(f) = p_2(f)$ ($p_1p$ et
$p_2p$ \'etant identifi\'es gr\^ace \`a $a$), il existe un
unique morphisme $g \colon x \to y$ tel que l'on ait $p(g)=f$.

\subsubsection{}
\label{XIII.1.10.2}
Consid\'erons le diagramme
\ifthenelse{\boolean{orig}}
{$$
\xymatrix{ \Phi \ar[r]^p & \Phi_1 \ar@<0.6ex>[r]^{p_1}
\ar@<-0.6ex>[r]_{p_2} & \Phi_2 \ar@<2ex>[r]^{p_{23}} \ar[r]^{p_{31}}
\ar@<-0.6ex>[r]_{p_{12}} & \Phi_3 } \text{,}
$$}
{$$
\xymatrix@C=.5cm{ \Phi \ar[r]^p & \Phi_1 \ar@<2pt>[rr]^{p_1,p_2}
\ar@<-2pt>[rr] && \Phi_2 \ar@<4pt>[rrr]^{p_{12},p_{23},p_{31}} \ar[rrr]
\ar@<-4pt>[rrr] &&& \Phi_3 } \text{,}
$$}
o\`u
\marginpar{358}
$\Phi$, $\Phi_i$, $1 \leqslant i \leqslant 3$, sont des cat\'egories
fibr\'ees sur~$E$ et les fl\`eches des morphismes de
cat\'egories fibr\'ees. Supposons donn\'es des isomorphismes de
foncteurs
$$
a \colon p_1 p \isomto p_2 p
$$
$$
a_1 \colon p_{31} p_2 \isomto p_{12} p_1 \text{,}\quad a_2 \colon p_{12}
p_2 \isomto p_{23} p_1 \text{,}\quad a_3 \colon p_{23} p_2 \isomto
p_{31}p_1
$$
tels que le diagramme suivant soit commutatif:
$$
\xymatrix@C=1.2cm{ p_{23} p_1 p \ar[r]^{\id.a} & p_{23} p_2 p \ar[r]^{a_3.\id}
& p_{31} p_1 p \ar[d]^{\id.a} \\ p_{12} p_2 p \ar[u]^{a_2.\id} &
p_{12} p_1 p \ar[l]^{\id.a} & p_{31} p_2 p \ar[l]^{a_1.\id} } \text{.}
$$
On identifie $p_1 p$ et $p_2 p$, $p_{31} p_2$ et $p_{12} p_1$, etc.

On dit que le diagramme ci-dessus est exact si les conditions
suivantes sont satisfaites:
\begin{enumerate}
\item[a)] Analogue \`a la condition a) de~\Ref{XIII.1.10.1}.
\item[b)] Pour tout objet $x_1$ de~$\Phi_1$ et pour tout isomorphisme
$u \colon p_1(x_1) \isomto p_2(x_1)$ tel que l'on ait
\begin{equation*}
\label{eq:XIII.1.10.2.1}
\tag{\thesubsubsection.1} { p_{23}(u) p_{31}(u) = p_{12}(u)^{-1}
\text{,} }
\end{equation*}
il existe un objet $x$ de~$\Phi$ tel que l'on ait un isomorphisme $i
\colon p(x) \isomto x_1$ rendant commutatif le diagramme
\begin{equation*}
\label{eq:XIII.1.10.2.2}
\tag{\thesubsubsection.2}
\begin{array}{c}
{ \xymatrix{ p_1p(x) \ar@{=}[r]
\ar[d]_{p_1(i)} & p_2p(x) \ar[d]^{p_2(i)} \\ p_1(x_1) \ar[r]^{\sim}_u
& p_2(x_1)
} }
\end{array}
\end{equation*}
\end{enumerate}

\subsubsection{}
\label{XIII.1.10.3}
On d\'efinit de fa\c con \'evidente la notion de morphisme de
diagrammes exacts de cat\'egories fibr\'ees au-dessus d'une
cat\'egorie $E$.
\end{definitions}

\subsubsection{}
\label{XIII.1.10.4}
Nous utiliserons
\marginpar{359}%
en particulier la notion de diagramme exact dans le
cas o\`u $E$ est un site et $\Phi$, $\Phi_i$, $1 \leqslant i
\leqslant 3$, des champs sur~$E$.

Soient $f \colon E \to E'$ un morphisme de sites et
\ifthenelse{\boolean{orig}}
{\begin{equation*}
\label{eq:XIII.1.10.4.1}
\tag{\thesubsubsection.1} { \xymatrix{ \Phi \ar[r]^p & \Phi_1
\ar@<0.6ex>[r]^{p_1} \ar@<-0.6ex>[r]_{p_2} & \Phi_2
\ar@<2ex>[r]^{p_{23}} \ar[r]^{p_{31}} \ar@<-0.6ex>[r]_{p_{12}} &
\Phi_3 } }
\end{equation*}}
{\begin{equation*}
\label{eq:XIII.1.10.4.1}
\tag{\thesubsubsection.1} { \xymatrix@C=.5cm{ \Phi \ar[r]^p & \Phi_1
\ar@<2pt>[rr]^{p_1,p_2} \ar@<-2pt>[rr] && \Phi_2
\ar@<4pt>[rrr]^{p_{12},p_{23},p_{31}} \ar[rrr] \ar@<-4pt>[rrr] &&&
\Phi_3 } }
\end{equation*}}
un diagramme exact de champs sur~$E$. On obtient par image directe un
diagramme
$$
\xymatrix@C=.5cm{ f_* \Phi \ar[r] & f_* \Phi_1 \ar@<2pt>[r]
\ar@<-2pt>[r] & f_* \Phi_2 \ar@<4pt>[r]\ar[r] \ar@<-4pt>[r] &
f_* \Phi_3 }
$$
qui est \'evidemment exact.

Si $h \colon E'' \to E$ est un morphisme de sites, on a de m\^eme un
diagramme exact
\ifthenelse{\boolean{orig}}
{$$
\xymatrix{ h^* \Phi \ar[r]^{p''} & h^* \Phi_1
\ar@<0.6ex>[r]^{p''_1} \ar@<-0.6ex>[r]_{p''_2} & h^* \Phi_2
\ar@<2ex>[r]^{p''_{23}} \ar[r]^{p''_{31}} \ar@<-0.6ex>[r]_{p''_{12}} &
h^* \Phi_3 }
$$}
{$$
\xymatrix@C=.5cm{ h^* \Phi \ar[r]^{p''} & h^* \Phi_1
\ar@<2pt>[rr]^{p''_1,p''_2} \ar@<-2pt>[rr] && h^* \Phi_2
\ar@<4pt>[rrr]^{p''_{12},p''_{23},p''_{31}} \ar[rrr] \ar@<-4pt>[rrr] &&&
h^* \Phi_3 }
$$}

V\'erifions d'abord la condition~a) de~\Ref{XIII.1.10.2}. Soient
$F''$ un objet de~$E''$, $x''$ et $y''$ deux objets de~$(h^*
\Phi)_{F''}$, $x''_1$ et $y''_1$ leurs images respectives dans
$h^* \Phi_1$ et $x''_2$, $y''_2$ leurs images dans $h^*
\Phi_2$. Soit $u''_1 \colon x''_1 \to y''_1$ un morphisme tel que l'on
ait $p''_1(u''_1) = p''_2(u''_1)$ et prouvons que $u''_1$ provient
d'un unique morphisme $u'' \colon x'' \to y''$. La question \'etant
locale sur~$F''$, on peut supposer qu'on a un objet $F_1$ de~$E$, un
morphisme de~$F''$ dans l'image inverse $F''_1$ de~$F_1$ par $h$, des
objets $x$, $y$ de~$\Phi_{F_1}$ dont les images inverses sur~$F''$
sont $x''$ et $y''$. Soient $x_1$, $y_1$ (\resp $x_2$, $y_2$) les
images de~$x$, $y$ dans $\Phi_1$ (\resp~$\Phi_2$). On peut supposer
que $u''_1$ provient d'un morphisme $u_1 \colon x_1 \to y_1$, tel que
l'on ait $p_1(u_1) = p_2(u_1)$. Vu l'exactitude de
\eqref{eq:XIII.1.10.4.1}, on obtient un unique morphisme $u \colon x
\to y$ dont l'image inverse par $h$ est le morphisme $u''$
cherch\'e.

La
\marginpar{360}
condition~b) de~\Ref{XIII.1.10.2} se v\'erifie de fa\c con
analogue. Soient $x''_1$ un objet de~$(h^* \Phi_1)_{F''}$, $u''
\colon p''_1(x''_1) \to p''_2(x''_1)$ un morphisme satisfaisant \`a
la relation
$$
p''_{23}(u'')p''_{31}(u'')=p''_{12}(u'')^{-1} \text{,}
$$
et prouvons qu'il existe un objet $x''$ de~$(h^* \Phi)_{F''}$ et
un isomorphisme $i'' \colon p''(x'') \simeq x''_1$ rendant commutatif
un diagramme analogue \`a \eqref{eq:XIII.1.10.2.2}. Comme la
question est locale sur~$F''$, on peut supposer qu'on a un objet
$F_1$, un morphisme $F'' \to F''_1$ comme ci-dessus, et un objet $x_1$
de~$(\Phi_1)_{F_1}$ dont l'image inverse dans $(h^* \Phi_1)_{F''}$
est $x''_1$. De m\^eme on peut supposer que $u''$ provient d'un
morphisme $u \colon p_1(x_1) \to p_2(x_1)$ satisfaisant \`a
\eqref{eq:XIII.1.10.2.2}. L'existence d'un objet $x$ de~$\Phi_{F_1}$,
dont l'image inverse par $h$ soit un \'el\'ement $x''$
r\'epondant \`a la question, r\'esulte alors de l'exactitude de
\eqref{eq:XIII.1.10.4.1}.

\begin{exemples}
\label{XIII.1.11}
\begin{enumerate}
\item[1)] Soit $f \colon X_1 \to X$ un \emph{morphisme de descente}
pour la cat\'egorie des faisceaux \'etales sur des sch\'emas
variables (par exemple un morphisme universellement submersif (SGA~4
VIII~9.3)). Soient $X_2 = X_1 \times_X X_1$, $g \colon X_2 \to X$ la
projection canonique et $F$ un faisceau d'ensembles sur~$X$. Il
r\'esulte alors de \loccit que l'on a une suite exacte de
faisceaux d'ensembles
\begin{equation*}
\label{eq:XIII.1.11.1}
\tag{\thesubsection.1} { \xymatrix@C=.5cm{ F \ar[r] & f_* f^* F
\ar@<2pt>[r] \ar@<-2pt>[r] & g_* g^* F \text{.}}
}
\end{equation*}
Si $\Phi$ est le champ en cat\'egories discr\`etes associ\'e
\`a $F$ et $\Phi_3$ le champ final sur~$X$, \ie le champ dont
toutes les fibres sont r\'eduites \`a un seul \'el\'ement
ayant pour seul morphisme le morphisme identique, dire que la suite
\eqref{eq:XIII.1.11.1} est exacte revient \`a dire qu'il en est
ainsi du diagramme de champs
$$
\xymatrix@C=.5cm{ \Phi \ar[r] & f_* f^* \Phi \ar@<2pt>[r]
\ar@<-2pt>[r] & g_* g^* \Phi \ar[r] \ar@<4pt>[r]
\ar@<-4pt>[r] & \Phi_3 } \text{.}
$$

\item[2)]
Soit
\marginpar{361}
$f \colon X_1 \to X$ un \emph{morphisme de descente effective}
pour la cat\'egorie des faisceaux \'etales sur des sch\'emas
variables (par exemple un morphisme propre surjectif, ou un morphisme
entier surjectif, ou un morphisme fid\`element plat localement de
pr\'esentation finie (SGA~4 VIII~9.4)). Soient $X_2 = X_1 \times_X
X_1$, $g \colon X_2 \to X$ la projection canonique, $X_3 = X_1
\times_X X_1 \times_X X_1$, $h \colon X_3 \to X$ le morphisme
canonique. Soient $\Phi$ un champ sur~$X$, $\Phi_1 = f_* f^* \Phi$,
$\Phi_2 = g_* g^* \Phi$, $\Phi_3 = h_* h^* \Phi$. On a alors un
diagramme exact
$$
\xymatrix@C=.5cm{ \Phi \ar[r] & \Phi_1 \ar@<2pt>[r] \ar@<-2pt>[r] &
\Phi_2 \ar[r] \ar@<4pt>[r] \ar@<-4pt>[r] & \Phi_3 \text{,} }
$$
o\`u les fl\`eches sont les morphismes canoniques associ\'es aux
projections.

Consid\'erons en effet $\Phi$ comme un champ sur la cat\'egorie
$\Sch_X$ des sch\'emas au-dessus de~$X$, munie de la topologie
\'etale. Alors, d'apr\`es \cite[VII 2.2.8]{XIII.2}, $\Phi$ est aussi un
champ pour la topologie la plus fine sur~$\Sch_X$ pour laquelle les
morphismes couvrants sont les morphismes de descente effective pour la
cat\'egorie des faisceaux \'etales. L'exactitude du diagramme
ci-dessus en r\'esulte aussit\^ot.
\end{enumerate}
\end{exemples}

\begin{proposition}
\label{XIII.1.12}
Soient $S$ un sch\'ema, $f\colon X \to Y$ un $S$-morphisme.
\begin{enumerate}
\item[1)] Soit
$$
\xymatrix@C=.5cm{ \Phi \ar[r] & \Phi_1 \ar@<2pt>[r] \ar@<-2pt>[r] &
\Phi_2 }
$$
un diagramme exact de champs sur~$X$. Supposons que $(\Phi_1, f)$ soit
cohomologique\-ment~propre relativement \`a~$S$ en dimension
$\leqslant 0$ et que $(\Phi_2, f)$ soit cohomologi\-quement~propre~relativement \`a~$S$ en dimension $\leqslant -1$. Alors $(\Phi,f)$
est cohomologiquement propre~relativement \`a~$S$ en dimension
$\leqslant 0$.

\item[2)] Soit
$$
\xymatrix@C=.5cm{ \Phi \ar[r] & \Phi_1 \ar@<2pt>[r] \ar@<-2pt>[r] &
\Phi_2 \ar@<4pt>[r] \ar[r] \ar@<-4pt>[r] & \Phi_3 }
$$
un diagramme exact de champs sur~$X$. Supposons que $(\Phi,f)$ soit
cohomologi\-quement~propre relativement \`a~$S$ en dimension
$\leqslant 1$, que $(\Phi_2, f)$ soit cohomologiquement~propre
\marginpar{362}
relativement \`a~$S$ en dimension $\leqslant 0$ et que
$(\Phi_3, f)$ soit cohomo\-logiquement~propre
\ifthenelse{\boolean{orig}}
{}
{relativement}
\`a~$S$ en dimension $\leqslant -1$. Alors $(\Phi,f)$ est
cohomologique\-ment~propre~relativement \`a~$S$ en dimension
$\leqslant 1$.
\end{enumerate}
\end{proposition}

Pour tout $S$-sch\'ema $S'$, on consid\`ere le diagramme
commutatif suivant dont tous les carr\'es sont cart\'esiens:
$$
\xymatrix{ X' \ar[r]^{f'} \ar[d]_h & Y' \ar[r] \ar[d]_g & S' \ar[d] \\
X \ar[r]^f & Y \ar[r] & S }
$$

D\'emontrons 2), la d\'emonstration de 1) \'etant
analogue. Comme les foncteurs image directe et image inverse
transforment diagramme exact de champs en diagramme exact
\eqref{XIII.1.10.4}, on a le morphisme de diagrammes exacts de champs
suivant:
$$
\xymatrix{ g^* f_* \Phi \ar[r]^\pi \ar[d]_{\varphi} & g^*
f_* \Phi_1 \ar@<2pt>[r]^{\pi_1} \ar@<-2pt>[r]_{\pi_2}
\ar[d]_{\varphi_1} & g^* f_* \Phi_2 \ar@<4pt>[r] \ar[r]
\ar@<-4pt>[r] \ar[d]_{\varphi_2} & g^* f_* \Phi_3
\ar[d]_{\varphi_3} \\ f'_* h^* \Phi \ar[r] & f'_* h^*
\Phi_1 \ar@<2pt>[r] \ar@<-2pt>[r] & f'_* h^* \Phi_2
\ar@<4pt>[r] \ar[r] \ar@<-4pt>[r] & f'_* h^* \Phi_3 }
\text{.}
$$
Par hypoth\`ese $\varphi_1$ est une \'equivalence de
cat\'egories, $\varphi_2$ est pleinement fid\`ele et $\varphi_3$
fid\`ele. Il r\'esulte donc du diagramme pr\'ec\'edent que
$\varphi$ est une \'equivalence.

\begin{proposition}
\label{XIII.1.13}
Soit $f \colon X \to Y$ un $S$-morphisme.
\begin{enumerate}
\item[1)] Soit
$$
\xymatrix@C=.5cm{ F \ar[r] & G \ar@<2pt>[r] \ar@<-2pt>[r] & H }
$$
un diagramme exact de faisceaux d'ensembles sur~$X$. Supposons que
$(G,f)$ soit cohomologiquement propre relativement \`a~$S$ en
dimension $\leqslant 0$ et que $(H,f)$ soit cohomologiquement propre
relativement \`a~$S$ en dimension $\leqslant -1$. Alors $(F,f)$ est
cohomologiquement propre relativement \`a~$S$ en dimension
$\leqslant 0$.
\item[2)]
Soit
\marginpar{363}
$F \to G$ un monomorphisme de faisceaux en groupes sur~$X$. Si
$Y_1$ est un sch\'ema \'etale sur~$Y$, on pose $X_1 = Y_1 \times_Y
X$, et on note $f_1$ (\resp $F_1$, \resp $G_1$) l'image inverse de~$f$
(\resp $F$, \resp $G$) sur~$Y_1$ (\cf \Ref{eq:XIII.1.0.2}). Supposons
que $(G,f)$ soit cohomologiquement propre relativement \`a~$S$ en
dimension $\leqslant 0$ (\resp en dimension $\leqslant 1$) et que,
pour tout sch\'ema $Y_1$ \'etale sur~$Y$ et pour tout torseur $Q$
sous $G_1$, $(Q/F_1, f_1)$ soit cohomologiquement propre relativement
\`a~$S$ en dimension $\leqslant -1$ (\resp en dimension $\leqslant
0$). Alors $(F,f)$ est cohomologiquement propre relativement \`a~$S$
en dimension $\leqslant 0$ (\resp en dimension $\leqslant 1$).
\item[3)] Soit $F \to G$ un monomorphisme de faisceaux en groupes sur~$X$. Supposons que $(F,f)$ soit cohomologiquement propre relativement
\`a~$S$ en dimension $\leqslant 1$ et que $(G,f)$ soit
cohomologiquement propre relativement \`a~$S$ en dimension
$\leqslant 0$. Alors, pour tout torseur $Q$ sous $G$, $(Q/F, f)$ est
cohomologiquement propre relativement \`a~$S$ en dimension
$\leqslant 0$.
\end{enumerate}
\end{proposition}

\subsubsection*{D\'emonstration}
1) Soit $\Phi$ (\resp $\Phi_1$, \resp $\Phi_2$) le champ en
cat\'egories discr\`etes associ\'e \`a $F$ (\resp $G$,
\resp $H$) et soit $\Phi_3$ le champ final sur~$X$. On a alors un
diagramme exact
$$
\xymatrix@C=.5cm{ \Phi \ar[r] & \Phi_1 \ar@<2pt>[r] \ar@<-2pt>[r] &
\Phi_2 \ar@<4pt>[r] \ar[r] \ar@<-4pt>[r] & \Phi_3 }
\text{.}
$$
Par hypoth\`ese $(\Phi_1, f)$ est cohomologiquement propre
relativement \`a~$S$ en dimension $\leqslant 1$ et $(\Phi_2, f)$ est
cohomologiquement propre relativement \`a~$S$ en dimension
$\leqslant 0$ \eqref{XIII.1.2}. Comme $(\Phi_3, f)$ est \'evidemment
cohomologiquement propre relativement \`a~$S$ en dimension
$\leqslant -1$, il r\'esulte de~\Ref{XIII.1.12} que $(\Phi, f)$ est
cohomologiquement propre relativement \`a~$S$ en dimension
$\leqslant 1$, \ie que $(F,f)$ est cohomologiquement propre
relativement \`a~$S$ en dimension $\leqslant 0$.

\smallskip
2) Montrons d'abord que, si $(G,f)$ est cohomologiquement
propre relativement \`a~$S$ en dimension $\leqslant 0$ et si les
$(Q/F_1, f_1)$ sont cohomologiquement
\marginpar{364}
propres relativement \`a~$S$ en dimension $\leqslant -1$, alors
$(F,f)$ est cohomologiquement propre relativement \`a~$S$ en
dimension $\leqslant 0$. D'apr\`es~\Ref{XIII.1.3.1} il suffit de
prouver que, pour tout sch\'ema $Y_1$ \'etale au-dessus de~$Y$ et
pour tout torseur $P$ sur~$X_1$ de groupe $F_1$, $(\leftexp{P\mkern-4mu}{F}_1,
f_1)$ est cohomologiquement propre relativement \`a~$S$ en dimension
$\leqslant 0$ quand on consid\`ere $\leftexp{P\mkern-4mu}{F}_1$ comme un
faisceau d'ensembles, et que le morphisme canonique
$$
d \colon g^*(\R^1 f_* F) \to \R^1 f'_* F'
$$
est injectif. La premi\`ere assertion r\'esulte aussit\^ot de 1)
car, si $Q$ d\'esigne le torseur d\'eduit de~$P$ par l'extension
$F_1 \to G_1$ du groupe structural, on a un isomorphisme
$\leftexp{Q}{G}_1 / \leftexp{P\mkern-4mu}{F}_1 \isomto Q/F_1$.

Montrons que $d$ est injectif. Il suffit de prouver que, si $Y_1$ est
un sch\'ema \'etale au-dessus de~$Y$, si $P$ et $\tilde{P}$ sont
deux torseurs sous $F_1$ dont les images inverses $P'$ et~$\tilde{P}'$
sur~$X'_1$ sont isomorphes, alors, quitte \`a faire une extension
\'etale surjective de~$Y_1$, $P$ et~$\tilde{P}$ deviennent
isomorphes. Choisissons un isomorphisme $p' \colon P' \isomto
\tilde{P}'$. Soient $Q$ (\resp~$\tilde{Q}$) le torseur d\'eduit de
$P$ (\resp $\tilde{P}$) par l'extension du groupe structural $F_1 \to
G_1$. Les images inverses $Q'$ (\resp $\tilde{Q}'$) de~$Q$
(\resp $\tilde{Q}$) sur~$X'_1$ se d\'eduisent de~$P'$
(\resp $\tilde{P}'$) par
\ifthenelse{\boolean{orig}}
{extension}
{l'extension}
du groupe structural $F'_1 \to G'_1$; soit $q' \colon Q' \isomto
\tilde{Q}'$ l'isomorphisme que l'on obtient de m\^eme \`a partir
de~$p'$. Comme $(G,f)$ est cohomologiquement propre relativement \`a~$S$ en dimension $\leqslant 0$, on peut supposer, quitte \`a faire
une extension \'etale surjective de~$Y_1$, que $q'$ est l'image d'un
isomorphisme $q \colon Q \isomto \tilde{Q}$. Au torseur~$P$
(\resp~$\tilde{P}$) est associ\'ee une section $x$ de~$Q/F_1$
(\resp une section $\tilde{x}$ de~$\tilde{Q}/F_1$), et, pour que $P$
et $\tilde{P}$ soient isomorphes, il faut et il suffit que l'on ait un
isomorphisme $Q \isomto \tilde{Q}$ tel que l'isomorphisme
$$
e \colon \H^0(X_1, Q/F_1) \to \H^0(X_1,\tilde{Q}/F_1)
$$
qu'on en d\'eduit, transforme $x$ en $\tilde{x}$. On prend
l'isomorphisme $q$. Les sections $e(x)$ et $\tilde{x}$ de
$\H^0(X_1,\tilde{Q}/F_1)$ ont m\^eme image dans
$\H^0(X'_1,\tilde{Q}'/F'_1)$. Comme $(\tilde{Q}/F_1,f_1)$ est
cohomologiquement propre relativement \`a~$S$ en dimension
$\leqslant -1$,
\marginpar{365}
quitte \`a faire une extension \'etale surjective de~$Y_1$, on a
bien $e(x) = \tilde{x}$, ce qui d\'emontre l'injectivit\'e de~$d$.

Pour achever la d\'emonstration, il reste \`a prouver que, si
$(G,f)$ est cohomologiquement propre relativement \`a~$S$ en
dimension $\leqslant 1$, et si, pour tout sch\'ema $Y_1$ \'etale
au-dessus de~$Y$ et tout torseur $Q$ sur~$X_1$ de groupe $F_1$,
$(Q/F_1, f_1)$ est cohomologiquement propre relativement \`a~$S$ en
dimension $\leqslant 0$, alors le morphisme $d$ est surjectif. Soient
$P'$ un torseur sur~$X'_1$ de groupe $F'_1$, $Q'$ le torseur sous
$G'_1$ obtenu \`a partir de~$P'$ par extension du groupe
structural. La donn\'ee de~$P'$ est \'equivalente \`a celle de
$Q'$ et d'une section $x'$ de~$\H^0(X'_1,Q'/F'_1)$. Il r\'esulte
alors de la surjectivit\'e du morphisme
$$
g^*(\R^1 f_* G) \to \R^1 f'_* G'
$$
que, quitte \`a faire une extension \'etale surjective de~$Y_1$,
il existe un torseur $Q$ sous~$G_1$, dont l'image inverse sur~$X'_1$
est isomorphe \`a $Q'$. Utilisant le fait que $(Q/F_1, f_1)$ est
cohomologiquement propre relativement \`a~$S$ en dimension
$\leqslant 0$, on peut de m\^eme supposer qu'il existe un
\'el\'ement $x$ de~$\H^0(X_1, Q/F_1)$ dont l'image dans
$\H^0(X'_1, Q'/F'_1)$ est $x'$. La donn\'ee de~$Q$ et de~$x$
d\'etermine un torseur $P$ sous $F_1$, dont l'image inverse sur~$X_1$ est isomorphe \`a $P'$, ce qui d\'emontre la
surjectivit\'e de~$d$.

\smallskip
3) Montrons que $(Q/F,f)$ est cohomologiquement propre
relativement \`a~$S$ en dimension $\leqslant -1$, \ie que, pour
tout $S$-sch\'ema $S'$, pour tout sch\'ema $Y_1$ \'etale
au-dessus de~$Y$, si $x$, $\tilde{x}$ sont deux \'el\'ements de
$\H^0(X_1,Q_1/F_1)$ dont les images $x'$, $\tilde{x}'$ dans
$\H^0(X'_1, Q'_1/F'_1)$ sont \'egales, alors, apr\`es extension
surjective de~$Y_1$, on a $x = \tilde{x}$. \`A~$x$ (\resp~$\tilde{x}$)
est associ\'e un torseur $P$ (\resp~$\tilde{P}$) sous $F_1$, tel que
$Q_1$ se d\'eduise de~$P$ (\resp~$\tilde{P}$) par l'extension $F_1
\to G_1$ du groupe structural. De la relation $x' = \tilde{x}'$
r\'esulte que l'on a un isomorphisme $u' \colon P' \isomto
\tilde{P}'$ tel que l'isomorphisme induit sur~$Q'_1$ par~$u'$ soit
l'identit\'e. Comme $(F,f)$ est cohomologiquement propre
relativement \`a~$S$ en dimension $\leqslant 0$, on en d\'eduit
que,
\marginpar{366}
apr\`es extension \'etale surjective de~$Y_1$, on a un
isomorphisme $u \colon P \to \tilde{P}$ relevant $u'$; le fait que
$(G,f)$ soit cohomologiquement propre relativement \`a~$S$ en
dimension $\leqslant -1$ entra\^ine alors que l'on a $x =
\tilde{x}$.

Montrons que $(Q/F,f)$ est cohomologiquement propre relativement \`a~$S$ en dimension $\leqslant 0$. Soient $Y''$ un sch\'ema \'etale
sur~$Y'$ et $x''$ un \'el\'ement de~$\H^0(X'',Q''/F'')$. \`A~$x''$
est associ\'e un torseur $P''$ sur~$X''$ de groupe $F''$. Comme
$(F,f)$ est cohomologiquement propre relativement \`a~$S$ en
dimension $\leqslant 1$, on peut trouver des morphismes \'etales
surjectifs $Y''_1 \to Y''$ et $Y_1 \to Y$, tels que l'on ait un
morphisme $Y''_1 \to Y'_1$, et un torseur $P$ sur~$X_1$ de groupe
$F_1$ dont l'image inverse sur~$X''_1$ soit isomorphe \`a l'image
inverse de~$P''$. Il r\'esulte alors du fait que $(G,f)$ est
cohomologiquement propre relativement \`a~$S$ en dimension
$\leqslant 0$ que l'on peut m\^eme choisir $Y''_1$ et $Y_1$ tels que
le torseur d\'eduit de~$P$ par extension du groupe structural $F_1
\to G_1$ soit isomorphe \`a $Q_1$; il correspond \`a~$P$ un
\'el\'ement $x$ de~$\H^0(X_1,Q_1/F_1)$, dont l'image dans
$\H^0(X''_1,Q''_1/F''_1)$ est isomorphe \`a l'image inverse de
$x''$, ce qui ach\`eve la d\'emonstration.

\begin{proposition}
\label{XIII.1.14}
Soient $f \colon X \to S$ un $S$-sch\'ema, $F$ un faisceau
d'ensembles ou de groupes sur~$X$ (\resp un faisceau de
ind-$\LL$-groupes, o\`u $\LL$ est un ensemble de nombres
premiers). Supposons $F$ localement constant, $(F,f)$
cohomologiquement propre en dimension $\leqslant 0$ (\resp en
dimension $\leqslant 1$) et $f$ localement $0$-acyclique
(\resp localement $1$-asph\'erique pour $\LL$) \textup{(SGA~4 XV~1.11)}.
Alors, pour toute sp\'ecialisation $\overline{s}_1 \to
\overline{s}_2$ de points g\'eom\'etriques de~$S$, le morphisme de
sp\'ecialisation \textup{(SGA~4 VIII~7.1)}
$$
a_0 \colon (f_* F)_{\overline{s}_2} \to (f_*
F)_{\overline{s}_1}
$$
est un isomorphisme, et, si $F$ est un faisceau en groupes, le
morphisme
$$
a_1 \colon (\R^1 f_* F)_{\overline{s}_2} \to (\R^1 f_*
F)_{\overline{s}_1}
$$
est injectif (\resp les morphismes $a_0$ et $a_1$ sont des
isomorphismes).
\end{proposition}

La
\marginpar{367}
d\'emonstration s'obtient en recopiant mot \`a mot celle de SGA~4
XVI~2.3, mais en y rempla\c cant l'expression \og propre\fg par
l'expression \og cohomologiquement propre\fg.

\ifthenelse{\boolean{orig}}
{}
{\enlargethispage{.5cm}}

\begin{corollaire}
\label{XIII.1.15}
Soient $f \colon X \to S$ un morphisme, $\Phi$ un champ sur~$X$, $\LL$
un ensemble de nombres premiers. Supposons que, pour tout sch\'ema
$X_1$ \'etale sur~$X$ et pour tout couple d'objets $x$, $y$ de
$\Phi_{X_1}$, le faisceau $\SheafHom_{X_1}(x,y)$ soit localement
constant, que le faisceau $\SheafAut_{X_1}(x)$ soit un
ind-$\LL$-groupe localement constant, et que le faisceau des
sous-gerbes maximales $S\Phi$ de~$\Phi$ \textup{\cite[III 2.1.7]{XIII.1}} soit localement
constant. Supposons que $(\Phi,f)$ soit cohomologiquement propre en
dimension $\leqslant 1$ et que $f$ soit localement $1$-asph\'erique
pour $\LL$. Alors, pour toute sp\'ecialisation $\overline{s}_1 \to
\overline{s}_2$ de points g\'eom\'etriques de~$S$, le morphisme de
sp\'ecialisation
$$
a \colon (f_* \Phi)_{\overline{s}_2} \to (f_*
\Phi)_{\overline{s}_1}
$$
est une \'equivalence de cat\'egories.
\end{corollaire}

Soient $\overline{S}_1$ (\resp $\overline{S}_2$) le localis\'e
strict de~$S$ en $\overline{s}_1$ (\resp le localis\'e strict de~$S$
en $\overline{s}_2$), $\overline{X}_2$, $\overline{\Phi}_2$
(\resp $\overline{X}_1$, $\overline{\Phi}_1$) les images inverses de
$X_2$, $\Phi_2$ sur~$\overline{S}_2$ (\resp de~$X_1$, $\Phi_1$ sur~$\overline{S}_1$) et consid\'erons le carr\'e cart\'esien
$$
\xymatrix{ \overline{X}_1 \ar[r]^h \ar[d]_{\overline{f}_1} &
\overline{X}_2 \ar[d]^{\overline{f}_2} \\ \overline{S}_1 \ar[r]^g &
\overline{S}_2 }
\text{.}
$$
On doit montrer que le foncteur
$$
\varphi \colon \overline{\Phi}_2(\overline{X}_2) \to \overline{\Phi}_1
(\overline{X}_1)
$$
est une \'equivalence. Le foncteur $\varphi$ est pleinement
fid\`ele; soient
\ifthenelse{\boolean{orig}}
{en}
{}
en
\marginpar{368}
effet $x$, $y$ deux objets de~$(\overline{\Phi}_2)_{\overline{X}_2}$;
le morphisme canonique
$$
\Hom_{\overline{X}_2}(x,y) \to
\Hom_{\overline{X}_1}(\varphi(x),\varphi(y))
$$
s'identifie au morphisme canonique
$$
\H^0(\overline{X}_2, \SheafHom_{\overline{X}_2}(x,y)) \to
\H^0(\overline{X}_1,h^*(\SheafHom_{\overline{X}_2}(x,y)) \text{.}
$$
Ce morphisme est un isomorphisme d'apr\`es~\Ref{XIII.1.14}.

Montrons que $\varphi$ est une \'equivalence. Soient $x_1$ un objet
de~$\overline{\Phi}_1(\overline{X}_1)$ et $G_1$ la sous-gerbe maximale
de~$\overline{\Phi}_1$ engendr\'ee par $x_1$. Le morphisme
$$
\H^0(\overline{X}_2, S\overline{\Phi}_2) \to \H^0(\overline{X}_1,
h^*(S \overline{\Phi}_2)) =
\H^0(\overline{X}_1,S\overline{\Phi}_1)
$$
est bijectif, et il existe donc une sous-gerbe maximale $G_2$ de
$\overline{\Phi}_2$ telle que $h^* G_2$ soit isomorphe \`a
$G_1$. Il suffit alors de prouver que le foncteur
$$
G_2 \to h_* h^* G_2
$$
est une \'equivalence. Mais, sous cette forme, la question est
locale pour la topologie \'etale sur~$\overline{X}_2$. On peut donc
supposer que $G_2$ est une gerbe de torseurs sous le groupe des
automorphismes d'un objet de~$G_2$, cas o\`u l'assertion r\'esulte
de~\Ref{XIII.1.14}.

\begin{corollaire}
\label{XIII.1.16}
Les hypoth\`eses sont celles de~\Ref{XIII.1.14}. Si l'on suppose de
plus que $F$ est un faisceau d'ensembles (\resp de ind-$\LL$-groupes)
et que $f_* F$ (\resp $\R^1 f_* F$) est constructible, alors
$f_* F$ (\resp $\R^1 f_* F$) est localement constant.
\end{corollaire}

Le corollaire r\'esulte de~\Ref{XIII.1.14} gr\^ace \`a SGA~4
IX~2.11.

\begin{remarque}
\label{XIII.1.17}
Rappelons que la condition $f$ localement $0$-acyclique est satisfaite
si $f$ est plat \`a fibres s\'eparables, $X$ et $Y$ localement
noeth\'eriens (SGA~4 XV~4.1), et que la condition $f$ localement
$1$-asph\'erique
\marginpar{369}
pour $\LL$ est satisfaite si $f$ est lisse, $\LL$ \'etant l'ensemble
des nombres premiers distincts des caract\'eristiques
r\'esiduelles de~$S$ (SGA~4 XV~2.1).
\end{remarque}

\section[Un cas particulier de propret\'e cohomologique]{Un cas particulier de propret\'e cohomologique: diviseurs
\`a croisements normaux relatifs}
\label{XIII.2}

\setcounter{subsection}{-1}

\subsection{}
\label{XIII.2.0}
Soient $R$ un anneau de valuation discr\`ete de corps des fractions
$K$ et $L$ une $K$-alg\`ebre \'etale; $L$ est alors produit direct
d'un nombre fini de corps $L_i$, o\`u $L_i$ est une extension
\'etale de~$K$. Si $L'_i$ d\'esigne l'extension galoisienne
engendr\'ee par $L_i$ dans une cl\^oture alg\'ebrique de~$L_i$,
on dit que $L$ est \emph{mod\'er\'ement ramifi\'ee}
\index{mod\'er\'ement ramifi\'ee (alg\`ebre)|hyperpage}%
sur~$R$ si les $L'_i$ sont des extensions mod\'er\'ement
ramifi\'ees au sens de~X~\Ref{X.3}, \ie si un groupe d'inertie~$I_i$ de~$L'_i|K$ est d'ordre premier \`a la caract\'eristique
r\'esiduelle $p$ de~$R$.

On sait que $I_i$ est en tous cas extension d'un groupe cyclique
d'ordre premier \`a $p$ par un $p$-groupe. (Cela r\'esulte de \cite[ch\ptbl IV prop\ptbl 7 cor\ptbl 4]{XIII.5} lorsque l'extension r\'esiduelle de~$R$ est
s\'eparable. La d\'emonstration donn\'ee dans
\loccit s'\'etend au cas g\'en\'eral de la fa\c con
suivante. Reprenons les hypoth\`eses et les notations de
\loccit mais sans supposer l'extension r\'esiduelle
s\'eparable. Soit $H_i$ le sous-groupe du groupe d'inertie $G_0$,
ensemble des \'el\'ements $s$ de~$G_0$ tels que l'on ait $s \pi /
\pi \in U^i$ pour toute uniformisante $\pi$ de~$A_L$. On v\'erifie
alors que $G_0/H_1$ est un groupe d'ordre premier \`a $p$ et que,
pour $i \geqslant 1$, les $H^i/H^{i+1}$ sont des $p$-groupes, d'o\`u
l'on d\'eduit le r\'esultat annonc\'e.)

\subsubsection{}
\label{XIII.2.0.1}
Dans le cas o\`u $R$ est strictement local, on a la
caract\'erisation simple suivante: la $K$-alg\`ebre $L$ est
mod\'er\'ement ramifi\'ee sur~$R$ si et seulement si les
$[L_i:K]$ sont premiers \`a $p$. De plus, si $L$ est
mod\'er\'ement ramifi\'ee sur~$R$, les $L_i$ sont des extensions
cycliques de~$K$. En effet, lorsque $R$ est strictement local, $I_i$
est \'egal au groupe de Galois
\marginpar{370}
de~$L'_i$ sur~$K$. Comme on vient de le rappeler $I_i$ est extension
d'un groupe cyclique d'ordre premier \`a $p$ par un $p$-groupe. Si
l'on suppose $L'_i$ mod\'er\'ement ramifi\'ee sur~$R$, $I_i$ est
alors un groupe cyclique d'ordre premier \`a $p$. Il en r\'esulte
que $[L_i:K]$ est premier \`a $p$ et que l'on a $L_i =
L'_i$. Inversement, si $[L_i:K]$ est premier \`a $p$, $I_i$ ne peut
contenir de~$p$-sous-groupe distingu\'e non trivial; $I_i$ est donc
un groupe cyclique d'ordre premier \`a $p$, ce qui prouve que $L_i$
est mod\'er\'ement ramifi\'ee sur~$R$.

\subsubsection{}
\label{XIII.2.0.2}
Soient $R$ un anneau de valuation discr\`ete de corps des fractions
$K$, $L$ une $K$-alg\`ebre \'etale, et soient $\overline{R}$ un
localis\'e strict de~$R$, $\overline{K}$ son corps de fractions,
$\overline{L} = L \otimes_K \overline{K}$. Alors, pour que $L$ soit
mod\'er\'ement ramifi\'ee sur~$R$, il faut et il suffit que
$\overline{L}$ soit mod\'er\'ement ramifi\'ee sur~$\overline{R}$. On se ram\`ene en effet au cas o\`u $L$ est un
corps. Soient alors $\overline{L} = \prod_i \overline{L}_i$, o\`u
les $\overline{L}_i$ sont des corps extensions de~$\overline{K}$; si
$L'$ est l'extension galoisienne engendr\'ee par $L$, et, si
\ifthenelse{\boolean{orig}}
{$L' = L' \otimes_K \overline{K}$,}
{$\overline{L'} = L' \otimes_K \overline{K}$,}
on a de m\^eme une d\'ecomposition de~$\overline{L'}$ en produit
de corps,
\ifthenelse{\boolean{orig}}
{$L' = \prod_j \overline{L'_j}$}
{$\overline{L'} = \prod_j \overline{L'_j}$}
et chaque $\overline{L}_i$ est sous-extension d'au moins l'un des
$\overline{L'_j}$. Comme $L'$ est une extension galoisienne de~$K$,
les $\overline{L'_j}$ sont des extensions galoisiennes de
$\overline{K}$. Supposons $L$ mod\'er\'ement ramifi\'ee sur~$R$;
comme le groupe de Galois de~$\overline{L'_j}|K$ est isomorphe au
groupe d'inertie de~$L'|K$, les $L'_j$ sont aussi mod\'er\'ement
ramifi\'ees sur~$\overline{R}$, et il en est donc de m\^eme des
$\overline{L}_i$. Inversement, supposons $\overline{L}$
mod\'er\'ement ramifi\'ee sur~$\overline{R}$. Pour chaque $j$,
soit $v_j$ la valuation discr\`ete de~$\overline{L'_j}$ qui prolonge
la valuation de~$\overline{K}$ et notons encore $v_j$ la valuation
induite sur~$L'$. Quand $j$ varie, $v_j$ parcourt l'ensemble des
valuations de~$L$ qui prolongent la valuation de~$K$. Soient $G =
\Gal(L'|K)$, $H=\Gal(L'|L)$, $I_j$ le groupe d'inertie de~$L'|K$ en
$v_j$, $J_j$ le groupe d'inertie de~$L'|L$ en $v_j$. Le groupe $I_i$
est extension d'un groupe cyclique d'ordre premier \`a $p$ par un
$p$-groupe $P_j$. Comme les $\overline{L}_i$ sont mod\'er\'ement
ramifi\'ees sur~$\overline{R}$, $I_j/J_j$ est d'ordre premier \`a
$p$, donc on a $P_j \subset J_j$. Par suite le groupe $H$ contient
tous les $P_j$ donc aussi le groupe $P$ engendr\'e par
\marginpar{371}
les $P_j$ pour $j$ variable. Mais le groupe~$P$ est invariant dans $G$
car un automorphisme int\'erieur de~$G$ transforme les~$I_j$ entre
eux donc aussi les $P_j$ entre eux. Il en r\'esulte que $P$ est un
sous-groupe de~$H$ distingu\'e dans~$G$, donc, puisque $L'$ est
l'extension galoisienne engendr\'ee par $L$, que l'on a $P=1$, ce
qui prouve que $L$ est mod\'er\'ement ramifi\'ee sur~$R$.

Soient plus g\'en\'eralement $R \to R'$ un morphisme d'anneaux de
valuation discr\`ete tel que l'image d'une uniformisante $\pi$ de
$R$ soit une uniformisante $\pi'$ de~$R'$ et que l'extension
r\'esiduelle $\kres(R')$ soit une extension s\'eparable de
$\kres(R)$. Soient $K$ le corps des fractions de~$R$, $K'$ le corps des
fractions de~$R'$, $L$ une $K$-alg\`ebre \'etale, $L' = L
\otimes_K K'$. Alors, pour que $L$ soit mod\'er\'ement
ramifi\'ee sur~$R$, il faut et il suffit que $L'$ soit
mod\'er\'ement ramifi\'ee sur~$R'$. On peut en effet supposer
$R$ et $R'$ strictement locaux. D'apr\`es~\Ref{XIII.2.0.1} il suffit
de prouver que, lorsque $L$ est un corps, il en est de m\^eme de~$L'$. Soient $\tilde{R}$ le normalis\'e de~$R$ dans $L$,
$\tilde{\pi}$ une uniformisante de~$\tilde{R}$, $\tilde{R}' =
\tilde{R} \otimes_R R'$. L'extension $\kres(\tilde{R})|\kres(R)$ \'etant
radicielle et l'extension $\kres(R')|\kres(R)$ \'etale,
$\kres(\tilde{R})\otimes_{\kres(R)}\kres(R')$ est un corps [EGA IV 4.3.2 et
4.3.5]. Ceci prouve que $R'$ est un anneau local, et, comme $\pi$ a
pour image $\pi'$ dans $R'$, on a $\kres(\tilde{R}')= R'/(\tilde{\pi})$;
par suite $R'$ est un anneau de valuation discr\`ete \cite[ch\ptbl I \S 2
prop\ptbl 2]{XIII.5} donc $L'$ est un corps.

\subsubsection{}
\label{XIII.2.0.3}
Par r\'eduction au cas strictement local, on voit qu'une
sous-alg\`ebre d'une alg\`ebre mod\'er\'ement ramifi\'ee est
mod\'er\'ement ramifi\'ee, que le produit tensoriel de deux
alg\`ebres mod\'er\'ement ramifi\'ees est mod\'er\'ement
ramifi\'e, qu'une alg\`ebre mod\'er\'ement ramifi\'ee le
reste apr\`es extension de l'anneau de valuation discr\`ete,
qu'une alg\`ebre qui devient mod\'er\'ement ramifi\'ee
apr\`es une extension mod\'er\'ement ramifi\'ee est
mod\'er\'ement ramifi\'ee.

\subsection{}
\label{XIII.2.1}
Soient
\marginpar{372}
$X$ un $S$-sch\'ema, $D$ un diviseur $\geq 0$ sur~$X$.
Rappelons (SGA~5 II~4.2) qu'on dit que $D$ est strictement \`a
\emph{croisements normaux relativement \`a}~$S$
\index{croisements normaux (diviseur \`a)|hyperpage}%
\index{diviseur \`a croisements normaux|hyperpage}%
\index{strictement \`a croisements normaux (diviseur)|hyperpage}%
s'il existe une famille finie $(f_i)_{i\in I}$ d'\'el\'ements de
$\Gamma(X,\cal{O}_X)$, telle que l'on ait $D=\sum_{i\in
I}\divisor(f_i)$ et que la condition suivante soit r\'ealis\'ee:

\setcounter{subsubsection}{-1}
\subsubsection{}
\label{XIII.2.1.0}
Pour tout point $x$ de~$\Supp D$, $X$ est lisse sur~$S$ en $x$, et, si
l'on note $I(x)$ l'ensemble des $i\in I$ tels que $f_i(x)=0$, le
sous-sch\'ema $V((f_i)_{i\in I(x)})$ est lisse sur~$S$ de
codimension
\ifthenelse{\boolean{orig}}
{$\mathrm{card}.\,I(x)$}
{$\card I(x)$}
dans~$X$.

Le diviseur $D$ est dit \`a \emph{croisements normaux relativement
\`a} $S$ si, localement sur~$X$ pour la topologie \'etale, il est
strictement \`a croisements normaux.

Soit $D$ un diviseur \`a croisements normaux relativement
\`a~$S$. On pose $Y=\Supp D$, $U=X-Y$, et on note $i\colon U\to X$
l'immersion canonique. Pour tout point g\'eom\'etrique
$\overline{s}$ de~$S$ et pour tout point maximal $y$ de la fibre
g\'eom\'etrique $Y_{\overline{s}}$, l'anneau
$R=\cal{O}_{X_{\overline{s}},y}$ est un anneau de valuation
discr\`ete.

Dans la suite de ce num\'ero, nous utiliserons la d\'efinition
technique suivante:

\begin{subdefinition}
\label{XIII.2.1.1}
Soit $F$ un faisceau d'ensembles sur~$U$. On dit que $F$ est
\emph{mod\'er\'ement ramifi\'e sur}~$X$
\index{mod\'er\'ement ramifi\'e (faisceau)|hyperpage}%
(\emph{le long de}~$D$) \emph{relativement \`a}~$S$ si, pour tout
point g\'eom\'etrique $\overline{s}$ de~$S$, la condition suivante
est satisfaite:

Pour tout point maximal $y$ de~$Y_{\overline{s}}$, la restriction de
$F$ au corps des fractions $K$ de~$\cal{O}_{X_{\overline{s}}}$ est
repr\'esentable par le spectre d'une $K$-alg\`ebre \'etale $L$,
mod\'er\'ement ramifi\'ee sur~$\cal{O}_{X_{\overline{s},y}}$.
\end{subdefinition}

Le plus souvent, quand il ne pourra en r\'esulter de confusion, nous
omettrons la mention de~$D$ dans la terminologie.

\begin{subdefinition}
\label{XIII.2.1.2}
Si
\marginpar{373}
$F$ est un faisceau en groupes sur~$U$, mod\'er\'ement
ramifi\'e sur~$X$ relativement \`a~$S$, on d\'esigne par
$$
\H^1_{\tame}(U,F)
$$
\label{indnot:mc}\oldindexnot{$\H^1_{\tame}(U,F)$|hyperpage}%
le sous-ensemble de~$\H^1(U,F)$ form\'e des classes de torseurs sous
$F$ qui sont mod\'er\'ement ramifi\'es sur~$X$ relativement
\`a~$S$.
\end{subdefinition}

Soit
$$
\xymatrix{ U \ar[r]^i \ar[d]_g & X \ar[dl]^f \\ T}
$$
un diagramme commutatif de~$S$-sch\'emas, avec $i$ comme
dans~\Ref{XIII.2.1}; on d\'esigne par
$$
\R^1_{\tame} g_* F
$$
\label{indnot:md}\oldindexnot{$\R^1_{\tame} g_* F$|hyperpage}%
le faisceau sur~$T$ associ\'e au pr\'efaisceau
$T'\mto\H^1_{\tame}(U',F)$, o\`u $T'$ parcourt les sch\'emas
\'etales au-dessus de~$T$ et o\`u $U'=U\times_T T'$; $\R^1_{\tame} g_*F$
est un sous-faisceau de~$\R^1g_*F$.

Notons que, si $g$ est coh\'erent, si $\overline{t}$ est un point
g\'eom\'etrique de~$T$, $\overline{T}$ le localis\'e strict de
$T$ en $\overline{t}$, $\overline{U}=U\times_T\overline{T}$, on a un
isomorphisme
\begin{equation*}
\label{eq:XIII.2.1.2.1}
\tag{\thesubsubsection.1} (\R^1_{\tame} g_*F)_{\overline{t}}\simeq
\H^1_{\tame}(\overline{U},\overline{F})
\end{equation*}

\subsubsection{}
\label{XIII.2.1.3}
Soit $C_t((U,X)/S)$ ou simplement $C_t$
\label{indnot:me}\oldindexnot{$C_t((U,X)/S)$ ou $C_t$|hyperpage}%
la cat\'egorie des rev\^etements \'etales de~$U$ qui sont
mod\'er\'ement ramifi\'es sur~$X$ relativement
\`a~$S$. Supposons $U$ connexe et soit $a$ un point
g\'eom\'etrique de~$U$. Soit $\Gamma_t$ le foncteur qui, \`a un
rev\^etement \'etale~$U'$ de~$U$ mod\'er\'ement ramifi\'e
sur~$X$ relativement \`a~$S$ fait correspondre l'ensemble des points
g\'eom\'etriques de~$U'$ au-dessus de~$a$. De~\Ref{XIII.2.0}
r\'esulte que le couple $(C_t,\Gamma_t)$ satisfait aux axiomes
($G_1$) \`a ($G_6$) de V~\Ref{V.4}. Par suite $\Gamma_t$ est
repr\'esentable par un pro-objet qu'on appelle le
\emph{rev\^etement universel mod\'er\'ement ramifi\'e
\marginpar{374}
de~$(U,X)$ relativement \`a~$S$ ponctu\'e en~$a$}.
\index{mod\'er\'ement ramifi\'e (rev\^etement universel)|hyperpage}%
\index{revetement universel moderement ramifie@rev\^etement universel mod\'er\'ement ramifi\'e|hyperpage}%
Le groupe oppos\'e au groupe des $U$-automorphismes du
rev\^etement universel mod\'er\'ement ramifi\'e est appel\'e
le \emph{groupe fondamental mod\'er\'ement ramifi\'e}
\index{groupe fondamental mod\'er\'ement ramifi\'e|hyperpage}%
\index{mod\'er\'ement ramifi\'e (groupe fondamental)|hyperpage}%
et not\'e
$$
\pi_1^{\tame}((U,X)/S,a)\text{ ou simplement }\pi_1^{\tame}(U,a)\text{ ou
m\^eme }\pi_1^{\tame}(U).
$$
\label{indnot:mf}\oldindexnot{$\pi_1^{\tame}((U,X)/S,a)$ ou $\pi_1^{\tame}(U,a)$ ou $\pi_1^{\tame}(U)$|hyperpage}%
C'est \'evidemment un quotient du groupe fondamental $\pi_1(U,a)$
(V~\Ref{V.6.9}).

\subsubsection{}
\label{XIII.2.1.4}
Soient $F$ un faisceau en groupes sur~$U$, $P$ un torseur \`a droite
de groupe $F$, $Q$ un torseur \`a gauche de groupe $F$, et supposons
$P$ et $Q$ mod\'er\'ement ramifi\'es sur~$X$ relativement \`a~$S$. Alors il en est de m\^eme de~$P\stackrel{F}{\wedge}Q$. On se
ram\`ene en effet \`a montrer que, si $R$ est un anneau de
valuation discr\`ete de corps des fractions $K$ et si~$F$ est un
sch\'ema en groupes \'etale fini sur~$K$, et, si $P$ et $Q$ sont
deux torseurs sous~$F$ mod\'er\'ement ramifi\'es sur~$R$, alors
il en est de m\^eme de~$P\stackrel{F}{\wedge} Q$. Or
$T=P\stackrel{F}{\wedge} Q$ est un quotient de~$P\times_K Q$. Si $L$,
$M$, $N$ d\'esignent les $K$-alg\`ebres repr\'esentant
respectivement $T$, $P$, $Q$, alors $L$ est une sous-alg\`ebre de
$M\otimes_K N$, et il r\'esulte de~\Ref{XIII.2.0.3} que $L$ est
mod\'er\'ement ramifi\'ee sur~$R$.

On d\'eduit de ce qui pr\'ec\`ede que, si $F$ est un faisceau en
groupes sur~$U$ et s'il existe un torseur de groupe $F$,
mod\'er\'ement ramifi\'e sur~$X$ relativement \`a~$S$, alors
$F$ est mod\'er\'ement ramifi\'e sur~$X$ relativement \`a~$S$. En effet le torseur $P^{\circ}$ oppos\'e de~$P$ est
mod\'er\'ement ramifi\'e sur~$X$ relativement \`a~$S$,
puisqu'il est isomorphe \`a $P$ en tant que faisceau d'ensembles. Si
$\leftexp{P\mkern-4mu}{F}$ est le groupe tordu de~$F$ par $P$, on a un
isomorphisme
$$
F\simeq P^{\circ}\stackrel{\leftexp{P\mkern-4mu}{F}}{\wedge}P
$$
et par suite $F$ est mod\'er\'ement ramifi\'e sur~$X$
relativement \`a~$S$.

On voit comme pr\'ec\'edemment que, si $F\to F'$ est un morphisme
de faisceaux en groupes sur~$U$, mod\'er\'ement ramifi\'es sur~$X$ relativement \`a~$S$, et si $P$ est un torseur sous $F$
mod\'er\'ement ramifi\'e sur~$X$ relativement
\marginpar{375}
\`a~$S$, alors le torseur $P'$ d\'eduit de~$P$ par l'extension du
groupe structural $F\to F'$ est mod\'er\'ement ramifi\'e sur~$X$
relativement \`a~$S$.

En particulier le morphisme canonique
$$
\H^1(U,F)\to \H^1(U,F')
$$
donne par restriction \`a $\H^1_{\tame}(U,F)$ un morphisme canonique
$$
\H^1_{\tame}(U,F)\to\H^1_{\tame}(U,F')
$$

\subsubsection{}
\label{XIII.2.1.5}
Soient $S'\to S$ un morphisme et notons $U'$ (\resp $X'$, etc.)
l'image inverse de~$U$ (\resp $X$, etc.) sur~$S'$. Si $F$ est un
faisceau d'ensembles sur~$U$ mod\'er\'ement ramifi\'e sur~$X$
relativement \`a~$S$, il r\'esulte de la
d\'efinition~\Ref{XIII.2.1.1} et de~\Ref{XIII.2.0.3} que $F'$ est
mod\'er\'ement ramifi\'e sur~$X'$ relativement \`a~$S'$.

Si maintenant $F$ est un faisceau en groupes sur~$U$, l'image inverse
sur~$S'$ d'un torseur sous $F$ mod\'er\'ement ramifi\'e sur~$X$
relativement \`a~$S$ est un torseur sous $F'$ mod\'er\'ement
ramifi\'e sur~$X'$ relativement \`a~$S'$. En particulier on a un
foncteur canonique
\begin{equation*}
\label{eq:XIII.2.1.5.1}
\tag{\thesubsubsection.1} C_t((U,X)/S)\to C_t((U',X')/S').
\end{equation*}
Supposons $U$ et $U'$ connexes et soient $a$ un point
g\'eom\'etrique de~$U$, $a'$ un point g\'eom\'etrique de~$U'$
au-dessus de~$a$; on d\'eduit de ce qui pr\'ec\`ede un morphisme
canonique
\begin{equation*}
\label{eq:XIII.2.1.5.2}
\tag{\thesubsubsection.2} \pi_1^{\tame}(U',a')\to\pi_1^{\tame}(U,a).
\end{equation*}

Si $S'\to S$ est un morphisme et $h\colon T'\to T$ la projection
canonique, le morphisme
$$
h^*(\R^1g_*F)\to\R^1g_*'F'
$$
donne par restriction un morphisme canonique
\begin{equation*}
\label{eq:XIII.2.1.5.3}
\tag{\thesubsubsection.3} h^*(\R^1_{\tame}g_*F)\to\R^1_{\tame}g_*'F'
\end{equation*}

\subsubsection{}
\label{XIII.2.1.6}
Soit
\marginpar{376}
$F$ un faisceau en groupes sur~$U$, mod\'er\'ement
ramifi\'e sur~$X$ relativement \`a~$S$. Les notations \'etant
celles de~\Ref{XIII.2.1.2}, on a des suites exactes canoniques:
\begin{equation*}
\label{eq:XIII.2.1.6.1}
\tag{\thesubsubsection.1}
\begin{array}{ccccccc}
1 & \longrightarrow & \H^1(X,i_*F)& \longrightarrow & \H^1_{\tame}(U,F) &
\longrightarrow & \H^0(X,\R^1_{\tame} i_*F)\\[10pt] 1 &
\longrightarrow & \R^1f_*(i_*F) & \longrightarrow & \R^1_{\tame}g_*F &
\longrightarrow & f_*(\R^1_{\tame}i_*F).\end{array}
\end{equation*}

La premi\`ere s'obtient \`a partir de la suite exacte
(SGA~4 III~3.2):
$$
1\to \H^1(X,i_*F)\to \H^1(U,F)\to
\H^0(X,\R^1i_*F).
$$
Il suffit en effet de montrer que l'image de~$\H^1(X,i_*F)$ dans
$\H^1(U,F)$ est en fait contenue dans $\H^1_{\tame}(U,F)$ et que l'image de
$\H^1_{\tame}(U,F)$ dans $\H^0(X,\R^1i_*F)$ est en fait contenue dans
$\H^0(X,\R^1_{\tame}i_*F)$. Or l'image inverse sur~$U$ d'un torseur sous
$i_*F$ est un torseur sous $i^*i_*F$ qui est \'evidemment
mod\'er\'ement ramifi\'e sur~$X$ relativement \`a~$S$, donc il
en est de m\^eme apr\`es l'extension du groupe structural
$i^*i_*F\to F$, ce qui prouve l'existence de la fl\`eche
$\H^1(X,i_*F)\to\H^1_{\tame}(U,F)$. Le fait que l'image de~$\H^1_{\tame}(U,F)$
dans $\H^0(X,\R^1i_*F)$ soit contenue dans $\H^0(X,\R^1_{\tame}i_*F)$
r\'esulte aussit\^ot de la d\'efinition de~$\R^1_{\tame}i_*F$. Ceci
prouve l'existence de la premi\`ere suite exacte, et la deuxi\`eme
s'en d\'eduit par localisation.

\subsection{}
\label{XIII.2.2}
On conserve les notations de~\Ref{XIII.2.1}. Nous allons d\'efinir
une notion d'objet mod\'er\'ement ramifi\'e d'un champ $\Phi$
sur~$U$ lorsque celui-ci est donn\'e, localement \emph{sur} $X$ et
sur~$S$ pour la topologie \'etale, comme \emph{image inverse d'un
champ} $\Psi$ \emph{sur}~$S$.

Soit d'abord $G$ une gerbe sur~$U$ et supposons donn\'e un morphisme
\'etale surjectif $S_1\to S$, un morphisme \'etale surjectif
$X_2\to X\times_S S_1$, une gerbe triviale $H$ sur~$S_1$ et un
isomorphisme
$$
G|{U_2}\to H|{U_2},
$$
o\`u
\marginpar{377}
$U_1=U\times_X X_1$, $U_2=U\times_X X_2$. Quand on choisit une
trivialisation de~$H|{X_2}$, l'isomorphisme ci-dessus identifie
$G|{U_2}$ au champ des torseurs sous un faisceau en groupes $F$. On
dit qu'un \'el\'ement $x$ de~$G_U$ est mod\'er\'ement
ramifi\'e sur~$X$ relativement
\ifthenelse{\boolean{orig}}
{\`a~$S$,}
{\`a~$S$}
si la restriction de~$x$
\`a $U_2$ est un torseur mod\'er\'ement ramifi\'e sur~$X$
relativement \`a~$S$. D'apr\`es~\Ref{XIII.2.1.4} cette notion ne
d\'epend pas de la fa\c con dont on a trivialis\'e~$H|X_2$.

Soit maintenant $\Phi$ un champ sur~$U$ et supposons donn\'es un
morphisme \'etale surjectif $S_1\to S$, un morphisme \'etale
surjectif $X_2\to X\times_S S_1$, un champ $\Psi$ sur~$S_1$ et un
isomorphisme
$$
i\colon \Phi|U_2\to \Psi|U_2.
$$
Soit $x$ un \'el\'ement de~$\Phi_U$, $G_x$ la sous-gerbe maximale
de~$\Phi$ engendr\'ee par $x$
\ifthenelse{\boolean{orig}}
{[III~2.1.7]}
{\cite[III~2.1.7]{XIII.1}}, $S\Phi$ le faisceau des
sous-gerbes maximales de~$\Phi$. L'isomorphisme $i$ induit un
isomorphisme
$$
S\Phi|U_2\to S\Psi|U_2.
$$
Il r\'esulte de~\Ref{XIII.5.7} que, quitte \`a remplacer $S_1$ par
une extension \'etale surjective, on a une unique sous-gerbe
maximale $H$ de~$\Psi$, que l'on peut supposer triviale, telle que $i$
d\'efinisse un isomorphisme
$$
G_x|U_2\to H|U_2.
$$
On dit que l'\'el\'ement $x$ est mod\'er\'ement ramifi\'e
sur~$X$ relativement \`a~$S$ s'il l'est en tant qu'\'el\'ement
de~$G_x$ munie de l'isomorphisme ci-dessus.

\subsubsection{}
\label{XIII.2.2.1}
Soit $\Phi$ un champ sur~$U$ donn\'e, localement sur~$X$ et $S$,
comme image inverse d'un champ sur~$S$, et soit
$$
\xymatrix{ U \ar[r]^i
\ar[d]_g & X \ar[dl]^f \\ T }
$$
un diagramme comme dans~\Ref{XIII.2.1.2}. Pour tout sch\'ema $T'$
\'etale sur~$T$, si $U'=U\times_T T'$, on consid\`ere le
sous-ensemble $(g_*^{\tame}\Phi)_{T'}$
\label{indnot:mg}\oldindexnot{$(g_*^{\tame}\Phi)_{T'}$|hyperpage}%
de~$(g_*\Phi)_{T'}=\Phi_{U'}$
form\'e
\marginpar{378}
des \'el\'ements de~$\Phi_{U'}$ qui sont mod\'er\'ement
ramifi\'es
\index{element moderement@\'el\'ement mod\'er\'ement ramifi\'e d'un champ|hyperpage}%
sur~$X$ relativement \`a~$S$. On appelle image directe
mod\'er\'ement ramifi\'ee
\index{image directe mod\'er\'ement ramifi\'ee d'un champ|hyperpage}%
\index{mod\'er\'ement ramifi\'ee (image directe d'un champ)|hyperpage}%
de~$\Phi$ par~$g$, et on note
$$
g_*^{\tame}\Phi
$$
la sous-cat\'egorie pleine de~$g_*\Phi$ dont les objets au-dessus
d'un sch\'ema $T'$ \'etale sur~$T$ sont les \'el\'ements de
$(g_*^{\tame}\Phi)_{T'}$. Il est clair que $g_*^{\tame}\Phi$ est un sous-champ de
$g_*\Phi$.

\subsubsection{}
\label{XIII.2.2.2}
Par r\'eduction au cas d'un champ de torseurs, on voit que, si
$\Phi$ est un champ sur~$U$ qui est localement pour la topologie
\'etale de~$S$ et $X$ image inverse d'un champ sur~$S$, le morphisme
canonique
$$
h^*(g_*\Phi)\to g_*'\Phi'
$$
donne par restriction un morphisme canonique
$$
h^*(g_*^{\tame}\Phi)\to g_*^{\prime\tame}\Phi'
$$

\begin{remarques}
\label{XIII.2.3}

a) Si $F$ est un faisceau d'ensembles localement constant
constructible sur~$U$, pour que $F$ soit mod\'er\'ement
ramifi\'e sur~$X$ relativement \`a~$S$, il suffit que la condition
de~\Ref{XIII.2.1.1} soit satisfaite pour les points
g\'eom\'etriques de~$S$ au-dessus des points maximaux de~$S$. Pour
le voir on peut supposer le diviseur $D$ strictement \`a croisements
normaux. Le faisceau $F$ est repr\'esentable par un rev\^etement
\'etale $V$ de~$U$. Si $\overline{s}$ est un point
g\'eom\'etrique de~$S$, $y$ un point maximal de
$Y_{\overline{s}}$, on note $\overline{S}$ le localis\'e strict de
$S$ en $\overline{s}$, $\overline{X}$ le localis\'e strict de~$X$ en
$\overline{y}$, $\overline{U}=U\times_X \overline{X}$,
$\overline{V}=V\times_X\overline{X}$. Si la condition
de~\Ref{XIII.2.1.1} est satisfaite aux points g\'eom\'etriques
au-dessus des points maximaux de~$S$, il r\'esulte de~\Ref{XIII.5.5}
plus bas que $\overline{V}$ est un rev\^etement de~$\overline{U}$
mod\'er\'ement ramifi\'e sur~$\overline{X}$ relativement \`a
$\overline{S}$; par suite $V$ est un rev\^etement \'etale de~$U$
mod\'er\'ement ramifi\'e sur~$X$ relativement \`a~$S$.

b)
\marginpar{379}
Soit $F$ un faisceau en groupes sur~$U$ mod\'er\'ement
ramifi\'e sur~$X$ relativement \`a~$S$. Si $\overline{s}$ est un
point g\'eom\'etrique de~$S$, $y$ un point maximal de
$Y_{\overline{s}}$, on note $K$ le corps des fractions de
$\cal{O}_{X_{\overline{s}},y}$. Supposons que, pour tout point
$\overline{s}$ et pour tout point $y$, la $K$-alg\`ebre $L$ dont le
spectre repr\'esente $F|K$ soit de rang premier \`a la
caract\'eristique r\'esiduelle $p$ de
$\cal{O}_{X_{\overline{s}}}$. On dira parfois, par abus de langage,
que $F$ est premier aux caract\'eristiques r\'esiduelles de~$S$.
\index{premier aux caract\'eristiques r\'esiduelles (faisceau)|hyperpage}%
Lorsqu'il en est ainsi, tout torseur $P$ sous $F$ est
mod\'er\'ement ramifi\'e sur~$X$ relativement \`a~$S$.

Soient en effet $\overline{R}$ le localis\'e strict de
$\cal{O}_{X_{\overline{s}},y}$ en $\overline{y}$, $\overline{K}$ son
corps des fractions, $\overline{F}$ l'image inverse de~$F$ sur~$\overline{K}$. Montrons que l'on peut supposer $F$ constant. Comme
$F$ est mod\'er\'ement ramifi\'e sur~$X$ relativement \`a~$S$,
$\overline{F}$ est repr\'esentable par le spectre d'une
$\overline{K}$-alg\`ebre $L=\prod L_i$, o\`u les $L_i$ sont des
extensions de~$\overline{K}$ de degr\'e premier \`a $p$. On peut
donc trouver une extension $K'$ de~$\overline{K}$ de degr\'e premier
\`a $p$ telle que $\overline{F}|K'$ soit un faisceau
constant. D'apr\`es~\Ref{XIII.2.0.3}, pour prouver que
$P|\overline{K}$ est mod\'er\'ement ramifi\'e sur~$\overline{R}$, il suffit de voir que $P|K'$ est mod\'er\'ement
ramifi\'e sur la cl\^oture int\'egrale de~$\overline{R}$ dans
$K'$, d'o\`u la r\'eduction au cas o\`u $\overline{F}$ est
constant. Supposons d\'esormais $\overline{F}$ constant. La
$\overline{K}$-alg\`ebre $H$ qui repr\'esente $P|\overline{K}$ est
alors produit d'extensions $H_i$ de~$\overline{K}$ isomorphes entre
elles. Comme le rang de~$H$ est premier \`a $p$, il en est de
m\^eme de~$[H_1:\overline{K}]$, ce qui prouve que $H$ est
mod\'er\'ement ramifi\'e sur~$X$ relativement \`a~$S$.

c) Soient $X$ un sch\'ema r\'egulier, $D$ un diviseur \`a
croisements normaux de~$X$ (SGA~5 I~3.1.5), $U=X-\Supp D$, $F$ un
faisceau d'ensembles sur~$U$. Si $y$ est un point maximal de~$\Supp
D$, on d\'esigne par $K$ le corps des fractions de
$\cal{O}_{X,y}$. On dit que $F$ est \emph{mod\'er\'ement
ramifi\'e relativement \`a}~$D$,
\index{mod\'er\'ement ramifi\'e (faisceau)|hyperpage}%
si, pour tout point maximal $y$ de~$\Supp D$, $F|K$ est
repr\'esentable par une $K$-alg\`ebre mod\'er\'ement
ramifi\'ee sur~$\cal{O}_{X,y}$.
\end{remarques}

\begin{theoreme}
\label{XIII.2.4}
Soient
\marginpar{380}
$f\colon X\to S$ un $S$-sch\'ema, $D$ un diviseur sur~$X$
\`a croisements normaux relativement \`a~$S$ \eqref{XIII.2.1},
$Y=\Supp D$, $U=X-Y$, $i\colon U\to X$ l'immersion canonique. Soit $F$
un faisceau d'ensembles (\resp de groupes) sur~$U$, satisfaisant \`a
l'une des conditions suivantes:
\begin{enumerate}
\item[a)] $F$ est localement pour la topologie \'etale \emph{sur}
$X$ et sur~$S$ l'image inverse d'un faisceau d'ensembles (\resp d'un
faisceau en groupes constructible) sur~$S$.
\item[b)] $F$ est localement constant constructible sur~$U$ et
mod\'er\'ement ramifi\'e sur~$X$ relativement \`a~$S$.
\end{enumerate}
Alors on a les conclusions suivantes:
\begin{enumerate}
\item[1)]$(F,i)$ est cohomologiquement propre relativement \`a~$S$
en dimension $\leq 0$ (\resp pour tout morphisme $h\colon S'\to S$, si
$i'\colon U'\to X'$ est l'image inverse de~$i$ sur~$S'$, si $F'=F|U'$
et si $k=h_{(X)}$, le morphisme canonique
$$
\Psi\colon k^*(\R^1_{\tame}i_*F)\to\R^1_{\tame}i_*'F'
$$
est un isomorphisme).

Si $F$ est un faisceau en groupes premier aux caract\'eristiques
r\'esiduelles de~$S$ (\Ref{XIII.2.3}.b)) (mod\'er\'ement
ramifi\'e sur~$X$ relativement \`a~$S$), alors $(F,i)$ est
cohomologiquement propre relativement \`a~$S$ en dimension $\leq 1$.
\item[2)]Si $F$ est un faisceau d'ensembles (\resp de groupes)
constructible, $i_*F$ (\resp $\R^1_{\tame}i_*F$) est constructible.
\end{enumerate}
\end{theoreme}

\subsubsection*{D\'emonstration}
Pour tout $S$-sch\'ema $S'$, on consid\`ere le diagramme suivant
dont tous les carr\'es sont cart\'esiens:
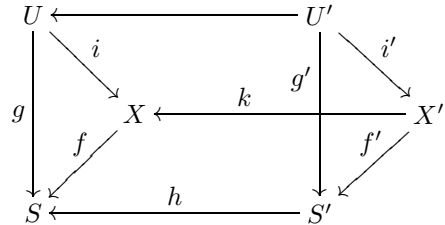
$$
\xymatrix{ U \ar[dd]_-g \ar[dr]^-i & & & U'\ar[lll] \ar[dd]_{\hbox{\raise10mm\hbox{$g'$}}}
\ar[dr]^-{i'}& \\ & X \ar[dl]_{f\mkern-8mu} & & & X'\ar[lll]_{k\hspace*{1cm}} \ar[dl]_{f'\mkern-10mu} \\ S
& & & S'\ar[lll]_-h & & }
$$
Comme
\marginpar{381}
la question est locale sur~$X$ pour la topologie \'etale, on peut
supposer que $D$ est un diviseur strictement \`a croisements normaux
relativement \`a~$S$ \eqref{XIII.2.1}; de plus, quitte \`a
restreindre $X$ \`a un voisinage de~$Y$, on peut supposer $X$ lisse
sur~$S$.

\Subsubsection*{D\'emonstration de~\Ref{XIII.2.4}~\textup{1)}}
\subsubsection{}
\label{XIII.2.4.1}
\emph{Cas d'un faisceau d'ensembles satisfaisant \`a} a). On peut
supposer que l'on~a $F=g^*G$, o\`u $G$ est un faisceau sur~$S$. Il
r\'esulte alors de (SGA~4 XVI~3.2) que le morphisme canonique
\begin{equation*}
\label{eq:XIII.2.4.1.1}
\tag{\thesubsubsection.1} f^*G\to i_*F
\end{equation*}
est un isomorphisme. Pour tout $S$-sch\'ema $h\colon S'\to S$, on a
de m\^eme un isomorphisme $f'^*G'\to i_*'F'$; par suite le morphisme
canonique
$$
\varphi\colon k^*(i_*F)\to i_*'F'
$$
s'identifie \`a l'isomorphisme naturel
$$
k^*f^*G\simeq f'^*G'
$$

\subsubsection{}
\label{XIII.2.4.2}
\emph{Cas d'un faisceau d'ensembles satisfaisant \`a} b). On doit
montrer que $\varphi$ est un isomorphisme et il suffit pour cela de
voir qu'il en est ainsi en chaque point g\'eom\'etrique
$\overline{x}'$ de~$X'$. Soit $\overline{S}$ (\resp $\overline{X}$,
\resp $\overline{S'}$, \resp $\overline{X'}$) le localis\'e strict
de~$S$ (\resp~$X$, \resp~$S'$, \resp $X'$) en $\overline{x}'$ et posons
$\overline{U}=U_{(\overline{X})}$,
$\overline{U'}=U_{(\overline{X'})}$, etc. Le morphisme
$\varphi_{\overline{x}}$ s'identifie au morphisme canonique
$$
\overline{\varphi}\colon \H^0(\overline{U},\overline{F})\to
\H^0(\overline{U'}, \overline{F'})
$$
On peut trouver un rev\^etement principal $V$ de~$\overline{U}$, du
type figurant dans~\Ref{XIII.5.4}, tel que l'image inverse de
$\overline{F}$ sur~$V$ soit un faisceau constant de valeur $C$. Si
$\Pi$ est le groupe de Galois de~$V$ sur~$\overline{U}$, $\Pi$
op\`ere sur~$\overline{F}|V$, et l'on a
\begin{equation*}
\label{eq:XIII.2.4.2.1}
\tag{\thesubsubsection.1}
\ifthenelse{\boolean{orig}}
{\H^0(U,F))\simeq \H^0(V,C_V)^{\Pi}}
{\H^0(\overline{U},\overline{F}))\simeq \H^0(V,C_V)^{\Pi}}
\end{equation*}
\label{indnot:mh}\oldindexnot{$\H^0(V,C_V)^{\Pi}$|hyperpage}%
o\`u le deuxi\`eme membre d\'esigne l'ensemble des
\'el\'ements de~$\H^0(V,C_V)$ invariants
\marginpar{382}
sous~$\Pi$. Comme $V'=V\times_{\overline{U}}\overline{U'}$ est une
rev\^etement principal de~$\overline{U'}$ de groupe de Galois
$\Pi'\simeq \Pi$, on voit que le morphisme $\overline{\varphi}$
s'obtient, en prenant les invariants sous $\Pi$, \`a partir du
morphisme canonique
$$
\H^0(V,C_V)\to \H^0(V',C_{V'}).
$$
Comme $V$ et $V'$ sont connexes \eqref{XIII.5.4}, ce morphisme, donc
aussi $\overline{\varphi}$, est un isomorphisme.

Notons que si de plus $F$ est un faisceau en groupes et si $P$ est un
torseur sur~$\overline{U}$ de groupe $\overline{F}$,
mod\'er\'ement ramifi\'e relativement \`a $D$, il r\'esulte
de la d\'emonstration pr\'ec\'edente et de~\Ref{XIII.2.2} que
$(\leftexp{P}{\overline{F}},\overline{i})$ est cohomologiquement
propre relativement \`a~$S$ en dimension $\leq 0$.

\subsubsection{}
\label{XIII.2.4.3}
\emph{Cas d'un faisceau en groupes.} Pour montrer que $\Psi$ est un
isomorphisme, il suffit de prouver que, pour tout point
g\'eom\'etrique $\overline{y}'$ de~$Y'$, le morphisme
$$
\Psi_{\overline{y}'}\colon
(k^*(\R^1_{\tame}i_*F))_{\overline{y}'}\to (\R^1_{\tame}
i_*'F')_{\overline{y}'}
$$
est un isomorphisme. Or, d'apr\`es~\Ref{XIII.2.1.2},
$\Psi_{\overline{y}'}$ s'identifie au morphisme canonique
$$
\overline{\Psi}\colon
\H^1_{\tame}(\overline{U},\overline{F})\to
\H^1_{\tame}(\overline{U'},\overline{F'}).
$$
Soient $\widetilde{U}$ le rev\^etement universel mod\'er\'ement
ramifi\'e de~$\overline{U}$ \eqref{XIII.2.1.3} et $\widetilde{F}$
l'image inverse de~$\overline{F}$ sur~$\widetilde{U}$. Il r\'esulte
de~\Ref{XIII.5.7} dans le cas a) et de~\Ref{XIII.5.5} dans le cas b),
que $\H^1_{\tame}(\overline{U},\overline{F})$ s'identifie au
sous-ensemble $\H^1(\overline{U},\overline{F})$ form\'e des classes
de~$\overline{F}$-torseurs dont l'image inverse sur~$\widetilde{U}$
est triviale. D'autre part un raisonnement
\ifthenelse{\boolean{orig}}
{classique}
{classique (\cf IX~\Ref{IX.5}, p\ptbl\pageref{page-300})}
montre que
l'ensemble des \'el\'ements de~$\H^1(\overline{U},\overline{F})$
dont l'image inverse sur~$\widetilde{U}$ est triviale s'identifie
\`a
$\H^1(\pi_1^{\tame}(\overline{U}),\H^0(\widetilde{U},\widetilde{F})).$
On obtient ainsi un isomorphisme canonique
\begin{equation*}
\label{eq:XIII.2.4.3.1} \tag{\thesubsubsection.1} {}
\H^1_{\tame}(\overline{U},\overline{F})\isomto
\H^1(\pi_1^{\tame}(\overline{U}),\H^0(\widetilde{U},\widetilde{F})).
\end{equation*}

Par suite le morphisme $\overline{\Psi}$ s'identifie au morphisme
canonique
$$
\H^1(\pi_1^{\tame}(\overline{U}),\H^0(\widetilde{U},\widetilde{F}))
\to
\H^1(\pi_1^{\tame}(\overline{U'}),\H^0(\widetilde{U}',\widetilde{F}')).
$$
Montrons
\marginpar{383}
que ce morphisme est un isomorphisme. Le morphisme
$\pi_1^{\tame}(\overline{U'})\to\pi_1^{\tame}(\overline{U})$
est un isomorphisme d'apr\`es~\Ref{XIII.5.6}, et il en est de
m\^eme du morphisme $\H^0(\widetilde{U},\widetilde{F})\to\allowbreak
\H^0(\widetilde{U}',\widetilde{F}')$. En effet, cela est \'evident
dans le cas b) car $\widetilde{F}$ est constant et $\widetilde{U}$ et
$\widetilde{U}'$ sont connexes. Dans le cas a), soit $\overline{G}$ un
faisceau en groupes constructible sur~$\overline{S}$ tel que l'on ait
$\overline{F}=\overline{g}^*\overline{G};$ comme les morphismes
$\widetilde{U}\to\overline{S}$ et
$\widetilde{U}'\to\overline{S'}$ sont $0$-acycliques
\eqref{XIII.5.7}, on a
$$
\H^0(\widetilde{U},\widetilde{F})\simeq
\H^0(\overline{S},\overline{G}) \simeq
\H^0(\overline{S'},\overline{G'})\simeq
\H^0(\widetilde{U}',\widetilde{F}'),
$$
ce qui entra\^ine que $\Psi$ est un isomorphisme. La derni\`ere
assertion de \Ref{XIII.2.4}~1) r\'esulte de ce qui pr\'ec\`ede,
compte tenu de \Ref{XIII.2.3}~b).

\subsubsection*{D\'emonstration de~\Ref{XIII.2.4}~\textup{2)}}

Le cas d'un faisceau d'ensembles constructible satisfaisant \`a a)
r\'esulte aussit\^ot de \eqref{eq:XIII.2.4.1.1}. Soit $F=g^*G$ un
faisceau en groupes satisfaisant \`a a), o\`u $G$ est un faisceau
constructible; on peut supposer $S$ affine; soient $(S_j)_{j\in J}$
une famille finie de sous-sch\'emas ferm\'es r\'eduits de~$S$
dont la r\'eunion recouvre $S$, tels que l'image inverse de~$G$ sur~$S_j$ soit un faisceau localement constant. Compte tenu de
\Ref{XIII.2.4}~1), il suffit, pour \'etablir que $\R^1_{\tame}i_*F$
est constructible, de voir qu'il en est ainsi apr\`es le changement
de base $S_j\to S$, pour chaque $j\in J$. On est donc
ramen\'e au cas~b) o\`u $F$ est localement constant.

On suppose d\'esormais que $F$ est un faisceau d'ensembles ou de
groupes satisfaisant \`a b). Comme la question est locale pour la
topologie \'etale sur~$X$, on peut supposer $X$ de pr\'esentation
finie sur~$S$, et, par passage \`a la limite, on peut supposer $X$
et $S$ noeth\'eriens.

Soit $D=\sum_{1\leq i\leq r}\divisor f_i$, o\`u, pour chaque point
$x$ de~$\Supp D$, si $I(x)$ est l'ensemble des $i$ tels que
$f_i(x)=0$, le sous-sch\'ema $V((f_i)_{i\in I(x)})$ est lisse sur~$S$ de codimension $\card I(x)$ dans~$X$. Soit $\mathcal{P}$
l'ensemble des parties de~$[ 1,r ]$ et, pour chaque $I\in\mathcal{P}$,
posons
$$
X_I=\textstyle\Big(\bigcap_{i\in I}V(f_i)\Big)\cap\Big(\bigcap_{i\not\in I}X_{f_i}\Big).
$$
Soit
\marginpar{384}
$z$ un point de~$X_I$. Quitte \`a se restreindre d'abord \`a un
voisinage \'etale de~$z$, on peut trouver un ouvert $W$ de~$X$
contenant $z$ et un rev\^etement principal $V$ de~$U\cap W$,
mod\'er\'ement ramifi\'e sur~$W$ relativement \`a~$S$, du type
consid\'er\'e dans~\Ref{XIII.5.6.1}, tel que l'image
r\'eciproque de~$F$ sur~$V$ soit un faisceau constant de
valeur~$C$. Soit $\pi$ le groupe de Galois du rev\^etement~$V$. Pour
tout point g\'eom\'etrique $\overline{x}$ de~$X_I$, on a alors
d'apr\`es \eqref{eq:XIII.2.4.2.1}
$$
\H^0(\overline{U},\overline{F})\simeq \H^0(V,C_V)^{\pi}.
$$
Il en r\'esulte que $i_*F|X_I\cap W$ est localement constant (SGA~4
IX~2.13) et par suite $i_*F$ est constructible.

Montrons enfin que si $F$ est un faisceau en groupes localement
constant, $\R^1_{\tame}i_*F|_{X_I}$ est constructible. Si
$\overline{x}$ est un point g\'eom\'etrique de~$X$, on a obtenu
dans \eqref{eq:XIII.2.4.3.1} l'expression
$$
\H^1_{\tame}(\overline{U},\overline{F})\isomto
\H^1(\pi^{\tame}_1(\overline{U}),\H^0(\widetilde{U},\widetilde{F})).
$$
Si $p$ est la caract\'eristique r\'esiduelle de~$\overline{X}$, on
a, d'apr\`es~\Ref{XIII.5.6},
$\pi_1^{\tame}(\overline{U})=\prod_{l\neq p }\ZZ_l(1)^{\card I}$. D\'esignons par $\LL$ l'ensemble des nombres premiers qui
divisent l'ordre du groupe fini $\H^0(\widetilde{U},\widetilde{F})$ et
soit $K=\prod_{l\in\LL -\{p\}\cap\LL }\ZZ_l(1)^{\card I}$. Il
r\'esulte de \cite[I \S5 ex.2]{XIII.4} que l'on a
$$
\H^1_{\tame}(\overline{U},\overline{F})\simeq
\H^1(K,\H^0(\widetilde{U},\widetilde{F})).
$$
Comme $K$ est topologiquement de type fini et
$\H^0(\widetilde{U},\widetilde{F})$ fini, on en d\'eduit tout
d'abord que les fibres du faisceau $\R^1_{\tame}i_*F|_{X_I}$ sont
finies. D'autre part, l'ensemble $\LL$ ne d\'epend pas du point
$\overline{x}$. Pour tout $q\in\LL$, soit $X_{I,q}$ le ferm\'e de
$X_I$ d'\'equation $q=0$ et soit $X_{I'}$ l'ouvert de~$X_I$
compl\'ementaire de la r\'eunion des $X_{I,q}$. Alors
$\R^1_{\tame}i_*F|_{X_{I,q}}$ et $\R^1_{\tame}i_*F|_{X_I'}$ sont
localement constants; en effet une fl\`eche de sp\'ecialisation de
points g\'eom\'etriques de~$X_{I,q}$ (\resp de~$X_{I'}$) induit un
isomorphisme sur les groupes $K$ \eqref{XIII.5.6.1} donc aussi sur les
ensembles $\H^1_{\tame}(\overline{U},\overline{F})$, et l'on peut appliquer
SGA~4 IX~2.13.

\begin{corollaire}
\label{XIII.2.5}
Soient
\marginpar{385}
$f\colon X\to S$ un morphisme, $D$ un diviseur sur~$X$
\`a croisements normaux relativement \`a~$S$ \eqref{XIII.2.1},
$Y=\Supp D$, $U=X- Y$, $i\colon U\to X$ l'immersion
canonique. Soit $\Phi$ un champ sur~$U$ et supposons donn\'es des
morphismes \'etales surjectifs $S_1\to S$ et $X_2\to
X\times_{S}S_1$, un champ $\Psi$ sur~$S_1$ et un isomorphisme
$\Phi|U_2\simeq\Psi|U_2$ (\cf \Ref{XIII.2.2}).

Alors, pour tout morphisme $h\colon S'\to S$, si $k=h_{(X)}$,
le foncteur canonique
$$
\varphi\colon k^*i_*\Phi\to i_*'\Phi'
$$
est pleinement fid\`ele. Si $\Psi$ est constructible, le foncteur
canonique
$$
\psi\colon k^*i_*^{\tame}\Phi\to i_*^{\prime\tame}\Phi'
$$
est une \'equivalence de cat\'egories.

De plus, si le champ $\Psi$ est constructible (\resp si $\Psi$ est
$1$-constructible \eqref{XIII.0} et $S$ localement noeth\'erien),
$i_*\Phi$ est constructible (\resp $i_*^t\Phi$ est
$1$-constructible).
\end{corollaire}

Montrons que $\varphi$ est pleinement fid\`ele. Il suffit de voir
que, pour tout point g\'eom\'etrique $\overline{x}'$ de~$X'$, il
en est ainsi du foncteur
$$
\overline{\varphi}\colon
\Phi(\overline{U})\to\Phi'(\overline{U'})
$$
(on a repris les notations de~\Ref{XIII.2.4.2}). Soient $a$, $b$ deux
\'el\'ements de~$\Phi(\overline{U})$, $a',\ b'$ leurs images par
$\overline{\varphi}$. Comme le morphisme
$\overline{U}\to\overline{S}$ est localement $0$-acyclique
\eqref{XIII.5.7}, on a un isomorphisme
$$
\H^0(\overline{U},\overline{S\Phi})
\overset{\sim}{\to}\H^0(S,\overline{S\Psi}).
$$
Par suite $a$ et $b$ proviennent par image inverse d'\'el\'ements
de~$\Psi$, et il en est donc de m\^eme de
$F=\SheafHom_{\overline{U}}(a,b)$. Comme le faisceau
$F'=\SheafHom_{\overline{U'}}(a',b')$ est l'image inverse sur~$\overline{U'}$ de~$F$, il r\'esulte de~\Ref{XIII.2.4.1} que le
morphisme canonique
$$
\H^0(\overline{U},F)\to \H^0(\overline{U'},F')
$$
est un isomorphisme, ce qui prouve que $\overline{\varphi}$ est
pleinement fid\`ele.

Montrons
\marginpar{386}
que, pour tout point g\'eom\'etrique $\overline{x}'$ de~$X'$, le
foncteur
$$
\overline{\psi}\colon i^{\tame}_*\Phi(\overline{X'})\to
i_*^{\prime\tame}\Phi'(\overline{X'})
$$
est une \'equivalence. D'apr\`es ce qui pr\'ec\`ede,
$\overline{\psi}$ est pleinement fid\`ele. Montrons que
$\overline{\psi}$ est essentiellement surjectif. Soit $a'$ un
\'el\'ement de~$\Phi'(\overline{U'})$ mod\'er\'ement
ramifi\'e sur~$X'$ relativement \`a~$S'$ \eqref{XIII.2.2} et
montrons qu'il est l'image d'un \'el\'ement mod\'er\'ement
ramifi\'e de~$\Phi(\overline{U})$. Il r\'esulte de~\Ref{XIII.2.4}
1) que le morphisme canonique
$$
\H^0(\overline{U},\overline{S}\overline{\Phi})\to
\H^0(\overline{U'},\overline{S}\overline{\Phi}')
$$
est un isomorphisme. Soit $G'$ la sous-gerbe maximale de~$\Phi'$
engendr\'ee par $a'$; il existe alors une sous-gerbe maximale $G$ de
$\overline{\Phi}$, image inverse d'une gerbe sur~$\overline{S}$, telle
que l'on ait
$$
\overline{m}^*G\simeq G',
$$
o\`u $m$ est le morphisme $\overline{U'}\to
\overline{U}$. Le foncteur canonique
\begin{equation*}
\label{eq:XIII.2.5.*}
\tag{$*$} {\overline{k}^*\overline{i}_*^{\tame}G\to
\overline{i}_*^{\prime\tame}G'}
\end{equation*}
est une \'equivalence, car $G$ s'identifie \`a une gerbe de
torseurs sous un faisceau en groupes constructible provenant de
$\overline{S}$, et l'on peut appliquer \Ref{XIII.2.4}~1). Il
r\'esulte alors de \eqref{eq:XIII.2.5.*} qu'il existe un
\'el\'ement $a$ de~$G(\overline{U})$, mod\'er\'ement
ramifi\'e sur~$X$ relativement \`a~$S$, dont l'image inverse sur~$\overline{U'}$ est $a'$, ce qui prouve que $\psi$ est une
\'equivalence.

Si $\Psi$ est constructible, il en est de m\^eme de~$i_*\Phi$; en
effet un objet $x$ de~$i_*\Phi$ est, localement pour la topologie
\'etale de~$S$, image inverse d'un objet $y$ de~$\Psi$; il
r\'esulte donc de~\Ref{eq:XIII.2.4.1.1} que $\SheafAut (x)$ est
l'image inverse de~$\SheafAut (y)$, donc est constructible. Enfin, si
$\Psi$ est $1$-constructible, il en est de m\^eme de
$i^{\tame}_*\Phi$ d'apr\`es~\Ref{XIII.6.3} ci-dessous.

\begin{corollaire}
\label{XIII.2.6}
Les
\marginpar{387}
notations sont celles de~\Ref{XIII.2.4}. Supposons que $S$ soit de
caract\'eristique nulle en tout point $s$ tel que l'on ait $Y_s\neq
\emptyset.$ Alors, si $\cal{ F}$ est un faisceau en groupes localement
constant constructible sur~$U$ (\resp un champ constructible sur~$U$ qui
est localement sur~$X$ et $S$ image inverse d'un champ constructible
sur~$S$), $(\cal{ F},i)$ est cohomologiquement propre relativement
\`a~$S$ en dimension $\leq 1.$
\end{corollaire}

Comme tout faisceau d'ensembles constructible sur~$U$ est
mod\'er\'ement ramifi\'e sur~$X$ relativement \`a~$S,$ le
corollaire r\'esulte de~\Ref{XIII.2.4} (\resp \Ref{XIII.2.5}).

\begin{corollaire}
\label{XIII.2.7}
Les notations sont celles de~\Ref{XIII.2.4}, mais on se donne de plus
un $S$-sch\'ema $T$, un morphisme propre $p\colon X\to T,$
et on suppose $X$ et $T$ de pr\'esentation finie sur~$S$; soit
$q=pi$. Soit $\cal{ F}$ un faisceau d'ensembles constructible sur~$U$
satisfaisant \`a l'une des conditions \textup{a)} ou \textup{b)} de~\Ref{XIII.2.4}
(\resp un faisceau en groupes satisfaisant \`a l'une des condition
\textup{a)}, \textup{b)} de~\Ref{XIII.2.4}, \resp un champ sur~$U$ qui est localement sur~$X$ et $S$ image inverse d'un champ constructible $G$ sur~$S$). Alors
on a les conclusions suivantes:
\begin{enumerate}
\item[1)] $(\cal{ F},q)$ est cohomologiquement propre relativement \`a~$S$ en dimension $\leq 0$ (\resp pour tout morphisme $h\colon
S^{'}\to S$, si $m=h_{(T)}$, le morphisme canonique:
$$
\Theta\colon m^*(\R^1_{\tame}q_*\cal{ F})\to
\R^1_{\tame}q_*^{'}\cal{ F}^{'}
$$
est un isomorphisme, \resp pour tout morphisme $S^{'}\to S$,
le morphisme canonique
$$
\xi \colon m^*(q_*^{\tame}\cal{ F})\to
{q'}_*^{{\tame}}\cal{F}^{'}
$$
est une \'equivalence).
\item[2)] le faisceau $q_*\cal{ F}$ (\resp le faisceau
$\R_{\tame}^1q_*\cal{ F}$, \resp le champ $q_*^{\tame}\cal{ F}$)
est constructible. Dans le dernier cas, si l'on suppose $S$ localement
noeth\'erien et $G$ $1$-constructible, il en est de m\^eme de
$q_*^{\tame}\cal{ F}.$
\end{enumerate}
\end{corollaire}

La
\marginpar{388}
premi\`ere partie r\'esulte aussit\^ot de~\Ref{XIII.2.4},
\Ref{XIII.2.5} et de la d\'emonstration
de~\Ref{XIII.1.8}. D\'emontrons~2). Si $\cal{ F}$ est un faisceau
d'ensembles constructible sur~$U$ satisfaisant \`a \Ref{XIII.2.4}~a)
ou \Ref{XIII.2.4}~b), il r\'esulte de \Ref{XIII.2.4}~2) que
$i_*\cal{ F}$ est constructible; il en est donc de m\^eme de
$q_*\cal{ F}=p_*(i_*\cal{ F})$ (SGA~4 XIV~1.1).

Soit $\cal{F}$ un faisceau en groupes constructible sur~$U$
satisfaisant \`a \Ref{XIII.2.4}~a) ou \Ref{XIII.2.4}~b) et prouvons
que $\R_{\tame}^1q_*\cal{ F}$ est constructible. Par passage \`a
la limite (EGA~IV 8.10.5 et 17.7.8) et en utilisant 1), on peut
supposer $S$ noeth\'erien. Soit alors $\Phi$ le champ sur~$X$ dont
la fibre en tout sch\'ema $X^{'}$ \'etale sur~$X$ est form\'ee
des torseurs sur~$U^{'}=U\times_{X}X^{'}$, de groupe $\cal{F}|U$, qui
sont mod\'er\'ement ramifi\'es sur~$X$ relativement \`a~$S$; on
a donc
$$
S(i_*^{\tame}\Phi)\simeq \R_{\tame}^1i_*\cal{ F}
$$
et ce faisceau est constructible d'apr\`es \Ref{XIII.2.4}~2).
Il r\'esulte donc de~\Ref{XIII.6.3} ci-dessous que $S(p_*\Phi)$
est constructible, \ie que $\R^1_{\tame}q_*\cal{ F}$ est
constructible.

Enfin, si $\cal{ F}$ est un champ sur~$U$ qui est localement sur~$X$ et
$S$ image inverse d'un champ constructible sur~$S,$
$i_*^{\tame}\cal{ F}$ est constructible, et il en est donc de
m\^eme de~$q_*^{\tame}\cal{ F}=p_*(i_*^{\tame}\cal{ F})$. Si
de plus $S$ est localement noeth\'erien et $SG$ constructible,
$S(i_*^{\tame}\cal{ F})$ est constructible d'apr\`es 6.3; il en
est donc de m\^eme de~$S(q_*i_*^{\tame}\cal{ F})$, \ie de
$S(q_*^{\tame}\cal{ F})$~\Ref{XIII.6.2}.

\begin{corollaire}
\label{XIII.2.8}
Soit
$$
\xymatrix{ U \ar[r]^{i} \ar[d]_{g} & X \ar[ld]^{f} \\ S }
$$
un diagramme commutatif de sch\'emas, dans lequel $U$ est l'ouvert
compl\'ementaire
\marginpar{389}
dans $X$ d'un diviseur \`a croisements normaux relativement \`a~$S$, $f$ un morphisme \emph {propre} de pr\'esentation finie. Soit
$\LL$ un ensemble de nombres premiers. On suppose~$g$ localement
$0$-acyclique (\resp localement $1$-asph\'erique pour $\LL$). Alors, si
$\cal{ F}$ est un faisceau d'ensembles sur~$U$ (\resp un faisceau de
$\LL$-groupes) sur~$U,$ localement constant constructible,
mod\'er\'ement ramifi\'e sur~$X$ relativement \`a~$S$,
$f_*\cal{ F}$ (\resp $\R^1_{\tame}f_*\cal{ F}$) est localement
constant constructible et $(\cal{ F},f)$ est cohomologiquement propre
relativement \`a~$S$ en dimension $\leq 0$ (\resp la formation de
$\R^1_{\tame}f_*\cal{ F}$ commute \`a tout changement de base
$S^{'}\to S$). Dans le cas non resp\'e, si $\cal{ F}$ est un
faisceau en groupes, pour toute sp\'ecialisation $\overline{s}_1
\to \overline{s}_2$ de points g\'eom\'etriques de~$S$, le
morphisme de sp\'ecialisation
$$
(\R^1_{\tame}f_*\cal{ F})_{\overline{s}_2}\to
(\R^1_{\tame}f_*\cal{ F})_{\overline{s}_1}
$$
est injectif.
\end{corollaire}
Le corollaire r\'esulte aussit\^ot de~\Ref{XIII.2.6} et
de~\Ref{XIII.1.14} (\resp de l'analogue de~\Ref{XIII.1.14} pour le
$\R^1_{\tame}f_*\cal{ F}$, lequel se d\'emontre comme \loccit).

\begin{corollaire} \label{XIII.2.9}
Soit
$$
\xymatrix{ U \ar[r]^{i} \ar[d]_{g} & X \ar[ld]^{f} \\ S }
$$
un diagramme commutatif de sch\'emas, dans lequel $U$ est l'ouvert
compl\'ementaire dans~$X$ d'un diviseur \`a croisements normaux
relativement \`a~$S$, $f$~un morphisme propre lisse de
pr\'esentation finie. Soit~$\LL $ l'ensemble des nombres premiers
distincts des caract\'eristiques r\'esiduelles de~$S$. Soit~$\cal{F}$ un faisceau de~$\LL $-groupes localement constant constructible
sur~$U$, mod\'er\'ement ramifi\'e sur~$X$ relativement
\`a~$S$. Alors $\R^1 f_*\cal{ F}$ est localement constant
constructible et $(\cal{ F},f)$ est cohomologiquement propre
relativement
\marginpar{390}
\`a~$S$ en dimension $\leq 1$.
\end{corollaire}
Le corollaire r\'esulte de~\Ref{XIII.2.8} et du fait que l'on a
$\R_{\tame}^1f_*\cal{ F}=\R^1f_*\cal{ F}$ (\Ref{XIII.2.3}~b)).

\subsection{}
\label{XIII.2.10}
Si $U$ est un sch\'ema connexe, $a$ un point g\'eom\'etrique
de~$U$, $\LL$~un ensemble de nombre premiers, on note
\begin{equation*}
\label{eq:XIII.2.10.0}
\tag{2.10.0} \pi_1^{\LL}(U,a)
\end{equation*}
\label{indnot:mi}\oldindexnot{$\pi_1^{\LL}(U,a)$|hyperpage}%
la limite projective des quotients finis de~$\pi_1(U,a)$ dont les
ordres ont tous leurs facteurs premiers dans~$\LL$.

Nous allons d\'efinir des morphismes de sp\'ecialisation pour le
groupe fondamental, g\'e\-n\'e\-ra\-li\-sant X.\Ref{X.2}.

Soit $g\colon U\to S$ un morphisme coh\'erent \`a fibres
g\'eom\'etriquement connexes (\resp un morphisme de la forme $g=
fi$, o\`u $f\colon X \to S$ est un morphisme propre de
pr\'esentation finie et o\`u $i\colon U\to X$ est une
immersion ouverte telle que $U$ soit le compl\'ementaire dans $X$
d'un diviseur \`a croisements normaux relativement \`a~$S$
(\cf \Ref{XIII.2.8})). Soit ${\LL}$ un ensemble de nombres premiers et
supposons, dans le cas non resp\'e, que, pour tout ${\LL}$-groupe
constant fini $C$, $(C_U,g)$ soit cohomologiquement propre
relativement \`a~$S$ en dimension $\leq 1$. Soient $\overline{s}_1
\to \overline{s}_2$ un morphisme de sp\'ecialisation de
points g\'eom\'etriques de~$S$, $\overline{S}$ le localis\'e
strict de~$S$ en $\overline{s}_2$,
$\overline{U}=U\times_S\overline{S}$. On a un diagramme commutatif
$$
\xymatrix{ U_{\bar{s}_1} \ar[r]^{h_1} \ar[d] & \bar{U}
\ar[d]^{\bar{g}} & U_{\bar{s}_1} \ar[l]_{h_2} \ar[d] \\ \bar{s}_1
\ar[r] & \bar{S} & \,\bar{s}_2\,. \ar[l] }
$$
Si $a_1$ est un point g\'eom\'etrique de~$U_{\bar{s}_1},\ a_2$ un
point g\'eom\'etrique de~$U_{\bar{s}_2}$ les morphismes $h_1$ et
$h_2$ d\'efinissent des morphismes canoniques
\marginpar{391}%
\begin{align*}
\pi_1\colon \pi_1^{\LL}(U_{\bar{s}_1},a_1)&\to
\pi_1^{\LL}(\bar{U},a_1) & \pi_2\colon
\pi_1^{\LL}(U_{\bar{s}_2},a_2)&\to \pi_1^{\LL}(\bar{U},a_2)\\
(\text{\resp~} \pi_1\colon \pi_1^\tame (U_{\bar{s}_1},a_1)&\to
\pi_1^\tame (\bar{U},a_1) & \pi_2\colon
\pi_1^\tame(U_{\bar{s}_2},a_2)&\to \pi_1^\tame(\bar{U},a_2))
\end{align*}
(V~\Ref{V.7} et~\Ref{eq:XIII.2.1.5.2}). Les hypoth\`eses de
propret\'e cohomologique (\resp \Ref{XIII.2.8}) prouvent que $\pi_2$
est un isomorphisme. Si l'on choisit une classe de chemins de~$a_1$
\`a $a_2$, on obtient un isomorphisme
$$
\pi_{12}\colon \pi_1^{\LL}(\bar{U},a_1)\isomto
\pi_1^{\LL}(\bar{U},a_2) \quad (\text{\resp~} \pi_{12}\colon
\pi_1^\tame (\bar{U},a_1)\isomto \pi_1^\tame(\bar{U},a_2))
$$
d'o\`u un morphisme $\pi=\pi_2^{-1}\pi_{12}\pi_1$
$$
\pi\colon \pi_1^{\LL}(U_{\bar{s}_1},a_1)\to
\pi_1^{\LL}(U_{\bar{s}_2},a_2) \quad (\text{\resp~}\pi\colon
\pi_1^\tame (U_{\bar{s}_1},a_1)\to \pi_1^\tame
(U_{\bar{s}_2},a_2)).
$$
Changer la classe de chemins de~$a_1$ \`a $a_2$ revient \`a
modifier $\pi$ par un automorphisme int\'erieur de
$\pi_1^{\LL}(X_{\bar{s}_2},a_2)$ (\resp de
$\pi_1^\tame(X_{\bar{s}_2},a_2)$). On appelle \emph{morphisme de
sp\'ecialisation pour le groupe fondamental}
\index{homomorphisme de sp\'ecialisation pour le groupe fondamental|hyperpage}%
\index{sp\'ecialisation pour le groupe fondamental (homomorphisme de)|hyperpage}%
associ\'e au morphisme $\bar{s}_1\to \bar{s}_2$ et on note
simplement
$$
\pi\colon \pi_1^{\LL}(X_{\bar{s}_1})\to
\pi_1^{\LL}(X_{\bar{s}_2}) \quad \text{(\resp~} \pi\colon
\pi_1^{\tame}(X_{\bar{s}_1})\to \pi_1^{\tame}(X_{\bar{s}_2}))
$$
l'un des morphismes d\'efini ci-dessus.

\begin{lemme}
\label{XIII.2.11}
Soient $f\colon X\to S$ un morphisme propre de
pr\'esentation finie, $D$ un diviseur sur~$X$ \`a croisements
normaux relativement \`a~$S$, $Y=\Supp D$, $U=X-Y$, $i\colon
U\to X$ le morphisme canonique, $\bar{s}_1\to
\bar{s}_2$ un morphisme de sp\'ecialisation de points
g\'eom\'etriques de~$S$, $y_1$ un point g\'eom\'etrique de
$Y_{\bar{s}_1}$, $y_2$ un point g\'eom\'etrique de
$Y_{\bar{s}_2}$, tel que la projection $z_1$ de~$y_1$ sur~$X$ soit une
g\'en\'erisation de la projection $z_2$ de~$y_2$. Soit $I_{y_1}$
un sous-groupe d'inertie de~$\pi_1^{\tame}(\bar{U}_{\bar{s}_1})$ en
$y_1.$ Alors l'image de~$I_{y_{1}}$ par le morphisme de
sp\'ecialisation
$$
\pi\colon \pi_1^{\tame}(U_{\bar{s}_1})\to
\pi_1^{\tame}(U_{\bar{s}_2})
$$
est
\marginpar{392}
un sous-groupe d'inertie de~$\pi_1^{\tame}(U_{\overline{s}_2})$
en~$y_2$.
\end{lemme}

Soient en effet $\overline{X}$ (\resp $\widetilde{X}$) le localis\'e
strict de~$X$ en $y_2$ (\resp en $y_1$), $\overline{U}=U \times_X
\overline{X}$ (\resp $\widetilde{U}=U\times_X \widetilde{X}$). On a un
morphisme canonique $\widetilde{U}\to \overline{U}$, et il
r\'esulte de~\Ref{XIII.1.10} que l'on a un diagramme commutatif
$$
\xymatrix{ \pi_1^{\tame} \left(\widetilde{U}_{\overline{s}_1} \right)
\ar[r]^{\pi '} \ar[d] & \pi_1^{\tame} \left(
\overline{U}_{\overline{s}_2} \right) \ar[d] \\ \pi_1^{\tame} \left(
U_{\overline{s}_1} \right) \ar[r]^{\pi } & \pi_1^{\tame} \left(
U_{\overline{s}_2} \right) }
$$
o\`u $\pi'$ est compos\'e du morphisme canonique
$\pi_1^{\tame}(\widetilde{U}_ {\overline{s}_1})\to
\pi_1^{\tame}(\overline{U}_{\overline{s}_1})$ et du morphisme de
sp\'ecialisation. Comme
$\pi_1^{\tame}(\widetilde{U}_{\overline{s}_1})$
(\resp $\pi_1^{\tame}(\overline{U}_{\overline{s}_1})$)
est un groupe d'inertie de
$\pi_1^{\tame}(U_{\overline{s}_1})$ en~$y_1$ et (\resp de~$\pi_1^{\tame}(U_{\overline{s}_2})$ en~$y_2$), il suffit de
prouver que $\pi'$ est surjectif. Mais cela r\'esulte de
l'expression obtenue dans~\Ref{XIII.5.6}.

\begin{corollaire}
\label{XIII.2.12}
Soit $X$ une courbe propre et lisse connexe de genre $g$ sur un corps
s\'e\-pa\-ra\-blement clos $k$ de caract\'eristique $p\geq
0$. Soit $U$ l'ouvert obtenu en enlevant \`a $X$ $n$ points
ferm\'es distincts $a_1,\dots,a_n$. Alors le groupe fondamental
mod\'er\'ement ramifi\'e $\pi_1^{\tame}(U)$ \eqref{XIII.2.1.3}
peut \^etre engendr\'e par $2g+n$ \'el\'ements $x_i,
y_i,\sigma_j$, avec $1\leq i\leq g$, $1\leq j\leq n$, tel que
$\sigma_j$ soit un g\'en\'erateur d'un groupe d'inertie
correspondant \`a $a_j$, et que l'on ait la relation
\begin{equation*}
\label{eq:XIII.2.12.*}
\tag{$*$} \prod_{1\leq i\leq
g}(x_iy_ix^{-1}_iy^{-1}_i)\cdot\prod_{1\leq j\leq n}\sigma_j=1.
\end{equation*}
Pour tout groupe fini $G$ \emph{d'ordre premier \`a} $p$,
engendr\'e par des \'el\'ements $\overline{x}_i,
\overline{y}_i,\overline{\sigma}_j$ satisfaisant \`a la relation
\eqref{eq:XIII.2.12.*}, il existe un rev\^etement \'etale de~$U$,
de groupe $G$, correspondant \`a un homomorphisme
$\pi_1^{\tame}(U)\to G$ qui
\marginpar{393}
envoie $x_i, y_i,\sigma_j$ sur~$\overline{x}_i,
\overline{y}_i,\overline{\sigma}_j$ respectivement. En d'autres
termes, si $p'$ d\'esigne l'ensemble des nombres premiers distincts
de~$p$, $\pi_1^{p'}(U)$ est le pro-$p'$-groupe engendr\'e par les
g\'en\'erateurs $x_i, y_i,\sigma_j$ li\'es par la seule relation~\eqref{eq:XIII.2.12.*}.
\end{corollaire}

\subsubsection*{D\'emonstration}
On peut supposer $k$ alg\'ebriquement clos. Supposons d'abord $k$ de
carac\-t\'eristique z\'ero. Il existe alors une sous-extension
alg\'ebriquement close $k'$ de~$k$, de degr\'e de transcendance
fini sur~$\QQ$, telle que $X$ provienne par extension des scalaires
d'une courbe propre et lisse $X'$ d\'efinie sur~$k'$, et l'on peut
supposer que les points $a_1,\dots,a_n$ proviennent de points rationnels
$a'_1,\dots,a'_n$ de~$X'$. Comme $k'$ est de degr\'e de transcendance
fini sur~$\QQ$, on peut trouver un plongement de~$k'$ dans le corps
des nombres complexes $\CC$; soit $\widetilde{U}=U'\times_{k'}\CC.$
Soit $k''$ une extension alg\'ebriquement close de~$k'$ telle que
l'on ait des $k'$-morphismes de~$k$ et de~$\CC$ dans $k''$. Si
$g'\colon U'\to k'$ est le morphisme structural, et si
$\cal{F}$ est un faisceau en groupes constant fini sur~$U''$, il
r\'esulte de~\Ref{XIII.2.9} que les morphismes de sp\'ecialisation
$$
(\R^1g'_*\cal{F}')_{\CC}\to(\R^1g'_*\cal{F}')_{k''}
\from(\R^1g'_*\cal{F}')_k
$$
sont des isomorphismes. En termes de groupes fondamentaux, cela montre
que l'on a un isomorphisme, d\'efini \`a automorphisme
int\'erieur pr\`es
$$
\pi_1(U)\to\pi_1(\widetilde{U}),
$$
et il est clair que cet isomorphisme transforme un groupe d'inertie
relatif \`a un point de~$X'-U'$ en un groupe d'inertie relatif au
m\^eme point. On peut donc supposer que l'on a $k=\CC$. Dans ce
dernier cas il r\'esulte du th\'eor\`eme d'existence de Riemann
(XII~\Ref{XII.5.2}) que le groupe fondamental $\pi_1(U)$ n'est autre
que le compl\'et\'e pour la topologie des sous-groupes d'indice
fini du groupe fondamental de l'espace analytique associ\'e \`a~$U$. Or ce dernier peut se calculer par voie transcendante \cite[ch\ptbl 7
\S 47]{XIII.3};
\marginpar{394}
il peut \^etre engendr\'e par $2g+n$ \'el\'ements $x_i,
y_i,\sigma_j$ tels que $\sigma_j$ soit l'image d'un g\'en\'erateur
du groupe fondamental local $\pi_1(D_j)$ d'un petit disque centr\'e
en $a_j$, \ie un g\'en\'erateur d'un groupe d'inertie
correspondant au point $a_j$, ces \'el\'ements satisfaisant \`a
la seule relation \eqref{eq:XIII.2.12.*}.

Si maintenant $k$ est de
caract\'eristique $p>0$, on peut trouver un anneau de valuation
discr\`ete complet $A$, de corps r\'esiduel $k$, de corps des
fractions $K$ de caract\'eristique z\'ero, et un sch\'ema
connexe $X_1$, propre et lisse sur~$S=\Spec A$, tel que l'on ait
$X_1\times_S \Spec k\simeq X$ (III~\Ref{III.7.4}). Les points $a_j$ se
rel\`event alors en des sections $s_j$ de~$X_1$ au-dessus de~$S$;
soit $Y_{1j}$ le sous-sch\'ema ferm\'e r\'eduit d'espace
sous-jacent $s_j(S)$, $Y_1$ la r\'eunion des $Y_{1j}$,
$U_1=X_1-Y_1$, $g_1\colon U_1\to S$ le morphisme
structural. Soient $\overline{K}$ une extension alg\'ebriquement
close de~$K$, $\overline{U}= U_1\times_S\overline{K}$. Si $C$ est un
groupe constant fini, il r\'esulte de~\Ref{XIII.2.8} que le
morphisme de sp\'ecialisation
$$
(\R^1g_{1*}C_{U_1})_k\to(\R^1g_{1*}C_{U_1})_{\overline{K}}
$$
est injectif et m\^eme bijectif si $C$ est d'ordre premier \`a
$p$. Or cela signifie, en termes de groupes fondamentaux, que le
morphisme de sp\'ecialisation \eqref{XIII.1.10}
$$
\pi \colon \pi_1(\overline{U})\to\pi^{\tame}_1(U)
$$
est surjectif, et que le morphisme de sp\'ecialisation
$$
\pi^{p'}_1(\overline{U})\to\pi^{p'}_1(U)
$$
est bijectif. Enfin, si $x_i, y_i,\sigma_j$ sont des
g\'en\'erateurs de~$\pi_1(\overline{U})$ tels que $\sigma_j$ soit
un g\'en\'erateur d'un groupe d'inertie correspondant au point
$b_j=Y_{1j}(\overline{K})$ de~$\overline{X}$, alors,
d'apr\`es~\Ref{XIII.1.11}, $\pi(\sigma_j)$ est un g\'en\'erateur
d'un groupe d'inertie correspondant \`a $a_j$, ce qui ach\`eve la
d\'emonstration.

\ifthenelse{\boolean{orig}}
{}
{\begin{remarqueMR}
\label{XIII.2.13}
Soient $k$ un corps alg\'ebriquement clos de caract\'eristique $p>0$, $X$~une courbe alg\'ebrique propre lisse sur~$k$, connexe, de genre $g$, $U$~un ouvert affine de~$X$, compl\'ementaire de~$r\geq1$ points rationnels
de~$X$. On dispose du groupe fondamental $\pi_1(U)$, de son quotient
$\pi_1^{\tame}(U)$ qui classifie les rev\^etements finis \'etales de~$U$,
mod\'er\'ement ramifi\'es aux points de~$X-U$, et du quotient
$\pi_1^{p'}(U)$ de~$\pi_1^{\tame}(U)$ qui classifie les rev\^etements \'etales
galoisiens de~$U$, de groupe de Galois d'ordre premier \`a~$p$.

Dans (Coverings of algebraic curves, Amer.\ J.~Math.\ \textbf{79}
(1957), p\ptbl825--856) S\ptbl Abhyankar formulait un certain nombre de conjectures sur la structure des groupes finis, groupes de Galois de rev\^etements finis
\'etales de~$U$.

Les conjectures concernant les groupes finis qui sont groupes de
Galois de rev\^etements \'etales connexes de~$U$ d'ordre premier \`a~$p$,
sont d\'emontr\'ees et pr\'ecis\'ees dans le corollaire \Ref{XIII.2.12}. Avant d'aborder les quotients finis de~$\pi_1(U)$, commen\c{c}ons par donner quelques
indications sur la \og taille\fg de ce groupe. Soit $X = \Spec(A)$. On
sait que $\Hom_{\mathrm{cont}}(\pi_1(U), \ZZ/p\ZZ)$ est d\'ecrit par
la th\'eorie d'Artin-Schreier (cor\ptbl XI~\Ref{XI.6.9}). Ce groupe est
isomorphe \`a $A/\wp(A)$ o\`u $\wp$ est l'application de~$A$ dans $A$, qui
envoie $a$ sur~$a^p-a$. Supposons d'abord que $U$ soit la droite
affine $\AA^1$, d'anneau $k[T]$. Soit~$E$ l'ensemble des
\'el\'ements de~$k[T]$ de la forme $\sum_i a_iT^i$, avec $a_i = 0$
lorsque $p$ divise~$i$. Il est imm\'ediat que l'application compos\'ee
$E\to k[T]\to k[T]/\wp k[T]$ est bijective. D\`es lors les coefficients
$a_i$, $(i,p)=1$, se comportent comme des coordonn\'ees d'un espace qui
param\`etre les rev\^etements de la droite affine cycliques de degr\'e
$p$. En particulier, on en d\'eduit que $\pi_1(\AA^1)$ n'est pas
topologiquement de type fini et que, si l'on effectue une extension
$k\to k'$ de corps alg\'ebriquement clos, le morphisme naturel $\pi_1(\AA^1\times_kk')\to\pi_1(\AA^1)$, qui est surjectif, n'est pas bijectif,
contrairement au cas propre. Dans le cas g\'en\'eral, la courbe $U$ se
r\'ealise comme sch\'ema fini sur la droite affine et on conclut que les
m\^emes ph\'enom\`enes se produisent pour $\pi_1(U)$ (et d'ailleurs, plus
g\'en\'eralement, pour $\pi_1(V)$ pour tout $k$-sch\'ema affine connexe $V$,
de type fini, de dimension $\geq1$).

Ceci \'etant, si $G$ est un groupe fini, notons $G^{(p')}$ le plus grand
groupe quotient de~$G$, d'ordre premier \`a~$p$. Pour qu'un groupe fini
$G$ soit un quotient topologique de~$\pi_1(U)$, il faut que $G^{(p')}$
soit un quotient topologique de~$\pi_1^{p'}(U)$, condition \`a laquelle
on sait en principe r\'epondre gr\^ace au corollaire \Ref{XIII.2.12}. Dans l'article
pr\'ecit\'e, Abhyankar conjecture que cette condition n\'ecessaire est
\'egalement suffisante. Cette conjecture a \'et\'e d\'emontr\'ee par M\ptbl Raynaud
dans le cas de la droite affine et par D\ptbl Harbater dans le cas g\'en\'eral
(Invent. Math.\ \textbf{116} (1994), p\ptbl425--462 et \textbf{117}, p\ptbl1--25). Par exemple,
dans le cas de la droite affine $\pi_1^{p'}(\AA^1)=1$ et on
conclut qu'un groupe fini $G$ est groupe de Galois d'un rev\^etement
\'etale connexe de~$\AA^1$ si et seulement si $G^{(p')} = 1$,
c'est-\`a-dire si et seulement si $G$ est engendr\'e par ses
$p$-sous-groupes de Sylow. Ainsi tout groupe fini simple d'ordre
multiple de~$p$ convient.
\end{remarqueMR}}

\section{Propret\'e cohomologique et locale acyclicit\'e g\'en\'erique}
\label{XIII.3}
\marginpar{395}

\begin{theoreme}
\label{XIII.3.1}
Soient $S$ un sch\'ema irr\'eductible de point
g\'en\'erique~$s$, $X$ et $Y$ deux $S$-sch\'emas de
pr\'esentation finie, $f\colon X \to Y$ un $S$-morphisme. Pour tout
$S$-sch\'ema~$S'$, on note $Y'$, $X'$, etc. l'image inverse de~$Y$,
$X$, etc. par le morphisme $S' \to S$. On a les propri\'et\'es
suivantes:

\begin{enumerate}
\item[1)]
{\upshape a)} On peut trouver un ouvert non vide~$S'$ de~$S$ tel que, pour
tout faisceau d'ensembles constant fini $F'$ sur~$X'$, $f'_*F'$ soit
constructible, et que $(F',f')$ soit cohomologiquement propre
relativement \`a~$S'$ en dimension $\le 0$.

{\upshape b)} Soit $F$ un faisceau d'ensembles constructible sur~$X$. Alors
on peut trouver un ouvert non vide~$S'$ de~$S$ (d\'ependant de~$F$)
tel que $f'_*F'$ soit constructible et que $(F',f')$ soit
cohomologiquement propre relativement \`a~$S'$ en dimension $\le 0$.

\item[2)] Supposons que les sch\'emas de type fini de dimension $\le
\dim X_s$ sur une cl\^oture alg\'ebrique $\overline{k}$ de~$\kres(s)$
soient fortement d\'esingularisables \textup{(SGA~5 I~3.1.5)}. Alors on a de
plus les propri\'et\'es suivantes:

{\upshape a)} On peut trouver un ouvert non vide~$S'$ de~$S$ tel que, pour
tout faisceau en groupes constant fini $F'$ sur~$X'$, d'ordre premier
aux caract\'eristiques r\'esiduelles de~$S$, si $\Phi'$ est le
champ des torseurs sous $F'$, $f'_*\Phi'$ soit $1$-constructible, et
que $(F',f')$ soit cohomologiquement propre relativement \`a~$S'$ en
dimension $\le 1$.

{\upshape b)} Soient $\LL$ l'ensemble des nombres premiers distincts des
caract\'eristiques r\'esiduelles de~$S$ et $\Phi$ un
ind-$\LL$-champ $1$-constructible sur~$X$ \eqref{XIII.0}, tel que,
pour tout sch\'ema~$X_1$ \'etale sur~$X$ et pour tout couple
d'objets $x, x_1$ de~$\Phi_{X_1}$, le faisceau $\SheafHom
_{X_1}(x,x_1)$ soit constructible. On suppose de plus $S$ localement
noeth\'erien. Alors on peut trouver un ouvert non vide~$S'$ de~$S$
tel que $f'_*\Phi'$ soit $1$-constructible, que, pour tout couple
d'objets $y, y_1$ d'une fibre $(f'_*\Phi')_{Y_1}$,
$\SheafHom_{Y_1}(y,y_1)$ soit constructible,
\marginpar{396}
et que $(\Phi',f')$ soit cohomologiquement propre relativement \`a~$S'$ en dimension $\le 1$.
\end{enumerate}
\end{theoreme}

\subsubsection*{D\'emonstration}
On peut suppose $S$ affine; d'apr\`es SGA~4 VIII~1.1, on peut
supposer $S$ int\`egre; enfin, par passage \`a la limite, on peut
supposer que $S$ est le spectre d'une alg\`ebre de type fini sur~$\ZZ$; en particulier $S$ est alors noeth\'erien. Comme la question
est locale sur~$Y$, on peut supposer $Y$ affine. De plus il suffit,
pour d\'emontrer le th\'eor\`eme, de le faire apr\`es
extension finie $S'\to S$, o\`u $S'$ est un sch\'ema int\`egre
et o\`u $S'\to S$ est compos\'e de morphismes \'etales et de
morphismes finis radiciels surjectifs.

\Subsection*{1) Cas des faisceaux d'ensembles constants}

1) 1. R\'eduction au cas o\`u $X$ est normal sur~$S$. Soit
$X_{1\overline{s}}$ le normalis\'e de~$(X_{\overline{s}})_{\red}$;
quitte \`a restreindre $S$ \`a un ouvert non vide et \`a faire
une extension radicielle de~$S$, on peut supposer que
$X_{1\overline{s}}$ provient d'un sch\'ema $X_1$ normal sur~$S$, et
que le morphisme $X_{1\overline{s}}\to X_{\overline{s}}$ provient d'un
morphisme fini surjectif $p\colon X_1\to X$ (EGA~IV 8.8.2 et
9.6.1). Supposons le th\'eor\`eme d\'emontr\'e pour
$fp$. Quitte \`a restreindre $S$ \`a un ouvert, on peut supposer
que, pour tout faisceau d'ensembles constant $F$ sur~$X$, $(p^*F,fp)$
est cohomologiquement propre relativement \`a~$S$ en dimension $\le
0$ et que $f_*p_*(p^*F)$ est
constructible. D'apr\`es~\Ref{XIII.1.9}, $(p_*p^*F,f)$ est alors
cohomologiquement propre relativement \`a~$S$ en dimension $\le
0$. Le morphisme
$$
F \to p_* p^* F=G
$$
est un monomorphisme. Il en r\'esulte d\'ej\`a, $f_*F$ \'etant
un sous-faisceau de~$f_*G$, que $f_*F$ est constructible (SGA~4
IX~2.9~(ii)) et que $(F,f)$ est cohomologiquement propre relativement
\`a~$S$ en dimension $\le -1$.

Soient $X_2=X_1 \times_{X} X_1$, $q \colon X_2 \to X$ le morphisme
canonique. D'apr\`es \Ref{XIII.1.11}~1),
\marginpar{397}
on a une suite exacte
$$
\xymatrix@C=.5cm{F \ar[r] & G \ar@<2pt>[r]\ar@<-2pt>[r] & q_*q^*F }.
$$

D'apr\`es ce que l'on vient de d\'emontrer, appliqu\'e \`a
$fq$ au lieu de~$f$, on peut supposer, quitte \`a restreindre $S$
\`a un ouvert non vide, que, pour tout faisceau d'ensembles constant
$F$ sur~$X$, $(q^*F,fq)$ est cohomologiquement propre relativement
\`a~$S$ en dimension $\le -1$, donc que $(q_*q^*F,f)$ est
cohomologiquement propre relativement \`a~$S$ en dimension $\le
-1$. Il r\'esulte alors de \Ref{XIII.1.13}~1) que $(F,f)$ est
cohomologiquement propre relativement \`a~$S$ en dimension $\le 0$.

1) 2. R\'eduction au cas o\`u $X$ est normal affine sur~$S$.

Soit $U_s$ un ouvert affine de~$X_s$ dense dans $X_s$. Quitte \`a
restreindre $S$ \`a un ouvert non vide, on peut supposer que $U_s
\to X_s$ se rel\`eve en une immersion ouverte $i\colon U\to X$,
sch\'ematiquement dominante relativement \`a~$S$ (EGA~IV
8.9.1). Comme le morphisme $X \to S$ est normal, on a d'apr\`es
SGA~2 XIV~1.18:
$$
\prof \et_{S- U} (X) \ge 2 \quoi ;
$$
par suite, pour tout faisceau constant $F$ sur~$X$, le morphisme
canonique
$$
F \to i_*i^*F
$$
est un isomorphisme. Il en r\'esulte que, si l'on suppose le
th\'eor\`eme d\'emontr\'e pour $fi$ et $i$, alors, apr\`es
restriction de~$S$ \`a un ouvert non vide, $(i^*F,i)$ et $(i^*F,fi)$
sont cohomologiquement propres relativement \`a~$S$ en dimension
$\le 0$. Il en est donc de m\^eme de~$(F,f)$ (\Ref{XIII.1.6}~2)). Comme de plus $f_* F=(fi)_*(i^*F)$ est constructible, ceci
ach\`eve la r\'eduction.

1) 3. Fin de la d\'emonstration.

On peut supposer $S$ normal (EGA~IV 7.8.3). On peut trouver une
compactification de~$X_s$:
\marginpar{398}
$$
\xymatrix{ X_s \ar[r]^{j_s} \ar[d]_{f_s} & P_s \ar[ld]^{g_s} \\ Y_s }
$$
o\`u $j_s$ est une immersion ouverte dominante et $g_s$ un morphisme
propre; quitte \`a faire une extension radicielle de~$\kres(s)$ et \`a
remplacer $P_s$ par son normalis\'e, ce qui ne change pas $X_s$, on
peut supposer $P_s$ g\'eom\'etriquement normal. Quitte \`a
restreindre $S$ \`a un ouvert non vide et \`a faire une extension
radicielle surjective, on peut supposer que le diagramme ci-dessus
provient d'un diagramme
$$
\xymatrix{ X \ar[r]^{j} \ar[d]_{f} & P \ar[ld]^{g} \\ Y }
$$
o\`u $P$ est un sch\'ema normal sur~$S$, $j$ une immersion ouverte
sch\'ematiquement dominante relativement \`a~$S$ et $g$ un
morphisme propre (EGA~IV 6.9.1, 9.9.4 et 9.6.1). Pour tout faisceau
d'ensembles constant fini $F$ sur~$X$ de valeur $C$, $j_*F$ est le
faisceau constant de valeur $C$ (SGA~4 2.14.1), et il en est de
m\^eme apr\`es tout changement de base $S' \to S$; il en
r\'esulte que $(F,j)$ est cohomologiquement propre relativement
\`a~$S$ en dimension $\le 0$. Il en est donc de m\^eme de~$(F,f)$,
puisque $g$ est propre \eqref{XIII.1.8}. Comme $g$ est propre
$f_*F=g_*C_P$ est constructible, ce qui ach\`eve la
d\'emonstration de 1)\kern2pt a).

\Subsection*{2) Cas d'un faisceau d'ensembles constructible}

Soit $F$ un faisceau d'ensembles constructible sur~$X$. D'apr\`es
SGA~4 IX 2.14~(ii), on peut trouver une famille finie de morphismes
$p_i \colon Z_i \to X$, et, sur chaque $Z_i$, un faisceau d'ensembles
constant fini $C_i$, de sorte que l'on ait un monomorphisme
\marginpar{399}
$$
j \colon F \to \prod_{i} {p_i}_* C_i =G \quoi.
$$
D'apr\`es~1)\kern2pt a),
on peut supposer, quitte \`a restreindre $S$ \`a un ouvert non
vide, que les $(C_i,fp_i)$ sont cohomologiquement propres relativement
\`a~$S$ en dimension $\le 0$, et que les $f_*{p_i}_*C_i$ sont
constructibles. On en conclut d\'ej\`a que $(G,f)$ est
cohomologiquement propre relativement \`a~$S$ en dimension $\le 0$
\eqref{XIII.1.9}, donc que $(F,f)$ est cohomologiquement propre
relativement \`a~$S$ en dimension $\le -1$, et que $f_*F$ est
constructible. Soit~$K$ la somme amalgam\'ee $K= G \coprod_{F} G$;
comme $F$ et $G$ sont constructibles, il en est de m\^eme de~$K$. On
conclut donc de ce qui pr\'ec\`ede que, quitte \`a restreindre
$S$ \`a un ouvert non vide, on peut supposer que $(K,f)$ est
cohomologiquement propre relativement \`a~$S$ en dimension $\le
-1$. Il r\'esulte alors de \Ref{XIII.1.13}~1) que $(F,f)$ est
cohomologiquement propre relativement \`a~$S$ en dimension $\le 0$.

\Subsection*{3) Cas des faisceaux en groupes constants}

Si $F$ est un faisceau en groupes constant sur~$X$, on note $\Phi$ le
champ des torseurs sous $F$.

3) 1. Montrons d'abord que, quitte \`a restreindre $S$ \`a un
ouvert non vide, pour tout faisceau en groupes constant $F$ sur~$X$,
d'ordre premier aux caract\'eristiques r\'esiduelles de~$S$,
$(F,f)$ est cohomologiquement propre relativement \`a~$S$ en
dimension $\le 0$, et que $f_*\Phi$ est constructible.

On se ram\`ene pour cela au cas o\`u $X$ est lisse sur~$S$. Quitte
\`a faire une extension finie de~$\kres(s)$, ce qui est loisible car on
peut la consid\'erer comme compos\'ee d'une extension \'etale et
d'une extension radicielle, on peut trouver un morphisme propre
surjectif $p_s \colon X_{1s}\to X_s$, o\`u $X_{1s}$ est un
sch\'ema lisse sur~$S$ de m\^eme dimension que $X_s$, et, quitte
\`a restreindre $S$ \`a un ouvert non vide, on peut supposer que
$p_s$ provient d'un morphisme propre surjectif $p \colon X_1 \to X$,
o\`u $X_1$ est un sch\'ema lisse sur~$S$
\marginpar{400}
(EGA~IV 9.6.1 et 12.1.6). Soient $X_2=X_1 \times_X X_1$, $q \colon
X_2 \to X$ le morphisme canonique. On a un diagramme exact de champs
sur~$X$
$$
\xymatrix@C=.5cm{\Phi \ar[r] & p_*p^* \Phi \ar@<2pt>[r]\ar@<-2pt>[r] &
q_*q^*\Phi }
$$
(\Ref{XIII.1.11}~2)). Le th\'eor\`eme \'etant suppos\'e
d\'emontr\'e dans le cas lisse, on voit tout d'abord que l'on peut
supposer $f_*p_*p^*\Phi$ constructible; il en est donc de m\^eme de
$f_*\Phi$ (\Ref{XIII.3.1.1} ci-dessous). De plus, d'apr\`es \Ref{XIII.1.6}~2),
on peut supposer $(p_*p^*\Phi,f)$ cohomologiquement propre
relativement \`a~$S$ en dimension $\le 0$; il en r\'esulte que
$(\Phi,f)$ est cohomologiquement propre relativement \`a~$S$ en
dimension $\le -1$. On peut donc supposer que $(q_*q^*\Phi,f)$ est
cohomologiquement propre relativement \`a~$S$ en dimension $\le -1$,
et cela entra\^ine que $(\Phi,f)$ est cohomologiquement propre
relativement \`a~$S$ en dimension $\le 0$ (\Ref{XIII.1.12}~1)).

On se ram\`ene ensuite comme dans 1)\kern2pt2 au cas o\`u $X$ est lisse
et affine sur~$S$. Soit alors
$$
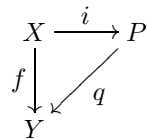
\xymatrix{ X \ar[r]^{i} \ar[d]_{f} & P \ar[ld]^{q} \\ Y }
$$
une compactification de~$X$, o\`u $i$ est une immersion ouverte
dominante et $q$ un morphisme propre. Comme on a $\dim P_s=\dim X_s$,
on peut appliquer l'hypoth\`ese de r\'esolution des
singularit\'es \`a $P_{\overline s}$. Quitte \`a faire une
extension \'etale et une extension radicielle de~$S$, on peut
trouver un morphisme propre $r \colon Z \to P$, o\`u $Z$ est lisse
sur~$S$, $r^{-1}(X)\simeq X$, et o\`u $r^{-1}(X)$ est le
compl\'ementaire dans $Z$ d'un diviseur \`a croisements normaux
relativement \`a~$S$. Tout torseur sous $F$ est alors
mod\'er\'ement ramifi\'e sur~$Z$ relativement \`a~$S$
(\Ref{XIII.2.3}~b)). Il r\'esulte donc de~\Ref{XIII.2.7} que $(F,f)$
est cohomologiquement propre relativement \`a~$S$ en dimension $\le
0$, ce qui d\'emontre notre assertion.

\ifthenelse{\boolean{orig}}
{3) \textbf{2. R\'eduction\marginpar{401} au cas o\`u $X$ est lisse sur~$S$}}%
{\pagebreak[4] 3) 2. R\'eduction\marginpar{401} au cas o\`u $X$ est lisse sur~$S$.}%
%
\par\nopagebreak
Quitte \`a faire une extension finie de~$\kres(s)$, on peut trouver un
morphisme propre surjectif $p_s\colon X_{1s}\to X$, o\`u~$X_{1s}$
est un sch\'ema lisse sur~$s$, et on peut supposer que $p_s$
provient d'un morphisme propre surjectif $p\colon X_1 \to X$, o\`u
$X_1$ est lisse sur~$S$. Supposons le th\'eor\`eme
d\'emontr\'e pour $fp$ et montrons-le pour $f$. Soit $F$ un
faisceau en groupes constant fini sur~$X$, d'ordre premier aux
caract\'eristiques r\'esiduelles de~$S$ et $\Phi$ le champ des
torseurs sous $F$. Soient $X_2=X_1 \times_X X_1$, $X_3=X_1 \times_X
X_1 \times_X X_1$, $q\colon X_2 \to X$, $r\colon X_3 \to X$ les
morphismes canoniques. D'apr\`es \Ref{XIII.1.11}~2), on a un
diagramme exact de champs
$$
\xymatrix@C=.5cm{\Phi \ar[r] & p_*p^*\Phi \ar@<2pt>[r]\ar@<-2pt>[r] &
q_*q^*\Phi \ar@<4pt>[r] \ar[r] \ar@<-4pt>[r] & r_*r^*\Phi}\quoi.
$$
Pour prouver que $(\Phi,f)$ est cohomologiquement propre relativement
\`a~$S$ en dimension $\le 1$, il suffit de montrer qu'il en est de
m\^eme de~$(p_*p^*\Phi,f)$, que $(q_*q^*\Phi, f )$ est
cohomologiquement propre relativement \`a~$S$ en dimension $\le 0$
et que $(r_*r^*\Phi,f)$ est cohomologiquement propre relativement
\`a~$S$ en dimension $\le -1$ (\Ref{XIII.1.12}~2)). D'apr\`es~3)\kern.3em 1
ci-dessus on peut supposer que, pour tout faisceau en groupes constant
fini~$F$, $(q^*\Phi,fq)$, $(q^*\Phi,q)$, $(r^*\Phi,fr)$, $(r^*\Phi,r)$
sont cohomologiquement propres relativement \`a~$S$ en dimension
$\le 0$. Il r\'esulte alors de \Ref{XIII.1.6}~2) que
$(q_*q^*\Phi,f)$ et $(r_*r^*\Phi,f)$ sont cohomologiquement propres
relativement \`a~$S$ en dimension $\le 0$. Le th\'eor\`eme
\'etant suppos\'e d\'emontr\'e dans le cas lisse, $(p^*\Phi,
fp)$ et $(p^*\Phi,p)$ donc aussi $(p_*p^*\Phi,f)$ (\Ref{XIII.1.6}~2))
sont cohomologiquement propres relativement \`a~$S$ en dimension
$\le 1$. Ceci montre bien que $(\Phi,f)$ est cohomologiquement propre
relativement \`a~$S$ en dimension $\le 1$.

De plus $f_*p_*p^*\Phi$ est $1$-constructible par hypoth\`ese;
d'apr\`es~\Ref{XIII.3.1} on peut supposer que $f_*q_*q^*\Phi$ est
constructible; il r\'esulte donc de~\Ref{XIII.3.1.1} ci-dessous que
$f_*\Phi$ est $1$-constructible.

3) 3. R\'eduction au cas o\`u $X$ est lisse affine sur~$S$.
\marginpar{402}

D'apr\`es 3)\kern2pt2, on peut supposer $X$ lisse sur~$S$. Soit $(U_i)_{i\in
I}$ un recouvrement fini de~$X$ par des ouverts affines, et soit $X_1$
la somme directe des $U_i$, $p\colon X_1 \to X$ le morphisme
canonique. Comme $p$ est un morphisme de descente effective pour la
cat\'egorie des faisceaux \'etales de type fini sur des
sch\'emas variables, on voit comme dans 3)\kern2pt2 que, si l'on suppose le
th\'eor\`eme d\'emontr\'e pour les $U_i$, \ie pour $X_1$, il
est aussi vrai pour~$X$.

3) 4. Cas o\`u $X$ est lisse affine sur~$S$.

On voit comme dans 3)1 que, quitte \`a restreindre $S$ \`a un
ouvert non vide et \`a faire une extension \'etale et une
extension radicielle surjectives de~$S$, on peut trouver un diagramme
commutatif
$$
\xymatrix{ X \ar[r]^{i} \ar[d]_{f} & P \ar[ld]^{q} \\ Y }
$$
o\`u $P$ est un sch\'ema lisse sur~$S$, $X$ le compl\'ementaire
dans $P$ d'un diviseur \`a croisements normaux relativement \`a~$S$ et $g$ un morphisme propre. Si $F$ est un faisceau constant
d'ordre premier aux caract\'eristiques r\'esiduelles de~$S$, tout
torseur sous $F$ est mod\'er\'ement ramifi\'e sur~$P$
relativement \`a~$S$ (\Ref{XIII.2.3}~b)). Le fait que $(F,f)$ soit
cohomologiquement propre relativement \`a~$S$ en dimension $\le 1$
et que $f_*\Phi$ soit constructible r\'esulte alors
de~\Ref{XIII.2.7}.

\ifthenelse{\boolean{orig}}
{4) \textbf{D\'emonstration de } 2)\kern2pt b)}
{\Subsection*{4) D\'emonstration de 2)\kern2pt b)}}

4) 1. Cas o\`u $\Phi$ est une gerbe.

On peut trouver un morphisme \'etale surjectif de type fini $p
\colon X_1 \to X$, tel que $p^*\Phi$ soit un gerbe triviale. Par
descente, comme dans 3)\kern2pt2, on voit qu'il suffit de prouver le
th\'eor\`eme pour $X_1$, $X_1\times_X X_1$ et
\marginpar{403}
$X_1\times_X X_1\times_X X_1$. On est donc ramen\'e au cas o\`u
$\Phi$ est le gerbe des torseurs sous un faisceau en groupes
constructible $F$, dont les fibres sont d'ordre premier aux
caract\'eristiques r\'esiduelles de~$S$.

D'apr\`es SGA~4 IX~2.14, on peut trouver une famille finie de
morphismes finis $p_i \colon Z_i \to X$, et, pour chaque $i$, un
faisceau en groupes constant fini $C_i$, d'ordre premier aux
caract\'eristiques r\'esiduelles de~$S$, de sorte que l'on ait un
\ifthenelse{\boolean{orig}}
{morphisme}
{monomorphisme}
$$
j\colon F \to \prod_{i} p_{i*}C_i=G \quoi.
$$
Soit $\Phi_i$ le champ de torseurs sous $C_i$ et $\Psi$ le champ des
torseurs sous $G$. Il r\'esulte de 2)\kern2pt a) que, quitte \`a
restreindre $S$ \`a un ouvert non vide, on peut supposer que les
$(C_i,fp_i)$ sont cohomologiquement propres relativement \`a~$S$ en
dimension $\le 1$, et que les champs $f_*p_{i*}\Phi_i$ sont
$1$-constructibles. Il r\'esulte alors de~\Ref{XIII.1.9} que les
$(p_{i*}C_i,f)$ sont cohomologiquement propres relativement \`a~$S$
en dimension $\le 1$; il en est donc de m\^eme de~$(G,f)$. De plus,
comme $p_{i*}\Phi_i$ est \'equivalent au champ des torseurs sous le
groupe $p_{i*}C_i$ (SGA~4 VIII~5.8), on voit que $f_*\Psi$ est
$1$-constructible.

Comme $\R^1f_*G$ est constructible, on peut trouver un faisceau
repr\'esentable par un $Y$-sch\'ema \'etale de type fini $T$ et
un \'epimorphisme
$$
a \colon T \to \R^1f_*G
$$
(SGA~4 IX~2.7); de plus on peut supposer que l'image de la section
identique de~$T(T)$ est d\'efinie par un torseur $Q$ sur~$X\times_Y
T=X_T$, de groupe $G|X_T$. Soient $f_T \colon X_T \to T$ le morphisme
canonique, $F_T=F|X_T$, etc. D'apr\`es 1)\kern2pt b) on peut supposer,
quitte \`a restreindre $S$ \`a un ouvert non vide, que
$(Q/F_T,f_T)$ est cohomologiquement propre relativement \`a~$S$ en
dimension $\le 0$ et que $f_{T*}(Q/F_T)$ est constructible. Il
r\'esulte alors de~\Ref{XIII.3.1.2} que $f_*\Psi$ est constructible.

Montrons que $(F,f)$ est cohomologiquement propre relativement \`a
\marginpar{404}
$S$ en dimension $\le 1$. D'apr\`es \Ref{XIII.1.13}~2) il suffit de
prouver que, pour tout sch\'ema $Y_1$ \'etale sur~$Y$ et pour tout
torseur $Q_1$ sur~$X_1=X\times_Y Y_1$, si $f_1 \colon X_1 \to Y_1$ est
le morphisme canonique, alors $(Q_1/F_1,f_1)$ est cohomologiquement
propre relativement \`a~$S$ en dimension $\le 0$. Or, par
d\'efinition de~$T$, $Q_1$ est, localement pour la topologie
\'etale de~$Y_1$, image inverse de~$Q$, ce qui d\'emontre notre
r\'eduction.

4) 2. Cas g\'en\'eral.\par\nopagebreak
On voit en utilisant le lemme~\Ref{XIII.6.1.1}, 4)\kern2pt1 et 1)\kern2pt a) que,
quitte \`a restreindre $S$ \`a un ouvert non vide, on peut
supposer $S(f_*\Phi)$ constructible et $(S\Phi,f)$ cohomologiquement
propre relativement \`a~$S$ en dimension $\le 0$. On peut alors
trouver un faisceau repr\'esentable par un $Y$-sch\'ema \'etale
de type fini $T$ et un \'epimorphisme
$$
a \colon T \to S(f_*\Phi)\quoi.
$$
On reprend les notations de 4)\kern2pt1, et on pose de plus $Z=T\times_Y T$,
$X_Z=X\times_Y Z$ et on note $f_Z \colon X_Z \to Z$ le morphisme
canonique. On peut supposer $T$ choisi de sorte que l'image $q$ de la
section identique de~$T(T)$ par $a$ soit d\'efinie par un objet~$p$
de~$(f_*\Phi)_T=\Phi_{X_T}$. Soient $p_1$ et $p_2$ (\resp~$q_1$ et
$q_2$) les images inverses de~$p$ (\resp~$q$) par les deux projections
de~$Z$ dans $T$. On peut supposer, quitte \`a restreindre $S$ \`a
un ouvert non vide, que $(\SheafAut_{X_T}(p),f_T)$ est
cohomologiquement propre relativement \`a~$S$ en dimension $\le 1$,
que $(\SheafHom_{X_Z}(p_1,p_2),f_Z)$ est cohomologiquement propre
relativement \`a~$S$ en dimension $\le 0$, et que
$f_{T*}(\SheafAut_{X_T}(p))$, $f_{Z*}(\SheafHom_{X_Z}(p_1,p_2))$ sont
constructibles ( a)\kern2pt1) et 4)\kern2pt1)).

On en d\'eduit d'abord que, pour tout sch\'ema $Y_1$ \'etale sur~$Y$ et pour tout couple d'objets $y$, $y_1$ de~$(f_*\Phi)_{Y_1}$, le
faisceau $\SheafAut_{Y_1}(y)$ (\resp $\SheafHom_{Y_1}(y,y_1)$) est
constructible, un tel faisceau \'etant, localement pour la topologie
\'etale de~$Y_1$, image inverse de~$f_{T*}(\SheafAut_{X_T}(p))$
(\resp de~$f_{Z*}(\SheafHom_{X_Z}(p_1,p_2))$).

Il reste \`a prouver que $(\Phi,f)$ est cohomologiquement propre
relativement
\marginpar{405}
\`a~$S$ en dimension $\le 1$. Il suffit pour cela de montrer que,
pour tout $S'$-sch\'ema $S$ et pour tout point g\'eom\'etrique
$\overline{y}'$ de~$Y$, si l'on note $\overline{Y}$ le localis\'e
strict de~$Y$ en $y'$, $\overline{X}=X \times_Y \overline{Y}$, etc.,
alors le foncteur canonique
$$
\overline{\varphi} \colon \Phi(\overline{X}) \to
\Phi'(\overline{X'})
$$
est une \'equivalence de cat\'egories.

Montrons que $\overline{\varphi}$ est pleinement fid\`ele. Soient
$x$, $y \in \Phi(\overline{X})$, $x'$, $y'$ leurs images dans
$\Phi'(\overline{X'})$ et montrons que le morphisme canonique
\begin{equation*}
\label{eq:XIII.4.*} \tag{$*$} {\Hom_{\overline{X}}(x,y)
\to \Hom_{\overline{X'}}(x',y') }
\end{equation*}
est bijectif. Par d\'efinition de~$T$ il existe deux morphismes de
$\overline{Y}$ dans $T$ tels que $x$ et~$y$ soient les images inverses
de~$p$ par ces deux morphismes. Cela revient \`a dire qu'il y a un
morphisme $\overline {Y} \to Z$, d'o\`u un morphisme $h\colon
\overline{X} \to X_Z$ tel que l'on ait
$$
h^*(p_1)=x \qquad h^*(p_2)=y.
$$
Par suite on a un isomorphisme canonique
$$
\SheafHom_{\overline{X}}(x,y)=h^*(\SheafHom_{X_Z}(p_1,p_2) )\quoi.
$$
mais, compte tenu du fait que $(\SheafHom_{X_Z}(p_1,p_2),f_Z)$ est
cohomologiquement propre relativement \`a~$S$ en dimension $\le 0$,
on voit qu'il en est de m\^eme de
$(\SheafHom_{\overline{X}}(x,y),\overline{f})$, ce qui prouve que le
morphisme \eqref{eq:XIII.4.*} est bijectif.

Montrons que $\overline{\varphi}$ est essentiellement surjectif. Soit
$x'\in \Phi'(\overline{X'})$. Comme $(S\Phi,f)$ est cohomologiquement
propre relativement \`a~$S$ en dimension $\le 0$, le morphisme
canonique
$$
\H^0(\overline{X},\overline{S\Phi}) \to
\H^0(\overline{X'},\overline{S\Phi}')
$$
est bijectif. Soit $G'$ la sous-gerbe maximale de~$\Phi'$
engendr\'ee par $x'$; il existe alors une sous-gerbe maximale $G$ de
$\overline{\Phi}$ dont l'image inverse sur~$\overline{X'}$ est
$G'$. D'apr\`es 4)\kern2pt1 $(G,\overline{f})$ est cohomologiquement propre
relativement \`a~$S$ en dimension $\le 1$; par suite le foncteur
canonique
$$
G(\overline{X}) \to G'(\overline{X'})
$$
\marginpar{406}%
est une \'equivalence de cat\'egories, ce qui prouve l'existence
d'un \'el\'ement $x$ de~$\Phi(\overline{X})$ dont l'image dans
$\Phi'(X')$ soit isomorphe \`a $x'$ et ach\`eve la
d\'emonstration du th\'eor\`eme.

\begin{sublemme}
\label{XIII.3.1.1}
Soient $S$ un sch\'ema localement noeth\'erien et
$$
\xymatrix@C=.5cm{ \Phi \ar[r]^p & \Phi_1 \ar@<2pt>[rr]^{p_1,p_2}
\ar@<-2pt>[rr] && \Phi_2 }
$$
un diagramme exact de champs sur~$S$ \eqref{XIII.1.10.1}. Si $\Phi_1$
est constructible, il en est de m\^eme de~$\Phi$. Si, pour tout
sch\'ema $S'$ \'etale sur~$S$ et pour tout couple
\ifthenelse{\boolean{orig}}
{d'objet}
{d'objets}
$x_1,~y_1$ de~$(\Phi_1)_{S'}$ le faisceau $\SheafHom_{S'}(x_1,y_1)$
est constructible, alors, pour tout couple d'objets $x,~y$ de
$\Phi_{S'}$, il en est de m\^eme de~$\SheafHom_{S'}(x,y)$. Supposons
que $\Phi_1$ soit $1$-constructible (\Ref{XIII.0}) et que $\Phi_2$ soit
constructible, alors $\Phi$ est $1$-constructible.
\end{sublemme}

Pour tout sch\'ema $S'$ \'etale sur~$S$ et pour tout objet $x$ de
$\Phi_{S'}$, on a un monomorphisme
$$
\SheafAut_{S'}(x)\to\SheafAut_{S'}(p(x)).
$$
Il r\'esulte donc de SGA~4 IX~2.9 que, si $\Phi_1$ est
constructible, il en est de m\^eme de~$\Phi$, et la deuxi\`eme
assertion du lemme se d\'emontre de la m\^eme fa\c con.

Supposons maintenant $\Phi_1$ $1$-constructible et $\Phi_2$
constructible. Le morphisme $p$ induit sur les faisceaux de
sous-gerbes maximales un morphisme
$$
\varphi\colon S\Phi\to S\Phi_1.
$$
Soit $G$ l'mage de~$S\Phi$ par $\varphi$; d'apr\`es SGA~4 IX~2.9,
$G$ est un faisceau constructible. On peut donc trouver un faisceau
repr\'esentable par un $S$-sch\'ema \'etale de type fini $T$ et
un \'epimorphisme
$$
a\colon T\to G
$$
(SGA~4 IX~2.7).
\marginpar{407}
De plus on peut choisir $T$ de sorte que l'image $y$ de section
identique de~$T(T)$ par $a$ soit d\'efinie par un objet $x_1$ de
$(\Phi_1)_T$ de la forme $x_1=p(x)$, o\`u $x\in\Phi_T$.

Il suffit de montrer que, pour tout point $s$ de~$S$, il existe un
ouvert non vide~$U$ de~$\overline{\{ s\} }$ tel que $S\Phi|U$ soit
localement constant constructible. Soient $s\in S$, $\overline{s}$ un
point g\'eom\'etrique au-dessus de~$s$ et $\overline{y}_1,\dots,\overline{y}_n$ les \'el\'ements de~$G_{\overline{s}}$. Par
d\'efinition de~$T$ il existe des morphismes
$h_i\colon\overline{s}\to T$ tels que l'on ait
$h_i^*(y)=\overline{y}_i$. Soit $S'$ le produit fibr\'e sur~$S$ de
$n$ sch\'emas isomorphes \`a $T$, $h\colon\overline{s}\to S'$ le
produit fibr\'e des $h_i$, $y_i$ (\resp $x_i$) l'image inverse de
$y$ (\resp $x$) par la $i$-\`eme projection de~$S'$ sur~$T$. Soit
$F_i$ le sous-faisceau de~$S\Phi |S'$ image r\'eciproque de~$y_i$
et montrons que les $F_i$ sont constructibles. Le faisceau $F_i$ est
un quotient du faisceau $F'_i$ tel que, pour tout sch\'ema $S''$
\'etale au-dessus de~$S'$, on ait
\begin{multline*}
F'_i(S'') = \{\text{classes mod. isomorphisme d'objet $z$ de
$\Phi_{S''}$ muni}\\
\ifthenelse{\boolean{orig}}
{\text{ d'un isomorphisme $i\colon p(z)\simeq p(x_i|S'')$}.}
{\text{ d'un isomorphisme $i\colon p(z)\simeq p(x_i|S'')$}\}.}
\end{multline*}
Il suffit de montrer que les $F_i'$ sont constructibles. Or, si l'on
pose $z_i=p_1p(x_i)$, $z_i'=p_2p(x_i)$, on a un monomorphisme
$$
\Psi_i\colon F'_i\to\SheafIsom_{S'}(z_i,z_i'),
$$
obtenu en associant, \`a tout sch\'ema $S''$ \'etale sur~$S'$ et
\`a tout objet $z$ de~$\Phi_{S''}$ tel que l'on ait un isomorphisme
$i\colon p(z)\simeq p(x_i|S'')$, l'isomorphisme de~$z_i$ dans $z_i'$
d\'efini par la condition que le diagramme
$$
\xymatrix@C=1.2cm{
p_1p(z) \ar[r]^-{p_1(i)}\ar[d]_-j& z_i\ar[d]\\ p_2p(z) \ar[r]^-{p_2(i)}& z'_i
}
$$
soit commutatif ($j$ est le morphisme canonique associ\'e au
diagramme exact).
\marginpar{408}

Le morphisme $\Psi_i$ est injectif car dire que deux objets $z,~z'$,
tels que l'on ait des isomorphismes $i\colon p(z)\isomto
p(x_i|S''),~i'\colon p(z')\isomto p(x_i|S'')$, d\'efinissent le
m\^eme \'el\'ement de~$\SheafIsom_{S''}(z_i,z'_i)$ revient \`a
dire que l'on a $p_1(i'^{-1}i)=p_2(i'^{-1}i)$, \ie que $i'^{-1}i$
provient d'un isomorphisme $z\to z'$. Le faisceau $F'_i$ \'etant un
sous-faisceau de~$\SheafIsom_{S'}(z_i,z'_i)$ est constructible.

Comme on peut trouver un ouvert non vide~$U$ de~$\overline{\{ s\}}$
tel que $y_1|U,\dots,y_n|U$ engendrent $G$, il r\'esulte du
lemme~\Ref{XIII.6.1.2} ci-dessous que $S\Phi |U$ est constructible,
donc, quitte \`a restreindre $U$, $S\Phi |U$ est localement constant
constructible.

\begin{sublemme}
\label{XIII.3.1.2}
Soient $S$ un sch\'ema localement noeth\'erien, $f\colon X\to S$
un morphisme, $F\to G$ un monomorphisme de faisceaux en groupes sur~$X$, $\Psi$ (\resp $\Psi_1$) le champ des torseurs sous $F$
(\resp sous $G$), $\Phi=f_*\Psi$, $\Phi_1=f_*\Psi_1$. On suppose
donn\'e un faisceau sur~$S$ repr\'esentable par un $S$-sch\'ema
\'etale de type fini $T$ et un morphisme surjectif
$$
a\colon T\to S\Phi_1\simeq\R^1f_*G,
$$
de sorte qu'il existe un torseur $Q$ sur~$X_T=X\times_ST$ de groupe
$G|X_T$ qui d\'efinisse dans $\R^1f_*G(T)$ l'image par $a$ de la
section identique de~$T(T)$. Soit $f_T\colon X_T\to T$ le morphisme
canonique et posons $F_T=F|X_T$. On suppose que $\Phi_1$ est
$1$-constructible et que $f_{T*}(Q/F_T)$ est constructible. Alors
$\Phi$ est $1$-constructible.
\end{sublemme}

Il suffit en effet de recopier la d\'emonstration
de~\Ref{XIII.3.1.1}, le fait que l'on ait des morphismes $\Psi_i$
\'etant remplac\'e par le fait que l'on a des isomorphismes
$$
F'_i\isomto f_*(Q/F_T)|S'.
$$

\begin{subremarque}
\label{XIII.3.1.3}
Supposons
\marginpar{409}
que $k=\kres(s)$ soit de caract\'eristique nulle; alors les
sch\'emas de type fini sur~$\overline{k}$ de dimension $\leq\dim X$
sont fortement d\'esingularisables et la d\'emonstration
de~\Ref{XIII.3.1} permet de prouver les r\'esultats suivants:
\begin{enumerate}
\item[a)] Il existe un ouvert non vide~$S_1$ de~$S$ tel que, pour tout
sch\'ema $S'$ au-dessus de~$S_1$ dont les points maximaux sont de
caract\'eristique nulle et pour tout faisceau d'ensembles localement
constant constructible $F$ sur~$X'=X\times_SS'$, $(F,f')$ soit
cohomologiquement propre relativement \`a~$S'$ en dimension $\leq
0$.

\item[b)] Si toutes les caract\'eristiques r\'esiduelles de~$S$ sont
nulles, il existe un ouvert non vide~$S_1$ de~$S$ tel que, pour tout
sch\'ema $S'$ au-dessus de~$S_1$ et pour tout faisceau en groupes
localement constant constructible $F$ sur~$X'$, $(F,f')$ soit
cohomologiquement propre relativement \`a~$S'$ en dimension $\leq1$.
\end{enumerate}

Il suffit en effet de recopier la d\'emonstration de
\Ref{XIII.3.1}.2)\kern2pt a). La proposition~\Ref{XIII.2.7} utilis\'ee dans
3)\kern2pt 4 s'applique au cas d'un faisceau localement constant $F$, car,
reprenant les notations de 3)\kern2pt4, tout torseur sous $F$ est
mod\'er\'ement ramifi\'e sur~$P$ relativement \`a~$S$ puisque
toutes les caract\'eristiques r\'esiduelles de~$S$ sont nulles.
\end{subremarque}

\begin{corollaire}
\label{XIII.3.2}
Soient $k$ un corps de caract\'eristique $p\geq 0$, $p'$ l'ensemble
des nombres premiers distincts de~$p$, $f\colon X\to k$ un morphisme
coh\'erent.
\begin{enumerate}
\item[1)] Pour tout faisceau d'ensembles $F$, $(F,f)$ est
cohomologiquement propre en dimension $\leq 0$.
\item[2)] Supposons satisfaite l'une des deux conditions suivantes:
\begin{enumerate}
\item[a)] $f$ est de type fini et les sch\'emas de type fini de
dimension $\leq\dim X$ sur une cl\^oture alg\'ebrique de~$k$ sont
fortement d\'e\-sin\-gu\-la\-ri\-sables.
\item[b)] Les sch\'emas de type fini sur une cl\^oture
alg\'ebrique de~$k$ sont
\marginpar{410}
fortement d\'esingularisables.
\end{enumerate}
Alors, pour tout faisceau de ind-$p'$-groupe $F$, $(F,f)$ est
cohomologiquement propre en dimension $\leq 1$.
\end{enumerate}
\end{corollaire}

Soit $F$ un faisceau d'ensembles (\resp de
ind-$p'$-groupes). D'apr\`es SGA~4 IX~2.7.2 on peut \'ecrire
$F$ comme limite inductive filtrante
$$
F=\varinjlim F_i,
$$
o\`u les $F_i$ sont des faisceaux d'ensembles (\resp de
ind-$p'$-groupes) constructibles. Comme $f$ est coh\'erent, $f_*$
(\resp $\R^1f_*$) commute aux limites inductives (SGA~4 VII~3.3). Si
l'on sait que les $(F_i,f)$ sont cohomologiquement propres en
dimension $\leq 0$ (\resp que tout faisceau d'ensemble est
cohomologiquement propre en dimension $\leq 0$ et que les $(F_i,f)$
sont cohomologiquement propres en dimension $\leq 1$), il en sera de
m\^eme de~$(F,f)$. On peut donc supposer $F$ constructible.

Si l'on suppose $f$ de type fini (\resp satisfaisant \`a a)), la
proposition r\'esulte de \Ref{XIII.3.1}~1)\kern2pt b) (\resp de
\Ref{XIII.3.1}~2)\kern2pt b)). Prouvons maintenant~\Ref{XIII.3.2} quand on ne
suppose plus $f$ de type fini. Pour tout sch\'ema $S'$ au-dessus de
$k$ et pour tout point g\'eom\'etrique $\overline{s}$ de~$S'$, on
note $\overline{k}$ (\resp $\overline{S'}$) le localis\'e strict de
$k$ en $\overline{s}$ (\resp de~$S'$ en $\overline{s}$),
$\overline{X}$ l'image inverse de~$X$ sur~$\overline{k}$, et on
consid\`ere le carr\'e cart\'esien
$$
\xymatrix@C=1.2cm{
\overline{X}\ar[d]_-{\overline{f}}&\ar[l]_g \overline{X'}\ar[d]^-{\overline{f}'}\\ \overline{k}&\ar[l]\overline{S'}
}
$$
Il suffit de prouver que, pour tout $S'$, pour tout $\overline{s}$, le
morphisme canonique
$$
\H^0(\overline{X},\overline{F})\to
\H^0(\overline{X'},\overline{F'})\qquad
\text{(\resp $\H^1(\overline{X},\overline{F})\to
\H^1(\overline{X'},\overline{F'})$)}
$$
est un isomorphisme. Il suffit de montrer que l'on a les relations
\marginpar{411}
\begin{equation*}
\label{eq:XIII.3.*}
\tag{$*$} \overline{F}\simeq g_*g^*\overline{F}\qquad
\ifthenelse{\boolean{orig}}
{\text{(\resp $\R^1g_*(g^*\overline{F}))=0$).}}
{\text{(\resp $\R^1g_*(g^*\overline{F})=0$).}}
\end{equation*}
Or, sous cette forme, la question est locale sur~$X$ pour la topologie
\'etale. On peut donc supposer $X$ affine; par passage \`a la
limite on peut supposer $X$ de type fini sur~$k$. On sait alors que
$(F,f)$ est cohomologiquement propre en dimension $\leq 0$
(\resp $\leq 1$) et qu'il en est de m\^eme quand on remplace $X$ par
un sch\'ema \'etale de type fini sur~$X$, ce qui prouve
\eqref{eq:XIII.3.*}.

\begin{theoreme}
\label{XIII.3.3}
Soient $S$ un sch\'ema irr\'eductible de point g\'en\'erique
$s$, $f\colon X\to S$ un morphisme de pr\'esentation
finie. Supposons que les sch\'emas de type fini de dimension
$\leqslant\dim X_s$ sur une cl\^oture alg\'ebrique $\bar{k}$ de
$k$ soient d\'esingularisables \textup{(EGA~IV 7.9.1)}. Alors, si $\LL$
d\'esigne l'ensemble des nombres premiers distincts des
caract\'eristiques r\'esiduelles de~$S$, on peut trouver un ouvert
non vide~$S_1$ de~$S$ tel que le morphisme $f|S_1$ soit
universellement localement $1$-asph\'erique pour~$\LL$.
\end{theoreme}

On peut supposer $S$ int\`egre et $X$ r\'eduit (SGA~4
VIII~1.1). Par passage \`a la limite on peut supposer $S$
noeth\'erien. De plus, pour d\'emontrer le th\'eor\`eme, il
suffit de le faire apr\`es extension finie $S_1\to S$, o\`u $S_1$
est un sch\'ema int\`egre et o\`u $S_1\to S$ est compos\'e
d'extensions \'etales et d'extensions radicielles surjectives.

Montrons d'abord que, quitte \`a restreindre $S$ \`a un ouvert non
vide, $f$ est universellement localement $0$-acyclique. Quitte \`a
faire une extension radicielle de~$\kres(s)$, on peut supposer que le
morphisme $(X_s)_\red\to s$ est s\'eparable; on peut donc supposer,
quitte \`a restreindre $S$ \`a un ouvert non vide et \`a faire
une extension radicielle surjective de~$S$, que le morphisme $f$ est
plat, \`a fibres g\'eom\'etriques s\'eparables (EGA~IV
12.1.1), ce qui entra\^ine que $f$ est universellement
$0$-acyclique (SGA~4 XV~4.1).

Montrons que, quitte \`a restreindre $S$ \`a un ouvert non vide,
$f$ est universellement localement $1$-asph\'erique pour
$\LL$. Quitte \`a faire une
\marginpar{412}
extension finie de~$\kres(s)$, ce qui est loisible car on peut la
consid\'erer comme compos\'ee d'une extension \'etale et d'une
extension radicielle, on peut trouver un morphisme propre surjectif
$p_s\colon Y_s\to X_s$, o\`u $Y_s$ est un sch\'ema lisse sur~$S$,
de m\^eme dimension que $X_s$, et on peut supposer que $p_s$
provient d'un morphisme propre surjectif $p\colon Y\to X$, o\`u $Y$
est un sch\'ema lisse sur~$S$ (EGA~IV 9.6.1 et 12.1.6). Il suffit de
montrer que, quitte \`a restreindre $S$ \`a un ouvert non vide,
pour tout diagramme \`a carr\'es cart\'esiens
$$
\xymatrix{ S''\ar[d]_i & X ''\ar[d]^j\ar[l]_{f''} \\ S'\ar[d] &
X'\ar[d]\ar[l]_{f'} \\ S & X, \ar[l]_f }
$$
o\`u $i$ est \'etale de pr\'esentation finie, et pour tout
faisceau de ind-$\LL$-groupes $F$ sur~$S''$, si $\Phi$ est le champ
des torseurs sous $F$, alors le morphisme canonique
$$
{f'}^*i_*\Phi \to j_*{f''}^*\Phi
$$
est une \'equivalence. Soient $Z=Y\times_{X}Y,
T=Y\times_{X}Y\times_{X}Y$. On a de fa\c con naturelle un diagramme
commutatif
$$
\xymatrix{ S'' \ar[d]_i & X''\ar[l]_{f''}\ar[d]_j &
Y''\ar[l]_{p''}\ar[d] & Z''\ar[d]\ar@<-.4ex>[l]\ar@<.4ex>[l] &
T''\ar[d]\ar[l]\ar@<-.8ex>[l]\ar@<.8ex>[l] \\ S' \ar[d] & X'
\ar[l]_{f'}\ar[d] & Y'\ar[d]\ar[l]_{p'}&
Z'\ar[d]\ar@<-.4ex>[l]\ar@<.4ex>[l] &
T'\ar[d]\ar[l]\ar@<-.8ex>[l]\ar@<.8ex>[l] \\ S & X\ar[l]_f & Y\ar[l]_p
& Z'\ar@<-.4ex>[l]\ar@<.4ex>[l] & T'\ar[l]\ar@<-.8ex>[l]\ar@<.8ex>[l]
}
$$
Soient
\marginpar{413}%
$q\colon Z\to X$ et $r\colon T\to X$ les morphismes canoniques,
les notations $q', r', q'', r''$ ayant un sens
\'evident. D'apr\`es \Ref{XIII.1.11}~2), on a le diagramme
essentiellement commutatif suivant, dont les lignes sont exactes:
$$
\ifthenelse{\boolean{orig}}
{\xymatrix{ f'^* i_*\Phi \ar[r]\ar[d]_a &
p'_*p'^*(f'^*i_x\Phi)\ar[d]_b\ar@<-.4ex>[r]\ar@<.4ex>[r] &
q'_*q'^*(f'^*i_*\Phi)\ar[d]_c\ar[r]\ar@<-.8ex>[r]\ar@<.8ex>[r]&
r'_*r'^*(f'^*i_*\Phi)\ar[d]_d \\ j_*f''^*\Phi \ar[r]&
j_*p''_*p''^*(f''^*\Phi)\ar@<-.4ex>[r]\ar@<.4ex>[r] &
j_*q''_*q''^*(f''^*\Phi)\ar[r]\ar@<-.8ex>[r]\ar@<.8ex>[r]&
j_*r''_*r''^*(f''^*\Phi) }}
{\xymatrix{ f'^* i_*\Phi \ar[r]\ar[d]_a &
p'_*p'^*(f'^*i_*\Phi)\ar[d]_b\ar@<-.4ex>[r]\ar@<.4ex>[r] &
q'_*q'^*(f'^*i_*\Phi)\ar[d]_c\ar[r]\ar@<-.8ex>[r]\ar@<.8ex>[r]&
r'_*r'^*(f'^*i_*\Phi)\ar[d]_d \\ j_*f''^*\Phi \ar[r]&
j_*p''_*p''^*(f''^*\Phi)\ar@<-.4ex>[r]\ar@<.4ex>[r] &
j_*q''_*q''^*(f''^*\Phi)\ar[r]\ar@<-.8ex>[r]\ar@<.8ex>[r]&
j_*r''_*r''^*(f''^*\Phi) }}
$$
Comme $Y$ est lisse sur~$S$, le morphisme $fp$ est universellement
localement $1$-asph\'erique pour $\LL$ (SGA~XV~2.1),
et il r\'esulte de \cite[VII 2.1.7]{XIII.2} que $b$ est une \'equivalence de
cat\'egories. D'autre part, quitte \`a restreindre $S$ \`a un
ouvert non vide, on peut supposer que les morphismes $Z\to S$ et $T\to
S$ sont universellement localement $0$-acycliques. Il en r\'esulte
que les foncteurs $c$ et $d$ sont pleinement fid\`eles, et le
diagramme ci-dessus montre alors que $a$ est une \'equivalence, ce
qui ach\`eve la d\'emonstration.
\begin{corollaire}
\label{XIII.3.4}
Soient $k$ un corps de caract\'eristique $p\geqslant0$, $p'$
l'ensemble des nombres premiers distincts de~$p$, $f\colon X\to k$ un
morphisme coh\'erent. Supposons satisfaite l'une des deux
conditions suivantes:
\begin{enumerate}
\item[a)] $f$ est de type fini et les sch\'emas de type fini de
dimension $\leqslant\dim X$ sur une cl\^oture alg\'ebrique de~$k$
sont d\'esingularisables.
\item[b)] Les sch\'emas de type fini sur une cl\^oture
alg\'ebrique de~$k$ sont d\'esingularisables.
\end{enumerate}
Alors $f$ est universellement localement $1$-asph\'erique pour $p'$.
\end{corollaire}

Le cas a) r\'esulte de~\Ref{XIII.3.3}. Dans le cas b), la question
\'etant locale sur~$X$, on peut supposer $X$ affine; par passage
\`a la limite (SGA~4 XV~1.3), on se ram\`ene au cas o\`u $X$ est
de type fini sur~$k$.
\begin{corollaire}
\label{XIII.3.5}
Soient
\marginpar{414}
$S$ un sch\'ema irr\'eductible de point g\'en\'erique
$s$, $f\colon X\to S$ un morphisme de pr\'esentation
finie. Supposons que les sch\'emas de type fini de dimension
$\leqslant\dim X_s$ sur une cl\^oture alg\'ebrique $\overline{k}$
de~$\kres(s)$ soient fortement d\'esingularisables \textup{(SGA~5 I~3.1.5)}. Si
$\LL$ d\'esigne l'ensemble des nombres premiers distincts des
caract\'eristiques r\'esiduelles de~$S$, on peut trouver un ouvert
non vide~$S_1$ de~$S$ tel que, pour toute sp\'ecialisation
$\bar{s}_1\to\bar{s}_2$ de points g\'eom\'etriques de~$S_1$, le
morphisme de sp\'ecialisation \eqref{XIII.2.10}
$$
\pi_1^{\LL}(X_{\bar{s}_1})\to\pi_1^{\LL}(X_{\bar{s}_2})
$$
soit bijectif.
\end{corollaire}

D'apr\`es~\Ref{XIII.3.1} et~\Ref{XIII.3.3}, on peut, quitte \`a
restreindre $S$ \`a un ouvert non vide, supposer que $f$ est
localement $1$-asph\'erique pour $\LL$, et que, pour tout faisceau
de~$\LL$-groupes constant fini $F$ sur~$X$, $(F,f)$ est
cohomologiquement propre en dimension $\leqslant 1$. Il r\'esulte
alors de~\Ref{XIII.1.14} que, pour toute sp\'ecialisation
$\bar{s}_1\to\bar{s}_2$ de points g\'eom\'etriques de~$S_1$, le
morphisme de sp\'ecialisation
$$
\left(\R^1 f_* F\right)_{\bar{s}_2} \to \left(\R^1 f_*
F\right)_{\bar{s}_1}
$$
est bijectif. Le corollaire n'est autre que la traduction de ce qui
pr\'ec\`ede en termes de groupes fondamentaux.

\section{Suites exactes d'homotopie}\label{XIII.4}\setcounter{subsection}{-1}
\subsection{}
\label{XIII.4.0}
Soient $X$ et $S$ deux sch\'emas connexes, $f\colon X\to S$ un
morphisme, $a$ un point g\'eom\'etrique de~$X$, $\LL$ un ensemble
de nombres premiers. Soit $K$ le noyau de l'homomorphisme canonique
$\pi_1(X,a)\to\pi_1(S,a)$ et $N$ le plus petit pro-sous-groupe
distingu\'e de~$K$ tel que $K/N$ soit un pro-$\LL$-groupe
$K^{\LL}$. Alors $N$ est distingu\'e dans $\pi_1(X,a)$ et on note
$$
\pi_1'(X,a)
$$
\label{indnot:mj}\oldindexnot{$\pi_1'(X,a)$|hyperpage}%
le
\marginpar{415}
quotient de~$\pi_1(X,a)$ par $N$. Si $a$ est un point
g\'eom\'etrique d'une fibre g\'eom\'etrique~$X_{\bar{s}}$, les
morphismes canoniques
$$
\pi_1(X_{\bar{s}},a)\to \pi_1(X,a) \to
\pi_1(S,a)
$$
permettent de d\'efinir des morphismes canoniques
$$
\pi_1^{\LL}(X_{\bar{s}},a)\lto{u} \pi_1'(X,a) \lto{v} \pi_1(S,a)
$$
On a $vu=0$.

\begin{proposition}
\label{XIII.4.1}
Soient $S$ un sch\'ema connexe, $f\colon X\to S$ un morphisme
localement $0$-acyclique \textup{(SGA~4 XV~1.11)}; supposons de plus $f$
$0$-acyclique (ce qui, lorsque $f$ est coh\'erent, revient \`a
dire que les fibres g\'eom\'etriques de~$f$ sont connexes \textup{(SGA~4
XV~1.16)}). Soit $\LL$ un ensemble de nombres premiers. Si $S'$ est un
sch\'ema \'etale sur~$S$, on note $X'$, $f'$ les images inverses
de~$X$, $f$ sur~$S'$. Supposons que, pour tout rev\^etement
\'etale $S'$ de~$S$ et pour tout rev\^etement \'etale $E$ de
$X'$, quotient d'un rev\^etement galoisien de groupe un
$\LL$-groupe, $(E,f')$ soit cohomologiquement propre relativement
\`a~$S'$ en dimension $\leqslant 0$ et $f'_* E$
constructible. Alors, si $\bar{s}$ est un point g\'eom\'etrique de
$S$ et $a$ un point g\'eom\'etrique de la fibre $X_{\bar{s}}$, la
suite d'homomorphismes de groupes
\begin{equation*}
\label{eq:XIII.4.1.1}
\tag{\thesubsection.1}
\pi_1^{\LL}(X_{\bar{s}},a)\lto{u} \pi_1'(X,a)\lto{v}\pi_1(S,a)\to1
\end{equation*}
est exacte
\end{proposition}

Cet \'enonc\'e g\'en\'eralise X.\Ref{X.1.4}, dont on va copier
la d\'emonstration.

Montrons d'abord que $v$ est surjectif. Il suffit de montrer que, pour
tout rev\^etement \'etale connexe $S'$ de~$S$, $X'$ est aussi
connexe (V~\Ref{V.6.9}). Soit $C$ un ensemble ayant au moins deux
\'el\'ements. Il r\'esulte du fait que $f$ est $0$-acyclique que
le morphisme canonique
$$
\H^0(S',C_{S'})\to \H^0(X',C_{X'})
$$
est bijectif, donc que $X'$ est connexe, d'o\`u la surjectivit\'e
de~$v$.

Par
\marginpar{416}
d\'efinition de~$K^{\LL}$ \eqref{XIII.4.0}, on a la suite exacte
$$
1\to K^{\LL} \to \pi_1'(X,a)\to \pi_1(S,a)\to 1
$$
Soient $\widetilde{S}$ le rev\^etement universel de~$S$ et
$\widetilde{X}=\widetilde{S}\times_SX$; le groupe $K^{\LL}$ classe les
rev\^etements galoisiens $P$ de groupe un $\LL$-groupe, tels qu'il
existe un rev\^etement \'etale $S'$ de~$S$ et un rev\^etement
galoisien $Q$ de~$X'=X\times_SS'$ tels que l'on ait un isomorphisme
$P\simeq Q\times_{X'}\widetilde{X}$. Pour que la suite
\eqref{eq:XIII.4.1.1} soit exacte, il faut et il suffit que le
morphisme canonique
$$
\pi_1^{\LL}(X_{\bar{s}},a)\to K^{\LL}
$$
soit surjectif. D'apr\`es l'interpr\'etation de~$K^{\LL}$ cela
revient \`a dire que, pour tout rev\^etement \'etale $S'$ de~$S$
et pour tout rev\^etement galoisien $Q$ de~$X'$ de groupe un
$\LL$-groupe, tel que $P=Q\times_{X'}\widetilde{X}$ soit connexe,
alors $Q|X_{\bar{s}}$ est connexe. Montrons que cette derni\`ere
condition est satisfaite. Soient en effet $S'$ un rev\^etement
\'etale de~$S$, $Q$ un rev\^etement galoisien de~$X'$ de groupe un
$\LL$-groupe $F$, tel que $Q|X_{\bar{s}}$ soit disconnexe et montrons
que, quitte \`a remplacer $S'$ par un rev\^etement \'etale, $Q$
devient disconnexe. Il existe un sous-groupe $G$ de~$F$ distinct de
$F$ et un torseur $R$ sous $G|X_{\bar{s}}$ tel que $Q|X_{\bar{s}}$
s'obtienne par extension du groupe structural $G\to F$ \`a partir de
$R$. Le rev\^etement \'etale $E=Q/G$ de~$X'$ est tel que
$E|X_{\bar{s}}$ ait une section. D'apr\`es~\Ref{XIII.1.16} $f_*E$
est localement constant constructible et, quitte \`a remplacer $S'$
par un rev\^etement \'etale, on peut m\^eme supposer que $f_*E$
est constant. Comme $(E,f')$ est cohomologiquement propre relativement
\`a~$S'$ en dimension $\leqslant0$ et comme $\H^0(X_{\bar{s}},
E|X_{\bar{s}})$ est non vide, on voit que $E$ a une section. Mais ceci
prouve que $Q$ est disconnexe, ce qui ach\`eve la d\'emonstration.

On d\'eduit de~\Ref{XIII.4.1} le lemme suivant, qui sera utilis\'e
dans~\Ref{XIII.4.6}.

\begin{lemme}
\label{XIII.4.2}
Soient $S$ un sch\'ema connexe, $f\colon X\to S$ un morphisme
localement $0$-acyclique et $0$-acyclique, $\LL$ un ensemble de
nombres premiers. Supposons
\marginpar{417}
que, pour tout faisceau de~$\LL$-groupes constant fini $F$ sur~$X$,
$(F,f)$ soit cohomologiquement propre relativement \`a~$S$ en
dimension $\leqslant1$, et que, pour tout rev\^etement \'etale
$S'$ de~$S$ et pour tout rev\^etement \'etale $E$ de~$X'$,
quotient d'un rev\^etement galoisien de groupe un $\LL$-groupe,
$f'_*E$ soit constructible. Alors, si $\bar{s}$ est un point
g\'eom\'etrique de~$S$ et $a$ un point g\'eom\'etrique de la
fibre $X_{\bar{s}}$, la suite d'homomorphismes de groupes
$$
\pi_1^{\LL}(X_{\bar{s}},a)\to\pi_1'(X,a)\to\pi_1(S,a)\to1
$$
est exacte.
\end{lemme}

les hypoth\`eses de~\Ref{XIII.4.1} sont satisfaites. Il r\'esulte
en effet de~\Ref{XIII.1.13} 3) que, pour tout sch\'ema $S'$
\'etale sur~$S$ et pour tout rev\^etement \'etale $E$ de~$X'$,
quotient d'un rev\^etement galoisien de groupe un $\LL$-groupe,
$(E,f')$ est cohomologiquement propre relativement \`a~$S'$ en
dimension $\leqslant 0$.
\begin{proposition}
\label{XIII.4.3}
Soient $S$ un sch\'ema connexe, $\LL$ un ensemble de nombres
premiers, $f\colon X\to S$ un morphisme $0$-acyclique, localement
$1$-asph\'erique pour $\LL$ \textup{(SGA~4 XV~1.11)}, $g\colon S\to X$ une
section de~$f$. Soient $\bar{s}$ un point g\'eom\'etrique de~$S$,
$a$ un point g\'eom\'etrique de la fibre $X_{\bar{s}}$. On suppose
que, pour tout faisceau de
\ifthenelse{\boolean{orig}}
{$\LL$-groupe}
{$\LL$-groupes}
constant~$F$, $(F,f)$ est cohomologiquement propre en dimension
$\leqslant1$, que l'image directe par~$f$ du champ des torseurs sous
$F$ est un champ $1$-constructible \eqref{XIII.0} et que, pour tout
rev\^etement \'etale $E$ de~$X'$, quotient d'un rev\^etement
galoisien de groupe un $\LL$-groupe, $f'_*E$ est constructible. Alors
la suite d'homomorphismes de groupes
\begin{equation*}
\label{eq:XIII.4.3.1}
\tag{\thesubsection.1}
1\to\pi_1^{\LL}(X_{\bar{s}},a)\lto{u} \pi_1'(X,a)\lto{v}\pi_1(S,a)\to1
\end{equation*}
est exacte.
\end{proposition}

Compte tenu de~\Ref{XIII.4.2}, il suffit de montrer l'injectivit\'e
du morphisme $u$, \ie de prouver que, pour tout rev\^etement
principal $\overline{Z}$ de~$X_{\bar{s}}$
\marginpar{418}
de groupe un $\LL$-groupe~$C$, il existe un rev\^etement \'etale
$Z$ de~$X$ et un morphisme d'une composante connexe de~$Z|X_{\bar{s}}$
dans~$\overline{Z}$ (V~\Ref{V.6.8}). Soient donc $\overline{Z}$ un
rev\^etement principal de~$X_{\bar{s}}$ de groupe un $\LL$-groupe~$C$ et~$\bar{z}$ sa classe dans
$\H^1(X_{\bar{s}},C_{X_{\bar{s}}})$. D'apr\`es \Ref{XIII.1.5}~d), on
a un isomorphisme canonique
$$
\left(\R^1f_*C_X\right)_{\bar{s}}\isomto\H^1(X_{\bar{s}},C_{X_{\bar{s}}}),
$$
et d'apr\`es~\Ref{XIII.1.16}, $\R^1f_*C_X$ est un faisceau
localement constant constructible. On peut donc trouver un
rev\^etement \'etale $S'$ de~$S$ tel que $\R^1f_*C_X|S'$ soit
constant. Si $\bar{s}\to S'$ est un point g\'eom\'etrique
au-dessus du point g\'eom\'etrique $\bar{s}\to S$, il existe un
\'el\'ement~$z$ de~$\H^0(S', \R^1f_*C_X)$ dont l'image dans
$\H^1(X_{\bar{s}},C_{X_{\bar{s}}})$ est $\bar{z}$. D'apr\`es le
lemme~\Ref{XIII.4.3.1} ci-dessous, on peut trouver un rev\^etement
\'etale $S'_1$ de~$S'$ et un torseur $P$ sur~$X'_1$ de groupe~$C$
dont l'image dans $\H^0(S'_1,\R^1f_*C)$ soit \'egale \`a la
restriction de~$z$. Le torseur~$P$ est repr\'esentable par un
rev\^etement \'etale $Z$ de~$X'_1=X\times_SS'_1$ tel que
$Z\times_{X'_1}X_{\bar{s}}$ soit isomorphe \`a $\overline{Z}$. Si
l'on consid\`ere $Z$ comme un rev\^etement \'etale de~$X$, on a
alors un morphisme de~$Z\times_XX_{\bar{s}}$ dans $\overline{Z}$, ce
qui ach\`eve la d\'emonstration.

\begin{sublemme}
\label{XIII.4.3.1}
Soient $f\colon X\to S$ un morphisme $0$-acyclique et localement
$0$-acyclique, $g$ une section de~$f$. Soit $C$ un groupe constant
fini tel que $(C_X,f)$ soit cohomologiquement propre en dimension
$\leqslant0$ et que l'image directe par $f$ du champ des torseurs sous
$C_X$ soit constructible. Alors, pour toute section $z$ de
$\H^0(S,\R^1f_*C_X)$, on peut trouver un rev\^etement \'etale
$S_1$ de~$S$ et, si $X_1=X\times_SS_1$, un \'el\'ement de
$\H^1(X_1,C_{X_1})$, dont l'image par le morphisme canonique
$$
\H^1(X_1,C_{X_1})\to\H^0(S_1,\R^1f_*C_{X_1})
$$
soit \'egale \`a la restriction de~$z$ \`a
$\H^0(S_1,\R^1f_*C_{X_1})$.
\end{sublemme}

Pour tout sch\'ema $S'$ \'etale sur~$S$, on pose $X'=X\times_SS'$,
et on note $g'$ (\resp $F'$, etc.) l'image inverse de~$g$ (\resp $F$,
etc.) par le morphisme
\marginpar{419}
$S'\to S$. Le pr\'efaisceau $G$ sur~$S$ d\'efini par
\newlength{\tmpRacinet} \setlength{\tmpRacinet}{\textwidth}
\addtolength{\tmpRacinet}{-7em}
$$
G(S') = \left\{\text{\parbox{\tmpRacinet}{classes mod. isomorphisme de
		torseurs $P$ sur~$X'$ de groupe $C_{X'}$, munis d'un
		isomorphisme $g'^*P\isomto C_{S'}$ }} \right\}
$$
est alors un faisceau. Cela r\'esulte en effet par descente du fait
qu'un isomorphisme d'un torseur $P$ sur~$X'$ est bien
d\'etermin\'e par sa restriction \`a $g'(S')$. De plus, on a un
morphisme surjectif
$$
G\to \R^1f_*F.
$$
Soient $z$ un \'el\'ement de~$\H^0(S,\R^1f_*C_X)$ et $H$ le
sous-faisceau de~$G$ image r\'eciproque de~$z$. Il suffit de montrer
que $H$ est un faisceau localement constant constructible. Or, cette
propri\'et\'e \'etant locale sur~$S$, on peut supposer que $z$
provient d'un \'el\'ement de~$\H^1(X,C_X)$ repr\'esent\'e par
un torseur $P$ tel que $g^*P$ soit isomorphe \`a $C_S$. Se donner un
isomorphisme $i\colon g^*P\isomto C_S$ revient \`a se donner une
section globale de~$\SheafAut_{C_S}(g^*P)$ et deux isomorphismes $i$
et $i'$ d\'efinissent le m\^eme \'el\'ement de~$G(X)$ si et
seulement si $ii'^{-1}$ est l'image d'un \'el\'ement de
$\Aut_{C_X}(P)$. Si l'on consid\`ere l'injection canonique
$$
f_*\SheafAut_{C_X}(P)\to\SheafAut_{C_S}(g^*P)\simeq C,
$$
$H$ s'identifie donc au quotient de~$\SheafAut_{C_S}(g^*P)$ par
$f_*\SheafAut_{C_X}(P)$. D'apr\`es~\Ref{XIII.1.16}
$f_*\SheafAut_{C_X}(P)$ est localement constant; il en est donc de
m\^eme de~$H$, ce qui ach\`eve la d\'emonstration.

\begin{exemples}
\label{XIII.4.4}
Notons que, si $S$ est un sch\'ema connexe, les hypoth\`eses
de~\Ref{XIII.4.1} sont satisfaites lorsque le morphisme $f$ est propre
plat de pr\'esentation finie, \`a fibres g\'eom\'etriques
s\'eparables connexes, $\LL$ \'etant quelconque
(\cf X~\Ref{X.1.3}). Les hypoth\`eses de~\Ref{XIII.4.3} sont
satisfaites si de plus $f$ est lisse et a une section, si l'on
d\'esigne par $\LL$ l'ensemble des nombres premiers distincts des
caract\'eristiques r\'esiduelles de~$S$ (SGA~4 XV~2.1 et XVI~5.2).

Les hypoth\`eses de~\Ref{XIII.4.1} sont aussi satisfaites si $S$ est connexe,
si l'on a un sch\'ema~$Z$ propre de pr\'esentation finie, plat sur~$S$, \`a
\marginpar{420}
fibres g\'eom\'etriques s\'eparables connexes, tel que $X$ soit
le compl\'ementaire dans $Z$ d'un diviseur \`a croisements normaux
relativement \`a~$S$, $\LL$ \'etant l'ensemble des nombres
premiers distincts des caract\'eristiques r\'esiduelles de~$S$
\eqref{XIII.2.9}. Les hypoth\`eses de~\Ref{XIII.4.3} sont
satisfaites si de plus $f$ est lisse et a une section.
\end{exemples}

\subsection{}
\label{XIII.4.5}
Reprenons les notations et les hypoth\`eses de~\Ref{XIII.4.3}. Si
$\bar{s}$ est un point g\'eom\'etrique de~$S$ et $a=g(\bar{s})$,
la section $g$ permet de d\'efinir un morphisme
$$
w\colon \pi_1(S,a)\to \pi_1'(X,a),
$$
de sorte que $\pi_1'(X,a)$ s'identifie au produit semi-direct de
$\pi_1(S,a)$ par $\pi_1(X_{\bar{s}},a)$. Le groupe pro-fini
$\pi_1(S,a)$ op\`ere donc sur~$\pi_1(X_{\bar{s}}, a)$. Comme
$\pi_1^\LL(X_{\bar{s}},a)$ est limite projective stricte de groupes
invariants par l'action de~$\pi_1(S,a)$, la donn\'ee de
$\pi_1^\LL(X_{\bar{s}},a)$ muni de cette action est \'equivalente
\`a la donn\'ee d'un syst\`eme projectif strict de sch\'emas
en groupes \'etales finis sur~$S$ que l'on note
$$
\pi_1^\LL(X/S,g,\bar{s})\quad\text{ou
simplement}\quad\pi_1^\LL(X/S,g).
$$
\label{indnot:mk}\oldindexnot{$\pi_1^\LL(X/S,g,\bar{s})$ ou $\pi_1^\LL(X/S,g)$|hyperpage}%

On a alors les propri\'et\'es suivantes:
\setcounter{subsubsection}{0}

\subsubsection{}
\label{XIII.4.5.1}
Pour tout sch\'ema en groupes \'etale fini $G$ sur~$S$, dont les
fibres sont des $\LL$-groupes, l'ensemble $E$ des classes de torseurs
$P$ sous l'image inverse $G_X$ de~$G$ sur~$X$, munis d'un isomorphisme
$g^*P\isomto G$, est canoniquement isomorphe \`a l'ensemble
$$
\Hom_S(\pi_1^\LL(X/S,g,\bar{s}),G) \quad\text{mod. automorphismes
int\'erieurs de~$G$.}
$$

\subsubsection{}
\label{XIII.4.5.2}
Pour tout sch\'ema en groupes \'etale fini $G$ sur~$S$ dont les
fibres sont des $\LL$-groupes, le faisceau $\R^1f_*G_X$ est
canoniquement isomorphe au faisceau associ\'e au pr\'efaisceau
$$
S'\mto \Hom_{S'}(\pi_1^\LL(X/S,g,\bar{s}),G)\
\text{mod. automorphismes int\'erieurs de~$G$.}
$$
($S'$ d\'esigne un sch\'ema \'etale sur~$S$).

\subsubsection{}
\label{XIII.4.5.3}
Soient
\marginpar{421}
$S'$ un
\ifthenelse{\boolean{orig}}
{$S$-ch\'ema}
{$S$-sch\'ema}
connexe, $\bar{s}$ un point g\'eom\'etrique de~$S'$, $X'$, $g'$
les images inverses respectives de~$X,g$ sur~$S'$. Alors
$\pi_1^\LL(X'/S',g',\bar{s})$ est canoniquement isomorphe \`a
l'image inverse de~$\pi_1^\LL(X/S,g,\bar{s})$ sur~$S'$. Pour tout
point g\'eom\'etrique $\xi$ de~$S$, la fibre
$\pi_1^\LL(X/S,g,\bar{s})_\xi$ est isomorphe \`a $\pi_1^\LL(X_\xi)$.

La donn\'ee de~$G$ est en effet \'equivalente \`a la donn\'ee
d'un $\LL$-groupe abstrait $\mathbf{G}$
sur lequel op\`ere $\pi_1(S,a)$, d'o\`u une action de
$\pi_1'(X,a)$ sur~$\mathbf{G}$. L'isomorphisme d\'efini
dans~\Ref{XIII.4.5.1} s'obtient alors par restriction au sous-ensemble
$E$ \`a partir du morphisme canonique
\[
\H^1(\pi_1'(X,a),\mathbf{G})\to
\H^1(\pi_1^\LL(X_{\bar{s}},a),\mathbf{G}) =
\Hom(\pi_1^\LL(X_{\bar{s}},a),\mathbf{G})/\text{aut. int. }G,
\]
l'ensemble $E$ s'envoyant bijectivement sur le sous-ensemble des
morphismes de $\pi_1^\LL(X_{\bar{s}},a)$ dans $\mathbf{G}$ qui sont
compatibles avec l'action de
$\pi_1(S,a)$. L'assertion~\Ref{XIII.4.5.3} r\'esulte alors de la
d\'efinition de~$\pi_1^\LL(X/S,g,\bar{s})$ compte-tenu de la suite
exacte d'homotopie \eqref{eq:XIII.4.3.1} et~\Ref{XIII.4.5.2} se
d\'eduit de \Ref{XIII.4.5.1} et~\Ref{XIII.4.5.3}.

\begin{proposition}[Formule de K\"unneth]
\label{XIII.4.6}
\index{Kunneth (formule de)@K\"unneth (formule de)|hyperpage}%
Soient $k$ un corps s\'eparablement clos de caract\'eristique
$p\geq 0$, $X$ et $Y$ deux $k$-sch\'emas connexes, $a$ un point
g\'eom\'etrique de~$X$, $b$ un point g\'eom\'etrique de~$Y$,
$c$ un point g\'eom\'etrique de~$X\times_k Y$ au-dessus de~$a$ et
$b$. On suppose satisfaite l'une des deux conditions suivantes:
\begin{enumerate}
\item[a)] $X$ est de type fini sur~$k$ et les sch\'emas de type fini
sur une cl\^oture alg\'ebrique $\bar{k}$ de~$k$, de dimension
$\leq\dim X$, sont fortement d\'esingularisables. \textup{(SGA~5 I~3.1.5)}.
\item[b)] $X$ est quasi-compact et quasi-s\'epar\'e et tout
sch\'ema de type fini sur~$\bar{k}$ est fortement
d\'esingularisable.
\end{enumerate}

Alors, si $p'$ est l'ensemble des nombres premiers distincts de~$p$,
le morphisme
\begin{equation*}
\label{eq:XIII.4.6.0} \tag{\thesubsection.0} {}
\pi_1^{p'}(X\times_k Y,c)\to\pi_1^{p'}(X,a)\times\pi_1^{p'}(Y,b),
\end{equation*}
d\'eduit
\marginpar{422}
des homomorphismes sur les groupes fondamentaux associ\'es aux
projections
$$
X\times_k Y\to X \quad \text{et} \quad X\times_k Y\to Y,
$$
est un isomorphisme.
\end{proposition}

On peut supposer $k$ alg\'ebriquement clos et $X$ r\'eduit
(SGA~4 VIII~1.1). Soient $Z=X\times_k Y$, $g\colon X\to k$ et $f\colon
Z\to Y$ les morphismes canoniques. Le morphisme $g$, donc aussi $f$,
est universellement localement $1$-asph\'erique pour $p'$
d'apr\`es~\Ref{XIII.3.5}. Comme $X$ est connexe, $f$ est
$0$-acyclique (SGA~4 XV~1.16). D'autre part il r\'esulte
de~\Ref{XIII.3.2} que, pour tout $p'$-groupe fini $C$, $(C_X,g)$ est
cohomologiquement propre relativement \`a $k$ en dimension $\leq
1$. Il en r\'esulte que $(C_Z,f)$ est cohomologiquement propre
relativement \`a $Y$ en dimension $\leq 1$ (\Ref{XIII.1.5}~c)) et
que $f_*C_Z$ et $\R^1f_*C_Z$ sont des faisceaux constants. Par suite
$g$ satisfait \`a toutes les hypoth\`eses de~\Ref{XIII.4.2}. On a
donc la suite exacte
$$
\pi_1^{p'}(X_b,c)\to\pi_1^{p'}(Z,c)\to\pi_1^{p'}(Y,b)\to 1.
$$
De plus le morphisme compos\'e
$$
\pi_1^{p'}(X_b,c)\to\pi_1^{p'}(Z,c)\to\pi_1^{p'}(Y,b)
$$
est un isomorphisme et l'on a donc la suite exacte
$$
1\to\pi_1^{p'}(X,a)\to\pi_1^{p'}(Z,c)\to\pi_1^{p'}(Y,b)\to 1.
$$
D'autre part le morphisme \eqref{eq:XIII.4.6.0} permet de d\'efinir
un morphisme de cette suite exacte dans la suite exacte
$$
1\to\pi_1^{p'}(X,a)\to
\pi_1^{p'}(X,a)\times\pi_1^{p'}(Y,b)\times\pi_1^{p'}(Y,b)\to 1,
$$
et il en r\'esulte que le morphisme \eqref{eq:XIII.4.6.0} est un
isomorphisme.

\subsection{}
\label{XIII.4.7} Soit
$$
\xymatrix@C=1.2cm{
X \ar[d]_-f&\ar[l] X_U \ar[d]_-{f_U}& \ar[l]X_{\bar{s}}\ar[d] \\
S&\ar[l] U &\ar[l]\bar{s}
}
$$
un
\marginpar{423}
diagramme dont les carr\'es sont cart\'esiens, o\`u $S$ est un
sch\'ema connexe par arcs (SGA~4 IX~2.12), $U$ un ouvert connexe de
$S$, $\bar{s}$ un point g\'eom\'etrique de~$U$. Soient $a$ un
point g\'eom\'etrique de~$X_{\bar{s}}$, $\LL$ un ensemble de
nombres premiers. Soit $g$ une section de~$f$ et supposons que les
conditions suivantes soient satisfaites:
\begin{enumerate}
\item[a)] Le morphisme $f$ est $0$-acyclique, localement $0$-acyclique,
et, pour tout rev\^etement \'etale $S'$ de~$S$ et tout
rev\^etement \'etale $E$ de~$X\times_S S'$ quotient d'un
rev\^etement galoisien de groupe un $\LL$-groupe, $(E,f_{(S')})$ est
cohomologiquement propre relativement \`a~$S'$ en dimension $\leq
0$.
\item[b)] Le morphisme $f_U$ est localement $1$-asph\'erique pour
$\LL$, et, pour tout faisceau de~$\LL$-groupe constant fini $F$ sur~$X_U$, $(F,f_U)$ est cohomologiquement propre en dimension $\leq 1$ et
les fibres de~$\R^1 f_{U*}F$ sont finies.
\end{enumerate}

On d\'eduit alors de~\Ref{XIII.4.1} et~\Ref{XIII.4.3} le diagramme
commutatif suivant dont les lignes sont exactes:
\begin{equation*}
\label{eq:XIII.4.7.0} \tag{\thesubsection.0} {}
\begin{array}{c}
\xymatrix{
1 \ar[r]& \pi_1^{\LL}(X_{\bar{s}},a) \ar[r]\ar@{=}[d]& \pi'_1(X_U,a) \ar[r]\ar[d]& \pi_1(U,a)\ar[r]\ar[d]& 1 \\
&\pi_1^{\LL}(X_{\bar{s}},a) \ar[r]&\pi'_1(X,a) \ar[r]& \pi_1(S,a) \ar[r]& 1
}
\end{array}
\end{equation*}

Gr\^ace \`a la section $g$, on a des morphismes
\ifthenelse{\boolean{orig}}
{$\pi_1(U,a)\to\pi'_1(X_U,a)$ $\pi_1(S,a)\to\pi'_1(X,a)$;}
{$$\pi_1(U,a)\to\pi'_1(X_U,a),\quad \pi_1(S,a)\to\pi'_1(X,a)\,\;$$}
on en
d\'eduit un morphisme de la somme amalgam\'ee de~$\pi_1(S,a)$ et
$\pi'_1(X_U,a)$ au-dessus de~$\pi_1(U,a)$ dans $\pi'_1(X,a)$:
\begin{equation*}
\label{eq:XIII.4.7.1} \tag{\thesubsection.1} {}
\varphi\colon\pi=\pi_1(S,a)\coprod_{\pi_1(U,a)}\pi'_1(X_U,a)\to\pi'_1(X,a).
\end{equation*}
Supposons satisfaite la condition suivante:
\begin{enumerate}
\item[c)] Si $T=S-U$, on a $\prof \et_T(S)\geq 2$ (SGA~2 XIV~1.1).
\end{enumerate}
Alors le foncteur qui, \`a un rev\^etement \'etale de~$S$, fait
correspondre sa restriction \`a $U$ est pleinement fid\`ele
(SGA~2 XVI~1.4). Il en r\'esulte que le morphisme
$\pi_1(U,a)\to\pi_1(S,a)$ est surjectif (V~\Ref{V.6.9}) et l'on
d\'eduit du
\marginpar{424}
diagramme \eqref{eq:XIII.4.7.0} qu'il en est de m\^eme du morphisme
$\pi'_1(X_U,a)\to\pi'_1(X,a)$; a fortiori $\varphi$ est un
\'epimorphisme. Soit
\begin{equation*}
\label{eq:XIII.4.7.2} \tag{\thesubsection.2} {}
K=\Ker(\pi_1(U,a)\to\pi_1(S,a)).
\end{equation*}
Le groupe $\pi$ de \eqref{eq:XIII.4.7.1} s'identifie au quotient de
$\pi'_1(X_U,a)$ par le sous-groupe invariant ferm\'e engendr\'e
par l'image $L$ de~$K$ dans $\pi'_1(X_U,a)$. Consid\'erons
$\pi'_1(X_U,a)$ comme produit semi-direct de~$\pi_1(U,a)$ par
$\pi_1^\LL (X_{\bar{s}},a)$. Le groupe $K$ op\`ere alors par
automorphismes int\'erieurs sur~$\pi_1^\LL (X_{\bar{s}},a)$, et le
quotient $\pi=\pi'_1(X_U,a)/L$ s'identifie au produit semi-direct de
$\pi_1(U,a)/K=\pi_1(S,a)$ par le groupe $\pi_1^\LL(X_{\bar{s}},a)_K$
des coinvariants de~$\pi_1^{\LL}(X_{\bar{s}},a)$ sous~$K$. On a
finalement un \'epimorphisme
\begin{equation*}
\label{eq:XIII.4.7.3} \tag{\thesubsection.3} {}
\varphi\colon\pi=\pi_1^\LL(X_{\bar{s}},a)_K\cdot\pi_1(S,a)\to\pi'_1(X,a).
\end{equation*}
\label{indnot:ml}\oldindexnot{$\pi_1^\LL(X_{\bar{s}},a)_K$|hyperpage}%

La proposition qui suit donne des conditions sous-lesquelles le
morphisme $\varphi$ est un isomorphisme.

\setcounter{subsubsection}{3}

\begin{subproposition}
\label{XIII.4.7.4} Les notations sont celles de~\Ref{XIII.4.7}. On
suppose que, en plus des conditions \textup{a), b), c)} les conditions
suivantes soient satisfaites:
\begin{enumerate}
\item[d)] Pour tout point $t$ de~$T=S-U$, le morphisme $f$ est
localement $1$-asph\'erique pour $\LL$ en $g(t)$.
\item[e)] Pour tout point $t\in T$, toute composante irr\'eductible
de la fibre $X_t$ contient $g(t)$ et, pour tout point $x$ de
$X_t-\{g(t)\}$ qui n'est pas maximal, on a
$$
\prof\hop_x(X)\geq 3 \quad\textup{(SGA~2 XIV~1.2)}
$$
et l'anneau $\cal{O}_{X,x}$ est noeth\'erien.
\end{enumerate}

Alors le morphisme \eqref{eq:XIII.4.7.3} est un isomorphisme.
\end{subproposition}

Comme on l'a dit pr\'ec\'edemment, le groupe $\pi$ s'identifie au
quotient de~$\pi'_1(X_U,a)$ par le sous-groupe invariant ferm\'e $L$
engendr\'e par l'image de~$K$ \eqref{eq:XIII.4.7.2} dans
$\pi'_1(X_U,a)$. Cela revient \`a dire que $\pi$ classe les
\marginpar{425}
rev\^etements principaux $Z$ de~$X_U$ tels que $g^{-1}_U(Z)$ se
prolonge en un rev\^etement \'etale de~$S$, et qui induisent sur~$X_{\bar{s}}$ un rev\^etement qui s'obtient par extension du groupe
structural \`a partir d'un rev\^etement principal de groupe un
$\LL$-groupe. Pour prouver que $\varphi$ est un isomorphisme, il
suffit de montrer qu'un tel rev\^etement $Z$ se prolonge \`a $X$
tout entier.

Montrons d'abord que $Z$ se prolonge en un rev\^etement \'etale
d'un ouvert contenant $X_U$ et $g(S)$. Soit $W$ un sch\'ema
\'etale sur~$X$ dont l'image contient $X_U$ et $g(S)$, et posons
$W_U=W\times_S U$. Du fait que le morphisme $W\to S$ est $0$-acyclique
et de la relation $\mathrm{prof}.\,\et_T S\geq 2$, r\'esulte que l'on a
$$
\mathrm{prof}.\,\et_{(W-W_U)} W\geq 2
$$
(SGA~2 XIV~1.13); par suite, si $Z|W_U$ se prolonge en un
rev\^etement \'etale de~$W$, ce prolongement est unique \`a
isomorphisme unique pr\`es. Il en r\'esulte que le probl\`eme de
prolonger $Z$ \`a un voisinage de~$g(S)\cup X_U$ est local pour la
topologie \'etale au voisinage des points de~$g(T)$. Si $t$ est un
point de~$T$, on pose $x=g(t)$, et on note $\bar{X}$ (\resp $\bar{S}$)
le localis\'e strict de~$X$ en $\bar{x}$ (\resp de~$S$ en
$\bar{t}$), $\bar{U}=U\times_S \bar{S}$, $\bar{X}_U=\bar{X}\times_X
X_U$, $\bar{g}\colon\bar{S}\to\bar{X}$ le morphisme d\'eduit de
$g$. Il suffit de montrer que, pour tout point $t$ de~$T$, l'image
inverse $\bar{Z}$ de~$Z$ sur~$\bar{X}_U$ se prolonge \`a $\bar{X}$
ou, ce qui revient au m\^eme, est triviale. Or, par d\'efinition
de~$Z$, l'image inverse de~$\bar{Z}$ sur~$\bar{U}$ est
triviale. Pour prouver que $\bar{Z}$ est trivial, il suffit de montrer
qu'il est de la forme $\bar{f}^*_U E$, o\`u $E$ est un
rev\^etement principal de~$\bar{U}$; on aura en effet
$E\simeq\bar{g}^*_U\bar{f}^*_U E \simeq\bar{g}^*_U\bar{Z}$, d'o\`u
le r\'esultat puisque $\bar{g}^*_U\bar{Z}$ est trivial. Or le
morphisme $\bar{f}_U$ \'etant $0$-acyclique et localement $0$-acyclique,
il suffit, pour prouver que $\bar{Z}$ provient de~$\bar{U}$, de
montrer que, pour tout point g\'eom\'etrique alg\'ebrique sur un
point de~$\bar{U}$, que l'on peut supposer \^etre le point
$\bar{s}$, $\bar{Z}|X_{\bar{s}}$ est trivial (SGA~4 XV~1.15). Mais
$\bar{Z}|X_{\bar{s}}$ \'etant obtenu par extension du groupe
structural \`a partir d'un rev\^etement principal de groupe un
$\LL$-groupe, cela r\'esulte du fait que le morphisme $\bar{f}$ est
$1$-asph\'erique pour $\LL$.

On
\marginpar{426}
a donc d\'emontr\'e qu'il existe un voisinage ouvert $V$ de
$g(S)\cup X_U$ tel que $Z$ se prolonge en un rev\^etement \'etale
$Z_V$ de~$V$. Montrons que $Z_V$ se prolonge \`a $X$ tout entier.
Il suffit de voir que, pour tout point $x$ de~$X-V$, on a
$$
\prof\hop_x X\geq 3.
$$
Or cela r\'esulte de l'hypoth\`ese e) et du fait qu'un point $x$
de~$X-V$ ne peut \^etre maximal dans sa fibre $X_t$ car, toute
composante irr\'eductible de~$X_t$ contenant $g(t)$, tout point
maximal de~$X_t$ appartient \`a $V$.

\begin{corollaire} \label{XIII.4.8}
Les hypoth\`eses sont celles de~\Ref{XIII.4.7.4} mais on suppose de
plus que l'on a $\pi_1(S,a)=1$. Alors on a un isomorphisme
\begin{equation*}
\label{eq:XIII.4.8.*} \tag{$*$}
{}\pi_1^\LL(X,a)\isomto\pi_1^\LL(X_{\bar{s}},a)_K.
\end{equation*}
En particulier, si $\pi_1^\LL(X_{\bar{s}},a)$ est topologiquement de
pr\'esentation finie et si $K$ op\`ere sur~$\pi_1^\LL(X_{\bar{s}},a)$ par l'interm\'ediaire d'un groupe de type
fini, alors $\pi_1^\LL(X,a)$ est topologiquement de pr\'esentation
finie.
\end{corollaire}

L'isomorphisme \eqref{eq:XIII.4.8.*} a \'et\'e d\'emontr\'e
dans~\Ref{XIII.4.7.4}. Supposons $\pi_1^\LL(X_{\bar{s}},a)$ quotient
du pro-$\LL$-groupe libre \`a $n$ g\'en\'erateurs
$L(x_1,\dots,x_n)$ par le sous-groupe invariant ferm\'e engendr\'e
par les \'el\'ements $y_1,\dots,y_p$ de~$L(x_1,\dots,x_n)$, et
supposons que $K$ agisse par l'interm\'ediaire d'un groupe
engendr\'e par des \'el\'ements $k_1,\dots,k_q$. Si pour tout
$i\in[1,n]$, $j\in[1,q]$, on note $z_{ij}$ un \'el\'ement de
$L(x_1,\dots,x_n)$ relevant l'\'el\'ement $(k_j\cdot
x_i)x_i^{-1}$, alors $\pi_1^\LL(X_{\bar{s}},a)_K$ est quotient de
$L(x_1,\dots,x_n)$ par le sous-groupe invariant ferm\'e engendr\'e
par les \'el\'ements $(y_i)_{i\in[1,p]}$,
$(z_{ij})_{i\in[1,n],j\in[1,q]}$.

\ifthenelse{\boolean{orig}}
{\setcounter{subsection}{5}}
{}

\begin{remarques}
\label{rem:XIII.4.6}
a) Les conditions a) \`a e) de~\Ref{XIII.4.7} sont
satisfaites lorsque $S$ est un sch\'ema normal connexe, $U$ un
ouvert dense r\'etrocompact de~$S$, $f$ un morphisme propre de
pr\'esentation finie, \`a fibres g\'eom\'etriquement connexes
et irr\'eductibles en tout point $t$ de~$T$, $f$ \'etant de plus
s\'eparable,
\marginpar{427}
lisse aux points de~$X_U\cup g(T)$, $\LL$ \'etant l'ensemble des
nombres premiers distincts des caract\'eristiques r\'esiduelles de
$S$, et $X$ \'etant r\'egulier en tout point de~$X_t$. La
condition a) r\'esulte en effet de SGA~4 XV~4.1 et 1.4; les
conditions b) et d) r\'esultent de~\Ref{XIII.1.4} et SGA~4 XV~2.1 et
XVI~5.2. Enfin e) r\'esulte de SGA~XIV~1.11.
\par\smallskip
b) Le corollaire~\Ref{XIII.4.5} s'applique pour calculer le
groupe fondamental $\pi_1^{p'}(X)$ d'une surface $X$ propre et lisse
sur un corps s\'eparablement clos $k$ de caract\'eristique $p$
($p'$~d\'esignant l'ensemble des nombres premiers distincts de
$p$). La m\'ethode nous a \'et\'e communiqu\'ee par
J.P\ptbl \textsc{Murre}; elle consiste \`a se ramener, en faisant \'eclater
$X$, au cas o\`u l'on a une fibration $X\to\PP^1_k$ et un ouvert $U$
de~$\PP^1_k$ satisfaisant aux hypoth\`eses de~\Ref{XIII.4.7} (voir
SGA~7 pour plus de d\'etails). La m\^eme m\'ethode peut \^etre
utilis\'ee plus g\'en\'eralement (\loccit) pour prouver que, si
$X$ est un $k$-sch\'ema connexe de type fini, et, si les sch\'emas
de type fini de dimension $\leq\dim X$ sur une cl\^oture
alg\'ebrique de~$k$ sont fortement d\'esingularisables
(SGA~5 I~3.1.5), alors $\pi_1^{p'}(X)$ est topologiquement de
pr\'esentation finie.
\end{remarques}

\section{Appendice I: Variations sur le lemme d'Abhyankar}
\label{XIII.5}
\index{lemme d'Abhyankar|hyperpage}%
Cet appendice contient diff\'erentes variantes du lemme d'Abhyankar.

\begin{proposition}
\label{XIII.5.1} Soient $X=\Spec A$ un sch\'ema local r\'egulier,
$D=\sum_{1\leq i\leq r} \divisor f_i$ un diviseur \`a croisements
normaux, o\`u les $f_i$ sont des \'el\'ements de l'id\'eal
maximal de~$A$ faisant partie d'un syst\`eme r\'egulier de
param\`etres. Soient $n_i$, $1\leq i\leq r$, des entiers $\geq 0$
et posons
$$
X'=X[T_1,\dots,T_r]/(T_1^{n_1}-f_1,\dots,T_r^{n_r}-f_r),
$$
$U'=U\times_X X'$. Alors $X'$ est r\'egulier et $U'$ est le
compl\'ementaire dans $X'$ du diviseur \`a croisements normaux
$\sum_{1\leq i\leq r}\divisor T_i$. Si les entiers $n_i$ sont premiers
\marginpar{428}
\`a la caract\'eristique r\'esiduelle $p$ de~$X$, $U'$ est un
rev\^etement \'etale connexe de~$U$, mod\'er\'ement
ramifi\'e relativement \`a $D$ \textup{(\Ref{XIII.2.3}~c))}.
\end{proposition}

En effet $X'$ est le spectre d'un anneau local $A'$ dont l'id\'eal
maximal est engendr\'e par $T_1,\dots,T_r$.
\ifthenelse{\boolean{orig}}
{\ignorespaces}
{(On peut supposer $r=\dim(A)$.)}
Comme $A'$ est fini et plat sur~$A$, donc de dimension $r$, $A'$ est
r\'egulier (EGA~$0_{\textup{IV}}$~17.1.1) et les $T_i$ forment un syst\`eme
r\'egulier de param\`etres de~$A'$. Supposons les $n_i$ premiers
\`a $p$. Comme tous les $f_i$ sont inversibles sur~$U$, le fait que
$U'$ soit \'etale sur~$U$ r\'esulte de I~\Ref{I.7.4}. De plus $U'$
est mod\'er\'ement ramifi\'e relativement \`a $D$; soit en
effet $x_i$ le point g\'en\'erique de~$V(f_i)$, $\bar{R}$ le
localis\'e strict de~$R=\cal{O}_{X,x_i}$, $\bar{K}$ le corps des
fractions de~$\bar{R}$. Alors la $\bar{K}$-alg\`ebre qui
repr\'esente $U'|\bar{K}$ s'obtient \`a partir du corps
$\bar{K}[T_i]/(T_i^{n_i}-f_i)$ en faisant une extension non
ramifi\'ee; elle est donc mod\'er\'ement ramifi\'ee sur~$\bar{R}$.

\begin{proposition}[Lemme d'Abhyankar absolu]
\label{XIII.5.2}
\index{Abhyankar absolu (lemme d')|hyperpage}%
Soit $X$ un sch\'ema local r\'egulier,
$$
D=\sum_{1\leq i\leq r} \divisor f_i
$$
un diviseur \`a croisements normaux comme dans~\Ref{XIII.5.1},
$Y=\Supp D$, $U=X-Y$. Soit~$V$ un rev\^etement \'etale de~$U$
mod\'er\'ement ramifi\'e relativement \`a $D$. Si $x_i$ est le
point g\'en\'erique du ferm\'e $V(f_i)$, $\cal{O}_{X,x_i}$ est
un anneau de valuation discr\`ete de corps des fractions $K_i$, et
on a $V|K_i=\Spec(\prod_{j\in J_i}L_j)$, o\`u les $L_j$ sont des
extensions finies s\'eparables de~$K_i$; notons $n_j$ l'ordre du
groupe d'inertie d'une extension galoisienne engendr\'ee par $L_j$
et $n_i$ le p.p.c.m. des $n_j$ quand $j$ parcourt $J_i$. Si l'on pose
$$
X'=X[T_1,\dots,T_r]/(T_1^{n_1}-f_1,\dots,T_r^{n_r}-f_r),
$$
$U'=U_{(X')}$, $V'=V_{(X')}$, etc., le rev\^etement \'etale $V'$
de~$U'$ se prolonge de mani\`ere unique \`a isomorphisme unique
pr\`es en un rev\^etement \'etale de~$X'$, et les $n_i$ sont
premiers \`a la caract\'eristique r\'esiduelle $p$ de~$X$.
\end{proposition}

L'unicit\'e r\'esulte du fait que $X'$ est normal
\eqref{XIII.5.1}; en effet
\marginpar{429}
un rev\^etement \'etale de~$X'$ qui prolonge $V'$ est isomorphe au
normalis\'e de~$X'$ dans la fibre de~$V'$ au point g\'en\'erique
de~$X'$ \eqref{I.10.2}. Si $\bar{x}'$ est un point g\'eom\'etrique
de~$Y'$, on note $\bar{X}'$ le localis\'e strict de~$X'$ en
$\bar{x}'$, $\bar{V}'=V'_{(\bar{X}')}$, etc. Par descente, compte tenu
de l'unicit\'e, il suffit de montrer que, pour tout point
g\'eom\'etrique $\bar{x}'$ de~$Y'$, le rev\^etement \'etale
$\bar{V}'$ de~$\bar{U}'$ se prolonge \`a $\bar{X}'$. \'Etant donn\'e
qu'un rev\^etement \'etale d'un ouvert du sch\'ema r\'egulier~$X'$ qui contient tous les points $x'$ tels que l'on ait $\dim
\cal{O}_{X',x'}\leq 1$ se prolonge \`a $X$ tout entier
(SGA~2 XIV~1.11), on peut m\^eme se borner aux points $\bar{x}'$ qui
se projettent sur un point maximal de~$Y'$. Or, en un tel point
$\bar{x}'$, le fait que $\bar{V}'$ se prolonge en un rev\^etement
\'etale de~$\bar{X}'$ r\'esulte de~\Ref{X.3.6}.

Montrons que les $n_i$ sont premiers \`a $p$. En effet, s'il n'en
\'etait pas ainsi, on aurait par exemple $p|n_1$. Quitte \`a
remplacer $X$ par
\ifthenelse{\boolean{orig}}
{$X[T_1,\dots,T_r]/(T_1^{n_1/p}-f_1,T_2^{n_2}-f_2,\dots,T_r^{n_r}-f_r)$,}
{$$X[T_1,\dots,T_r]/(T_1^{n_1/p}-f_1,T_2^{n_2}-f_2,\dots,T_r^{n_r}-f_r),$$}
on se ram\`ene au cas o\`u l'on~a $X'=X[T_1]/(T_1^p-f_1)$. Il suffit
de montrer que $V$ se prolonge en un rev\^etement \'etale de~$X$,
car on aura alors $n_1=1$ contrairement \`a l'hypoth\`ese. On peut
supposer pour cela que $X$ est strictement local. Soit $Z$ le
sous-sch\'ema ferm\'e de~$X$ d'\'equation $p=0$ et $Z_1=Z\cap
X_{f_1}$; $Z_1$ est un ouvert non vide de~$Z$. D'apr\`es ce qui
pr\'ec\`ede le rev\^etement \'etale $V'$ de~$U'$ se prolonge
en un rev\^etement \'etale $W'$ de~$X'$. Soient $W''_1$ et $W''_2$
les images inverses de~$W'$ par les deux projections $X''=X'\times_X
X'\rightrightarrows X'$, et montrons que l'isomorphisme de descente
$u\colon W''_1|U''\to W''_2|U''$ se prolonge en un $X''$-morphisme
$W''_1\to W''_2$ qui sera n\'ecessairement une donn\'ee de
descente sur~$W'$ relativement \`a $X'\to X$. Soit $Z''$
(\resp $Z''_1$) l'image inverse de~$Z$ (\resp $Z_1$) dans $X''$. Comme
le morphisme $Z''\to Z$ est radiciel, il existe un isomorphisme
$v\colon W''_1|Z''\to W''_2|Z''$ qui prolonge l'isomorphisme
$u|Z''_1$. Mais, comme $X$ est hens\'elien, on a une bijection
$$
\Hom_{X''}(W''_1,W''_2)\simeq\Hom_{Z''}(W''_1|Z'',W_2''|Z''),
$$
d'o\`u un morphisme $w\colon W''_1\to W''_2$ relevant $v$. Le
sous-sch\'ema \`a la fois ouvert
\marginpar{430}
et ferm\'e de~$X''$ au-dessus duquel $u$ et $w$ co\"incident
contient $Z''_1$, donc \'egal \`a~$X''$, d'o\`u le fait que $V$
se prolonge \`a~$X$.

\setcounter{subsection}{3} \setcounter{subsubsection}{-1}

\subsubsection{}
\label{XIII.5.3.0} Reprenons les hypoth\`eses et les notations
de~\Ref{XIII.5.2} en supposant de plus
\ifthenelse{\boolean{orig}}
{$S$}
{$X$}
strictement local. Alors
il r\'esulte de \loccit que tout rev\^etement \'etale connexe
de~$U$ mod\'er\'ement ramifi\'e relativement \`a $D$ est
quotient d'un rev\^etement mod\'er\'ement ramifi\'e de la
forme
$$
U'=U[T_1,\dots,T_r]/(T_1^{n_1}-f_1,\dots,T_r^{n_r}-f_r),
$$
o\`u les $n_i$ sont des entiers premiers \`a $p$. Soit
$\bbmu_n$ le groupe des racines $n$-i\`emes de l'unit\'e de
$U$. Le groupe des $U$-automorphismes de~$U'$ n'est autre que le
groupe $\bbmu_{n_1}\times\dots\times\bbmu_{n_r}$, une racine
$n_i$-i\`eme de l'unit\'e $\xi_i$ op\'erant sur~$U'$ en
transformant $T_i$ en $\xi_i T_i$. On a donc l'\'enonc\'e suivant:

\setcounter{subsection}{2}

\begin{corollaire}
\label{XIII.5.3} Soit $X$ un sch\'ema strictement local r\'egulier
de caract\'eristique r\'esiduelle $p\geq 0$, $D=\sum_{1\leq i\leq
n} \divisor f_i$ un diviseur \`a croisements normaux sur~$X$,
$U=X-\Supp D$. Posons
$$
\tilde U=\varprojlim_{(n_i)}
U[T_1,\dots,T_r]/(T_1^{n_1}-f_1,\dots,T_r^{n_r}-f_r),
$$
la limite projective \'etant prise suivant l'ensemble ordonn\'e
filtrant (pour la relation de divisibilit\'e) des familles d'entiers
$n_i>0$, premiers \`a $p$. Alors $\tilde U$ est un rev\^etement
universel mod\'er\'ement ramifi\'e de~$U$. Par suite le groupe
fondamental mod\'er\'ement ramifi\'e de~$U$ est
$$
\pi_1^{\tame}(U)\simeq \prod_{\ell\neq p} \ZZ_\ell [1]^r
\qquad \text{(isomorphisme canonique)}
$$
o\`u l'on a pos\'e
$\ZZ_\ell[1]=\varprojlim_{n>0}\bbmu_{\ell^n}$. Le groupe
$\pi_1^{\tame}(U)$ est non canoniquement isomorphe \`a $\prod_{\ell\neq p}
\ZZ_\ell^r$.
\label{indnot:mm}\oldindexnot{$\ZZ_\ell[1]$|hyperpage}%
\end{corollaire}

\begin{proposition}
\label{XIII.5.4} Soient $f\colon X\to S$ un morphisme de
sch\'emas, $D=\sum_{1\leq i\leq r} \divisor f_i$ un diviseur \`a
croisements normaux relativement \`a~$S$ \eqref{XIII.2.1}, o\`u,
\marginpar{431}
pour chaque point $x$ de~$Y=\Supp D$, si $I(x)\subset[1,r]$ est
l'ensemble des $i$ tels que l'on ait $f_i(x)=0$, le sous-sch\'ema
$V((f_i)_{i\in I(x)})$ est lisse sur~$S$ de codimension
\ifthenelse{\boolean{orig}}
{$\mathrm{card}.\,I(x)$}
{$\card I(x)$}
dans $X$. Soit $U=X-Y$. Soient $x$ un point de
$Y$, $X_1=\Spec\cal{O}_{X,x}$, $U_1=U\times_X X_1$, $n_i$, $i\in I(x)$
des entiers et
$$
X'=X[T_i]_{i\in I(x)}/(T_i^{n_i}-f_i).
$$
Alors, si $x'$ est le point de~$X'$ au-dessus de~$x$, $X'$ est lisse
sur~$S$ en $x'$. Si les entiers~$n_i$ sont premiers \`a la
caract\'eristique $p$ de~$\kres(x)$, $U'_1=U_1\times_X X'$ est un
rev\^etement \'etale connexe de~$U_1$ mod\'er\'ement
ramifi\'e sur~$X_1$ relativement \`a~$S_1$ \eqref{XIII.2.1.1}.
\end{proposition}

Si $s=f(x)$, la fibre g\'eom\'etrique $X'_{\bar{s}}$ est
r\'eguli\`ere au point $x'$ \eqref{XIII.5.1}; comme $X'$ est plat
sur~$S$ au voisinage de~$Y$, cela prouve que $X'$ est lisse sur~$S$ en
$x'$ (EGA~IV~12.1.6). Si les entiers $n_i$ sont premiers \`a $p$,
$U'_1$ est un rev\^etement \'etale de~$U_1$ \eqref{I.7.4}; il est
mod\'er\'ement ramifi\'e sur~$X_1$ relativement \`a~$S$ car il
en est ainsi sur les fibres g\'eom\'etriques en chaque point de~$S$
\eqref{XIII.5.1}. Enfin le fait que $U'_1$ soit connexe r\'esulte de
(SGA~4 XVI~3.2).

\begin{proposition}[Lemme d'Abhyankar relatif]
\label{XIII.5.5}
\index{Abhyankar relatif (lemme d')|hyperpage}%
Soient $X$ un $S$-sch\'ema, $D$ un
diviseur \`a croisements normaux relativement \`a~$S$, comme
dans~\eqref{XIII.5.4}. Soient $Y=X-\Supp D$,
$U=X- Y$, $x$ un point de~$Y$, $X_{1}$ le localis\'e
strict de~$X$ en un point g\'eom\'etrique au-dessus de~$x$,
$U_{1}=U\times_{X}X_{1}$, $V_{1}$ un rev\^etement \'etale de
$U_{1}$. On suppose que, pour tout point maximal $s$ de~$S$,
$V_{1\overline{s}}$ est mod\'er\'ement ramifi\'e sur~$X_{1\overline{s}}$ relativement \`a $\overline{s}$. Alors on peut
trouver des entiers $n_{i}$ premiers \`a la caract\'eristique $p$
de~$\kres(x)$, avec $i\in I(x)$, tels que, si l'on pose
$$
X'_{1}=X_{1}[T_{i}]_{i\in I(x)}\big/ (T_{i}^{n_{i}}-f_{i}) \quoi,
$$
$U'_{1}=U_{1}\times_{X_{1}}X'_{1}$, etc., le rev\^etement \'etale
$V'_{1}$ de~$U'_{1}$ se prolonge de mani\`ere unique \`a
isomorphisme unique pr\`es en un rev\^etement \'etale de
$X'_{1}$. En particulier
\marginpar{432}
$V_{1}$ est mod\'er\'ement ramifi\'e sur~$X_{1}$ relativement
\`a~$S$.
\end{proposition}

On peut supposer $S$ local noeth\'erien de point ferm\'e~$f(x)$.
Pour chaque point maximal $s$ de~$S$ et pour chaque $i\in I(x)$, soit
$x_{i}$ le point g\'en\'erique du ferm\'e $V(f_{i})$ de la fibre
$X_{1\overline{s}}$. L'anneau local $(\cal{O}_{X_{1},x_{i}})_{\red}$
est un anneau de valuation discr\`ete de corps de fractions $K_{i}$
et l'on a
\ifthenelse{\boolean{orig}}
{$V_{1}|K_{i}=\Spec(\prod_{j\in J(x_{i})}L_{j})$,}
{$V_{1}|K_{i}=\Spec(\prod_{j\in I(x_{i})}L_{j})$,}
o\`u $L_{j}$ est une extension finie s\'eparable de~$K$; on note
$n_{j}$ l'ordre du groupe d'inertie d'une extension galoisienne
engendr\'ee par $L_{j}$ et $n_{i}$ le p.p.c.m. des $n_{j}$ quand
$s$ parcourt les points maximaux de~$S$ et
\ifthenelse{\boolean{orig}}
{$j\in J(x_{i})$.}
{$j\in I(x_{i})$.}

Les $n_{i}$ \'etant ainsi choisis, nous allons montrer que $V'_{1}$
se prolonge de fa\c con unique en un rev\^etement \'etale de
$X'_{1}$. L'unicit\'e r\'esulte du fait que, $X'$ \'etant lisse
sur~$S$ aux points de~$Y$, on a $\prof\et_{Y'_{1}}(X'_{1})\geqslant 2$
(SGA~4 XVI~3.2 ou SGA~2 XIV~1.19). Soient $x'_{1}$ un point de
$Y'_{1}$, $\overline{x}'_{1}$ un point g\'eom\'etrique au-dessus
de~$x'_{1}$, et notons $\overline{X'_{1}}$ le localis\'e strict de
$X'_{1}$ en $\overline{x}'_{1}$,
$\overline{U'_{1}}=U'_{1(\overline{X'_{1}})}$,~etc. Par descente,
compte tenu de l'unicit\'e, il suffit de montrer que
$\overline{V}'_{1}$ se prolonge \`a $\overline{X'_{1}}$. De plus on
peut se borner \`a prendre pour $x'_{1}$ les point maximaux de
$Y'_{1}$; en effet on aura alors un prolongement de~$V'_{1}$ sur un
ouvert $W'_{1}$ de~$X'_{1}$ contenant les points maximaux de~$Y'_{1}$;
or, si $Z'_{1}=X'_{1}- W'_{1}$, on a
$\codim(Z'_{1s},X'_{1s})\geqslant 2$ si $s$ est un point maximal de
$S$ et $\codim(Z'_{1s},X'_{1s})\geqslant 1$, $\prof\et_{s}(S)\geqslant
1$ si $s$ est un point de~$S$ qui n'est pas maximal; le fait que
$V'_{1}$ se prolonge \`a $X'_{1}$ tout entier r\'esulte alors
de~(SGA~2 XIV~1.20). Mais, en un point g\'eom\'etrique
$\overline{x}'_{1}$ au-dessus d'un point maximal de~$Y'_{1}$,
$X'_{1\red}$ est le spectre d'un anneau de valuation discr\`ete, et
le fait que $\overline{V'_{1}}$ se prolonge \`a $\overline{X'_{1}}$
r\'esulte alors de~X~\Ref{X.3.6}.

Montrons que les $n_{i}$ sont premiers \`a $p$. En effet, s'il n'en
\'etait pas ainsi, on aurait un indice $i_{0}\in I(x)$ tel que $p$
divise $n_{i_{0}}$. Quitte \`a remplacer $X$ par
$X_{1}[T_{i_{0}},T_{i}]_{i\in I(x)}\big/
(T_{i}^{n_{i_{0}}/p}-f_{i_{0}},T_{i}^{n_{i}}-f_{i})$, on se ram\`ene
au cas o\`u $X'_{1}=X_{1}[T]/T^{p}-f_{i_{0}}$. D'apr\`es ce qui
pr\'ec\`ede le rev\^etement
\marginpar{433}
\'etale $V'_{1}$ de~$U'_{1}$ se prolonge en un rev\^etement
\'etale $E'_{1}$ de~$X'_{1}$. Soit $\eta$ le point ferm\'e de
$S$; comme le morphisme $X'_{1\eta}\to X_{1\eta}$ est radiciel,
$V_{1\eta}$ se prolonge en un rev\^etement \'etale $E_{1\eta}$ de
$X_{1\eta}$. On en d\'eduit alors, comme dans~\Ref{XIII.5.2}, que
$E'_{1}$ est muni d'une donn\'ee de descente relativement au
morphisme $X'_{1}\to X_{1}$, qui prolonge la donn\'ee de descente
naturelle que l'on a sur~$E'_{1}|U'_{1}$. Il en r\'esulte que
$V_{1}$ se prolonge \`a $X_{1}$; mais ceci entra\^ine que l'on a
$n_{i_{0}}=1$ contrairement \`a l'hypoth\`ese $n_{i_{0}}=p$.

\begin{corollaire}
\label{XIII.5.6}
Soient $X$ un $S$-sch\'ema, $D=\sum_{1\leqslant i\leqslant
r}\divisor f_{i}$ un diviseur \`a croisements normaux relativement
\`a~$S$, comme dans~\Ref{XIII.5.4}. Soient $\overline{x}$ un point
g\'eom\'etrique de~$X$, $\overline{X}$ le localis\'e strict de
$X$ en $\overline{x}$, $\overline{Y}=Y_{(\overline{X})}$,
$\overline{U}=\overline{X}-\overline{Y}$ et
$$
\widetilde{U}= \varprojlim_{(n_{i})}
\overline{U}[T_{i}]_{i\in I(x)}\big/ (T_{i}^{n_{i}}-f_{i}) \quoi,
$$
la limite projective \'etant prise suivant l'ensemble filtrant des
familles d'entiers $n_{i}>0$, premiers \`a la caract\'eristique
$p$ de~$\kres(x)$. Alors $\widetilde{U}$ est un rev\^etement universel
mod\'er\'ement ramifi\'e de~$\overline{U}$ relativement \`a~$S$. Par suite le groupe fondamental mod\'er\'ement ramifi\'e
de~$\overline{U}$ est
$$
\pi_{1}^{\tame}(\overline{U})\simeq \prod_{\ell\neq
p}\ZZ_{\ell}[1]^{I(x)} \qquad \textrm{(isomorphisme canonique)\quoi.}
$$
Le groupe $\pi_{1}^{\tame}(\overline{U})$ est isomorphe non
canoniquement \`a $\prod_{\ell\neq p}\ZZ_{\ell}^{I(x)}$.
\end{corollaire}

\begin{subremarque}
\label{XIII.5.6.1}
Soient $X$ un $S$-sch\'ema, $D=\sum_{1\leqslant i\leqslant
r}\divisor f_{i}$ un diviseur \`a croisements normaux relativement
\`a~$S$, comme dans~\Ref{XIII.5.4}, $U=X-\Supp D$.
Pour toute partie $I\subset[1,r]$ soit
$$
X_{I}= \Big(\tbigcap_{i\in I}V(f_{i})\Big) \cap \Big(\tbigcap_{i\in\complement
I}X_{f_{i}}\Big)\quoi .
$$
Soit $p$ un entier premier ou nul et soit $Z$ un sous-ensemble de
$X_{I}$ dont tous les points sont de caract\'eristique~$p$. Soit
$$
\widetilde{U}_{I}= \varprojlim_{(n_{i})}
\ifthenelse{\boolean{orig}}
{U[T_{i}]}
{U[T_{i}]_{i\in I}}
\big/ (T_{i}^{n_{i}}-f_{i}) \quoi,
$$
o\`u
\marginpar{434}
la limite projective est prise suivant l'ensemble filtrant des
familles d'entiers $n_{i}>0$, premiers \`a $p$. Alors, pour tout
point g\'eom\'etrique $\overline{x}$ de~$Z$, l'image inverse de
$\widetilde{U}_{I}$ sur~$\overline{U}$ s'identifie au rev\^etement
universel mod\'er\'ement ramifi\'e de~$\overline{U}$.
\end{subremarque}

\begin{corollaire}
\label{XIII.5.7}
Les notations sont celles de~\Ref{XIII.5.6}. Soient $\overline{S}$ le
localis\'e strict de~$S$ en~$\overline{x}$,
$$
\overline{g}\colon\overline{U}\to\overline{S} \quad \text{et} \quad
\tilde{g}\colon\widetilde{U}\to\widetilde{S}
$$
les morphismes canoniques. Alors les morphismes $\overline{g}$ et
$\tilde g$ sont $0$-acycliques \textup{(SGA~4 XV~1.3)}. Soient $G$ un faisceau
en groupes constructible sur~$\overline{S}$, $F=\overline{g}^*G$, $P$
un torseur sous $F$. Alors, pour que $P$ soit mod\'er\'ement
ramifi\'e sur~$\overline{X}$ relativement \`a $\overline{S}$, il
faut et il suffit que son image inverse $\widetilde{P}$ sur~$\widetilde{U}$ soit triviale.
\end{corollaire}

En effet, pour tout sch\'ema
$\overline{X'}=\overline{X}[T_{i}]_{i\in I(x)}\big/
(T_{i}^{n_{i}}-f_{i})$, o\`u les $n_{i}$ sont des entiers~$>0$
premiers \`a $p$, le morphisme
$\overline{f}'\colon\overline{X'}\to\overline{S}$ est $0$-acyclique.
Les fibres g\'eom\'etriques de~$\overline{f}'$ aux diff\'erents
points de~$\overline{S}$ sont donc connexes et m\^eme
irr\'eductibles. Il en est donc de m\^eme des fibres
g\'eom\'etriques des morphismes
$\overline{g}'\colon\overline{U'}\to\overline{S}$, ce qui prouve que
les $\overline{g}'$, donc aussi $\tilde g$, sont
$0$-acycliques (SGA~4 XV~1.16).

Il est clair qu'un torseur $P$ sur~$\overline{U}$ de groupe $F$, dont
l'image inverse sur~$\widetilde{U}$ est triviale, est
mod\'er\'ement ramifi\'e sur~$\overline{X}$ relativement \`a
$\overline{S}$. Montrons que r\'eciproquement, si $P$ est
mod\'er\'ement ramifi\'e sur~$\overline{X}$ relativement \`a
$\overline{S}$, son image inverse sur~$\widetilde{U}$ est triviale.

Il r\'esulte de (SGA~4 IX~2.14~(ii)) que l'on peut trouver un
morphisme fini $n\colon S_{1}\to\overline{S}$, un faisceau en groupes
constant $C$ sur~$S_{1}$, un monomorphisme $G\to n_*C$.
Consid\'erons le diagramme commutatif suivant form\'e de
carr\'es cart\'esiens:
$$
\xymatrix{ \widetilde{U}_{1} \ar[r] \ar[d]_{r} &U_{1} \ar[r]^{g_{1}}
\ar[d]_{q} &S_{1} \ar[d]_{n} \cr \widetilde{U} \ar[r] &\overline{U}
\ar[r]^{\overline{g}} &\overline{S} \rlap{\quoi.}\cr }
$$
Soient
\marginpar{435}
$C_{1}$ (\resp $\widetilde{C}_{1}$) l'image inverse de~$C$ sur~$U_{1}$
(\resp sur~$\widetilde{U}_{1}$). On a un diagramme commutatif, dans
lequel $i$ et $j$ sont des isomorphismes (SGA~4 VIII~5.8):
\begin{equation*}
\label{eq:XIII.5.7.*}
\tag{$*$}
\begin{array}{c}
\xymatrix{ \H^{1}(\overline{U},q_*C_{1}) \ar[r]^{i}
\ar[d] &\H^{1}(U_{1},C_{1}) \ar[d] \cr
\H^{1}(\widetilde{U},r_*\widetilde{C}_{1}) \ar[r]^{j}
&\H^{1}(\widetilde{U}_{1},\widetilde{C}_{1}) \rlap{\quoi.}\cr }
\end{array}
\end{equation*}
Soit $Q$ le torseur sous $q_*C_{1}$ d\'eduit de~$P$ par
l'extension du groupe structural $F\to q_*C_{1}$.
D'apr\`es~\Ref{XIII.2.1.4}, $Q$ est mod\'er\'ement ramifi\'e
sur~$\overline{X}$ relativement \`a $\overline{S}$. Au torseur $Q$
correspond, gr\^ace \`a $i$, un torseur $Q_{1}$ sous $C_{1}$, et
il est clair que $Q_{1}$ est mod\'er\'ement ramifi\'e sur~$X_{1}=\overline{X}\times_{\overline{S}}S_{1}$ relativement \`a~$S_{1}$. Il r\'esulte donc de~\Ref{XIII.5.6} que l'image inverse
$\widetilde{Q}_{1}$ de~$Q_{1}$ sur~$\widetilde{U}_{1}$ est triviale,
et le diagramme~\eqref{eq:XIII.5.7.*} montre alors que l'image inverse
$\widetilde{Q}$ de~$Q$ sur~$\widetilde{U}$ est triviale.

Consid\'erons le diagramme commutatif suivant, dont la deuxi\`eme
ligne est exacte (SGA~4 XII~3.1):
\begin{equation*}
\label{eq:XIII.5.7.**}
\tag{$**$}
\begin{array}{c}
\xymatrix{ \H^{0}(\overline{S},n_*C/G) \ar[r]
\ar[d]_{k} &\H^{1}(\overline{S},G)\rlap{$~=~1$} \ar[d] \cr
\H^{0}(\widetilde{U},r_*\widetilde{C}_{1}/\widetilde{F}) \ar[r]
&\H^{1}(\widetilde{U},\widetilde{F}) \ar[r]
&\H^{1}(\widetilde{U},r_*\widetilde{C}_{1}) \rlap{\quoi.} \cr }
\end{array}
\end{equation*}
Comme le morphisme $\widetilde{U}\to\overline{S}$ est $0$-acyclique,
$k$ est un isomorphisme. Le fait que $\widetilde{P}$ soit trivial
r\'esulte alors de~\eqref{eq:XIII.5.7.**}.

\section[Appendice II: finitude pour les images
directes des champs]{Appendice II: th\'eor\`eme de finitude pour les images
directes des champs}
\label{XIII.6}

\begin{proposition}
\label{XIII.6.1}
Soient $S$ un sch\'ema localement noeth\'erien, $f\colon X\to S$
un morphisme. Si $S'$ est un $S$-sch\'ema, on note $X'$
(\resp $f'$, etc.) l'image inverse de~$X$ (\resp $f$, etc.) par le
morphisme $S'\to S$. Supposons que, pour tout
\marginpar{436}
sch\'ema $S'$ \'etale sur~$S$, pour tout faisceau d'ensembles
constructible $F$ sur~$X'$, $f'_*F$ soit constructible, et que, pour
tout faisceau en groupes constructible $F$ sur~$X'$, $\R^{1}f'_*F$
soit constructible. Soit $\Phi$ un champ $1$-constructible sur~$X$~\eqref{XIII.0}. Alors $f_*\Phi$ est $1$-constructible.
\end{proposition}

Pour tout sch\'ema $S'$ \'etale sur~$S$ et pour tout objet $x$ de
$(f_*\Phi)_{S'}$, on a un isomorphisme
$$
\SheafAut_{S'}(x)\simeq f'_*\SheafAut_{X'}(x) \quoi,
$$
o\`u, dans le deuxi\`eme membre de l'\'egalit\'e, $x$ est
consid\'er\'e comme objet de~$\Phi_{X'}$. Les hypoth\`eses
faites entra\^inent donc que $f_*\Phi$ est constructible.
Soit $S\Phi$ le faisceau des sous-gerbes maximales de
$\Phi$~\cite[III~2.1.7]{XIII.1}. Comme $f_*(S\Phi)$ est constructible, on
peut lui appliquer (SGA~4 IX~2.7), et le fait que $f_*\Phi$ soit
$1$-constructible r\'esulte alors du lemme qui suit.

\begin{sublemme}
\label{XIII.6.1.1}
Soient $S$ un sch\'ema localement noeth\'erien, $f\colon X\to S$
un morphisme, $\Phi$ un champ sur~$X$. On suppose donn\'e un faisceau
sur~$S$, repr\'esentable par un $S$-sch\'ema \'etale de type
fini $T$, un morphisme surjectif
$$
a\colon T\to f_*(S\Phi )
$$
et un objet $p$ de la fibre $\Phi_{X_T}$ (o\`u $X_T=X\times_ST$),
d\'efinissant dans $f_*(S\Phi )(T)=\allowbreak S\Phi (X_T)$ un \'el\'ement
\'egal \`a l'image $q$ par $a$ de la section identique de
$T(T)$. Soit $f_T\colon X_T\to T$ le morphisme canonique et supposons
que le faisceau $\R^1f_{T*}(\SheafAut_{X_T}(p))$ soit constructible;
alors il en est de m\^eme de~$S(f_*\Phi )$.
\end{sublemme}

Le morphisme canonique $f^*f_*\Phi\to\Phi$ donne un morphisme
$$
S(f^*f_*\Phi )\simeq f^*(S(f_*\Phi ))\to S\Phi ,
$$
d'o\`u un morphisme canonique
$$
\varphi\colon S(f_*\Phi )\to f_*(S\Phi ).
$$
Soient
\marginpar{437}
$F=S(f_*\Phi )$ et $G$ l'image de~$F$ par $\varphi$; d'apr\`es
(SGA~4 IX,~2.9) $G$ est un faisceau constructible.

Il suffit de montrer que, pour tout point $s$ de~$S$, il existe un
ouvert non vide~$U$ de~$s$ tel que $F|U$ soit localement constant
constructible. Soient $s\in S$, $\overline{s}$ un point
g\'eom\'etrique au-dessus de~$s$, $\overline{q}_1,\dots,\overline{q}_n$ les \'el\'ements de~$G_{\overline{s}}$. Par
d\'efinition de~$T$, il existe des $S$-morphismes
$h_i\colon\overline{s}\to S'$ tels que l'on ait
$\overline{q}_i=h_i^*(q)$. Soient $S'$ le produit fibr\'e sur~$S$ de
$n$ sch\'emas isomorphes \`a $T$, $\overline{s}\to S'$ le produit
fibr\'e des $h_i$, $X'=X\times_SS'$, $q_i$ (\resp $p_i$) l'image
inverse de~$q$ (\resp $p$) par la $i$-\`eme projection de~$S'$ sur~$T$. Si $\Psi_i$ est la sous-gerbe maximale de~$\Phi|X'$
engendr\'ee par $p_i$, le faisceau
$F_i=\R^1f_*'(\SheafAut_{X'}(p_i))$ n'est autre que le faisceau
$S(f_*'\Psi_i)$ des sous-gerbes maximales de~$f_*'\Psi_i$. En
particulier l'injection canonique $\Psi_i\to\Phi |X'$ donne un
morphisme
$$
\alpha_i\colon F_i\to F|S'.
$$

Nous allons montrer que $\alpha_i$ est une bijection de~$F_i$ sur
l'image inverse de~$q_i$ dans $F|S'$. Pour tout sch\'ema $S''$
\'etale sur~$S'$ toute section $y$ de~$F_i(S'')$ a pour image
$q_i|S''$ dans $F(S'')$, car, localement pour la topologie \'etale
sur~$S''$, $y$ est d\'efini par un objet~$x$ de~$\Phi_{X''}$ qui est
isomorphe \`a $p_i|X''$. R\'eciproquement, si $y\in F(S'')$ a pour
image $q_i|S''$ dans $F(S'')$, localement pour la topologie \'etale
sur~$S''$, $y$ est d\'efini par un objet~$x$ de~$\Phi_{X''}$ qui est
isomorphe \`a $p_i$; par suite $x$ est un objet de~$\Psi_{iX''}$ et
par suite $y\in F_i(S'')$.

La d\'emonstration s'ach\`eve en utilisant~\Ref{XIII.6.1.2}
ci-dessous. On peut en effet trouver un voisinage ouvert $U'$ de~$s$
tel que $q_1|U',\dots,q_n|U'$ soient des sections de~$G(U')$ et
engendrent ce faisceau. Comme les $F_i|U'$ et $G|U'$ sont
constructibles, il en est de m\^eme de~$F|U'$
d'apr\`es~\Ref{XIII.6.1.2}; quitte \`a remplacer $U$ par un ouvert
plus petit, $F|U$ est localement constant, ce qui ach\`eve la
d\'emonstration.

\begin{sublemme}
\label{XIII.6.1.2}
Soient
\marginpar{438}
$S$ un sch\'ema localement noeth\'erien, $F\to G$ un
morphisme surjectif de faisceaux en groupes sur~$S$. Soient $q_i$ une
famille finie de sections de~$G$ sur~$X$ qui engendrent $G$, et, pour
chaque $i$, soit $F_i$ le sous-faisceau de~$F$ image inverse de
$q_i$. Alors, si $G$ et les $F_i$ sont constructibles, il en est de
m\^eme de~$F$.
\end{sublemme}

Pour prouver que $F$ est constructible, il suffit de montrer que, pour
tout point~$s$ de~$S$, il existe un voisinage ouvert $U$ de~$s$ tel
que $F|U$ soit localement constant constructible. Soit donc $s$ un
point de~$S$. Comme les faisceaux $F_i$ et $G$ sont constructibles, on
peut trouver un voisinage ouvert $U$ de~$s$ tel que $F_i|U$ et $G|U$
soient localement constants. Montrons alors que $F|U$ est localement
constant. D'apr\`es (SGA~4 IX~2.13~(i)), il suffit de voir que, si
$\overline{s}$ est un point g\'eom\'etrique au-dessus de~$s$,
$\tilde{s}$ un point g\'eom\'etrique de~$U$ et
$\overline{s}\to\tilde{s}$ un morphisme de sp\'ecialisation, le
morphisme canonique
$$
F_{\tilde{s}}\to F_{\overline{s}}
$$
est bijectif.

On consid\`ere les diagrammes commutatifs
$$
\xymatrix{ (F_i)_{\tilde{s}} \ar[r]^-{\tilde{a}_i} \ar[d]^{\wr} &
F_{\tilde{s}} \ar[r]^-{\tilde{a}} \ar[d]_b & G_{\tilde{s}}
\ar[d]^{\wr}\\ (F_i)_{\overline{s}} \ar[r]^-{\overline{a}_i} &
F_{\overline{s}} \ar[r]^-{\overline{a}} & G_{\overline{s}} }
$$
Soit $\overline{q}_i$ (\resp $\tilde{q}_i$) l'image inverse de~$q_i$
dans $G_{\overline{s}}$ (\resp $G_{\tilde{s}}$); les morphismes
$\overline{a}$ et $\tilde{a}$ sont surjectifs, et le morphisme
$\overline{a}_i$ (\resp $\tilde{a}_i$) induit une bijection de
$(F_i)_{\overline{s}}$ (\resp $(F_i)_{\tilde{s}}$) sur~$\overline{a}^{-1}(q_i)$ (\resp sur~$\tilde{a}^{-1}(q_i)$). Il r\'esulte
donc du diagramme ci-dessus que $b$ est un isomorphisme.

\begin{corollaire}
\label{XIII.6.2}
Soient
\marginpar{439}
$S$ un sch\'ema localement noeth\'erien, $f\colon X\to S$
un morphisme propre. Soit $\Phi$ un champ $1$-constructible sur~$X$,
alors $f_*\Phi$ est un champ \hbox{$1$-constructible}.
\end{corollaire}

La d\'emonstration de~\Ref{XIII.6.1} prouve aussi le r\'esultat
suivant, compte tenu de \Ref{XIII.2.4}~2).

\begin{corollaire}
\label{XIII.6.3}
Soient $S$ un sch\'ema localement noeth\'erien, $f\colon X\to S$ un
morphisme, $D$ un diviseur sur~$X$ \`a croisements normaux
relativement \`a~$S$ \eqref{XIII.2.1}, $Y=\Supp D$, $U=X-Y$,
$i\colon U\to X$ l'immersion canonique. Soit $\Phi$ un champ sur~$U$
donn\'e, localement pour la topologie \'etale sur~$X$ et $S$,
comme image inverse d'un champ $\Psi$ $1$-constructible sur~$S$. Alors
le champ $i_*^{\tame}\Phi$ est $1$-constructible.
\end{corollaire}


\backmatter
\def\indexname{Index terminologique}
\printindex

\def\indexname{Index des notations}

\begin{theindex}

\medskip\item $\Delta_{X/Y}$ ou simplement $\Delta $\dotfill\pageref{indnot:ab}
\medskip\item $\mathit {\Omega }^1_{X/Y}$\dotfill\pageref{indnot:ac}
\medskip\item $\cal {P}_{X/Y}^n$\dotfill\pageref{indnot:ad}, \pageref{indnot:ad1}
\medskip\item $\Delta_{X/Y}^n$\dotfill\pageref{indnot:ae}
\medskip\item $\goth {d}_{X/Y}$\dotfill\pageref{indnot:af}
\medskip\item $d_{X/S}^n$\dotfill\pageref{indnot:bb}
\medskip\item $\goth {g}_{X/S}$\dotfill\pageref{indnot:cb}
\medskip\item $\cal {C}(\pi )$\dotfill\pageref{indnot:eb}
\medskip\item $\pi $\dotfill\pageref{indnot:ec}
\medskip\item $\Prodash \cal {C}(\pi)$\dotfill\pageref{indnot:ed}
\medskip\item $\Gamma $\dotfill\pageref{indnot:ef}
\medskip\item $\pi_1(S,a)$\dotfill\pageref{indnot:eg}
\medskip\item $\pi_1(S;a,a')$\dotfill\pageref{indnot:eh}
\medskip\item $\cal {C}(S)$\dotfill\pageref{indnot:ei}
\medskip\item $\pi_1(f;a')$\dotfill\pageref{indnot:ej}
\medskip\item $\Ens $\dotfill\pageref{indnot:fb}
\medskip\item $\Cat $\dotfill\pageref{indnot:fc}
\medskip\item $\Ob (\cal {C})$\dotfill\pageref{indnot:fd}
\medskip\item $\Fl (\cal {C})$\dotfill\pageref{indnot:fe}
\medskip\item $\SheafHom (\cal {C},\cal {C}')$\dotfill\pageref{indnot:ff}
\medskip\item $\cal {C}^\circ $\dotfill\pageref{indnot:fg}
\medskip\item $\Cat_{/\cal {E}}$\dotfill\pageref{indnot:fh}
\medskip\item $\SheafHom_{\cal {E}/-}(\cal {F},\cal {G})$\dotfill\pageref{indnot:fi}
\medskip\item $v*u$\dotfill\pageref{indnot:fj}
\medskip\item $\cal {F}\times_{\cal {E}}\cal {G}$\dotfill\pageref{indnot:fk}
\medskip\item $\lambda^*\colon \Cat_{/\cal {E}}\to \Cat_{/\cal {E}'}$\dotfill\pageref{indnot:fl}
\medskip\item $\bf {\Gamma }(\cal {G}/\cal {E})$ et $\Gamma (\cal {G}/\cal {E})$\dotfill\pageref{indnot:fm}
\medskip\item $\cal {F}_S$\dotfill\pageref{indnot:fn}
\medskip\item $f^*_{\cal {F}}(\xi )$ ou $f^*(\xi )$\dotfill\pageref{indnot:fo}
\medskip\item $\alpha_f(\xi )$\dotfill\pageref{indnot:fp}
\medskip\item $\SheafHom_\cart (\cal {F},\cal {G})$\dotfill\pageref{indnot:fq}
\medskip\item $\Cat^\cart_{/\cal {E}}$\dotfill\pageref{indnot:fr}
\medskip\item $\varprojLim \cal {F}/\cal {E}$\dotfill\pageref{indnot:fs}
\medskip\item $f_*^{\cal {F}}$ ou~$f_*$\dotfill\pageref{indnot:ft}
\medskip\item $\Sch $\dotfill\pageref{indnot:hb}
\medskip\item $\bbmu_{n\,S}$\dotfill\pageref{indnot:kb}
\medskip\item $X^\an $\dotfill\pageref{indnot:lb}
\medskip\item $f^\an $\dotfill\pageref{indnot:lc}
\medskip\item $SF$ ou $S(F)$\dotfill\pageref{indnot:mb}
\medskip\item $\H^1_{\tame}(U,F)$\dotfill\pageref{indnot:mc}
\medskip\item $\R^1_{\tame} g_* F$\dotfill\pageref{indnot:md}
\medskip\item $C_t((U,X)/S)$ ou $C_t$\dotfill\pageref{indnot:me}
\medskip\item $\pi_1^{\tame}((U,X)/S,a)$ ou $\pi_1^{\tame}(U,a)$ ou $\pi_1^{\tame}(U)$\dotfill\pageref{indnot:mf}
\medskip\item $(g_*^{\tame}\Phi )_{T'}$\dotfill\pageref{indnot:mg}
\medskip\item $\H^0(V,C_V)^{\Pi }$\dotfill\pageref{indnot:mh}
\medskip\item $\pi_1^{\LL }(U,a)$\dotfill\pageref{indnot:mi}
\medskip\item $\pi_1'(X,a)$\dotfill\pageref{indnot:mj}
\medskip\item $\pi_1^\LL (X/S,g,\bar {s})$ ou $\pi_1^\LL (X/S,g)$\dotfill\pageref{indnot:mk}
\medskip\item $\pi_1^\LL (X_{\bar {s}},a)_K$\dotfill\pageref{indnot:ml}
\medskip\item $\ZZ_\ell [1]$\dotfill\pageref{indnot:mm}
\par
\end{theindex}
\end{document}